\documentclass[letterpaper]{article}
\usepackage{exscale,relsize}
\usepackage[dvips]{graphicx,color}
\usepackage{fancybox}

\usepackage{amsmath}
\usepackage{amssymb}
\usepackage{bm}
\usepackage[mathscr]{eucal}
\usepackage{amsfonts}
\usepackage{latexsym}
\usepackage{amsxtra}
\usepackage{amsbsy}
\usepackage{amsthm}
\usepackage{amscd}
\usepackage{amsopn}
\usepackage{amstext}
\usepackage{upref}
\usepackage{oldgerm}
\usepackage[varumlaut]{yfonts}
\newcounter{z0}

\newtheorem{bbbb}{Proposition}[subsection]
\newtheorem{cccc}{Lemma}[subsection]
\newtheorem{dddd}{Theorem}[subsection]
\newtheorem{ffff}{Corollary}[subsection]
\newtheorem{bbbbb}{Proposition}[section]
\newtheorem{ccccc}{Lemma}[section]
\newtheorem{ddddd}{Theorem}[section]

\theoremstyle{definition}

\theoremstyle{definition}

\theoremstyle{definition}

\theoremstyle{definition}

\theoremstyle{definition}
\newtheorem{eeee}{Remark}[subsection]

\theoremstyle{definition}
\newtheorem{eeeee}{Remark}[section]

\theoremstyle{definition}

\newcommand{\me}{\mathrm{e}}
\newcommand{\mi}{\mathrm{i}}
\newcommand{\md}{\mathrm{d}}
\renewcommand{\Im}{\mathrm{Im}}
\renewcommand{\Re}{\mathrm{Re}}
\newcommand{\id}{\pmb{\mi \md}}
\numberwithin{equation}{section}
\allowdisplaybreaks[4]

\newenvironment{lcase}{\left\lbrace \begin{aligned}}{\end{aligned} \right.}
\usepackage{stmaryrd}
\DeclareMathOperator*{\ad}{ad}
\DeclareMathOperator{\tr}{tr}
\DeclareMathOperator{\diag}{diag}
\let \Im \relax
\let \Re \relax
\DeclareMathOperator{\Im}{Im}
\DeclareMathOperator{\Re}{Re}
\usepackage{epic}
\usepackage{eepic}
\begin{document}
\fontsize{10pt}{11pt}\selectfont
\fontencoding{T1}\selectfont
\baselineskip=11pt
\frenchspacing
\flushbottom
\raggedbottom
\title{Trans-Series Asymptotics of Solutions to the Degenerate Painlev\'{e} III Equation: A Case Study}
\author{A.~Vartanian}
\date{\today}
\maketitle
\begin{abstract}
\noindent 
A one-parameter family of trans-series asymptotics as $\tau \! \to \! \pm \infty$ and $\tau \! \to \! \pm \mi \infty$ 
for solutions of the degenerate Painlev\'{e} III equation (DP3E), $u^{\prime \prime}(\tau) \! = \! \frac{(u^{\prime}
(\tau))^{2}}{u(\tau)} \! - \! \frac{u^{\prime}(\tau)}{\tau} \! + \! \frac{1}{\tau}(-8 \varepsilon (u(\tau))^{2} \! + \! 2ab) \! 
+ \!\frac{b^{2}}{u(\tau)}$, where $\varepsilon \! \in \! \lbrace \pm 1 \rbrace$, $a \! \in \! \mathbb{C}$, and $b \! \in 
\! \mathbb{R} \setminus \lbrace 0 \rbrace$, are parametrised in terms of the monodromy data of an associated 
first-order $2 \times 2$ matrix linear ODE via the isomonodromy deformation approach: trans-series asymptotics 
for the associated Hamiltonian and principal auxiliary functions and the solution of one of the $\sigma$-forms of 
the DP3E are also obtained. The actions of various Lie-point symmetries for the DP3E are derived.

\vspace{0.40cm}

\textbf{2020 Mathematics Subject Classification.} 33E17, 34M35, 34M40, 34M50, 34M55, 34M56, 

34M60

\vspace{0.23cm}


\vspace{0.23cm}

\textbf{Key Words.} Isomonodromy deformations, Stokes phenomena, symmetries
\end{abstract}
\section{Introduction} \label{sec1} 
In this section, which is partitioned into five inter-dependent subsections, the reader is given a concise overview of the 
information subsumed in the text: (i) in Subsection~\ref{sec1a}, the degenerate Painlev\'{e} III equation (DP3E) is 
introduced, and the qualitative behaviours of the asymptotic results the reader can expect to excise from this work are 
delineated; (ii) in Subsection~\ref{sec1b}, the DP3E's associated Hamiltonian and principal auxiliary functions, as well 
as one of its $\sigma$-forms, are introduced; (iii) in Subsection~\ref{sec1c}, pre- and post-gauge-transformed Lax pairs 
giving rise to isomonodromy deformations are reviewed; (iv) in Subsection~\ref{sec1d}, canonical asymptotics of the 
post-gauge-transformed Lax-pair solution matrix is presented in conjunction with the corresponding monodromy data; 
and (v) in Subsection~\ref{sec1e}, the monodromy manifold and the direct and inverse problems of monodromy theory 
are introduced, and a synopsis of the organisation of this work is given.
\subsection{The Degenerate Painlev\'{e} III Equation (DP3E)} \label{sec1a} 
This paper continues the studies in \cite{avlkv,a1,av2,avkavint,avkavalgeb,kitavzap2023} of the DP3E,
\begin{equation} \label{eq1.1} 
u^{\prime \prime}(\tau) \! = \! \frac{(u^{\prime}(\tau))^{2}}{u(\tau)} \! - \! \frac{u^{\prime}(\tau)}{\tau} \! + \! \frac{1}{\tau} \! 
\left(-8 \varepsilon (u(\tau))^{2} \! + \! 2 ab \right) \! + \! \frac{b^{2}}{u(\tau)}, \quad \varepsilon \! \in \! \lbrace \pm 1 \rbrace,
\end{equation}
where the prime denotes differentiation with respect to $\tau$, $\mathbb{C} \! \ni \! a$ is the parameter of formal 
monodromy, and $\mathbb{R} \setminus \lbrace 0 \rbrace \! \ni \! b$ is a parameter (see, also, \cite{vgilss}, 
Chapter 7, Section 33); in fact, making the formal change of variables $\tau \! \to \! t^{1/2}$, $u(\tau) \! \to \! 
\tilde{\eta}_{0}^{2}t^{-1/2} \tilde{\lambda}(t)$, $a \! \to \! \mp \mi \tilde{c}_{0} \tilde{\eta}_{0}$, and $b \! \to \! \pm \mi 2 
\tilde{\eta}_{0}^{3}$, where $\tilde{c}_{0} \! \in \! \mathbb{C}$ and $\mi \tilde{\eta}_{0} \! \in \! \mathbb{R} \setminus 
\lbrace 0 \rbrace$, and setting $\varepsilon \! = \! +1$, one shows that the DP3E \eqref{eq1.1} transforms into, in the 
classification scheme of \cite{a3}, the degenerate third Painlev\'{e} equation of type $D_{7}$,
\begin{equation} \label{dpee3d7}
(P_{\mathrm{III}^{\prime}})_{D_{7}}: \quad \quad \frac{\md^{2} \tilde{\lambda}}{\md t^{2}} \! = \! \frac{1}{\tilde{\lambda}} 
\! \left(\frac{\md \tilde{\lambda}}{\md t} \right)^{2} \! - \! \frac{1}{t} \frac{\md \tilde{\lambda}}{\md t} \! + \! \tilde{\eta}_{0}^{2} 
\! \left(-2 \frac{\tilde{\lambda}^{2}}{t^{2}} \! + \! \frac{\tilde{c}_{0}}{t} \! - \! \frac{1}{\tilde{\lambda}} \right).
\end{equation}
It is known that, in the complex plane of the independent variable, Painlev\'{e} equations admit, in open sectors near 
the point at infinity containing one special ray, pole-free solutions that are characterised by divergent asymptotic 
expansions: such solutions, called \emph{tronqu\'{e}e} solutions by Boutroux, usually contain free parameters 
manifesting in exponentially small terms for large values of the independent variable. There also exist pole-free solutions 
that are void of parameters in larger open sectors near the point at infinity containing three special rays: such solutions are 
called \emph{tritronqu\'{e}e} solutions (see, for example, \cite{edelere}, Chapter 3). In contrast to the asymptotic results of 
\cite{a1,av2}, this work entails an analysis of one-parameter families of \emph{trans-series\/} (\cite{edelere}, Chapter 5) 
asymptotic (as $\tau \! \to \! \pm \infty$ and $\tau \! \to \! \pm \mi \infty$) solutions related to the underlying quasi-linear 
Stokes phenomenon associated with the DP3E \eqref{eq1.1}: such solutions are also referred to as instanton-type solutions 
in the physics literature \cite{sgaiakmm} (see, also, \cite{aritaapa,a66,sachokap1,sachokap2}, and Chapter 11 of \cite{a5}); 
in particular, tronqu\'{e}e solutions that are free of poles not only on the real and the imaginary axes of $\tau$, but also in 
open sectors about the point at infinity, are considered.\footnote{The terms trans-series \cite{avsmv4,edga} and tronqu\'{e}e 
are used interchangeably in this work.} The existence of one-parameter tronqu\'{e}e solutions for a scaled version of the 
DP3E \eqref{eq1.1} was proved in \cite{ylddpt} via direct asymptotic analysis. A review of recent manifestations of the DP3E 
\eqref{eq1.1} and $(P_{\mathrm{III}^{\prime}})_{D_{7}}$ \eqref{dpee3d7} in variegated mathematical and physical settings 
such as, for example, non-linear optics, number theory, asymptotics, non-linear waves, random matrix theory, and differential 
geometry, is presented in Appendix \ref{litsurvdp3e}. 

An effectual approach for studying the asymptotic behaviour of solutions (in particular, the connection formulae for their 
asymptotics) of the Painlev\'{e} equations $\mathrm{P} \mathrm{I},\mathrm{P} \mathrm{II},\dotsc,\mathrm{P} \mathrm{VI}$ 
is the Isomonodromic Deformation Method (IDM) \cite{a5,a18,a2,a22,sachokap3}: specific features of the IDM as applied, 
in particular, to the DP3E~\eqref{eq1.1} can be located in Sections~1 and~2 of \cite{a1}. It is imperative, within the IDM 
framework, to mention the seminal r\^{o}le played by the recent monograph \cite{a5}, as it summarizes and reflects not only 
the key technical and theoretical developments and advances of the IDM since the appearance of \cite{a2}, but also of 
an equivalent, technically distinct approach based on the Deift-Zhou non-linear steepest descent analysis of the associated 
RHP \cite{a6}. The methodological paradigm adopted in this paper is the IDM. Even though the DP3E \eqref{eq1.1} 
resembles one of the canonical variants of the Painlev\'{e} equations $\mathrm{P} \mathrm{I},\mathrm{P} \mathrm{II},
\dotsc,\mathrm{P} \mathrm{VI}$, the associated asymptotic analysis of its solutions via the IDM subsumes additional 
technical complications due to the necessity of having to extract the explicit functional dependencies of the contributing 
error terms rather than merely estimating them, which requires a considerably more detailed study of the error functions. 
By studying the isomonodromic deformations of a first-order $3 \times 3$ matrix linear ODE (see, also, Section 8 of 
\cite{DP}) with two irregular singular points, asymptotics as $\tau \! \to \! \infty$ and $\tau \! \to \! 0$ of solutions to the 
DP3E \eqref{eq1.1} for the case $a \! = \! 0$, as well as the corresponding connection formulae, were obtained in 
\cite{aveekit} via the IDM. As observed in \cite{avaev}, though, there is an alternative first-order $2 \times 2$ matrix linear 
ODE whose isomonodromy deformations are described, for arbitrary $a \! \in \! \mathbb{C}$, by the DP3E \eqref{eq1.1}: 
it is this latter $2 \times 2$ ODE system that is adopted in this work.

In order to eschew a flood of superfluous notation and to motivate, in as succinct a manner as possible, the qualitative 
behaviour of the solution of the DP3E \eqref{eq1.1} that the reader will encounter in this work, consider, for example, 
asymptotics as $\tau \! \to \! +\infty$ with $\varepsilon b \! > \! 0$ of $u(\tau)$. As is well known \cite{afvkav,avsmv3,
edelere,a5,mlaud,avsmv2,shimo2,shimo,shimo3,avsmv1,xxia}, the Painlev\'{e} equations admit a one-parameter family of 
trans-series solutions of the form ``(power series) $+$ (exponentially small terms)''. As noted in Remark \ref{remlndaip34} 
below, $u(\tau)$ admits the `complete' asymptotic trans-series representation $u(\tau) \genfrac{}{}{0pt}{3}{=}{\tau \to 
+\infty} c_{0,k}(\tau^{1/3} \! + \! v_{0,k}(\tau))$, $k \! \in \! \lbrace \pm 1 \rbrace$, where $c_{0,k} \! := \! \tfrac{1}{2} 
\varepsilon (\varepsilon b)^{2/3} \me^{-\mi 2 \pi k/3}$, and $v_{0,k}(\tau) \! := \! \tau^{-1/3} 
\mathrm{u}_{{\scriptscriptstyle \mathrm{R},k}}(\tau) \! + \! \mathrm{u}_{{\scriptscriptstyle \mathrm{E},k}}(\tau)$, 
with $\mathbb{C} \llbracket \tau^{-1/3} \rrbracket \! \ni \! \mathrm{u}_{{\scriptscriptstyle \mathrm{R},k}}(\tau) 
\! = \! \sum_{n=0}^{\infty} \upsilon_{n,k}({\scriptscriptstyle \mathscr{M}})(\tau^{-1/3})^{n}$ and 
$\mathrm{u}_{{\scriptscriptstyle \mathrm{E},k}}(\tau) \! = \! \sum_{m=1}^{\infty} \sum_{j=0}^{\infty} 
\mathfrak{v}_{m,j,k}({\scriptscriptstyle \mathscr{M}})(\tau^{-1/3})^{j} \big(\me^{-\frac{3 \sqrt{\smash[b]{3}}}{2}
(\sqrt{\smash[b]{3}}+ \mi k)(\varepsilon b)^{1/3} \tau^{2/3}} \big)^{m}$, and where the monodromy-data-dependent expansion 
coefficients $\upsilon_{n,k}({\scriptscriptstyle \mathscr{M}})$ and $\mathfrak{v}_{m,j,k}({\scriptscriptstyle \mathscr{M}})$ can 
be determined recursively provided that certain leading coefficients are known \emph{a priori}. The purpose of this work, 
though, is not to address the complete asymptotic trans-series representation stated above, but, rather, to determine the 
coefficient of the leading-order exponentially small correction term to the asymptotics of solutions of the DP3E \eqref{eq1.1}, 
which is, to the best of the author's knowledge as at the time of the presents, the decidedly non-trivial task within the IDM 
paradigm, in which case, the asymptotic trans-series representation for $u(\tau)$ reads:\footnote{The notation $\lambda_{1}(t) 
\genfrac{}{}{0pt}{3}{=}{t \to +\infty} \mathcal{O}(\lambda_{2}(t))$ means that there exists $\widehat{C} \! > \! 0$ and sufficiently 
small $\hat{\epsilon} \! > \! 0$ such that $\lvert \lambda_{1}(t)/\lambda_{2}(t) \rvert \! \leqslant \! \widehat{C}$ for all 
$t \! > \! 1/\hat{\epsilon}$.}
\begin{equation} \label{mainyoo} 
u(\tau) \underset{\tau \to +\infty}{=} c_{0,k} \! \left(\tau^{1/3} \! + \! \sum_{m=0}^{\infty} 
\frac{\mathfrak{u}_{m}(k)}{(\tau^{1/3})^{m+1}} \! + \! \mathrm{A}_{k} \me^{-\frac{3 \sqrt{\smash[b]{3}}}{2}(\sqrt{\smash[b]{3}}
+\mi k)(\varepsilon b)^{1/3} \tau^{2/3}} \big(1 \! + \! \mathcal{O}(\tau^{-1/3}) \big) \right), \quad k \! \in \! \lbrace \pm 1 \rbrace.
\end{equation}
While the expansion coefficients $\lbrace \mathfrak{u}_{m}(k) \rbrace_{m=0}^{\infty}$, $k \! \in \! \lbrace \pm 1 
\rbrace$, can be determined (not always uniquely!) by substituting the trans-series representation~\eqref{mainyoo} 
into the DP3E \eqref{eq1.1} and solving a system of recurrence relations for the $\mathfrak{u}_{m}(k)$'s, the 
monodromy-data-dependent expansion coefficients, $\mathrm{A}_{k}$, $k \! \in \! \lbrace \pm 1 \rbrace$, can not, 
and must, therefore, be determined independently; in fact, the principal technical accomplishment of this work is the 
determination, via the IDM, of the explicit dependence of the coefficients $\mathrm{A}_{k}$, $k \! \in \! \lbrace \pm 1 
\rbrace$, on the Stokes multiplier $s_{0}^{0}$ (see, in particular, Section \ref{finalsec}, equations \eqref{geek109} 
and \eqref{geek111}). Even though the motivational discussion above for the introduction of the 
monodromy-data-dependent expansion coefficients $\mathrm{A}_{k}$, $k \! \in \! \lbrace \pm 1 \rbrace$, relies on the 
asymptotics of $u(\tau)$ as $\tau \! \to \! +\infty$ for $\varepsilon b \! > \! 0$, it must be emphasized that, in this work, 
the coefficients $\mathrm{A}_{k}$, $k \! \in \! \lbrace \pm 1 \rbrace$, and their analogues, corresponding to trans-series 
asymptotics of $u(\tau)$, the associated Hamiltonian and principal auxiliary functions, and one of the $\sigma$-forms of 
the DP3E~\eqref{eq1.1} as $\tau \! \to \! +\infty \me^{\mi \pi \varepsilon_{1}}$ for $\varepsilon b \! = \! \lvert \varepsilon 
b \rvert \me^{\mi \pi \varepsilon_{2}}$, $\varepsilon_{1},\varepsilon_{2} \! \in \! \lbrace 0,\pm 1 \rbrace$, and $\tau 
\! \to \! +\infty \me^{\mi \pi \hat{\varepsilon}_{1}/2}$ for $\varepsilon b \! = \! \lvert \varepsilon b \rvert \me^{\mi \pi 
\hat{\varepsilon}_{2}}$, $\hat{\varepsilon}_{1} \! \in \! \lbrace \pm 1 \rbrace$ and $\hat{\varepsilon}_{2} \! \in \! \lbrace 
0,\pm 1 \rbrace$, are obtained (see, in particular, Section \ref{sec2}, Theorems \ref{theor2.1} and \ref{appen}, 
respectively).
\begin{eeee} \label{remlndaip34} 
In the seminal work \cite{ylddpt}, the authors consider, in particular, the existence and uniqueness of tronqu\'{e}e 
solutions of the $\mathrm{P} \mathrm{III}$ equation with parameters $(1,\beta,0,-1)$, denoted by 
$\mathrm{P}^{(\mathrm{ii})}_{\mathrm{III}}$ in equation~(1.5) of \cite{ylddpt}: $v^{\prime \prime}(x) \! = \! 
\tfrac{(v^{\prime}(x))^{2}}{v(x)} \! - \! \tfrac{v^{\prime}(x)}{x} \! + \! \tfrac{1}{x}((v(x))^{2} \! + \! \beta) \! - \! \tfrac{1}{v(x)}$, 
where $\mathbb{C} \! \ni \! \beta$ is arbitrary; $\mathrm{P}^{(\mathrm{ii})}_{\mathrm{III}}$ can be derived {}from the 
DP3E \eqref{eq1.1} via the mapping $\mathscr{S}_{\scriptscriptstyle \varepsilon} \colon (\tau,u(\tau),a,b) \! \to \! 
(\alpha x,\gamma v(x),\tfrac{\beta}{2} \me^{-\mi (2m+1) \pi/2},b)$, $\varepsilon \! = \! \pm 1$, $m \! = \! 0,1$, where 
$\alpha \! := \! 2^{-3/2}b^{-1/2} \me^{\mi (2+ \varepsilon) \pi/4} \me^{\mi (2m^{\prime}+m) \pi/2}$, and $\gamma \! := 
\! -2^{-3/2} \varepsilon b^{1/2} \me^{-\mi (2+ \varepsilon) \pi/4} \me^{-\mi (2m^{\prime}+m) \pi/2}$, $m^{\prime} \! = \! 0,1$. 
In Theorem~2 of \cite{ylddpt}, the authors prove that, in any open sector of angle less than $3 \pi/2$, there exist 
one-parameter solutions of $\mathrm{P}^{(\mathrm{ii})}_{\mathrm{III}}$ with asymptotic expansion $v(x) \! \sim \! 
v_{\scriptscriptstyle f}^{\scriptscriptstyle (m_{1})}(x) \! := \! x^{1/3} \sum_{n=0}^{\infty}a_{n}^{\scriptscriptstyle (m_{1})}
(x^{-2/3})^{n}$ for $S_{\scriptscriptstyle k}^{\scriptscriptstyle (m_{1})} \! \ni \! x \! \to \! \infty$, $m_{1} \! = \! 0,1,2$, 
where the sectors $S_{\scriptscriptstyle k}^{\scriptscriptstyle (m_{1})}$, $k \! = \! 0,1,2,3$, are defined in equation 
(1.10) of \cite{ylddpt}, $a_{0}^{\scriptscriptstyle (m_{1})} \! := \! \exp (\mi 2 \pi m_{1}/3)$, and the ($x$-independent) 
coefficients $a_{n}^{\scriptscriptstyle (m_{1})}$, $n \! \in \! \mathbb{N}$, solve the recursion relations (1.12) of 
\cite{ylddpt}; moreover, the authors prove that, for any branch of $x^{1/3}$, there exists a unique solution of 
$\mathrm{P}^{(\mathrm{ii})}_{\mathrm{III}}$ in $\mathbb{C} \setminus \pmb{\leftthreetimes}$ with asymptotic expansion 
$v_{\scriptscriptstyle f}^{\scriptscriptstyle (m_{1})}(x)$, where $\pmb{\leftthreetimes}$ is an arbitrary branch cut connecting 
the singular points $0$ and $\infty$ (they also address the existence of the exponentially small correction term(s) of 
the tronqu\'{e}e solution of $\mathrm{P}^{(\mathrm{ii})}_{\mathrm{III}}$). This crucial result of \cite{ylddpt}, in 
conjunction with the invertibility of the mapping $\mathscr{S}_{\scriptscriptstyle \varepsilon}$, implies the existence 
and the uniqueness of the asymptotic (as $\tau \! \to \! +\infty$ with $\varepsilon b \! > \! 0$) trans-series representation 
\eqref{mainyoo}. \hfill $\blacksquare$
\end{eeee} 
\subsection{Hamiltonian Structure, Auxiliary Functions, and the $\sigma$-Form} \label{sec1b} 
Herewith follows a brief synopsis of select results {}from \cite{a1} that are relevant for the present work; for complete 
details, see, in particular, Sections 1, 2, and 6 of \cite{a1}, and \cite{avkavint}.

An important formal property of the DP3E \eqref{eq1.1} is its associated Hamiltonian structure; in fact, as shown in 
Proposition 1.3 of \cite{a1}, upon setting
\begin{equation} \label{hamk1} 
\mathcal{H}_{\epsilon_{1}}(\hat{p}(\tau),\hat{q}(\tau);\tau) \! := \! (\hat{p}(\tau) \hat{q}(\tau))^{2} \tau^{-1} \! - \! 2 
\epsilon_{1} \hat{p}(\tau) \hat{q}(\tau)(\mi a \! + \! 1/2) \tau^{-1} \! + \! 4 \varepsilon \hat{q}(\tau) \! + \! \mi b \hat{p}
(\tau) \! + \! \frac{1}{2}(\mi a \! + \! 1/2)^{2} \tau^{-1},
\end{equation} 
where the functions $\hat{p}(\tau)$ and $\hat{q}(\tau)$ are the generalised impulse and co-ordinate, respectively, 
$\epsilon_{1} \! \in \! \lbrace \pm 1 \rbrace$, and $\epsilon_{1}^{2} \! = \! \varepsilon^{2} \! = \! 1$, Hamilton's equations, that is,
\begin{equation} \label{hamk2} 
\hat{p}^{\prime}(\tau) \! = \! -\dfrac{\partial \mathcal{H}_{\epsilon_{1}}(\hat{p}(\tau),\hat{q}(\tau);\tau)}{\partial \hat{q}} 
\, \qquad \, \text{and} \, \qquad \, \hat{q}^{\prime}(\tau) \! = \! \dfrac{\partial \mathcal{H}_{\epsilon_{1}}(\hat{p}(\tau),
\hat{q}(\tau);\tau)}{\partial \hat{p}},
\end{equation}
are equivalent to either one of the degenerate $\mathrm{P} \mathrm{III}$ equations
\begin{gather}
\hat{p}^{\prime \prime}(\tau) \! = \! \dfrac{(\hat{p}^{\prime}(\tau))^{2}}{\hat{p}(\tau)} \! - \! \dfrac{\hat{p}^{\prime}(\tau)}{\tau} 
\! + \! \dfrac{1}{\tau} \! \left(-2 \mi b(\hat{p}(\tau))^{2} \! + \! 8 \varepsilon (\mi a \epsilon_{1} \! + \! (\epsilon_{1} \! - \! 1)/2) 
\right) \! - \! \dfrac{16}{\hat{p}(\tau)}, \label{hamk3} \\
\hat{q}^{\prime \prime}(\tau) \! = \! \dfrac{(\hat{q}^{\prime}(\tau))^{2}}{\hat{q}(\tau)} \! - \! \dfrac{\hat{q}^{\prime}(\tau)}{\tau} 
\! + \! \dfrac{1}{\tau} \! \left(-8 \varepsilon (\hat{q}(\tau))^{2} \! - \! b(2a \epsilon_{1} \! - \! \mi (1 \! + \! \epsilon_{1})) \right) 
\! + \! \dfrac{b^{2}}{\hat{q}(\tau)}; \label{hamk4}
\end{gather}
it was also noted during the proof of the above-mentioned result that the Hamiltonian system \eqref{hamk2} can be rewritten as
\begin{equation} \label{hamk5} 
\hat{p}(\tau) \! = \! \dfrac{\tau (\hat{q}^{\prime}(\tau) \! - \! \mi b)}{2(\hat{q}(\tau))^{2}} \! + \! \dfrac{\mi \epsilon_{1}
(a \! - \! \mi/2)}{\hat{q}(\tau)} \qquad \text{and} \, \qquad \, \hat{q}(\tau) \! = \! -\dfrac{\tau (\hat{p}^{\prime}(\tau) 
\! + \! 4 \varepsilon)}{2(\hat{p}(\tau))^{2}} \! + \! \dfrac{\mi \epsilon_{1}(a \! - \! \mi/2)}{\hat{p}(\tau)}.
\end{equation}
As shown in Section 2 of \cite{a1}, the \emph{Hamiltonian function}, $\mathcal{H}(\tau)$, is defined as follows:
\begin{equation} \label{hamk10} 
\mathcal{H}(\tau) \! := \! \left. \mathcal{H}_{\epsilon_{1}}(\hat{p}(\tau),\hat{q}(\tau);\tau) \vphantom{M^{M}} 
\right\vert_{\epsilon_{1}=-1},
\end{equation}
where $\hat{p}(\tau)$ is calculated {}from the first (left-most) relation of equations \eqref{hamk5} with $\hat{q}(\tau) \! 
= \! u(\tau)$; moreover, as shown in Section 2 of \cite{a1}, the definition \eqref{hamk10} implies the following explicit 
expression for $\mathcal{H}(\tau)$ in terms of $u(\tau)$:
\begin{equation} \label{eqh1} 
\mathcal{H}(\tau) \! = \! (a \! - \! \mi/2) \dfrac{b}{u(\tau)} \! + \! \dfrac{1}{2 \tau}(a \! - \! \mi/2)^{2} \! + \! 
\dfrac{\tau}{4(u(\tau))^{2}} \! \left((u^{\prime}(\tau))^{2} \! + \! b^{2} \right) \! + \! 4 \varepsilon u(\tau).
\end{equation}

It was shown in Section 1 of \cite{a1} that the function $\sigma (\tau)$ defined by
\begin{align} \label{hamk7} 
\sigma (\tau) \! :=& \, \tau \mathcal{H}_{\epsilon_{1}}(\hat{p}(\tau),\hat{q}(\tau);\tau) \! + \! \hat{p}(\tau) \hat{q}(\tau) 
\! + \! \frac{1}{2}(\mi a \! + \! 1/2)^{2} \! - \! \epsilon_{1} (\mi a \! + \! 1/2) \! + \! \frac{1}{4} \nonumber \\
=& \, \left(\hat{p}(\tau) \hat{q}(\tau) \! - \! \epsilon_{1}(\mi a \! + \! (1 \! - \! \epsilon_{1})/2) \right)^{2} \! + \! \tau 
(4 \varepsilon \hat{q}(\tau) \! + \! \mi b \hat{p}(\tau))
\end{align}
satisfies the second-order non-linear ODE (related to the DP3E~\eqref{eq1.1})
\begin{equation} \label{hamk9} 
(\tau \sigma^{\prime \prime}(\tau) \! - \! \sigma^{\prime}(\tau))^{2} \! = \! 2(2 \sigma (\tau) \! - \! \tau \sigma^{\prime}(\tau))
(\sigma^{\prime}(\tau))^{2} \! - \! 32 \mi \varepsilon b \tau \left(((1 \! - \! \epsilon_{1})/2 \! - \! \mi a \epsilon_{1}) 
\sigma^{\prime}(\tau) \! + \! 2 \mi \varepsilon b \tau \right).
\end{equation}
Equation \eqref{hamk9} is referred to as the $\sigma$-form of the DP3E \eqref{eq1.1}. Motivated by the definition 
\eqref{hamk10} for the Hamiltonian function, setting $\epsilon_{1} \! = \! -1$, letting the generalised co-ordinate 
$\hat{q}(\tau) \! = \! u(\tau)$, and using the first (left-most) relation of equations~\eqref{hamk5} to calculate the 
generalised impulse, it suffices, for the purposes of the present work, to define the function (cf. definition \eqref{hamk7}) 
$\sigma (\tau)$ and the second-order non-linear ODE it satisfies as follows:
\begin{equation} \label{thmk23} 
\sigma (\tau) \! := \! \tau \mathcal{H}(\tau) \! + \! \dfrac{\tau (u^{\prime}(\tau) \! - \! \mi b)}{2u(\tau)} \! + \! 
\frac{1}{2}(\mi a \! + \! 1/2)^{2} \! + \! \frac{1}{4},
\end{equation}
and
\begin{equation} \label{thmk22} 
(\tau \sigma^{\prime \prime}(\tau) \! - \! \sigma^{\prime}(\tau))^{2} \! = \! 2(2 \sigma (\tau) \! - \! \tau \sigma^{\prime}(\tau))
(\sigma^{\prime}(\tau))^{2} \! - \! 32 \mi \varepsilon b \tau ((1 \! + \! \mi a) \sigma^{\prime}(\tau) \! + \! 2 \mi \varepsilon b \tau).
\end{equation}

Via the B\"{a}cklund transformations given in Subsection~6.1 of \cite{a1}, let
\begin{gather} 
u_{-}(\tau) \! := \! \frac{\mi \varepsilon b}{8(u(\tau))^{2}} \! \left(\tau (u^{\prime}(\tau) \! - \! \mi b) \! + \! (1 \! - \! 2 \mi a_{-})
u(\tau) \right), \label{yoominus1} \\
u_{+}(\tau) \! := \! -\frac{\mi \varepsilon b}{8(u(\tau))^{2}} \! \left(\tau (u^{\prime}(\tau) \! + \! \mi b) \! + \! (1 \! + \! 2 \mi 
a_{+})u(\tau) \right), \label{yooplus1}
\end{gather}
where $u(\tau)$ denotes any solution of the DP3E \eqref{eq1.1}, and $a_{\pm} \! := \! a \! \pm \! \mi$; in fact, as shown 
in Subsection 6.1 of \cite{a1}, $u_{-}(\tau)$ (resp., $u_{+}(\tau)$) solves the DP3E~\eqref{eq1.1} for $a \! = \! a_{-}$ 
(resp., $a \! = \! a_{+}$). {}From the results of \cite{avkavint}, define the two \emph{principal auxiliary functions} 
\begin{gather}
f_{-}(\tau) \! := \! -\frac{2 \mi}{\varepsilon b}u(\tau)u_{-}(\tau), \label{eqf} \\
f_{+}(\tau) \! := \! u(\tau)u_{+}(\tau), \label{yooplus2}
\end{gather}
where $f_{-}(\tau)$ solves the second-order non-linear ODE {}\footnote{This is a consequence of the ODE for the function 
$f(\tau)$ presented on p.~1168 of \cite{a1} upon making the notational change $f(\tau) \! \to \! f_{-}(\tau)$ and setting 
$\epsilon_{1} \! = \! -1$.} 
\begin{equation} \label{thmk13} 
\tau^{2} \! \left(f_{-}^{\prime \prime}(\tau) \! + \! 4 \mi \varepsilon b \right)^{2} \! - \! \left(4f_{-}(\tau) \! + \! 2 \mi a \! + \! 1 
\right)^{2} \left((f_{-}^{\prime}(\tau))^{2} \! + \! 8 \mi \varepsilon b f_{-}(\tau) \right) \! = \! 0,
\end{equation}
and $f_{+}(\tau)$ solves the second-order non-linear ODE {}\footnote{See equation~(2) in \cite{avkavint}.} 
\begin{equation} \label{yooplus3} 
(\varepsilon b \tau)^{2} \! \left(f_{+}^{\prime \prime}(\tau) \! - \! 2(\varepsilon b)^{2} \right)^{2} \! + \! \left(8f_{+}(\tau) \! + \! \mi 
\varepsilon b (2 \mi a \! - \! 1) \right)^{2} \left((f_{+}^{\prime}(\tau))^{2} \! - \! 4(\varepsilon b)^{2}f_{+}(\tau) \right) \! = \! 0.
\end{equation}
It follows {}from the definitions \eqref{yoominus1}--\eqref{yooplus2} that the functions $f_{\pm}(\tau)$ possess the 
alternative representations
\begin{gather}
2f_{-}(\tau) \! = \! -\mi (a \! - \! \mi/2) \! + \! \frac{\tau (u^{\prime}(\tau) \! - \! \mi b)}{2u(\tau)}, \label{yoominus3} \\
\frac{4 \mi}{\varepsilon b}f_{+}(\tau) \! = \! \mi (a \! + \! \mi/2) \! + \! \frac{\tau (u^{\prime}(\tau) \! + \! \mi b)}{2u(\tau)}; 
\label{yooplus4}
\end{gather}
incidentally, equations \eqref{yoominus3} and \eqref{yooplus4} imply the corollary
\begin{equation} \label{yoominus4} 
\frac{4 \mi}{\varepsilon b}f_{+}(\tau) \! = \! 2f_{-}(\tau) \! + \! \mi \tau \! \left(\frac{2a}{\tau} \! + \! \frac{b}{u(\tau)} \right).
\end{equation}
For the monodromy data considered in \cite{avlkv}, preliminary asymptotics as $\tau \! \to \! +\infty$ with $\varepsilon 
b \! > \! 0$ for $\int_{0}^{\tau} \xi^{-1}f_{+}(\xi) \, \md \xi$ have been presented in \cite{avkavint}.
\subsection{Lax Pairs and Isomonodromic Deformations} \label{sec1c} 
In this subsection, the reader is reminded about some basic facts regarding the isomonodromy deformation theory for 
the DP3E~\eqref{eq1.1}.
\begin{eeee} \label{pregauge} 
Pre-gauge-transformed Lax-pair-associated functions are denoted with `hats', whilst post-gauge-transformed 
Lax-pair-associated functions are not; in some cases, these functions are equal, and in others, they are not (see 
the discussion below). \hfill $\blacksquare$
\end{eeee}
The study of the DP3E \eqref{eq1.1} is based on the following pre-gauge-transformed Lax pair (see Proposition 2.1 
of \cite{a1}, with notational amendments):
\begin{equation} \label{eqFGmain}
\partial_{\mu} \widehat{\Psi}(\mu,\tau) \! = \! \widehat{\mathscr{U}}(\mu,\tau) \widehat{\Psi}(\mu,\tau), \qquad 
\qquad \partial_{\tau} \widehat{\Psi}(\mu,\tau) \! = \! \widehat{\mathscr{V}}(\mu,\tau) \widehat{\Psi}(\mu,\tau),
\end{equation}
where
\begin{align}
\widehat{\mathscr{U}}(\mu,\tau) =& \, -2 \mi \tau \mu \sigma_{3} \! + \! 2 \tau \! 
\begin{pmatrix}
0 & \frac{2 \mi \hat{A}(\tau)}{\sqrt{\smash[b]{-\hat{A}(\tau) \hat{B}(\tau)}}} \\
-\hat{D}(\tau) & 0
\end{pmatrix} 
\! - \! \dfrac{1}{\mu} \! \left(\mi a \! + \! \dfrac{1}{2} \! + \! \dfrac{2 \tau \hat{A}(\tau) \hat{D}(\tau)}{
\sqrt{\smash[b]{-\hat{A}(\tau) \hat{B}(\tau)}}} \right) \! \sigma_{3} \nonumber \\
+& \, \dfrac{1}{\mu^{2}} \! 
\begin{pmatrix}
0 & \hat{\alpha}(\tau) \\
\mi \tau \hat{B}(\tau) & 0
\end{pmatrix}, \label{mpeea1} \\
\widehat{\mathscr{V}}(\mu,\tau) =& \, -\mi \mu^{2} \sigma_{3} \! + \! \mu \! 
\begin{pmatrix}
0 & \frac{2 \mi \hat{A}(\tau)}{\sqrt{\smash[b]{-\hat{A}(\tau) \hat{B}(\tau)}}} \\
-\hat{D}(\tau) & 0
\end{pmatrix} 
\! + \! \left(\dfrac{\mi a}{2 \tau} \! - \! \dfrac{\hat{A}(\tau) \hat{D}(\tau)}{
\sqrt{\smash[b]{-\hat{A}(\tau) \hat{B}(\tau)}}} \right) \! \sigma_{3} \nonumber \\
-& \, \frac{1}{\mu} \frac{1}{2 \tau} \! 
\begin{pmatrix}
0 & \hat{\alpha}(\tau) \\
\mi \tau \hat{B}(\tau) & 0
\end{pmatrix}, \label{mpeea2} 
\end{align}
with $\sigma_{3} \! = \! \diag (1,-1)$,
\begin{equation} \label{firstintegral} 
\hat{\alpha}(\tau) \! := \! -2(\hat{B}(\tau))^{-1} \! \left(\mi a \sqrt{
\smash[b]{-\hat{A}(\tau) \hat{B}(\tau)}} + \! \tau (\hat{A}(\tau) \hat{D}
(\tau) \! + \! \hat{B}(\tau) \hat{C}(\tau)) \right),
\end{equation}
and where the differentiable, scalar-valued functions $\hat{A}(\tau)$, $\hat{B}(\tau)$, $\hat{C}(\tau)$, and 
$\hat{D}(\tau)$ satisfy the system of isomonodromy deformations
\begin{equation} \label{eq1.4}
\begin{gathered}
\hat{A}^{\prime}(\tau) \! = \! 4 \hat{C}(\tau) \sqrt{\smash[b]{-\hat{A}(\tau) \hat{B}(\tau)}}, \, \qquad \, \quad 
\, \hat{B}^{\prime}(\tau) \! = \! -4 \hat{D}(\tau) \sqrt{\smash[b]{-\hat{A}(\tau) \hat{B}(\tau)}}, \\
(\tau \hat{C}(\tau))^{\prime} \! = \! 2 \mi a \hat{C}(\tau) \! - \! 2 \tau \hat{A}(\tau), \, \quad \, \qquad \, (\tau 
\hat{D}(\tau))^{\prime} \! = \! -2 \mi a \hat{D}(\tau) \! + \! 2 \tau \hat{B}(\tau), \\
\left(\sqrt{\smash[b]{-\hat{A}(\tau) \hat{B}(\tau)}} \, \right)^{\prime} \! = \! 2(\hat{A}(\tau) \hat{D}(\tau) \! - \! 
\hat{B}(\tau) \hat{C}(\tau)).
\end{gathered}
\end{equation}
(Note: the isomonodromy deformations~\eqref{eq1.4} are, for arbitrary values of $\mu \! \in \! \mathbb{C}$, the 
Frobenius compatibility condition for the system \eqref{eqFGmain}.)
\begin{eeee} \label{alphwave} 
In fact, $-\mi \hat{\alpha}(\tau) \hat{B}(\tau) \! = \! \varepsilon b$, $\varepsilon \! = \! \pm 1$, so that the definition 
\eqref{firstintegral} is the \emph{first integral\/} of the system \eqref{eq1.4} (see Lemma {\rm 2.1} of \cite{a1}, 
with notational amendments). \hfill $\blacksquare$
\end{eeee}
\begin{eeee} \label{xalp} 
With conspicuous changes in notation (cf. system~(4) in \cite{a1}), whilst transforming {}from the original Lax pair
\begin{gather*}
\partial_{\lambda} \Phi (\lambda,\tau) \! = \! \tau \! \left(-\mi \sigma_{3} \! - \! 
\dfrac{1}{\lambda} \dfrac{\mi a}{2 \tau} \sigma_{3} \! - \! \dfrac{1}{\lambda} \! 
\begin{pmatrix}
0 & \hat{C}(\tau) \\
\hat{D}(\tau) & 0
\end{pmatrix} \! + \! \dfrac{1}{\lambda^{2}} \dfrac{\mi}{2} \! 
\begin{pmatrix}
\sqrt{\smash[b]{-\hat{A}(\tau) \hat{B}(\tau)}} & \hat{A}(\tau) \\
\hat{B}(\tau) & -\sqrt{\smash[b]{-\hat{A}(\tau) \hat{B}(\tau)}} 
\end{pmatrix} \right) \! \Phi (\lambda,\tau), \\
\partial_{\tau} \Phi (\lambda,\tau) \! = \! \left(-\mi \lambda \sigma_{3} 
\! + \! \dfrac{\mi a}{2 \tau} \sigma_{3} \! - \! 
\begin{pmatrix}
0 & \hat{C}(\tau) \\
\hat{D}(\tau) & 0
\end{pmatrix} \! - \! \dfrac{1}{\lambda} \dfrac{\mi}{2} \! 
\begin{pmatrix}
\sqrt{\smash[b]{-\hat{A}(\tau) \hat{B}(\tau)}} & \hat{A}(\tau) \\
\hat{B}(\tau) & -\sqrt{\smash[b]{-\hat{A}(\tau) \hat{B}(\tau)}} 
\end{pmatrix} \right) \! \Phi (\lambda,\tau),
\end{gather*}
to the Fuchs-Garnier pair \eqref{eqFGmain}, the Fabry-type transformation (cf. Proposition~2.1 in \cite{a1})
\begin{equation*}
\lambda \! = \! \mu^{2} \qquad \text{and} \qquad \Phi (\lambda,\tau) \! := \! \sqrt{\smash[b]{\mu}} 
\left(
\begin{pmatrix}
1 & 0 \\
0 & 0
\end{pmatrix} \! + \! \dfrac{1}{\mu} \! 
\begin{pmatrix}
0 & -\frac{\hat{A}(\tau)}{\sqrt{\smash[b]{-\hat{A}(\tau) \hat{B}(\tau)}}} \\
0 & 1
\end{pmatrix} \right) \! \widehat{\Psi}(\mu,\tau)
\end{equation*}
was used; if, instead, one applies the slightly more general transformation
\begin{equation*}
\Phi (\lambda,\tau) \! := \! \sqrt{\smash[b]{\mu}} \! 
\left(
\begin{pmatrix}
1 & 0 \\
0 & 0
\end{pmatrix} \! + \! \dfrac{1}{\mu} \! 
\begin{pmatrix}
-\frac{\hat{A}(\tau) \mathbb{P}^{\ast}}{\sqrt{\smash[b]{-\hat{A}(\tau) \hat{B}(\tau)}}} 
& -\frac{\hat{A}(\tau)}{\sqrt{\smash[b]{-\hat{A}(\tau) \hat{B}(\tau)}}} \\
\mathbb{P}^{\ast} & 1
\end{pmatrix} \right) \! \widehat{\Psi}(\mu,\tau)
\end{equation*}
for some constant or $\tau$-dependent $\mathbb{P}^{\ast}$, then, in lieu of, say, the $\mu$-part of the Fuchs-Garnier 
pair \eqref{eqFGmain}, that is, $\partial_{\mu} \widehat{\Psi}(\mu,\tau) \! = \! \widehat{\mathscr{U}}(\mu,\tau) 
\widehat{\Psi}(\mu,\tau)$, one arrives at
\begin{equation*}
\partial_{\mu} \widehat{\Psi}(\mu,\tau) \! = \! \left(\widehat{\mathcal{L}}_{-1} \mu \! + \! \widehat{\mathcal{L}}_{0} 
\! + \! \widehat{\mathcal{L}}_{1} \mu^{-1} \! + \! \widehat{\mathcal{L}}_{2} \mu^{-2} \right) \! \widehat{\Psi}(\mu,\tau),
\end{equation*}
where
\begin{gather*}
\widehat{\mathcal{L}}_{-1} \! = \! -2 \mi \tau \! 
\begin{pmatrix}
1 & 0 \\
-2 \mathbb{P}^{\ast} & -1
\end{pmatrix}, \, \qquad \, \widehat{\mathcal{L}}_{0} \! = \! -2 \tau \! 
\begin{pmatrix}
0 & 0 \\
\hat{D}(\tau) & 0
\end{pmatrix} \! - \! \dfrac{4 \mi \tau \hat{A}(\tau)}{\sqrt{\smash[b]{-\hat{A}(\tau) \hat{B}(\tau)}}} \! 
\begin{pmatrix}
-\mathbb{P}^{\ast} & -1 \\
(\mathbb{P}^{\ast})^{2} & \mathbb{P}^{\ast}
\end{pmatrix}, \\
\widehat{\mathcal{L}}_{1} \! = \! \left(\mi a \! + \! \dfrac{1}{2} \! + \! \dfrac{2 \tau \hat{A}
(\tau) \hat{D}(\tau)}{\sqrt{\smash[b]{-\hat{A}(\tau) \hat{B}(\tau)}}} \right) \! 
\begin{pmatrix}
-1 & 0 \\
2 \mathbb{P}^{\ast} & 1
\end{pmatrix}, \, \quad \, \widehat{\mathcal{L}}_{2} \! = \! \mi \tau \! 
\begin{pmatrix}
0 & 0 \\
\hat{B}(\tau) & 0
\end{pmatrix} \! + \! \hat{\alpha}(\tau) \! 
\begin{pmatrix}
\mathbb{P}^{\ast} & 1 \\
-(\mathbb{P}^{\ast})^{2} & -\mathbb{P}^{\ast}
\end{pmatrix},
\end{gather*}
with $\hat{\alpha}(\tau)$ defined by equation \eqref{firstintegral}. Setting $\mathbb{P}^{\ast} \! = \! 0$, one arrives at the 
Fuchs-Garnier pair stated in Proposition~2.1 of \cite{a1}, system~(1.4) of \cite{av2}, and system \eqref{eqFGmain} of 
the present work. \hfill $\blacksquare$
\end{eeee}

A relation between the Fuchs-Garnier pair~\eqref{eqFGmain} and the DP3E~\eqref{eq1.1} is given by (see, in particular, 
Proposition~1.2 of \cite{a1}, with notational amendments)
\begin{bbbb}[{\rm \cite{a1,av2}}] \label{prop1.2}
Let $\hat{u} \! = \! \hat{u}(\tau)$ and $\hat{\varphi} \! = \! \hat{\varphi}(\tau)$ solve the system
\begin{equation} \label{eq1.5}
\begin{gathered}
\hat{u}^{\prime \prime}(\tau) \! = \! \dfrac{(\hat{u}^{\prime}(\tau))^{2}}{\hat{u}(\tau)} \! - \! 
\dfrac{\hat{u}^{\prime}(\tau)}{\tau} \! + \! \dfrac{1}{\tau} \! \left(-8 \varepsilon (\hat{u}(\tau))^{2} 
\! + \! 2ab \right) \! + \! \dfrac{b^{2}}{\hat{u}(\tau)}, \quad \quad \hat{\varphi}^{\prime}(\tau) 
\! = \! \dfrac{2a} \tau \! + \! \dfrac{b}{\hat{u}(\tau)},
\end{gathered}
\end{equation}
where $\varepsilon \! = \! \pm 1$, and $a,b \! \in \! \mathbb{C}$ are independent of $\tau$$;$ then,
\begin{equation} \label{eq:ABCD} 
\begin{gathered}
\hat{A}(\tau) \! := \! \frac{\hat{u}(\tau)}{\tau} \me^{\mi \hat{\varphi}(\tau)}, \qquad \qquad 
\hat{B}(\tau) \! := \! -\frac{\hat{u}(\tau)}{\tau} \me^{-\mi \hat{\varphi}(\tau)}, \\
\hat{C}(\tau) \! := \! \dfrac{\varepsilon \tau \hat{A}^{\prime}(\tau)}{4 \hat{u}(\tau)} \! = \! 
\frac{\varepsilon \me^{\mi \hat{\varphi}(\tau)}}{2 \tau} \! \left(\mi (a \! + \! \mi/2) \! + \! 
\frac{\tau (\hat{u}^{\prime}(\tau) \! + \! \mi b)}{2 \hat{u}(\tau)} \right), \\
\hat{D}(\tau) \! := \! -\dfrac{\varepsilon \tau \hat{B}^{\prime}(\tau)}{4 \hat{u}(\tau)} \! = \! 
-\frac{\varepsilon \me^{-\mi \hat{\varphi}(\tau)}}{2 \tau} \! \left(\mi (a \! - \! \mi/2) \! - \! 
\frac{\tau (\hat{u}^{\prime}(\tau) \! - \! \mi b)}{2 \hat{u}(\tau)} \right)
\end{gathered}
\end{equation}
solve the system~\eqref{eq1.4}. Conversely, let $\hat{A}(\tau) \! \not\equiv \! 0$, $\hat{B}(\tau) \! \not\equiv \! 0$, 
$\hat{C}(\tau)$, and $\hat{D}(\tau)$ solve the system~\eqref{eq1.4}, and define
\begin{equation} \label{tempeq}
\hat{u}(\tau) \! := \! \varepsilon \tau \sqrt{\smash[b]{-\hat{A}(\tau) \hat{B}(\tau)}}, 
\quad \hat{\varphi}(\tau) \! := \! -\dfrac{\mi}{2} \ln \! \left(-\hat{A}(\tau)/\hat{B}
(\tau) \right), \quad b \! := \! \hat{u}(\tau) \! \left(\hat{\varphi}^{\prime}(\tau) 
\! - \! 2a \tau^{-1} \right);
\end{equation}
then, $b$ is independent of $\tau$, and $\hat{u}(\tau)$ and $\hat{\varphi}(\tau)$ solve the system~\eqref{eq1.5}.
\end{bbbb}
\begin{bbbb} \label{equzforhatf} 
Let (cf. equation \eqref{yoominus3}$)$
\begin{equation} \label{hatsoff1} 
2 \hat{f}_{-}(\tau) \! := \! -\mi (a \! - \! \mi/2) \! + \! \frac{\tau}{2} \! \left(\frac{\hat{u}^{\prime}(\tau) \! - \! \mi b}{\hat{u}(\tau)} \right),
\end{equation}
and (cf. equation \eqref{yooplus4}$)$
\begin{equation} \label{pga1} 
\frac{4 \mi}{\varepsilon b} \hat{f}_{+}(\tau) \! := \! \mi (a \! + \! \mi/2) \! + \! \frac{\tau}{2} \! \left(\frac{\hat{u}^{\prime}(\tau) 
\! + \! \mi b}{\hat{u}(\tau)} \right).
\end{equation}
Then, for $\varepsilon \! \in \! \lbrace \pm 1 \rbrace$,
\begin{gather}
2 \hat{f}_{-}(\tau) \! = \! \frac{2 \varepsilon \tau^{2} \hat{A}(\tau) \hat{D}(\tau)}{\hat{u}(\tau)} \! = \! \frac{\tau}{2} 
\frac{\md}{\md \tau} \! \left(\ln \! \left(\frac{\hat{u}(\tau)}{\tau} \right) \! - \! \mi \hat{\varphi}(\tau) \right), \label{hatsoff2} \\
\intertext{and} 
\frac{4 \mi}{\varepsilon b} \hat{f}_{+}(\tau) \! = \! -\frac{2 \varepsilon \tau^{2} \hat{B}(\tau) \hat{C}(\tau)}{\hat{u}(\tau)} \! 
= \! \frac{\tau}{2} \frac{\md}{\md \tau} \! \left(\ln \! \left(\frac{\hat{u}(\tau)}{\tau} \right) \! + \! \mi \hat{\varphi}(\tau) \right); 
\label{hatsoff3}
\end{gather}
furthermore,
\begin{equation} \label{pga2} 
\frac{4 \mi}{\varepsilon b} \hat{f}_{+}(\tau) \! = \! 2 \hat{f}_{-}(\tau) \! + \! \mi \tau \hat{\varphi}^{\prime}(\tau) \! = \! 2 
\hat{f}_{-}(\tau) \! + \! \mi \tau \! \left(\frac{2a}{\tau} \! + \! \frac{b}{\hat{u}(\tau)} \right).
\end{equation}
\end{bbbb}

\emph{Proof}. Without loss of generality, consider, say, the proof for the function $\hat{f}_{-}(\tau)$: the proof for the 
function $\hat{f}_{+}(\tau)$ is analogous. One commences by establishing the following relation:
\begin{equation} \label{hatsoff4} 
\frac{\hat{u}^{\prime}(\tau) \! - \! \mi b}{\hat{u}(\tau)} \! = \! \frac{2}{\tau} \! \left(\frac{2 \tau \hat{A}(\tau) 
\hat{D}(\tau)}{\sqrt{\smash[b]{-\hat{A}(\tau) \hat{B}(\tau)}}} \! + \! (\mi a \! + \! 1/2) \right).
\end{equation}
{}From definition \eqref{firstintegral}, the system of isomonodromy deformations \eqref{eq1.4}, Remark \ref{alphwave}, and 
the definition of the function $\hat{u}(\tau)$ given by the first (left-most) member of equations \eqref{tempeq}, it follows via 
differentiation that
\begin{align*} 
\frac{\hat{u}^{\prime}(\tau) \! - \! \mi b}{\hat{u}(\tau)} =& \, \frac{2 \tau (\hat{A}
(\tau) \hat{D}(\tau) \! - \! \hat{B}(\tau) \hat{C}(\tau)) \! + \! \sqrt{\smash[b]{-\hat{A}
(\tau) \hat{B}(\tau)}}}{\tau \sqrt{\smash[b]{-\hat{A}(\tau) \hat{B}(\tau)}}} \! - \! 
\frac{\mi (\varepsilon b)}{\varepsilon \hat{u}(\tau)} \\
=& \, \frac{2 \tau (\hat{A}(\tau) \hat{D}(\tau) \! - \! \hat{B}(\tau) \hat{C}(\tau)) \! + 
\! \sqrt{\smash[b]{-\hat{A}(\tau) \hat{B}(\tau)}}}{\tau \sqrt{\smash[b]{-\hat{A}(\tau) 
\hat{B}(\tau)}}} \! - \! \frac{\hat{\alpha}(\tau) \hat{B}(\tau)}{\varepsilon \hat{u}(\tau)} \\
=& \, \frac{2 \tau (\hat{A}(\tau) \hat{D}(\tau) \! - \! \hat{B}(\tau) \hat{C}(\tau)) \! + 
\! \sqrt{\smash[b]{-\hat{A}(\tau) \hat{B}(\tau)}}}{\tau \sqrt{\smash[b]{-\hat{A}(\tau) 
\hat{B}(\tau)}}} \\
+& \, \frac{2 \tau (\hat{A}(\tau) \hat{D}(\tau) \! + \! \hat{B}(\tau) \hat{C}(\tau)) \! + \! 
2 \mi a \sqrt{\smash[b]{-\hat{A}(\tau) \hat{B}(\tau)}}}{\tau \sqrt{\smash[b]{-\hat{A}
(\tau) \hat{B}(\tau)}}} \\
=& \, \frac{2}{\tau} \! \left(\frac{2 \tau \hat{A}(\tau) \hat{D}(\tau)}{\sqrt{\smash[b]{
-\hat{A}(\tau) \hat{B}(\tau)}}} \! + \! (\mi a \! + \! 1/2) \right);
\end{align*} 
conversely, {}from the system of isomonodromy deformations \eqref{eq1.4}, the system \eqref{eq1.5}, and the definitions 
\eqref{eq:ABCD} and \eqref{tempeq}, it follows that
\begin{align*}
\frac{4 \hat{A}(\tau) \hat{D}(\tau)}{\sqrt{\smash[b]{-\hat{A}(\tau) \hat{B}(\tau)}}} =& \, \frac{4 \varepsilon \tau \hat{A}(\tau) 
\hat{D}(\tau)}{\hat{u}(\tau)} \! = \! \frac{4 \varepsilon \tau}{\hat{u}(\tau)} \! \left(-\frac{\varepsilon}{4} \hat{B}^{\prime}(\tau) 
\me^{\mi \hat{\varphi}(\tau)} \right) \! = \! \frac{\tau \me^{\mi \hat{\varphi}(\tau)}}{\hat{u}(\tau)} \frac{\md}{\md \tau} \! \left(
\frac{\hat{u}(\tau)}{\tau} \me^{-\mi \hat{\varphi}(\tau)} \right) \\
=& \, \frac{\tau}{\hat{u}(\tau)} \! \left(-\mi \hat{\varphi}^{\prime}(\tau) \frac{\hat{u}(\tau)}{\tau} \! - \! \frac{\hat{u}(\tau)}{\tau^{2}} 
\! + \! \frac{\hat{u}^{\prime}(\tau)}{\tau} \right) \! = \! \frac{\tau}{\hat{u}(\tau)} \! \left(-\frac{\hat{u}(\tau)}{\tau} \! \left(
\frac{2 \mi a}{\tau} \! + \! \frac{\mi b}{\hat{u}(\tau)} \right) \! - \! \frac{\hat{u}(\tau)}{\tau^{2}} \! + \! \frac{\hat{u}^{\prime}
(\tau)}{\tau} \right) \\
=& \, \frac{\hat{u}^{\prime}(\tau) \! - \! \mi b}{\hat{u}(\tau)} \! - \! \frac{2}{\tau}(\mi a \! + \! 1/2),
\end{align*}
whence
\begin{equation*}
\frac{2}{\tau} \! \left(\frac{2 \tau \hat{A}(\tau) \hat{D}(\tau)}{\sqrt{\smash[b]{-\hat{A}(\tau) \hat{B}(\tau)}}} \! + \! 
(\mi a \! + \! 1/2) \right) \! = \! \frac{\hat{u}^{\prime}(\tau) \! - \! \mi b}{\hat{u}(\tau)},
\end{equation*}
which establishes equation \eqref{hatsoff4}. Via definition~\eqref{hatsoff1} and equation~\eqref{hatsoff4}, one shows that
\begin{equation} \label{hatsoff5} 
\hat{f}_{-}(\tau) \! = \! \frac{\tau \hat{A}(\tau) \hat{D}(\tau)}{\sqrt{\smash[b]{-\hat{A}(\tau) \hat{B}(\tau)}}},
\end{equation}
hence, via the definition for $\hat{u}(\tau)$ given by the first (left-most) member of equations~\eqref{tempeq}, one arrives at 
the first (left-most) relation of equation \eqref{hatsoff2}; moreover, it follows {}from the ODE for the function $\hat{\varphi}(\tau)$ 
given in the system \eqref{eq1.5} and definition \eqref{hatsoff1} that
\begin{equation*}
\tau^{-1} \hat{f}_{-}(\tau) \! = \! \frac{1}{4} \! \left(\frac{\hat{u}^{\prime}(\tau)}{\hat{u}(\tau)} \! + \! \frac{2 \mi a}{\tau} \! - \! 
\mi \hat{\varphi}^{\prime}(\tau) \right) \! - \! \frac{1}{2 \tau}(\mi a \! + \! 1/2) \! = \! \frac{1}{4} \! \left(\frac{\md}{\md \tau} 
\ln \! \left(\frac{\hat{u}(\tau)}{\tau} \right) \! - \! \mi \hat{\varphi}^{\prime}(\tau) 
\right),
\end{equation*}
which implies the second (right-most) relation of equation \eqref{hatsoff2}. Equations \eqref{hatsoff2} and \eqref{hatsoff3} 
imply the Corollary \eqref{pga2}, which is consistent with, and can also be derived {}from, the definition \eqref{firstintegral} 
and the first integral of system~\eqref{eq1.4} (cf. Remark~\ref{alphwave}). \hfill $\qed$

Herewith follows the post-gauge-transformed Fuchs-Garnier pair.
\begin{bbbb} \label{newlax1} 
Let $\widehat{\Psi}(\mu,\tau)$ be a fundamental solution of the system~\eqref{eqFGmain}. Set
\begin{equation} \label{newlax2}
\begin{gathered}
A(\tau) \! := \! \hat{A}(\tau) \tau^{-\mi a}, \, \quad \, B(\tau) \! := \! \hat{B}(\tau) \tau^{\mi a}, \, \quad \, C(\tau) \! 
:= \! \hat{C}(\tau) \tau^{-\mi a}, \quad D(\tau) \! := \! \hat{D}(\tau) \tau^{\mi a}, \\
\alpha (\tau) \! := \! \hat{\alpha}(\tau) \tau^{-\mi a}, \, \quad \, \widehat{\Psi}(\mu,\tau) \! := \! \tau^{\frac{\mi a}{2} 
\sigma_{3}} \Psi (\mu,\tau).
\end{gathered}
\end{equation}
Then$:$ {\rm (i)} $\Psi (\mu,\tau)$ is a fundamental solution of
\begin{equation} \label{newlax3}
\partial_{\mu} \Psi(\mu,\tau) \! = \! \widetilde{\mathscr{U}}(\mu,\tau) \Psi (\mu,\tau), \qquad \quad \partial_{\tau} 
\Psi(\mu,\tau) \! = \! \widetilde{\mathscr{V}}(\mu,\tau) \Psi(\mu,\tau),
\end{equation}
where
\begin{align}
\widetilde{\mathscr{U}}(\mu,\tau) =& \, -2 \mi \tau \mu \sigma_{3} \! + \! 2 \tau \! 
\begin{pmatrix}
0 & \frac{2 \mi A(\tau)}{\sqrt{\smash[b]{-A(\tau)B(\tau)}}} \\
-D(\tau) & 0
\end{pmatrix} 
\! - \! \dfrac{1}{\mu} \! \left(\mi a \! + \! \dfrac{1}{2} \! + \! \dfrac{2 \tau A(\tau) D(\tau)}{\sqrt{\smash[b]{-A(\tau)B(\tau)}}} 
\right) \! \sigma_{3} \nonumber \\
+& \, \dfrac{1}{\mu^{2}} \! 
\begin{pmatrix}
0 & \alpha (\tau) \\
\mi \tau B(\tau) & 0
\end{pmatrix}, \label{nlxa} \\
\widetilde{\mathscr{V}}(\mu,\tau) =& \, -\mi \mu^{2} \sigma_{3} \! + \! \mu \! 
\begin{pmatrix}
0 & \frac{2 \mi A(\tau)}{\sqrt{\smash[b]{-A(\tau)B(\tau)}}} \\
-D(\tau) & 0
\end{pmatrix} 
\! - \! \dfrac{A(\tau)D(\tau)}{\sqrt{\smash[b]{-A(\tau) B(\tau)}}} \sigma_{3} \! 
- \! \frac{1}{\mu} \frac{1}{2 \tau} \! 
\begin{pmatrix}
0 & \alpha (\tau) \\
\mi \tau B(\tau) & 0
\end{pmatrix}, \label{nlxb}
\end{align}
with
\begin{equation} \label{aphnovij} 
\alpha (\tau) \! := \! -2(B(\tau))^{-1} \! \left(\mi a \sqrt{\smash[b]{-A(\tau)B(\tau)}} + \! \tau (A(\tau)D(\tau) \! + \! 
B(\tau)C(\tau)) \right);
\end{equation}
and {\rm (ii)} if the coefficient functions $\hat{A}(\tau)$, $\hat{B}(\tau)$, $\hat{C}(\tau)$, and $\hat{D}(\tau)$ satisfy 
the system of isomonodromy deformations \eqref{eq1.4} and the functions $A(\tau)$, $B(\tau)$, $C(\tau)$, and 
$D(\tau)$ are defined by equations \eqref{newlax2}, then the Frobenius compatibility condition of the system 
\eqref{newlax3}, for arbitrary values of $\mu \! \in \! \mathbb{C}$, is that the differentiable, scalar-valued functions 
$A(\tau)$, $B(\tau)$, $C(\tau)$, and $D(\tau)$ satisfy the corresponding system of isomonodromy deformations
\begin{equation} \label{newlax8}
\begin{gathered}
A^{\prime}(\tau) \! = \! -\frac{\mi a}{\tau}A(\tau) \! + \! 4C(\tau) \sqrt{\smash[b]{-A(\tau)B(\tau)}}, \quad \, \, 
\quad B^{\prime}(\tau) \! = \! \frac{\mi a}{\tau}B(\tau) \! - \! 4D(\tau) \sqrt{\smash[b]{-A(\tau)B(\tau)}}, \\
(\tau C(\tau))^{\prime} \! = \! \mi aC(\tau) \! - \! 2 \tau A(\tau), \, \, \quad \quad \, \, (\tau D(\tau))^{\prime} \! = \! 
-\mi aD(\tau) \! + \! 2 \tau B(\tau), \\
\left(\sqrt{\smash[b]{-A(\tau)B(\tau)}} \, \right)^{\prime} \! = \! 2(A(\tau)D(\tau) \! - \! B(\tau)C(\tau)).
\end{gathered}
\end{equation}
\end{bbbb}

\emph{Proof}. If $\widehat{\Psi}(\mu,\tau)$ is a fundamental solution of the system~\eqref{eqFGmain}, then it follows 
{}from the isomonodromy deformations \eqref{eq1.4} and the definitions \eqref{newlax2} that $\Psi (\mu,\tau)$ solves 
the system \eqref{newlax3}, and that the functions $A(\tau)$, $B(\tau)$, $C(\tau)$, and $D(\tau)$ satisfy the corresponding 
isomonodromy deformations \eqref{newlax8}. One verifies the Frobenius compatibility condition for the system 
\eqref{newlax3} by showing that, $\forall \, \mu \! \in \! \mathbb{C}$, $\partial_{\tau} \widetilde{\mathscr{U}}(\mu,\tau) \! - \! 
\partial_{\mu} \widetilde{\mathscr{V}}(\mu,\tau) \! + \! [\widetilde{\mathscr{U}}(\mu,\tau),\widetilde{\mathscr{V}}(\mu,\tau)] \! = \! 
\left(
\begin{smallmatrix}
0 & 0 \\
0 & 0
\end{smallmatrix}
\right)$, where, for $\mathfrak{X},\mathfrak{Y} \! \in \! \mathrm{M}_{2}(\mathbb{C})$, $[\mathfrak{X},\mathfrak{Y}] \! 
:= \! \mathfrak{X} \mathfrak{Y} \! - \! \mathfrak{Y} \mathfrak{X}$ is the matrix commutator. \hfill $\qed$
\begin{eeee} \label{newlax6} 
Definitions \eqref{firstintegral}, \eqref{newlax2}, and \eqref{aphnovij}, and Remark \ref{alphwave} imply that 
$-\mi \alpha (\tau)B(\tau) \! = \! \varepsilon b$, $\varepsilon \! = \! \pm 1$. \hfill $\blacksquare$
\end{eeee}
\begin{bbbb} \label{propuf} 
Let $u(\tau)$ and $\varphi(\tau)$ solve the system
\begin{equation} \label{equu17}
\begin{gathered}
u^{\prime \prime}(\tau) \! = \! \dfrac{(u^{\prime}(\tau))^{2}}{u(\tau)} \! - \! \dfrac{u^{\prime}(\tau)}{\tau} \! + \! 
\dfrac{1}{\tau} \! \left(-8 \varepsilon (u(\tau))^{2} \! + \! 2ab \right) \! + \! \dfrac{b^{2}}{u(\tau)}, \quad \quad 
\varphi^{\prime}(\tau) \! = \! \dfrac{a} \tau \! + \! \dfrac{b}{u(\tau)},
\end{gathered}
\end{equation}
where $\varepsilon \! = \! \pm 1$, and $a,b \! \in \! \mathbb{C}$ are independent of $\tau$$;$ then,
\begin{equation} \label{equu18} 
\begin{gathered}
A(\tau) \! := \! \frac{u(\tau)}{\tau} \me^{\mi \varphi(\tau)}, \qquad \quad B(\tau) \! := \! -\frac{u(\tau)}{\tau} 
\me^{-\mi \varphi(\tau)}, \\
C(\tau) \! := \! \dfrac{\varepsilon \tau}{4u(\tau)} \! \left(A^{\prime}(\tau) \! + \! \frac{\mi a}{\tau}A(\tau) \right) \! = \! 
\frac{\varepsilon \me^{\mi \varphi (\tau)}}{2 \tau} \! \left(\mi (a \! + \! \mi/2) \! + \! \frac{\tau (u^{\prime}(\tau) \! + 
\! \mi b)}{2u(\tau)} \right), \\
D(\tau) \! := \! -\dfrac{\varepsilon \tau}{4u(\tau)} \! \left(B^{\prime}(\tau) \! - \! \frac{\mi a}{\tau}B(\tau) \right) \! = \! 
-\frac{\varepsilon \me^{-\mi \varphi (\tau)}}{2 \tau} \! \left(\mi (a \! - \! \mi/2) \! - \! \frac{\tau (u^{\prime}(\tau) \! - \! 
\mi b)}{2u(\tau)} \right)
\end{gathered}
\end{equation}
solve the system \eqref{newlax8}. Conversely, let $A(\tau) \! \not\equiv \! 0$, $B(\tau) \! \not\equiv \! 0$, $C(\tau)$, 
and $D(\tau)$ solve the system \eqref{newlax8}, and define
\begin{equation} \label{equu19} 
u(\tau) \! := \! \varepsilon \tau \sqrt{\smash[b]{-A(\tau)B(\tau)}}, \quad \varphi(\tau) \! := \! -\frac{\mi}{2} \ln \! 
\left(-A(\tau)/B(\tau) \right), \quad b \! := \! u(\tau) \! \left(\varphi^{\prime}(\tau) \! - \! a \tau^{-1} \right);
\end{equation}
then, $b$ is independent of $\tau$, and $u(\tau)$ and $\varphi (\tau)$ solve the system \eqref{equu17}.
\end{bbbb}

\emph{Proof}. Via the definition of $\hat{u}(\tau)$ given by the first (left-most) member of equations \eqref{tempeq} 
and the definitions \eqref{newlax2}, one arrives at the definition for $u(\tau)$ given by the first (left-most) member 
of equations \eqref{equu19}; in particular, it follows that $u(\tau) \! = \! \hat{u}(\tau)$, and, {}from the first equation 
of system \eqref{eq1.5}, $u(\tau)$ solves the DP3E \eqref{eq1.1} (see the first equation of the system 
\eqref{equu17}). Let $\varphi (\tau)$ be defined as in equations \eqref{equu19}, that is, $\varphi (\tau) \! = \! 
-\mi \ln (\sqrt{\smash[b]{-A(\tau)B(\tau)}}/B(\tau))$; then, via differentiation, the definition \eqref{aphnovij}, and the 
corresponding system of isomonodromy deformations \eqref{newlax8}, it follows that
\begin{align*}
\varphi^{\prime}(\tau) =& \, -\mi \! \left(\frac{1}{\sqrt{\smash[b]{-A(\tau)B(\tau)}}} 
\left(\sqrt{\smash[b]{-A(\tau)B(\tau)}} \, \right)^{\prime} \! - \! \frac{B^{\prime}
(\tau)}{B(\tau)} \right) \nonumber \\
=& \, -\mi \! \left(\frac{2(A(\tau)D(\tau) \! - \! B(\tau)C(\tau))}{\sqrt{\smash[b]{
-A(\tau)B(\tau)}}} \! - \! \frac{1}{B(\tau)} \left(\frac{\mi a}{\tau}B(\tau) \! - \! 
4D(\tau) \sqrt{\smash[b]{-A(\tau)B(\tau)}} \right) \right) \nonumber \\
=& \, -\frac{a}{\tau} \! + \! \frac{2 \mi}{\sqrt{\smash[b]{-A(\tau)B(\tau)}}}
(A(\tau)D(\tau) \! + \! B(\tau)C(\tau)) \nonumber \\
=& \, -\frac{a}{\tau} \! + \! \frac{2 \mi}{\sqrt{\smash[b]{-A(\tau)B(\tau)}}} \! 
\left(-\frac{\mi \varepsilon b}{2 \tau} \! - \! \frac{\mi a}{\tau} \sqrt{\smash[b]{
-A(\tau)B(\tau)}} \right) \! = \! \frac{a}{\tau} \! + \! \frac{b}{u(\tau)},
\end{align*}
that is, $\varphi (\tau)$ solves the ODE given by the second (right-most) member of the system \eqref{equu17}; 
moreover, it also follows {}from the definitions \eqref{tempeq}, \eqref{newlax2}, and \eqref{equu19} that
\begin{equation} \label{equu20} 
\varphi (\tau) \! = \! \hat{\varphi}(\tau) \! - \! a \ln \tau.
\end{equation}
The definitions~\eqref{equu18} for the functions $A(\tau)$, $B(\tau)$, $C(\tau)$, and $D(\tau)$ are a consequence 
of the definitions \eqref{eq:ABCD} and \eqref{newlax2}, the fact that $u(\tau) \! = \! \hat{u}(\tau)$, and equation 
\eqref{equu20}. A series of lengthy, but otherwise straightforward, differentiation arguments completes the proof. 
\hfill $\qed$
\begin{eeee} \label{phitophi} 
It also follows {}from the ODE satisfied by $\hat{\varphi}(\tau)$ given in the system \eqref{eq1.5} and equation 
\eqref{equu20} that $\varphi (\tau)$ solves the corresponding ODE given in the system \eqref{equu17}. \hfill $\blacksquare$
\end{eeee}
\begin{bbbb} \label{equzforeff} 
Let 
\begin{gather}
2f_{-}(\tau) \! := \! -\mi (a \! - \! \mi/2) \! + \! \frac{\tau}{2} \! \left(\frac{u^{\prime}(\tau) \! - \! \mi b}{u(\tau)} \right), 
\label{hatsoff7} \\
\intertext{and} 
\frac{4 \mi}{\varepsilon b}f_{+}(\tau) \! := \! \mi (a \! + \! \mi/2) \! + \! \frac{\tau}{2} \! \left(\frac{u^{\prime}(\tau) \! 
+ \! \mi b}{u(\tau)} \right). \label{pga3} 
\end{gather}
Then, for $\varepsilon \! \in \! \lbrace \pm 1 \rbrace$,
\begin{gather} 
2f_{-}(\tau) \! = \! \frac{2 \varepsilon \tau^{2}A(\tau)D(\tau)}{u(\tau)} \! = \! \frac{\tau}{2} \frac{\md}{\md \tau} \! 
\left(\ln \! \left(\frac{u(\tau)}{\tau} \right) \! - \! \mi (\varphi (\tau) \! + \! a \ln \tau) \right), \label{hatsoff8} \\
\intertext{and} 
\frac{4 \mi}{\varepsilon b}f_{+}(\tau) \! = \! -\frac{2 \varepsilon \tau^{2}B(\tau)C(\tau)}{u(\tau)} \! = \! \frac{\tau}{2} 
\frac{\md}{\md \tau} \! \left(\ln \! \left(\frac{u(\tau)}{\tau} \right) \! + \! \mi (\varphi (\tau) \! + \! a \ln \tau) \right); 
\label{pga4}
\end{gather}
furthermore,
\begin{equation} \label{pga5} 
\frac{4 \mi}{\varepsilon b}f_{+}(\tau) \! = \! 2f_{-}(\tau) \! + \! \mi \tau \frac{\md}{\md \tau} \! \left(\varphi (\tau) \! 
+ \! a \ln \tau \right) \! = \! 2f_{-}(\tau) \! + \! \mi \tau \! \left(\frac{2a}{\tau} \! + \! \frac{b}{u(\tau)} \right).
\end{equation}
\end{bbbb}

\emph{Proof}. Via definition \eqref{aphnovij}, the system \eqref{equu17}, the corresponding system of isomonodromy 
deformations \eqref{newlax8}, Remark \ref{newlax6}, and the definitions \eqref{equu18} and \eqref{equu19}, one 
establishes the veracity of the relation
\begin{equation} \label{hatsoff9} 
\frac{u^{\prime}(\tau) \! - \! \mi b}{u(\tau)} \! = \! \frac{2}{\tau} \! \left(\frac{2 \tau A(\tau)D(\tau)}{\sqrt{\smash[b]{-A(\tau)
B(\tau)}}} \! + \! (\mi a \! + \! 1/2) \right),
\end{equation}
and then proceeds, \emph{mutatis mutandis}, as in the proof of Proposition \ref{equzforhatf}. The Corollary \eqref{pga5} 
follows {}from, and is consistent with, the definition \eqref{aphnovij} and the first integral of system \eqref{newlax8} (cf. 
Remark~\ref{newlax6}). \hfill $\qed$
\begin{eeee} \label{efftohateff} 
One deduces {}from the definitions \eqref{newlax2}, equation \eqref{equu20}, and Propositions \ref{equzforhatf} 
and \ref{equzforeff} that $f_{\pm}(\tau) \! = \! \hat{f}_{\pm}(\tau)$. \hfill $\blacksquare$
\end{eeee}
\begin{eeee} \label{hamian} 
A lengthy algebraic exercise reveals that, in terms of the coefficient functions $A(\tau)$, $B(\tau)$, $C(\tau)$, and 
$D(\tau)$ satisfying the corresponding isomonodromy deformations \eqref{newlax8}, the Hamiltonian function (cf. 
equation \eqref{eqh1}$)$ reads:
\begin{equation*} 
\mathcal{H}(\tau) \! = \! \frac{1}{2 \tau} \! \left(\mi a \! + \! \dfrac{1}{2} \! + \! \dfrac{2 \tau A(\tau)D(\tau)}{
\sqrt{\smash[b]{-A(\tau)B(\tau)}}} \right)^{2} \! + \! 4 \tau \sqrt{\smash[b]{-A(\tau)B(\tau)}} \! - \! \dfrac{\mi (\varepsilon b)
D(\tau)}{B(\tau)} \! + \! 2 \tau C(\tau)D(\tau) \! + \! \dfrac{A(\tau)D(\tau)}{\sqrt{\smash[b]{-A(\tau)B(\tau)}}}. \tag*{$\blacksquare$}
\end{equation*}
\end{eeee}
\begin{eeee} \label{rem1.1}
Hereafter, all explicit $\tau$ dependencies are suppressed, except where imperative. \hfill $\blacksquare$
\end{eeee}
\subsection{Canonical Solutions and the Monodromy Data} \label{sec1d} 
A succinct discussion of the monodromy data associated with the system \eqref{newlax3} is presented in this 
subsection (see, in particular, \cite{a1,av2}).

For $\mu \! \in \! \mathbb{C}$, the system \eqref{newlax3} has two irregular singular points, one being the point 
at infinity ($\mu \! = \! \infty$) and the other being the origin ($\mu \! = \! 0$). For $\delta_{\infty},\delta_{0} \! > \! 0$ 
and $m \! \in \! \mathbb{Z}$, define the (sectorial) neighbourhoods $\Omega_{m}^{\infty}$ and $\Omega_{m}^{0}$, 
respectively, of these singular points:
\begin{align*}
\Omega_{m}^{\infty} &:= \left\{\mathstrut \mu \! \in \! \mathbb{C}; \, \lvert \mu \rvert \! > \! \delta_{\infty}^{-1}, \, 
\frac{\pi}{2}(m \! - \! 1) \! < \! \arg (\mu) \! + \! \frac{1}{2} \arg (\tau) \! < \! \dfrac{\pi}{2}(m \! + \! 1) \right\}, \\
\Omega_{m}^{0} &:= \left\{\mathstrut \mu \! \in \! \mathbb{C}; \, \lvert \mu \rvert \! < \! \delta_{0}, \, \pi (m \! - \! 
1) \! < \! \arg (\mu) \! - \! \frac{1}{2} \arg (\tau) \! - \! \frac{1}{2} \arg (\varepsilon b) \! < \! \pi (m \! + \! 1) \right\}.
\end{align*}
\begin{bbbb}[{\rm \cite{a1,av2}}] \label{prop1.4} 
There exist solutions $\mathbb{Y}_{m}^{\infty}(\mu) \! = \! \mathbb{Y}_{m}^{\infty}(\mu,\tau)$ and $\mathbb{X}_{m}^{0}
(\mu) \! = \! \mathbb{X}_{m}^{0}(\mu,\tau)$, $m \! \in \! \mathbb{Z}$, of the system \eqref{newlax3} that are uniquely 
defined by the following asymptotic expansions:
\begin{align*}
\mathbb{Y}_{m}^{\infty}(\mu) \underset{\Omega_{m}^{\infty} \ni \mu \to \infty}{:=}& \, \left(\mathrm{I} \! + \! \Psi^{(1)} 
\mu^{-1} \! + \! \Psi^{(2)} \mu^{-2} \! + \! \dotsb \right) \exp \! \left(-\mi \! \left(\tau \mu^{2} \! + \! \left(a \! - \! \mi/2 
\right) \ln \mu \right) \! \sigma_{3} \right), \\
\mathbb{X}_{m}^{0}(\mu) \underset{\Omega_{m}^{0} \ni \mu \to 0}{:=}& \, \Psi_{0} \! \left(\mathrm{I} \! + \! 
\hat{\mathcal{Z}}_{1} \mu \! + \! \dotsb \right) \exp \! \left(-\mi \sqrt{\tau \varepsilon b} \, \mu^{-1} \sigma_{3} \right),
\end{align*}
where $\mathrm{I} \! = \! \diag (1,1)$, $\ln \mu \! := \! \ln \vert \mu \vert \! + \! \mi \arg \mu$,
\begin{gather*}
\Psi^{(1)} \! = \! 
\begin{pmatrix}
0 & \frac{A(\tau)}{\sqrt{\smash[b]{-A(\tau)B(\tau)}}} \\
-\mi D(\tau)/2 & 0
\end{pmatrix}, \, \quad \, \quad \, \Psi^{(2)} \! = \! 
\begin{pmatrix}
\psi^{(2)}_{11} & 0 \\
0 & \psi^{(2)}_{22}
\end{pmatrix}, \\
\psi^{(2)}_{11} \! := \! -\dfrac{\mi}{2} \! \left(\tau \sqrt{\smash[b]{-A(\tau)B(\tau)}}
+ \! \tau C(\tau)D(\tau) \! + \! \dfrac{A(\tau)D(\tau)}{\sqrt{\smash[b]{-A(\tau)
B(\tau)}}} \right), \\
\psi^{(2)}_{22} \! := \! \dfrac{\mi \tau}{2} \! \left(\sqrt{\smash[b]{-A(\tau)B(\tau)}} 
+ \! C(\tau)D(\tau) \right), \\ 
\Psi_{0} \! = \! \dfrac{\mi}{\sqrt{2}} \! \left(\dfrac{(\varepsilon b)^{1/4}}{
\tau^{1/4} \sqrt{\smash[b]{B(\tau)}}} \right)^{\sigma_{3}} \left(\sigma_{1} 
\! + \! \sigma_{3} \right), \, \quad \, \quad \, \hat{\mathcal{Z}}_{1} \! = \! 
\begin{pmatrix}
z_{1}^{(11)} & z_{1}^{(12)} \\
-z_{1}^{(12)} & - z_{1}^{(11)}
\end{pmatrix}, \\
z_{1}^{(11)} \! := -\dfrac{\mi \! \left(\mi a \! + \! \frac{1}{2} \! + \! \frac{2 \tau 
A(\tau)D(\tau)}{\sqrt{\smash[b]{-A(\tau)B(\tau)}}} \right)^{2}}{2 \sqrt{\smash[b]{
\tau \varepsilon b}}} \! - \! \dfrac{2 \mi \tau^{3/2} \sqrt{\smash[b]{-A(\tau)
B(\tau)}}}{\sqrt{\smash[b]{\varepsilon b}}} \! - \! \dfrac{D(\tau) \sqrt{\smash[b]{
\tau \varepsilon b}}}{B(\tau)}, \\
z_{1}^{(12)} \! := \! -\frac{\mi \! \left(\mi a \! + \! \frac{1}{2} \! + \! \frac{2 \tau 
A(\tau)D(\tau)}{\sqrt{\smash[b]{-A(\tau)B(\tau)}}} \right)}{2 \sqrt{\smash[b]{\tau 
\varepsilon b}}},
\end{gather*}
and $\sigma_{1} \! = \! 
\left(
\begin{smallmatrix}
0 & 1 \\
1 & 0
\end{smallmatrix}
\right)$.
\end{bbbb}
\begin{eeee} \label{newrem12} 
The canonical solutions $\mathbb{X}_{m}^{0}(\mu)$, $m \! \in \! \mathbb{Z}$, are defined uniquely provided that 
the branch of $(B(\tau))^{1/2}$ is fixed: hereafter, the branch of $(B(\tau))^{1/2}$ is not fixed; therefore, the set 
of canonical solutions $\lbrace \mathbb{X}_{m}^{0}(\mu) \rbrace_{m \in \mathbb{Z}}$ is defined up to a sign. This 
ambiguity doesn't affect the definition of the Stokes multipliers $($see equations \eqref{eqnstokmult} below); rather, 
it results in a sign discrepancy in the definition of the connection matrix, $G$ (see equation \eqref{eqdefg} below). 
\hfill $\blacksquare$
\end{eeee}

The \emph{canonical solutions}, $\mathbb{Y}_{m}^{\infty}(\mu)$ and $\mathbb{X}_{m}^{0}(\mu)$, $m \! \in \! 
\mathbb{Z}$, enable one to define the \emph{Stokes matrices}, $S_{m}^{\infty}$ and $S_{m}^{0}$, respectively:
\begin{equation}
\mathbb{Y}_{m+1}^{\infty}(\mu) \! = \! \mathbb{Y}_{m}^{\infty}(\mu)S_{m}^{\infty}, \, \quad \, \quad \, 
\mathbb{X}_{m+1}^{0}(\mu) \! = \! \mathbb{X}_{m}^{0}(\mu)S_{m}^{0}. \label{eqnstokmult}
\end{equation}
The Stokes matrices are independent of $\mu$ and $\tau$, and have the following structures:
\begin{equation*}
S_{2m}^{\infty} \! = \! 
\begin{pmatrix}
1 & 0 \\
s_{2m}^{\infty} & 1
\end{pmatrix}, \quad S_{2m+1}^{\infty} \! = \! 
\begin{pmatrix}
1 & s_{2m+1}^{\infty} \\
0 & 1
\end{pmatrix}, \quad S_{2m}^{0} \! = \! 
\begin{pmatrix}
1 & s_{2m}^{0} \\
0 & 1
\end{pmatrix}, \quad S_{2m+1}^{0} \! = \! 
\begin{pmatrix}
1 & 0 \\
s_{2m+1}^{0} & 1
\end{pmatrix}.
\end{equation*}
The parameters $s_{m}^{\infty}$ and $s_{m}^{0}$ are called the \emph{Stokes multipliers}: it can be shown that
\begin{equation} \label{eq1.8} 
S_{m+4}^{\infty} \! = \! \me^{-2 \pi (a-\mi/2) \sigma_{3}}S_{m}^{\infty} 
\me^{2 \pi (a-\mi/2) \sigma_{3}}, \, \quad \, \quad \, S_{m+2}^{0} \! = \! S_{m}^{0}.
\end{equation}
Equations \eqref{eq1.8} imply that the number of independent Stokes multipliers does not exceed six; for example, 
$s_{0}^{0}$, $s_{1}^0$, $s_{0}^{\infty}$, $s_{1}^{\infty}$, $s_{2}^{\infty}$, and $s_{3}^{\infty}$. Furthermore, due to 
the special structure of the system \eqref{newlax3}, that is, the coefficient matrices of odd (resp., even) powers of 
$\mu$ in $\widetilde{\mathscr{U}}(\mu,\tau)$ are diagonal (resp., off-diagonal) and \emph{vice-versa} for 
$\widetilde{\mathscr{V}}(\mu,\tau)$, one can deduce the following relations for the Stokes matrices:
\begin{equation} \label{eq1.9} 
S_{m+2}^{\infty} \! = \! \sigma_{3} \me^{-\pi (a-\mi/2) \sigma_{3}}S_{m}^{\infty} \me^{\pi (a-\mi/2) \sigma_{3}} 
\sigma_{3}, \, \quad \, \quad \, S_{m+1}^{0} = \! \sigma_{1}S_{m}^{0} \sigma_{1}.
\end{equation}
Equations \eqref{eq1.9} reduce the number of independent Stokes multipliers by two, that is, all Stokes multipliers 
can be expressed in terms of $s_{0}^{0}$, $s_{0}^{\infty}$, $s_{1}^{\infty}$, and $a$. There is one more relation 
between the Stokes multipliers that follows {}from the so-called cyclic relation (see equation \eqref{cycrel} below). 
Define the monodromy matrix at the point at infinity, $M^{\infty}$, and the monodromy matrix at the origin, $M^{0}$, 
via the following relations:
\begin{gather*}
\mathbb{Y}_{0}^{\infty}(\mu \me^{-2 \pi \mi}) \! := \! \mathbb{Y}_{0}^{\infty}(\mu)M^{\infty}, \, \quad \, \quad \, 
\mathbb{X}_{0}^{0}(\mu \me^{-2 \pi \mi}) \! := \! \mathbb{X}_{0}^{0}(\mu)M^{0}.
\end{gather*}
Since $\mathbb{Y}_{0}^{\infty}(\mu)$ and $\mathbb{X}_{0}^{0}(\mu)$ are solutions of the system \eqref{newlax3}, 
they differ by a right-hand (matrix) factor $G$:
\begin{equation}
\mathbb{Y}_{0}^{\infty}(\mu) \! := \! \mathbb{X}_{0}^{0}(\mu)G, \label{eqdefg}
\end{equation}
where $G$ is called the \emph{connection matrix}. As matrices relating fundamental solutions of the system 
\eqref{newlax3}, the monodromy, connection, and Stokes matrices are independent of $\mu$ and $\tau$; moreover, 
since $\mathrm{tr}(\widetilde{\mathscr{U}}(\mu,\tau)) \! = \! \mathrm{tr}(\widetilde{\mathscr{V}}(\mu,\tau)) \! = \! 0$, 
it follows that $\det (M^{\infty}) \! = \! \det (M^{0}) \! = \! \det (G) \! = \! 1$. {}From the definition of the monodromy 
and connection matrices, one deduces the following \emph{cyclic relation}:
\begin{equation} \label{cycrel} 
GM^{\infty} \! = \! M^{0}G.
\end{equation} 
The monodromy matrices can be expressed in terms of the Stokes matrices:
\begin{equation*}
M^{\infty} \! = \! S_{0}^{\infty}S_{1}^{\infty}S_{2}^{\infty}S_{3}^{\infty} \me^{-2 \pi (a-\mi/2) \sigma_{3}}, \, 
\quad \, \quad \, M^{0} \! = \! S_{0}^{0}S_{1}^{0}.
\end{equation*}
The Stokes multipliers, $s_{0}^{0}$, $s_{0}^{\infty}$, and $s_{1}^{\infty}$, the elements of the connection matrix, 
$(G)_{ij} \! =: \! g_{ij}$, $i,j \! \in \! \lbrace 1,2 \rbrace$, and the parameter of formal monodromy, $a$, are called the 
\emph{monodromy data}.
\subsection{The Monodromy Manifold, the Direct and Inverse Problems of Monodromy Theory, and Organisation 
of the Paper} \label{sec1e} 
In this subsection, the monodromy manifold is introduced, the direct and inverse problems of monodromy theory are 
discussed (see, for example, \cite{bolisachp,a5,a2,a20}, and Section 2 of \cite{kitp2ellip}), and the contents of this 
work are delineated.

Consider $\mathbb{C}^{8}$ with co-ordinates $(a,s_{0}^{0},s_{0}^{\infty},s_{1}^{\infty},g_{11},g_{12},g_{21},g_{22})$. 
The---algebraic---variety defined by $\det (G) \! = \! 1$ and the \emph{semi-cyclic relation}
\begin{equation} \label{semcyc} 
G^{-1}S_{0}^{0} \sigma_{1}G \! = \! S_{0}^{\infty}S_{1}^{\infty} \sigma_{3} 
\me^{-\pi (a- \mi/2) \sigma_{3}}
\end{equation}
are called the \emph{manifold of the monodromy data}, $\mathscr{M}$. Since only three of the four equations in the 
semi-cyclic relation \eqref{semcyc} are independent, it follows that $\mathrm{dim}_{\mathbb{C}}(\mathscr{M}) \! = \! 
4$; more specifically, the system of algebraic equations defining $\mathscr{M}$ reads:\footnote{In these equations, 
$\me^{\pi a}$ is considered to be a parameter.}
\begin{equation} \label{monoeqns}
\begin{gathered} 
s_{0}^{\infty}s_{1}^{\infty} \! = \! -1 \! - \! \me^{-2 \pi a} \! - \! \mi s_{0}^{0} \me^{-\pi a}, \, \quad \, \quad \, 
g_{21}g_{22} \! - \! g_{11}g_{12} \! + \! s_{0}^{0}g_{11}g_{22} \! = \! \mi \me^{-\pi a}, \\
g_{11}^{2} \! - \! g_{21}^{2} \! - \! s_{0}^{0} g_{11} g_{21} \! = \! \mi s_{0}^{\infty} \me^{-\pi a}, \, \quad \, 
g_{22}^{2} \! - \! g_{12}^{2} \! + \! s_{0}^{0} g_{12} g_{22} \! = \! \mi s_{1}^{\infty} \me^{\pi a}, \, \quad \, 
g_{11}g_{22} \! - \! g_{12} g_{21} \! = \! 1.
\end{gathered}
\end{equation}
\begin{eeee} \label{newrem13}
To achieve a one-to-one correspondence between the coefficients of the system \eqref{newlax3} and the points on 
$\mathscr{M}$, one has to factorize $\mathscr{M}$ by the involution $G \! \to \! -G$ (cf. Remark \ref{newrem12}), 
that is, $G \! \in \! \mathrm{PSL}(2,\mathbb{C})$. \hfill $\blacksquare$
\end{eeee}
As shown in Section 2 of \cite{a1}, equations \eqref{monoeqns} defining $\mathscr{M}$ are equivalent to one of the 
following three systems: \textbf{(i)}\footnote{This case does not exclude the possibility that $g_{12} \! = \! 0$ or 
$g_{21} \! = \! 0$. There is a misprint in Section~2, p. 1172 of \cite{a1}: in item (1), below equations (33), the formula 
for the Stokes multiplier $s_{1}^{\infty}$ must be changed to $s_{1}^{\infty} \! = \! -\tfrac{\mi (g_{22}+ \mi g_{12} 
\me^{-\pi a}) \me^{-\pi a}}{g_{11}}$.} $g_{11}g_{22} \! \neq \! 0$ $\Rightarrow$
\begin{gather} \label{monok1} 
s_{0}^{\infty} \! = \! -\dfrac{(g_{21} \! + \! \mi \me^{\pi a}g_{11})}{g_{22}}, \, \quad \, s_{1}^{\infty} \! = \! 
-\dfrac{\mi (g_{22} \! + \! \mi g_{12} \me^{-\pi a}) \me^{-\pi a}}{g_{11}}, \, \quad \, s_{0}^{0} \! = \! 
\dfrac{\mi \me^{-\pi a} \! + \! g_{11}g_{12} \! - \! g_{21}g_{22}}{g_{11}g_{22}};
\end{gather}
\textbf{(ii)} $g_{11} \! \neq \! 0$ and $g_{22} \! = \! 0$, in which case the parameters are $s_{0}^{0}$ and $g_{11}$, and
\begin{gather} \label{monok2} 
g_{12} \! = \! -\dfrac{\mi \me^{-\pi a}}{g_{11}}, \quad g_{21} \! = \! -\mi \me^{\pi a}
g_{11}, \quad s_{0}^{\infty} \! = \! -\mi g_{11}^{2}(1 \! + \! \me^{2 \pi a} \! + \! 
\mi s_{0}^{0} \me^{\pi a}) \me^{\pi a}, \quad s_{1}^{\infty} \! = \! -\dfrac{\mi 
\me^{-3 \pi a}}{g_{11}^{2}};
\end{gather}
and \textbf{(iii)} $g_{11} \! = \! 0$ and $g_{22} \! \neq \! 0$, in which case the parameters are $s_{0}^{0}$ and $g_{22}$, and
\begin{gather} \label{monok3} 
g_{12} \! = \! \mi \me^{\pi a}g_{22}, \, \quad \, g_{21} \! = \! \frac{\mi \me^{-\pi a}}{g_{22}}, 
\, \quad \, s_{0}^{\infty} \! = \! -\dfrac{\mi \me^{-\pi a}}{g_{22}^{2}}, \, \quad \, 
s_{1}^{\infty} \! = \! -\mi g_{22}^{2}(1 \! + \! \me^{2 \pi a} \! + \! \mi s_{0}^{0} 
\me^{\pi a}) \me^{-\pi a}.
\end{gather}
Asymptotics as $\tau \! \to \! \pm 0$ and $\tau \! \to \! \pm \mi 0$ (resp., as $\tau \! \to \! \pm \infty$ and $\tau \! \to \! 
\pm \mi \infty$) of the general (resp., general regular) solution of the DP3E \eqref{eq1.1}, and its associated 
Hamiltonian function, $\mathcal{H}(\tau)$, parametrised in terms of the proper open subset of $\mathscr{M}$ 
corresponding to case \textbf{(i)} were presented in \cite{a1},\footnote{Asymptotics as $\tau \! \to \! \pm 0$ and 
$\tau \! \to \! \pm \mi 0$ for the corresponding $\pmb{\pmb{\tau}}$-function, but without the `constant term', were 
also conjectured in \cite{a1}.} and asymptotics as $\tau \! \to \! \pm \infty$ and $\tau \! \to \! \pm \mi \infty$ of general 
regular and singular solutions of the DP3E \eqref{eq1.1}, and its associated Hamiltonian and auxiliary functions, 
$\mathcal{H}(\tau)$ and $f_{-}(\tau)$,\footnote{Denoted as $f(\tau)$ in \cite{av2}.} respectively, parametrised in 
terms of the proper open subset of $\mathscr{M}$ corresponding to case \textbf{(i)} were obtained in \cite{av2}; 
furthermore, three-real-parameter families of solutions of the DP3E \eqref{eq1.1} that possess infinite sequences 
of poles and zeros asymptotically located along the imaginary and real axes were identified, and the asymptotic 
distribution of these poles and zeros were also derived. The purpose of the present work, therefore, is to close the 
aforementioned gaps, and to continue to cover $\mathscr{M}$ by deriving asymptotics (as $\tau \! \to \! \pm \infty$ 
and $\tau \! \to \! \pm \mi \infty$) of $u(\tau)$, and the related functions $f_{\pm}(\tau)$, $\mathcal{H}(\tau)$, and 
$\sigma (\tau)$, that are parametrised in terms of the complementary proper open subsets of $\mathscr{M}$ 
corresponding to cases \textbf{(ii)} and \textbf{(iii)}.\footnote{Asymptotics as $\tau \! \to \! \pm 0$ and $\tau \! \to \! 
\pm \mi 0$ for $u(\tau)$, $\mathcal{H}(\tau)$, $f_{\pm}(\tau)$, and $\sigma (\tau)$ corresponding to cases \textbf{(ii)} 
and \textbf{(iii)} will be presented elsewhere.} For notational consistency with the main body of the text, cases 
\textbf{(ii)} and \textbf{(iii)} for $\mathscr{M}$ will, henceforth, be referred to via the integer index $k \! \in \! \lbrace 
\pm 1 \rbrace$; more specifically, case \textbf{(ii)}, that is, $g_{11} \! \neq \! 0$, $g_{22} \! = \! 0$, and $g_{12}g_{21} 
\! = \! -1$, will be designated by $k \! = \! +1$, and case \textbf{(iii)}, that is, $g_{11} \! = \! 0$, $g_{22} \! \neq \! 0$, 
and $g_{12}g_{21} \! = \! -1$, will be designated by $k \! = \! -1$.

Without loss of generality, and with a slight, temporary amendment of the notation, reconsider, for given $a \! \in \! 
\mathbb{C}$, $b \! \in \! \mathbb{R} \setminus \lbrace 0 \rbrace$, and $\varepsilon \! \in \! \lbrace \pm 1 \rbrace$, the 
first-order linear matrix ODE that constitutes the $\mu$-part of the post-gauge-transformed Fuchs-Garnier pair given 
in the system \eqref{newlax3},\footnote{One merely makes the purely notational change $\widetilde{\mathscr{U}}
(\mu,\tau) \! \to \! \widetilde{\mathscr{U}}(\mu,\tau;\vec{\mathbf{y}})$ in equation \eqref{nlxa}. Analogous statements 
can be made regarding the $\mu$-part of the pre-gauge-transformed Fuchs-Garnier pair presented in the system 
\eqref{eqFGmain}.}
\begin{equation} \label{fordmpimp} 
\partial_{\mu} \Psi(\mu,\tau) \! = \! \widetilde{\mathscr{U}}(\mu,\tau;\vec{\mathbf{y}}) \Psi (\mu,\tau),
\end{equation}
where $\mu,\tau \! \in \! \mathbb{C}$, $\mathbb{C}^{5} \! \ni \! \vec{\mathbf{y}} \! := \! (A(\tau),B(\tau),C(\tau),D(\tau),
\sqrtsign{\smash[b]{-A(\tau)B(\tau)}})$ is a vector-valued function constructed {}from the matrix elements of the 
coefficient matrices in the decomposition of (cf. equation \eqref{nlxa}) $\mathrm{M}_{2}(\mathbb{C}) \! \ni \! 
\widetilde{\mathscr{U}}(\mu,\tau;\vec{\mathbf{y}})$ into partial fractions, $\widetilde{\mathscr{U}}(\mu,\tau;\vec{\mathbf{y}})$ 
is a rational function with respect to the spectral parameter $\mu$ with poles that are independent of $\tau$, and $\tr 
(\widetilde{\mathscr{U}}(\mu,\tau;\vec{\mathbf{y}})) \! = \! 0$. The \emph{direct problem of monodromy theory\/} (DMP) 
can be stated as follows: using the tuple of coefficients $(\tau,A(\tau),B(\tau),C(\tau),D(\tau),\sqrtsign{\smash[b]{-A(\tau)
B(\tau)}})$, find the monodromy data $\mathfrak{M} \! := \! (a,s_{0}^{0},s_{0}^{\infty},s_{1}^{\infty},g_{11},g_{12},g_{21},
g_{22}) \! \in \! \mathscr{M}$ (recall that the monodromy data are not independent and are related via the algebraic 
equations \eqref{monoeqns}, which define the complex manifold $\mathscr{M} \! \in \! \mathbb{C}^{8}$ called the 
manifold of the monodromy data), or, in other words, it is a correspondence $(\tau,A(\tau),B(\tau),C(\tau),D(\tau),
\sqrtsign{\smash[b]{-A(\tau)B(\tau)}})$ $\to$ system~\eqref{fordmpimp} $\to$ $\mathfrak{M} \! \in \! \mathscr{M}$. The 
\emph{inverse problem of monodromy theory\/} (IMP) can be stated as follows: using the data set $\lbrace \tau,
\mathfrak{M} \rbrace$, find $\vec{\mathbf{y}} \! \in \! \mathbb{C}^{5}$ such that the system~\eqref{fordmpimp} constructed 
with the help of the co-ordinate (or coefficient) functions of $\vec{\mathbf{y}}$ has the monodromy data $\mathfrak{M} 
\! \in \! \mathscr{M}$, or, in other words, it is the inverse map $\lbrace \tau,\mathfrak{M} \rbrace \! \to \! (\tau,A(\tau),
B(\tau),C(\tau),D(\tau),\sqrtsign{\smash[b]{-A(\tau)B(\tau)}})$.\footnote{If there exists a solution of the IMP, then it is 
unique \cite{bolisachp,a5,a2,a20,kitp2ellip}.} Thus, if one fixes the collection of the monodromy data $\mathfrak{M} 
\! \in \! \mathscr{M}$ and denotes by $\mathscr{T} \subset \mathbb{C}$ the set of all $\tau$ for which the IMP is solvable, 
then the functions $A(\tau),B(\tau),C(\tau),D(\tau),\sqrtsign{\smash[b]{-A(\tau)B(\tau)}} \colon \mathscr{T} \! \to \! 
\mathbb{C}$ are determined, and thus, via Proposition \ref{propuf}, the $2$-tuple $(u(\tau),\varphi (\tau))$ solves the 
system \eqref{equu17}.\footnote{As long as the monodromy data is given, the function $\varphi (\tau)$ is fixed modulo 
$2 \pi l$, $l \! \in \! \mathbb{Z}$, or, alternatively, the constant of integration in the system \eqref{equu17} is defined via 
the monodromy data modulo $2 \pi l$. The function $\varphi (\tau)$ belongs to the class of functions defined by the 
equivalence relation $\varphi \! \equiv \! \varphi \! + \! 2 \pi l$, $l \! \in \! \mathbb{Z}$.} The complete set of the monodromy 
data corresponding to the system \eqref{fordmpimp} (equivalently, the system \eqref{newlax3}) depends, in general, on 
both $\tau$ and $\vec{\mathbf{y}}$, and will be denoted by $\mathfrak{M}(\tau;\vec{\mathbf{y}})$. As a consequence of 
the requirement that the monodromy data be independent of $\tau$ and $\vec{\mathbf{y}}$, that is, $\mathfrak{M}
(\tau;\vec{\mathbf{y}}) \! = \! \mathrm{const.}$, it is necessary that $\vec{\mathbf{y}} \! = \! \vec{\mathbf{y}}(\tau)$ satisfy 
the system of isomonodromy deformations (non-linear ODEs) \eqref{newlax8}, which can be presented in the form 
$\tfrac{\md}{\md \tau} \vec{\mathbf{y}}(\tau) \! = \! \big(-\tfrac{\mi a}{\tau}A(\tau) \! + \! 4C(\tau) 
\sqrtsign{\smash[b]{-A(\tau)B(\tau)}},\tfrac{\mi a}{\tau}B(\tau) \! - \! 4D(\tau) \sqrtsign{\smash[b]{-A(\tau)B(\tau)}},
\tfrac{(\mi a-1)}{\tau}C(\tau) \! - \! 2A(\tau),-\tfrac{(\mi a+1)}{\tau}D(\tau) \! + \! 2B(\tau),2(A(\tau)D(\tau) \! - \! B(\tau)
C(\tau)) \big)$. Clearly, $\mathfrak{M}(\tau;\vec{\mathbf{y}}) \! \in \! \mathscr{M}$. Denote by $\mathbb{M}_{3}$ the 
collection of monodromy data for which the IMP is explicitly solvable: for other $\mathfrak{M}(\tau;\vec{\mathbf{y}}) \! 
\in \! \mathscr{M}$, it is possible to solve the IMP asymptotically (as $\tau \to \! +\infty$, say); this leads to, for example, 
asymptotic formulae for solutions of the DP3E \eqref{eq1.1}. Let $\mathcal{D} \subset \mathscr{M} \setminus 
\mathbb{M}_{3}$ be a domain. The IMP is said to be \emph{asymptotically solvable\/} (as $\tau \! \to \! +\infty$, say) 
if, for any $\mathfrak{M} \! \in \! \mathcal{D}$ representing the monodromy data, there exists an asymptotically locally 
uniform {}\footnote{A function $f(\tau,\lambda)$ is said to be \emph{asymptotically locally uniform\/} (as $\tau \! \to \! 
+\infty$, say) if, for any point $\lambda$ in the domain of definition of $f(\tau,\lambda)$, there exist functions 
$h_{1}(\tau,\lambda)$ and $h_{2}(\tau,\lambda)$ such that, for any $\tilde{\epsilon}_{\ast} \! > \! 0$, there exist 
numbers $T$ and $\tilde{\delta}_{\ast} \! = \! \tilde{\delta}_{\ast}(\lambda,\tilde{\epsilon}_{\ast}) \! > \! 0$ such that, for 
any $(T,+\infty) \! \ni \! \tau$ and for all $\tilde{\lambda} \! \in \! \mathbb{B}_{\tilde{\delta}_{\ast}}(\lambda) \! := \! \lbrace 
\mathstrut \tilde{\lambda}; \, \lvert \tilde{\lambda} \! - \! \lambda \rvert \! < \! \tilde{\delta}_{\ast} \rbrace$ (the open ball of 
radius $\tilde{\delta}_{\ast}$ centred at $\lambda$), the inequality $h_{1}(\tau,\lambda)(1 \! - \! \tilde{\epsilon}_{\ast}) \! 
< \! \lvert f(\tau,\tilde{\lambda}) \rvert \! < \! h_{2}(\tau,\lambda)(1 \! + \! \tilde{\epsilon}_{\ast})$ is satisfied; furthermore, 
if $h_{1}(\tau,\lambda),h_{2}(\tau,\lambda) \! \to \! 0$ (as $\tau \! \to \! +\infty$, say) in the latter inequality, then 
$f(\tau,\lambda)$ is said to be a \emph{locally uniformly decreasing\/} function \cite{a20}.} vector-valued function 
$\vec{\mathbf{y}}^{\scriptscriptstyle \clubsuit} \! = \! \vec{\mathbf{y}}^{\scriptscriptstyle \clubsuit}(\tau;\mathfrak{M}) \! 
:= \! (A(\tau;\mathfrak{M}),B(\tau;\mathfrak{M}),C(\tau;\mathfrak{M}),D(\tau;\mathfrak{M}),\sqrtsign{\smash[b]{-A
(\tau;\mathfrak{M})B(\tau;\mathfrak{M})}}) \! \in \! \mathbb{C}^{5}$ constructed {}from the matrix elements of the 
$\mathrm{M}_{2}(\mathbb{C})$-coefficients of the system~\eqref{fordmpimp} that is analytic in $(T,+\infty) \times 
\mathcal{D}$ and invertible with respect to $\mathfrak{M}$, and the monodromy data 
$\mathfrak{M}^{\scriptscriptstyle \clubsuit}(\tau;\mathfrak{M})$ corresponding to 
$\vec{\mathbf{y}}^{\scriptscriptstyle \clubsuit}(\tau;\mathfrak{M})$ can be represented as 
$\mathfrak{M}^{\scriptscriptstyle \clubsuit}(\tau;\mathfrak{M}) \! = \! \mathfrak{M} \! + \! 
\mathfrak{G}(\tau;\mathfrak{M})$, where $\mathfrak{G}(\tau;\mathfrak{M})$ is a locally uniformly decreasing vector-valued 
function, that is, $\lvert \lvert \mathfrak{M}^{\scriptscriptstyle \clubsuit}(\tau;\mathfrak{M}) \! - \! \mathfrak{M} \rvert \rvert \! 
= \! \lvert \lvert \mathfrak{G}(\tau;\mathfrak{M}) \rvert \rvert \! < \! \mathrm{C} \lvert \tau \rvert^{-\delta_{\ast}}$ as $\tau \! 
\to \! +\infty$,\footnote{$\lvert \lvert \pmb{\cdot} \rvert \rvert$ is any norm in $\mathbb{C}^{8}$.} where $\mathrm{C} \! > \! 0$ 
and $\delta_{\ast} \! > \! 0$ are the same for all $\mathfrak{M}^{\scriptscriptstyle \clubsuit}(\tau;\mathfrak{M})$ 
\cite{a20,kitp2ellip}.\footnote{There are also asymptotics obtained via the IDM for which the vector-valued function(s) 
$\vec{\mathbf{y}}^{\scriptscriptstyle \clubsuit} \! = \! \vec{\mathbf{y}}^{\scriptscriptstyle \clubsuit}(\tau;\mathfrak{M})$ have 
poles for certain $\mathfrak{M} \! \in \! \mathcal{D}$ with $\infty$ (the point at infinity) being an accumulation point of the 
poles (see, for example, \cite{av2}). In such cases, $(T,+\infty)$ must be replaced by $\cup_{m=0}^{\infty}(T_{2m},
T_{2m+1})$, with $T_{m} \! \nearrow \infty$, where the poles lie in the intervals (lacunae) $(T_{2m+1},T_{2m+2})$, and 
where the ratio of the lengths of the intervals containing the poles to the lengths of the intervals devoid of poles must tend 
to zero, that is, $\tfrac{\lvert T_{2m+2}-T_{2m+1} \rvert}{\lvert T_{2m+1}-T_{2m} \rvert} \! \to \! 0$ as $\mathbb{N} \! \ni \! 
m \! \to \! \infty$ (see \cite{a20} for technical details). In such cases, $\cup_{m=0}^{\infty}(T_{2m},T_{2m+1}) \! \times \! 
\mathcal{D}$ should be regarded as the domain of definition for $\vec{\mathbf{y}}^{\scriptscriptstyle \clubsuit}(\tau;\mathfrak{M})$, 
and the IDM enables one to prove the existence of an analytic solution for $\tau \! \in \! \mathbb{C}$ whose asymptotic behaviour 
on $\cup_{m=0}^{\infty}(T_{2m},T_{2m+1})$ is determined by $\vec{\mathbf{y}}^{\scriptscriptstyle \clubsuit}(\tau;\mathfrak{M})$ 
and with poles in the intervals $(T_{2m+1},T_{2m+2})$ \cite{a20}. For complexified $\tau$ with $\lvert \tau \rvert \! \to \! +\infty$, 
$(T,+\infty)$ must be replaced by a Swiss-cheese-like, multiply-connected strip domain (see, for example, \cite{av2}).} In fact, 
according to the \textsc{Theorem} in \cite{a20}, if the IMP is solvable for the domain $\mathcal{D}$, then, for any 
$\mathfrak{M}_{0} \! \in \! \mathcal{D}$ representing the monodromy data for the system \eqref{fordmpimp}, there exists a 
\emph{unique\/} vector-valued function $\vec{\mathbf{y}} \! = \! \vec{\mathbf{y}}(\tau;\mathfrak{M}_{0}) \! := \! (A(\tau;\mathfrak{M}_{0}),
B(\tau;\mathfrak{M}_{0}),C(\tau;\mathfrak{M}_{0}),D(\tau;\mathfrak{M}_{0}),\sqrtsign{\smash[b]{-A(\tau;\mathfrak{M}_{0})
B(\tau;\mathfrak{M}_{0})}}) \! \in \! \mathbb{C}^{5}$ formed by the matrix elements of the $\mathrm{M}_{2}(\mathbb{C})$-coefficients 
of the system~\eqref{fordmpimp} that is analytic in $(T,+\infty) \times \mathcal{D}$ such that the monodromy data $\mathfrak{M}
(\tau;\mathfrak{M}_{0})$ corresponding to $\vec{\mathbf{y}}(\tau;\mathfrak{M}_{0})$ coincides with $\mathfrak{M}_{0}$ for all 
$\tau \! \in \! (T,+\infty)$, namely, $\lvert \lvert \mathfrak{M}(\tau;\mathfrak{M}_{0}) \! - \! \mathfrak{M}_{0} \rvert \rvert \! = \! 
o(\tau^{-\delta_{\ast}})$ uniformly as $\tau \! \to \! +\infty$, $\delta_{\ast} \! > \! 0$.
\begin{eeee} \label{Koefabcdab2uvpup} 
The explication above of the DMP and IMP for the $\mu$-part of the system \eqref{newlax3} was formulated within the 
framework of the $\mathbb{C}$-valued functions $A(\tau)$, $B(\tau)$, $C(\tau)$, $D(\tau)$, and $\sqrtsign{\smash[b]{-A(\tau)
B(\tau)}}$ (solving the system of isomonodromy deformations \eqref{newlax8}) which appear as matrix elements of the 
$\mathrm{M}_{2}(\mathbb{C})$-coefficients of (cf. equation \eqref{nlxa}) $\widetilde{\mathscr{U}}(\mu,\tau)$ in its partial 
fraction decomposition with respect to the spectral parameter $\mu$. Equivalently, via the definition~\eqref{aphnovij}, 
Remark~\ref{newlax6}, and Proposition \ref{propuf}, one may eschew the $\mathbb{C}$-valued functions 
$A(\tau)$, $B(\tau)$, $C(\tau)$, $D(\tau)$, and $\sqrtsign{\smash[b]{-A(\tau)B(\tau)}}$ altogether and re-express 
$\widetilde{\mathscr{U}}(\mu,\tau) \! \in \! \mathrm{M}_{2}(\mathbb{C})$ solely in terms of the $3$-tuple of 
$\mathbb{C}$-valued functions $(u(\tau),\varphi (\tau),u^{\prime}(\tau))$, where, in particular, the $2$-tuple 
$(u(\tau),\varphi (\tau))$ solves the system~\eqref{equu17}, that is,
\begin{align} \label{alteqdmpimp2} 
\widetilde{\mathscr{U}}(\mu,\tau) \! =& -2 \mi \tau \mu \sigma_{3} \! + \! 2 \tau \! 
\begin{pmatrix}
0 & 2 \mi \varepsilon \me^{\mi \varphi (\tau)} \\
\frac{\varepsilon \me^{-\mi \varphi (\tau)}}{2 \tau} \! \left(\mi (a \! - \! \tfrac{\mi}{2}) \! - \! 
\frac{\tau (u^{\prime}(\tau)-\mi b)}{2u(\tau)} \right) & 0
\end{pmatrix} \nonumber \\
-& \, \frac{1}{\mu} \frac{\tau (u^{\prime}(\tau) \! - \! \mi b)}{2u(\tau)} \sigma_{3} \! + \! \dfrac{1}{\mu^{2}} \! 
\begin{pmatrix}
0 & -\frac{\mi \varepsilon b \tau}{u(\tau)} \me^{\mi \varphi (\tau)} \\
-\mi u(\tau) \me^{-\mi \varphi (\tau)} & 0
\end{pmatrix},
\end{align}
and regurgitate \emph{verbatim\/} the above discussion of the DMP and IMP in terms of the $\mathbb{C}$-valued functions 
$u(\tau)$, $\varphi (\tau)$, and $u^{\prime}(\tau)$; but, since the former, and not the latter, approach has been adopted in 
the present work, this matter will not be addressed further. \hfill $\blacksquare$
\end{eeee}
The contents of this paper, the main body of which is devoted to the asymptotic analysis (as $\tau \! \to \! +\infty$ for 
$\varepsilon b \! > \! 0$) of $u(\tau)$ and the related, auxiliary functions $f_{\pm}(\tau)$, $\mathcal{H}(\tau)$, and 
$\sigma (\tau)$, are now described. In Section \ref{sec2}, the main asymptotic results as $\tau \! \to \! \pm \infty$ and 
$\tau \! \to \! \pm \mi \infty$ for $u(\tau)$, $f_{\pm}(\tau)$, $\mathcal{H}(\tau)$, and $\sigma (\tau)$ parametrised in terms 
of the monodromy data corresponding to the cases designated by the index $k \! \in \! \lbrace \pm 1 \rbrace$ (see the 
discussion above) are stated. In Section \ref{sec3}, the asymptotic (as $\tau \! \to \! +\infty$ for $\varepsilon b \! > \! 0$) 
solution of the DMP for the $\mu$-part of the system \eqref{newlax3}, under certain tempered restrictions on its coefficient 
functions (in some class(es) of functions) that are consistent with the monodromy data corresponding to $k \! \in \! \lbrace 
\pm 1 \rbrace$, is presented; in particular, with the coefficient functions satisfying the asymptotic conditions \eqref{iden5}, 
the asymptotic representation for the connection matrix, $G$, corresponding to $k \! \in \! \lbrace \pm 1 \rbrace$ stated in 
Theorem \ref{theor3.3.1} is obtained, and, in conjunction with the parametrisations \eqref{monok2} and \eqref{monok3}, the 
complete asymptotic representation for the monodromy data is derived. The latter analysis is predicated on focusing the 
principal emphasis on the study of the global asymptotic properties of the fundamental solution of the system \eqref{newlax3} 
via the possibility of `matching' different local asymptotic expansions of $\Psi (\mu,\tau)$ at singular and turning points, namely, 
matching WKB-asymptotics of the fundamental solution of the system \eqref{newlax3} with its parametrix represented in terms 
of parabolic-cylinder functions in open neighbourhoods of double-turning points. In Section \ref{finalsec}, the asymptotic results 
derived in Section \ref{sec3} are inverted in order to solve the IMP for the $\mu$-part of the system \eqref{newlax3}, that is, 
explicit asymptotics for the coefficient functions of the $\mu$-part of the system \eqref{newlax3} are parametrised in terms 
of the monodromy data corresponding to $k \! \in \! \lbrace \pm 1 \rbrace$; in particular, via the inversion of the asymptotic 
representation for the connection matrix corresponding to $k \! \in \! \lbrace \pm 1 \rbrace$, explicit asymptotic expressions 
for the coefficient functions parametrised in terms of points on $\mathscr{M}$ are obtained. Under the permanency of the 
isomonodromy condition on the corresponding connection matrices, namely, the monodromy data are constant and satisfy 
certain conditions, one deduces that the asymptotics obtained via inversion represent an asymptotic solution of the IMP and 
satisfy all the restrictions imposed in Section \ref{sec3}; however, since it is not immediately apparent that an asymptotic 
solution of the IMP represents an asymptotic expansion of the functions in the systems \eqref{newlax8} and \eqref{equu17}, 
because the asymptotic solution of the corresponding monodromy problem was obtained via the IDM, one can use the justification 
scheme presented in \cite{a20} (see, also, \cite{bolisachp,a5,a22}) to prove solvability of the corresponding monodromy 
problem, {}from which it follows, therefore, that there exist---exact---solutions of the system of isomonodromy deformations 
\eqref{newlax8} whose asymptotics coincide with those obtained in this section. In order to extend the results derived in 
Sections \ref{sec3} and \ref{finalsec} for asymptotics of $u(\tau)$, $f_{\pm}(\tau)$, $\mathcal{H}(\tau)$, and $\sigma (\tau)$ 
on the positive semi-axis $(\tau \! \to \! +\infty)$ for $\varepsilon b \! > \! 0$ to asymptotics on the negative semi-axis 
$(\tau \! \to \! -\infty)$ and on the imaginary axis $(\tau \! \to \! \pm \mi \infty)$ for both positive and negative values of 
$\varepsilon b$, one applies the (group) action of the Lie-point symmetries changing $\tau \! \to \! -\tau$, $\tau \! \to \! \tau$, 
$a \! \to \! -a$, and $\tau \! \to \! \pm \mi \tau$ derived in Appendix \ref{sectonsymm} on the proper open subsets of 
$\mathscr{M}$ corresponding to $k \! \in \! \lbrace \pm 1 \rbrace$. Finally, in Appendix \ref{feetics}, asymptotics as 
$\tau \! \to \! \pm \infty$ and $\tau \! \to \! \pm \mi \infty$ with $\pm (\varepsilon b) \! > \! 0$ for the multi-valued function 
$\hat{\varphi}(\tau)$ are presented.
\section{Summary of Results} \label{sec2} 
In this work, the detailed analysis of asymptotics as $\tau \! \to \! +\infty$ for $\varepsilon b \! > \! 0$ of $u(\tau)$ and 
the associated functions $f_{\pm}(\tau)$, $\mathcal{H}(\tau)$, $\sigma (\tau)$, and $\hat{\varphi}(\tau)$ is presented. 
In order to arrive at the corresponding asymptotics of $u(\tau)$, $f_{\pm}(\tau)$, $\mathcal{H}(\tau)$, $\sigma (\tau)$, 
and $\hat{\varphi}(\tau)$ for positive, negative, and pure-imaginary values of $\tau$ for both positive and negative values 
of $\varepsilon b$, one applies the actions of the Lie-point symmetries changing $\tau \! \to \! -\tau$, $\tau \! \to \! \tau$, 
$a \! \to \! -a$, and $\tau \! \to \! \pm \mi \tau$ on $\mathscr{M}$ (see Appendices \ref{sectonsymmt}--\ref{sectonsymmtit}, 
respectively). The `composed' symmetries of these actions on $\mathscr{M}$ are presented in Appendix 
\ref{sectonsymmcomp} in terms of two auxiliary mappings, both of which are isomorphisms on $\mathscr{M}$, denoted by 
$\mathscr{F}^{\scriptscriptstyle \lbrace \ell \rbrace}_{\scriptscriptstyle \varepsilon_{1},\varepsilon_{2},m(\varepsilon_{2})}$, 
which is relevant for real $\tau$, and $\hat{\mathscr{F}}^{\scriptscriptstyle \lbrace \hat{\ell} \rbrace}_{\scriptscriptstyle 
\hat{\varepsilon}_{1},\hat{\varepsilon}_{2},\hat{m}(\hat{\varepsilon}_{2})}$, which is relevant for pure-imaginary $\tau$; 
more precisely, {}from Appendix \ref{sectonsymmcomp},\footnote{Due to the involution $G \! \to \! -G$ (cf. Remarks 
\ref{newrem12} and \ref{newrem13}), it suffices to take $\tilde{l} \! = \! l^{\prime} \! = \! +1$ in equations 
\eqref{laxhat76}--\eqref{laxhat121}.}
\begin{align} \label{newpam1} 
\mathscr{F}^{\scriptscriptstyle \lbrace \ell \rbrace}_{\scriptscriptstyle \varepsilon_{1},\varepsilon_{2},m(\varepsilon_{2})} 
\colon \mathscr{M} \! \to \! \mathscr{M}, \, \, &(a,s_{0}^{0},s_{0}^{\infty},s_{1}^{\infty},g_{11},g_{12},g_{21},g_{22}) \! 
\mapsto \! \left((-1)^{\varepsilon_{2}}a,s_{0}^{0}(\varepsilon_{1},\varepsilon_{2},m(\varepsilon_{2}) \vert \ell), \right. 
\nonumber \\
&\left. \, s_{0}^{\infty}(\varepsilon_{1},\varepsilon_{2},m(\varepsilon_{2}) \vert \ell),s_{1}^{\infty}(\varepsilon_{1},
\varepsilon_{2},m(\varepsilon_{2}) \vert \ell),g_{11}(\varepsilon_{1},\varepsilon_{2},m(\varepsilon_{2}) \vert \ell), 
\right. \nonumber \\
&\left. \, g_{12}(\varepsilon_{1},\varepsilon_{2},m(\varepsilon_{2}) \vert \ell),g_{21}(\varepsilon_{1},\varepsilon_{2},
m(\varepsilon_{2}) \vert \ell),g_{22}(\varepsilon_{1},\varepsilon_{2},m(\varepsilon_{2}) \vert \ell) \right),
\end{align}
where $\varepsilon_{1},\varepsilon_{2} \! \in \! \lbrace 0,\pm 1 \rbrace$, 
$m(\varepsilon_{2}) \! = \! 
\left\{
\begin{smallmatrix}
0, \, \, \varepsilon_{2}=0, \\
\pm \varepsilon_{2}, \, \, \varepsilon_{2} \in \lbrace \pm 1 \rbrace,
\end{smallmatrix} 
\right.$ $\ell \! \in \! \lbrace 0,1 \rbrace$, and the explicit expressions for $s_{0}^{0}
(\varepsilon_{1},\varepsilon_{2},m(\varepsilon_{2}) \vert \ell)$, $s_{0}^{\infty}
(\varepsilon_{1},\varepsilon_{2},m(\varepsilon_{2}) \vert \ell)$, $s_{1}^{\infty}
(\varepsilon_{1},\varepsilon_{2},m(\varepsilon_{2}) \vert \ell)$, and $g_{ij}
(\varepsilon_{1},\varepsilon_{2},m(\varepsilon_{2}) \vert \ell)$, $i,j \! \in \! \lbrace 
1,2 \rbrace$, are given in equations~\eqref{laxhat76}--\eqref{laxhat90} 
and~\eqref{laxhat99}--\eqref{laxhat113}, and
\begin{align} \label{newpam3} 
\hat{\mathscr{F}}^{\scriptscriptstyle \lbrace \hat{\ell} \rbrace}_{\scriptscriptstyle 
\hat{\varepsilon}_{1},\hat{\varepsilon}_{2},\hat{m}(\hat{\varepsilon}_{2})} \colon 
\mathscr{M} \! \to \! \mathscr{M}, \, \, &(a,s_{0}^{0},s_{0}^{\infty},s_{1}^{\infty},
g_{11},g_{12},g_{21},g_{22}) \! \mapsto \! \left((-1)^{1+\hat{\varepsilon}_{2}}a,
\hat{s}_{0}^{0}(\hat{\varepsilon}_{1},\hat{\varepsilon}_{2},\hat{m}
(\hat{\varepsilon}_{2}) \vert \hat{\ell}), \right. \nonumber \\
&\left. \, \hat{s}_{0}^{\infty}(\hat{\varepsilon}_{1},\hat{\varepsilon}_{2},\hat{m}
(\hat{\varepsilon}_{2}) \vert \hat{\ell}),\hat{s}_{1}^{\infty}(\hat{\varepsilon}_{1},
\hat{\varepsilon}_{2},\hat{m}(\hat{\varepsilon}_{2}) \vert \hat{\ell}),\hat{g}_{11}
(\hat{\varepsilon}_{1},\hat{\varepsilon}_{2},\hat{m}(\hat{\varepsilon}_{2}) \vert 
\hat{\ell}), \right. \nonumber \\
&\left. \, \hat{g}_{12}(\hat{\varepsilon}_{1},\hat{\varepsilon}_{2},\hat{m}
(\hat{\varepsilon}_{2}) \vert \hat{\ell}),\hat{g}_{21}(\hat{\varepsilon}_{1},
\hat{\varepsilon}_{2},\hat{m}(\hat{\varepsilon}_{2}) \vert \hat{\ell}),\hat{g}_{22}
(\hat{\varepsilon}_{1},\hat{\varepsilon}_{2},\hat{m}(\hat{\varepsilon}_{2}) 
\vert \hat{\ell}) \right),
\end{align}
where $\hat{\varepsilon}_{1} \! \in \! \lbrace \pm 1 \rbrace$, $\hat{\varepsilon}_{2} \! \in \! \lbrace 0,\pm 1 \rbrace$, 
$\hat{m}(\hat{\varepsilon}_{2}) \! = \! 
\left\{
\begin{smallmatrix}
0, \, \, \hat{\varepsilon}_{2} \in \lbrace \pm 1 \rbrace, \\
\pm \hat{\varepsilon}_{1}, \, \, \hat{\varepsilon}_{2}=0,
\end{smallmatrix}
\right.$ $\hat{\ell} \! \in \! \lbrace 0,1 \rbrace$, and the expressions for $\hat{s}_{0}^{0}(\hat{\varepsilon}_{1},
\hat{\varepsilon}_{2},\hat{m}(\hat{\varepsilon}_{2}) \vert \hat{\ell})$, $\hat{s}_{0}^{\infty}(\hat{\varepsilon}_{1},
\hat{\varepsilon}_{2},\hat{m}(\hat{\varepsilon}_{2}) \vert \hat{\ell})$, $\hat{s}_{1}^{\infty}(\hat{\varepsilon}_{1},
\hat{\varepsilon}_{2},\hat{m}(\hat{\varepsilon}_{2}) \vert \hat{\ell})$, and $\hat{g}_{ij}(\hat{\varepsilon}_{1},
\hat{\varepsilon}_{2},\hat{m}(\hat{\varepsilon}_{2}) \vert \hat{\ell})$, $i,j \! \in \! \lbrace 1,2 \rbrace$, are 
given in equations \eqref{laxhat91}--\eqref{laxhat98} and \eqref{laxhat114}--\eqref{laxhat121}.
\begin{eeeee} \label{stokmxx} 
It is worth noting that $s_{0}^{0}(\varepsilon_{1},\varepsilon_{2},m(\varepsilon_{2}) \vert \ell) \! = \! s_{0}^{0} \! 
= \! \hat{s}_{0}^{0}(\hat{\varepsilon}_{1},\hat{\varepsilon}_{2},\hat{m}(\hat{\varepsilon}_{2}) \vert \hat{\ell})$; 
furthermore, it follows that $\operatorname{card} \lbrace (\varepsilon_{1},\varepsilon_{2},m(\varepsilon_{2}) 
\vert \ell) \rbrace \! = \! 30$ and $\operatorname{card} \lbrace (\hat{\varepsilon}_{1},\hat{\varepsilon}_{2},
\hat{m}(\hat{\varepsilon}_{2}) \vert \hat{\ell}) \rbrace \! = \! 16$, that is, for $\ell,\hat{\ell} \! \in \! \lbrace 0,1 \rbrace$,
\begin{equation*} 
(\varepsilon_{1},\varepsilon_{2},m(\varepsilon_{2}) \vert \ell) \! = \! 
\begin{lcase}
(0,0,0 \vert \ell), \\
(-1,0,0 \vert \ell), \\
(1,0,0 \vert \ell), \\
(0,-1,-1 \vert \ell), \\
(0,-1,1 \vert \ell), \\
(0,1,-1 \vert \ell), \\
(0,1,1 \vert \ell), \\
(-1,-1,-1 \vert \ell), \\
(1,-1,-1 \vert \ell), \\
(-1,-1,1 \vert \ell), \\
(1,-1,1 \vert \ell), \\
(-1,1,-1 \vert \ell), \\
(1,1,-1 \vert \ell), \\
(-1,1,1 \vert \ell), \\
(1,1,1 \vert \ell),
\end{lcase} \, \, \qquad \, \, \text{and} \, \, \qquad \, \, (\hat{\varepsilon}_{1},\hat{\varepsilon}_{2},
\hat{m}(\hat{\varepsilon}_{2}) \vert \hat{\ell}) \! = \! 
\begin{lcase}
(1,1,0 \vert \hat{\ell}), \\
(1,-1,0 \vert \hat{\ell}), \\
(-1,1,0 \vert \hat{\ell}), \\
(-1,-1,0 \vert \hat{\ell}), \\
(1,0,-1 \vert \hat{\ell}), \\
(-1,0,-1 \vert \hat{\ell}), \\
(1,0,1 \vert \hat{\ell}), \\
(-1,0,1 \vert \hat{\ell}).
\end{lcase} \tag*{$\blacksquare$}
\end{equation*}
\end{eeeee}
\noindent
Via the above-defined notation(s) and Remark \ref{stokmxx}, asymptotics as $\tau \! \to \! \pm \infty$ (resp., $\tau 
\! \to \! \pm \mi \infty$) for $\pm (\varepsilon b) \! > \! 0$ of $u(\tau)$, $f_{\pm}(\tau)$, $\mathcal{H}(\tau)$, and 
$\sigma (\tau)$ are presented in Theorem \ref{theor2.1} (resp., Theorem \ref{appen}) below, whilst asymptotics as 
$\tau \! \to \! \pm \infty$ (resp., $\tau \! \to \! \pm \mi \infty$) for $\pm (\varepsilon b) \! > \! 0$ of $\hat{\varphi}(\tau)$ 
are presented in Appendix \ref{feetics}, Theorem \ref{pfeetotsa} (resp., Theorem \ref{pfeetotsb}).
\begin{eeeee} \label{rem2.1} 
The roots and fractional powers of positive quantities are assumed positive, whilst the branches of the roots of 
complex quantities can be taken arbitrarily, unless stated otherwise; moreover, it is assumed that, for negative 
real $z$, the following branches are always taken: $z^{1/3} \! := \! -\lvert z \rvert^{1/3}$ and $z^{2/3} \! := \! 
(z^{1/3})^{2}$. \hfill $\blacksquare$
\end{eeeee}
\begin{eeeee} \label{noteaboutreal} 
If one is only interested in the asymptotics as $\tau \! \to \! +\infty$ for $\varepsilon b \! > \! 0$ of the functions 
$u(\tau)$, $f_{\pm}(\tau)$, $\mathcal{H}(\tau)$, and $\sigma (\tau)$, then, in Theorem \ref{theor2.1} below, one 
sets $(\varepsilon_{1},\varepsilon_{2},m(\varepsilon_{2}) \vert \ell) \! = \! (0,0,0 \vert 0)$ and uses the fact that 
(see Appendix \ref{sectonsymmcomp}, the identity map \eqref{laxhat76}) $s_{0}^{0}(0,0,0 \vert 0) \! = \! s_{0}^{0}$, 
$s_{0}^{\infty}(0,0,0 \vert 0) \! = \! s_{0}^{\infty}$, $s_{1}^{\infty}(0,0,0 \vert 0) \! = \! s_{1}^{\infty}$, and $g_{ij}
(0,0,0 \vert 0) \! = \! g_{ij}$, $i,j \! \in \! \lbrace 1,2 \rbrace$. \hfill $\blacksquare$
\end{eeeee}
\begin{ddddd} \label{theor2.1} 
Let $u(\tau)$ be a solution of the {\rm DP3E} \eqref{eq1.1} and $\hat{\varphi}(\tau)$ be the general solution of 
the {\rm ODE} $\hat{\varphi}^{\prime}(\tau) \! = \! 2a \tau^{-1} \! + \! b(u(\tau))^{-1}$ for $\varepsilon b \! > \! 0$ 
corresponding to the monodromy data $(a,s^{0}_{0},s^{\infty}_{0},s^{\infty}_{1},g_{11},g_{12},g_{21},g_{22})$. 
Let $\varepsilon_{1},\varepsilon_{2} \! \in  \! \lbrace 0,\pm 1 \rbrace$, 
$m(\varepsilon_{2}) \! = \! 
\left\{
\begin{smallmatrix}
0, \, \, \varepsilon_{2}=0, \\
\pm \varepsilon_{2}, \, \, \varepsilon_{2} \in \lbrace \pm 1 \rbrace,
\end{smallmatrix} 
\right.$ $\ell \! \in \! \lbrace 0,1 \rbrace$, and $\varepsilon b \! = \! \vert \varepsilon b \vert \me^{\mi \pi 
\varepsilon_{2}}$.\footnote{See Remark \ref{eps2nonzero} below.} For $k \! = \! +1$, let
\begin{equation*}
g_{11}(\varepsilon_{1},\varepsilon_{2},m(\varepsilon_{2}) \vert \ell)g_{12}
(\varepsilon_{1},\varepsilon_{2},m(\varepsilon_{2}) \vert \ell)g_{21}(\varepsilon_{1},
\varepsilon_{2},m(\varepsilon_{2}) \vert \ell) \! \neq \! 0 \quad \text{and} \quad 
g_{22}(\varepsilon_{1},\varepsilon_{2},m(\varepsilon_{2}) \vert \ell) \! = \! 0,
\end{equation*}
and, for $k \! = \! -1$, let
\begin{equation*}
g_{11}(\varepsilon_{1},\varepsilon_{2},m(\varepsilon_{2}) \vert \ell) \! = \! 0 \quad 
\text{and} \quad g_{12}(\varepsilon_{1},\varepsilon_{2},m(\varepsilon_{2}) \vert 
\ell)g_{21}(\varepsilon_{1},\varepsilon_{2},m(\varepsilon_{2}) \vert \ell)g_{22}
(\varepsilon_{1},\varepsilon_{2},m(\varepsilon_{2}) \vert \ell) \! \neq \! 0.
\end{equation*}
Then, for $s_{0}^{0}(\varepsilon_{1},\varepsilon_{2},m(\varepsilon_{2}) \vert \ell) \! \neq \! \mi \me^{(-1)^{1
+\varepsilon_{2}} \pi a}$,\footnote{For $s_{0}^{0}(\varepsilon_{1},\varepsilon_{2},m(\varepsilon_{2}) \vert \ell) \! = \! 
\mi \me^{(-1)^{1+\varepsilon_{2}} \pi a}$, the exponentially small correction terms in the asymptotics \eqref{thmk11}, 
\eqref{thmk17}, \eqref{efhpls1}, \eqref{thmk21}, and \eqref{thmk25} are absent.}
\begin{align} \label{thmk11} 
u(\tau) \underset{\tau \to +\infty \me^{\mi \pi \varepsilon_{1}}}{=} u_{0,k}^{\ast}(\tau) 
- &\dfrac{(-1)^{\varepsilon_{1}} \mi \varepsilon (\varepsilon b \me^{-\mi \pi \varepsilon_{2}})^{1/2} 
\me^{\mi \pi k/4}(s_{0}^{0}(\varepsilon_{1},\varepsilon_{2},m(\varepsilon_{2}) \vert \ell) \! - \! 
\mi \me^{(-1)^{1+\varepsilon_{2}} \pi a})}{\sqrt{\pi} \, 2^{3/2}3^{1/4}(2 \! + \! \sqrt{3})^{\mi 
k(-1)^{1+\varepsilon_{2}}a}} \nonumber \\
\times& \,  \me^{-(\beta (\tau)+\mi k \vartheta (\tau))} \! \left(1 \! + \! \mathcal{O} \big(\tau^{-1/3} \big) \right), \quad 
k \! \in \! \lbrace \pm 1 \rbrace,
\end{align}
where
\begin{equation} \label{thmk1} 
u_{0,k}^{\ast}(\tau) \! = \! c_{0,k} \tau^{1/3} \! \left(1 \! + \! \tau^{-2/3} \sum_{m=0}^{\infty} 
\dfrac{\mathfrak{u}_{m}(k)}{((-1)^{\varepsilon_{1}} \tau^{1/3})^{m}} \right),
\end{equation}
with
\begin{gather}
c_{0,k} \! := \! \frac{\varepsilon (\varepsilon b)^{2/3}}{2} \me^{-\mi 2 \pi k/3}, \label{thmk2} \\
\mathfrak{u}_{0}(k) \! = \! \dfrac{a \me^{-\mi 2 \pi k/3}}{3(\varepsilon b)^{1/3}} 
\! = \! \dfrac{a}{6 \alpha_{k}^{2}}, \, \quad \, \quad \, \mathfrak{u}_{1}(k) \! = \! 
\mathfrak{u}_{2}(k) \! = \! \mathfrak{u}_{3}(k) \! = \! \mathfrak{u}_{5}(k) \! = \! 
\mathfrak{u}_{7}(k) \! = \! \mathfrak{u}_{9}(k) \! = \! 0, \label{thmk3} \\
\mathfrak{u}_{4}(k) \! = \! -\dfrac{a(a^{2} \! + \! 1)}{3^{4}(\varepsilon b)}, 
\qquad \mathfrak{u}_{6}(k) \! = \! \dfrac{a^{2}(a^{2} \! + \! 1) 
\me^{-\mi 2 \pi k/3}}{3^{5}(\varepsilon b)^{4/3}}, \qquad \mathfrak{u}_{8}(k) 
\! = \! \dfrac{a(a^{2} \! + \! 1) \me^{\mi 2 \pi k/3}}{3^{5}(\varepsilon b)^{5/3}}, 
\label{thmk4}
\end{gather}
where
\begin{equation} \label{thmk5} 
\alpha_{k} \! := \! 2^{-1/2}(\varepsilon b)^{1/6} \me^{\mi \pi k/3},
\end{equation}
and, for $m \! \in \! \mathbb{Z}_{+} \! := \! \lbrace 0 \rbrace \cup \mathbb{N}$,
\begin{align} \label{thmk6} 
\mathfrak{u}_{2(m+5)}(k) \! =& \, \dfrac{1}{27} \! \left(\frac{c_{0,k}}{b} \right)^{2} \! 
\left(\vphantom{M^{M^{M^{M^{M}}}}} \mathfrak{w}_{2(m+3)}(k) \! - \! 2 
\mathfrak{u}_{0}(k) \mathfrak{w}_{2(m+2)}(k) \! + \! \eta_{2(m+2)}(k) \! - \! 
\mathfrak{u}_{0}(k) \eta_{2(m+1)}(k) \right. \nonumber \\
+&\left. \, \sum_{p=0}^{2m} \eta_{p}(k) \mathfrak{w}_{2(m+1)-p}(k) \right) \! - \! 
\dfrac{1}{3} \sum_{p=0}^{2(m+4)}(\mathfrak{u}_{p}(k) \! + \! \mathfrak{w}_{p}(k)) 
\mathfrak{u}_{2(m+4)-p}(k) \nonumber \\
-& \, \frac{1}{3} \! \left(\frac{c_{0,k}}{b} \right)^{2} \! \left(\frac{2m \! + \! 7}{3} 
\right)^{2} \mathfrak{u}_{2(m+3)}(k),
\end{align}
\begin{equation} \label{thmk7}
\mathfrak{u}_{2(m+5)+1}(k) \! = \! 0,
\end{equation}
where
\begin{equation} \label{thmk8} 
\mathfrak{w}_{0}(k) \! = \! -\mathfrak{u}_{0}(k), \quad \quad \mathfrak{w}_{1}(k) \! = \! 0, 
\quad \quad \mathfrak{w}_{n+2}(k) \! = \! -\mathfrak{u}_{n+2}(k) \! - \! \sum_{p=0}^{n} 
\mathfrak{w}_{p}(k) \mathfrak{u}_{n-p}(k), \quad n \! \in \! \mathbb{Z}_{+},
\end{equation}
with
\begin{gather} \label{thmk10} 
\eta_{j}(k) \! := \! -2(j \! + \! 3) \mathfrak{u}_{j+2}(k) \! + \! \sum_{p=0}^{j}(p \! + \! 1)
(j \! - \! p \! + \! 1) \mathfrak{u}_{p}(k) \mathfrak{u}_{j-p}(k), \quad j \! \in \! \mathbb{Z}_{+},
\end{gather}
and
\begin{equation} \label{thmk12} 
\vartheta (\tau) \! := \! \dfrac{3 \sqrt{3}}{2}(-1)^{\varepsilon_{2}}(\varepsilon b)^{1/3} \tau^{2/3} \,, 
\qquad \qquad \beta (\tau) \! := \! \dfrac{9}{2}(-1)^{\varepsilon_{2}}(\varepsilon b)^{1/3} \tau^{2/3}.
\end{equation}

Let the auxiliary function $f_{-}(\tau)$ (corresponding to $u(\tau)$ above) defined by equation \eqref{hatsoff7} solve 
the {\rm ODE} \eqref{thmk13}, and let the auxiliary function $f_{+}(\tau)$ (corresponding to $u(\tau)$ above) defined 
by equation \eqref{pga3} solve the {\rm ODE} \eqref{yooplus3}. Then, for $s_{0}^{0}(\varepsilon_{1},\varepsilon_{2},
m(\varepsilon_{2}) \vert \ell) \! \neq \! \mi \me^{(-1)^{1+\varepsilon_{2}} \pi a}$,
\begin{align} \label{thmk17} 
2f_{-}(\tau) \underset{\tau \to +\infty \me^{\mi \pi \varepsilon_{1}}}{=} f_{0,k}^{\ast}(\tau) 
- &\dfrac{(-1)^{\varepsilon_{1}}k(\varepsilon b \me^{-\mi \pi \varepsilon_{2}})^{1/6} 
\me^{\mi \pi k/4} \me^{\mi \pi k/3}(s_{0}^{0}(\varepsilon_{1},\varepsilon_{2},
m(\varepsilon_{2}) \vert \ell) \! - \! \mi \me^{(-1)^{1+\varepsilon_{2}} \pi a})}{\sqrt{\pi} \, 
2^{k/2}3^{1/4}(\sqrt{3} \! + \! 1)^{-k}(2 \! + \! \sqrt{3})^{\mi k(-1)^{1+\varepsilon_{2}}a}} 
\nonumber \\
\times& \, \tau^{1/3} \me^{-(\beta (\tau)+\mi k \vartheta (\tau))} \! \left(1 \! + \! \mathcal{O} \big(\tau^{-1/3} \big) 
\right), \quad k \! \in \! \lbrace \pm 1 \rbrace,
\end{align}
where
\begin{equation} \label{thmk14} 
f_{0,k}^{\ast}(\tau) \! = \! -\mi \! \left((-1)^{\varepsilon_{2}}a \! - \! \mi/2 \right) 
\! + \! \dfrac{\mi (-1)^{\varepsilon_{2}}(\varepsilon b)^{1/3} \me^{\mi 2 \pi k/3}}{2} 
\tau^{2/3} \! \left(-2 \! + \! \tau^{-2/3} \sum_{m=0}^{\infty} \dfrac{\mathfrak{r}_{m}
(k)}{((-1)^{\varepsilon_{1}} \tau^{1/3})^{m}} \right),
\end{equation}
and
\begin{align} \label{efhpls1} 
\frac{4 \mi (-1)^{\varepsilon_{2}}}{\varepsilon b}f_{+}(\tau) \underset{\tau \to +\infty \me^{\mi \pi \varepsilon_{1}}}{=} 
\mathfrak{f}_{0,k}^{\ast}(\tau) + &\dfrac{(-1)^{\varepsilon_{1}}(\varepsilon b \me^{-\mi \pi \varepsilon_{2}})^{1/6} 
\me^{\mi \pi k/4} \me^{\mi \pi k/3}(2^{(k+1)/2} \! - \! k(\sqrt{3} \! + \! 1)^{k})}{\sqrt{\pi} \, 2^{k/2}3^{1/4}(2 \! + \! 
\sqrt{3})^{\mi k(-1)^{1+\varepsilon_{2}}a}} \nonumber \\
\times& \, \big(s_{0}^{0}(\varepsilon_{1},\varepsilon_{2},m(\varepsilon_{2}) \vert \ell) \! - \! \mi 
\me^{(-1)^{1+\varepsilon_{2}} \pi a} \big) \tau^{1/3} \me^{-(\beta (\tau)+\mi k \vartheta (\tau))} \nonumber \\
\times& \left(1 \! + \! \mathcal{O} \big(\tau^{-1/3} \big) \right), \quad k \! \in \! \lbrace \pm 1 \rbrace,
\end{align}
where
\begin{equation} \label{efhpls2} 
\mathfrak{f}_{0,k}^{\ast}(\tau) \! = \! \mi \! \left((-1)^{\varepsilon_{2}}a \! + \! \mi/2 \right) \! + \! \mi (-1)^{\varepsilon_{2}}
(\varepsilon b)^{1/3} \me^{\mi 2 \pi k/3} \tau^{2/3} \! \left(1 \! + \! \tau^{-2/3} \sum_{m=0}^{\infty} \dfrac{(\frac{1}{2} 
\mathfrak{r}_{m}(k) \! + \! 2 \mathfrak{w}_{m}(k))}{((-1)^{\varepsilon_{1}} \tau^{1/3})^{m}} \right),
\end{equation}
with
\begin{equation} \label{thmk15}
\mathfrak{r}_{0}(k) \! = \! \frac{a \! - \! \mi (-1)^{\varepsilon_{2}}/2}{3 \alpha_{k}^{2}}, 
\, \quad \, \mathfrak{r}_{1}(k) \! = \! 0, \, \quad \, \mathfrak{r}_{2}(k) \! = \! \frac{\mi 
(-1)^{\varepsilon_{2}}a(1 \! + \! \mi (-1)^{\varepsilon_{2}}a)}{18 \alpha_{k}^{4}}, \, 
\quad \, \mathfrak{r}_{3}(k) \! = \! 0,
\end{equation}
\begin{align} \label{thmk16} 
\mi 2 \alpha_{k}^{2} \mathfrak{r}_{m+4}(k) \! =& \, \sum_{p=0}^{m} \left(
\mi 4 \alpha_{k}^{2}(\mathfrak{u}_{m+2-p}(k) \! - \! \mathfrak{u}_{0}(k) 
\mathfrak{u}_{m-p}(k)) \! - \! \dfrac{(-1)^{\varepsilon_{2}}}{3}(m \! - \! p \! 
+ \! 2) \mathfrak{u}_{m-p}(k) \right) \! \mathfrak{w}_{p}(k) \nonumber \\
+& \, \mi 4 \alpha_{k}^{2}(\mathfrak{u}_{m+4}(k) \! - \! \mathfrak{u}_{0}(k) 
\mathfrak{u}_{m+2}(k)) \! - \! \dfrac{(-1)^{\varepsilon_{2}}}{3}(m \! + \! 4) 
\mathfrak{u}_{m+2}(k), \quad m \! \in \! \mathbb{Z}_{+}.
\end{align}

Let the Hamiltonian function $\mathcal{H}(\tau)$ (corresponding to $u(\tau)$ above) be defined by equation 
\eqref{eqh1}. Then, for $s_{0}^{0}(\varepsilon_{1},\varepsilon_{2},m(\varepsilon_{2}) \vert \ell) \! \neq \! \mi 
\me^{(-1)^{1+\varepsilon_{2}} \pi a}$,
\begin{align} \label{thmk21} 
\mathcal{H}(\tau) \underset{\tau \to +\infty \me^{\mi \pi \varepsilon_{1}}}{=} \mathcal{H}_{0,k}^{\ast}(\tau) 
- &\dfrac{(-1)^{\varepsilon_{1}}(\varepsilon b \me^{-\mi \pi \varepsilon_{2}})^{1/6} \me^{\mi \pi k/4} \me^{\mi \pi k/3}
(s_{0}^{0}(\varepsilon_{1},\varepsilon_{2},m(\varepsilon_{2}) \vert \ell) \! - \! \mi \me^{(-1)^{1+\varepsilon_{2}} 
\pi a})}{\sqrt{\pi} \, 2^{k/2}3^{3/4}(\sqrt{3} \! + \! 1)^{-k} (2 \! + \! \sqrt{3})^{\mi k(-1)^{1+\varepsilon_{2}}a}} 
\nonumber \\
\times& \, \tau^{-2/3} \me^{-(\beta (\tau)+\mi k \vartheta (\tau))} \! \left(1 \! + \! \mathcal{O} \big(\tau^{-1/3} \big) 
\right), \quad k \! \in \! \lbrace \pm 1 \rbrace,
\end{align}
where
\begin{align} \label{thmk18} 
\mathcal{H}_{0,k}^{\ast}(\tau) \! =& \, 3(\varepsilon b)^{2/3} \me^{-\mi 2 \pi k/3} 
\tau^{1/3} \! + \! 2(\varepsilon b)^{1/3} \me^{\mi 2 \pi k/3}(a \! - \! \mi 
(-1)^{\varepsilon_{2}}/2) \tau^{-1/3} \! + \! \dfrac{1}{6} \! \left((a \! - \! \mi 
(-1)^{\varepsilon_{2}}/2)^{2} \right. \nonumber \\
-&\left. \, 1/3 \right) \! \tau^{-1} \! + \! \alpha_{k}^{2}(\tau^{-1/3})^{5} 
\sum_{m=0}^{\infty} \left(\vphantom{M^{M^{M^{M}}}} \! -4(a \! - \! \mi 
(-1)^{\varepsilon_{2}}/2) \mathfrak{u}_{m+2}(k) \! + \! \alpha_{k}^{2} 
\mathfrak{d}_{m}(k) \right. \nonumber \\
+&\left. \, \sum_{p=0}^{m} \left(\tilde{\mathfrak{h}}_{p}(k) \! - \! 4(a \! - \! \mi (-1)^{\varepsilon_{2}}/2) 
\mathfrak{u}_{p}(k) \right) \! \mathfrak{w}_{m-p}(k) \right) \! \big((-1)^{\varepsilon_{1}} \tau^{-1/3} \big)^{m},
\end{align}
with
\begin{align} \label{thmk19} 
\mathfrak{d}_{m}(k) \! :=& \, \sum_{p=0}^{m+2}(8 \mathfrak{u}_{p}(k) \mathfrak{u}_{m+2-p}(k) \! + \! 
(4 \mathfrak{u}_{p}(k) \! - \! \mathfrak{r}_{p}(k)) \mathfrak{r}_{m+2-p}(k)) \nonumber \\
-& \, \sum_{p_{1}=0}^{m} \sum_{m_{1}=0}^{p_{1}} \mathfrak{r}_{m_{1}}(k) \mathfrak{r}_{p_{1}-m_{1}}(k) 
\mathfrak{u}_{m-p_{1}}(k), \quad m \! \in \! \mathbb{Z}_{+},
\end{align}
and
\begin{equation} \label{thmk20} 
\tilde{\mathfrak{h}}_{0}(k) \! = \! -\dfrac{(12a^{2} \! + \! 1) \me^{\mi \pi k/3}}{18(\varepsilon b)^{1/3}}, \quad \quad 
\tilde{\mathfrak{h}}_{1}(k) \! = \! 0, \quad \quad \tilde{\mathfrak{h}}_{m+2}(k) \! = \! \alpha_{k}^{2} \mathfrak{d}_{m}(k).
\end{equation}

Let the auxiliary function $\sigma (\tau)$ (corresponding to $u(\tau)$ above) defined by equation \eqref{thmk23} 
solve the {\rm ODE} \eqref{thmk22}. Then, for $s_{0}^{0}(\varepsilon_{1},\varepsilon_{2},m(\varepsilon_{2}) \vert \ell) 
\! \neq \! \mi \me^{(-1)^{1+\varepsilon_{2}} \pi a}$,
\begin{align} \label{thmk25} 
\sigma (\tau) \underset{\tau \to +\infty \me^{\mi \pi \varepsilon_{1}}}{=} 
\sigma_{0,k}^{\ast}(\tau) - &\dfrac{(-1)^{\varepsilon_{1}}(\varepsilon b 
\me^{-\mi \pi \varepsilon_{2}})^{1/6} \me^{\mi \pi k/4} \me^{\mi \pi k/3}
(s_{0}^{0}(\varepsilon_{1},\varepsilon_{2},m(\varepsilon_{2}) \vert \ell) \! - \! 
\mi \me^{(-1)^{1+\varepsilon_{2}} \pi a})}{\sqrt{\pi} \, 2^{k/2}3^{3/4}(\sqrt{3} 
\! + \! 1)^{-k}(1 \! + \! k \sqrt{3})^{-1}(2 \! + \! \sqrt{3})^{\mi k(-1)^{1+\varepsilon_{2}}a}} \nonumber \\
\times& \, \tau^{1/3} \me^{-(\beta (\tau)+\mi k \vartheta (\tau))} \! \left(1 \! + \! \mathcal{O} 
\big(\tau^{-1/3} \big) \right), \quad k \! \in \! \lbrace \pm 1 \rbrace,
\end{align}
where
\begin{align} \label{thmk24} 
\sigma_{0,k}^{\ast}(\tau) \! =& \, 3(\varepsilon b)^{2/3} \me^{-\mi 2 \pi k/3} \tau^{4/3} 
\! - \! \mi (-1)^{\varepsilon_{2}}2(\varepsilon b)^{1/3} \me^{\mi 2 \pi k/3}(1 \! + 
\! \mi (-1)^{\varepsilon_{2}}a) \tau^{2/3} \! + \! \dfrac{1}{3} \! \left((1 \! + \! \mi 
(-1)^{\varepsilon_{2}}a)^{2} \right. \nonumber \\
+&\left. 1/3 \right) \! + \! \alpha_{k}^{2} \tau^{-2/3} \sum_{m=0}^{\infty} \left(
-4(a \! - \! \mi (-1)^{\varepsilon_{2}}/2) \mathfrak{u}_{m+2}(k) \! + \! \alpha_{k}^{2} 
\mathfrak{d}_{m}(k) \! + \! \sum_{p=0}^{m}(\tilde{\mathfrak{h}}_{p}(k) \right. \nonumber \\
-&\left. 4(a \! - \! \mi (-1)^{\varepsilon_{2}}/2) \mathfrak{u}_{p}(k)) \mathfrak{w}_{m-p}
(k) \vphantom{M^{M^{M^{M}}}} \! + \! \mi (-1)^{\varepsilon_{2}} \mathfrak{r}_{m+2}(k) \right) \! 
\big((-1)^{\varepsilon_{1}} \tau^{-1/3} \big)^{m}.
\end{align}
\end{ddddd}
\begin{eeeee} \label{eps2nonzero} 
To be unequivocally clear, the first two sentences of the formulation of Theorem \ref{theor2.1} do not imply that 
$\varepsilon_{2} \! = \! 0$ (similar comments apply, \emph{mutatis mutandis\/}, to Theorems \ref{appen}, \ref{pfeetotsa}, 
and \ref{pfeetotsb}). The first sentence of Theorem \ref{theor2.1} states that $u(\tau)$ is a solution of the {\rm DP3E} 
\eqref{eq1.1} and $\hat{\varphi}(\tau)$ is the general solution of the {\rm ODE} $\hat{\varphi}^{\prime}(\tau) \! = \! 2a 
\tau^{-1} \! + \! b(u(\tau))^{-1}$ for $\varepsilon b \! > \! 0$ corresponding to the monodromy data $(a,s^{0}_{0},
s^{\infty}_{0},s^{\infty}_{1},g_{11},g_{12},g_{21},g_{22})$. Taking into account Remarks \ref{stokmxx} and 
\ref{noteaboutreal}, these monodromy co-ordinates are ascribed a clearer notational sense, namely, $s_{0}^{0} \! := \! 
s_{0}^{0}(0,0,0 \vert 0)$, $s_{0}^{\infty} \! := \! s_{0}^{\infty}(0,0,0 \vert 0)$, $s_{1}^{\infty} \! := \! s_{1}^{\infty}(0,0,0 \vert 0)$, 
and $g_{ij} \! := \! g_{ij}(0,0,0 \vert 0)$, $i,j \! \in \! \lbrace 1,2 \rbrace$. This means that one first solves the DP3E \eqref{eq1.1} 
for $u(\tau)$ as $\tau \! \to \! +\infty$ ($\varepsilon_{1} \! = \! 0$) for $\varepsilon b \! > \! 0$ ($\varepsilon_{2} \! = \! 0$) 
corresponding to the monodromy data satisfying the restrictions (take, say, the case $k \! = \! +1$) $g_{11}(0,0,0 \vert 0)
g_{12}(0,0,0 \vert 0)g_{21}(0,0,0 \vert 0) \! = \! g_{11}g_{12}g_{21} \! \neq \! 0$ and $g_{22}(0,0,0 \vert 0) \! = \! g_{22} 
\! = \! 0$, that is,
\begin{equation} \label{tempasu} 
u(\tau) \underset{\tau \to +\infty}{=} u_{0,1}^{\ast}(\tau) \! + \! \frac{\varepsilon (\varepsilon b)^{1/2} 
\me^{-\frac{\mi \pi}{4}}(2 \! + \! \sqrt{3})^{\mi a}(s_{0}^{0} \! - \! \mi \me^{-\pi a})}{\sqrtsign{\pi} \, 2^{3/2}3^{1/4}} \, 
\me^{-(\beta (\tau)+\mi \vartheta (\tau))} \! \left(1 \! + \! \mathcal{O} \big(\tau^{-1/3} \big) \right),
\end{equation}
and then use this $(\varepsilon_{1},\varepsilon_{2},m(\varepsilon_{2}) \vert \ell) \! = \! (0,0,0 \vert 0)$ asymptotics 
\eqref{tempasu} as a ``base'',  ``seed'', or ``germ'' solution to which Lie-point symmetries can be applied (akin to 
Darboux transformations in the theory of solitons); for example, if one wants the solution $u(\tau)$ of the DP3E 
\eqref{eq1.1} as $\tau \! \to \! -\infty$ for $\varepsilon b \! < \! 0$, which corresponds to any one of the parameter 
values $(\varepsilon_{1},\varepsilon_{2},m(\varepsilon_{2}) \vert \ell) \! = \! (-1,1,\pm 1 \vert \ell)$, $\ell \! = \! 0,1$, 
provided that the corresponding monodromy data satisfy the restrictions
\begin{equation*}
g_{22}(-1,1,\pm 1 \vert \ell) \! = \! 0 \quad \text{and} \quad g_{11}(-1,1,\pm 1 \vert \ell)g_{12}(-1,1,\pm 1 \vert \ell)
g_{21}(-1,1,\pm 1 \vert \ell) \! \neq \! 0, \quad \ell \! = \! 0,1,
\end{equation*}
where explicit expressions for $s_{0}^{0}(-1,1,\pm 1 \vert \ell)$, $s_{0}^{\infty}(-1,1,\pm 1 \vert \ell)$, $s_{1}^{\infty}
(-1,1,\pm 1 \vert \ell)$, and $g_{ij}(-1,1,\pm 1 \vert \ell)$, $i,j \! \in \! \lbrace 1,2 \rbrace$, in terms of $s_{0}^{0}$, 
$s_{0}^{\infty}$, $s_{1}^{\infty}$, and $g_{ij}$, are given in Appendix \ref{sectonsymmcomp}, equations 
\eqref{laxhat87}, \eqref{laxhat89}, \eqref{laxhat110}, and \eqref{laxhat112}, one makes the changes $\tau \! \to \! 
\tau \me^{\mi \pi}$ ($\varepsilon_{1} \! = \! -1)$ and $\varepsilon b \! \to \! \lvert \varepsilon b \rvert \me^{\mi \pi}$ 
($\varepsilon_{2} \! = \! +1)$ in equation \eqref{tempasu}, and, taking into account Remark \ref{rem2.1}, arrives 
at the asymptotics of $u(\tau)$ as $\tau \! \to \! -\infty$ for $\varepsilon b \! < \! 0$. \hfill $\blacksquare$
\end{eeeee}
\begin{eeeee} \label{remgrom} 
For $\mi a \! \in \! \mathbb{Z}$, a separate analysis based on B\"{a}cklund transformations is required in order to 
generate the analogue of the sequence of $\mathbb{C}$-valued expansion coefficients $\lbrace \mathfrak{u}_{m}(k) 
\rbrace$, $m \! \in \! \mathbb{Z}_{+}$, $k \! = \! \pm 1$, and the corresponding function $u_{0,k}^{\ast}(\tau)$; this 
comment applies, \emph{mutatis mutandis}, to the $\mathbb{C}$-valued expansion coefficients $\lbrace 
\hat{\mathfrak{u}}_{m}(k) \rbrace$ and the corresponding function $\hat{u}_{0,k}^{\ast}(\tau)$ given in Theorem 
\ref{appen} below (see, also, Theorems \ref{pfeetotsa} and \ref{pfeetotsb}). In fact, as discussed in Section 1 of 
\cite{a1}, for fixed values of $\mi a \! = \! n \! \in \! \mathbb{Z}$, $\varepsilon$, and $b$, there is only one algebraic 
solution (rational function of $\tau^{1/3}$) of the DP3E \eqref{eq1.1} which is a multi-valued function with three 
branches (see, also, \cite{yomura}): this solution can be derived via the $\lvert n \rvert$-fold iteration of the B\"{a}cklund 
transformations given in Subsection 6.1 of \cite{a1} to the simplest solution of the DP3E \eqref{eq1.1} (for $a \! = \! 0$), 
namely, $u(\tau) \! = \! \tfrac{1}{2} \varepsilon (\varepsilon b)^{2/3} \tau^{1/3}$. The case $\mi a \! \in \! \mathbb{Z}$ 
will be considered elsewhere. In this context, it must be mentioned that an expansive analysis, based on the RHP 
approach, of algebraic solutions of the $\mathrm{P} \mathrm{III}$ equation of $\mathrm{D}7$ type has recently 
appeared in \cite{BuckinghamMiller2022}; in particular, the authors present a study of algebraic solutions of the 
DP3E \eqref{eq1.1} for the parameter values $\varepsilon \! = \! -1$, $b \! = \! \mi$, and $a \! = \! -\mi n$, $n \! \in 
\! \mathbb{Z}$. \hfill $\blacksquare$
\end{eeeee}
\begin{eeeee} \label{unboundstrip} 
Define the simply-connected strip domain 
\begin{equation}  \label{hazd68} 
\mathfrak{D}_{u}^{\blacktriangledown} \! := \! \left\lbrace \mathstrut \tau \! \in \! \mathbb{C}; \mathstrut \, 
\Re (\theta^{\ddagger}(\tau)) \! > \! d_{1,\ast}^{\diamond}, \, \lvert \Im (\theta^{\ddagger}(\tau)) \rvert \! < 
\! d_{2,\ast}^{\diamond} \right\rbrace,
\end{equation}
where $\theta^{\ddagger}(\tau) \! = \! 3^{3/2}(-1)^{\varepsilon_{2}}(\varepsilon b)^{1/3} \tau^{2/3}$, and 
$d_{1,\ast}^{\diamond},d_{2,\ast}^{\diamond} \! > \! 0$ are some ($\tau$-independent) constants. The 
asymptotics of $u(\tau)$, $f_{\pm}(\tau)$, $\mathcal{H}(\tau)$, and $\sigma (\tau)$ stated in Theorem 
\ref{theor2.1} are actually valid in $\mathfrak{D}_{u}^{\blacktriangledown}$. \hfill $\blacksquare$
\end{eeeee}
\begin{ddddd} \label{appen} 
Let $u(\tau)$ be a solution of the {\rm DP3E} \eqref{eq1.1} and $\hat{\varphi}(\tau)$ be the general solution of 
the {\rm ODE} $\hat{\varphi}^{\prime}(\tau) \! = \! 2a \tau^{-1} \! + \! b(u(\tau))^{-1}$ for $\varepsilon b \! > \! 0$ 
corresponding to the monodromy data $(a,s^{0}_{0},s^{\infty}_{0},s^{\infty}_{1},g_{11},g_{12},g_{21},g_{22})$. 
Let $\hat{\varepsilon}_{1} \! \in \! \lbrace \pm 1 \rbrace$, $\hat{\varepsilon}_{2} \! \in  \! \lbrace 0,\pm 1 \rbrace$, 
$\hat{m}(\hat{\varepsilon}_{2}) \! = \! 
\left\{
\begin{smallmatrix}
0, \, \, \hat{\varepsilon}_{2} \in \lbrace \pm 1 \rbrace, \\
\pm \hat{\varepsilon}_{1}, \, \, \hat{\varepsilon}_{2}=0,
\end{smallmatrix}
\right.$ $\hat{\ell} \! \in \! \lbrace 0,1 \rbrace$, and $\varepsilon b \! = \! \vert \varepsilon b \vert 
\me^{\mi \pi \hat{\varepsilon}_{2}}$. For $k \! = \! +1$, let
\begin{equation*}
\hat{g}_{11}(\hat{\varepsilon}_{1},\hat{\varepsilon}_{2},\hat{m}(\hat{\varepsilon}_{2}) 
\vert \hat{\ell}) \hat{g}_{12}(\hat{\varepsilon}_{1},\hat{\varepsilon}_{2},\hat{m}
(\hat{\varepsilon}_{2}) \vert \hat{\ell}) \hat{g}_{21}(\hat{\varepsilon}_{1},
\hat{\varepsilon}_{2},\hat{m}(\hat{\varepsilon}_{2}) \vert \hat{\ell}) \! \neq \! 0 \quad 
\text{and} \quad \hat{g}_{22}(\hat{\varepsilon}_{1},\hat{\varepsilon}_{2},\hat{m}
(\hat{\varepsilon}_{2}) \vert \hat{\ell}) \! = \! 0,
\end{equation*}
and, for $k \! = \! -1$, let
\begin{equation*}
\hat{g}_{11}(\hat{\varepsilon}_{1},\hat{\varepsilon}_{2},\hat{m}(\hat{\varepsilon}_{2}) 
\vert \hat{\ell}) \! = \! 0 \quad \text{and} \quad \hat{g}_{12}(\hat{\varepsilon}_{1},
\hat{\varepsilon}_{2},\hat{m}(\hat{\varepsilon}_{2}) \vert \hat{\ell}) \hat{g}_{21}
(\hat{\varepsilon}_{1},\hat{\varepsilon}_{2},\hat{m}(\hat{\varepsilon}_{2}) \vert 
\hat{\ell}) \hat{g}_{22}(\hat{\varepsilon}_{1},\hat{\varepsilon}_{2},\hat{m}
(\hat{\varepsilon}_{2}) \vert \hat{\ell}) \! \neq \! 0.
\end{equation*}
Then, for $\hat{s}_{0}^{0}(\hat{\varepsilon}_{1},\hat{\varepsilon}_{2},\hat{m}(\hat{\varepsilon}_{2}) \vert \hat{\ell}) \! \neq 
\! \mi \me^{(-1)^{\hat{\varepsilon}_{2}} \pi a}$,\footnote{For $\hat{s}_{0}^{0}(\hat{\varepsilon}_{1},\hat{\varepsilon}_{2},
\hat{m}(\hat{\varepsilon}_{2}) \vert \hat{\ell}) \! = \! \mi \me^{(-1)^{\hat{\varepsilon}_{2}} \pi a}$, the exponentially small 
correction terms in the asymptotics \eqref{appen9}, \eqref{appen14}, \eqref{efhpls3}, \eqref{appen18}, and \eqref{appen20} 
are absent.}
\begin{align} \label{appen9} 
u(\tau) \underset{\tau \to +\infty \me^{\mi \pi \hat{\varepsilon}_{1}/2}}{=} 
\hat{u}_{0,k}^{\ast}(\tau) - &\dfrac{\mi \me^{-\mi \pi \hat{\varepsilon}_{1}/2} 
\varepsilon (\varepsilon b \me^{-\mi \pi \hat{\varepsilon}_{2}})^{1/2} 
\me^{\mi \pi k/4}(\hat{s}_{0}^{0}(\hat{\varepsilon}_{1},\hat{\varepsilon}_{2},\hat{m}
(\hat{\varepsilon}_{2}) \vert \hat{\ell}) \! - \! \mi \me^{(-1)^{\hat{\varepsilon}_{2}} 
\pi a})}{\sqrt{\pi} \, 2^{3/2}3^{1/4}(2 \! + \! \sqrt{3})^{\mi k(-1)^{\hat{\varepsilon}_{2}}a}} 
\nonumber \\
\times& \, \me^{-(\hat{\beta}(\tau_{\ast})+\mi k \hat{\vartheta}(\tau_{\ast}))} \! \left(1 \! + \! \mathcal{O} 
\big(\tau^{-1/3} \big) \right), \quad k \! \in \! \lbrace \pm 1 \rbrace,
\end{align}
where
\begin{equation} \label{appen1} 
\hat{u}_{0,k}^{\ast}(\tau) \! = \! \me^{-\mi \pi \hat{\varepsilon}_{1}/2}c_{0,k} \tau_{\ast}^{1/3} \! \left(1 \! + \! 
\tau_{\ast}^{-2/3} \sum_{m=0}^{\infty} \dfrac{\hat{\mathfrak{u}}_{m}(k)}{\big(\tau_{\ast}^{1/3} \big)^{m}} \right),
\end{equation}
with $c_{0,k}$ defined by equation \eqref{thmk2},
\begin{gather}
\tau_{\ast} \! := \! \tau \me^{-\mi \pi \hat{\varepsilon}_{1}/2}, \label{teas} \\
\hat{\mathfrak{u}}_{0}(k) \! = \! -\dfrac{a \me^{-\mi 2 \pi k/3}}{3(\varepsilon 
b)^{1/3}} \! = \! -\dfrac{a}{6 \alpha_{k}^{2}}, \qquad \hat{\mathfrak{u}}_{1}
(k) \! = \! \hat{\mathfrak{u}}_{2}(k) \! = \! \hat{\mathfrak{u}}_{3}(k) \! 
= \! \hat{\mathfrak{u}}_{5}(k) \! = \! \hat{\mathfrak{u}}_{7}(k) \! = \! 
\hat{\mathfrak{u}}_{9}(k) \! = \! 0, \label{appen2} \\
\hat{\mathfrak{u}}_{4}(k) \! = \! \dfrac{a(a^{2} \! + \! 1)}{3^{4}(\varepsilon b)}, 
\qquad \hat{\mathfrak{u}}_{6}(k) \! = \! \dfrac{a^{2}(a^{2} \! + \! 1) 
\me^{-\mi 2 \pi k/3}}{3^{5}(\varepsilon b)^{4/3}}, \qquad \hat{\mathfrak{u}}_{8}
(k) \! = \! -\dfrac{a(a^{2} \! + \! 1) \me^{\mi 2 \pi k/3}}{3^{5}(\varepsilon b)^{5/3}}, 
\label{appen3}
\end{gather}
where $\alpha_{k}$ is defined by equation \eqref{thmk5}, and, for $m \! \in \! \mathbb{Z}_{+}$, 
\begin{align} \label{appen4} 
\hat{\mathfrak{u}}_{2(m+5)}(k) \! =& \, \frac{1}{27} \! \left(\frac{c_{0,k}}{b} 
\right)^{2} \! \left(\vphantom{M^{M^{M^{M}}}} \hat{\mathfrak{w}}_{2(m+3)}(k) 
\! - \! 2 \hat{\mathfrak{u}}_{0}(k) \hat{\mathfrak{w}}_{2(m+2)}(k) \! + \! 
\hat{\eta}_{2(m+2)}(k) \! - \! \hat{\mathfrak{u}}_{0}(k) \hat{\eta}_{2(m+1)}(k) 
\right. \nonumber\\
+&\left. \, \sum_{p=0}^{2m} \hat{\eta}_{p}(k) \hat{\mathfrak{w}}_{2(m+1)-p}(k) 
\right) \! - \! \frac{1}{3} \sum_{p=0}^{2(m+4)}(\hat{\mathfrak{u}}_{p}(k) \! 
+ \! \hat{\mathfrak{w}}_{p}(k)) \hat{\mathfrak{u}}_{2(m+4)-p}(k)\nonumber \\
-& \, \frac{1}{3} \! \left(\frac{c_{0,k}}{b} \right)^{2} \! \left(\dfrac{2m \! + \! 7}{3} 
\right)^{2} \hat{\mathfrak{u}}_{2(m+3)}(k),
\end{align}
\begin{equation} \label{appen5}
\hat{\mathfrak{u}}_{2(m+5)+1}(k) \! = \! 0,
\end{equation}
where 
\begin{equation} \label{appen6} 
\hat{\mathfrak{w}}_{0}(k) \! = \! -\hat{\mathfrak{u}}_{0}(k), \quad \quad 
\hat{\mathfrak{w}}_{1}(k) \! = \! 0, \quad \quad \hat{\mathfrak{w}}_{n+2}(k) \! 
= \! -\hat{\mathfrak{u}}_{n+2}(k) \! - \! \sum_{p=0}^{n} \hat{\mathfrak{w}}_{p}
(k) \hat{\mathfrak{u}}_{n-p}(k), \quad n \! \in \! \mathbb{Z}_{+},
\end{equation}
with
\begin{gather}
\hat{\eta}_{j}(k) \! := \! -2(j \! + \! 3) \hat{\mathfrak{u}}_{j+2}(k) \! + \! \sum_{p=0}^{j}
(p \! + \! 1)(j \! - \! p \! + \! 1) \hat{\mathfrak{u}}_{p}(k) \hat{\mathfrak{u}}_{j-p}(k), 
\quad j \! \in \! \mathbb{Z}_{+}, \label{appen8}
\end{gather}
and
\begin{equation} \label{thmk12q} 
\hat{\vartheta}(\tau) \! := \! \dfrac{3 \sqrt{3}}{2}(-1)^{\hat{\varepsilon}_{2}}(\varepsilon b)^{1/3} 
\tau^{2/3} \,, \qquad \qquad \hat{\beta}(\tau) \! := \! \dfrac{9}{2}(-1)^{\hat{\varepsilon}_{2}}
(\varepsilon b)^{1/3} \tau^{2/3}.
\end{equation}

Let the auxiliary function $f_{-}(\tau)$ (corresponding to $u(\tau)$ above) defined by equation \eqref{hatsoff7} solve the 
{\rm ODE} \eqref{thmk13}, and let the auxiliary function $f_{+}(\tau)$ (corresponding to $u(\tau)$ above) defined by 
equation \eqref{pga3} solve the {\rm ODE}~\eqref{yooplus3}. Then, for $\hat{s}_{0}^{0}(\hat{\varepsilon}_{1},
\hat{\varepsilon}_{2},\hat{m}(\hat{\varepsilon}_{2}) \vert \hat{\ell}) \! \neq \! \mi \me^{(-1)^{\hat{\varepsilon}_{2}} \pi a}$,
\begin{align} \label{appen14} 
2f_{-}(\tau) \underset{\tau \to +\infty \me^{\mi \pi \hat{\varepsilon}_{1}/2}}{=} 
\hat{f}_{0,k}^{\ast}(\tau) - &\dfrac{k(\varepsilon b \me^{-\mi \pi 
\hat{\varepsilon}_{2}})^{1/6} \me^{\mi \pi k/4} \me^{\mi \pi k/3}
(\hat{s}_{0}^{0}(\hat{\varepsilon}_{1},\hat{\varepsilon}_{2},\hat{m}
(\hat{\varepsilon}_{2}) \vert \hat{\ell}) \! - \! \mi \me^{(-1)^{
\hat{\varepsilon}_{2}} \pi a})}{\sqrt{\pi} \, 2^{k/2}3^{1/4}(\sqrt{3} \! + \! 1)^{-k}
(2 \! + \! \sqrt{3})^{\mi k(-1)^{\hat{\varepsilon}_{2}}a}} \nonumber \\
\times& \, \tau_{\ast}^{1/3} \me^{-(\hat{\beta}(\tau_{\ast})+\mi k \hat{\vartheta}(\tau_{\ast}))} \! 
\left(1 \! + \! \mathcal{O} \big(\tau^{-1/3} \big) \right), \quad k \! \in \! \lbrace \pm 1 \rbrace,
\end{align}
where
\begin{equation} \label{appen11} 
\hat{f}_{0,k}^{\ast}(\tau) \! = \! -\mi \! \left((-1)^{1+\hat{\varepsilon}_{2}}a \! - \! \mi/2 \right) \! + \! 
\dfrac{\mi (-1)^{\hat{\varepsilon}_{2}}(\varepsilon b)^{1/3} \me^{\mi 2 \pi k/3}}{2} \tau_{\ast}^{2/3} \! \left(-2 \! + 
\! \tau_{\ast}^{-2/3} \sum_{m=0}^{\infty} \dfrac{\hat{\mathfrak{r}}_{m}(k)}{\big(\tau_{\ast}^{1/3} \big)^{m}} \right),
\end{equation}
and
\begin{align} \label{efhpls3} 
\frac{4 \mi (-1)^{\hat{\varepsilon}_{2}}}{\varepsilon b}f_{+}(\tau) \underset{\tau \to 
+\infty \me^{\mi \pi \hat{\varepsilon}_{1}/2}}{=} \hat{\mathfrak{f}}_{0,k}^{\ast}(\tau) + 
&\dfrac{(\varepsilon b \me^{-\mi \pi \hat{\varepsilon}_{2}})^{1/6} \me^{\mi \pi k/4} 
\me^{\mi \pi k/3}(2^{(k+1)/2} \! - \! k(\sqrt{3} \! + \! 1)^{k})}{\sqrt{\pi} \, 2^{k/2}
3^{1/4}(2 \! + \! \sqrt{3})^{\mi k(-1)^{\hat{\varepsilon}_{2}}a}} \nonumber \\
\times& \, \big(\hat{s}_{0}^{0}(\hat{\varepsilon}_{1},\hat{\varepsilon}_{2},\hat{m}
(\hat{\varepsilon}_{2}) \vert \hat{\ell}) \! - \! \mi \me^{(-1)^{\hat{\varepsilon}_{2}} \pi a} \big) 
\tau_{\ast}^{1/3} \me^{-(\hat{\beta}(\tau_{\ast})+\mi k \hat{\vartheta}(\tau_{\ast}))} \nonumber \\
\times& \, \left(1 \! + \! \mathcal{O} \big(\tau^{-1/3} \big) \right), \quad k \! \in \! \lbrace \pm 1 \rbrace,
\end{align}
where
\begin{equation} \label{efhpls4} 
\hat{\mathfrak{f}}_{0,k}^{\ast}(\tau) \! = \! \mi \! \left((-1)^{1+\hat{\varepsilon}_{2}}a \! + \! \mi/2 \right) \! + \! \mi 
(-1)^{\hat{\varepsilon}_{2}}(\varepsilon b)^{1/3} \me^{\mi 2 \pi k/3} \tau_{\ast}^{2/3} \! \left(1 \! + \! \tau_{\ast}^{-2/3} 
\sum_{m=0}^{\infty} \dfrac{(\frac{1}{2} \hat{\mathfrak{r}}_{m}(k) \! + \! 2 \hat{\mathfrak{w}}_{m}(k))}{\big(
\tau_{\ast}^{1/3} \big)^{m}} \right),
\end{equation}
with
\begin{equation} \label{appen12}
\hat{\mathfrak{r}}_{0}(k) \! = \! -\dfrac{(a \! + \! \mi (-1)^{\hat{\varepsilon}_{2}}/2)}{
3 \alpha_{k}^{2}}, \, \quad \, \hat{\mathfrak{r}}_{1}(k) \! = \! 0, \, \quad \, 
\hat{\mathfrak{r}}_{2}(k) \! = \! \dfrac{\mi a((-1)^{1+\hat{\varepsilon}_{2}} \! 
+ \! \mi a)}{18 \alpha_{k}^{4}}, \, \quad \, \hat{\mathfrak{r}}_{3}(k) \! = \! 0,
\end{equation}
\begin{align} \label{appen13} 
\mi 2 \alpha_{k}^{2} \hat{\mathfrak{r}}_{m+4}(k) \! =& \, \sum_{p=0}^{m} \left(
\mi 4 \alpha_{k}^{2}(\hat{\mathfrak{u}}_{m+2-p}(k) \! - \! \hat{\mathfrak{u}}_{0}
(k) \hat{\mathfrak{u}}_{m-p}(k)) \! - \! \dfrac{(-1)^{\hat{\varepsilon}_{2}}}{3}(m \! - \! 
p \! + \! 2) \hat{\mathfrak{u}}_{m-p}(k) \right) \! \hat{\mathfrak{w}}_{p}(k) \nonumber \\
+& \, \mi 4 \alpha_{k}^{2}(\hat{\mathfrak{u}}_{m+4}(k) \! - \! \hat{\mathfrak{u}}_{0}
(k) \hat{\mathfrak{u}}_{m+2}(k)) \! - \! \dfrac{(-1)^{\hat{\varepsilon}_{2}}}{3}(m \! + 
\! 4) \hat{\mathfrak{u}}_{m+2}(k), \quad m \! \in \! \mathbb{Z}_{+}.
\end{align}

Let the Hamiltonian function $\mathcal{H}(\tau)$ (corresponding to $u(\tau)$ above) be defined by equation 
\eqref{eqh1}. Then, for $\hat{s}_{0}^{0}(\hat{\varepsilon}_{1},\hat{\varepsilon}_{2},\hat{m}(\hat{\varepsilon}_{2}) 
\vert \hat{\ell}) \! \neq \! \mi \me^{(-1)^{\hat{\varepsilon}_{2}} \pi a}$,
\begin{align} \label{appen18} 
\mathcal{H}(\tau) \underset{\tau \to +\infty \me^{\mi \pi \hat{\varepsilon}_{1}/2}}{=} 
\hat{\mathcal{H}}_{0,k}^{\ast}(\tau) - &\dfrac{\me^{-\mi \pi \hat{\varepsilon}_{1}/2}
(\varepsilon b \me^{-\mi \pi \hat{\varepsilon}_{2}})^{1/6} \me^{\mi \pi k/4} 
\me^{\mi \pi k/3}(\hat{s}_{0}^{0}(\hat{\varepsilon}_{1},\hat{\varepsilon}_{2},\hat{m}
(\hat{\varepsilon}_{2}) \vert \hat{\ell}) \! - \! \mi \me^{(-1)^{\hat{\varepsilon}_{2}} 
\pi a})}{\sqrt{\pi} \, 2^{k/2}3^{3/4}(\sqrt{3} \! + \! 1)^{-k}(2 \! + \! \sqrt{3})^{\mi k
(-1)^{\hat{\varepsilon}_{2}}a}} \nonumber \\
\times& \, \tau_{\ast}^{-2/3} \me^{-(\hat{\beta}(\tau_{\ast})+\mi k \hat{\vartheta}(\tau_{\ast}))} 
\! \left(1 \! + \! \mathcal{O} \big(\tau^{-1/3} \big) \right), \quad k \! \in \! \lbrace \pm 1 \rbrace,
\end{align}
where
\begin{align} \label{appen15} 
\hat{\mathcal{H}}_{0,k}^{\ast}(\tau) \! =& \, \me^{-\mi \pi \hat{\varepsilon}_{1}/2} 
\! \left(3(\varepsilon b)^{2/3} \me^{-\mi 2 \pi k/3} \tau_{\ast}^{1/3} \! + \! 
(-1)^{\hat{\varepsilon}_{2}}2(\varepsilon b)^{1/3} \me^{\mi 2 \pi k/3} \! \left(
(-1)^{1+\hat{\varepsilon}_{2}}a \! - \! \mi/2 \right) \! \tau_{\ast}^{-1/3} \right. 
\nonumber \\
+&\left. \, \dfrac{1}{6} \! \left(((-1)^{1+\hat{\varepsilon}_{2}}a \! - \! \mi/2)^{2} \! - 
\! 1/3 \right) \! \tau_{\ast}^{-1} \! + \! (-1)^{\hat{\varepsilon}_{2}} \alpha_{k}^{2}
(\tau_{\ast}^{-1/3})^{5} \sum_{m=0}^{\infty} \left(-4((-1)^{1+\hat{\varepsilon}_{2}}
a \! - \! \mi/2) \right. \right. \nonumber \\
\times&\left. \left. \hat{\mathfrak{u}}_{m+2}(k) \! + \! (-1)^{\hat{\varepsilon}_{2}} 
\alpha_{k}^{2} \hat{\mathfrak{d}}_{m}(k) \! + \! \sum_{p=0}^{m} \left(\hat{\mathfrak{h}}_{p}^{\ast}(k) \! - \! 
4((-1)^{1+\hat{\varepsilon}_{2}} a \! - \! \mi/2) \hat{\mathfrak{u}}_{p}(k) \right) \! \hat{\mathfrak{w}}_{m-p}(k) 
\right) \! \big(\tau_{\ast}^{-1/3} \big)^{m} \right),
\end{align}
with
\begin{align} \label{appen16} 
\hat{\mathfrak{d}}_{m}(k) \! :=& \, \sum_{p=0}^{m+2}(8 \hat{\mathfrak{u}}_{p}
(k) \hat{\mathfrak{u}}_{m+2-p}(k) \! + \! (4 \hat{\mathfrak{u}}_{p}(k) \! - \! 
\hat{\mathfrak{r}}_{p}(k)) \hat{\mathfrak{r}}_{m+2-p}(k)) \nonumber \\
-& \, \sum_{p_{1}=0}^{m} \sum_{m_{1}=0}^{p_{1}} \hat{\mathfrak{r}}_{m_{1}}(k) 
\hat{\mathfrak{r}}_{p_{1}-m_{1}}(k) \hat{\mathfrak{u}}_{m-p_{1}}(k), \quad m \! 
\in \! \mathbb{Z}_{+},
\end{align}
and
\begin{equation} \label{appen17} 
\hat{\mathfrak{h}}_{0}^{\ast}(k) \! = \! \dfrac{(-1)^{1+\hat{\varepsilon}_{2}}
(12a^{2} \! + \! 1) \me^{\mi \pi k/3}}{18(\varepsilon b)^{1/3}}, 
\quad \quad \hat{\mathfrak{h}}_{1}^{\ast}(k) \! = \! 0, \quad \quad 
\hat{\mathfrak{h}}_{m+2}^{\ast}(k) \! = \! (-1)^{\hat{\varepsilon}_{2}} 
\alpha_{k}^{2} \hat{\mathfrak{d}}_{m}(k).
\end{equation}

Let the auxiliary function $\sigma (\tau)$ (corresponding to $u(\tau)$ above) defined by equation \eqref{thmk23} 
solve the {\rm ODE} \eqref{thmk22}. Then, for $\hat{s}_{0}^{0}(\hat{\varepsilon}_{1},\hat{\varepsilon}_{2},\hat{m}
(\hat{\varepsilon}_{2}) \vert \hat{\ell}) \! \neq \! \mi 
\me^{(-1)^{\hat{\varepsilon}_{2}} \pi a}$,
\begin{align} \label{appen20} 
\sigma (\tau) \underset{\tau \to +\infty \me^{\mi \pi \hat{\varepsilon}_{1}/2}}{=} 
\hat{\sigma}_{0,k}^{\ast}(\tau) - &\dfrac{(\varepsilon b \me^{-\mi \pi 
\hat{\varepsilon}_{2}})^{1/6} \me^{\mi \pi k/4} \me^{\mi \pi k/3}(\hat{s}_{0}^{0}
(\hat{\varepsilon}_{1},\hat{\varepsilon}_{2},\hat{m}(\hat{\varepsilon}_{2}) \vert 
\hat{\ell}) \! - \! \mi \me^{(-1)^{\hat{\varepsilon}_{2}} \pi a})}{\sqrt{\pi} \, 2^{k/2}
3^{3/4}(\sqrt{3} \! + \! 1)^{-k}(1 \! + \! k \sqrt{3})^{-1}(2 \! + \! \sqrt{3})^{\mi 
k(-1)^{\hat{\varepsilon}_{2}}a}} \nonumber \\
\times& \, \tau_{\ast}^{1/3} \me^{-(\hat{\beta}(\tau_{\ast})+\mi k \hat{\vartheta}(\tau_{\ast}))} 
\! \left(1 \! + \! \mathcal{O} \big(\tau^{-1/3} \big) \right), \quad k \! \in \! \lbrace \pm 1 \rbrace,
\end{align}
where
\begin{align} \label{appen19} 
\hat{\sigma}_{0,k}^{\ast}(\tau) \! =& \, 3(\varepsilon b)^{2/3} \me^{-\mi 2 \pi k/3} \tau_{\ast}^{4/3} \! - \! \mi 
(-1)^{\hat{\varepsilon}_{2}}2(\varepsilon b)^{1/3} \me^{\mi 2 \pi k/3}(1 \! + \! \mi (-1)^{1+\hat{\varepsilon}_{2}}a) 
\tau_{\ast}^{2/3} \nonumber \\
+& \, \dfrac{1}{3} \! \left((1 \! + \! \mi (-1)^{1+\hat{\varepsilon}_{2}}a)^{2} \! + \! 1/3 \right) \! + \! 
(-1)^{\hat{\varepsilon}_{2}} \alpha_{k}^{2} \tau_{\ast}^{-2/3} \sum_{m=0}^{\infty} \left(-4((-1)^{1+
\hat{\varepsilon}_{2}}a \! - \! \mi/2) \hat{\mathfrak{u}}_{m+2}(k) \right. \nonumber \\
+&\left. (-1)^{\hat{\varepsilon}_{2}} \alpha_{k}^{2} \hat{\mathfrak{d}}_{m}(k) \! + \! \sum_{p=0}^{m} \left(
\hat{\mathfrak{h}}_{p}^{\ast}(k) \! - \! 4((-1)^{1+\hat{\varepsilon}_{2}}a \! - \! \mi/2) \hat{\mathfrak{u}}_{p}(k) 
\right) \! \hat{\mathfrak{w}}_{m-p}(k) \right. \nonumber \\
+&\left. \mi \hat{\mathfrak{r}}_{m+2}(k) \vphantom{M^{M^{M^{f}}}} \right) \! \big(\tau_{\ast}^{-1/3} \big)^{m}.
\end{align}
\end{ddddd}
\begin{eeeee} \label{unboundstripcompl} 
Define the simply-connected strip domain
\begin{gather} 
\widehat{\mathfrak{D}}_{u}^{\blacktriangle} \! := \! \left\lbrace \mathstrut \tau \! \in \! \mathbb{C}; \mathstrut 
\, \Re (\hat{\theta}^{\ddagger}(\tau \me^{-\mi \pi \hat{\varepsilon}_{1}/2})) \! > \! \hat{d}_{1,\ast}^{\diamond}, 
\, \lvert \Im (\hat{\theta}^{\ddagger}(\tau \me^{-\mi \pi \hat{\varepsilon}_{1}/2})) \rvert \! < \! 
\hat{d}_{2,\ast}^{\diamond} \right\rbrace, \label{hazd668}
\end{gather}
where $\hat{\theta}^{\ddagger}(\tau) \! = \! 3^{3/2}(-1)^{\hat{\varepsilon}_{2}}(\varepsilon b)^{1/3} \tau^{2/3}$, 
and $\hat{d}_{1,\ast}^{\diamond},\hat{d}_{2,\ast}^{\diamond} \! > \! 0$ are some ($\tau$-independent) constants. 
The asymptotics of $u(\tau)$, $f_{\pm}(\tau)$, $\mathcal{H}(\tau)$, and $\sigma (\tau)$ stated in Theorem 
\ref{appen} are actually valid in $\widehat{\mathfrak{D}}_{u}^{\blacktriangle}$. \hfill $\blacksquare$
\end{eeeee}
\section{Asymptotic Solution of the Direct Problem of Monodromy Theory} \label{sec3} 
In this section, the monodromy data introduced in Subsection~\ref{sec1d} is calculated as $\tau \! \to \! +\infty$ 
for $\varepsilon b \! > \! 0$ (corresponding to $(\varepsilon_{1},\varepsilon_{2},m(\varepsilon_{2}) \vert \ell) \! = 
\! (0,0,0 \vert 0)$; cf. Section \ref{sec2}): this constitutes the first step towards the proof of the results stated in 
Theorems \ref{theor2.1}, \ref{appen}, \ref{pfeetotsa}, and \ref{pfeetotsb}. 

The aforementioned calculation consists of three components: (i) the matrix WKB analysis for the $\mu$-part of the 
system \eqref{newlax3}, that is,
\begin{equation} \label{eq3.1} 
\partial_{\mu} \Psi (\mu) \! = \! \widetilde{\mathscr{U}}(\mu,\tau) \Psi (\mu),
\end{equation}
where $\Psi (\mu) \! = \! \Psi (\mu,\tau)$ (see Subsection \ref{subsec3.1}); (ii) the approximation of $\Psi (\mu)$ in the 
neighbourhoods of the turning points (see Subsection \ref{sec3.2}); and (iii) the matching of these asymptotics (see 
Subsection \ref{sec3.3}).

Before commencing the asymptotic analysis, the notation used throughout this work is introduced:
\begin{enumerate}
\item[(1)] $\mathrm{I} \! = \! \diag (1,1)$ is the $2 \times 2$ identity 
matrix, $\sigma_{1} \! = \! 
\left(
\begin{smallmatrix}
0 & 1 \\
1 & 0
\end{smallmatrix}
\right)$, $\sigma_{2} \! = \! 
\left(
\begin{smallmatrix}
0 & -\mi \\
\mi & 0
\end{smallmatrix}
\right)$, and $\sigma_{3} \! = \! 
\left(
\begin{smallmatrix}
1 & 0 \\
0 & -1
\end{smallmatrix}
\right)$ are the Pauli matrices, $\sigma_{\pm} \! := \! \tfrac{1}{2}(\sigma_{1} \! \pm \! \mi \sigma_{2})$, 
$\mathbb{R}_{\pm} \! := \! \lbrace \mathstrut x \! \in \! \mathbb{R}; \, \pm x \! > \! 0 \rbrace$, and 
$\mathbb{C}_{\pm} \! := \! \lbrace \mathstrut z \! \in \! \mathbb{C}; \, \pm \Im (z) \! > \! 0 \rbrace$;
\item[(2)] for $(\varsigma_{1},\varsigma_{2}) \! \in \! \mathbb{R} \times\mathbb{R}$, the function 
$(z \! - \! \varsigma_{1})^{\mi \varsigma_{2}} \colon \mathbb{C} \setminus (-\infty,\varsigma_{1}] \! \to \! 
\mathbb{C}$, $z \! \mapsto \! \exp (\mi \varsigma_{2} \ln (z \! - \! \varsigma_{1}))$, with the branch cut 
taken along $(-\infty,\varsigma_{1}]$ and the principal branch of the logarithm chosen (that is, $\arg 
(z \! - \! \varsigma_{1}) \! \in \! (-\pi,\pi])$;
\item[(3)] for $\omega_{o} \! \in \! \mathbb{C}$ and $\widehat{\Upsilon} \! \in \! \mathrm{M}_{2}
(\mathbb{C})$, $\omega_{o}^{\ad (\sigma_{3})} \widehat{\Upsilon } \! := \! \omega_{o}^{\sigma_{3}} 
\widehat{\Upsilon} \omega_{o}^{-\sigma_{3}}$;
\item[(4)] for $\mathrm{M}_{2}(\mathbb{C}) \! \ni \! \mathfrak{I}(z)$, $(\mathfrak{I}(z))_{ij}$ or 
$\mathfrak{I}_{ij}(z)$, $i,j \! \in \! \lbrace 1,2 \rbrace$, denotes the $(i \, j)$-element of $\mathfrak{I}(z)$;
\item[(5)] $\hat{w}(t) \genfrac{}{}{0pt}{3}{=}{t \to +\infty} o(1)$ means there exists $C_{1} \! > \! 0$ and 
$\epsilon_{1} \! > \! 0$ such that $\lvert \hat{w}(t) \rvert \! \leqslant \! C_{1} \lvert t \rvert^{-\epsilon_{1}}$;
\item[(6)] for $\mathrm{M}_{2}(\mathbb{C}) \! \ni \! \hat{\mathfrak{Y}}(z)$, $\hat{\mathfrak{Y}}(z) 
\genfrac{}{}{0pt}{3}{=}{z \to z_{0}} \mathcal{O}(\pmb{\ast})$ (resp., $o(\pmb{\ast}))$ means 
$\hat{\mathfrak{Y}}_{ij}(z) \genfrac{}{}{0pt}{3}{=}{z \to z_{0}} \mathcal{O}(\pmb{\ast}_{ij})$ (resp., 
$o(\pmb{\ast}_{ij}))$, $i,j \! \in \! \lbrace 1,2 \rbrace$;
\item[(7)] for $\mathrm{M}_{2}(\mathbb{C}) \! \ni \! \hat{\mathfrak{B}}(z)$, $\lvert \lvert \hat{\mathfrak{B}}(\pmb{\cdot}) 
\rvert \rvert \! := \! \left(\sum_{i,j=1}^{2} \hat{\mathfrak{B}}_{ij}(\pmb{\cdot}) \overline{\hat{\mathfrak{B}}_{ij}
(\pmb{\cdot})} \, \right)^{1/2}$ denotes the Hilbert-Schmidt norm, where $\overline{\pmb{\star}}$ denotes 
complex conjugation of $\pmb{\star}$; and
\item[(8)] for some $\delta_{\ast} \! > \! 0$, $\mathscr{O}_{\delta_{\ast}}(z_{0})$ denotes the (open) 
$\delta_{\ast}$-neighbourhood of the point $z_{0}$, that is, for $z_{0} \! \in \! \mathbb{C}$, 
$\mathscr{O}_{\delta_{\ast}}(z_{0}) \! := \! \lbrace \mathstrut z \! \in \! \mathbb{C}; \, \lvert z \! - \! z_{0} 
\rvert \! < \! \delta_{\ast} \rbrace$, and, for $z_{0}$ the point at infinity, $\mathscr{O}_{\delta_{\ast}}
(\infty) \! := \! \lbrace \mathstrut z \! \in \! \mathbb{C}; \, \lvert z \rvert \! > \! \delta_{\ast}^{-1} \rbrace$.
\end{enumerate}
\subsection{Matrix WKB Analysis} \label{subsec3.1} 
This subsection is devoted to the WKB analysis of equation \eqref{eq3.1} as $\tau \! \to \! +\infty$ for 
$\varepsilon b \! > \! 0$.

In order to transform equation \eqref{eq3.1} into a form amenable to WKB analysis, one uses the result of Proposition 
4.1.1 in \cite{a1} (see, also, Proposition 3.2.1 in \cite{av2}), which is summarised here for the reader's convenience.
\begin{bbbb}[\textrm{\cite{a1,av2}}] \label{prop3.1.1} 
In the system \eqref{newlax3}, let
\begin{equation} \label{eq3.2}
\begin{gathered}
A(\tau) \! = \! a(\tau) \tau^{-2/3}, \, \quad \, B(\tau) \! = \! b(\tau) \tau^{-2/3}, 
\, \quad \, C(\tau) \! = \! c(\tau) \tau^{-1/3}, \, \quad \, D(\tau) \! = \! d(\tau) 
\tau^{-1/3}, \\
\widetilde{\mu} \! = \! \mu \tau^{1/6}, \, \quad \, \quad \, \widetilde{\Psi}
(\widetilde{\mu}) \! := \! \tau^{-\frac{1}{12} \sigma_{3}} \Psi (\widetilde{\mu} 
\tau^{-1/6}),
\end{gathered}
\end{equation}
where $\widetilde{\Psi}(\widetilde{\mu}) \! = \! \widetilde{\Psi}(\widetilde{\mu},\tau)$. Then, the $\mu$-part of the system 
\eqref{newlax3} transforms as follows:
\begin{equation} \label{eq3.3}
\partial_{\widetilde{\mu}} \widetilde{\Psi}(\widetilde{\mu}) \! = \! \tau^{2/3} 
\mathcal{A}(\widetilde{\mu},\tau) \widetilde{\Psi}(\widetilde{\mu}),
\end{equation}
where
\begin{equation} \label{eq3.4} 
\mathcal{A}(\widetilde{\mu},\tau) \! := \! -2 \mi \widetilde{\mu} \sigma_{3} \! + \! 
\begin{pmatrix}
0 & -\frac{4 \mi \sqrt{\smash[b]{-a(\tau)b(\tau)}}}{b(\tau)} \\
-2d(\tau) & 0
\end{pmatrix} \! - \! \dfrac{1}{\widetilde{\mu}} \dfrac{\mi r(\tau)(\varepsilon 
b)^{1/3}}{2} \sigma_{3} \! + \! \dfrac{1}{\widetilde{\mu}^{2}} \! 
\begin{pmatrix}
0 & \frac{\mi (\varepsilon b)}{b(\tau)} \\
\mi b(\tau) & 0
\end{pmatrix},
\end{equation}
with
\begin{equation} \label{eq3.5} 
\dfrac{\mi r(\tau)(\varepsilon b)^{1/3}}{2} \! = \! \mi (a \! - \! \mi/2) \tau^{-2/3} \! + \! 
\dfrac{2a(\tau)d(\tau)}{\sqrt{\smash[b]{-a(\tau)b(\tau)}}}.
\end{equation}
\end{bbbb}

As in Subsection 3.2 of \cite{av2}, define the functions $h_{0}(\tau)$, $\hat{r}_{0}(\tau)$, and $\hat{u}_{0}
(\tau)$ via the relations
\begin{gather}
\sqrt{\smash[b]{-a(\tau)b(\tau)}} + \! c(\tau)d(\tau) \! + \! \dfrac{a(\tau)d(\tau) 
\tau^{-2/3}}{2 \sqrt{\smash[b]{-a(\tau)b(\tau)}}} \! - \! \dfrac{1}{4}(a \! - \! \mi/2)^{2} 
\tau^{-4/3} \! = \! \dfrac{3}{4}(\varepsilon b)^{2/3} \! - \! h_{0}(\tau) \tau^{-2/3}, 
\label{iden1} \\
r(\tau) \! = \! -2 \! + \! \hat{r}_{0}(\tau), \label{iden3oldr} \\
\sqrt{\smash[b]{-a(\tau)b(\tau)}} = \! \dfrac{(\varepsilon b)^{2/3}}{2}
(1 \! + \! \hat{u}_{0}(\tau)). \label{iden4oldu}
\end{gather}
As follows {}from the first integral \eqref{aphnovij} (cf. Remark \ref{newlax6}), the functions $a(\tau)$, $b(\tau)$, 
$c(\tau)$, and $d(\tau)$ are related via the formula
\begin{equation} \label{iden6} 
a(\tau)d(\tau) \! + \! b(\tau)c(\tau) \! + \! \mi a \sqrt{\smash[b]{-a(\tau)b(\tau)}} 
\tau^{-2/3} \! = \! -\mi \varepsilon b/2, \quad \varepsilon \! \in \! \lbrace \pm 
1 \rbrace.
\end{equation}
It is worth noting that equations \eqref{iden1}--\eqref{iden6} are self-consistent; in fact, a calculation reveals 
that they are equivalent to
\begin{align}
a(\tau)d(\tau) \! =& \, \dfrac{(\varepsilon b)^{2/3}}{2}(1 \! + \! \hat{u}_{0}(\tau)) \! 
\left(-\dfrac{\mi (\varepsilon b)^{1/3}}{2} \! + \! \dfrac{\mi (\varepsilon b)^{1/3} 
\hat{r}_{0}(\tau)}{4} \! - \! \dfrac{\mi}{2}(a \! - \! \mi/2) \tau^{-2/3} \right), 
\label{iden7} \\
b(\tau)c(\tau) \! =& \, \dfrac{(\varepsilon b)^{2/3}}{2}(1 \! + \! \hat{u}_{0}(\tau)) 
\! \left(-\dfrac{\mi (\varepsilon b)^{1/3}}{2} \! + \! \mi (\varepsilon b)^{1/3} \! 
\left(\dfrac{\hat{u}_{0}(\tau)}{1 \! + \! \hat{u}_{0}(\tau)} \! - \! \dfrac{\hat{r}_{0}
(\tau)}{4} \right) \! - \! \dfrac{\mi}{2}(a \! + \! \mi/2) \tau^{-2/3} \right), 
\label{iden8} \\
-h_{0}(\tau) \tau^{-2/3} \! =& \, \dfrac{(\varepsilon b)^{2/3}}{2} \! \left(
\dfrac{(\hat{u}_{0}(\tau))^{2} \! + \! \frac{1}{2} \hat{u}_{0}(\tau) \hat{r}_{0}
(\tau)}{1 \! + \! \hat{u}_{0}(\tau)} \! - \! \dfrac{(\hat{r}_{0}(\tau))^{2}}{8} \right) 
\! + \! \dfrac{(\varepsilon b)^{1/3}(a \! - \! \mi/2) \tau^{-2/3}}{2(1 \! + \! 
\hat{u}_{0}(\tau))}; \label{iden9}
\end{align}
moreover, via equations \eqref{iden4oldu}, \eqref{iden7}, and \eqref{iden8}, one shows that
\begin{align}
-c(\tau)d(\tau) =& \, \left(\dfrac{\mi (\varepsilon b)^{1/3}}{2} \! - \! \mi 
(\varepsilon b)^{1/3} \! \left(\dfrac{\hat{u}_{0}(\tau)}{1 \! + \! \hat{u}_{0}
(\tau)} \! - \! \dfrac{\hat{r}_{0}(\tau)}{4} \right) \! + \! \dfrac{\mi}{2}
(a \! + \! \mi/2) \tau^{-2/3} \right) \nonumber \\
\times& \, \left(\dfrac{\mi (\varepsilon b)^{1/3}}{2} \! - \! \dfrac{\mi 
(\varepsilon b)^{1/3} \hat{r}_{0}(\tau)}{4} \! + \! \dfrac{\mi}{2}
(a \! - \! \mi/2) \tau^{-2/3} \right). \label{iden10}
\end{align}

In this work, in lieu of the functions $h_{0}(\tau)$, $\hat{r}_{0}(\tau)$, and $\hat{u}_{0}(\tau)$, it is more 
convenient to work with the functions $\hat{h}_{0}(\tau)$, $\tilde{r}_{0}(\tau)$, and $v_{0}(\tau)$, respectively, 
which are defined as follows: for $k \! = \! \pm 1$,
\begin{gather}
h_{0}(\tau) \! := \! \left(\dfrac{3(\varepsilon b)^{2/3}}{4} \! \left(1 \! - \! \me^{-\mi 2 \pi k/3} \right) 
\! + \! \hat{h}_{0}(\tau) \right) \! \tau^{2/3}, \label{iden2} \\
-2 \! + \! \hat{r}_{0}(\tau) \! := \! \me^{\mi 2 \pi k/3} \! \left(-2 \! + \! \tilde{r}_{0}(\tau) \tau^{-1/3} \right), 
\label{iden3} \\
1 \! + \! \hat{u}_{0}(\tau) \! := \! \me^{-\mi 2 \pi k/3} \! \left(1 \! + \! v_{0}(\tau) \tau^{-1/3} \right). \label{iden4}
\end{gather}
The WKB analysis of equation \eqref{eq3.3} is predicated on the assumption that the functions $\hat{h}_{0}
(\tau)$, $\tilde{r}_{0}(\tau)$, and $v_{0}(\tau)$ satisfy the---asymptotic---conditions
\begin{equation}
\begin{gathered}
\lvert \hat{h}_{0}(\tau) \rvert \underset{\tau \to +\infty}{=} \mathcal{O}(\tau^{-2/3}), \qquad \lvert \tilde{r}_{0}
(\tau) \rvert \underset{\tau \to +\infty}{=} \mathcal{O}(\tau^{-1/3}), \qquad \lvert v_{0}(\tau) \rvert 
\underset{\tau \to +\infty}{=} \mathcal{O}(\tau^{-1/3}). \label{iden5}
\end{gathered}
\end{equation}
\begin{eeee} \label{rempowlzerowz} 
Some solutions $u(\tau)$ of the DP3E \eqref{eq1.1} may, and in fact do, have poles and zeros located 
on the positive real line. In order to be able to study such solutions, one must consider a slightly more 
general, complex domain $\widetilde{\mathfrak{D}}_{u}$; however, since, \emph{a priori\/}, one does not 
know the solutions $u(\tau)$ which possess such poles and zeros, nor their exact locations, it is necessary 
to introduce a formal definition for $\widetilde{\mathfrak{D}}_{u}$. Denote by $\mathscr{P}_{u}$ and 
$\mathscr{Z}_{u}$, respectively, the countable sets of poles and zeros of the function $u(\tau)$. As a 
consequence of the Painlev\'{e} property, these sets may have accumulation points at the origin and at 
the point at infinity. Define neighbourhoods of $\mathscr{P}_{u}$ and $\mathscr{Z}_{u}$, respectively, as 
follows:\footnote{There is a misprint in Subsection 3.1 of \cite{av2}: in the definitions (3.2) and (3.3), the 
inequality $>$ must be changed to $<$.} for some $\epsilon_{\ast} \! > \! 0$, let
\begin{gather*} 
\mathscr{P}_{u}(\epsilon_{\ast}) \! := \! \left\{\mathstrut \tau \! \in \! \mathbb{C}; \mathstrut \, \lvert 
\theta^{\ddagger}(\tau) \! - \! \theta^{\ddagger}(\tau_{\mathfrak{p}}) \rvert \! < \! C_{\ast} \lvert 
\tau_{\mathfrak{p}} \rvert^{-\epsilon_{\ast}}, \, \tau_{\mathfrak{p}} \! \in \! \mathscr{P}_{u} \right\}, \\
\mathscr{Z}_{u}(\epsilon_{\ast}) \! := \! \left\{\mathstrut \tau \! \in \! \mathbb{C}; \mathstrut \, \lvert 
\theta^{\ddagger}(\tau) \! - \! \theta^{\ddagger}(\tau_{\mathfrak{z}}) \rvert \! < \! C_{\ast} \lvert 
\tau_{\mathfrak{z}} \rvert^{-\epsilon_{\ast}}, \, \tau_{\mathfrak{z}} \! \in \! \mathscr{Z}_{u} \right\},
\end{gather*}
where $\theta^{\ddagger}(\tau)$ is given in Remark~\ref{unboundstrip}, and $C_{\ast} \! > \! 0$ is some 
($\tau$-independent) constant. Now, define the Swiss-cheese-like, multiply-connected domain 
$\widetilde{\mathfrak{D}}_{u} \! := \! \mathfrak{D}_{u}^{\blacktriangledown} \setminus (\mathscr{P}_{u}
(\epsilon_{\ast}) \cup \mathscr{Z}_{u}(\epsilon_{\ast}))$, where the simply-connected strip domain 
$\mathfrak{D}_{u}^{\blacktriangledown}$ is defined by equation \eqref{hazd68}. Theoretically speaking, 
therefore, it is to be understood that the asymptotic analysis is undertaken in the sense that 
$\widetilde{\mathfrak{D}}_{u} \! \ni \! \tau$ and $\Re (\tau) \! \to \! +\infty$ (with $\varepsilon b \! > \! 0$); 
however, due to the---asymptotic---conditions \eqref{iden5}, which reflect the sought-after class(es) of functions 
analysed herein, it turns out that $\mathscr{P}_{u}(\epsilon_{\ast}) \! = \! \mathscr{Z}_{u}(\epsilon_{\ast}) \! 
= \! \emptyset$ (see \cite{av2}, Section 4), in which case $\epsilon_{\ast}$ is vacuous and may be set equal 
to zero, and $\widetilde{\mathfrak{D}}_{u} \! = \! \mathfrak{D}_{u}^{\blacktriangledown}$. Henceforth, in the 
asymptotics of all expressions, formulae, etc., depending on $u(\tau)$, the `notation' $\tau \! \to \! +\infty$ 
means $\mathfrak{D}_{u}^{\blacktriangledown} \! \ni \! \tau$ and $\Re (\tau) \! \to \! +\infty$. \hfill $\blacksquare$
\end{eeee}
\begin{eeee} \label{abeych} 
The function $\hat{h}_{0}(\tau)$ defined by equation \eqref{iden2} plays a prominent r\^{o}le in the 
asymptotic estimates of this work; for further reference, therefore, a compact expression for it, which 
simplifies several of the ensuing estimates, is presented here: via equation~\eqref{iden9} and the 
definition~\eqref{iden2}, one shows that
\begin{equation} \label{expforeych} 
\hat{h}_{0}(\tau) \! = \! \alpha_{k}^{2} \tau^{-2/3} \! \left(\dfrac{\varkappa_{0}^{2}(\tau)}{4} \! - \! 
\dfrac{(a \! - \! \mi/2)}{1 \! + \! v_{0}(\tau) \tau^{-1/3}} \right), \quad k \! = \! \pm 1,
\end{equation}
where $\alpha_{k}$ is defined by equation \eqref{thmk5}, and the function $\varkappa_{0}^{2}(\tau)$ has 
the following equivalent representations:
\begin{align} \label{expforkapp} 
\left(\dfrac{\varkappa_{0}(\tau)}{\tau^{1/3}} \right)^{2} =& \, \left(2 \alpha_{k} \! + \! \dfrac{(\varepsilon b)^{1/3}
r(\tau)}{2 \alpha_{k}} \right)^{2} \! + \! \left(\dfrac{1}{\alpha_{k}^{2}} \! + \! \dfrac{r(\tau)}{(\varepsilon b)^{1/3}
(1 \! + \! \hat{u}_{0}(\tau))} \right) \! \left(-2(\varepsilon b)^{2/3}(1 \! + \! \hat{u}_{0}(\tau)) \! + \! 
\dfrac{(\varepsilon b)}{\alpha_{k}^{2}} \right) \nonumber \\
=& \, -\dfrac{\varepsilon b}{8 \alpha_{k}^{4}} \! \left(\dfrac{(8v_{0}^{2}(\tau) \! + \! 4 \tilde{r}_{0}(\tau)v_{0}
(\tau) \! - \! (\tilde{r}_{0}(\tau))^{2}) \tau^{-2/3} \! - \! (\tilde{r}_{0}(\tau))^{2}v_{0}(\tau) \tau^{-1}}{1 \! + \! 
v_{0}(\tau) \tau^{-1/3}} \right) \nonumber \\
=& \, -\left(2 \alpha_{k} \! - \! \dfrac{(\varepsilon b)^{1/3}r(\tau)}{2 \alpha_{k}} \right) \! \left(2 \alpha_{k} 
\! + \! \dfrac{(\varepsilon b)^{1/3}r(\tau)}{2 \alpha_{k}} \right) \nonumber \\
+& \, \dfrac{1}{\alpha_{k}^{2}} \! \left(\dfrac{2 \varepsilon b}{\alpha_{k}^{2}} \! + \! (\varepsilon b)^{2/3} \! 
\left(-2(1 \! + \! \hat{u}_{0}(\tau)) \! + \! \dfrac{r(\tau)}{1 \! + \! \hat{u}_{0}(\tau)} \right) \right).
\end{align}
It follows {}from the conditions \eqref{iden5} that $\lvert \varkappa_{0}^{2}(\tau) \rvert \genfrac{}{}{0pt}{3}{=}{
\tau \to +\infty} \mathcal{O}(\tau^{-2/3})$. \hfill $\blacksquare$
\end{eeee}

From Proposition \ref{prop1.2}, the definitions~\eqref{newlax2}, equations \eqref{eq3.2}, equation \eqref{iden4oldu}, 
and the definition \eqref{iden4}, one deduces that, in terms of the function $v_{0}(\tau)$, the solution of the DP3E 
\eqref{eq1.1} is given by
\begin{equation} \label{iden217} 
u(\tau) \! = \! c_{0,k} \tau^{1/3} \big(1 \! + \! \tau^{-1/3}v_{0}(\tau) \big), \quad k \! = \! \pm 1,
\end{equation}
where $c_{0,k}$ is defined by equation \eqref{thmk2}. As per the argument at the end of Subsection \ref{sec1a} 
regarding the particular form of the asymptotics for $u(\tau)$ as $\tau \! \to \! +\infty$ with $\varepsilon b \! > \! 0$ 
(cf. equation \eqref{mainyoo} and Remark \ref{remlndaip34}), it follows that, in conjunction with the representation 
\eqref{iden217}, the function $v_{0}(\tau)$ can be presented in the form
\begin{equation} \label{tr1} 
v_{0}(\tau) \! := \! v_{0,k}(\tau) \! \underset{\tau \to +\infty}{=} \sum_{m=0}^{\infty} \dfrac{\mathfrak{u}_{m}
(k)}{(\tau^{1/3})^{m+1}} \! + \! \mathrm{A}_{k} \me^{-(\beta (\tau)+\mi k \vartheta (\tau))} \! \left(1 \! + \! \mathcal{O}
(\tau^{-1/3}) \right), \quad k \! = \! \pm 1,
\end{equation}
where the sequence of $\mathbb{C}$-valued expansion coefficients $\lbrace \mathfrak{u}_{m}(k) \rbrace_{m=0}^{\infty}$ 
are determined in Proposition \ref{recursys} below, $\vartheta (\tau)$ and $\beta (\tau)$ are defined in equations 
\eqref{thmk12}, and, in the course of the ensuing analysis, it will be established that $\mathrm{A}_{k}$ depends on the 
Stokes multiplier $s_{0}^{0}$ (see Section \ref{finalsec}, equations \eqref{geek109} and \eqref{geek111}). {}From equation 
\eqref{iden217} and the expansion \eqref{tr1}, it follows that the associated solution of the DP3E \eqref{eq1.1} has asymptotics
\begin{equation} \label{recur15} 
u(\tau) \underset{\tau \to +\infty}{=} c_{0,k} \tau^{1/3} \! \left(1 \! + \! \sum_{m=0}^{\infty} \dfrac{\mathfrak{u}_{m}
(k)}{(\tau^{1/3})^{m+2}} \! + \! \mathrm{A}_{k} \tau^{-1/3} \me^{-(\beta (\tau)+\mi k \vartheta (\tau))} \! \left(1 \! 
+ \! \mathcal{O} \big(\tau^{-1/3} \big) \right) \right), \quad k \! = \! \pm 1.
\end{equation}
\begin{bbbb} \label{recursys} 
For $u(\tau)$ the corresponding solution of the {\rm DP3E} \eqref{eq1.1}, let the function $v_{0}(\tau) \! := \! v_{0,k}
(\tau)$, $k \! = \! \pm 1$, have the asymptotic expansion stated in equation \eqref{tr1}$;$ then, the expansion 
coefficients $\mathfrak{u}_{m}(k)$, $m \! \in \! \mathbb{Z}_{+}$, are given in equations \eqref{thmk2}--\eqref{thmk10}.
\end{bbbb}

\emph{Proof}. See Appendix \ref{apprecuru}. \hfill $\qed$ 

It follows {}from equations \eqref{hatsoff9}, \eqref{eq3.2}, \eqref{eq3.5}, and \eqref{iden3oldr} that
\begin{align} \label{textfeqn1} 
\dfrac{u^{\prime}(\tau) \! - \! \mi b}{u(\tau)} =& \, \frac{2}{\tau^{1/3}} \! \left(\frac{2a(\tau)d(\tau)}{\sqrt{\smash[b]{-a(\tau)
b(\tau)}}} \! + \! \tau^{-2/3}(\mi a \! + \! 1/2) \right) \! = \! \mi (\varepsilon b)^{1/3} \tau^{-1/3}(-2 \! + \! \hat{r}_{0}(\tau));
\end{align}
thus, via the definition \eqref{iden3}, it follows that
\begin{equation} \label{tr2} 
\tilde{r}_{0}(\tau) \! = \! 2 \tau^{1/3} \! - \! \dfrac{\mi \me^{-\mi 2 \pi k/3} \tau^{2/3}}{(\varepsilon b)^{1/3}} 
\! \left(\dfrac{u^{\prime}(\tau) \! - \! \mi b}{u(\tau)} \right), \quad k \! = \! \pm 1.
\end{equation}
\begin{bbbb} \label{proprr} 
For $u(\tau)$ the corresponding solution of the {\rm DP3E} \eqref{eq1.1} having the differentiable asymptotics 
\eqref{recur15}, with $\mathfrak{u}_{m}(k)$, $m \! \in \! \mathbb{Z}_{+}$, $k \! = \! \pm 1$, given in Proposition 
\ref{recursys}, let the function $\tilde{r}_{0}(\tau)$ be given by equation \eqref{tr2}$;$ then, $\tilde{r}_{0}(\tau)$ 
has the following asymptotic expansion:
\begin{equation} \label{tr3} 
\tilde{r}_{0}(\tau) \! := \! \tilde{r}_{0,k}(\tau) \underset{\tau \to +\infty}{=} \sum_{m=0}^{\infty} \dfrac{\mathfrak{r}_{m}
(k)}{(\tau^{1/3})^{m+1}} \! + \! 2(1 \! + \! k\sqrt{3}) \mathrm{A}_{k} \me^{-(\beta (\tau)+\mi k \vartheta (\tau))} \! 
\left(1 \! + \! \mathcal{O}(\tau^{-1/3}) \right), \quad k \! = \! \pm 1,
\end{equation}
where the expansion coefficients $\mathstrut \mathfrak{r}_{m}(k)$, $m \! \in \! \mathbb{Z}_{+}$, are given in 
equations \eqref{thmk15} and \eqref{thmk16}.
\end{bbbb}

\emph{Proof}. Substituting the differentiable asymptotics \eqref{recur15} for $u(\tau)$ into equation \eqref{tr2} and 
using the expressions for the coefficients $c_{0,k}$, $\mathfrak{u}_{m}(k)$, and $\mathfrak{w}_{m}(k)$, $k \! = \! 
\pm 1$, $m \! \in \! \mathbb{Z}_{+}$, given in the proof of Proposition \ref{recursys} (cf. Appendix \ref{apprecuru}), 
one arrives at, after a lengthy, but otherwise straightforward, algebraic calculation, the asymptotics for $\tilde{r}_{0}
(\tau) \! := \! \tilde{r}_{0,k}(\tau)$ stated in the proposition. \hfill $\qed$
\begin{eeee} \label{remforkay} 
Hereafter, explicit $k$ dependencies for the subscripts of the functions $v_{0}(\tau)$ and $\tilde{r}_{0}(\tau)$ (cf. 
equations \eqref{tr1} and \eqref{tr3}, respectively) will be suppressed, except where absolutely necessary and/or 
where confusion may arise. \hfill $\blacksquare$
\end{eeee}

In certain domains of the complex $\widetilde{\mu}$-plane (see the discussion below), the leading term of asymptotics 
(as $\tau \! \to \! +\infty$ for $\varepsilon b \! > \! 0)$ of a fundamental solution of equation \eqref{eq3.3} is given by the 
following matrix WKB formula (see, for example, Chapter 5 of \cite{F}),\footnote{Hereafter, for simplicity of notation, 
explicit $\tau$ dependencies will be suppressed, except where absolutely necessary.}
\begin{equation} \label{eq3.16} 
T(\widetilde{\mu}) \exp \! \left(-\sigma_{3} \mi \tau^{2/3} \int_{}^{\widetilde{\mu}}l(\xi) \, \md \xi \! - \! 
\int_{}^{\widetilde{\mu}} \text{diag} \! \left((T(\xi))^{-1} \partial_{\xi} T(\xi) \right) \md \xi \right) \! := 
\widetilde{\Psi}_{\scriptscriptstyle \mathrm{WKB}}(\widetilde{\mu}),
\end{equation}
where
\begin{equation}
l(\widetilde{\mu}) \! := \! \sqrtsign{\det (\mathcal{A}(\widetilde{\mu}))}, \label{eq3.17}
\end{equation}
and the matrix $T(\widetilde{\mu})$, which diagonalizes $\mathcal{A}(\widetilde{\mu})$, that is, $(T(\widetilde{\mu}))^{-1} 
\mathcal{A}(\widetilde{\mu})T(\widetilde{\mu}) \! = \! -\mi l(\widetilde{\mu}) \sigma_{3}$, is given by
\begin{equation}
T(\widetilde{\mu}) \! = \! \dfrac{\mi}{\sqrt{\smash[b]{2 \mi l(\widetilde{\mu})(\mathcal{A}_{11}(\widetilde{\mu}) 
\! - \! \mi l(\widetilde{\mu}))}}} \left(\mathcal{A}(\widetilde{\mu}) \! - \! \mi l(\widetilde{\mu}) \sigma_{3} \right) 
\sigma_{3}. \label{eq3.18}
\end{equation}
\begin{bbbb}[\textrm{\cite{av2}}] \label{prop3.1.2}
Let $T(\widetilde{\mu})$ be given in equation \eqref{eq3.18}, with $\mathcal{A}(\widetilde{\mu})$ and 
$l(\widetilde{\mu})$ defined by equations \eqref{eq3.4} and \eqref{eq3.17}, respectively. Then, $\det 
(T(\widetilde{\mu})) \! = \! 1$ and $\tr ((T(\widetilde{\mu}))^{-1} \partial_{\widetilde{\mu}}T(\widetilde{\mu})) 
\! = \! 0$$;$ moreover,
\begin{equation}
\diag \! \left((T(\widetilde{\mu}))^{-1} \partial_{\widetilde{\mu}}T(\widetilde{\mu}) \right) \! = \! -\dfrac{1}{2} \! \left(
\dfrac{\mathcal{A}_{12}(\widetilde{\mu}) \partial_{\widetilde{\mu}} \mathcal{A}_{21}(\widetilde{\mu}) \! - \! 
\mathcal{A}_{21}(\widetilde{\mu}) \partial_{\widetilde{\mu}} \mathcal{A}_{12}(\widetilde{\mu})}{2l(\widetilde{\mu})
(\mi \mathcal{A}_{11}(\widetilde{\mu}) \! + \! l(\widetilde{\mu}))} \right) \! \sigma_{3}. \label{eq3.23}
\end{equation}
\end{bbbb}
\begin{ffff} \label{cor3.1.1}
Let $\widetilde{\Psi}_{\scriptscriptstyle \mathrm{WKB}}(\widetilde{\mu})$ be defined by equation \eqref{eq3.16}, with 
$l(\widetilde{\mu})$ defined by equation \eqref{eq3.17} and $T(\widetilde{\mu})$ given in equation \eqref{eq3.18}$;$ 
then, $\det (\widetilde{\Psi}_{\scriptscriptstyle \mathrm{WKB}}(\widetilde{\mu})) \! = \! 1$.
\end{ffff}

\begin{figure}[hptb]
\begin{center}
\includegraphics[bb=30mm 25mm 250mm 255mm,scale=0.92]{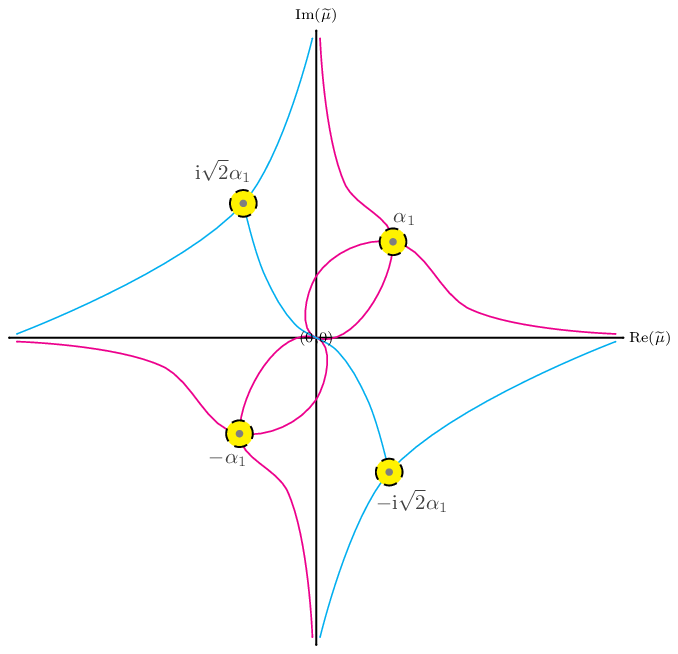}
\vspace{-11cm}
\caption{The Stokes graph for $k \! = \! +1$.}
\label{stograf1}
\end{center}
\end{figure}

\begin{figure}[hptb]
\begin{center}
\includegraphics[bb=30mm 25mm 250mm 255mm,scale=0.92]{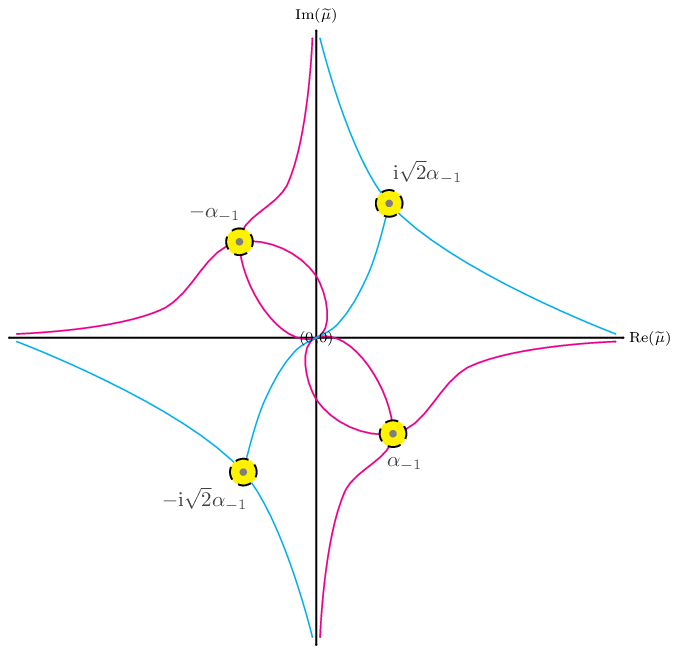}
\vspace{-11cm}
\caption{The Stokes graph for $k \! = \! -1$.}
\label{stograf2}
\end{center}
\end{figure}

The domains in the complex $\widetilde{\mu}$-plane where equation \eqref{eq3.16} gives the---leading---asymptotic 
approximation of solutions to equation \eqref{eq3.3} are defined in terms of the \emph{Stokes graph} (see, for 
example, \cite{F,mlaud,W}). The vertices of the Stokes graph are the singular points of equation \eqref{eq3.3}, that 
is, $\widetilde{\mu} \! = \! 0$ and $\widetilde{\mu} \! = \! \infty$, and the \emph{turning points\/}, which are the roots 
of the equation $l^{2}(\widetilde{\mu}) \! = \! 0$. The edges of the Stokes graph are the \emph{Stokes curves\/}, defined 
as $\Im (\int_{\widetilde{\mu}_{\scriptscriptstyle \mathrm{TP}}}^{\widetilde{\mu}}l(\xi) \, \md \xi) \! = \! 0$, where 
$\widetilde{\mu}_{\scriptscriptstyle \mathrm{TP}}$ denotes a turning point. \emph{Canonical domains\/} are those domains in 
the complex $\widetilde{\mu}$-plane containing one, and only one, Stokes curve and bounded by two adjacent Stokes curves. 
(Note that the restriction of any branch of $l(\widetilde{\mu})$ to a canonical domain is a single-valued function.) In each 
canonical domain, for any choice of the branch of $l(\widetilde{\mu})$, there exists a fundamental solution of equation 
\eqref{eq3.3} which has asymptotics whose leading term is given by equation \eqref{eq3.16}. {}From the definition of 
$l(\widetilde{\mu})$ given by equation \eqref{eq3.17}, one arrives at
\begin{equation} \label{eq3.19} 
l^{2}(\widetilde{\mu}) \! := \! l_{k}^{2}(\widetilde{\mu}) \! = \! \frac{4}{\widetilde{\mu}^{4}} \! \left(\left(\widetilde{\mu}^{2} 
\! - \! \alpha_{k}^{2} \right)^{2} \left(\widetilde{\mu}^{2} \! + \! 2 \alpha_{k}^{2} \right) \! + \! \widetilde{\mu}^{2} \hat{h}_{0}
(\tau) \! + \! \widetilde{\mu}^{4} \! \left(a \! - \! \mi/2 \right) \tau^{-2/3} \right), \quad 
k \! = \! \pm 1,
\end{equation}
where $\alpha_{k}$ is defined by equation~\eqref{thmk5}. It follows {}from equation~\eqref{eq3.19} that there are six 
turning points. For $k \! = \! \pm 1$, the conditions \eqref{iden5} imply that one pair of turning points coalesce at 
$\alpha_{k}$ with asymptotics $\mathcal{O}(\tau^{-1/3})$, another pair has asymptotics $-\alpha_{k} \! + \! \mathcal{O}
(\tau^{-1/3})$, and the two remaining turning points have the asymptotic behaviour $\pm \mi \sqrt{2} \alpha_{k} \! + \! 
\mathcal{O}(\tau^{-2/3})$. For simplicity of notation, denote by $\widetilde{\mu}_{1}(k)$ any one of the turning points 
coalescing at $\alpha_{k}$, and denote by $\widetilde{\mu}_{2}(k)$ the turning point approaching $\mi k \sqrt{2} \alpha_{k}$. 
Let $\mathscr{G}_{\scriptscriptstyle \mathbb{S}}(k)$, $k \! = \! \pm 1$, be the part of the Stokes graph that consists 
of the vertices $0,\infty,\widetilde{\mu}_{1}(k)$ and $\widetilde{\mu}_{2}(k)$, and the union of the---oriented---edges 
$\operatorname{arc}(\mi k \infty,\widetilde{\mu}_{2}(k))$, $\operatorname{arc}(\widetilde{\mu}_{2}(k),0)$ and 
$\operatorname{arc}(\widetilde{\mu}_{2}(k),-\infty)$, and $\operatorname{arc}(\mi k \infty,\widetilde{\mu}_{1}(k))$, 
$\operatorname{arc}(\widetilde{\mu}_{1}(k),0)$, $\operatorname{arc}(0,\widetilde{\mu}_{1}(k))$ and $\operatorname{arc}
(\widetilde{\mu}_{1}(k),+\infty)$; the complete Stokes graph is given by $\mathscr{G}_{\scriptscriptstyle \mathbb{S}}(k) 
\cup \me^{\mi \pi} \mathscr{G}_{\scriptscriptstyle \mathbb{S}}(k)$ (see Figure \ref{stograf1} (resp., Figure \ref{stograf2}) 
for the case $k \! = \! +1$ (resp., $k \! = \! -1))$.
\begin{bbbb} \label{prop3.1.3} 
Let $l_{k}^{2}(\widetilde{\mu})$, $k \! = \! \pm 1$, be given in equation \eqref{eq3.19}$;$ then,
\begin{equation}
\int_{\widetilde{\mu}_{0,k}}^{\widetilde{\mu}}l_{k}(\xi) \, \md \xi \underset{\tau \to +\infty}{=} \varUpsilon_{k}
(\widetilde{\mu}) \! - \! \varUpsilon_{k}(\widetilde{\mu}_{0,k}) \! + \! \mathcal{O}(\mathscr{E}_{k}(\widetilde{\mu})) 
\! + \! \mathcal{O}(\mathscr{E}_{k}(\widetilde{\mu}_{0,k})), \label{eq3.21}
\end{equation}
where, for $\delta \! > \! 0$, $\widetilde{\mu},\widetilde{\mu}_{0,k} \! \in \! \mathbb{C} \setminus 
(\mathscr{O}_{\tau^{-1/3+ \delta}}(\pm \alpha_{k}) \cup \mathscr{O}_{\tau^{-2/3+2 \delta}}(\pm \mi \sqrt{2} 
\alpha_{k}) \cup \lbrace 0,\infty \rbrace)$ and the path of integration lies in the corresponding 
canonical domain,
\begin{align}
\varUpsilon_{k}(\xi) :=& \, (\xi \! + \! 2 \alpha_{k}^{2} \xi^{-1})(\xi^{2} \! + \! 2 \alpha_{k}^{2})^{1/2} \! + \! 
\tau^{-2/3}(a \! - \! \mi/2) \ln \! \left(\xi \! + \! (\xi^{2} \! + \! 2 \alpha_{k}^{2})^{1/2} \right) \nonumber \\
+& \, \dfrac{\tau^{-2/3}}{2 \sqrt{3}} \! \left((a \! - \! \mi/2) \! + \! \frac{\tau^{2/3}}{\alpha_{k}^{2}} \hat{h}_{0}(\tau) 
\right) \! \ln \! \left(\left(\dfrac{3^{1/2}(\xi^{2} \! + \! 2 \alpha_{k}^{2})^{1/2} \! - \! \xi \! + \! 2 \alpha_{k}}{3^{1/2}
(\xi^{2} \! + \! 2 \alpha_{k}^{2})^{1/2} \! + \! \xi \! + \! 2 \alpha_{k}} \right) \! \left(\dfrac{\xi \! - \! \alpha_{k}}{\xi \! 
+ \! \alpha_{k}} \right) \right), \label{equpsi}
\end{align}
and
\begin{equation} \label{eqeel} 
\tau^{4/3} \mathscr{E}_{k}(\xi) \! := \! 
\begin{cases}
{\fontsize{9pt}{11pt}\selectfont 
\begin{aligned}[b]
&\frac{\big((a \! - \! \mi/2) \! + \! \frac{\tau^{2/3}}{\alpha_{k}^{2}} \hat{h}_{0}(\tau) \big)^{2}}{192 \sqrt{3}
(\xi \! \mp \! \alpha_{k})^{2}} \! + \! \mathcal{O} \! \left(\frac{c_{1,k} \! + \! c_{2,k} \tau^{2/3} \hat{h}_{0}(\tau) 
\! + \! c_{3,k}(\tau^{2/3} \hat{h}_{0}(\tau))^{2}}{\xi \! \mp \! \alpha_{k}} \right), 
\end{aligned}} &\text{$\xi \! \in \! \mathbb{U}_{k}^{1}$,} \\
{\fontsize{8pt}{11pt}\selectfont 
\begin{aligned}[b]
&\frac{\big((a \! - \! \mi/2) \! - \! \frac{\tau^{2/3}}{2 \alpha_{k}^{2}} \hat{h}_{0}(\tau) \big)^{2}}{d_{0,k}(\xi \! \mp \! \mi 
\sqrtsign{2} \alpha_{k})^{1/2}} \! + \! \mathcal{O} \! \left((\xi \! \mp \! \mi \sqrt{2} \alpha_{k})^{1/2} \big(c_{4,k} \! + \! 
c_{5,k} \tau^{2/3} \hat{h}_{0}(\tau) \! + \! c_{6,k}(\tau^{2/3} \hat{h}_{0}(\tau))^{2} \big) \right),
\end{aligned}} &\text{$\xi \! \in \! \mathbb{U}_{k}^{2}$,} \\
{\fontsize{10pt}{11pt}\selectfont 
\begin{aligned}[b]
&\mathfrak{f}_{1,k}(\xi^{-1}) \! + \! \tau^{2/3} \hat{h}_{0}(\tau) \mathfrak{f}_{2,k}(\xi^{-1}) \! + \! (\tau^{2/3} \hat{h}_{0}
(\tau))^{2} \mathfrak{f}_{3,k}(\xi^{-1}),
\end{aligned}} &\text{$\xi \! \to \! \infty$,} \\
{\fontsize{10pt}{11pt}\selectfont 
\begin{aligned}[b]
\mathfrak{f}_{4,k}(\xi) \! + \! \tau^{2/3} \hat{h}_{0}(\tau) \mathfrak{f}_{5,k}(\xi) \! + \! (\tau^{2/3} \hat{h}_{0}(\tau))^{2} 
\mathfrak{f}_{6,k}(\xi),
\end{aligned}} &\text{$\xi \! \to \! 0$,} 
\end{cases}
\end{equation}
where $\mathbb{U}_{k}^{1} \! := \! \mathscr{O}_{\tau^{-1/3 + \delta_{k}}}(\pm \alpha_{k})$, $\mathbb{U}_{k}^{2} \! 
:= \! \mathscr{O}_{\tau^{-2/3 +2 \delta_{k}}}(\pm \mi \sqrt{2} \alpha_{k})$, the parameter $\delta_{k}$ satisfies 
(see Corollary \ref{cor3.1.2}$)$ $0 \! < \! \delta \! < \! \delta_{k} \! < \! 1/9$, $d_{0,k}^{-1} \! := \! 2^{-1/4} 
\me^{\mp \mi 3 \pi/4} \alpha_{k}^{-3/2}/27$, $\mathfrak{f}_{j,k}(z)$, $j \! = \! 1,2,\dotsc,6$, are analytic functions 
of $z$ in a neighbourhood of $z \! = \! 0$ given in equations \eqref{alpbet1}--\eqref{alpbet6}, and $c_{m,k}$, 
$m \! = \! 1,2,\dotsc,6$, are $\mathcal{O}(1)$.
\end{bbbb}

\emph{Proof}. Let $l_{k}^{2}(\widetilde{\mu})$, $k \! = \! \pm 1$, be given in equation \eqref{eq3.19}, with $\alpha_{k}$ 
defined by equation \eqref{thmk5}. Recalling {}from the conditions \eqref{iden5} that $\lvert \hat{h}_{0}(\tau) \rvert 
\genfrac{}{}{0pt}{3}{=}{\tau \to +\infty} \mathcal{O}(\tau^{-2/3})$, set
\begin{equation}
l_{k,\infty}^{2}(\widetilde{\mu}) \! = \! 4 \widetilde{\mu}^{-4}(\widetilde{\mu}^{2} \! - \! \alpha_{k}^{2})^{2}(\widetilde{\mu}^{2} 
\! + \! 2 \alpha_{k}^{2}). \label{eqlsquared}
\end{equation}
Define
\begin{equation}
\Delta_{k,\tau}(\widetilde{\mu}) \! := \! \dfrac{l_{k}^{2}(\widetilde{\mu}) 
\! - \! l_{k,\infty}^{2}(\widetilde{\mu})}{l_{k,\infty}^{2}(\widetilde{\mu})} 
\! = \! \dfrac{\widetilde{\mu}^{2} \hat{h}_{0}(\tau) \! + \! \widetilde{\mu}^{4}
(a \! - \! \mi/2) \tau^{-2/3}}{(\widetilde{\mu}^{2} \! - \! \alpha_{k}^{2})^{2}
(\widetilde{\mu}^{2} \! + \! 2 \alpha_{k}^{2})}; \label{eqdellt}
\end{equation}
hence, presenting $l_{k}(\widetilde{\mu})$ as $l_{k}(\widetilde{\mu}) \! = \! l_{k,\infty}(\widetilde{\mu})
(1 \! + \! \Delta_{k,\tau}(\widetilde{\mu}))^{1/2}$, a straightforward calculation, via the conditions \eqref{iden5}, 
shows that, for $k \! = \! \pm 1$,
\begin{align}
l_{k}(\widetilde{\mu}) \underset{\tau \to +\infty}{=}& \, l_{k,\infty}(\widetilde{\mu}) \! \left(1 \! + \! \Delta_{k,\tau}
(\widetilde{\mu})/2  \! + \! \mathcal{O}(-(\Delta_{k,\tau}(\widetilde{\mu}))^{2}/8) \right) \nonumber \\
\underset{\tau \to +\infty}{=}& \, 2(1 \! - \! \alpha_{k}^{2}/\widetilde{\mu}^{2})(\widetilde{\mu}^{2} \! + \! 2 
\alpha_{k}^{2})^{1/2} \! + \! \dfrac{\hat{h}_{0}(\tau) \! + \! \widetilde{\mu}^{2}(a \! - \! \mi/2) \tau^{-2/3}}{
(\widetilde{\mu}^{2} \! - \! \alpha_{k}^{2})(\widetilde{\mu}^{2} \! + \! 2 \alpha_{k}^{2})^{1/2}} \nonumber \\
+& \, \mathcal{O} \! \left(-\dfrac{\widetilde{\mu}^{2} \big(\hat{h}_{0}(\tau) \! + \! \widetilde{\mu}^{2}(a \! - \! \mi/2) 
\tau^{-2/3} \big)^{2}}{4(\widetilde{\mu}^{2} \! - \! \alpha_{k}^{2})^{3}(\widetilde{\mu}^{2} \! + \! 2 \alpha_{k}^{2})^{3/2}} 
\right). \label{eqforl}
\end{align}
Integration of the two terms in the second line of equation \eqref{eqforl} gives rise to the leading term of asymptotics in 
equation \eqref{eq3.21}, and integration of the error term in the third line of equation \eqref{eqforl} leads to an explicit 
expression for the error function, $\mathscr{E}_{k}(\pmb{\cdot})$, whose asymptotics at the turning and the singular points 
read: (i) for $\xi \! \in \! \mathscr{O}_{\tau^{-1/3+ \delta_{k}}}(\pm \alpha_{k})$, $0 \! < \! \delta \! < \! \delta_{k} \! < \! 1/9$,
\begin{align} \label{terrbos1} 
\tau^{4/3} \mathscr{E}_{k}(\xi) \underset{\tau \to +\infty}{=}& \, \dfrac{\big((a \! - \! \mi/2) \! + \! \alpha_{k}^{-2} 
\tau^{2/3} \hat{h}_{0}(\tau) \big)^{2}}{192 \sqrt{3} (\xi \! \mp \! \alpha_{k})^{2}} \! + \! \dfrac{\hat{d}_{-1,k}(\tau)}{\xi 
\! \mp \! \alpha_{k}} \! + \! \hat{d}_{0,k}(\tau) \ln (\xi \! \mp \! \alpha_{k}) \nonumber \\
+& \, \sum_{m \in \mathbb{Z}_{+}} \hat{d}_{m+1,k}(\tau)(\xi \! \mp \! \alpha_{k})^{m+1},
\end{align}
where
\begin{equation*}
\hat{d}_{m,k}(\tau) \! := \! \hat{\mathfrak{c}}^{\flat}_{m,k} \! + \! \hat{\mathfrak{c}}^{\natural}_{m,k} \tau^{2/3} 
\hat{h}_{0}(\tau) \! + \! \hat{\mathfrak{c}}^{\sharp}_{m,k}(\tau^{2/3} \hat{h}_{0}(\tau))^{2}, \quad m \! \in \! 
\lbrace -1 \rbrace \cup \mathbb{Z}_{+},
\end{equation*}
and $\hat{\mathfrak{c}}^{r}_{m,k}$, $r \! \in \! \lbrace \flat,\natural,\sharp \rbrace$, are $\mathcal{O}(1)$, and thus, retaining only 
the first two terms of the expansion \eqref{terrbos1}, one arrives at the representation for $\mathscr{E}_{k}(\xi)$ stated in the 
first line of equation \eqref{eqeel}; (ii) for $\xi \! \in \! \mathscr{O}_{\tau^{-2/3+ 2 \delta_{k}}}(\pm \mi \sqrt{2} \alpha_{k})$,
\begin{equation} \label{terrbos2} 
\tau^{4/3} \mathscr{E}_{k}(\xi) \underset{\tau \to +\infty}{=} \dfrac{2^{-1/4} \big((a \! - \! \mi/2) \! - \! \tau^{2/3} \hat{h}_{0}
(\tau)/2 \alpha_{k}^{2} \big)^{2}}{27 \alpha_{k}^{3/2} \me^{\pm \mi 3 \pi/4}(\xi \! \mp \! \mi \sqrt{2} \alpha_{k})^{1/2}} \! + \! 
(\xi \! \mp \! \mi \sqrt{2} \alpha_{k})^{1/2} \sum_{m \in \mathbb{Z}_{+}} \tilde{d}_{m,k}(\tau)(\xi \! \mp \! \mi \sqrt{2} \alpha_{k})^{m},
\end{equation}
where
\begin{equation*}
\tilde{d}_{m,k}(\tau) \! := \! \tilde{\mathfrak{c}}^{\flat}_{m,k} \! + \! \tilde{\mathfrak{c}}^{\natural}_{m,k} \tau^{2/3} 
\hat{h}_{0}(\tau) \! + \! \tilde{\mathfrak{c}}^{\sharp}_{m,k}(\tau^{2/3} \hat{h}_{0}(\tau))^{2}, \quad m \! \in \! \mathbb{Z}_{+},
\end{equation*}
and $\tilde{\mathfrak{c}}^{r}_{m,k}$, $r \! \in \! \lbrace \flat,\natural,\sharp \rbrace$, are $\mathcal{O}(1)$, and thus, keeping only 
the first two terms of the expansion \eqref{terrbos2}, one arrives at the representation for $\mathscr{E}_{k}(\xi)$ stated in the 
second line of equation \eqref{eqeel}; (iii) as $\xi \! \to \! \infty$, one arrives at the representation for $\mathscr{E}_{k}(\xi)$ stated 
in the third line of equation \eqref{eqeel}, where
\begin{gather}
\mathfrak{f}_{1,k}(z) \! = \! \frac{(a \! - \! \mi/2)^{2}}{8}z^{2} \! + \! (a \! - \! \mi/2)^{2}z^{6} \sum_{m \in \mathbb{Z}_{+}} 
\hat{\mathfrak{c}}^{\circ,1}_{m,k}z^{2m}, \label{alpbet1} \\
\mathfrak{f}_{2,k}(z) \! = \! \dfrac{(a \! - \! \mi/2)}{8}z^{4} \! + \! (a \! - \! \mi/2)z^{8} \sum_{m \in \mathbb{Z}_{+}} 
\hat{\mathfrak{c}}^{\circ,2}_{m,k}z^{2m}, \label{alpbet2} \\
\mathfrak{f}_{3,k}(z) \! = \! \dfrac{1}{24}z^{6} \! + \! z^{10} \sum_{m \in \mathbb{Z}_{+}} \hat{\mathfrak{c}}^{\circ,3}_{m,k}
z^{2m}, \label{alpbet3}
\end{gather}
and $\hat{\mathfrak{c}}^{\circ,r}_{m,k}$, $r \! = \! 1,2,3$, $m \! \in \! \mathbb{Z}_{+}$, are $\mathcal{O}(1)$; and (iv) as $\xi \! 
\to \! 0$, one arrives at the representation for $\mathscr{E}_{k}(\xi)$ stated in the fourth line of equation \eqref{eqeel}, where
\begin{gather}
\mathfrak{f}_{4,k}(z) \! = \! \dfrac{(a \! - \! \mi/2)^{2}}{56 \sqrt{2} \, \alpha_{k}^{9}}z^{7} \! + \! (a \! - \! \mi/2)^{2}z^{9} 
\sum_{m \in \mathbb{Z}_{+}} \tilde{d}^{\circ,4}_{m,k}z^{2m}, \label{alpbet4} \\
\mathfrak{f}_{5,k}(z) \! = \! \dfrac{(a \! - \! \mi/2)}{20 \sqrt{2} \, \alpha_{k}^{9}}z^{5} \! + \! (a \! - \! \mi/2)z^{7} 
\sum_{m \in \mathbb{Z}_{+}} \tilde{d}^{\circ,5}_{m,k}z^{2m}, \label{alpbet5} \\
\mathfrak{f}_{6,k}(z) \! = \! \dfrac{1}{24 \sqrt{2} \, \alpha_{k}^{9}}z^{3} \! + \! z^{5} \sum_{m \in \mathbb{Z}_{+}} 
\tilde{d}^{\circ,6}_{m,k}z^{2m}, \label{alpbet6}
\end{gather}
and $\tilde{d}^{\circ,r}_{m,k}$, $r \! = \! 4,5,6$, $m \! \in \! \mathbb{Z}_{+}$, are $\mathcal{O}(1)$. \hfill $\qed$
\begin{ffff} \label{cor3.1.2} 
Set $\widetilde{\mu}_{0,k} \! = \! \alpha_{k} \! + \! \tau^{-1/3} \widetilde{\Lambda}$, $k \! = \! \pm 1$, where $\widetilde{\Lambda} 
\genfrac{}{}{0pt}{3}{=}{\tau \to +\infty} \mathcal{O}(\tau^{\delta_{k}})$, $0 \! < \! \delta \! < \! \delta_{k} \! < \! 1/9;$ then,
\begin{align} \label{eq3.22} 
\int_{\widetilde{\mu}_{0,k}}^{\widetilde{\mu}}l_{k}(\xi) \, \md \xi 
\underset{\tau \to +\infty}{=}& \, \varUpsilon_{k}(\widetilde{\mu}) \! + \! 
\varUpsilon_{k}^{\sharp} \! + \! \mathcal{O}(\mathscr{E}_{k}(\widetilde{\mu})) 
\! + \! \mathcal{O}(\tau^{-1} \widetilde{\Lambda}^{3}) \! + \! \mathcal{O}
(\tau^{-1} \widetilde{\Lambda}) \nonumber \\
+& \, \mathcal{O} \! \left(\dfrac{\tau^{-1}}{\widetilde{\Lambda}} \left(
\mathfrak{c}_{1,k} \! + \! \mathfrak{c}_{2,k} \tau^{2/3} \hat{h}_{0}(\tau) 
\! + \! \mathfrak{c}_{3,k}(\tau^{2/3} \hat{h}_{0}(\tau))^{2} \right) \right),
\end{align}
where $\varUpsilon_{k}(\widetilde{\mu})$ and $\mathscr{E}_{k}(\widetilde{\mu})$ are defined by equations 
\eqref{equpsi} and \eqref{eqeel}, respectively,
\begin{align}
\varUpsilon_{k}^{\sharp} :=& \, \mp 3 \sqrt{3} \alpha_{k}^{2} \! \mp \! 2 \sqrt{3} 
\tau^{-2/3} \widetilde{\Lambda}^{2} \! - \! \tau^{-2/3}(a \! - \! \mi/2) \ln \!  
\left((\sqrt{3} \pm \! 1) \alpha_{k} \me^{\mi \pi (1 \mp 1)/2} \right) \nonumber \\
\mp& \, \dfrac{\tau^{-2/3}}{2 \sqrt{3}} \! \left((a \! - \! \mi/2) \! + \! \alpha_{k}^{-2} 
\tau^{2/3} \hat{h}_{0}(\tau) \right) \! \left(\ln \widetilde{\Lambda} \! - \! 
\dfrac{1}{3} \ln \tau \! - \! \ln (3 \alpha_{k}) \right), \label{equpsisharp}
\end{align}
with the upper (resp., lower) signs taken according to the branch of the square-root function $\lim_{\xi^{2} \to +\infty}
(\xi^{2} \! + \! 2 \alpha_{k}^{2})^{1/2} \! = \! +\infty$ (resp., $\lim_{\xi^{2} \to +\infty}(\xi^{2} \! + \! 2 \alpha_{k}^{2})^{1/2} 
\! = \! -\infty)$, and $\mathfrak{c}_{m,k}$, $m \! = \! 1,2,3$, are $\mathcal{O}(1)$.
\end{ffff}

\emph{Proof}. Substituting $\widetilde{\mu}_{0,k}$, as given in the corollary, for the argument of the functions 
$\varUpsilon_{k}(\xi)$ and $\mathscr{E}_{k}(\xi)$ (cf. equation \eqref{equpsi} and the first line of equation \eqref{eqeel}, 
respectively) and expanding with respect to the ``small parameter'' $\tau^{-1/3} \widetilde{\Lambda}$, one arrives at the 
following estimates:
\begin{equation} \label{terrbos3} 
-\varUpsilon_{k}(\widetilde{\mu}_{0,k}) \underset{\tau \to +\infty}{=} 
\varUpsilon_{k}^{\sharp} \! + \! \mathcal{O}(\tau^{-1} \widetilde{
\Lambda}^{3}) \! + \! \mathcal{O}(\tau^{-1} \widetilde{\Lambda}) \! + \! 
\mathcal{O} \! \left(\tau^{-1} \widetilde{\Lambda}((a \! - \! \mi/2) \! + \! 
\alpha_{k}^{-2} \tau^{2/3} \hat{h}_{0}(\tau)) \right),
\end{equation}
where $\varUpsilon_{k}^{\sharp}$ is defined by equation \eqref{equpsisharp},
\begin{align} \label{terrbos4} 
\mathcal{O}(\mathscr{E}_{k}(\widetilde{\mu}_{0,k})) \underset{\tau \to +\infty}{=}& 
\, \mathcal{O} \! \left(\dfrac{\tau^{-2/3}}{\widetilde{\Lambda}^{2}}((a \! - \! \mi/2) 
\! + \! \alpha_{k}^{-2} \tau^{2/3} \hat{h}_{0}(\tau))^{2} \right) \nonumber \\
+& \, \mathcal{O} \! \left(\dfrac{\tau^{-1}}{\widetilde{\Lambda}} \left(c_{1,k} \! + \! 
c_{2,k} \tau^{2/3} \hat{h}_{0}(\tau) \! + \! c_{3,k}(\tau^{2/3} \hat{h}_{0}(\tau))^{2} 
\right) \right),
\end{align}
and $c_{m,k}$, $m \! = \! 1,2,3$, are $\mathcal{O}(1)$. {}From equations \eqref{iden9}, \eqref{iden2}, \eqref{iden3}, and 
\eqref{iden4}, one shows that
\begin{equation} \label{terrbos5} 
-\tau^{2/3} \hat{h}_{0}(\tau) \! = \! \dfrac{\alpha_{k}^{2}(a \! - \! \mi/2)}{1 \! + \! v_{0}(\tau) \tau^{-1/3}} \! + \! 
\dfrac{\alpha_{k}^{4}(8v_{0}^{2}(\tau) \! + \! 4 \tilde{r}_{0}(\tau)v_{0}(\tau) \! - \! (\tilde{r}_{0}(\tau))^{2} \! - \! 
v_{0}(\tau)(\tilde{r}_{0}(\tau))^{2} \tau^{-1/3})}{4(1 \! + \! v_{0}(\tau) \tau^{-1/3})},
\end{equation}
whence, via the conditions \eqref{iden5},
\begin{align} \label{terrbos6}
(a \! - \! \mi/2) \! + \! \frac{\tau^{2/3}}{\alpha_{k}^{2}} \hat{h}_{0}(\tau) 
\underset{\tau \to +\infty}{=}& \, -\dfrac{\alpha_{k}^{2}}{4} \left(8v_{0}^{2}
(\tau) \! + \! 4v_{0}(\tau) \tilde{r}_{0}(\tau) \! - \! (\tilde{r}_{0}(\tau))^{2} 
\right) \! + \! (a \! - \! \mi/2) v_{0}(\tau) \tau^{-1/3} \nonumber \\
+& \, \mathcal{O}((2v_{0}^{2}(\tau) \! + \! v_{0}(\tau) \tilde{r}_{0}(\tau))
v_{0}(\tau) \tau^{-1/3}) \! + \! \mathcal{O}(v_{0}^{2}(\tau) \tau^{-2/3}).
\end{align}
Note {}from the conditions \eqref{iden5} and the expansion \eqref{terrbos6} that
\begin{equation*}
(a \! - \! \mi/2) \! + \! \frac{\tau^{2/3}}{\alpha_{k}^{2}} \hat{h}_{0}(\tau) \underset{\tau \to +\infty}{=} \mathcal{O}
(\tau^{-2/3}) \quad \text{and} \! \quad \! c_{1,k} \! + \! c_{2,k} \tau^{2/3} \hat{h}_{0}(\tau) \! + \! c_{3,k}
(\tau^{2/3} \hat{h}_{0}(\tau))^{2} \underset{\tau \to +\infty}{=} \! \mathcal{O}(1):
\end{equation*}
{}from the expansions \eqref{terrbos3} and \eqref{terrbos4} and the latter two estimates, it follows that
\begin{gather}
-\varUpsilon_{k}(\widetilde{\mu}_{0,k}) \underset{\tau \to +\infty}{=} 
\varUpsilon_{k}^{\sharp} \! + \! \mathcal{O}(\tau^{-1} \widetilde{
\Lambda}^{3}) \! + \! \mathcal{O}(\tau^{-1} \widetilde{\Lambda}) \! 
+ \! \mathcal{O}(\tau^{-5/3} \widetilde{\Lambda}), \label{terrbos7} \\
\mathcal{O}(\mathscr{E}_{k}(\widetilde{\mu}_{0,k})) \underset{\tau \to +\infty}{=} 
\mathcal{O}(\tau^{-1} \widetilde{\Lambda}^{-1}) \! + \! \mathcal{O}
(\tau^{-2} \widetilde{\Lambda}^{-2}), \label{terrbos8}
\end{gather}
whence, introducing the inequality $0 \! < \! \delta \! < \! \delta_{k} \! < \! 1/9$ in order to guarantee that the error estimates in 
the expansions \eqref{terrbos7} and \eqref{terrbos8} are $o(1)$ after multiplication by the ``large parameter'' $\tau^{2/3}$ (cf. 
equation \eqref{eq3.16}), retaining only leading-order contributions, one arrives at
\begin{align*}
-\varUpsilon_{k}(\widetilde{\mu}_{0,k}) \! + \! \mathcal{O}(\mathscr{E}_{k}
(\widetilde{\mu}_{0,k})) \underset{\tau \to +\infty}{=}& \, \varUpsilon_{k}^{\sharp} 
\! + \! \mathcal{O} \! \left(\dfrac{\tau^{-1}}{\widetilde{\Lambda}} \left(
c_{1,k} \! + \! c_{2,k} \tau^{2/3} \hat{h}_{0}(\tau) \! + \! c_{3,k}
(\tau^{2/3} \hat{h}_{0}(\tau))^{2} \right) \right) \nonumber \\
+& \, \mathcal{O}(\tau^{-1} \widetilde{\Lambda}^{3}) \! + \! \mathcal{O}
(\tau^{-1} \widetilde{\Lambda}),
\end{align*}
which, via equation \eqref{eq3.21}, implies the result stated in the corollary. \hfill $\qed$
\begin{ffff} \label{cor3.1.4} 
Let the conditions stated in Corollary {\rm \ref{cor3.1.2}} be valid$;$ then, for the branch of $l_{k}(\xi)$, $k \! = \! \pm 1$, 
that is positive for large and small positive $\xi$,
\begin{align} \label{eq3.37} 
-\mi \tau^{2/3} \int_{\widetilde{\mu}_{0,k}}^{\widetilde{\mu}}l_{k}(\xi) \, \md \xi 
\underset{\underset{\widetilde{\mu} \to \infty}{\tau \to +\infty}}{=}& \, -\mi \big(\tau^{2/3} \widetilde{\mu}^{2} \! + \! (a \! - \! \mi/2) 
\ln \widetilde{\mu} \big) \! + \! \mi 3(\sqrt{3}- \! 1) \alpha_{k}^{2} \tau^{2/3} \! + \! \mi 2 \sqrt{3} \, \widetilde{\Lambda}^{2} \! + \! 
C^{\scriptscriptstyle \mathrm{WKB}}_{\infty,k} \nonumber \\
-& \, \dfrac{\mi}{2 \sqrt{3}} \big((a \! - \! \mi/2) \! + \! \alpha_{k}^{-2} \tau^{2/3} \hat{h}_{0}(\tau) \big) \! \left(\dfrac{1}{3} \ln 
\tau \! - \! \ln \widetilde{\Lambda} \! + \! \ln \left(\dfrac{6 \alpha_{k}}{(\sqrt{3} \! + \! 1)^{2}} \right) \right) \nonumber \\
+& \, \mathcal{O} \! \left(\dfrac{\tau^{-1/3}}{\widetilde{\Lambda}} \big(\mathfrak{c}_{1,k} \! + \! \mathfrak{c}_{2,k} \tau^{2/3} 
\hat{h}_{0}(\tau) \! + \! \mathfrak{c}_{3,k}(\tau^{2/3} \hat{h}_{0}(\tau))^{2} \big) \right) \nonumber \\
+& \, \mathcal{O}(\tau^{-1/3} \widetilde{\Lambda}^{3}) \! + \! \mathcal{O}(\tau^{-1/3} \widetilde{\Lambda}) \! + \! 
\mathcal{O}(\tau^{-2/3} \widetilde{\mu}^{-3}),
\end{align}
where
\begin{equation} \label{eq3.38} 
C^{\scriptscriptstyle \mathrm{WKB}}_{\infty,k} \! := \! \mi (a \! - \! \mi/2) \ln ((\sqrtsign{3} \! + \! 1) \alpha_{k}/2),
\end{equation}
and
\begin{align} \label{eq3.39} 
-\mi \tau^{2/3} \int_{\widetilde{\mu}_{0,k}}^{\widetilde{\mu}}l_{k}(\xi) \, \md \xi 
\underset{\underset{\widetilde{\mu} \to 0}{\tau \to +\infty}}{=}& \, \dfrac{1}{\widetilde{\mu}} \mi 2 \sqrt{2} \alpha_{k}^{3} 
\tau^{2/3} \! - \! \mi 3 \sqrt{3} \alpha_{k}^{2} \tau^{2/3} \! - \! \mi 2 \sqrt{3} \, \widetilde{\Lambda}^{2} \! + \! 
\dfrac{\mi}{2 \sqrt{3}} \! \left((a \! - \! \mi/2) \vphantom{M^{M^{m}}} \right. \nonumber \\
+&\left. \, \alpha_{k}^{-2} \tau^{2/3} \hat{h}_{0}(\tau) \right) \! \left(\dfrac{1}{3} \ln \tau \! - \! \ln \widetilde{\Lambda} 
\! + \! \ln (3 \alpha_{k} \me^{-\mi \pi k}) \right) \! + \! C^{\scriptscriptstyle \mathrm{WKB}}_{0,k} \nonumber \\
+& \, \mathcal{O} \! \left(\dfrac{\tau^{-1/3}}{\widetilde{\Lambda}} \big(\mathfrak{c}_{4,k} \! + \! \mathfrak{c}_{5,k} 
\tau^{2/3} \hat{h}_{0}(\tau) \! + \! \mathfrak{c}_{6,k}(\tau^{2/3} \hat{h}_{0}(\tau))^{2} \big) \right) \nonumber \\
+& \, \mathcal{O}(\tau^{-1/3} \widetilde{\Lambda}^{3}) \! + \! \mathcal{O}(\tau^{-1/3} \widetilde{\Lambda}) \! + \! 
\mathcal{O} \big(\tau^{2/3}(\hat{h}_{0}(\tau))^{2} \widetilde{\mu}^{3} \big),
\end{align}
where
\begin{equation} \label{eq3.40} 
C^{\scriptscriptstyle \mathrm{WKB}}_{0,k} \! := \! -\mi (a \! - \! \mi/2) \ln ((\sqrtsign{3} \! + \! 1)/\sqrtsign{2}),
\end{equation}
and $\mathfrak{c}_{m,k}$, $m \! = \! 1,2,\dotsc,6$, are $\mathcal{O}(1)$.
\end{ffff}

\emph{Proof}. Consequence of Corollary \ref{cor3.1.2}, equation \eqref{eq3.22}, upon choosing consistently 
the corresponding br\-a\-n\-c\-h\-e\-s in equations \eqref{equpsi} and \eqref{equpsisharp} and taking the limits 
$\widetilde{\mu} \! \to \! \infty$ and $\widetilde{\mu} \! \to \! 0$: the error estimate $\mathcal{O}(\mathscr{E}_{k}
(\xi))$ in equation \eqref{eq3.22} is given in equation \eqref{eqeel}; in particular, {}from the last two lines of 
equation \eqref{eqeel},
\begin{equation*}
\mathcal{O} \big(\tau^{2/3} \mathscr{E}_{k}(\widetilde{\mu}) \big) 
\underset{\underset{\widetilde{\mu} \to \infty}{\tau \to +\infty}}{=} \mathcal{O} \big(\tau^{-2/3} \widetilde{\mu}^{-3} \big) 
\, \qquad \, \text{and} \, \quad \, \mathcal{O} \big(\tau^{2/3} \mathscr{E}_{k}(\widetilde{\mu}) \big) 
\underset{\underset{\widetilde{\mu} \to 0}{\tau \to +\infty}}{=} \mathcal{O} \big(\tau^{2/3}(\hat{h}_{0}(\tau))^{2} 
\widetilde{\mu}^{3} \big),
\end{equation*}
which implies the results stated in the corollary. \hfill $\qed$
\begin{bbbb} \label{prop3.1.4} 
Let $T(\widetilde{\mu})$ be given in equation~\eqref{eq3.18}, with $\mathcal{A}(\widetilde{\mu})$ defined by equation 
\eqref{eq3.4} and $l_{k}^{2}(\widetilde{\mu})$, $k \! = \! \pm 1$, given in equation \eqref{eq3.19}$;$ then,
\begin{equation}
\int_{\widetilde{\mu}_{0,k}}^{\widetilde{\mu}} \diag \! \left((T(\xi))^{-1} \partial_{\xi}T(\xi) \right) \md \xi 
\underset{\tau \to +\infty}{=} \left(\mathcal{I}_{\tau,k}(\widetilde{\mu}) \! + \! \mathcal{O}
(\mathcal{E}_{\scriptscriptstyle T,k}(\widetilde{\mu})) \! + \! \mathcal{O}(\mathcal{E}_{\scriptscriptstyle T,k}
(\widetilde{\mu}_{0,k})) \right) \! \sigma_{3}, \label{eqintTminus1T}
\end{equation}
where, for $\delta \! > \! 0$, $\widetilde{\mu},\widetilde{\mu}_{0,k} \! \in \! \mathbb{C} \setminus 
(\mathscr{O}_{\tau^{-1/3+ \delta}}(\pm \alpha_{k}) \cup \mathscr{O}_{\tau^{-2/3+2 \delta}}(\pm \mi \sqrt{2} 
\alpha_{k}) \cup \lbrace 0,\infty \rbrace)$ and the path of integration lies in the corresponding 
canonical domain,
\begin{equation} \label{eqItee} 
\mathcal{I}_{\tau,k}(\widetilde{\mu}) \! = \! \mathfrak{p}_{k}(\tau)(\digamma_{\tau,k}(\widetilde{\mu}) \! - \! 
\digamma_{\tau,k}(\widetilde{\mu}_{0,k})),
\end{equation}
with
\begin{equation}
\mathfrak{p}_{k}(\tau) \! := \! \dfrac{\alpha_{k}^{2} \left(-2 \! + \! \tilde{r}_{0}(\tau) \tau^{-1/3} \! + \! 2(1 \! + \! v_{0}(\tau) 
\tau^{-1/3})^{2} \right) \! - \! (a \! - \! \mi/2) \tau^{-2/3}}{8(-2 \! + \! \tilde{r}_{0}(\tau) \tau^{-1/3})(1 \! + \! v_{0}(\tau) 
\tau^{-1/3})}, \label{eqpeetee}
\end{equation}
\begin{equation}
\digamma_{\tau,k}(\xi) \! := \! \dfrac{2}{\xi^{2} \! - \! \alpha_{k}^{2}} \! + \! \dfrac{2}{3 \sqrt{3} \alpha_{k}^{2}} \ln \! \left(
\left(\dfrac{3^{1/2}(\xi^{2} \! + \! 2 \alpha_{k}^{2})^{1/2} \! - \! \xi \! + \! 2 \alpha_{k}}{3^{1/2}(\xi^{2} \! + \! 2 \alpha_{k}^{2})^{1/2} 
\! + \! \xi \! + \! 2 \alpha_{k}} \right) \! \left(\dfrac{\xi \! - \! \alpha_{k}}{\xi \! + \! \alpha_{k}} \right) \right) \! - \! \dfrac{2}{3 
\alpha_{k}^{2}} \dfrac{\xi (\xi^{2} \! + \! 2 \alpha_{k}^{2})^{1/2}}{\xi^{2} \! - \! \alpha_{k}^{2}}, \label{eqFtee}
\end{equation}
and
\begin{equation} \label{eqee2} 
\mathcal{E}_{\scriptscriptstyle T,k}(\xi) \! := \! 
\begin{cases}
\mathfrak{p}_{k}(\tau) \! \left(\frac{\mathfrak{c}_{1,k}^{\blacklozenge} \tilde{r}_{0}(\tau) \tau^{-1/3} + 
\mathfrak{c}_{2,k}^{\blacklozenge} \hat{\mathfrak{f}}_{1,k}(\tau)}{(\xi \mp \alpha_{k})^{2}} \! + \! \frac{
\mathfrak{c}_{3,k}^{\blacklozenge} \tilde{r}_{0}(\tau) \tau^{-1/3}}{\xi \mp \alpha_{k}} \right), & \text{$\xi \! \in \! \mathbb{U}_{k}^{1}$,} \\
\mathfrak{p}_{k}(\tau) \hat{\mathfrak{f}}_{3,k}(\tau) \! \left(\frac{\mathfrak{c}_{4,k}^{\blacklozenge}}{(\xi \mp \mi \sqrt{2} 
\alpha_{k})^{1/2}} \! + \! \mathfrak{c}_{5,k}^{\blacklozenge} \ln (\xi \! \mp \! \mi \sqrt{2} \alpha_{k}) \right), & \text{$\xi \! \in \! 
\mathbb{U}_{k}^{2}$,} \\
\mathfrak{p}_{k}(\tau) \xi^{-4} \! \left(\mathfrak{c}_{6,k}^{\blacklozenge} \tilde{r}_{0}(\tau) \tau^{-1/3} \! + \! \mathcal{O} 
\big((\mathfrak{c}_{7,k}^{\blacklozenge} \tilde{r}_{0}(\tau) \tau^{-1/3} \! + \! \mathfrak{c}_{8,k}^{\blacklozenge} \tau^{-2/3}) 
\xi^{-2} \big) \right), & \text{$\xi \! \to \! \infty$,} \\
\mathfrak{p}_{k}(\tau) \tilde{r}_{0}(\tau) \tau^{-1/3} \xi^{2} \big(\mathfrak{c}_{9,k}^{\blacklozenge} \! + \! \mathcal{O}(\xi) \big), 
& \text{$\xi \! \to \! 0$,}
\end{cases}
\end{equation}
where $\mathbb{U}_{k}^{1} \! := \! \mathscr{O}_{\tau^{-1/3+ \delta_{k}}}(\pm \alpha_{k})$, $\mathbb{U}_{k}^{2} \! := \! 
\mathscr{O}_{\tau^{-2/3+2 \delta_{k}}}(\pm \mi \sqrt{2} \alpha_{k})$, the parameter $\delta_{k}$ satisfies (cf. Corollary 
\ref{cor3.1.2}$)$ $0 \! < \! \delta \! < \! \delta_{k} \! < \! 1/9$, the functions $\hat{\mathfrak{f}}_{1,k}(\tau)$ and 
$\hat{\mathfrak{f}}_{3,k}(\tau)$ are given in equation \eqref{mfkeff}, and $\mathfrak{c}_{m,k}^{\blacklozenge}$, 
$m \! = \! 1,2,\dotsc,9$, are $\mathcal{O}(1)$.
\end{bbbb}

\emph{Proof}. {}From equations \eqref{eq3.4}, \eqref{iden3}, and \eqref{eqlsquared}--\eqref{eqforl}, one shows that
\begin{align}
2l_{k}(\xi)(\mi \mathcal{A}_{11}(\xi) \! + \! l_{k}(\xi)) \underset{\tau \to +\infty}{=}& 
\, \mathcal{P}_{\infty,k}(\xi) \! + \! \mathcal{P}_{1,k}(\xi) \Delta_{k,\tau}(\xi) \! + 
\! \mathcal{O} \! \left(l_{k,\infty}^{2}(\xi) \Delta^{2}_{k,\tau}(\xi) \right) \nonumber \\
+& \, \mathcal{O} \! \left(l_{k,\infty}(\xi) \Delta^{2}_{k,\tau}(\xi) \! \left(2 \xi 
\! + \! \dfrac{(\varepsilon b)^{1/3}}{2 \xi}(-2 \! + \! \hat{r}_{0}(\tau)) \right) 
\right), \label{eqdenomTminus1T}
\end{align}
where
\begin{align}
\mathcal{P}_{\infty,k}(\xi) :=& \, 2l_{k,\infty}^{2}(\xi) \! + \! 2l_{k,\infty}(\xi) 
\! \left(2 \xi \! + \! \dfrac{(\varepsilon b)^{1/3}}{2 \xi}(-2 \! + \! \hat{r}_{0}
(\tau)) \right), \label{eq3.24} \\
\mathcal{P}_{1,k}(\xi) :=& \, 2l_{k,\infty}^{2}(\xi) \! + \! l_{k,\infty}(\xi) \! 
\left(2 \xi \! + \! \dfrac{(\varepsilon b)^{1/3}}{2 \xi}(-2 \! + \! \hat{r}_{0}
(\tau)) \right), \label{eq3.25}
\end{align}
and, via equations \eqref{eq3.4}, \eqref{iden7}, \eqref{iden3}, and \eqref{iden4},
\begin{align} \label{eqnumTminus1T} 
\mathcal{A}_{12}(\xi) \partial_{\xi} \mathcal{A}_{21}(\xi) \! - \! \mathcal{A}_{21}(\xi) \partial_{\xi} \mathcal{A}_{12}
(\xi)=& \, -\frac{4(\varepsilon b)^{2/3}}{\xi^{3}} \left(\dfrac{2(1 \! + \! \hat{u}_{0}(\tau))^{2} \! + \! (-2 \! + \! 
\hat{r}_{0}(\tau))}{2(1 \! + \! \hat{u}_{0}(\tau))} \right) \nonumber \\
+& \, \frac{4(\varepsilon b)^{1/3}(a \! - \! \mi/2) \tau^{-2/3}}{\xi^{3}(1 \! + \! \hat{u}_{0}(\tau))}.
\end{align}
Substituting equations \eqref{eqdenomTminus1T} and \eqref{eqnumTminus1T} into equation \eqref{eq3.23} and 
expanding $(2l_{k}(\xi)(\mi \mathcal{A}_{11}(\xi) \! + \! l_{k}(\xi)))^{-1}$ into a series of powers of $\Delta_{k,\tau}(\xi)$, 
one arrives at (cf. equation \eqref{eq3.16})
\begin{align} \label{iden31} 
\int_{\widetilde{\mu}_{0,k}}^{\widetilde{\mu}} \diag \! \left((T(\xi))^{-1} \partial_{\xi}T(\xi) \right) \md \xi 
\underset{\tau \to +\infty}{=}& \, \left(\varkappa_{k}(\tau) \int_{\widetilde{\mu}_{0,k}}^{\widetilde{\mu}} \dfrac{1}{\xi^{3} 
\mathcal{P}_{\infty,k}(\xi)} \, \md \xi \right. \nonumber \\
+&\left. \, \mathcal{O} \! \left(\varkappa_{k}(\tau) \int_{\widetilde{\mu}_{0,k}}^{\widetilde{\mu}} \dfrac{\xi^{3} 
\mathcal{P}_{1,k}(\xi) \Delta_{k,\tau}(\xi)}{(\xi^{3} \mathcal{P}_{\infty,k}(\xi))^{2}} \, \md \xi \right) \right) \! \sigma_{3},
\end{align}
where
\begin{equation} \label{iden30} 
\varkappa_{k}(\tau) \! := \! (\varepsilon b)^{2/3} \! \left(\dfrac{2(1 \! + 
\! \hat{u}_{0}(\tau))^{2} \! + \! (-2 \! + \! \hat{r}_{0}(\tau))}{1 \! + \! 
\hat{u}_{0}(\tau)} \right) \! - \! \dfrac{2(\varepsilon b)^{1/3}(a \! - \! 
\mi/2) \tau^{-2/3}}{1 \! + \! \hat{u}_{0}(\tau)}.
\end{equation}
Via equations \eqref{eqlsquared} and \eqref{eq3.24}, a calculation reveals that
\begin{equation}
\dfrac{\varkappa_{k}(\tau)}{\xi^{3} \mathcal{P}_{\infty,k}(\xi)} \! = \! 
\mathfrak{p}_{k}(\tau) \! \left(\dfrac{\xi \left(\xi (4 \xi^{2} \! + \! 
(\varepsilon b)^{1/3}(-2 \! + \! \hat{r}_{0}(\tau))) \! - \! 4(\xi^{2} \! - \! 
\alpha_{k}^{2})(\xi^{2} \! + \! 2 \alpha_{k}^{2})^{1/2} \right)}{(\xi^{2} \! 
- \!\alpha_{k}^{2})(\xi^{2} \! + \! 2 \alpha_{k}^{2})^{1/2}(\xi^{2} \! + \! 
\hat{\mathfrak{z}}_{k}^{+}(\tau))(\xi^{2} \! + \! \hat{\mathfrak{z}}_{k}^{-}
(\tau))} \right), \label{iden32}
\end{equation}
where $\mathfrak{p}_{k}(\tau)$ is defined by equation \eqref{eqpeetee}, and
\begin{equation} \label{iden33} 
\hat{\mathfrak{z}}_{k}^{\pm}(\tau) \! := \! \dfrac{(\varepsilon b)^{1/3}}{4
(-2 \! + \! \hat{r}_{0}(\tau))} \! \left(\left(\dfrac{-2 \! + \! \hat{r}_{0}(\tau)}{2} 
\right)^{2} \! - \! 3 \me^{\mi \pi k/3} \! \mp \! \sqrt{\left(\left(\dfrac{-2 \! 
+ \! \hat{r}_{0}(\tau)}{2} \right)^{2} \! - \! 3 \me^{\mi \pi k/3} \right)^{2} \! 
+ \! 8(-2 \! + \! \hat{r}_{0}(\tau))} \, \right).
\end{equation}
One shows {}from equations \eqref{iden3} and \eqref{iden4}, the conditions \eqref{iden5}, and the definition \eqref{iden33} that
\begin{align} \label{iden34} 
\hat{\mathfrak{z}}_{k}^{\pm}(\tau) \underset{\tau \to +\infty}{=}& \, \dfrac{(\varepsilon b)^{1/3} \me^{-\mi \pi k/3}}{2} 
\left(1 \! + \! \left(\dfrac{1 \! \pm \! \sqrt{3}}{4} \right) \tilde{r}_{0}(\tau) \tau^{-1/3} \! + \! \left(\dfrac{3\sqrt{3} \! \pm \! 
5}{16 \sqrt{3}} \right) \big(\tilde{r}_{0}(\tau) \tau^{-1/3} \big)^{2} \right. \nonumber \\
+&\left. \, \mathcal{O} \big((\tilde{r}_{0}(\tau) \tau^{-1/3})^{3} \big) \right),
\end{align}
whence, via equation \eqref{iden32}, the first term on the right-hand side of equation \eqref{iden31} can be presented as follows:
\begin{equation} \label{iden35} 
\varkappa_{k}(\tau) \int_{\widetilde{\mu}_{0,k}}^{\widetilde{\mu}} \dfrac{1}{
\xi^{3} \mathcal{P}_{\infty,k}(\xi)} \, \md \xi \underset{\tau \to +\infty}{=} 
\mathcal{I}_{\tau,k}(\widetilde{\mu}) \! + \! \mathcal{I}_{{\scriptscriptstyle A},k}
(\widetilde{\mu}) \! + \! \mathcal{O}(\mathcal{I}_{{\scriptscriptstyle B},k}
(\widetilde{\mu})),
\end{equation}
where
\begin{gather}
\mathcal{I}_{\tau,k}(\widetilde{\mu}) \! := \! \mathfrak{p}_{k}(\tau) \int_{\widetilde{\mu}_{0,k}}^{\widetilde{\mu}} 
\left(\dfrac{4 \xi^{2}(\xi^{2} \! + \! 2 \alpha_{k}^{2})^{1/2}}{(\xi^{2} \! + \! 2 \alpha_{k}^{2})(\xi^{2} \! - \! 
\alpha_{k}^{2})^{2}} \! - \! \dfrac{4 \xi}{(\xi^{2} \! - \! \alpha_{k}^{2})^{2}} \right) \md \xi, \label{iden36} \\
\mathcal{I}_{{\scriptscriptstyle A},k}(\widetilde{\mu}) \! := \! \mathfrak{p}_{k}(\tau) \tilde{r}_{0}(\tau) \tau^{-1/3} 
\int_{\widetilde{\mu}_{0,k}}^{\widetilde{\mu}} \left(\dfrac{4 \alpha_{k}^{2} \xi^{2}(\xi^{2} \! + \! 2 \alpha_{k}^{2}
)^{1/2}}{(\xi^{2} \! + \! 2 \alpha_{k}^{2})(\xi^{2} \! - \! \alpha_{k}^{2})^{3}} \! - \! \dfrac{2 \alpha_{k}^{2} \xi}{(\xi^{2} 
\! - \! \alpha_{k}^{2})^{3}} \right) \md \xi, \label{iden37} \\
\mathcal{I}_{{\scriptscriptstyle B},k}(\widetilde{\mu}) \! := \! \mathfrak{p}_{k}(\tau) \big(\tilde{r}_{0}(\tau) \tau^{-1/3} 
\big)^{2} \int_{\widetilde{\mu}_{0,k}}^{\widetilde{\mu}} \left(\dfrac{\alpha_{k}^{4} \xi^{2}(\xi^{2} \! + \! 2 \alpha_{k}^{2}
)^{1/2}}{(\xi^{2} \! + \! 2 \alpha_{k}^{2})(\xi^{2} \! - \! \alpha_{k}^{2})^{4}} \! - \!  \dfrac{4 \xi^{3}}{(\xi^{2} \! - \! 
\alpha_{k}^{2})^{4}} \! + \! \dfrac{4 \xi^{4}(\xi^{2} \! + \! 2 \alpha_{k}^{2})^{1/2}}{(\xi^{2} \! + \! 2 \alpha_{k}^{2})
(\xi^{2} \! - \! \alpha_{k}^{2})^{4}} \right) \md \xi. \label{iden38}
\end{gather}
A partial fraction decomposition shows that
\begin{equation} \label{iden39} 
\dfrac{\xi^{2}}{(\xi^{2} \! + \! 2 \alpha_{k}^{2})(\xi^{2} \! - \! \alpha_{k}^{2})^{2}} \! = \! \dfrac{\alpha_{k}^{-3}}{36} 
\dfrac{1}{\xi \! - \! \alpha_{k}} \! + \! \dfrac{\alpha_{k}^{-2}}{12} \dfrac{1}{(\xi \! - \! \alpha_{k})^{2}} \! - \! 
\dfrac{\alpha_{k}^{-3}}{36} \dfrac{1}{\xi \! + \! \alpha_{k}} \! + \! \dfrac{\alpha_{k}^{-2}}{12} \dfrac{1}{(\xi \! + \! 
\alpha_{k})^{2}} \! - \! \dfrac{2 \alpha_{k}^{-2}}{9} \dfrac{1}{\xi^{2} \! + \! 2 \alpha_{k}^{2}};
\end{equation}
substituting equation \eqref{iden39} into equation \eqref{iden36} and integrating, one arrives at equations 
\eqref{eqItee}--\eqref{eqFtee}.

Equations \eqref{iden37} and \eqref{iden38} contribute to the error function, $\mathcal{E}_{\scriptscriptstyle T,k}
(\mathbf{\cdot})$, in equation \eqref{eqintTminus1T}; therefore, only its asymptotics at the turning and the singular 
points are requisite. Evaluating the integrals in equations \eqref{iden37} and \eqref{iden38}, one shows that
\begin{equation} \label{iden40} 
\mathcal{I}_{{\scriptscriptstyle A},k}(\widetilde{\mu}) 
\underset{\tau \to +\infty}{=} 
\begin{cases}
\mathfrak{p}_{k}(\tau) \tilde{r}_{0}(\tau) \tau^{-1/3} \big(\hat{\mathfrak{h}}_{1,k}(\widetilde{\mu}) \! - \! \hat{\mathfrak{h}}_{1,k}
(\widetilde{\mu}_{0,k}) \big), & \text{$\widetilde{\mu} \! \in \! 
\mathscr{O}_{\tau^{-1/3+ \delta_{k}}}(\pm \alpha_{k})$,} \\\mathfrak{p}_{k}(\tau) \tilde{r}_{0}(\tau) \tau^{-1/3} \big(
\hat{\mathfrak{h}}_{2,k}(\widetilde{\mu}) \! - \! \hat{\mathfrak{h}}_{2,k}(\widetilde{\mu}_{0,k}) \big), & 
\text{$\widetilde{\mu} \! \in \! \mathscr{O}_{\tau^{-2/3+ 2 \delta_{k}}}(\pm \mi \sqrt{2} \alpha_{k})$,} \\
\mathfrak{p}_{k}(\tau) \tilde{r}_{0}(\tau) \tau^{-1/3} \big(\hat{\mathfrak{h}}_{3,k}(\widetilde{\mu}) \! - \! \hat{\mathfrak{h}}_{3,k}
(\widetilde{\mu}_{0,k}) \big), & \text{$\widetilde{\mu} \! \to \! \infty$,} \\
\mathfrak{p}_{k}(\tau) \tilde{r}_{0}(\tau) \tau^{-1/3} \big(\hat{\mathfrak{h}}_{4,k}(\widetilde{\mu}) \! - \! \hat{\mathfrak{h}}_{4,k}
(\widetilde{\mu}_{0,k}) \big), & \text{$\widetilde{\mu} \! \to \! 0$,}
\end{cases}
\end{equation}
where
\begin{gather*}
\hat{\mathfrak{h}}_{1,k}(\xi) \! := \! c_{1,k}^{\flat}(\xi \! \mp \! \alpha_{k})^{-2} \! + \! c_{2,k}^{\flat}(\xi \! \mp  \! \alpha_{k})^{-1} 
\! + \! c_{3,k}^{\flat} \ln (\xi \! \mp \! \alpha_{k}) \! + \! \sum_{m \in \mathbb{Z}_{+}}d_{m,k}^{\flat}(\xi \! \mp \! \alpha_{k})^{m}, \\
\hat{\mathfrak{h}}_{2,k}(\xi) \! := \! (\xi \! \mp \! \mi \sqrt{2} \alpha_{k})^{1/2} \sum_{m \in \mathbb{Z}_{+}}c_{m,k}^{\natural}(\xi 
\! \mp \! \mi \sqrt{2} \alpha_{k})^{m} \! + \! \sum_{m \in \mathbb{Z}_{+}}d_{m,k}^{\natural}(\xi \! \mp \! \mi \sqrt{2} \alpha_{k})^{m}, \\
\hat{\mathfrak{h}}_{3,k}(\xi) \! := \! \xi^{-4} \sum_{m \in \mathbb{Z}_{+}}c_{m,k}^{\sharp,\infty} \xi^{-2m}, \, \quad \, \quad \, 
\hat{\mathfrak{h}}_{4,k}(\xi) \! := \! \xi^{2} \sum_{m \in \mathbb{Z}_{+}}c_{m,k}^{\sharp,0} \xi^{m},
\end{gather*}
and $c_{1,k}^{\flat}$, $c_{2,k}^{\flat}$, $c_{3,k}^{\flat}$, $d_{m,k}^{\flat}$, $c_{m,k}^{\natural}$, $d_{m,k}^{\natural}$, 
$c_{m,k}^{\sharp,\infty}$, and $c_{m,k}^{\sharp,0}$ are $\mathcal{O}(1)$, and
\begin{equation} \label{iden41} 
\mathcal{I}_{{\scriptscriptstyle B},k}(\widetilde{\mu}) \underset{\tau \to +\infty}{=} 
\begin{cases}
\mathfrak{p}_{k}(\tau)(\tilde{r}_{0}(\tau) \tau^{-1/3})^{2} \big(\hat{\mathfrak{h}}_{5,k}(\widetilde{\mu}) \! - \! \hat{\mathfrak{h}}_{5,k}
(\widetilde{\mu}_{0,k}) \big), & \text{$\widetilde{\mu} \! \in \! \mathscr{O}_{\tau^{-1/3+ \delta_{k}}}(\pm \alpha_{k})$,} \\
\mathfrak{p}_{k}(\tau)(\tilde{r}_{0}(\tau) \tau^{-1/3})^{2} \big(\hat{\mathfrak{h}}_{6,k}(\widetilde{\mu}) \! - \! \hat{\mathfrak{h}}_{6,k}
(\widetilde{\mu}_{0,k}) \big), & \text{$\widetilde{\mu} \! \in \! \mathscr{O}_{\tau^{-2/3+ 2 \delta_{k}}}(\pm \mi \sqrt{2} \alpha_{k})$,} \\
\mathfrak{p}_{k}(\tau)(\tilde{r}_{0}(\tau) \tau^{-1/3})^{2} \big(\hat{\mathfrak{h}}_{7,k}(\widetilde{\mu}) \! - \! \hat{\mathfrak{h}}_{7,k}
(\widetilde{\mu}_{0,k}) \big), & \text{$\widetilde{\mu} \! \to \! \infty$,} \\
\mathfrak{p}_{k}(\tau)(\tilde{r}_{0}(\tau) \tau^{-1/3})^{2} \big(\hat{\mathfrak{h}}_{8,k}(\widetilde{\mu}) \! - \! \hat{\mathfrak{h}}_{8,k}
(\widetilde{\mu}_{0,k}) \big), & \text{$\widetilde{\mu} \! \to \! 0$,}
\end{cases}
\end{equation}
where
\begin{gather*}
\hat{\mathfrak{h}}_{5,k}(\xi) \! := \! \hat{c}_{1,k}^{\flat}(\xi \! \mp \! \alpha_{k})^{-3} \! + \! \hat{c}_{2,k}^{\flat}
(\xi \! \mp  \! \alpha_{k})^{-2} \! + \! \hat{c}_{3,k}^{\flat}(\xi \! \mp  \! \alpha_{k})^{-1} \! + \! \hat{c}_{4,k}^{\flat} 
\ln (\xi \! \mp \! \alpha_{k}) \! + \! \sum_{m \in \mathbb{Z}_{+}} \hat{d}_{m,k}^{\flat}(\xi \! \mp \! \alpha_{k})^{m}, \\
\hat{\mathfrak{h}}_{6,k}(\xi) \! := \! (\xi \! \mp \! \mi \sqrt{2} \alpha_{k})^{1/2} \sum_{m \in \mathbb{Z}_{+}} 
\hat{c}_{m,k}^{\natural}(\xi \! \mp \! \mi \sqrt{2} \alpha_{k})^{m} \! + \! \sum_{m \in \mathbb{Z}_{+}} 
\hat{d}_{m,k}^{\natural}(\xi \! \mp \! \mi \sqrt{2} \alpha_{k})^{m}, \\
\hat{\mathfrak{h}}_{7,k}(\xi) \! := \! \xi^{-6} \sum_{m \in \mathbb{Z}_{+}} \hat{c}_{m,k}^{\sharp,\infty} \xi^{-2m}, \, \quad 
\, \quad \, \hat{\mathfrak{h}}_{8,k}(\xi) \! := \! \xi^{3} \sum_{m \in \mathbb{Z}_{+}} \hat{c}_{m,k}^{\sharp,0} \xi^{m},
\end{gather*}
and $\hat{c}_{1,k}^{\flat}$, $\hat{c}_{2,k}^{\flat}$, $\hat{c}_{3,k}^{\flat}$, $\hat{c}_{4,k}^{\flat}$, $d_{m,k}^{\flat}$, 
$\hat{c}_{m,k}^{\natural}$, $\hat{d}_{m,k}^{\natural}$, $\hat{c}_{m,k}^{\sharp,\infty}$, and $\hat{c}_{m,k}^{\sharp,0}$ 
are $\mathcal{O}(1)$.

One now estimates the second term on the right-hand side of equation \eqref{iden31}. {}From equations 
\eqref{eqlsquared}---\eqref{eqforl}, it follows, after simplification, that
\begin{align} \label{iden43} 
\int_{\widetilde{\mu}_{0,k}}^{\widetilde{\mu}} \dfrac{\xi^{3} 
\mathcal{P}_{1,k}(\xi) \Delta_{k,\tau}(\xi)}{(\xi^{3} \mathcal{P}_{\infty,k}
(\xi))^{2}} \, \md \xi \! = \! \int_{\widetilde{\mu}_{0,k}}^{\widetilde{\mu}}& 
\, \dfrac{\xi \left(\xi (4 \xi^{2} \! + \! (\varepsilon b)^{1/3}(-2 \! + \! 
\hat{r}_{0}(\tau))) \! + \! 8(\xi^{2} \! - \! \alpha_{k}^{2})(\xi^{2} \! + \! 2 
\alpha_{k}^{2})^{1/2} \right)}{\left(\xi (4 \xi^{2} \! + \! (\varepsilon b)^{1/3}
(-2 \! + \! \hat{r}_{0}(\tau))) \! + \! 4(\xi^{2} \! - \! \alpha_{k}^{2})(\xi^{2} 
\! + \! 2 \alpha_{k}^{2})^{1/2} \right)^{2}} \nonumber \\
\times& \, \dfrac{(\xi^{2} \hat{h}_{0}(\tau) \! + \! \xi^{4}(a \! - \! \mi/2) 
\tau^{-2/3})}{4(\xi^{2} \! - \! \alpha_{k}^{2})^{3}(\xi^{2} \! + \! 2 
\alpha_{k}^{2})^{3/2}} \, \md \xi.
\end{align}
Evaluating the integral in equation \eqref{iden43}, a lengthy calculation shows that its asymptotics at the turning and the singular 
points are given by
\begin{equation} \label{iden47} 
\varkappa_{k}(\tau) \int_{\widetilde{\mu}_{0,k}}^{\widetilde{\mu}} 
\dfrac{\xi^{3} \mathcal{P}_{1,k}(\xi) \Delta_{k,\tau}(\xi)}{(\xi^{3} 
\mathcal{P}_{\infty,k}(\xi))^{2}} \, \md \xi \underset{\tau \to +\infty}{=} 
\begin{cases}
\hat{\mathfrak{h}}_{9,k}(\widetilde{\mu}) \! - \! \hat{\mathfrak{h}}_{9,k}
(\widetilde{\mu}_{0,k}), & \text{$\widetilde{\mu} \! \in \! 
\mathscr{O}_{\tau^{-1/3+ \delta_{k}}}(\pm \alpha_{k})$,} \\
\hat{\mathfrak{h}}_{10,k}(\widetilde{\mu}) \! - \! \hat{\mathfrak{h}}_{10,k}
(\widetilde{\mu}_{0,k}), & \text{$\widetilde{\mu} \! \in \! 
\mathscr{O}_{\tau^{-2/3+ 2 \delta_{k}}}(\pm \mi \sqrt{2} \alpha_{k})$,} \\
\hat{\mathfrak{h}}_{11,k}(\widetilde{\mu}) \! - \! \hat{\mathfrak{h}}_{11,k}
(\widetilde{\mu}_{0,k}), & \text{$\widetilde{\mu} \! \to \! \infty$,} \\
\hat{\mathfrak{h}}_{12,k}(\widetilde{\mu}) \! - \! \hat{\mathfrak{h}}_{12,k}
(\widetilde{\mu}_{0,k}), & \text{$\widetilde{\mu} \! \to \! 0$,}
\end{cases}
\end{equation}
where
\begin{align*}
\hat{\mathfrak{h}}_{9,k}(\xi) :=& \, \tilde{c}_{1,k}^{\sharp} \mathfrak{p}_{k}(\tau) 
\hat{\mathfrak{f}}_{1,k}(\tau)(\xi \! \mp \! \alpha_{k})^{-2} \! + \! 
\mathfrak{p}_{k}(\tau)(\tilde{c}_{2,k}^{\sharp} \hat{\mathfrak{f}}_{2,k}
(\tau) \! + \! \tilde{c}_{3,k}^{\sharp} \tilde{r}_{0}(\tau) \tau^{-1/3} 
\hat{\mathfrak{f}}_{1,k}(\tau))(\xi \! \mp \! \alpha_{k})^{-3}, \\
\hat{\mathfrak{h}}_{10,k}(\xi) :=& \, \tilde{c}_{4,k}^{\sharp} \mathfrak{p}_{k}(\tau) 
\hat{\mathfrak{f}}_{3,k}(\tau)(\xi \! \mp \! \mi \sqrt{2} \alpha_{k})^{-1/2} \! 
+ \! \tilde{c}_{5,k}^{\sharp} \mathfrak{p}_{k}(\tau) \hat{\mathfrak{f}}_{3,k}
(\tau) \ln (\xi \! \mp \! \mi \sqrt{2} \alpha_{k}), \\
\hat{\mathfrak{h}}_{11,k}(\xi) :=& \, \mathfrak{p}_{k}(\tau) \tau^{-2/3} \xi^{-6} 
\! \left(\tilde{c}_{6,k}^{\sharp} \! + \! \xi^{-2}(\tilde{c}_{7,k}^{\sharp} \! + \! 
\tilde{c}_{8,k}^{\sharp} \tau^{2/3} \hat{h}_{0}(\tau)) \right. \\
+&\left. \mathcal{O} \! \left(\tilde{r}_{0}(\tau) \tau^{-1/3}(\tilde{c}_{9,k}^{\sharp} 
\! + \! \xi^{-2}(\tilde{c}_{10,k}^{\sharp} \! + \! \tilde{c}_{11,k}^{\sharp} \tau^{2/3} 
\hat{h}_{0}(\tau))) \right) \right), \\
\hat{\mathfrak{h}}_{12,k}(\xi) :=& \, \mathfrak{p}_{k}(\tau) \tau^{-2/3} \xi^{4} 
\! \left(\tilde{c}_{12,k}^{\sharp} \tau^{2/3} \hat{h}_{0}(\tau) \! + \! \xi 
\tilde{c}_{13,k}^{\sharp} \tau^{2/3} \hat{h}_{0}(\tau) \! + \! \xi^{2}
(\tilde{c}_{14,k}^{\sharp} \! + \! \tilde{c}_{15,k}^{\sharp} \tau^{2/3} 
\hat{h}_{0}(\tau)) \right. \\
+&\left. \mathcal{O} \! \left(\tilde{r}_{0}(\tau) \tau^{-1/3}
(\tilde{c}_{16,k}^{\sharp} \tau^{2/3} \hat{h}_{0}(\tau) \! + \! \xi 
\tilde{c}_{17,k}^{\sharp} \tau^{2/3} \hat{h}_{0}(\tau) \! + \! \xi^{2}
(\tilde{c}_{18,k}^{\sharp} \! + \! \tilde{c}_{19,k}^{\sharp} \tau^{2/3} 
\hat{h}_{0}(\tau))) \right) \right),
\end{align*}
and $\tilde{c}_{m,k}^{\sharp}$, $m \! = \! 1,2,\dotsc,19$, are $\mathcal{O}(1)$, and
\begin{equation} \label{mfkeff} 
\hat{\mathfrak{f}}_{j,k}(\tau) \! = \! \left((a \! - \! \mi/2) \! + \! \dfrac{2 \hat{s}(j) \hat{h}_{0}(\tau) 
\tau^{2/3}}{(3 \! + \! (-1)^{j+1}) \alpha_{k}^{2}} \right) \! \tau^{-2/3}, \quad j \! = \! 1,2,3,
\end{equation}
where $\hat{s}(1) \! = \! \hat{s}(2) \! = \! +1$ and $\hat{s}(3) \! = \! -1$. Thus, assembling the error estimates 
\eqref{iden40}, \eqref{iden41}, and \eqref{iden47}, and retaining only leading-order terms, one arrives at the 
error function defined by equation \eqref{eqee2}. \hfill $\qed$
\begin{ffff} \label{cor3.1.3} 
Set $\widetilde{\mu}_{0,k} \! = \! \alpha_{k} \! + \! \tau^{-1/3} \widetilde{\Lambda}$, $k \! = \! \pm 1$, where 
$\widetilde{\Lambda} \genfrac{}{}{0pt}{3}{=}{\tau \to +\infty} \mathcal{O}(\tau^{\delta_{k}})$, $0 \! < \! \delta 
\! < \! \delta_{k} \! < \! 1/9;$ then,
\begin{align}
\int_{\widetilde{\mu}_{0,k}}^{\widetilde{\mu}} \diag \! \left((T(\xi))^{-1} \partial_{\xi}T(\xi) \right) \md \xi 
\underset{\tau \to +\infty}{=}& \, \left(\mathfrak{p}_{k}(\tau)(\digamma_{\tau,k}(\widetilde{\mu}) \! + \! 
\digamma_{\tau,k}^{\sharp}(\tau)) \! + \! \mathcal{O}(\mathcal{E}_{\scriptscriptstyle T,k}(\widetilde{\mu})) 
\right. \nonumber \\
+&\left. \, \mathcal{O} \! \left((\mathfrak{c}_{3,k} \tau^{-1/3} \! + \! \mathfrak{c}_{4,k}(\tilde{r}_{0}(\tau) 
\! + \! 4v_{0}(\tau))) \right. \right. \nonumber \\
\times&\left. \left. \left(\dfrac{\mathfrak{c}_{1,k} \tau^{-1/3} \! + \! \mathfrak{c}_{2,k} \tilde{r}_{0}(\tau)}{
\widetilde{\Lambda}^{2}} \right) \right) \right) \! \sigma_{3}, \label{eq3.36}
\end{align}
where $\mathfrak{p}_{k}(\tau)$, $\digamma_{\tau,k}(\xi)$ and $\mathcal{E}_{\scriptscriptstyle T,k}(\xi)$ are 
defined by equations \eqref{eqpeetee}, \eqref{eqFtee}, and \eqref{eqee2}, respectively,
\begin{equation}
\digamma_{\tau,k}^{\sharp}(\tau) \! := \! -\dfrac{\tau^{1/3}}{\alpha_{k} 
\widetilde{\Lambda}} \! \left(\dfrac{\sqrt{3} \! \mp \! 1}{\sqrt{3}} \right) \! 
\mp \! \dfrac{2}{3 \sqrt{3} \alpha_{k}^{2}} \! \left(-\dfrac{1}{3} \ln \tau 
\! + \! \ln \widetilde{\Lambda} \right) \! \pm \! \dfrac{(5 \! \pm \! 
3 \sqrt{3})}{6 \sqrt{3} \alpha_{k}^{2}} \! \pm \! \dfrac{2}{3 \sqrt{3} 
\alpha_{k}^{2}} \ln (3 \alpha_{k}), \label{eqfteesharp}
\end{equation}
with the upper (resp., lower) signs taken according to the branch of the square-root function $\lim_{\xi^{2} \to +\infty}
(\xi^{2} \! + \! 2 \alpha_{k}^{2})^{1/2} \! = \! +\infty$ (resp., $\lim_{\xi^{2} \to +\infty}(\xi^{2} \! + \! 2 \alpha_{k}^{2})^{1/2} 
\! = \! -\infty)$, and $\mathfrak{c}_{m,k}$, $m \! = \! 1,2,3,4$, are $\mathcal{O}(1)$.
\end{ffff}

\emph{Proof}. Substituting $\widetilde{\mu}_{0,k}$, as given in the corollary, for the argument of the functions 
$\digamma_{\tau,k}(\xi)$ and $\mathcal{E}_{\scriptscriptstyle T,k}(\xi)$ (cf. equation \eqref{eqFtee} and the first 
line of equation \eqref{eqee2}, respectively) and expanding with respect to the small parameter $\tau^{-1/3} 
\widetilde{\Lambda}$, one arrives at the following estimates:
\begin{equation} \label{asympforf1} 
-\digamma_{\tau,k}(\widetilde{\mu}_{0,k}) \underset{\tau \to +\infty}{=} 
\digamma_{\tau,k}^{\sharp}(\tau) \! + \! \mathcal{O}(\tau^{-1/3} 
\widetilde{\Lambda}),
\end{equation}
where $\digamma_{\tau,k}^{\sharp}(\tau)$ is defined by equation \eqref{eqfteesharp}, and
\begin{equation} \label{asympforf2} 
\mathcal{O}(\mathcal{E}_{\scriptscriptstyle T,k}(\widetilde{\mu}_{0,k})) 
\underset{\tau \to +\infty}{=} \mathcal{O} \! \left(\dfrac{\mathfrak{p}_{k}
(\tau) \tilde{r}_{0}(\tau)}{\tau^{-1/3} \widetilde{\Lambda}^{2}} \right) \! + 
\! \mathcal{O} \! \left(\dfrac{\mathfrak{p}_{k}(\tau) \hat{\mathfrak{f}}_{1,k}
(\tau)}{\tau^{-2/3} \widetilde{\Lambda}^{2}} \right) \! + \! \mathcal{O} \! 
\left(\dfrac{\mathfrak{p}_{k}(\tau) \tilde{r}_{0}(\tau)}{\widetilde{\Lambda}} \right).
\end{equation}
{}From the conditions \eqref{iden5} and the definitions \eqref{eqpeetee} and \eqref{mfkeff} (for $j \! = \! 1)$, one shows that
\begin{equation} \label{asympforf3} 
\mathfrak{p}_{k}(\tau) \underset{\tau \to +\infty}{=} \mathfrak{p}_{k}^{\infty}
(\tau) \! + \! \mathcal{O}((\tilde{r}_{0}(\tau) \! - \! 2v_{0}(\tau)) \tau^{-1}) \! + 
\! \mathcal{O}(((\tilde{r}_{0}(\tau) \! - \! 2v_{0}(\tau))(\tilde{r}_{0}(\tau) \! + \! 
4v_{0}(\tau)) \! + \! 4v_{0}^{2}(\tau)) \tau^{-2/3}),
\end{equation}
where
\begin{equation} \label{peekayity} 
\mathfrak{p}_{k}^{\infty}(\tau) \! := \! \frac{\tau^{-1/3}}{16} \! \left(-\alpha_{k}^{2}(\tilde{r}_{0}(\tau) \! + \! 4v_{0}(\tau)) 
\! + \! (a \! - \! \mi/2) \tau^{-1/3} \right),
\end{equation}
and
\begin{align} \label{asympforf4} 
\hat{\mathfrak{f}}_{1,k}(\tau) \underset{\tau \to +\infty}{=} \tau^{-2/3} \! \left(\dfrac{1}{2}(a \! - \! \mi/2) \! + \! \mathcal{O}
(v_{0}(\tau) \tau^{-1/3}) \! + \! \mathcal{O} \! \left(8v_{0}^{2}(\tau) \! + \! 4v_{0}(\tau) \tilde{r}_{0}(\tau) \! - \! 
(\tilde{r}_{0}(\tau))^{2} \right) \right);
\end{align}
thus, {}from the conditions \eqref{iden5} and the asymptotics \eqref{asympforf2}--\eqref{asympforf4}, it follows that, for 
$c_{m,k}$, $m \! = \! 1,2,\dotsc,6$, that are $\mathcal{O}(1)$,
\begin{align} \label{asympforf5} 
\mathcal{O}(\mathcal{E}_{\scriptscriptstyle T,k}(\widetilde{\mu}_{0,k})) 
\underset{\tau \to +\infty}{=}& \,  \mathcal{O} \! \left(\left(\dfrac{c_{1,k} 
\tau^{-1/3} \! + \! c_{2,k} \tilde{r}_{0}(\tau)}{\widetilde{\Lambda}^{2}} 
\right) \! \left(c_{3,k} \tau^{-1/3} \! + \! c_{4,k}(\tilde{r}_{0}(\tau) 
\! + \! 4v_{0}(\tau)) \right) \right) \nonumber \\
+& \, \mathcal{O} \! \left(\dfrac{\tau^{-1/3}}{\widetilde{\Lambda}} \! 
\left(c_{5,k} \tilde{r}_{0}(\tau) \tau^{-1/3} \! + \! c_{6,k} \tilde{r}_{0}
(\tau)(\tilde{r}_{0}(\tau) \! + \! 4v_{0}(\tau)) \right) \right) \nonumber \\
\underset{\tau \to +\infty}{=}& \, \mathcal{O}(\tau^{-2/3} 
\widetilde{\Lambda}^{-2}) \! + \! \mathcal{O}(\tau^{-1} 
\widetilde{\Lambda}^{-1}).
\end{align}
{}From the conditions \eqref{iden5}, equation \eqref{eqItee}, and the asymptotics \eqref{asympforf1} and \eqref{asympforf3}, 
it follows that
\begin{equation} \label{asympforf6} 
\mathcal{I}_{\tau,k}(\widetilde{\mu}) \underset{\tau \to +\infty}{=} \mathfrak{p}_{k}(\tau)(\digamma_{\tau,k}(\widetilde{\mu}) 
\! + \! \digamma_{\tau,k}^{\sharp}(\tau)) \! + \! \mathcal{O}((\tilde{r}_{0}(\tau) \! + \! 4v_{0}(\tau)) \tau^{-2/3} \widetilde{\Lambda}) 
\! + \! \mathcal{O}(\tau^{-1} \widetilde{\Lambda}).
\end{equation}
Therefore, via the asymptotic estimates \eqref{asympforf5} and \eqref{asympforf6}, and the fact that 
$\widetilde{\Lambda} \genfrac{}{}{0pt}{3}{=}{\tau \to +\infty} \mathcal{O}(\tau^{\delta_{k}})$, $0 \! < \! 
\delta \! < \! \delta_{k} \! < \! 1/9$, the result stated in the corollary (cf. equation \eqref{eq3.36}) is a 
consequence of Proposition \ref{prop3.1.4} (cf. equation \eqref{eqintTminus1T}), upon retaining only 
leading-order contributions. \hfill $\qed$
\begin{ffff} \label{cor3.1.5}
Let the conditions stated in Corollary {\rm \ref{cor3.1.3}} be valid$;$ then, for the branch of $l_{k}(\xi)$, $k \! = \! \pm 1$, that is 
positive for large and small positive $\xi$,
\begin{align} \label{eq3.41} 
\int_{\widetilde{\mu}_{0,k}}^{\widetilde{\mu}} \diag \! \left((T(\xi))^{-1} \partial_{\xi}T(\xi) \right) \md \xi 
\underset{\underset{\widetilde{\mu} \to \infty}{\tau \to +\infty}}{=}& \, \left(\mathfrak{p}_{k}(\tau) 
\digamma_{\tau,k}^{\sharp,\infty}(\tau) \! + \! \mathcal{O} \! \left(\left(\dfrac{\mathfrak{c}_{1,k} \tau^{-1/3} 
\! + \! \mathfrak{c}_{2,k} \tilde{r}_{0}(\tau)}{\widetilde{\Lambda}^{2}} \right) \right. \right. \nonumber \\
\times&\left. \left. (\mathfrak{c}_{3,k} \tau^{-1/3} \! + \! \mathfrak{c}_{4,k}(\tilde{r}_{0}(\tau) \! + \! 4v_{0}(\tau))) 
\right) \right. \nonumber \\
+&\left. \, \mathcal{O}(\widetilde{\mu}^{-2} \tau^{-1/3}(\mathfrak{c}_{5,k} \tau^{-1/3} \! + \! \mathfrak{c}_{6,k}
(\tilde{r}_{0}(\tau) \! + \! 4v_{0}(\tau)))) \right) \! \sigma_{3},
\end{align}
where $\mathfrak{p}_{k}(\tau)$ is defined by equation \eqref{eqpeetee},
\begin{equation} \label{asympforf7} 
\digamma_{\tau,k}^{\sharp,\infty}(\tau) \! := \! -\dfrac{(\sqrt{3} \! - \! 1) \tau^{1/3}}{\sqrt{3} \alpha_{k} \widetilde{\Lambda}} 
\! - \! \dfrac{2}{3 \sqrt{3} \alpha_{k}^{2}} \! \left(-\dfrac{1}{3} \ln \tau \! + \! \ln \widetilde{\Lambda} \right) \! + \! \dfrac{5 \! 
- \! \sqrt{3}}{6 \sqrt{3} \alpha_{k}^{2}} \! + \! \dfrac{2}{3 \sqrt{3} \alpha_{k}^{2}} \ln (3(2 \! - \! \sqrt{3}) \alpha_{k}),
\end{equation}
and
\begin{align} \label{eq3.43} 
\int_{\widetilde{\mu}_{0,k}}^{\widetilde{\mu}} \diag ((T(\xi))^{-1} \partial_{\xi}T(\xi)) \, \md \xi 
\underset{\underset{\widetilde{\mu} \to 0}{\tau \to +\infty}}{=}& \, \left(\mathfrak{p}_{k}(\tau) \digamma_{\tau,k}^{\sharp,0}
(\tau) \! + \! \mathcal{O} \! \left(\left(\dfrac{\mathfrak{c}_{7,k} \tau^{-1/3} \! + \! \mathfrak{c}_{8,k} 
\tilde{r}_{0}(\tau)}{\widetilde{\Lambda}^{2}} \right) \right. \right. \nonumber \\
\times&\left. \left. (\mathfrak{c}_{9k} \tau^{-1/3} \! + \! \mathfrak{c}_{10,k}(\tilde{r}_{0}(\tau) \! + \! 4v_{0}(\tau))) \right) \right. 
\nonumber \\
+&\left. \, \mathcal{O}(\widetilde{\mu}^{2} \tau^{-1/3}(\mathfrak{c}_{11,k} \tau^{-1/3} \! + \! \mathfrak{c}_{12,k}(\tilde{r}_{0}
(\tau) \! + \! 4v_{0}(\tau)))) \right) \! \sigma_{3},
\end{align}
where
\begin{equation} \label{asympforf8} 
\digamma_{\tau,k}^{\sharp,0}(\tau) \! := \! -\dfrac{(\sqrt{3} \! + \! 1) \tau^{1/3}}{\sqrt{3} \alpha_{k} \widetilde{\Lambda}} 
\! + \! \dfrac{2}{3 \sqrt{3} \alpha_{k}^{2}} \! \left(-\dfrac{1}{3} \ln \tau \! + \! \ln \widetilde{\Lambda} \right) \! - \! \dfrac{(5 
\! + \! 9 \sqrt{3})}{6 \sqrt{3} \alpha_{k}^{2}} \! + \! \dfrac{2}{3 \sqrt{3} \alpha_{k}^{2}} \ln (\me^{\mi k \pi}/3 \alpha_{k}),
\end{equation}
and $\mathfrak{c}_{m,k}$, $m \! = \! 1,2,\ldots,12$, are $\mathcal{O}(1)$.
\end{ffff}

\emph{Proof}. Choosing consistently the corresponding branches in equations \eqref{eqFtee} and \eqref{eqfteesharp}, and 
via the third and fourth lines of equation \eqref{eqee2}, respectively, one shows, via the conditions \eqref{iden5} and the 
asymptotics \eqref{asympforf3}, that (cf. equation \eqref{eq3.36})
\begin{align}
\digamma_{\tau,k}(\widetilde{\mu}) \underset{\underset{\widetilde{\mu} 
\to \infty}{\tau \to +\infty}}{=}& \, -\dfrac{2}{3 \alpha_{k}^{2}} \! + \! 
\dfrac{2}{3 \sqrt{3} \alpha_{k}^{2}} \ln (2 \! - \! \sqrt{3}) \! + \! 
\mathcal{O}(\widetilde{\mu}^{-2}), \label{iden48} \\
\digamma_{\tau,k}(\widetilde{\mu}) \underset{\underset{\widetilde{\mu} 
\to 0}{\tau \to +\infty}}{=}& \, -\dfrac{2}{\alpha_{k}^{2}} \! + \! \dfrac{2}{3 
\sqrt{3} \alpha_{k}^{2}} \ln (\me^{\mi k \pi}) \! + \! \mathcal{O}
(\widetilde{\mu}^{2}), \label{iden49} \\
\mathcal{O}(\mathcal{E}_{\scriptscriptstyle T,k}(\widetilde{\mu})) \underset{
\underset{\widetilde{\mu} \to \infty}{\tau \to +\infty}}{=}& \, \mathcal{O}
(\widetilde{\mu}^{-4} \tilde{r}_{0}(\tau)(\tilde{r}_{0}(\tau) \! + \! 4v_{0}
(\tau)) \tau^{-2/3}) \! + \! \mathcal{O} \! \left(\widetilde{\mu}^{-4} 
\tilde{r}_{0}(\tau) \tau^{-1} \right), \label{iden50} \\
\mathcal{O}(\mathcal{E}_{\scriptscriptstyle T,k}(\widetilde{\mu})) \underset{
\underset{\widetilde{\mu} \to 0}{\tau \to +\infty}}{=}& \, \mathcal{O}
(\widetilde{\mu}^{2} \tilde{r}_{0}(\tau)(\tilde{r}_{0}(\tau) \! + \! 4v_{0}(\tau)) 
\tau^{-2/3}) \! + \! \mathcal{O} \! \left(\widetilde{\mu}^{2} \tilde{r}_{0}
(\tau) \tau^{-1} \right). \label{iden51}
\end{align}
Via the conditions \eqref{iden5}, equation \eqref{eqfteesharp}, and the asymptotics \eqref{asympforf3} and 
\eqref{iden48}--\eqref{iden51}, it follows that (cf. equation \eqref{eq3.36})
\begin{gather}
\mathfrak{p}_{k}(\tau)(\digamma_{\tau,k}(\widetilde{\mu}) \! + \! 
\digamma_{\tau,k}^{\sharp}(\tau)) \underset{\underset{\widetilde{\mu} 
\to \infty}{\tau \to +\infty}}{=} \mathfrak{p}_{k}(\tau) \digamma_{
\tau,k}^{\sharp,\infty}(\tau) \! + \! \mathcal{O}(\widetilde{\mu}^{-2}
(\tilde{r}_{0}(\tau) \! + \! 4v_{0}(\tau)) \tau^{-1/3}) \! + \! \mathcal{O}
(\widetilde{\mu}^{-2} \tau^{-2/3}), \label{asympforf9} \\
\mathfrak{p}_{k}(\tau)(\digamma_{\tau,k}(\widetilde{\mu}) \! + \! 
\digamma_{\tau,k}^{\sharp}(\tau)) \underset{\underset{\widetilde{\mu} 
\to 0}{\tau \to +\infty}}{=} \mathfrak{p}_{k}(\tau) \digamma_{\tau,k}^{
\sharp,0}(\tau) \! + \! \mathcal{O}(\widetilde{\mu}^{2}(\tilde{r}_{0}(\tau) 
\! + \! 4v_{0}(\tau)) \tau^{-1/3}) \! + \! \mathcal{O}(\widetilde{\mu}^{2} 
\tau^{-2/3}), \label{asympforf10}
\end{gather}
where $\digamma_{\tau,k}^{\sharp,\infty}(\tau)$ and $\digamma_{\tau,k}^{\sharp,0}(\tau)$ are defined by 
equations \eqref{asympforf7} and \eqref{asympforf8}, respectively. The results stated in the corollary are 
now a consequence of the conditions \eqref{iden5}, equation \eqref{eq3.36}, and the asymptotic expansions 
\eqref{iden50}--\eqref{asympforf10}, upon retaining only leading-order terms. \hfill $\qed$
\begin{bbbb} \label{prop3.1.5} 
Let $T(\widetilde{\mu})$ be given in equation~\eqref{eq3.18}, with $\mathcal{A}(\widetilde{\mu})$ defined by 
equation \eqref{eq3.4} and $l_{k}^{2}(\widetilde{\mu})$, $k \! = \! \pm 1$, given in equation \eqref{eq3.19}, with 
the branches defined as in Corollary {\rm \ref{cor3.1.4}}$;$ then,
\begin{align} \label{eq3.51} 
T(\widetilde{\mu}) \underset{\underset{\widetilde{\mu} \to \infty}{\tau \to +\infty}}{=}& \, (b(\tau))^{-\frac{1}{2} 
\ad (\sigma_{3})} \! \left(\mathrm{I} \! + \! \dfrac{1}{\widetilde{\mu}} 
\begin{pmatrix}
0 & -\frac{(\varepsilon b)^{2/3}}{2}(1 \! + \! \hat{u}_{0}(\tau)) \\
\frac{2(a-\mi/2) \tau^{-2/3} - (\varepsilon b)^{1/3}(-2+\hat{r}_{0}
(\tau))}{4(\varepsilon b)^{2/3}(1+ \hat{u}_{0}(\tau))} & 0
\end{pmatrix} \right. \nonumber \\
+&\left. \, \mathcal{O} \! \left(\dfrac{1}{\widetilde{\mu}^{2}} \! 
\begin{pmatrix}
\mathfrak{c}_{1}(\tau) & 0 \\
0 & \mathfrak{c}_{1}(\tau)
\end{pmatrix} \right) \right),
\end{align}
and
\begin{align} \label{eq3.52} 
T(\widetilde{\mu}) \underset{\underset{\widetilde{\mu} \to 0}{\tau \to +\infty}}{=}& \, \frac{1}{\sqrt{2}} \! 
\left(\dfrac{b(\tau)}{\sqrt{\smash[b]{\varepsilon b}}} \right)^{-\frac{1}{2} \ad (\sigma_{3})} \! \left(
\begin{pmatrix}
1 & -1 \\
1 & 1
\end{pmatrix} \! + \! \widetilde{\mu} \, \dfrac{(-2 \! + \! \hat{r}_{0}(\tau))}{4(\varepsilon b)^{1/6}} 
\begin{pmatrix}
-1 & -1 \\
1 & -1
\end{pmatrix} \right. \nonumber \\
+&\left. \mathcal{O} \! \left(\widetilde{\mu}^{2} \! 
\begin{pmatrix}
\mathfrak{c}_{2}(\tau) & \mathfrak{c}_{3}(\tau) \\
\mathfrak{c}_{4}(\tau) & \mathfrak{c}_{2}(\tau)
\end{pmatrix}
\right) \right),
\end{align}
where $\mathfrak{c}_{1}(\tau)$, $\mathfrak{c}_{2}(\tau)$, $\mathfrak{c}_{3}(\tau)$, and $\mathfrak{c}_{4}(\tau)$, 
respectively, are defined by equations \eqref{tempsea1}--\eqref{tempsea4}.
\end{bbbb}

\emph{Proof}. The proof is presented for the asymptotics \eqref{eq3.51}. Let the conditions stated in the proposition 
be valid. Then, via equations \eqref{iden7}, \eqref{iden3}, and \eqref{iden4}, and the conditions \eqref{iden5}, one 
shows that
\begin{align}
l_{k}(\widetilde{\mu}) \underset{\underset{\widetilde{\mu} \to \infty}{\tau \to 
+\infty}}{=}& \, 2 \widetilde{\mu} \! + \! \dfrac{1}{\widetilde{\mu}}(a \! - \! \mi/2) 
\tau^{-2/3} \! + \! \mathcal{O}(\widetilde{\mu}^{-3} \hat{\lambda}_{1}(\tau)), 
\label{iden52} \\
\mi (\mathcal{A}(\widetilde{\mu}) \! - \! \mi l_{k}(\widetilde{\mu}) \sigma_{3}) 
\sigma_{3} \underset{\underset{\widetilde{\mu} \to \infty}{\tau \to +\infty}}{=}& 
\, 4 \widetilde{\mu} \, \mathrm{I} \! + \! 
\begin{pmatrix}
0 & -\frac{4 \sqrt{\smash[b]{-a(\tau)b(\tau)}}}{b(\tau)} \\
-\mi 2d(\tau) & 0
\end{pmatrix} \! + \! \dfrac{1}{\widetilde{\mu}} \mathfrak{d}_{0,0}^{\lozenge}
(\tau) \mathrm{I} \nonumber \\
+& \, \dfrac{1}{\widetilde{\mu}^{2}} \! 
\begin{pmatrix}
0 & \frac{(\varepsilon b)}{b(\tau)} \\
-b(\tau) & 0
\end{pmatrix} \! + \! \mathcal{O} \! \left(\widetilde{\mu}^{-3} 
\hat{\lambda}_{1}(\tau) \!  
\begin{pmatrix}
c_{1,k} & 0 \\
0 & c_{2,k}
\end{pmatrix} \right), \label{iden53} \\
\dfrac{1}{\sqrt{\smash[b]{2 \mi l_{k}(\widetilde{\mu})(\mathcal{A}_{11}
(\widetilde{\mu}) \! - \! \mi l_{k}(\widetilde{\mu}))}}} 
\underset{\underset{\widetilde{\mu} \to \infty}{\tau \to +\infty}}{=}& \, 
\dfrac{1}{4 \widetilde{\mu}} \! \left(1 \! - \! \dfrac{1}{\widetilde{\mu}^{2}} 
\dfrac{\mathfrak{d}_{1,0}^{\lozenge}(\tau)}{8} \! + \! \mathcal{O}
(\widetilde{\mu}^{-4} \hat{\lambda}_{2}(\tau)) \right), \label{iden54}
\end{align}
where
\begin{gather*}
\mathfrak{d}_{m,j}^{\lozenge}(\tau) \! := \! \dfrac{(\varepsilon b)^{1/3}}{2}
(-2 \! + \! \hat{r}_{0}(\tau)) \! + \! (-1)^{j}(2m \! + \! 1)(a \! - \! \mi/2) 
\tau^{-2/3}, \quad m,j \! \in \! \lbrace 0,1 \rbrace, \\
\hat{\lambda}_{1}(\tau) \! := \! -3 \alpha_{k}^{4} \! + \! \hat{h}_{0}
(\tau) \! - \! \dfrac{1}{4} \left(a \! - \! \mi/2 \right)^{2} \tau^{-4/3}, \\
\hat{\lambda}_{2}(\tau) \! := \! c_{3,k} \hat{\lambda}_{1}(\tau) \! + \! 
c_{4,k}(\mathfrak{d}_{1,0}^{\lozenge}(\tau))^{2} \! + \! c_{5,k} 
\tau^{-2/3} \mathfrak{d}_{0,0}^{\lozenge}(\tau),
\end{gather*}
and $c_{m,k}$, $m \! = \! 1,2,\dotsc,5$, are $\mathcal{O}(1)$; thus, via the conditions \eqref{iden5}, equation \eqref{eq3.18}, 
and the expansions \eqref{iden52}--\eqref{iden54}, one arrives at the asymptotics \eqref{eq3.51}, where
\begin{equation} \label{tempsea1} 
\mathfrak{c}_{1}(\tau) \! := \! \mathfrak{d}_{0,1}^{\lozenge}(\tau)/8.
\end{equation}
Proceeding analogously, one arrives at the asymptotics \eqref{eq3.52}, where
\begin{gather}
\mathfrak{c}_{2}(\tau) \! := \! -\dfrac{(-2 \! + \! \hat{r}_{0}(\tau))^{2}}{32
(\varepsilon b)^{1/3}}, \label{tempsea2} \\
\mathfrak{c}_{3}(\tau) \! := \! \dfrac{-3 \alpha_{k}^{4} \! + \! \hat{h}_{0}(\tau)}{4 
\alpha_{k}^{6}} \! - \! \dfrac{3(-2 \! + \! \hat{r}_{0}(\tau))^{2}}{32(\varepsilon 
b)^{1/3}} \! + \! \dfrac{2(1 \! + \! \hat{u}_{0}(\tau))}{(\varepsilon b)^{1/3}}, 
\label{tempsea3} \\
\mathfrak{c}_{4}(\tau) \! := \! \dfrac{3 \alpha_{k}^{4} \! - \! \hat{h}_{0}(\tau)}{4 
\alpha_{k}^{6}} \! + \! \dfrac{3(-2 \! + \! \hat{r}_{0}(\tau))^{2}}{32(\varepsilon 
b)^{1/3}} \! + \! \dfrac{2 \mathfrak{d}_{0,1}^{\lozenge}(\tau)}{(\varepsilon 
b)^{2/3}(1 \! + \! \hat{u}_{0}(\tau))}, \label{tempsea4}
\end{gather}
with $\mathfrak{d}_{0,1}^{\lozenge}(\tau)$ defined above. \hfill $\qed$
\begin{bbbb} \label{prop3.1.6} 
Let $T(\widetilde{\mu})$ be given in equation \eqref{eq3.18}, with $\mathcal{A}(\widetilde{\mu})$ defined by 
equation \eqref{eq3.4} and $l_{k}^{2}(\widetilde{\mu})$, $k \! = \! \pm 1$, given in equation \eqref{eq3.19}. Set 
$\widetilde{\mu}_{0,k} \! = \! \alpha_{k} \! + \! \tau^{-1/3} \widetilde{\Lambda}$, where $\widetilde{\Lambda} 
\genfrac{}{}{0pt}{3}{=}{\tau \to +\infty} \mathcal{O}(\tau^{\delta_{k}})$, $0 \! < \! \delta \! < \! \delta_{k} \! < \! 1/9;$ then,
{\fontsize{9pt}{10pt}\selectfont
\begin{align} \label{iden55} 
T(\widetilde{\mu}_{0,k}) \underset{\tau \to +\infty}{=}& \, \dfrac{(b(\tau))^{
-\frac{1}{2} \ad (\sigma_{3})}}{(2 \sqrt{3}(\varpi \! + \! \sqrt{3}))^{1/2}} \left(
\begin{pmatrix}
\varpi \! + \! \sqrt{3} & (2 \varepsilon b)^{1/2} \varpi \\
-\frac{\sqrt{2} \varpi}{(\varepsilon b)^{1/2}} & \varpi \! + \! \sqrt{3}
\end{pmatrix} \! + \! 
\begin{pmatrix}
\frac{\varpi}{3 \alpha_{k}} & -\frac{(2 \varepsilon b)^{1/2}
(2 \varpi + \sqrt{3}) \varpi}{3(\varpi + \sqrt{3}) \alpha_{k}} \\
\frac{\sqrt{2}(2 \varpi + \sqrt{3}) \varpi}{3(\varepsilon b)^{1/2}
(\varpi + \sqrt{3}) \alpha_{k}} & \frac{\varpi}{3 \alpha_{k}} 
\end{pmatrix} \! \tau^{-1/3} \widetilde{\Lambda} \right. 
\nonumber \\
+&\left. \, 
\begin{pmatrix}
\mathbb{T}_{11,k}(\varpi;\tau) & \mathbb{T}_{12,k}(\varpi;\tau) \\
\mathbb{T}_{21,k}(\varpi;\tau) & \mathbb{T}_{22,k}(\varpi;\tau)
\end{pmatrix} \! \dfrac{1}{\widetilde{\Lambda}} \! + \! 
\mathcal{O} \! \left(
\begin{pmatrix}
\mathfrak{c}_{1,k} & \mathfrak{c}_{2,k} \\
\mathfrak{c}_{3,k} & \mathfrak{c}_{1,k} 
\end{pmatrix} \! (\tau^{-1/3} \widetilde{\Lambda})^{2} \right) \right),
\end{align}}
where
\begin{align}
\mathbb{T}_{11,k}(\varpi;\tau) \! =& \, \mathbb{T}_{22,k}(\varpi;\tau) \! 
:= \! \dfrac{\varpi}{4} \! \left(\dfrac{\alpha_{k} \tilde{r}_{0}(\tau)}{2} \! - \! 
\frac{\tau^{-1/3} \hat{\mathfrak{g}}^{\ast}_{k}(\tau)}{3 \alpha_{k}} \right), 
\label{unikay} \\
\mathbb{T}_{12,k}(\varpi;\tau) \! :=& \left(\frac{\varepsilon b}{2} \right)^{1/2} 
\! \left(\varpi \alpha_{k}v_{0}(\tau) \! - \! \dfrac{\alpha_{k} \tilde{r}_{0}(\tau)}{4
(\varpi \! + \! \sqrt{3})} \! - \! \dfrac{(1 \! + \! 2 \sqrt{3} \varpi) \tau^{-1/3} 
\hat{\mathfrak{g}}^{\ast}_{k}(\tau)}{6(\varpi \! + \! \sqrt{3}) \alpha_{k}} \right), 
\label{unikby} \\
\mathbb{T}_{21,k}(\varpi;\tau) \! :=& \, \dfrac{\varpi}{(2 \varepsilon b)^{1/2}} 
\! \left(\dfrac{(\varepsilon b)^{1/3}(\tilde{r}_{0}(\tau) \! + \! 2v_{0}(\tau)) \! + \! 
2(a \! - \! \mi/2) \me^{\mi \pi k/3} \tau^{-1/3}}{2^{3/2}(\varepsilon b)^{1/6} 
\me^{-\mi \pi k/3}(1 \! + \! v_{0}(\tau) \tau^{-1/3})} \right. \nonumber \\
+&\left. \, 
\dfrac{\alpha_{k} 
\tilde{r}_{0}(\tau) \! + \! \frac{2(1+2 \sqrt{3} \varpi) \tau^{-1/3} 
\hat{\mathfrak{g}}^{\ast}_{k}(\tau)}{3 \alpha_{k}}}{4(\varpi \! + \! \sqrt{3}) 
\varpi} \right), \label{unikcy}
\end{align}
with $\hat{\mathfrak{g}}^{\ast}_{k}(\tau) \! := \! \tau^{2/3} \hat{\mathfrak{f}}_{1,k}(\tau)$, where $\hat{\mathfrak{f}}_{1,k}
(\tau)$ is given in equation \eqref{mfkeff} (for $j \! = \! 1)$, $\big(\widetilde{\Lambda}^{2} \big)^{1/2} \! := \! \varpi 
\widetilde{\Lambda}$, $\varpi \! = \! \pm 1$, and $\mathfrak{c}_{m,k}$, $m \! = \! 1,2,3$, are $\mathcal{O}(1)$.
\end{bbbb}

\emph{Proof}. Set $T(\widetilde{\mu}) \! = \! (T(\widetilde{\mu}))_{i,j=1,2}$. {}From the formula for $T(\widetilde{\mu})$ 
given in equation \eqref{eq3.18}, with $\mathcal{A}(\widetilde{\mu})$ defined by equation \eqref{eq3.4} and 
$l_{k}^{2}(\widetilde{\mu})$, $k \! = \! \pm 1$, given in equation \eqref{eq3.19}, one shows that
\begin{equation} \label{iden56} 
\begin{gathered}
T_{11}(\widetilde{\mu}) \! = \! T_{22}(\widetilde{\mu}) \! = \! \dfrac{\mi (\mathcal{A}_{11}(\widetilde{\mu}) \! - \! \mi l_{k}
(\widetilde{\mu}))}{\sqrt{\smash[b]{2 \mi l_{k}(\widetilde{\mu})(\mathcal{A}_{11}(\widetilde{\mu}) \! - \! \mi l_{k}
(\widetilde{\mu}))}}}, \, \qquad \, T_{12}(\widetilde{\mu}) \! = \! -\dfrac{\mi \mathcal{A}_{12}(\widetilde{\mu})}{
\sqrt{\smash[b]{2 \mi l_{k}(\widetilde{\mu})(\mathcal{A}_{11}(\widetilde{\mu}) \! - \! \mi l_{k}(\widetilde{\mu}))}}}, \\
T_{21}(\widetilde{\mu}) \! = \! \dfrac{\mi \mathcal{A}_{21}(\widetilde{\mu})}{\sqrt{\smash[b]{2 \mi l_{k}(\widetilde{\mu})
(\mathcal{A}_{11}(\widetilde{\mu}) \! - \! \mi l_{k}(\widetilde{\mu}))}}}.
\end{gathered}
\end{equation}
{}From equations \eqref{eq3.4}, \eqref{iden7}, \eqref{iden3}, and \eqref{iden4}, the conditions \eqref{iden5}, and 
equation \eqref{mfkeff} for $\hat{\mathfrak{f}}_{1,k}(\tau)$ (with associated asymptotics \eqref{asympforf4}), one 
shows, upon taking $\widetilde{\mu}_{0,k}$ as stated in the proposition, that
\begin{align}
\dfrac{1}{\sqrt{\smash[b]{2 \mi l_{k}(\widetilde{\mu}_{0,k})(\mathcal{A}_{11}
(\widetilde{\mu}_{0,k}) \! - \! \mi l_{k}(\widetilde{\mu}_{0,k}))}}} 
\underset{\tau \to +\infty}{=}& \, \dfrac{(\varpi \tau^{-1/3} 
\widetilde{\Lambda})^{-1}}{4(2 \sqrt{3}(\varpi \! + \! \sqrt{3}))^{1/2}} 
\! \left(1 \! + \! \dfrac{(5 \varpi \! + \! 7 \sqrt{3})}{6(\varpi \! + \! 
\sqrt{3}) \alpha_{k}} \tau^{-1/3} \widetilde{\Lambda} \right. \nonumber \\
-&\left. \, \left(\dfrac{\alpha_{k} \tilde{r}_{0}(\tau) \! + \! 2(1 \! + \! 2 \sqrt{3} 
\varpi)(3 \alpha_{k})^{-1} \hat{\mathfrak{g}}^{\ast}_{k}(\tau) \tau^{-1/3}}{8 
\varpi (\varpi \! + \! \sqrt{3})} \right) \! \dfrac{1}{\widetilde{\Lambda}} \right. 
\nonumber \\
+&\left. \, \mathcal{O} \big((\tau^{-1/3} \widetilde{\Lambda})^{2} \big) \vphantom{M^{M^{M^{M^{M^{M}}}}}} \right), 
\label{iden59}
\end{align}
\begin{align}
\mi \mathcal{A}_{11}(\widetilde{\mu}_{0,k}) \! + \! l_{k}(\widetilde{\mu}_{0,k}) 
\underset{\tau \to +\infty}{=}& \, 4 \varpi (\varpi \! + \! \sqrt{3}) \tau^{-1/3} 
\widetilde{\Lambda} \! \left(1 \! - \! \dfrac{\sqrt{3}(7 \! + \! \sqrt{3} \varpi)}{
6(\varpi \! + \! \sqrt{3}) \alpha_{k}} \tau^{-1/3} \widetilde{\Lambda} \! + \! 
\mathcal{O}((\tau^{-1/3} \widetilde{\Lambda})^{2}) \right. \nonumber \\
+&\left. \, \left(\dfrac{\alpha_{k} \tilde{r}_{0}(\tau) \! + \! 2 \varpi (\sqrt{3} 
\alpha_{k})^{-1} \hat{\mathfrak{g}}^{\ast}_{k}(\tau) \tau^{-1/3}}{4 \varpi 
(\varpi \! + \! \sqrt{3})} \right) \! \dfrac{1}{\widetilde{\Lambda}} 
\vphantom{M^{M^{M^{M^{M^{M}}}}}} \right), \label{iden60} \\
-\mi \mathcal{A}_{12}(\widetilde{\mu}_{0,k}) \underset{\tau \to +\infty}{=}& 
\, (b(\tau))^{-1} \! \left(-2(\varepsilon b) \alpha_{k}^{-3} \tau^{-1/3} 
\widetilde{\Lambda} \! + \! 3(\varepsilon b) \alpha_{k}^{-4}(\tau^{-1/3} 
\widetilde{\Lambda})^{2} \! + \! \mathcal{O}((\tau^{-1/3} \widetilde{\Lambda})^{3}) 
\right. \nonumber \\
-&\left. \, 2(\varepsilon b)^{2/3} \me^{-\mi 2 \pi k/3}v_{0}(\tau) \tau^{-1/3} 
\vphantom{M^{M^{M^{M^{M}}}}} \right), \label{iden61} \\
\mi \mathcal{A}_{21}(\widetilde{\mu}_{0,k}) \underset{\tau \to +\infty}{=}& 
\, b(\tau) \! \left(2 \alpha_{k}^{-3} \tau^{-1/3} \widetilde{\Lambda} \! - \! 3 
\alpha_{k}^{-4}(\tau^{-1/3} \widetilde{\Lambda})^{2} \! + \! \mathcal{O}
((\tau^{-1/3} \widetilde{\Lambda})^{3}) \right. \nonumber \\
+&\left. \, \dfrac{\me^{\mi \pi k/3} \tau^{-1/3} \left((\varepsilon b)^{1/3}
(\tilde{r}_{0}(\tau) \! + \! 2v_{0}(\tau)) \! + \! 2(a \! - \! \mi/2) \me^{\mi \pi k/3} 
\tau^{-1/3} \right)}{(\varepsilon b)^{2/3}(1 \! + \! v_{0}(\tau) \tau^{-1/3})} 
\right), \label{iden62}
\end{align}
where $\hat{\mathfrak{g}}^{\ast}_{k}(\tau)$ and $\varpi$ are defined in the proposition. Substituting the expansions 
\eqref{iden59}---\eqref{iden62} into equations \eqref{iden56} (with $\widetilde{\mu} \! = \! \widetilde{\mu}_{0,k})$, 
one arrives at the asymptotics for $T(\widetilde{\mu}_{0,k})$ stated in the proposition. \hfill $\qed$
\subsection{Parametrix Near the Double-Turning Points} \label{sec3.2}
The matrix WKB formula (cf. equation \eqref{eq3.16}) doesn't provide an approximation for solutions of equation 
\eqref{eq3.3} in shrinking (as $\tau \! \to \! +\infty$ with $\varepsilon b \! > \! 0$) neighbourhoods of the turning points, 
where a more refined approximation must be constructed. There are two simple turning points approaching 
$\pm \mi \sqrt{2} \alpha_{k}$, $k \! = \! \pm 1$: the approximate solution of equation \eqref{eq3.3} in the 
neighbourhoods of these turning points is representable in terms of Airy functions (see, for example, \cite{a5,a2}, 
\textbf{Riemann-Hilbert Problem~4} in \cite{thombothpdmill}, \cite{BuckinghamMiller2022}, and Subsections 3.5 
and 3.6 in \cite{shimpeefive}). There are, additionally, two pairs of double-turning points, one pair coalescing at 
$-\alpha_{k}$, and another pair coalescing at $\alpha_{k}$: in neighbourhoods of $\pm \alpha_{k}$, the approximate 
solution of equation \eqref{eq3.3} is expressed in terms of parabolic-cylinder functions (see, for example, 
\cite{F,a5,a18,a2,W}). In order to obtain asymptotics for $u(\tau)$ and the associated, auxiliary functions $f_{\pm}
(\tau)$, $\mathcal{H}(\tau)$, $\sigma (\tau)$, and $\hat{\varphi}(\tau)$, it is sufficient to study a subset of the complete 
set of the monodromy data, which can be calculated via the approximation of the general solution of equation 
\eqref{eq3.3} in a neighbourhood of the double-turning point $\alpha_{k}$, because the remaining monodromy data 
can be calculated via equations \eqref{monoeqns}, which define the monodromy manifold.\footnote{More precisely, 
equations \eqref{monok2} (resp., equations \eqref{monok3}) for $k \! = \! +1$ (resp., $k \! = \! -1$).} For the 
asymptotic conditions \eqref{iden5} on the functions $\hat{h}_{0}(\tau)$, $\tilde{r}_{0}(\tau)$, and $v_{0}(\tau)$, 
this parametrix (approximation) is given in Lemma \ref{nprcl} below.
\begin{cccc} \label{nprcl} 
Set
\begin{equation} \label{prpr1} 
\nu (k) \! + \! 1 \! := \! -\dfrac{p_{k}(\tau)q_{k}(\tau)}{2 \mu_{k}(\tau)}, \quad k \! = \! \pm 1,
\end{equation}
where $\mu_{k}(\tau)$, $p_{k}(\tau)$, and $q_{k}(\tau)$ are defined by equations \eqref{prcy54}, \eqref{prcy57}, and 
\eqref{prcy58}, respectively,\footnote{See, also, the corresponding definitions \eqref{prcy10}, \eqref{prcy15}--\eqref{prcy20}, 
\eqref{prcy22}, \eqref{prcy33}--\eqref{prcy35}, \eqref{prcy40}, \eqref{prcy45}, \eqref{prcy46}, and \eqref{prcy53}.} 
and let $\widetilde{\mu} \! = \! \widetilde{\mu}_{0,k} \! = \! \alpha_{k} \! + \! \tau^{-1/3} \widetilde{\Lambda}$, where 
$\widetilde{\Lambda} \genfrac{}{}{0pt}{3}{=}{\tau \to +\infty} \mathcal{O}(\tau^{\delta_{k}})$, $0 \! < \! \delta \! < \! \delta_{k} 
\! < \! 1/9$. Concomitant with equations \eqref{iden1}--\eqref{iden6}, the definitions \eqref{iden2}--\eqref{iden4}, and the 
conditions \eqref{iden5}, impose the following restrictions:
\begin{equation} \label{pc4} 
\begin{gathered}
0 \underset{\tau \to +\infty}{<} \Re (\nu (k) \! + \! 1) \underset{\tau \to +\infty}{<} 1, \, \, \qquad \, \, \Im (\nu (k) \! + \! 1) 
\underset{\tau \to +\infty}{\leqslant} \mathcal{O}(1), \\
0 \underset{\tau \to +\infty}{<} \delta_{k} \underset{\tau \to +\infty}{<} \dfrac{1}{6(3 \! + \! \Re (\nu (k) \! + \! 1))}, \quad 
k \! = \! \pm 1.
\end{gathered}
\end{equation}
Then, there exists a fundamental solution of equation \eqref{eq3.3}, $\widetilde{\Psi}(\widetilde{\mu}) \! = \! 
\widetilde{\Psi}_{k}(\widetilde{\mu},\tau)$, $k \! = \! \pm 1$, with asymptotics
\begin{align} \label{prpr2} 
\widetilde{\Psi}_{k}(\widetilde{\mu},\tau) \underset{\tau \to +\infty}{=}& \, (b(\tau))^{-\frac{1}{2} \sigma_{3}} 
\mathcal{G}_{0,k} \mathfrak{B}_{k}^{\frac{1}{2} \sigma_{3}} 
\begin{pmatrix}
1 & 0 \\
(\frac{4 \mi \sqrt{3} \mathcal{Z}_{k}}{\chi_{k}(\tau)} \! - \! 1) \mathfrak{A}_{k} & 1
\end{pmatrix} \! \left(\mathrm{I} \! + \! \gimel_{{\scriptscriptstyle A,k}}(\tau) \widetilde{\Lambda} \! + \! 
\gimel_{\scriptscriptstyle{B,k}}(\tau) \widetilde{\Lambda}^{2} \right) \nonumber \\
\times& \left(\mathrm{I} \! + \! \mathcal{O} \! \left(\tilde{\mathfrak{C}}_{k}(\tau) \frac{\lvert \nu (k) \! + \! 1 
\rvert^{2}}{\lvert p_{k}(\tau) \rvert^{2}} \tau^{-\big(\frac{1}{3}-2(3+ \Re (\nu (k)+1)) \delta_{k} \big)} \right) \right) \! 
\Phi_{M,k}(\widetilde{\Lambda}),
\end{align}
where
\begin{gather} 
\gimel_{{\scriptscriptstyle A,k}}(\tau) \! := \! 
\begin{pmatrix}
\frac{4 \mi \sqrt{3} \mathcal{Z}_{k} \mathfrak{A}_{k} \ell_{0,k}^{+}}{\chi_{k}(\tau)} & \ell_{0,k}^{+} \\
(\frac{4 \sqrt{3} \mathcal{Z}_{k} \mathfrak{A}_{k}}{\chi_{k}(\tau)})^{2} \ell_{0,k}^{+} \! + \! \ell_{1,k}^{+} 
\! + \! \ell_{2,k}^{+} & -\frac{4 \mi \sqrt{3} \mathcal{Z}_{k} \mathfrak{A}_{k} \ell_{0,k}^{+}}{\chi_{k}(\tau)}
\end{pmatrix}, \label{prpr3} \\
\gimel_{{\scriptscriptstyle B,k}}(\tau) \! := \! \ell_{0,k}^{+}(\ell_{1,k}^{+} \! + \! \ell_{2,k}^{+}) \! 
\begin{pmatrix}
1 & 0 \\
-\frac{4 \mi \sqrt{3} \mathcal{Z}_{k} \mathfrak{A}_{k}}{\chi_{k}(\tau)} & 0
\end{pmatrix}, \label{prpr4}
\end{gather}
with $\mathcal{G}_{0,k}$, $\mathcal{Z}_{k}$, $\mathfrak{A}_{k}$, $\mathfrak{B}_{k}$, $\ell_{0,k}^{+}$, $\ell_{1,k}^{+}$, 
$\chi_{k}(\tau)$, and $\ell_{2,k}^{+}$ defined by equations \eqref{prcy9}, \eqref{prcy10}, \eqref{prcy15}, \eqref{prcy16}, 
\eqref{prcy40}, \eqref{prcy45}, \eqref{prcy46}, and \eqref{prcy53}, respectively,\footnote{See, also, the corresponding 
definition \eqref{prcy5}.} $\mathrm{M}_{2}(\mathbb{C}) \! \ni \! \tilde{\mathfrak{C}}_{k}(\tau) 
\genfrac{}{}{0pt}{3}{=}{\tau \to +\infty} \mathcal{O}(1)$, and $\Phi_{M,k}(\widetilde{\Lambda})$ is a fundamental solution of
\begin{equation} \label{prcy1} 
\dfrac{\partial \Phi_{M,k}(\widetilde{\Lambda})}{\partial \widetilde{\Lambda}} \! = \! \left(\mu_{k}(\tau) \widetilde{\Lambda} 
\sigma_{3} \! + \! p_{k}(\tau) \sigma_{+} \! + \! q_{k}(\tau) \sigma_{-} \right) \! \Phi_{M,k}(\widetilde{\Lambda}):
\end{equation}
$\Phi_{M,k}(\widetilde{\Lambda})$ has the explicit representation
\begin{equation} \label{prcy2} 
\Phi_{M,k}(\widetilde{\Lambda}) \! = \! 
\begin{pmatrix}
D_{-\nu (k)-1}(\mi (2 \mu_{k}(\tau))^{1/2} \widetilde{\Lambda}) & D_{\nu (k)}((2 \mu_{k}(\tau))^{1/2} \widetilde{\Lambda}) \\
\mathbb{D}_{k}^{\ast}(\tau,\widetilde{\Lambda})D_{-\nu (k)-1}(\mi (2 \mu_{k}(\tau))^{1/2} \widetilde{\Lambda}) & 
\mathbb{D}_{k}^{\ast}(\tau,\widetilde{\Lambda})D_{\nu (k)}((2 \mu_{k}(\tau))^{1/2} \widetilde{\Lambda})
\end{pmatrix},
\end{equation}
where $\mathbb{D}_{k}^{\ast}(\tau,\widetilde{\Lambda}) \! := \! \frac{1}{p_{k}(\tau)} \! \left(\frac{\partial}{\partial 
\widetilde{\Lambda}} \! - \! \mu_{k}(\tau) \widetilde{\Lambda} \right)$, and $D_{\pmb{\ast}}(\boldsymbol{\cdot})$ is the 
parabolic-cylinder function {\rm \cite{a24}}.
\end{cccc}

\emph{Proof}. The derivation of the parametrix \eqref{prpr2} for a fundamental solution of equation \eqref{eq3.3} consists 
of applying the sequence of invertible linear transformations $\mathfrak{F}_{j}$, $j \! = \! 1,2,\dotsc,11$; for $k \! = \! \pm 1$,
\begin{align*}
\text{\pmb{(i)}} \quad \mathfrak{F}_{1} \colon& \operatorname{SL}_{2}
(\mathbb{C}) \! \ni \! \widetilde{\Psi}(\widetilde{\mu}) \! \mapsto \! 
\widetilde{\Psi}_{k}(\widetilde{\Lambda}) \! := \! \widetilde{\Psi}
(\alpha_{k} \! + \! \tau^{-1/3} \widetilde{\Lambda}), \\
\text{\pmb{(ii)}} \quad \mathfrak{F}_{2} \colon& \operatorname{SL}_{2}
(\mathbb{C}) \! \ni \! \widetilde{\Psi}_{k}(\widetilde{\Lambda}) \! \mapsto 
\! \tilde{\Phi}_{k}(\widetilde{\Lambda}) \! := \! (b(\tau))^{\frac{1}{2} 
\sigma_{3}} \widetilde{\Psi}_{k}(\widetilde{\Lambda}), \\
\text{\pmb{(iii)}} \quad \mathfrak{F}_{3} \colon& \operatorname{SL}_{2}
(\mathbb{C}) \! \ni \! \tilde{\Phi}_{k}(\widetilde{\Lambda}) \! \mapsto 
\! \Phi_{k}^{\sharp}(\widetilde{\Lambda}) \! := \! \mathcal{G}_{0,k}^{-1} 
\tilde{\Phi}_{k}(\widetilde{\Lambda}), \\
\text{\pmb{(iv)}} \quad \mathfrak{F}_{4} \colon& \operatorname{SL}_{2}
(\mathbb{C}) \! \ni \! \Phi_{k}^{\sharp}(\widetilde{\Lambda}) \! \mapsto 
\! \hat{\Phi}_{k}(\widetilde{\Lambda}) \! := \! \mathcal{G}_{1,k}^{-1} 
\Phi_{k}^{\sharp}(\widetilde{\Lambda}), \\
\text{\pmb{(v)}} \quad \mathfrak{F}_{5} \colon& \operatorname{SL}_{2}
(\mathbb{C}) \! \ni \! \hat{\Phi}_{k}(\widetilde{\Lambda}) \! \mapsto \! 
\hat{\Phi}_{0,k}(\widetilde{\Lambda}) \! := \! \tau^{-\frac{1}{6} \sigma_{3}} 
\hat{\Phi}_{k}(\widetilde{\Lambda}), \\
\text{\pmb{(vi)}} \quad \mathfrak{F}_{6} \colon& \operatorname{SL}_{2}
(\mathbb{C}) \! \ni \! \hat{\Phi}_{0,k}(\widetilde{\Lambda}) \! \mapsto \! 
\Phi_{0,k}(\widetilde{\Lambda}) \! := \! (\mathrm{I} \! + \! \mi \omega_{0,k} 
\sigma_{-}) \hat{\Phi}_{0,k}(\widetilde{\Lambda}), \\
\text{\pmb{(vii)}} \quad \mathfrak{F}_{7} \colon& \operatorname{SL}_{2}
(\mathbb{C}) \! \ni \! \Phi_{0,k}(\widetilde{\Lambda}) \! \mapsto \! 
\Phi_{0,k}^{\flat}(\widetilde{\Lambda}) \! := \! (\mathrm{I} \! - \! \ell_{0,k} 
\widetilde{\Lambda} \sigma_{+}) \Phi_{0,k}(\widetilde{\Lambda}), \\
\text{\pmb{(viii)}} \quad \mathfrak{F}_{8} \colon& \operatorname{SL}_{2}
(\mathbb{C}) \! \ni \! \Phi_{0,k}^{\flat}(\widetilde{\Lambda}) \! \mapsto \! 
\Phi_{0,k}^{\sharp}(\widetilde{\Lambda}) \! := \! (\mathrm{I} \! - \! \ell_{1,k} 
\widetilde{\Lambda} \sigma_{-}) \Phi_{0,k}^{\flat}(\widetilde{\Lambda}), \\
\text{\pmb{(ix)}} \quad \mathfrak{F}_{9} \colon& \operatorname{SL}_{2}
(\mathbb{C}) \! \ni \! \Phi_{0,k}^{\sharp}(\widetilde{\Lambda}) \! \mapsto 
\! \Phi_{0,k}^{\natural}(\widetilde{\Lambda}) \! := \! \mathcal{G}_{2,k}^{-1} 
\Phi_{0,k}^{\sharp}(\widetilde{\Lambda}), \\
\text{\pmb{(x)}} \quad \mathfrak{F}_{10} \colon& \operatorname{SL}_{2}
(\mathbb{C}) \! \ni \! \Phi_{0,k}^{\natural}(\widetilde{\Lambda}) \! \mapsto 
\! \Phi_{k}^{\ast}(\widetilde{\Lambda}) \! := \! (\mathrm{I} \! - \! \ell_{2,k} 
\widetilde{\Lambda} \sigma_{-}) \Phi_{0,k}^{\natural}(\widetilde{\Lambda}), \\
\text{\pmb{(xi)}} \quad \mathfrak{F}_{11} \colon& \operatorname{SL}_{2}
(\mathbb{C}) \! \ni \! \Phi_{k}^{\ast}(\widetilde{\Lambda}) \! \mapsto \! 
\Phi_{M,k}(\widetilde{\Lambda}) \! := \! \hat{\chi}_{k}^{-1}(\widetilde{\Lambda}) 
\Phi_{k}^{\ast}(\widetilde{\Lambda}) \! \in \! \mathrm{M}_{2}(\mathbb{C}),
\end{align*}
where the $\mathrm{M}_{2}(\mathbb{C})$-valued, $\tau$-dependent functions $\mathcal{G}_{0,k}$, 
$\mathcal{G}_{1,k}$, $\mathrm{I} \! + \! \mi \omega_{0,k} \sigma_{-}$, $\mathcal{G}_{2,k}$, and 
$\hat{\chi}_{k}(\widetilde{\Lambda})$, and the $\tau$-dependent parameters $\ell_{0,k}$, $\ell_{1,k}$, 
and $\ell_{2,k}$ are described in steps \pmb{(iii)}, \pmb{(iv)}, \pmb{(vi)}, \pmb{(ix)}, \pmb{(xi)}, \pmb{(vii)}, 
\pmb{(viii)}, and \pmb{(x)}, respectively, below, and $\mathrm{M}_{2}(\mathbb{C}) \! \ni \! \Phi_{M,k}
(\widetilde{\Lambda})$ is given in equation~\eqref{prcy2}.

\pmb{(i)} The gist of this step is to simplify the system \eqref{eq3.3} in a proper neighbourhood of the 
(coalescing) double-turning point $\alpha_{k}$, $k \! \in \! \lbrace \pm 1 \rbrace$. Let $\widetilde{\Psi}
(\widetilde{\mu})$ solve equation \eqref{eq3.3}; then, using equations \eqref{iden3oldr}, \eqref{iden4oldu}, 
\eqref{iden7}, \eqref{iden3}, and \eqref{iden4}, the conditions \eqref{iden5}, and applying the transformation 
$\mathfrak{F}_{1}$, one shows that, for $k \! = \! \pm 1$,
\begin{equation} \label{prcy3} 
\dfrac{\partial \widetilde{\Psi}_{k}(\widetilde{\Lambda})}{\partial \widetilde{\Lambda}} \underset{\tau \to +\infty}{=} 
(b(\tau))^{-\frac{1}{2} \ad (\sigma_{3})} \! \left(\hat{\mathcal{P}}_{0,k}(\tau) \! + \! \hat{\mathcal{P}}_{1,k}(\tau) 
\widetilde{\Lambda} \! + \! \hat{\mathcal{P}}_{2,k}(\tau) \widetilde{\Lambda}^{2} \! + \! \mathcal{O} \big(
\hat{\mathbb{E}}_{k}(\tau) \widetilde{\Lambda}^{3} \big) \right) \! \widetilde{\Psi}_{k}(\widetilde{\Lambda}),
\end{equation}
where
\begin{align}
\hat{\mathcal{P}}_{0,k}(\tau) =& 
\begin{pmatrix}
-\mi \alpha_{k} \tilde{r}_{0}(\tau) & -2 \mi (\varepsilon b)^{2/3} \me^{-\mi 2 \pi k/3}v_{0}(\tau) \\
-\frac{(\mi (\varepsilon b)^{1/3} \me^{\mi \pi k/3}(\tilde{r}_{0}(\tau)+2v_{0}(\tau))+2 \mi (a - \mi/2) 
\me^{\mi 2 \pi k/3} \tau^{-1/3})}{(\varepsilon b)^{2/3}(1+v_{0}(\tau) \tau^{-1/3})} & \mi \alpha_{k} \tilde{r}_{0}(\tau)
\end{pmatrix} \nonumber \\
=:& 
\begin{pmatrix}
\hat{\mathcal{A}}_{0} & \hat{\mathcal{B}}_{0} \\
\hat{\mathcal{C}}_{0} & -\hat{\mathcal{A}}_{0}
\end{pmatrix}, \label{prcy4} \\
\hat{\mathcal{P}}_{1,k}(\tau) =& 
\begin{pmatrix}
\mi (-4 \! + \! \tilde{r}_{0}(\tau) \tau^{-1/3}) & 4 \mi \sqrt{2}(\varepsilon b)^{1/2} \\
4 \mi \sqrt{2}(\varepsilon b)^{-1/2} & -\mi (-4 \! + \! \tilde{r}_{0}(\tau) \tau^{-1/3})
\end{pmatrix} \! =: \! 
\begin{pmatrix}
\hat{\mathcal{A}}_{1} & \hat{\mathcal{B}}_{1} \\
\hat{\mathcal{C}}_{1} & -\hat{\mathcal{A}}_{1}
\end{pmatrix}, \label{prcy5} \\
\hat{\mathcal{P}}_{2,k}(\tau) =& 
\begin{pmatrix}
\frac{\mi \sqrt{2} \me^{\mi 2 \pi k/3}}{(\varepsilon b)^{1/6}}(-2 \! + \! \tilde{r}_{0}(\tau) \tau^{-1/3}) \tau^{-1/3} & 
-12 \mi (\varepsilon b)^{1/3} \me^{-\mi \pi k/3} \tau^{-1/3} \\
-12 \mi (\varepsilon b)^{-2/3} \me^{-\mi \pi k/3} \tau^{-1/3} & -\frac{\mi \sqrt{2} \me^{\mi 2 \pi k/3}}{(\varepsilon b)^{1/6}}
(-2 \! + \! \tilde{r}_{0}(\tau) \tau^{-1/3}) \tau^{-1/3}
\end{pmatrix} \nonumber \\
=:& 
\begin{pmatrix}
\hat{\mathcal{A}}_{2} & \hat{\mathcal{B}}_{2} \\
\hat{\mathcal{C}}_{2} & -\hat{\mathcal{A}}_{2}
\end{pmatrix}, \label{prcy6}
\end{align}
and
\begin{equation} \label{prcy7} 
\hat{\mathbb{E}}_{k}(\tau) \! = \! 
\begin{pmatrix}
\mi \alpha_{k}^{-2}(-2 \! + \! \tilde{r}_{0}(\tau) \tau^{-1/3}) \tau^{-2/3} & -32 \mi \alpha_{k} \tau^{-2/3} \\
-4 \mi \alpha_{k}^{-5} \tau^{-2/3} & -\mi \alpha_{k}^{-2}(-2 \! + \! \tilde{r}_{0}(\tau) \tau^{-1/3}) \tau^{-2/3}
\end{pmatrix}.
\end{equation}
Observe that $\tr (\hat{\mathcal{P}}_{0,k}(\tau)) \! = \! \tr (\hat{\mathcal{P}}_{1,k}(\tau)) \! = \! \tr (\hat{\mathcal{P}}_{2,k}
(\tau)) \! = \! \tr (\hat{\mathbb{E}}_{k}(\tau)) \! = \! 0$.

\pmb{(ii)} This intermediate step removes the scalar-valued function $b(\tau)$ {}from equation \eqref{prcy3}. Let 
$\widetilde{\Psi}_{k}(\widetilde{\Lambda})$ solve equation \eqref{prcy3}; then, applying the transformation 
$\mathfrak{F}_{2}$, one shows that, for $k \! = \! \pm 1$,
\begin{equation} \label{prcy8} 
\dfrac{\partial \tilde{\Phi}_{k}(\widetilde{\Lambda})}{\partial \widetilde{\Lambda}} \underset{\tau \to +\infty}{=} 
\left(\hat{\mathcal{P}}_{0,k}(\tau) \! + \! \hat{\mathcal{P}}_{1,k}(\tau) \widetilde{\Lambda} \! + \! \hat{\mathcal{P}}_{2,k}
(\tau) \widetilde{\Lambda}^{2} \! + \! \mathcal{O} \big(\hat{\mathbb{E}}_{k}(\tau) \widetilde{\Lambda}^{3} \big) \right) 
\! \tilde{\Phi}_{k}(\widetilde{\Lambda}).
\end{equation}

\pmb{(iii)} The essence of this step is to transform the coefficient matrix $\hat{\mathcal{P}}_{1,k}(\tau)$ (cf. definition 
\eqref{prcy5}) into diagonal form. Let $\tilde{\Phi}_{k}(\widetilde{\Lambda})$ be a solution of equation \eqref{prcy8}; 
then, applying the transformation $\mathfrak{F}_{3}$, where
\begin{equation} \label{prcy9} 
\mathcal{G}_{0,k} \! := \! \left(\dfrac{\hat{\mathcal{C}}_{1}}{2 \lambda^{\ast}_{1}
(k)} \right)^{1/2} 
\begin{pmatrix}
\frac{\hat{\mathcal{A}}_{1}+ \lambda^{\ast}_{1}(k)}{\hat{\mathcal{C}}_{1}} & 
\frac{\hat{\mathcal{A}}_{1}- \lambda^{\ast}_{1}(k)}{\hat{\mathcal{C}}_{1}} \\
1 & 1
\end{pmatrix}, \quad k \! = \! \pm 1,
\end{equation}
with $\hat{\mathcal{A}}_{1}$ and $\hat{\mathcal{C}}_{1}$ given in equation \eqref{prcy5}, and
\begin{equation} \label{prcy10} 
\lambda^{\ast}_{1}(k) \! := \! \mi 4 \sqrt{3} \mathcal{Z}_{k} \! = \! \mi 4 \sqrt{3} \left(1 \! - \! \dfrac{1}{6} \tilde{r}_{0}(\tau) 
\tau^{-1/3} \! + \! \dfrac{1}{48} \big(\tilde{r}_{0}(\tau) \tau^{-1/3} \big)^{2} \right)^{1/2},
\end{equation}
one shows that
\begin{equation} \label{prcy11} 
\dfrac{\partial \Phi_{k}^{\sharp}(\widetilde{\Lambda})}{\partial \widetilde{\Lambda}} \underset{\tau \to +\infty}{=} 
\left(\mathcal{P}_{0,k}^{\vartriangle}(\tau) \! + \! \mathcal{P}_{1,k}^{\vartriangle}(\tau) \widetilde{\Lambda} \! + \! 
\mathcal{P}_{2,k}^{\vartriangle}(\tau) \widetilde{\Lambda}^{2} \! + \! \mathcal{O} \big(\mathcal{G}_{0,k}^{-1} 
\hat{\mathbb{E}}_{k}(\tau) \mathcal{G}_{0,k} \widetilde{\Lambda}^{3} \big) \right) \! \Phi_{k}^{\sharp}(\widetilde{\Lambda}),
\end{equation}
where
\begin{gather}
\mathcal{P}_{0,k}^{\vartriangle}(\tau) \! := \! \mathcal{G}_{0,k}^{-1} \hat{\mathcal{P}}_{0,k}(\tau) \mathcal{G}_{0,k} 
\! = \! \mathfrak{A}_{k} \sigma_{3} \! + \! \mathfrak{B}_{k} \sigma_{+} \! + \! \mathfrak{C}_{k} \sigma_{-}, \label{prcy12} \\
\mathcal{P}_{1,k}^{\vartriangle}(\tau) \! := \! \mathcal{G}_{0,k}^{-1} \hat{\mathcal{P}}_{1,k}(\tau) \mathcal{G}_{0,k} 
\! = \! \mi 4 \sqrt{3} \mathcal{Z}_{k} \sigma_{3}, \label{prcy13} \\
\mathcal{P}_{2,k}^{\vartriangle}(\tau) \! := \! \mathcal{G}_{0,k}^{-1} \hat{\mathcal{P}}_{2,k}(\tau) \mathcal{G}_{0,k} 
\! = \! \mathfrak{A}_{0,k}^{\sharp} \sigma_{3} \! + \! \mathfrak{B}_{0,k}^{\sharp} \sigma_{+} \! + \! 
\mathfrak{C}_{0,k}^{\sharp} \sigma_{-}, \label{prcy14}
\end{gather}
with
\begin{align}
\mathfrak{A}_{k} \! =& \, \dfrac{1}{(6 \varepsilon b)^{1/2} \mathcal{Z}_{k}} \! \left(-\dfrac{\mi \alpha_{k}
(\varepsilon b)^{1/2}}{2 \sqrt{2}} \tilde{r}_{0}(\tau)(-4 \! + \! \tilde{r}_{0}(\tau) \tau^{-1/3}) \! - \! 2 \mi 
(\varepsilon b)^{2/3} \me^{-\mi 2 \pi k/3}v_{0}(\tau) \right. \nonumber \\
-&\left. \mi (\varepsilon b)^{1/3} \! \left(\dfrac{(\varepsilon b)^{1/3} \me^{\mi \pi k/3}(\tilde{r}_{0}(\tau) \! + \! 2v_{0}(\tau)) 
\! + \! 2(a \! - \! \mi/2) \me^{\mi 2 \pi k/3} \tau^{-1/3}}{1 \! + \! v_{0}(\tau) \tau^{-1/3}} \right) \right), \label{prcy15} \\
\mathfrak{B}_{k} \! =& \, \dfrac{1}{(6 \varepsilon b)^{1/2} \mathcal{Z}_{k}} \! \left(-\dfrac{\mi \alpha_{k}
(\varepsilon b)^{1/2}}{2 \sqrt{2}} \tilde{r}_{0}(\tau)(-4 \! + \! \tilde{r}_{0}(\tau) \tau^{-1/3} \! - \! 4 \sqrt{3} \mathcal{Z}_{k}) 
\! - \! 2 \mi (\varepsilon b)^{2/3} \me^{-\mi 2 \pi k/3}v_{0}(\tau) \right. \nonumber \\
+&\left. \mi (\varepsilon b)^{1/3} \! \left(\dfrac{(\varepsilon b)^{1/3} \me^{\mi \pi k/3}(\tilde{r}_{0}(\tau) \! + \! 2v_{0}(\tau)) 
\! + \! 2(a \! - \! \mi/2) \me^{\mi 2 \pi k/3} \tau^{-1/3}}{1 \! + \! v_{0}(\tau) \tau^{-1/3}} \right) \right. \nonumber \\
\times&\left. \left(1 \! + \! \frac{1}{16}(-4 \! + \! \tilde{r}_{0}(\tau) \tau^{-1/3})(-4 \! + \! \tilde{r}_{0}(\tau) \tau^{-1/3} 
\! - \! 4 \sqrt{3} \mathcal{Z}_{k}) \right) \right), \label{prcy16} \\
\mathfrak{C}_{k} \! =& \, \dfrac{1}{(6 \varepsilon b)^{1/2} \mathcal{Z}_{k}} \! \left(\dfrac{\mi \alpha_{k}
(\varepsilon b)^{1/2}}{2 \sqrt{2}} \tilde{r}_{0}(\tau)(-4 \! + \! \tilde{r}_{0}(\tau) \tau^{-1/3} \! + \! 4 \sqrt{3} \mathcal{Z}_{k}) 
\! + \! 2 \mi (\varepsilon b)^{2/3} \me^{-\mi 2 \pi k/3}v_{0}(\tau) \right. \nonumber \\
-&\left. \mi (\varepsilon b)^{1/3} \! \left(\dfrac{(\varepsilon b)^{1/3} \me^{\mi \pi k/3}(\tilde{r}_{0}(\tau) \! + \! 2v_{0}(\tau)) 
\! + \! 2(a \! - \! \mi/2) \me^{\mi 2 \pi k/3} \tau^{-1/3}}{1 \! + \! v_{0}(\tau) \tau^{-1/3}} \right) \right. \nonumber \\
\times&\left. \left(1 \! + \! \dfrac{1}{16}(-4 \! + \! \tilde{r}_{0}(\tau) \tau^{-1/3})(-4 \! + \! \tilde{r}_{0}(\tau) \tau^{-1/3} 
\! + \! 4 \sqrt{3} \mathcal{Z}_{k}) \right) \right), \label{prcy17} \\
\mathfrak{A}_{0,k}^{\sharp} \! =& \, -\dfrac{\mi (\varepsilon b)^{1/3} \me^{-\mi \pi k/3} \tau^{-1/3}}{2(6 \varepsilon b)^{1/2} 
\mathcal{Z}_{k}} \! \left(48 \! + \! (-2 \! + \! \tilde{r}_{0}(\tau) \tau^{-1/3})(-4 \! + \! \tilde{r}_{0}(\tau) \tau^{-1/3}) \right), 
\label{prcy18} \\
\mathfrak{B}_{0,k}^{\sharp} \! =& \, \dfrac{\mi (\varepsilon b)^{1/3} \me^{-\mi \pi k/3} \tau^{-1/3}}{2(6 \varepsilon b)^{1/2} 
\mathcal{Z}_{k}}(-4 \! + \! \tilde{r}_{0}(\tau) \tau^{-1/3} \! - \! 4 \sqrt{3} \mathcal{Z}_{k})(-4 \! + \! \tfrac{1}{2} \tilde{r}_{0}(\tau) 
\tau^{-1/3}), \label{prcy19} \\
\mathfrak{C}_{0,k}^{\sharp} \! =& \, -\dfrac{\mi (\varepsilon b)^{1/3} \me^{-\mi \pi k/3} \tau^{-1/3}}{2(6 \varepsilon b)^{1/2} 
\mathcal{Z}_{k}}(-4 \! + \! \tilde{r}_{0}(\tau) \tau^{-1/3} \! + \! 4 \sqrt{3} \mathcal{Z}_{k})(-4 \! + \! \tfrac{1}{2} \tilde{r}_{0}(\tau) 
\tau^{-1/3}). \label{prcy20}
\end{align}
Observe that $\tr (\mathcal{P}_{0,k}^{\vartriangle}(\tau)) \! = \! \tr (\mathcal{P}_{1,k}^{\vartriangle}(\tau)) \! = \! \tr 
(\mathcal{P}_{2,k}^{\vartriangle}(\tau)) \! = \! \tr (\mathcal{G}_{0,k}^{-1} \hat{\mathbb{E}}_{k}(\tau) \mathcal{G}_{0,k}) \! = \! 0$. 

\pmb{(iv)} The idea behind the transformation for equation \eqref{prcy11} that is subsumed in this step is to put the 
coefficient matrix $\mathcal{P}_{0,k}^{\vartriangle}(\tau)$ (cf. definition \eqref{prcy12}) into Jordan canonical form, 
namely, to find a unimodular, $\tau$-dependent function $\mathcal{G}_{1,k}$ such that
\begin{equation} \label{prcy21} 
\mathcal{G}_{1,k}^{-1} \mathcal{P}_{0,k}^{\vartriangle}(\tau) \mathcal{G}_{1,k} \! = \! \mi \omega_{0,k} \sigma_{3} 
\! + \! \tau^{1/3} \sigma_{+}, \quad k \! = \! \pm 1,
\end{equation}
where (cf. equations \eqref{expforeych}, \eqref{expforkapp}, and \eqref{prcy15}--\eqref{prcy17})
\begin{align} \label{prcy22} 
\omega_{0,k}^{2} \! := \! \det (\mathcal{P}_{0,k}^{\vartriangle}(\tau)) 
\! =& \, \varkappa_{0}^{2}(\tau) \! + \! \dfrac{4(a \! - \! \mi/2)v_{0}(\tau) 
\tau^{-1/3}}{1 \! + \! v_{0}(\tau) \tau^{-1/3}} \! = \! 4 \! \left((a \! - \! \mi/2) 
\! + \! \alpha_{k}^{-2} \tau^{2/3} \hat{h}_{0}(\tau) \right) \nonumber \\
=& \, -\alpha_{k}^{2} \! \left(\dfrac{8v_{0}^{2}(\tau) \! + \! 4v_{0}(\tau) 
\tilde{r}_{0}(\tau) \! - \! (\tilde{r}_{0}(\tau))^{2} \! - \! v_{0}(\tau)(\tilde{r}_{0}
(\tau))^{2} \tau^{-1/3}}{1 \! + \! v_{0}(\tau) \tau^{-1/3}} \right) \nonumber \\
+& \, \dfrac{4(a \! - \! \mi/2)v_{0}(\tau) \tau^{-1/3}}{1 \! + \! v_{0}(\tau) 
\tau^{-1/3}};
\end{align}
the following lower-triangular solution for $\mathcal{G}_{1,k}$ is chosen:
\begin{equation} \label{prcy23} 
\mathcal{G}_{1,k} \! = \! \mathfrak{B}_{k}^{\frac{1}{2} \sigma_{3}} 
\tau^{-\frac{1}{6} \sigma_{3}} \! \left(\mathrm{I} \! + \! (\mi \omega_{0,k} 
\! - \! \mathfrak{A}_{k}) \tau^{-1/3} \sigma_{-} \right), \quad k \! = \! \pm 1.
\end{equation}
Let $\Phi_{k}^{\sharp}(\widetilde{\Lambda})$ solve equation~\eqref{prcy11}; then, applying the transformation 
$\mathfrak{F}_{4}$, one shows that
\begin{equation} \label{prcy24} 
\dfrac{\partial \hat{\Phi}_{k}(\widetilde{\Lambda})}{\partial \widetilde{\Lambda}} \underset{\tau \to +\infty}{=} 
\left(\mathcal{P}_{0,k}^{\triangledown}(\tau) \! + \! \mathcal{P}_{1,k}^{\triangledown}(\tau) \widetilde{\Lambda} 
\! + \! \mathcal{P}_{2,k}^{\triangledown}(\tau) \widetilde{\Lambda}^{2} \! + \! \mathcal{O} \big(\mathcal{G}_{1,k}^{-1} 
\mathcal{G}_{0,k}^{-1} \hat{\mathbb{E}}_{k}(\tau) \mathcal{G}_{0,k} \mathcal{G}_{1,k} \widetilde{\Lambda}^{3} \big) 
\right) \! \hat{\Phi}_{k}(\widetilde{\Lambda}),
\end{equation}
where
\begin{align}
\mathcal{P}_{0,k}^{\triangledown}(\tau) \! :=& \, \mathcal{G}_{1,k}^{-1} \mathcal{P}_{0,k}^{\vartriangle}(\tau) 
\mathcal{G}_{1,k} \! = \! \mi \omega_{0,k} \sigma_{3} \! + \! \tau^{1/3} \sigma_{+}, \label{prcy25} \\
\mathcal{P}_{1,k}^{\triangledown}(\tau) \! :=& \, \mathcal{G}_{1,k}^{-1} \mathcal{P}_{1,k}^{\vartriangle}(\tau) 
\mathcal{G}_{1,k} \! = \! \mi 4 \sqrt{3} \mathcal{Z}_{k} \sigma_{3} \! - \! \mi 8 \sqrt{3}(\mi \omega_{0,k} \! - \! 
\mathfrak{A}_{k}) \mathcal{Z}_{k} \tau^{-1/3} \sigma_{-}, \label{prcy26} \\
\mathcal{P}_{2,k}^{\triangledown}(\tau) \! :=& \, \mathcal{G}_{1,k}^{-1} \mathcal{P}_{2,k}^{\vartriangle}(\tau) 
\mathcal{G}_{1,k} \nonumber \\
=& 
\begin{pmatrix}
\mathfrak{A}_{0,k}^{\sharp} \! + \! \frac{(\mi \omega_{0,k}-\mathfrak{A}_{k}) \mathfrak{B}_{0,k}^{\sharp}}{\mathfrak{B}_{k}} 
& \frac{\mathfrak{B}_{0,k}^{\sharp} \tau^{1/3}}{\mathfrak{B}_{k}} \\
\frac{2(\mi \omega_{0,k}-\mathfrak{A}_{k})(\mathfrak{A}_{k} \mathfrak{B}_{0,k}^{\sharp}-\mathfrak{B}_{k} 
\mathfrak{A}_{0,k}^{\sharp})+(\mathfrak{B}_{k} \mathfrak{C}_{0,k}^{\sharp}-\mathfrak{C}_{k} \mathfrak{B}_{0,k}^{\sharp}) 
\mathfrak{B}_{k}}{\mathfrak{B}_{k} \tau^{1/3}} & -(\mathfrak{A}_{0,k}^{\sharp} \! + \! \frac{(\mi \omega_{0,k}
-\mathfrak{A}_{k}) \mathfrak{B}_{0,k}^{\sharp}}{\mathfrak{B}_{k}})
\end{pmatrix}. \label{prcy27}
\end{align}
Note that, at this stage, the matrix $\mathcal{P}_{1,k}^{\triangledown}(\tau)$ is not diagonal; instead, it now contains an 
additional, lower off-diagonal contribution.

\pmb{(v)} This step entails a straightforward $\tau$-dependent scaling. Let $\hat{\Phi}_{k}(\widetilde{\Lambda})$ solve 
equation \eqref{prcy24}; then, applying the transformation $\mathfrak{F}_{5}$, one shows that, for $k \! = \! \pm 1$,
\begin{align} \label{prcy28} 
\dfrac{\partial \hat{\Phi}_{0,k}(\widetilde{\Lambda})}{\partial \widetilde{
\Lambda}} \underset{\tau \to +\infty}{=}& \, \left(\tilde{\mathcal{P}}_{0,k}^{
\blacktriangle}(\tau) \! + \! \tilde{\mathcal{P}}_{1,k}^{\blacktriangle}(\tau) 
\widetilde{\Lambda} \! + \! \tilde{\mathcal{P}}_{2,k}^{\blacktriangle}(\tau) 
\widetilde{\Lambda}^{2} \right. \nonumber \\
+&\left. \, \mathcal{O} \! \left(\tau^{-\frac{1}{6} \sigma_{3}} \mathcal{G}_{1,
k}^{-1} \mathcal{G}_{0,k}^{-1} \hat{\mathbb{E}}_{k}(\tau) \mathcal{G}_{0,k} 
\mathcal{G}_{1,k} \tau^{\frac{1}{6} \sigma_{3}} \widetilde{\Lambda}^{3} 
\right) \right) \! \hat{\Phi}_{0,k}(\widetilde{\Lambda}),
\end{align}
where
\begin{align}
\tilde{\mathcal{P}}_{0,k}^{\blacktriangle}(\tau) \! :=& \, \tau^{-\frac{1}{6} 
\sigma_{3}} \mathcal{P}_{0,k}^{\triangledown}(\tau) \tau^{\frac{1}{6} 
\sigma_{3}} \! = \! \mi \omega_{0,k} \sigma_{3} \! + \! \sigma_{+}, 
\label{prcy29} \\
\tilde{\mathcal{P}}_{1,k}^{\blacktriangle}(\tau) \! :=& \, \tau^{-\frac{1}{6} 
\sigma_{3}} \mathcal{P}_{1,k}^{\triangledown}(\tau) \tau^{\frac{1}{6} 
\sigma_{3}} \! = \! \mi 4 \sqrt{3} \mathcal{Z}_{k} \sigma_{3} \! - \! \mi 
8 \sqrt{3}(\mi \omega_{0,k} \! - \! \mathfrak{A}_{k}) \mathcal{Z}_{k} 
\sigma_{-}, \label{prcy30} \\
\tilde{\mathcal{P}}_{2,k}^{\blacktriangle}(\tau) \! :=& \, \tau^{-\frac{1}{6} 
\sigma_{3}} \mathcal{P}_{2,k}^{\triangledown}(\tau) \tau^{\frac{1}{6} 
\sigma_{3}} \nonumber \\
=& 
\begin{pmatrix}
\mathfrak{A}_{0,k}^{\sharp} \! + \! \frac{(\mi \omega_{0,k}-
\mathfrak{A}_{k}) \mathfrak{B}_{0,k}^{\sharp}}{\mathfrak{B}_{k}} & 
\frac{\mathfrak{B}_{0,k}^{\sharp}}{\mathfrak{B}_{k}} \\
\frac{2(\mi \omega_{0,k}-\mathfrak{A}_{k})(\mathfrak{A}_{k} 
\mathfrak{B}_{0,k}^{\sharp}-\mathfrak{B}_{k} \mathfrak{A}_{0,k}^{\sharp}) 
+(\mathfrak{B}_{k} \mathfrak{C}_{0,k}^{\sharp}-\mathfrak{C}_{k} 
\mathfrak{B}_{0,k}^{\sharp}) \mathfrak{B}_{k}}{\mathfrak{B}_{k}} & 
-(\mathfrak{A}_{0,k}^{\sharp} \! + \! \frac{(\mi \omega_{0,k}-
\mathfrak{A}_{k}) \mathfrak{B}_{0,k}^{\sharp}}{\mathfrak{B}_{k}})
\end{pmatrix}. \label{prcy31}
\end{align}

\pmb{(vi)} The purpose of this step is to transform the coefficient matrix $\tilde{\mathcal{P}}_{0,k}^{\blacktriangle}
(\tau)$ (cf. equation \eqref{prcy29}) into off-diagonal form. Let $\hat{\Phi}_{0,k}(\widetilde{\Lambda})$ solve 
equation \eqref{prcy28}; then, applying the transformation $\mathfrak{F}_{6}$, one shows that, for $k \! = \! \pm 1$,
\begin{align} \label{prcy32} 
\dfrac{\partial \Phi_{0,k}(\widetilde{\Lambda})}{\partial \widetilde{\Lambda}} \underset{\tau \to +\infty}{=}& \, 
\left(
\begin{pmatrix}
0 & 1 \\
-\omega_{0,k}^{2} & 0
\end{pmatrix} \! + \! 
\begin{pmatrix}
\mi 4 \sqrt{3} \mathcal{Z}_{k} & 0 \\
\mi 8 \sqrt{3} \mathcal{Z}_{k} \mathfrak{A}_{k} & -\mi 4 \sqrt{3} \mathcal{Z}_{k}
\end{pmatrix} \! \widetilde{\Lambda} \! + \! 
\begin{pmatrix}
\mathfrak{P}_{0,k}^{\ast} & \mathfrak{Q}_{0,k}^{\ast} \\
\mathfrak{R}_{0,k}^{\ast} & -\mathfrak{P}_{0,k}^{\ast}
\end{pmatrix} \! \widetilde{\Lambda}^{2} \right. \nonumber \\
+&\left. \, \mathcal{O} \big(\mathbb{E}_{k}^{\ast}(\tau) \widetilde{\Lambda}^{3} \big) \right) \! 
\Phi_{0,k}(\widetilde{\Lambda}),
\end{align}
where
\begin{gather}
\mathfrak{P}_{0,k}^{\ast} \! := \! \mathfrak{A}_{0,k}^{\sharp} \! - \! \mathfrak{B}_{0,k}^{\sharp} \mathfrak{A}_{k} 
\mathfrak{B}_{k}^{-1}, \label{prcy33} \\
\mathfrak{Q}_{0,k}^{\ast} \! := \! \mathfrak{B}_{0,k}^{\sharp} \mathfrak{B}_{k}^{-1}, \label{prcy34} \\
\mathfrak{R}_{0,k}^{\ast} \! := \! -\mathfrak{B}_{0,k}^{\sharp} \mathfrak{A}_{k}^{2} \mathfrak{B}_{k}^{-1} \! + \! 2 
\mathfrak{A}_{k} \mathfrak{A}_{0,k}^{\sharp} \! + \! \mathfrak{B}_{k} \mathfrak{C}_{0,k}^{\sharp}, \label{prcy35}
\end{gather}
and
\begin{gather} \label{prcy36} 
\mathbb{E}_{k}^{\ast}(\tau) \! := \! \left(\mathrm{I} \! + \! \mi \omega_{0,k} \sigma_{-} \right) \tau^{-\frac{1}{6} 
\sigma_{3}} \mathcal{G}_{1,k}^{-1} \mathcal{G}_{0,k}^{-1} \hat{\mathbb{E}}_{k}(\tau) \mathcal{G}_{0,k} 
\mathcal{G}_{1,k} \tau^{\frac{1}{6} \sigma_{3}} \left(\mathrm{I} \! - \! \mi \omega_{0,k} \sigma_{-} \right).
\end{gather}

\pmb{(vii)} This step, in conjunction with steps \pmb{(viii)} and \pmb{(x)} below, is precipitated by the fact that, in 
order to derive a---canonical---model problem solvable in terms of parabolic-cylinder functions (see step \pmb{(xi)} 
below), one must eliminate the coefficient matrix of the $\widetilde{\Lambda}^{2}$ term {}from equation \eqref{prcy32}; 
in particular, this step focuses on the excision of the $(1 \, 2)$-element. Let $\Phi_{0,k}(\widetilde{\Lambda})$ solve 
equation \eqref{prcy32}; then, applying the transformation $\mathfrak{F}_{7}$, with $\tau$-dependent parameter 
$\ell_{0,k}$, one shows, via the conditions \eqref{iden5}, that, for $k \! = \! \pm 1$,
\begin{align} \label{prcy37} 
\dfrac{\partial \Phi_{0,k}^{\flat}(\widetilde{\Lambda})}{\partial \widetilde{\Lambda}} \underset{\tau \to +\infty}{=}& 
\, \left(
\begin{pmatrix}
0 & -\ell_{0,k} \! + \! 1 \\
-\omega_{0,k}^{2} & 0
\end{pmatrix} \! + \! 
\begin{pmatrix}
\mi 4 \sqrt{3} \mathcal{Z}_{k} \! + \! \omega_{0,k}^{2} \ell_{0,k} & 0 \\
\mi 8 \sqrt{3} \mathcal{Z}_{k} \mathfrak{A}_{k} & -\mi 4 \sqrt{3} \mathcal{Z}_{k} \! - \! \omega_{0,k}^{2} \ell_{0,k}
\end{pmatrix} \! \widetilde{\Lambda} \right. \nonumber \\
+&\left. \, 
\begin{pmatrix}
-\mi 8 \sqrt{3} \mathcal{Z}_{k} \mathfrak{A}_{k} \ell_{0,k} \! + \! 
\mathfrak{P}_{0,k}^{\ast} & \omega_{0,k}^{2} \ell_{0,k}^{2} \! + \! \mi 
8 \sqrt{3} \mathcal{Z}_{k} \ell_{0,k} \! + \! \mathfrak{Q}_{0,k}^{\ast} \\
\mathfrak{R}_{0,k}^{\ast} & \mi 8 \sqrt{3} \mathcal{Z}_{k} 
\mathfrak{A}_{k} \ell_{0,k} \! - \! \mathfrak{P}_{0,k}^{\ast}
\end{pmatrix} \! \widetilde{\Lambda}^{2} \right. \nonumber \\
+&\left. \, \mathcal{O} \big(\mathbb{E}^{\triangledown}_{k}(\ell_{0,k};\tau) \widetilde{\Lambda}^{3} \big) \right) 
\! \Phi_{0,k}^{\flat}(\widetilde{\Lambda}),
\end{align}
where
\begin{equation} \label{prcy38} 
\mathbb{E}^{\triangledown}_{k}(\ell_{0,k};\tau) \! := \!  \mathbb{E}^{\ast}_{k}(\tau) \! + \! 
\begin{pmatrix}
-\mathfrak{R}_{0,k}^{\ast} \ell_{0,k} & -\mi 8 \sqrt{3} \mathcal{Z}_{k} \mathfrak{A}_{k} \ell_{0,k}^{2} 
\! + \! 2 \mathfrak{P}_{0,k}^{\ast} \ell_{0,k} \\0 & \mathfrak{R}_{0,k}^{\ast} \ell_{0,k}
\end{pmatrix},
\end{equation}
with $\mathbb{E}_{k}^{\ast}(\tau)$ defined by equation \eqref{prcy36}. One now chooses $\ell_{0,k}$ so that the 
$(1 \, 2)$-element of the coefficient matrix of the $\widetilde{\Lambda}^{2}$ term in equation \eqref{prcy37} is equal 
to zero, that is, $\omega_{0,k}^{2} \ell_{0,k}^{2} \! + \! \mi 8 \sqrt{3} \mathcal{Z}_{k} \ell_{0,k} \! + \! \mathfrak{Q}_{0,k}^{\ast} 
\! = \! 0$; the roots are given by
\begin{equation} \label{prcy39} 
\ell_{0,k}^{\pm} \! = \! \dfrac{-\mi 8 \sqrt{3} \mathcal{Z}_{k} \! \pm \! \sqrtsign{(\mi 8 \sqrt{3} \mathcal{Z}_{k})^{2} \! - \! 
4 \omega_{0,k}^{2} \mathfrak{Q}_{0,k}^{\ast}}}{2 \omega_{0,k}^{2}}, \quad k \! = \! \pm 1.
\end{equation}
Noting {}from the conditions \eqref{iden5}, the asymptotics \eqref{tr1} and \eqref{tr3}, equations \eqref{prcy16} 
and \eqref{prcy19}, and the definitions \eqref{prcy10}, \eqref{prcy22}, and \eqref{prcy34} that 
$\mathcal{Z}_{k} \genfrac{}{}{0pt}{3}{=}{\tau \to +\infty} 1 \! + \! \mathcal{O}(\tau^{-2/3})$, $\omega_{0,k}^{2} 
\genfrac{}{}{0pt}{3}{=}{\tau \to +\infty} \mathcal{O}(\tau^{-2/3})$, and $\mathfrak{Q}_{0,k}^{\ast} 
\genfrac{}{}{0pt}{3}{=}{\tau \to +\infty} \mathcal{O}(1)$, it follows that, for the class of functions consistent with 
the conditions \eqref{iden5}, the `$+$-root' in equation \eqref{prcy39} is chosen:
\begin{equation} \label{prcy40} 
\ell_{0,k} \! := \! \ell_{0,k}^{+} \! = \! \dfrac{-\mi 8 \sqrt{3} \mathcal{Z}_{k} \! + \! \sqrtsign{(\mi 8 \sqrt{3} \mathcal{Z}_{k})^{2} 
\! - \! 4 \omega_{0,k}^{2} \mathfrak{Q}_{0,k}^{\ast}}}{2 \omega_{0,k}^{2}}.
\end{equation}
Via the formula for the $\tau$-dependent parameter $\ell_{0,k} \! := \! \ell_{0,k}^{+}$ given in equation \eqref{prcy40}, one 
rewrites equation \eqref{prcy37} as follows: for $k \! = \! \pm 1$,
\begin{align} \label{prcy41} 
\dfrac{\partial \Phi_{0,k}^{\flat}(\widetilde{\Lambda})}{\partial \widetilde{\Lambda}} \underset{\tau \to +\infty}{=}& 
\, \left(
\begin{pmatrix}
0 & -\ell_{0,k}^{+} \! + \! 1 \\
-\omega_{0,k}^{2} & 0
\end{pmatrix} \! + \! 
\begin{pmatrix}
\mi 4 \sqrt{3} \mathcal{Z}_{k} \! + \! \omega_{0,k}^{2} \ell_{0,k}^{+} & 0 \\
\mi 8 \sqrt{3} \mathcal{Z}_{k} \mathfrak{A}_{k} & -\mi 4 \sqrt{3} \mathcal{Z}_{k} \! - \! \omega_{0,k}^{2} \ell_{0,k}^{+}
\end{pmatrix} \! \widetilde{\Lambda} \right. \nonumber \\
+&\left. \, 
\begin{pmatrix}
-\mi 8 \sqrt{3} \mathcal{Z}_{k} \mathfrak{A}_{k} \ell_{0,k}^{+} \! + \! 
\mathfrak{P}_{0,k}^{\ast} & 0 \\
\mathfrak{R}_{0,k}^{\ast} & \mi 8 \sqrt{3} \mathcal{Z}_{k} \mathfrak{A}_{k} 
\ell_{0,k}^{+} \! - \! \mathfrak{P}_{0,k}^{\ast}
\end{pmatrix} \! \widetilde{\Lambda}^{2} \! + \! \mathcal{O} \big(\mathbb{E}_{k}^{\triangledown}
(\ell_{0,k}^{+};\tau) \widetilde{\Lambda}^{3} \big) \right) \! \Phi_{0,k}^{\flat}(\widetilde{\Lambda}).
\end{align}

\pmb{(viii)} This step focuses on the excision of the $(2 \, 1)$-element {}from the coefficient matrix of the 
$\widetilde{\Lambda}^{2}$ term in equation \eqref{prcy41}. Let $\Phi_{0,k}^{\flat}(\widetilde{\Lambda})$ 
solve equation \eqref{prcy41}; then, under the action of the transformation $\mathfrak{F}_{8}$, with 
$\tau$-dependent parameter $\ell_{1,k}$, one shows that, for $k \! = \! \pm 1$,
\begin{align} \label{prcy42} 
\dfrac{\partial \Phi_{0,k}^{\sharp}(\widetilde{\Lambda})}{\partial \widetilde{\Lambda}} \underset{\tau \to +\infty}{=}& 
\, \left(
\begin{pmatrix}
0 & -\ell_{0,k}^{+} \! + \! 1 \\
-\omega_{0,k}^{2} \! - \! \ell_{1,k} & 0
\end{pmatrix} \! + \! 
\left((\mi 4 \sqrt{3} \mathcal{Z}_{k} \! + \! \omega_{0,k}^{2} \ell_{0,k}^{+} \! + \! \ell_{1,k}(-\ell_{0,k}^{+} \! + \! 1)) 
\sigma_{3} \right. \right. \nonumber \\
+&\left. \left. \mi 8 \sqrt{3} \mathcal{Z}_{k} \mathfrak{A}_{k} \sigma_{-} \right) \! \widetilde{\Lambda} \! + \! 
\left((\mathfrak{R}_{0,k}^{\ast} \! - \! 2(\mi 4 \sqrt{3} \mathcal{Z}_{k} \! + \! \omega_{0,k}^{2} \ell_{0,k}^{+}) 
\ell_{1,k} \! - \! \ell_{1,k}^{2}(-\ell_{0,k}^{+} \! + \! 1)) \sigma_{-} \right. \right. \nonumber \\
-&\left. \left. (\mi 8 \sqrt{3} \mathcal{Z}_{k} \mathfrak{A}_{k} \ell_{0,k}^{+} \! - \! \mathfrak{P}_{0,k}^{\ast}) 
\sigma_{3} \right) \! \widetilde{\Lambda}^{2} \! + \! \mathcal{O} 
\big(\overset{\ast}{\mathbb{E}}_{k}^{\raise-6.75pt\hbox{$\scriptstyle \blacktriangle$}}(\ell_{0,k}^{+},\ell_{1,k};\tau) 
\widetilde{\Lambda}^{3} \big) \right) \! \Phi_{0,k}^{\sharp}(\widetilde{\Lambda}),
\end{align}
where
\begin{equation} \label{prcy43} 
\overset{\ast}{\mathbb{E}}_{k}^{\raise-6.75pt\hbox{$\scriptstyle \blacktriangle$}}(\ell_{0,k}^{+},\ell_{1,k};\tau) 
\! := \! \mathbb{E}_{k}^{\triangledown}(\ell_{0,k}^{+};\tau) \! + \! 2 \ell_{1,k}(-\mathfrak{P}_{0,k}^{\ast} \! + \! 
\mi 8 \sqrt{3} \mathcal{Z}_{k} \mathfrak{A}_{k} \ell_{0,k}^{+}) \sigma_{-}.
\end{equation}
One now chooses $\ell_{1,k}$ so that the $(2 \, 1)$-element of the coefficient matrix of the $\widetilde{\Lambda}^{2}$ 
term in equation \eqref{prcy42} vanishes, that is, $(-\ell_{0,k}^{+} \! + \! 1) \ell_{1,k}^{2} \! + \! 2(\mi 4 \sqrt{3} 
\mathcal{Z}_{k} \! + \! \omega_{0,k}^{2} \ell_{0,k}^{+}) \ell_{1,k} \! - \! \mathfrak{R}_{0,k}^{\ast} \! = \! 0$; the roots 
are given by
\begin{equation} \label{prcy44} 
\ell_{1,k}^{\pm} \! = \! \dfrac{-(\mi 4 \sqrt{3} \mathcal{Z}_{k} \! + \! \omega_{0,k}^{2} \ell_{0,k}^{+}) \! \pm \! 
\sqrt{(\mi 4 \sqrt{3} \mathcal{Z}_{k} \! + \! \omega_{0,k}^{2} \ell_{0,k}^{+})^{2} \! + \! \mathfrak{R}_{0,k}^{\ast}
(-\ell_{0,k}^{+} \! + \! 1)}}{-\ell_{0,k}^{+} \! + \! 1}, \quad k \! = \! \pm 1.
\end{equation}
Noting {}from the conditions \eqref{iden5}, the asymptotics \eqref{tr1} and \eqref{tr3}, equations 
\eqref{prcy15}--\eqref{prcy20}, and the definition \eqref{prcy35} that $\mathfrak{R}_{0,k}^{\ast} 
\genfrac{}{}{0pt}{3}{=}{\tau \to +\infty} \mathcal{O}(\tau^{-2/3})$, and, recalling ({}from step \pmb{(vii)} above) 
the asymptotics $\mathcal{Z}_{k} \genfrac{}{}{0pt}{3}{=}{\tau \to +\infty} 1 \! + \! \mathcal{O}(\tau^{-2/3})$, 
$\omega_{0,k}^{2} \genfrac{}{}{0pt}{3}{=}{\tau \to +\infty} \mathcal{O}(\tau^{-2/3})$, and $\mathfrak{Q}_{0,k}^{\ast} 
\genfrac{}{}{0pt}{3}{=}{\tau \to +\infty} \mathcal{O}(1)$, it follows {}from the definition~\eqref{prcy40} for 
$\ell_{0,k}^{+}$ that, for the class of functions consistent with the conditions \eqref{iden5}, the `$+$-root' in 
equation \eqref{prcy44} is taken:
\begin{equation} \label{prcy45} 
\ell_{1,k} \! := \! \ell_{1,k}^{+} \! = \! \dfrac{-(\mi 4 \sqrt{3} \mathcal{Z}_{k} \! + \! \omega_{0,k}^{2} \ell_{0,k}^{+}) 
\! + \! \chi_{k}(\tau)}{-\ell_{0,k}^{+} \! + \! 1},
\end{equation}
where
\begin{equation} \label{prcy46} 
\chi_{k}(\tau) \! := \! \sqrtsign{\big(\mi 4 \sqrt{3} \mathcal{Z}_{k} \! + \! \omega_{0,k}^{2} \ell_{0,k}^{+} \big)^{2} \! + \! 
\mathfrak{R}_{0,k}^{\ast}(-\ell_{0,k}^{+} \! + \! 1)}.
\end{equation}
Via the formula for the $\tau$-dependent parameter $\ell_{1,k} \! := \! \ell_{1,k}^{+}$ defined by equations \eqref{prcy45} 
and \eqref{prcy46}, one rewrites equation \eqref{prcy42} as follows: for $k \! = \! \pm 1$,
\begin{align} \label{prcy47} 
\dfrac{\partial \Phi_{0,k}^{\sharp}(\widetilde{\Lambda})}{\partial \widetilde{\Lambda}} \underset{\tau \to +\infty}{=}& 
\, \left(
\begin{pmatrix}
0 & -\ell_{0,k}^{+} \! + \! 1 \\
-\omega_{0,k}^{2} \! - \! \ell_{1,k}^{+} & 0
\end{pmatrix} \! + \! \left(\chi_{k}(\tau) \sigma_{3} \! + \! \mi 8 \sqrt{3} \mathcal{Z}_{k} \mathfrak{A}_{k} 
\sigma_{-} \right) \! \widetilde{\Lambda} \right. \nonumber \\
+&\left. \, \left(\mathfrak{P}_{0,k}^{\ast} \! - \! \mi 8 \sqrt{3} \mathcal{Z}_{k} \mathfrak{A}_{k} \ell_{0,k}^{+} \right) 
\! \widetilde{\Lambda}^{2} \sigma_{3} \! + \! \mathcal{O} \big(
\overset{\ast}{\mathbb{E}}_{k}^{\raise-6.75pt\hbox{$\scriptstyle \blacktriangle$}}(\ell_{0,k}^{+},\ell_{1,k}^{+};\tau) 
\widetilde{\Lambda}^{3} \big) \right) \! \Phi_{0,k}^{\sharp}(\widetilde{\Lambda}).
\end{align}

\pmb{(ix)} This step is necessitated by the fact that the coefficient matrix of the $\widetilde{\Lambda}$ term in 
equation \eqref{prcy47} remains to be re-diagonalised. Let $\Phi_{0,k}^{\sharp}(\widetilde{\Lambda})$ solve 
equation \eqref{prcy47}; then, under the action of the transformation $\mathfrak{F}_{9}$, where
\begin{equation} \label{prcy48} 
\mathcal{G}_{2,k} \! := \! 
\begin{pmatrix}
1 & 0 \\
\frac{\mi 4 \sqrt{3} \mathcal{Z}_{k} \mathfrak{A}_{k}}{\chi_{k}(\tau)} & 1
\end{pmatrix}, \quad k \! = \! \pm 1,
\end{equation}
with $\mathcal{Z}_{k}$, $\mathfrak{A}_{k}$, and $\chi_{k}(\tau)$ defined by equations \eqref{prcy10}, 
\eqref{prcy15}, and \eqref{prcy46}, respectively, one shows that
\begin{align} \label{prcy49} 
\dfrac{\partial \Phi_{0,k}^{\natural}(\widetilde{\Lambda})}{\partial \widetilde{\Lambda}} \underset{\tau \to +\infty}{=}& 
\, \left(
\begin{pmatrix}
\frac{\mi 4 \sqrt{3} \mathcal{Z}_{k} \mathfrak{A}_{k}}{\chi_{k}(\tau)}(-\ell_{0,k}^{+} \! + \! 1) & -\ell_{0,k}^{+} \! + \! 1 \\
-(\frac{\mi 4 \sqrt{3} \mathcal{Z}_{k} \mathfrak{A}_{k}}{\chi_{k}(\tau)})^{2}(-\ell_{0,k}^{+} \! + \! 1) \! - \! \ell_{1,k}^{+} 
\! - \! \omega_{0,k}^{2} & -\frac{\mi 4 \sqrt{3} \mathcal{Z}_{k} \mathfrak{A}_{k}}{\chi_{k}(\tau)}(-\ell_{0,k}^{+} \! + \! 1)
\end{pmatrix} \right. \nonumber \\
+&\left. \, \chi_{k}(\tau) \widetilde{\Lambda} \sigma_{3} \! + \! (\mathfrak{P}_{0,k}^{\ast} \! - \! \mi 8 \sqrt{3} 
\mathcal{Z}_{k} \mathfrak{A}_{k} \ell_{0,k}^{+}) 
\begin{pmatrix}
1 & 0 \\
-\frac{\mi 8 \sqrt{3} \mathcal{Z}_{k} \mathfrak{A}_{k}}{\chi_{k}(\tau)} & -1
\end{pmatrix} \! \widetilde{\Lambda}^{2} \right. \nonumber \\
+&\left. \, \mathcal{O} \big(\mathcal{G}_{2,k}^{-1} 
\overset{\ast}{\mathbb{E}}_{k}^{\raise-6.75pt\hbox{$\scriptstyle \blacktriangle$}}(\ell_{0,k}^{+},\ell_{1,k}^{+};\tau) 
\mathcal{G}_{2,k} \widetilde{\Lambda}^{3} \big) \right) \! \Phi_{0,k}^{\natural}(\widetilde{\Lambda}).
\end{align}

\pmb{(x)} This penultimate step focuses on the annihilation of the nilpotent coefficient sub-matrix of the 
$\widetilde{\Lambda}^{2}$ term in equation \eqref{prcy49}. Let $\Phi_{0,k}^{\natural}(\widetilde{\Lambda})$ solve 
equation \eqref{prcy49}; then, under the action of the transformation $\mathfrak{F}_{10}$, with $\tau$-dependent 
parameter $\ell_{2,k}$, one shows that, for $k \! = \! \pm 1$,
\begin{align} \label{prcy50} 
\dfrac{\partial \Phi_{k}^{\ast}(\widetilde{\Lambda})}{\partial \widetilde{\Lambda}} \underset{\tau \to +\infty}{=}& 
\, \left(
\begin{pmatrix}
\frac{\mi 4 \sqrt{3} \mathcal{Z}_{k} \mathfrak{A}_{k}}{\chi_{k}(\tau)}
(-\ell_{0,k}^{+} \! + \! 1) & -\ell_{0,k}^{+} \! + \! 1 \\
-(\frac{\mi 4 \sqrt{3} \mathcal{Z}_{k} \mathfrak{A}_{k}}{\chi_{k}(\tau)})^{2}
(-\ell_{0,k}^{+} \! + \! 1) \! - \! \ell_{1,k}^{+} \! - \! \ell_{2,k} \! - \! 
\omega_{0,k}^{2} & -\frac{\mi 4 \sqrt{3} \mathcal{Z}_{k} \mathfrak{A}_{k}}{
\chi_{k}(\tau)}(-\ell_{0,k}^{+} \! + \! 1)
\end{pmatrix}  \right. \nonumber \\
+&\left. 
\begin{pmatrix}
\chi_{k}(\tau) \! + \! \ell_{2,k}(-\ell_{0,k}^{+} \! + \! 1) & 0 \\
-\frac{\mi 8 \sqrt{3} \mathcal{Z}_{k} \mathfrak{A}_{k}}{\chi_{k}(\tau)} 
\ell_{2,k}(-\ell_{0,k}^{+} \! + \! 1) & -(\chi_{k}(\tau) \! + \! \ell_{2,k}
(-\ell_{0,k}^{+} \! + \! 1))
\end{pmatrix} \! \widetilde{\Lambda} \right. \nonumber \\
+&\left. \left( \left(-\ell_{2,k}^{2}(-\ell_{0,k}^{+} \! + \! 1) \! - \! 2 \ell_{2,k} \chi_{k}(\tau) \! - \! \frac{\mi 8 \sqrt{3} \mathcal{Z}_{k} 
\mathfrak{A}_{k}}{\chi_{k}(\tau)}(\mathfrak{P}_{0,k}^{\ast} \! - \! \mi 8 \sqrt{3} \mathcal{Z}_{k} \mathfrak{A}_{k} \ell_{0,k}^{+}) 
\right) \! \sigma_{-} \right. \right. \nonumber \\
+&\left. \left. (\mathfrak{P}_{0,k}^{\ast} \! - \! \mi 8 \sqrt{3} \mathcal{Z}_{k} \mathfrak{A}_{k} \ell_{0,k}^{+}) \sigma_{3} 
\right) \! \widetilde{\Lambda}^{2} \! + \! \mathcal{O} \big(
\overset{\ast}{\mathbb{E}}_{k}^{\raise-6.75pt\hbox{$\scriptstyle \ast$}}(\ell_{0,k}^{+},\ell_{1,k}^{+},
\ell_{2,k};\tau) \widetilde{\Lambda}^{3} \big) \right) \! \Phi_{k}^{\ast}(\widetilde{\Lambda}),
\end{align}
where
\begin{equation} \label{prcy51} 
\overset{\ast}{\mathbb{E}}_{k}^{\raise-6.75pt\hbox{$\scriptstyle \ast$}}(\ell_{0,k}^{+},\ell_{1,k}^{+},\ell_{2,k};\tau) \! := \! 
\mathcal{G}_{2,k}^{-1} \overset{\ast}{\mathbb{E}}_{k}^{\raise-6.75pt\hbox{$\scriptstyle \blacktriangle$}}
(\ell_{0,k}^{+},\ell_{1,k}^{+};\tau) \mathcal{G}_{2,k} \! - \! 2 \ell_{2,k}(\mathfrak{P}_{0,k}^{\ast} \! - \! \mi 8 \sqrt{3} 
\mathcal{Z}_{k} \mathfrak{A}_{k} \ell_{0,k}^{+}) \sigma_{-}.
\end{equation}
One now chooses $\ell_{2,k}$ so that the $(2 \, 1)$-element of the nilpotent coefficient matrix of the $\widetilde{\Lambda}^{2}$ 
terms in equation \eqref{prcy50} is equal to zero, that is, $(-\ell_{0,k}^{+} \! + \! 1) \ell_{2,k}^{2} \! + \! 2 \chi_{k}(\tau) \ell_{2,k} \! 
+ \! \mi 8 \sqrt{3} \mathcal{Z}_{k} \mathfrak{A}_{k} \chi_{k}^{-1}(\tau)(\mathfrak{P}_{0,k}^{\ast} \! - \! \mi 8 \sqrt{3} \mathcal{Z}_{k} 
\mathfrak{A}_{k} \ell_{0,k}^{+}) \! = \! 0$; the roots are given by
\begin{equation} \label{prcy52} 
\ell_{2,k}^{\pm} \! = \! \dfrac{-\chi_{k}(\tau) \! \pm \! \sqrt{\chi_{k}^{2}(\tau) \! - \! \mi 8 \sqrt{3} \mathcal{Z}_{k} 
\mathfrak{A}_{k} \chi_{k}^{-1}(\tau)(-\ell_{0,k}^{+} \! + \! 1)(\mathfrak{P}_{0,k}^{\ast} \! - \! \mi 8 \sqrt{3} 
\mathcal{Z}_{k} \mathfrak{A}_{k} \ell_{0,k}^{+})}}{-\ell_{0,k}^{+} \! + \! 1}, \quad k \! = \! \pm 1.
\end{equation}
Arguing as in steps \pmb{(vii)} and \pmb{(viii)} above, for the class of functions consistent with the conditions 
\eqref{iden5}, the `$+$-root' in equation~\eqref{prcy52} is taken:
\begin{equation} \label{prcy53} 
\ell_{2,k} \! := \! \ell_{2,k}^{+} \! = \! \dfrac{-\chi_{k}(\tau) \! + \! \mu_{k}(\tau)}{-\ell_{0,k}^{+} \! + \! 1},
\end{equation} 
where
\begin{equation} \label{prcy54} 
\mu_{k}(\tau) \! := \! \sqrtsign{\chi_{k}^{2}(\tau) \! - \! \mi 8 \sqrt{3} \mathcal{Z}_{k} \mathfrak{A}_{k} \chi_{k}^{-1}(\tau)
(-\ell_{0,k}^{+} \! + \! 1)(\mathfrak{P}_{0,k}^{\ast} \! - \! \mi 8 \sqrt{3} \mathcal{Z}_{k} \mathfrak{A}_{k} \ell_{0,k}^{+})},
\end{equation}
with $\chi_{k}(\tau)$ defined by equation \eqref{prcy46}. Via the formula for the $\tau$-dependent parameter 
$\ell_{2,k} \! := \! \ell_{2,k}^{+}$ defined by equations \eqref{prcy53} and \eqref{prcy54}, one simplifies equation 
\eqref{prcy50} to read
\begin{equation} \label{prcy55} 
\dfrac{\partial \Phi_{k}^{\ast}(\widetilde{\Lambda})}{\partial \widetilde{\Lambda}} \underset{\tau \to +\infty}{=} 
\left(\daleth_{k}(\tau,\widetilde{\Lambda}) \! + \! \mathcal{O} \big(\beth_{k}(\tau,\widetilde{\Lambda}) \big) 
\right) \! \Phi_{k}^{\ast}(\widetilde{\Lambda}), \quad k \! = \! \pm 1,
\end{equation}
where
\begin{equation} \label{prcy56} 
\daleth_{k}(\tau,\widetilde{\Lambda}) \! := \! \mu_{k}(\tau) \widetilde{\Lambda} \sigma_{3} \! + \! p_{k}(\tau) 
\sigma_{+} \! + \! q_{k}(\tau) \sigma_{-},
\end{equation}
with
\begin{gather}
p_{k}(\tau) \! := \! -\ell_{0,k}^{+} \! + \! \hat{\mathbb{L}}_{k}(\tau) \! + \! 1, \label{prcy57} \\
q_{k}(\tau) \! := \! \big(4 \sqrt{3} \mathcal{Z}_{k} \mathfrak{A}_{k} \chi_{k}^{-1}(\tau) \big)^{2}
(-\ell_{0,k}^{+} \! + \! 1) \! - \! \ell_{1,k}^{+} \! - \! \ell_{2,k}^{+} \! - \! \omega_{0,k}^{2}, \label{prcy58}
\end{gather}
and
\begin{align} \label{prcy59} 
\beth_{k}(\tau,\widetilde{\Lambda}) \! :=& \, 
\dfrac{\mi 4 \sqrt{3} \mathcal{Z}_{k} \mathfrak{A}_{k}}{\chi_{k}(\tau)}(-\ell_{0,k}^{+} \! + \! 1) \sigma_{3} \! - \! 
\hat{\mathbb{L}}_{k}(\tau) \sigma_{+} \! - \! \dfrac{\mi 8 \sqrt{3} \mathcal{Z}_{k} \mathfrak{A}_{k}}{\chi_{k}(\tau)} 
\ell_{2,k}^{+}(-\ell_{0,k}^{+} \! + \! 1) \widetilde{\Lambda} \sigma_{-} \nonumber \\
+& \, \big(\mathfrak{P}_{0,k}^{\ast} \! - \! \mi 8 \sqrt{3} \mathcal{Z}_{k} \mathfrak{A}_{k} \ell_{0,k}^{+} \big) 
\widetilde{\Lambda}^{2} \sigma_{3} \! + \!  \overset{\ast}{\mathbb{E}}_{k}^{\raise-6.75pt\hbox{$\scriptstyle \ast$}}
(\ell_{0,k}^{+},\ell_{1,k}^{+},\ell_{2,k}^{+};\tau) \widetilde{\Lambda}^{3},
\end{align}
where the yet-to-be-determined scalar function $\hat{\mathbb{L}}_{k}(\tau)$ is chosen in the proof of Lemma 
\ref{ginversion}.\footnote{It will be shown that $\hat{\mathbb{L}}_{k}(\tau) \genfrac{}{}{0pt}{3}{=}{\tau \to +\infty} 
\mathcal{O}(\tau^{-2/3})$, $k \! \in \! \lbrace \pm 1 \rbrace$: this fact will be used throughout the remainder of the proof.} 

\pmb{(xi)} The rationale for this---final---step is to transform equation \eqref{prcy55} into a `model' matrix linear ODE 
describing the coalescence of turning points. Let $\Phi_{M,k}(\widetilde{\Lambda})$, $k \! = \! \pm 1$, be a fundamental 
solution of equation \eqref{prcy1}; then, changing variables according to $\widetilde{\Lambda} \! = \! \widetilde{\Lambda}(z) 
\! = \! a_{k}^{\ast}(\tau)b^{\ast}z$, where $a_{k}^{\ast}(\tau) \! := \! (\mi 4 \sqrt{3}/ \mu_{k}(\tau))^{1/2}$ and $b^{\ast} \! := \! 
2^{-3/2}3^{-1/4} \me^{-\mi \pi/4}$, and defining $\phi_{M,k}(z) \! := \! \Phi_{M,k}(\widetilde{\Lambda}(z))$, one shows that 
$\phi_{M,k}(z)$ solves the canonical matrix ODE
\begin{equation} \label{prcy60} 
\partial_{z} \phi_{M,k}(z) \! = \! \left(\dfrac{z}{2} \sigma_{3} \! + \! P_{k}
(\tau) \sigma_{+} \! + \! Q_{k}(\tau) \sigma_{-} \right) \! \phi_{M,k}(z), 
\quad k \! = \! \pm 1,
\end{equation}
where $P_{k}(\tau) \! := \! a_{k}^{\ast}(\tau)b^{\ast}p_{k}(\tau)$ and $Q_{k}(\tau) \! := \! a_{k}^{\ast}(\tau)b^{\ast}
q_{k}(\tau)$, with fundamental solution expressed in terms of the parabolic-cylinder function, $D_{\star}
(\pmb{\pmb{\cdot}})$,\footnote{See, for example, \cite{a5,a18,a2}.}
\begin{equation} \label{prcy61} 
\phi_{M,k}(z) \! = \! 
\begin{pmatrix}
D_{-\nu (k)-1}(\mi z) & D_{\nu (k)}(z) \\
\frac{1}{P_{k}(\tau)}(\frac{\partial}{\partial z} \! - \! \frac{z}{2})D_{-\nu (k)-1}(\mi z) & 
\frac{1}{P_{k}(\tau)}(\frac{\partial}{\partial z} \! - \! \frac{z}{2})D_{\nu (k)}(z)
\end{pmatrix},
\end{equation}
where $-(\nu (k) \! + \! 1) \! := \! P_{k}(\tau)Q_{k}(\tau)$. Inverting the dependent- and independent-variable linear 
transformations given above, one arrives at the formula for the parameter $\nu (k) \! + \! 1$ defined by equation 
\eqref{prpr1} and the representation for $\Phi_{M,k}(\widetilde{\Lambda})$ given in equation \eqref{prcy2}.

Finally, in order to establish the asymptotic representation \eqref{prpr2}, one has to estimate the unimodular 
function $\hat{\chi}_{k}(\widetilde{\Lambda})$ defined in the transformation $\mathfrak{F}_{11}$. Under the 
action of the transformation $\mathfrak{F}_{11}$, one rewrites equation \eqref{prcy55} as follows:
\begin{equation} \label{prcy62} 
\dfrac{\partial \hat{\chi}_{k}(\widetilde{\Lambda})}{\partial \widetilde{\Lambda}} \underset{\tau \to +\infty}{=} 
\beth_{k}(\tau,\widetilde{\Lambda}) \hat{\chi}_{k}(\widetilde{\Lambda}) \! + \! \left[\daleth_{k}
(\tau,\widetilde{\Lambda}),\hat{\chi}_{k}(\widetilde{\Lambda}) \right], \quad k \! = \! \pm 1,
\end{equation}
where $\daleth_{k}(\tau,\widetilde{\Lambda})$ is defined by equations \eqref{prcy56}--\eqref{prcy58}, and 
$\beth_{k}(\tau,\widetilde{\Lambda})$ is defined by equation \eqref{prcy59}. The normalised solution of 
equation \eqref{prcy62}, that is, the one for which $\hat{\chi}_{k}(0) \! = \! \mathrm{I}$, is given by
\begin{equation} \label{prcy63} 
\hat{\chi}_{k}(\widetilde{\Lambda}) \! = \! \mathrm{I} \! + \! \int_{0}^{\widetilde{\Lambda}} \Phi_{M,k}
(\widetilde{\Lambda}) \Phi_{M,k}^{-1}(\xi) \beth_{k}(\tau,\xi) \hat{\chi}_{k}(\xi) \Phi_{M,k}(\xi) \Phi_{M,k}^{-1}
(\widetilde{\Lambda}) \, \md \xi, \quad k \! = \! \pm 1.
\end{equation}
In order to prove the required estimate for $\hat{\chi}_{k}(\widetilde{\Lambda})$, one uses the method of 
successive approximations, namely,
\begin{equation*}
\hat{\chi}_{k}^{(m)}(\widetilde{\Lambda}) \! = \! \mathrm{I} \! + \! \int_{0}^{\widetilde{\Lambda}} \Phi_{M,k}
(\widetilde{\Lambda}) \Phi_{M,k}^{-1}(\xi) \beth_{k}(\tau,\xi) \hat{\chi}_{k}^{(m-1)}(\xi) \Phi_{M,k}(\xi) 
\Phi_{M,k}^{-1}(\widetilde{\Lambda}) \, \md \xi, \quad k \! = \! \pm 1, \quad m \! \in \! 
\mathbb{N},
\end{equation*}
with $\hat{\chi}^{(0)}_{k}(\widetilde{\Lambda}) \! \equiv \! \mathrm{I}$, to construct a Neumann series solution for 
$\hat{\chi}_{k}(\widetilde{\Lambda})$ ($\hat{\chi}_{k}(\widetilde{\Lambda}) \! := \! \lim_{m \to \infty} \hat{\chi}^{(m)}_{k}
(\widetilde{\Lambda})$); in this instance, however, it suffices to estimate the matrix norm of the associated resolvent 
kernel. Via the above iteration argument, a calculation shows that, for $k \! = \! \pm 1$,
\begin{equation} \label{prcy64} 
\lvert \lvert \hat{\chi}_{k}(\widetilde{\Lambda}) \! - \! \mathrm{I} \rvert \rvert \underset{\tau \to +\infty}{\leqslant} 
\exp \! \left(\int_{0}^{\widetilde{\Lambda}} \lvert \lvert \Phi_{M,k}(\widetilde{\Lambda}) \rvert \rvert \lvert \lvert 
\Phi_{M,k}^{-1}(\xi) \rvert \rvert \lvert \lvert \beth_{k}(\tau,\xi) \rvert \rvert \lvert \lvert \Phi_{M,k}(\xi) \rvert \rvert 
\lvert \lvert \Phi_{M,k}^{-1}(\widetilde{\Lambda}) \rvert \rvert \lvert \md \xi \rvert \right) \! - \! 1,
\end{equation}
where $\lvert \md \xi \rvert$ denotes integration with respect to arc length. Noting that (see Remark \ref{leedasyue})  
$\det (\Phi_{M,k}(z)) \! = \! -\me^{-\mi \pi (\nu (k)+1)/2}(2 \mu_{k}(\tau))^{1/2}p^{-1}_{k}(\tau)$, it follows {}from the 
estimate \eqref{prcy64} that, for $k \! = \! \pm 1$,
\begin{equation} \label{prcy65} 
\lvert \lvert \hat{\chi}_{k}(\widetilde{\Lambda}) \! - \! \mathrm{I} \rvert \rvert \underset{\tau \to +\infty}{\leqslant} 
\exp \! \left(\dfrac{\lvert p_{k}(\tau) \rvert^{2} \lvert \lvert \Phi_{M,k}(\widetilde{\Lambda}) \rvert \rvert^{2}}{
\lvert 2 \mu_{k}(\tau) \rvert (\me^{\pi \Im (\nu (k)+1)/2})^{2}} \int_{0}^{\widetilde{\Lambda}} \lvert \lvert \Phi_{M,k}
(\xi) \rvert \rvert^{2} \lvert \lvert \beth_{k}(\tau,\xi) \rvert \rvert \, \lvert \md \xi \rvert \right) \! - \! 1.
\end{equation}
One now proceeds to estimate the respective norms in equation \eqref{prcy65}.

One commences with the estimation of the norm $\lvert \lvert \beth_{k}(\tau,\xi) \rvert \rvert$, $k \! = \! \pm 1$, appearing 
in equation \eqref{prcy65}. Via equations \eqref{prcy7}, \eqref{prcy11}, \eqref{prcy24}, \eqref{prcy28}, \eqref{prcy32}, 
\eqref{prcy36}, \eqref{prcy37}, \eqref{prcy38}, \eqref{prcy41}, \eqref{prcy42}, \eqref{prcy43}, \eqref{prcy47}, 
\eqref{prcy49}, \eqref{prcy50}, \eqref{prcy51}, and \eqref{prcy59}, one shows that, for $k \! = \! \pm 1$, in terms of the 
composition of the linear transformations $\mathfrak{F}_{j}$, $j \! = \! 1,2,\dotsc,11$,
\begin{align} \label{prcy81} 
\beth_{k}(\tau,\widetilde{\Lambda}) \! :=& \, \left(\mathfrak{F}_{11} \circ \mathfrak{F}_{10} \circ \mathfrak{F}_{9} 
\circ \mathfrak{F}_{8} \circ \mathfrak{F}_{7} \circ \mathfrak{F}_{6} \circ \mathfrak{F}_{5} \circ \mathfrak{F}_{4} 
\circ \mathfrak{F}_{3} \circ \mathfrak{F}_{2} \circ \mathfrak{F}_{1} \right) \! (\widetilde{\Psi}(\widetilde{\mu},\tau) 
\! - \! \widetilde{\Psi}_{k}(\widetilde{\mu},\tau)) \nonumber \\
=& \, \dfrac{\mi 4 \sqrt{3} \mathcal{Z}_{k} \mathfrak{A}_{k}}{\chi_{k}(\tau)}(-\ell_{0,k}^{+} \! + \! 1) \sigma_{3} \! - 
\! \hat{\mathbb{L}}_{k}(\tau) \sigma_{+} \! - \! \dfrac{\mi 8 \sqrt{3} \mathcal{Z}_{k} \mathfrak{A}_{k}}{\chi_{k}(\tau)} 
\ell_{2,k}^{+}(-\ell_{0,k}^{+} \! + \! 1) \widetilde{\Lambda} \sigma_{-} \nonumber \\
+& \, (\mathfrak{P}_{0,k}^{\ast} \! - \! \mi 8 \sqrt{3} \mathcal{Z}_{k} \mathfrak{A}_{k} \ell_{0,k}^{+}) 
\widetilde{\Lambda}^{2} \sigma_{3} \! + \! \left(-2 \ell_{2,k}^{+}(\mathfrak{P}_{0,k}^{\ast} \! - \! \mi 8 \sqrt{3} 
\mathcal{Z}_{k} \mathfrak{A}_{k} \ell_{0,k}^{+}) \sigma_{-} \right. \nonumber \\
+&\left. \mathcal{G}_{2,k}^{-1} 
\left(
\begin{pmatrix}
1 & 0 \\
\mi \omega_{0,k} & 1
\end{pmatrix} \tau^{-\frac{1}{6} \sigma_{3}} \mathcal{G}_{1,k}^{-1} \mathcal{G}_{0,k}^{-1} \hat{\mathbb{E}}_{k}
(\tau) \mathcal{G}_{0,k} \mathcal{G}_{1,k} \tau^{\frac{1}{6} \sigma_{3}} 
\begin{pmatrix}
1 & 0 \\
-\mi \omega_{0,k} & 1
\end{pmatrix} \right. \right. 
\nonumber \\
+&\left. \left. 
\begin{pmatrix}
-\mathfrak{R}_{0,k}^{\ast} \ell_{0,k}^{+} & \ell_{0,k}^{+}(2 \mathfrak{P}_{0,k}^{\ast} \! - \! \mi 8 \sqrt{3} \mathcal{Z}_{k} 
\mathfrak{A}_{k} \ell_{0,k}^{+}) \\
-2 \ell_{1,k}^{+}(\mathfrak{P}_{0,k}^{\ast} \! - \! \mi 8 \sqrt{3} \mathcal{Z}_{k} \mathfrak{A}_{k} \ell_{0,k}^{+}) & 
\mathfrak{R}_{0,k}^{\ast} \ell_{0,k}^{+} 
\end{pmatrix} \right) \! \mathcal{G}_{2,k} \right) \! \widetilde{\Lambda}^{3},
\end{align}
whence, via the definitions \eqref{prcy9}, \eqref{prcy23}, \eqref{prcy33}--\eqref{prcy35}, \eqref{prcy40}, and 
\eqref{prcy48}, and a matrix-multipli\-c\-a\-t\-i\-o\-n argument, one arrives at, for $k \! = \! \pm 1$,
\begin{align} \label{prcy82} 
\beth_{k}(\tau,\widetilde{\Lambda}) \! =& \, \dfrac{\mi 4 \sqrt{3} \mathcal{Z}_{k} 
\mathfrak{A}_{k}}{\chi_{k}(\tau)}(-\ell_{0,k}^{+} \! + \! 1) \sigma_{3} \! - \! 
\hat{\mathbb{L}}_{k}(\tau) \sigma_{+} \! - \! \dfrac{\mi 8 \sqrt{3} \mathcal{Z}_{k} 
\mathfrak{A}_{k}}{\chi_{k}(\tau)} \ell_{2,k}^{+}(-\ell_{0,k}^{+} \! + \! 1) 
\widetilde{\Lambda} \sigma_{-} \nonumber \\
+& \, (\mathfrak{A}_{0,k}^{\sharp} \! + \! \mathfrak{A}_{k} \omega_{0,k}^{2}
(\ell_{0,k}^{+})^{2}) \widetilde{\Lambda}^{2} \sigma_{3} \! + \! 
\begin{pmatrix}
\mathcal{N}_{11}^{\ast}(\tau) \! + \! \mathcal{M}_{11}^{\ast}(\tau) & 
\mathcal{N}_{12}^{\ast}(\tau) \! + \! \mathcal{M}_{12}^{\ast}(\tau) \\
\mathcal{N}_{21}^{\ast}(\tau) \! + \! \mathcal{M}_{21}^{\ast}(\tau) & 
-(\mathcal{N}_{11}^{\ast}(\tau) \! + \! \mathcal{M}_{11}^{\ast}(\tau))
\end{pmatrix} \! \widetilde{\Lambda}^{3},
\end{align}
where
{\fontsize{10pt}{11pt}\selectfont
\begin{align}
\mathcal{N}_{11}^{\ast}(\tau) \! :=& \, \ell_{0,k}^{+} \mathfrak{A}_{k} \! 
\left(\dfrac{\mathfrak{A}_{k} \mathfrak{B}_{0,k}^{\sharp}}{\mathfrak{B}_{k}} 
\! - \! 2 \mathfrak{A}_{0,k}^{\sharp} \right) \! \left(1 \! - \! \dfrac{\mi 
4 \sqrt{3} \mathcal{Z}_{k}}{\chi_{k}(\tau)} \right) \! - \! \ell_{0,k}^{+} \! 
\left(\mathfrak{B}_{k} \mathfrak{C}_{0,k}^{\sharp} \! - \! \dfrac{\mi 
4 \sqrt{3} \mathcal{Z}_{k}}{\chi_{k}(\tau)} \mathfrak{A}_{k}^{2} 
\omega_{0,k}^{2}(\ell_{0,k}^{+})^{2} \right), \label{prcy83} \\
\mathcal{N}_{12}^{\ast}(\tau) \! :=& \, \ell_{0,k}^{+} \! \left(2 
\mathfrak{A}_{0,k}^{\sharp} \! + \! \mathfrak{A}_{k} \omega_{0,k}^{2}
(\ell_{0,k}^{+})^{2} \! - \! \dfrac{\mathfrak{A}_{k} \mathfrak{B}_{0,k}^{
\sharp}}{\mathfrak{B}_{k}} \right), \label{prcy84} \\
\mathcal{N}_{21}^{\ast}(\tau) \! :=& \, -\dfrac{\mi 4 \sqrt{3} \mathcal{Z}_{k} 
\mathfrak{A}_{k}}{\chi_{k}(\tau)} \! \left(\ell_{0,k}^{+} \mathfrak{A}_{k} \! 
\left(\dfrac{\mathfrak{A}_{k} \mathfrak{B}_{0,k}^{\sharp}}{\mathfrak{B}_{k}} 
\! - \! 2 \mathfrak{A}_{0,k}^{\sharp} \right) \! \left(2 \! - \! \dfrac{\mi 4 
\sqrt{3} \mathcal{Z}_{k}}{\chi_{k}(\tau)} \right) \! - \! \ell_{0,k}^{+} \! 
\left(2 \mathfrak{B}_{k} \mathfrak{C}_{0,k}^{\sharp} \! - \! \dfrac{
\mi 4 \sqrt{3} \mathcal{Z}_{k}}{\chi_{k}(\tau)} \right. \right. \nonumber \\
\times&\left. \left. \, \mathfrak{A}_{k}^{2} \omega_{0,k}^{2}(\ell_{0,k}^{+})^{2} 
\right) \right) \! - \! 2(\mathfrak{A}_{0,k}^{\sharp} \! + \! \mathfrak{A}_{k} 
\omega_{0,k}^{2}(\ell_{0,k}^{+})^{2})(\ell_{1,k}^{+} \! + \! \ell_{2,k}^{+}), 
\label{prcy85} \\
\mathcal{M}_{11}^{\ast}(\tau) \! :=& \, \dfrac{\hat{\mathcal{C}}_{1}}{2 
\lambda_{1}^{\ast}(k) \mathfrak{B}_{k}} \! \left((\hat{\mathbb{E}}_{k}(\tau))_{11} 
\! \left(\hat{\mathfrak{g}}_{11} \mathfrak{B}_{k} \! + \! \hat{\mathfrak{g}}_{12} 
\mathfrak{A}_{k} \! \left(\dfrac{\mi 4 \sqrt{3} \mathcal{Z}_{k}}{\chi_{k}(\tau)} 
\! - \! 1 \right) \! + \! \hat{\mathfrak{g}}_{12} \! \left(\mathfrak{B}_{k} \! + \! 
\mathfrak{A}_{k} \! \left(\dfrac{\mi 4 \sqrt{3} \mathcal{Z}_{k}}{\chi_{k}(\tau)} \! 
- \! 1 \right) \right) \right) \right. \nonumber \\
+&\left. \, (\hat{\mathbb{E}}_{k}(\tau))_{12} \! \left(\mathfrak{B}_{k} \! + \! 
\mathfrak{A}_{k} \! \left(\dfrac{\mi 4 \sqrt{3} \mathcal{Z}_{k}}{\chi_{k}(\tau)} 
\! - \! 1 \right) \right) \! - \! (\hat{\mathbb{E}}_{k}(\tau))_{21} 
\hat{\mathfrak{g}}_{12} \! \left(\hat{\mathfrak{g}}_{11} \mathfrak{B}_{k} \! 
+ \! \hat{\mathfrak{g}}_{12} \mathfrak{A}_{k} \! \left(\dfrac{\mi 4 \sqrt{3} 
\mathcal{Z}_{k}}{\chi_{k}(\tau)} \! - \! 1 \right) \right) \right), \label{prcy86} \\
\mathcal{M}_{12}^{\ast}(\tau) \! :=& \, \dfrac{\hat{\mathcal{C}}_{1}}{2 
\lambda_{1}^{\ast}(k) \mathfrak{B}_{k}} \left(2(\hat{\mathbb{E}}_{k}(\tau))_{11} 
\hat{\mathfrak{g}}_{12} \! + \! (\hat{\mathbb{E}}_{k}(\tau))_{12} \! - \! 
(\hat{\mathbb{E}}_{k}(\tau))_{21}(\hat{\mathfrak{g}}_{12})^{2} \right), 
\label{prcy87} \\
\mathcal{M}_{21}^{\ast}(\tau) \! :=& \, \dfrac{\hat{\mathcal{C}}_{1}}{2 
\lambda_{1}^{\ast}(k) \mathfrak{B}_{k}} \! \left(-2(\hat{\mathbb{E}}_{k}
(\tau))_{11} \! \left(\mathfrak{B}_{k} \! + \! \mathfrak{A}_{k} \! \left(
\dfrac{\mi 4 \sqrt{3} \mathcal{Z}_{k}}{\chi_{k}(\tau)} \! - \! 1 \right) \right) \! 
\left(\hat{\mathfrak{g}}_{11} \mathfrak{B}_{k} \! + \! \hat{\mathfrak{g}}_{12} 
\mathfrak{A}_{k} \! \left(\dfrac{\mi 4 \sqrt{3} \mathcal{Z}_{k}}{\chi_{k}(\tau)} 
\! - \! 1 \right) \right) \right. \nonumber \\
-&\left. \, (\hat{\mathbb{E}}_{k}(\tau))_{12} \! \left(\mathfrak{B}_{k} \! + \! 
\mathfrak{A}_{k} \! \left(\dfrac{\mi 4 \sqrt{3} \mathcal{Z}_{k}}{\chi_{k}(\tau)} 
\! - \! 1 \right) \right)^{2} \! + \! (\hat{\mathbb{E}}_{k}(\tau))_{21} \! \left(
\hat{\mathfrak{g}}_{11} \mathfrak{B}_{k} \! + \! \hat{\mathfrak{g}}_{12} 
\mathfrak{A}_{k} \! \left(\dfrac{\mi 4 \sqrt{3} \mathcal{Z}_{k}}{\chi_{k}(\tau)} 
\! - \! 1 \right) \right)^{2} \right), \label{prcy88} 
\end{align}}
with
\begin{equation} \label{prcy89} 
\hat{\mathfrak{g}}_{11} \! := \! \dfrac{\hat{\mathcal{A}}_{1} \! + \! 
\lambda_{1}^{\ast}(k)}{\hat{\mathcal{C}}_{1}} \, \quad \, \text{and} \, 
\quad \, \hat{\mathfrak{g}}_{12} \! := \! \dfrac{\hat{\mathcal{A}}_{1} 
\! - \! \lambda_{1}^{\ast}(k)}{\hat{\mathcal{C}}_{1}}.
\end{equation}
In order to realise the $\tau \! \to \! +\infty$ asymptotics for $\beth_{k}(\tau,\widetilde{\Lambda})$, $k \! = \! \pm 1$, 
via equation \eqref{prcy82} (cf. definitions \eqref{prcy83}--\eqref{prcy89}), and subsequently estimate the norm 
$\lvert \lvert \beth_{k}(\tau,\xi) \rvert \rvert$, the $\tau \! \to \! +\infty$ asymptotics of the functions $\mathcal{Z}_{k}$, 
$\mathcal{G}_{0,k}$, $\mathfrak{A}_{k}$, $\mathfrak{B}_{k}$, $\mathfrak{C}_{k}$, $\mathfrak{A}_{0,k}^{\sharp}$, 
$\mathfrak{B}_{0,k}^{\sharp}$, $\mathfrak{C}_{0,k}^{\sharp}$, $\omega_{0,k}^{2}$, $\chi_{k}(\tau)$, and $\mu_{k}
(\tau)$, and the $\tau$-dependent parameters $\ell_{0,k}^{+}$, $\ell_{1,k}^{+}$, and $\ell_{2,k}^{+}$ are required: for 
the reader's convenience, they are presented in Appendix \ref{formforleadasymps}. Substituting the asymptotics 
\eqref{tr1}, \eqref{tr3}, and \eqref{prcyzeek1}--\eqref{prcyell2k1} into the definitions \eqref{prcy83}--\eqref{prcy89}, 
recalling the definitions \eqref{prcy33} and \eqref{prcy35}, and using equation \eqref{prcy82}, one arrives at the estimate
\begin{equation} \label{estbet1} 
\lvert \lvert \beth_{k}(\tau,\xi) \rvert \rvert \underset{\tau \to +\infty}{=} \mathcal{O} \big(\tau^{-1/3} 
\lvert \xi \rvert^{3} \big), \quad k \! = \! \pm 1.
\end{equation}

There remains the matter of estimating the norm of the unimodular function $\hat{\chi}_{k}(\xi)$, $k \! = \! \pm 1$. In 
order to do so, one has to derive a uniform approximation for $\hat{\chi}_{k}(\xi)$ on $\mathbb{R} \cup \mi \mathbb{R}$. 
Towards this goal, one uses the following integral representation for the parabolic-cylinder function (see, for example, 
\cite{EMOT}); for $k \! = \! \pm 1$,
\begin{equation} \label{prcy66} 
D_{\nu (k)}(z) \! = \! \dfrac{2^{\frac{\nu (k)}{2}} \me^{-\frac{z^{2}}{4}}}{\Gamma (-\nu (k)/2)} \int_{0}^{+\infty} 
\me^{-\frac{\xi z^{2}}{2}} \xi^{-\frac{\nu (k)}{2}-1}(1 \! + \! \xi)^{\frac{\nu (k)-1}{2}} \, \md \xi, \quad \Re 
(\nu (k)) \! < \! 0, \quad \lvert \arg (z) \rvert \! \leqslant \! \pi/4,
\end{equation}
where $\Gamma (\pmb{\cdot})$ is the gamma function. As the integral representation \eqref{prcy66} will be applied 
simultaneously to the entries of the $\mathrm{M}_{2}(\mathbb{C})$-valued function (cf. equation \eqref{prcy2}) 
$\Phi_{M,k}(\xi)$ in order to arrive at a uniform approximation for $\hat{\chi}_{k}(\xi)$ on the Stokes rays $\arg (\xi) 
\! = \! 0,\pm \pi/2,\pm \pi,\dotsc$, $0 \! \leqslant \! \lvert \xi \rvert \! < \! +\infty$, it implies the restrictions \eqref{pc4} 
on $\nu (k) \! + \! 1$; in fact, for the purposes of this proof, it suffices to have a uniform approximation for $\hat{\chi}_{k}
(\xi)$ on, say, the Stokes rays $\pmb{\widehat{\mathcal{S}}} \! := \! \lbrace \mathstrut \xi \! \in \! \mathbb{C}; \, 0 \! 
\leqslant \! \lvert \xi \rvert \! < \! +\infty, \, \arg (\xi) \! = \! 0,-\pi/2,-\pi,-3 \pi/2 \rbrace$. Using the following functional 
relations and values for the gamma function (see, for example, \cite{a24}),
\begin{gather*}
\Gamma (z \! + \! 1) \! = \! z \Gamma (z), \, \quad \, \Gamma (z) \Gamma (1 \! - \! z) \! = \! \dfrac{\pi}{\sin (\pi z)}, 
\, \quad \, \sqrt{\pi} \, \Gamma (2z) \! = \! 2^{2z-1} \Gamma (z) \Gamma (z \! + \! 1/2), \\
\Gamma (1/2) \! = \! \sqrt{\pi}, \, \quad \, \int_{0}^{+\infty} \dfrac{t^{x-1}}{(1 \! + \! t)^{x+y}} \, \md t \! = \! 
\dfrac{\Gamma (x) \Gamma (y)}{\Gamma (x \! + \! y)}, \quad \Re (x),\Re (y) \! > \! 0,
\end{gather*}
the linear relations relating any three of the four parabolic-cylinder functions (cf. equation \eqref{prcy61}) 
$D_{-\nu (k)-1}(\pm \mi z)$ and $D_{\nu (k)}(\pm z)$,
\begin{align*}
\sqrt{2 \pi}D_{\nu (k)}(z) =& \, \Gamma (\nu (k) \! + \! 1) \! \left(\me^{\mi \pi \nu (k)/2}D_{-\nu (k)-1}(\mi z) 
\! + \! \me^{-\mi \pi \nu (k)/2}D_{-\nu (k)-1}(-\mi z) \right), \\
D_{\nu (k)}(z) =& \, \me^{-\mi \pi \nu (k)}D_{\nu (k)}(-z) \! + \! \dfrac{\sqrt{2 \pi} \me^{-\mi \pi (\nu (k)+1)/2}}{
\Gamma (-\nu (k))}D_{-\nu (k)-1}(\mi z), \\
D_{\nu (k)}(z) =& \, \me^{\mi \pi \nu (k)}D_{\nu (k)}(-z) \! + \! \dfrac{\sqrt{2 \pi} \me^{\mi \pi (\nu (k)+1)/2}}{
\Gamma (-\nu (k))}D_{-\nu (k)-1}(-\mi z),
\end{align*}
and the fact that (see the asymptotics \eqref{geek13} below) $\nu (k) \! + \! 1 \! \to \! 0$ as $\tau \! \to \! +\infty$, 
one arrives at, via the restrictions \eqref{pc4} on $\nu (k) \! + \! 1$, equation \eqref{prcy2}, and the integral 
representation \eqref{prcy66}, estimates for the moduli $\lvert (\Phi_{M,k}(\xi))_{ij} \rvert$, $k \! = \! \pm 1$, 
$i,j \! = \! 1,2$, on the Stokes rays $\pmb{\widehat{\mathcal{S}}}$: for the convenience of the reader, they 
are stated in Appendix \ref{appenformoduli}. To eschew technical redundancies, consider, say, the case 
$k \! = \! +1$, and, without loss of generality, $\arg (\widetilde{\Lambda}) \! = \! \pm \pi/2$:\footnote{The pair of 
values $\arg (\widetilde{\Lambda}) \! = \! \pm \pi/2$ on the Stokes rays are chosen for illustrative purposes only, 
in order to present the general scheme of the calculations: for any of the remaining $\binom{4}{2} \! - \! 1 \! = 
\! 5$ pairs of values of $\arg (\widetilde{\Lambda})$ on the Stokes rays, one arrives at the same estimate (see 
equation \eqref{prcy100}) for $\lvert \lvert \hat{\chi}_{k}(\widetilde{\Lambda}) \! - \! \mathrm{I} \rvert \rvert$, 
$k \! = \! \pm 1$, but with different $\mathcal{O}(1)$ constants.} the case $k \! = \! -1$ is analogous. Using the 
asymptotic expansions for the parabolic-cylinder functions (see Remark \ref{leedasyue}), one shows that: 
\textbf{(a)} for $\arg (\widetilde{\Lambda}) \genfrac{}{}{0pt}{3}{=}{\tau \to +\infty} \pi/2 \! + \! \mathcal{O} 
\big(\tau^{-2/3} \big)$,
\begin{gather} \label{prcy75} 
\begin{split}
\lvert (\Phi_{M,1}(\widetilde{\Lambda}))_{11} \rvert \underset{\tau \to +\infty}{=}& \mathcal{O} \! \left(\tilde{\rho}_{0} 
\lvert \widetilde{\Lambda} \rvert^{-\Re (\nu (1)+1)} \right), \quad \lvert (\Phi_{M,1}(\widetilde{\Lambda}))_{12} \rvert 
\underset{\tau \to +\infty}{=} \mathcal{O} \! \left(\tilde{\rho}_{1} \lvert \widetilde{\Lambda} \rvert^{-\Re (\nu (1)+1)} \right), \\
\lvert (\Phi_{M,1}(\widetilde{\Lambda}))_{21} \rvert \underset{\tau \to +\infty}{=}& \mathcal{O} \! \left(\tilde{\rho}_{2} 
\dfrac{\lvert \widetilde{\Lambda} \rvert^{\Re (\nu (1)+1)}}{\lvert p_{1}(\tau) \rvert} \right), \quad \lvert (\Phi_{M,1}
(\widetilde{\Lambda}))_{22} \rvert \underset{\tau \to +\infty}{=} \mathcal{O} \! \left(\tilde{\rho}_{3} \dfrac{\lvert 
\widetilde{\Lambda} \rvert^{\Re (\nu (1)+1)}}{\lvert p_{1}(\tau) \rvert} \right),
\end{split}
\end{gather}
where
\begin{gather*}
\tilde{\rho}_{0} \! = \! \eta_{+} \me^{-3 \pi \Im (\nu (1)+1)/2}, \quad \quad \tilde{\rho}_{3} \! = \! 2^{3/2}3^{1/4}/\eta_{+}, \\
\tilde{\rho}_{1} \! = \! \dfrac{2^{3/2} \eta_{+}}{\sqrt{\pi}} \me^{-\pi \Im (\nu (1)+1)} \lvert \cos (\tfrac{\pi}{2}(\nu (1) \! + \! 1)) 
\rvert \lvert \sin (\tfrac{\pi}{2}(\nu (1) \! + \! 1)) \rvert \Gamma (\Re (\nu (1) \! + \! 1)), \\
\tilde{\rho}_{2} \! = \! \dfrac{8 \cdot 3^{1/4}}{\sqrt{\pi} \eta_{+}} \me^{\pi \Im (\nu (1)+1)/2} \lvert \cos (\tfrac{\pi}{2}
(\nu (1) \! + \! 1)) \rvert \lvert \sin (\tfrac{\pi}{2}(\nu (1) \! + \! 1)) \rvert \Gamma (-\Re (\nu (1))),
\end{gather*}
with $\eta_{+} \! := \! (2^{3/2}3^{1/4})^{-\Re (\nu (1)+1)} \me^{3 \pi \Im (\nu (1)+1)/4}$; and \textbf{(b)} for $\arg 
(\widetilde{\Lambda}) \genfrac{}{}{0pt}{3}{=}{\tau \to +\infty} -\pi/2 \! + \! \mathcal{O} \big(\tau^{-2/3} \big)$,
\begin{gather} \label{prcy76} 
\begin{split}
\lvert (\Phi_{M,1}(\widetilde{\Lambda}))_{11} \rvert \underset{\tau \to +\infty}{=}& \mathcal{O} \! \left(\hat{\rho}_{0} 
\lvert \widetilde{\Lambda} \rvert^{-\Re (\nu (1)+1)} \right), \quad \lvert (\Phi_{M,1}(\widetilde{\Lambda}))_{12} \rvert 
\underset{\tau \to +\infty}{=} \mathcal{O} \! \left(\hat{\rho}_{1} \dfrac{\lvert \widetilde{\Lambda} \rvert^{\Re (\nu (1)
+1)}}{\lvert \widetilde{\Lambda} \rvert} \right), \\
\lvert (\Phi_{M,1}(\widetilde{\Lambda}))_{21} \rvert \underset{\tau \to +\infty}{=}& \mathcal{O} \! \left(\hat{\rho}_{2} 
\dfrac{\lvert \widetilde{\Lambda} \rvert^{-\Re (\nu (1)+1)}}{\lvert p_{1}(\tau) \rvert \lvert \widetilde{\Lambda} \rvert} 
\right), \quad \lvert (\Phi_{M,1}(\widetilde{\Lambda}))_{22} \rvert \underset{\tau \to +\infty}{=} \mathcal{O} \! \left(
\hat{\rho}_{3} \dfrac{\lvert \widetilde{\Lambda} \rvert^{\Re (\nu (1)+1)}}{\lvert p_{1}(\tau) \rvert} \right),
\end{split}
\end{gather}
where
\begin{equation*}
\hat{\rho}_{0} \! = \! \eta_{-} \me^{\pi \Im (\nu (1)+1)/2}, \! \! \quad \! \! \hat{\rho}_{1} \! = \! 2^{-3/2}3^{-1/4}/\eta_{-}, 
\! \! \quad \! \! \hat{\rho}_{2} \! = \! \eta_{-} \me^{\pi \Im (\nu (1)+1)/2} \lvert \nu (1) \! + \! 1 \rvert, \! \! \quad \! \! 
\hat{\rho}_{3} \! = \! 2^{3/2}3^{1/4}/\eta_{-},
\end{equation*}
with $\eta_{-} \! := \! (2^{3/2}3^{1/4})^{-\Re (\nu (1)+1)} \me^{-\pi \Im (\nu (1)+1)/4}$. Hence, via the elementary 
inequalities $\lvert \Re (\nu (1) \! + \! 1) \rvert \! \leqslant \! \lvert \nu (1) \! + \! 1 \rvert$ and $\lvert \Im (\nu (1) 
\! + \! 1) \rvert \! \leqslant \! \lvert \nu (1) \! + \! 1 \rvert$, it follows {}from the estimates \eqref{prcy75} and 
\eqref{prcy71}--\eqref{prcy74} that, for $\arg (\widetilde{\Lambda}) \genfrac{}{}{0pt}{3}{=}{\tau \to +\infty} \pi/2 
\! + \! \mathcal{O} \big(\tau^{-2/3} \big)$,
\begin{gather}
\lvert \lvert \Phi_{M,1}(\xi) \rvert \rvert^{2} \underset{\tau \to +\infty}{=} \mathcal{O}(\tilde{\mathfrak{c}}_{M}^{\sharp}) 
\! + \! \mathcal{O} \! \left(\dfrac{\tilde{\mathfrak{c}}_{M}^{\sharp} \lvert \nu (1) \! + \! 1 \rvert^{2} \lvert \xi \rvert^{2}}{\lvert 
p_{1}(\tau) \rvert^{2}} \right), \label{prcy77} \\
\lvert \lvert \Phi_{M,1}(\widetilde{\Lambda}) \rvert \rvert^{2} \underset{\tau \to +\infty}{=} \mathcal{O} \! \left(\lvert 
\widetilde{\Lambda} \rvert^{2 \Re (\nu (1)+1)} \! \left(\dfrac{\tilde{\mathfrak{c}}_{M}}{\lvert p_{1}(\tau) \rvert^{2}} 
\! + \! \mathcal{O} \! \left(\dfrac{\tilde{\mathfrak{c}}_{M}}{\lvert \widetilde{\Lambda} \rvert^{4 \Re (\nu (1)+1)}} 
\right) \right) \right), 
\label{prcy78}
\end{gather}
where $\tilde{\mathfrak{c}}_{M}^{\sharp} \! := \! 2 \max_{m=0,1,2,3} \lbrace (\tilde{\varrho}_{m}(1))^{2} \rbrace$, 
and $\tilde{\mathfrak{c}}_{M} \! := \! 2 \max_{m=0,1,2,3} \lbrace \tilde{\rho}_{m}^{2} \rbrace$, and, {}from the 
estimates \eqref{prcy76} and \eqref{prcy67}--\eqref{prcy70}, it follows that, for $\arg (\widetilde{\Lambda}) 
\genfrac{}{}{0pt}{3}{=}{\tau \to +\infty} -\pi/2 \! + \! \mathcal{O} \big(\tau^{-2/3} \big)$,
\begin{gather}
\lvert \lvert \Phi_{M,1}(\xi) \rvert \rvert^{2} \underset{\tau \to +\infty}{=} \mathcal{O}(\hat{\mathfrak{c}}_{M}^{\sharp}) 
\! + \! \mathcal{O} \! \left(\dfrac{\hat{\mathfrak{c}}_{M}^{\sharp} \lvert \nu (1) \! + \! 1 \rvert^{2} \lvert \xi \rvert^{2}}{\lvert 
p_{1}(\tau) \rvert^{2}} \right), \label{prcy79} \\
\lvert \lvert \Phi_{M,1}(\widetilde{\Lambda}) \rvert \rvert^{2} \underset{\tau \to +\infty}{=} \mathcal{O} \! \left(\lvert 
\widetilde{\Lambda} \rvert^{2 \Re (\nu (1)+1)} \! \left(\dfrac{\hat{\mathfrak{c}}_{M}}{\lvert p_{1}(\tau) \rvert^{2}} \! + \! 
\mathcal{O} \! \left(\dfrac{\hat{\mathfrak{c}}_{M}}{\lvert \widetilde{\Lambda} \rvert^{2 \min \lbrace 1,2 \Re (\nu (1)+1) 
\rbrace}} \right) \right) \right), \label{prcy80}
\end{gather}
where $\hat{\mathfrak{c}}_{M}^{\sharp} \! := \! 2 \max_{m=0,1,2,3} \lbrace (\hat{\varrho}_{m}(1))^{2} \rbrace$, and 
$\hat{\mathfrak{c}}_{M} \! := \! \max_{m=0,1,2,3} \lbrace \hat{\rho}_{m}^{2} \rbrace$. Assembling the asymptotics 
\eqref{prcy77}--\eqref{prcy80} and invoking the restriction \eqref{pc4} on $\delta_{k}$ (for $k \! = \! +1)$, one deduces 
{}from the asymptotics \eqref{prcy65} and \eqref{estbet1} that, for $\arg (\widetilde{\Lambda}) 
\genfrac{}{}{0pt}{3}{=}{\tau \to +\infty} \pm \pi/2 \! + \! \mathcal{O}(\tau^{-2/3})$,
\begin{equation} \label{prcy100} 
\lvert \lvert \hat{\chi}_{k}(\widetilde{\Lambda}) \! - \! \mathrm{I} \rvert \rvert \underset{\tau \to +\infty}{\leqslant} 
\mathcal{O} \! \left(\mathfrak{c}_{k}^{\Ydown}(\tau) \lvert \nu (k) \! + \! 1 \rvert^{2} \lvert p_{k}(\tau) \rvert^{-2} 
\tau^{-\big(\frac{1}{3}-2(3+ \Re (\nu (k)+1)) \delta_{k} \big)} \right), \quad k \! = \! +1,
\end{equation}
where, for $\arg (\widetilde{\Lambda}) \genfrac{}{}{0pt}{3}{=}{\tau \to +\infty} \tfrac{\pi}{2} \! + \! \mathcal{O}(\tau^{-2/3})$, 
$\mathfrak{c}_{1}^{\Ydown}(\tau) \! := \! \tilde{\mathfrak{c}}_{M}^{\sharp} \tilde{\mathfrak{c}}_{M}(2^{3/2}3^{1/4} 
\me^{\pi \Im (\nu (1)+1)/2})^{-2}$, and, for $\arg (\widetilde{\Lambda}) \genfrac{}{}{0pt}{3}{=}{\tau \to +\infty} 
-\tfrac{\pi}{2} \! + \! \mathcal{O}(\tau^{-2/3})$, $\mathfrak{c}_{1}^{\Ydown}(\tau) \! := \! \hat{\mathfrak{c}}_{M}^{\sharp} 
\hat{\mathfrak{c}}_{M}(2^{3/2}3^{1/4} \me^{\pi \Im (\nu (1)+1)/2})^{-2}$ (see Remark \ref{aspvals}). Via an analogous 
series of calculations, one arrives at a similar estimate (cf. asymptotics \eqref{prcy100}) for the case $k \! = \! -1$.

Forming the composition of the inverses of the linear transformations $\mathfrak{F}_{j}$, $j \! = \! 1,2,\dotsc,11$, that is,
\begin{align}
\widetilde{\Psi}_{k}(\widetilde{\mu},\tau) \! :=& \, \left(\mathfrak{F}_{1}^{-1} \circ \mathfrak{F}_{2}^{-1} \circ 
\mathfrak{F}_{3}^{-1} \circ \mathfrak{F}_{4}^{-1} \circ \mathfrak{F}_{5}^{-1} \circ \mathfrak{F}_{6}^{-1} \circ 
\mathfrak{F}_{7}^{-1} \circ \mathfrak{F}_{8}^{-1} \circ \mathfrak{F}_{9}^{-1} \circ \mathfrak{F}_{10}^{-1} \circ 
\mathfrak{F}_{11}^{-1} \right) \! \Phi_{M,k}(\widetilde{\Lambda}) \nonumber \\
=& \, (b(\tau))^{-\frac{1}{2} \sigma_{3}} \mathcal{G}_{0,k} \mathcal{G}_{1,k} \tau^{\frac{1}{6} \sigma_{3}} 
\begin{pmatrix}
1 & 0 \\
-\mi \omega_{0,k} & 1
\end{pmatrix} \! 
\begin{pmatrix}
1 & \ell_{0,k}^{+} \widetilde{\Lambda} \\
0 & 1
\end{pmatrix} \! 
\begin{pmatrix}
1 & 0 \\
\ell_{1,k}^{+} \widetilde{\Lambda} & 1
\end{pmatrix} \! \mathcal{G}_{2,k} \! 
\begin{pmatrix}
1 & 0 \\
\ell_{2,k}^{+} \widetilde{\Lambda} & 1
\end{pmatrix} \nonumber \\
\times& \, \hat{\chi}_{k}(\widetilde{\Lambda}) \Phi_{M,k}(\widetilde{\Lambda}), 
\quad k \! = \! \pm 1,
\label{prcy101}
\end{align}
one arrives at the asymptotic representation for $\widetilde{\Psi}_{k}(\widetilde{\mu},\tau)$ given in equation 
\eqref{prpr2}. \hfill $\qed$
\begin{eeee} \label{delrestsn} 
Heretofore, it was assumed that (cf. Corollaries \ref{cor3.1.2}--\ref{cor3.1.5}) $0 \! < \! \delta \! < \! \delta_{k} 
\! < \! 1/9$, $k \! = \! \pm 1$; however, the set of restrictions \eqref{pc4} implies the following, more stringent 
restriction on $\delta_{k}$:\footnote{Note: $18 \genfrac{}{}{0pt}{3}{<}{\tau \to +\infty} 6(3 \! + \! \Re (\nu (k) \! 
+ \! 1)) \genfrac{}{}{0pt}{3}{<}{\tau \to +\infty} 24$.}
\begin{equation} \label{restr1} 
0 \underset{\tau \to +\infty}{<} \delta_{k} \underset{\tau \to +\infty}{<} 1/24, 
\quad k \! = \! \pm 1.
\end{equation}
Since $(0,1/24) \subset (0,1/9)$, the latter restriction \eqref{restr1} on $\delta_{k}$ implies, and is consistent with, 
the earlier one; henceforth, the restriction \eqref{restr1} on $\delta_{k}$ will be enforced. \hfill $\blacksquare$
\end{eeee}
\begin{eeee} \label{aspvals} 
Using the fact that (see the asymptotics \eqref{geek13} below) $\nu (k) \! + \! 1 \! \to \! 0$ as $\tau \to +\infty$, 
$k \! = \! \pm 1$, one shows, via the expansion for the gamma function \cite{a24}
\begin{equation*}
\dfrac{1}{\Gamma (z \! + \! 1)} \! = \! \sum_{j=0}^{\infty} \mathfrak{d}_{j}^{\ast}z^{j}, \quad \lvert z \rvert \! < \! 1,
\end{equation*}
where $\mathfrak{d}_{0}^{\ast} \! = \! 1$ and $\mathfrak{d}_{n+1}^{\ast} \! = \! \tfrac{1}{(n+1)} \sum_{j=0}^{n}(-1)^{j}
s_{j+1} \mathfrak{d}_{n-j}^{\ast}$, $n \! \in \! \mathbb{Z}_{+}$, with $s_{1} \! = \! -\psi (1)$ the Euler-Mascheroni 
constant,\footnote{$-\psi (1) \! = \! 0.577215664901532860606512 \ldots$.} and $s_{m} \! = \! \zeta (m)$, 
$\mathbb{N} \! \ni \! m \! \geqslant \! 2$, where $\zeta (z)$ is the Riemann Zeta function, and well-known inequalities 
for complex-valued trigonometric functions, that the auxiliary parameters introduced in step~\pmb{(xi)} of the 
proof of Lemma \ref{nprcl} have (for the case $k \! = \! +1)$ the following asymptotics: \pmb{(1)} for $\arg 
(\widetilde{\Lambda}) \genfrac{}{}{0pt}{3}{=}{\tau \to +\infty} \tfrac{\pi}{2} \! + \! \mathcal{O} \big(\tau^{-2/3} \big)$,
\begin{gather*}
(\tilde{\varrho}_{0}(1))^{2} \underset{\tau \to +\infty}{=} (2 \! + \! \vert \sec 
\theta \rvert)^{2}(1 \! + \! \mathcal{O}(\lvert \nu (1) \! + \! 1 \rvert)), \\
(\tilde{\varrho}_{1}(1))^{2} \underset{\tau \to +\infty}{=} \dfrac{\pi}{2}(1 \! + 
\! 2 \sec^{2} \theta)^{2}(1 \! + \! \mathcal{O}(\lvert \nu (1) \! + \! 1 \rvert)), \\
(\tilde{\varrho}_{2}(1))^{2} \underset{\tau \to +\infty}{=} 192(2 \sqrt{\pi} \! + \! 
\vert \sec \theta \rvert)^{2}(1 \! + \! \mathcal{O}(\lvert \nu (1) \! + \! 1 \rvert)), \\
(\tilde{\varrho}_{3}(1))^{2} \underset{\tau \to +\infty}{=} 96 \pi (1 \! + \! 2 
\sec^{2} \theta)^{2}(1 \! + \! \mathcal{O}(\lvert \nu (1) \! + \! 1 \rvert)), \\
\tilde{\rho}_{0}^{2} \underset{\tau \to +\infty}{=} 1 \! + \! \mathcal{O}
(\lvert \nu (1) \! + \! 1 \rvert), \quad \quad \tilde{\rho}_{1}^{2} 
\underset{\tau \to +\infty}{=} 2 \pi \sec^{2}(\theta)(1 \! + \! \mathcal{O}
(\lvert \nu (1) \! + \! 1 \rvert)), \\
\tilde{\rho}_{2}^{2} \underset{\tau \to +\infty}{=} 16 \sqrt{3} \pi \lvert \nu 
(1) \! + \! 1 \rvert^{2}(1 \! + \! \mathcal{O}(\lvert \nu (1) \! + \! 1 \rvert)), 
\quad \quad \tilde{\rho}_{3}^{2} \underset{\tau \to +\infty}{=} 8 \sqrt{3}
(1 \! + \! \mathcal{O}(\lvert \nu (1) \! + \! 1 \rvert)),
\end{gather*}
where $\theta \! := \! \arg (\nu (1) \! + \! 1)$, whence $\tilde{\mathfrak{c}}_{M}^{\sharp} \genfrac{}{}{0pt}{3}{=}{\tau \to +\infty} 
\mathcal{O}(1)$ and $\tilde{\mathfrak{c}}_{M} \genfrac{}{}{0pt}{3}{=}{\tau \to +\infty} \mathcal{O}(1)$ $\Rightarrow$ 
$\mathfrak{c}_{1}^{\Ydown}(\tau) \genfrac{}{}{0pt}{3}{=}{\tau \to +\infty} \mathcal{O}(1)$ (as claimed); and \pmb{(2)} 
for $\arg (\widetilde{\Lambda}) \genfrac{}{}{0pt}{3}{=}{\tau \to +\infty} -\tfrac{\pi}{2} \! + \! \mathcal{O} \big(\tau^{-2/3} \big)$,
\begin{gather*}
(\hat{\varrho}_{0}(1))^{2} \underset{\tau \to +\infty}{=} \sec^{2}(\theta)
(1 \! + \! \mathcal{O}(\lvert \nu (1) \! + \! 1 \rvert)), \quad \quad 
(\hat{\varrho}_{1}(1))^{2} \underset{\tau \to +\infty}{=} \dfrac{\pi}{2}
(1 \! + \! \mathcal{O}(\lvert \nu (1) \! + \! 1 \rvert)), \\
(\hat{\varrho}_{2}(1))^{2} \underset{\tau \to +\infty}{=} 192 \sec^{2}
(\theta)(1 \! + \! \mathcal{O}(\lvert \nu (1) \! + \! 1 \rvert)), \quad 
\quad (\hat{\varrho}_{3}(1))^{2} \underset{\tau \to +\infty}{=} 96 
\pi (1 \! + \! \mathcal{O}(\lvert \nu (1) \! + \! 1 \rvert)), \\
\hat{\rho}_{0}^{2} \underset{\tau \to +\infty}{=} 1 \! + \! \mathcal{O}
(\lvert \nu (1) \! + \! 1 \rvert), \quad \quad \hat{\rho}_{1}^{2} 
\underset{\tau \to +\infty}{=} \dfrac{1}{8 \sqrt{3}}(1 \! + \! \mathcal{O}
(\lvert \nu (1) \! + \! 1 \rvert)), \\
\hat{\rho}_{2}^{2} \underset{\tau \to +\infty}{=} \mathcal{O}
(\lvert \nu (1) \! + \! 1 \rvert^{2}), \quad \quad \hat{\rho}_{3}^{2} 
\underset{\tau \to +\infty}{=} 8 \sqrt{3}(1 \! + \! \mathcal{O}
(\lvert \nu (1) \! + \! 1 \rvert)),
\end{gather*}
whence $\hat{\mathfrak{c}}_{M}^{\sharp} \genfrac{}{}{0pt}{3}{=}{\tau \to +\infty} \mathcal{O}(1)$ and 
$\hat{\mathfrak{c}}_{M} \genfrac{}{}{0pt}{3}{=}{\tau \to +\infty} \mathcal{O}(1)$ $\Rightarrow$ 
$\mathfrak{c}_{1}^{\Ydown}(\tau) \genfrac{}{}{0pt}{3}{=}{\tau \to +\infty} \mathcal{O}(1)$ (as claimed). 
The analysis for the case $k \! = \! -1$ is analogous. \hfill $\blacksquare$
\end{eeee}
\begin{eeee} \label{leedasyue} 
In Lemma \ref{nprcl} and hereafter, the function $\Phi_{M,k}(\pmb{\cdot})$ plays a crucial r\^{o}le; therefore, 
its asymptotics are presented here: for $m \! \in \! \lbrace -1,0,1,2 \rbrace$ and $k \! = \! \pm 1$,
\begin{equation*}
\Phi_{M,k}(z) \underset{\underset{\scriptscriptstyle \arg (z)=\frac{m \pi}{2}+\frac{\pi}{4}-\frac{1}{2} \arg 
(\mu_{k}(\tau))}{\scriptscriptstyle \mathbb{C} \ni z \to \infty}}{=} \left(\mathrm{I} \! + \! \sum_{j=1}^{\infty} 
\hat{\psi}_{j,k}(\tau)z^{-j} \right) \! \me^{\left(\frac{1}{2} \mu_{k}(\tau)z^{2} - (\nu (k)+1) \ln ((2 \mu_{k}(\tau))^{1/2}z) 
\right) \sigma_{3}} \mathcal{R}_{m}(k),
\end{equation*}
where
\begin{gather*}
\mathcal{R}_{-1}(k) \! := \! 
\begin{pmatrix}
\me^{-\mi \pi (\nu (k)+1)/2} & 0 \\
0 & -\frac{(2 \mu_{k}(\tau))^{1/2}}{p_{k}(\tau)}
\end{pmatrix}, \\
\mathcal{R}_{0}(k) \! := \! 
\begin{pmatrix}
\me^{-\mi \pi (\nu (k)+1)/2} & 0 \\
-\frac{\mi \sqrt{2 \pi}(2 \mu_{k}(\tau))^{1/2} \me^{-\mi \pi (\nu (k)+1)/2}}{p_{k}(\tau) \Gamma (\nu (k)+1)} 
& -\frac{(2 \mu_{k}(\tau))^{1/2}}{p_{k}(\tau)}
\end{pmatrix}, \\
\mathcal{R}_{1}(k) \! := \! 
\begin{pmatrix}
\me^{\mi 3 \pi (\nu (k)+1)/2} & \frac{\sqrt{2 \pi} \me^{\mi \pi (\nu (k)+1)}}{\Gamma (-\nu (k))} \\
-\frac{\mi \sqrt{2 \pi}(2 \mu_{k}(\tau))^{1/2} \me^{-\mi \pi (\nu (k)+1)/2}}{p_{k}(\tau) \Gamma (\nu (k)+1)} 
& -\frac{(2 \mu_{k}(\tau))^{1/2}}{p_{k}(\tau)}
\end{pmatrix}, \\
\mathcal{R}_{2}(k) \! := \! 
\begin{pmatrix}
\me^{\mi 3 \pi (\nu (k)+1)/2} & \frac{\sqrt{2 \pi} \me^{\mi \pi (\nu (k)+1)}}{\Gamma (-\nu (k))} \\
0 &  -\frac{(2 \mu_{k}(\tau))^{1/2} \me^{-2 \pi \mi (\nu (k)+1)}}{p_{k}(\tau)}
\end{pmatrix},
\end{gather*}
and $\hat{\psi}_{j,k}(\tau)$, $j \! \in \! \mathbb{N}$, are off-diagonal (resp., diagonal) $\mathrm{M}_{2}(\mathbb{C})$-valued 
functions for $j$ odd (resp., $j$ even); for example, 
\begin{gather*}
\hat{\psi}_{1,k}(\tau) \! = \! -\dfrac{1}{2 \mu_{k}(\tau)} \! 
\begin{pmatrix}
0 & p_{k}(\tau) \\
-q_{k}(\tau) & 0
\end{pmatrix}, \, \, \qquad \, \, \hat{\psi}_{2,k}(\tau) \! = \! 
\dfrac{(\nu (k) \! + \! 1)}{4 \mu_{k}(\tau)} \! 
\begin{pmatrix}
1 \! + \! (\nu (k) \! + \! 1) & 0 \\
0 & 1 \! - \! (\nu (k) \! + \! 1)
\end{pmatrix}, \\
\hat{\psi}_{3,k}(\tau) \! = \! \dfrac{1}{8(\mu_{k}(\tau))^{2}} \! 
\begin{pmatrix}
0 & (1 \! - \! (\nu (k) \! + \! 1))(2 \! - \! (\nu (k) \! + \! 1))p_{k}(\tau) \\
(1 \! + \! (\nu (k) \! + \! 1))(2 \! + \! (\nu (k) \! + \! 1))q_{k}(\tau) & 0
\end{pmatrix}.
\end{gather*}
These asymptotics can be derived {}from the asymptotics of the parabolic-cylinder functions \cite{EMOT}. \hfill $\blacksquare$
\end{eeee}
\subsection{Asymptotic Matching} \label{sec3.3} 
In this subsection, the connection matrix is calculated asymptotically (as $\tau \! \to \! +\infty$ with $\varepsilon b \! > \! 0)$ 
in terms of the matrix elements of the function $\mathcal{A}(\widetilde{\mu},\tau)$ (cf. equation \eqref{eq3.4}) that are 
defined in terms of the set of functions $\hat{h}_{0}(\tau)$, $\tilde{r}_{0}(\tau)$, $v_{0}(\tau)$, and $b(\tau)$ concomitant 
with the conditions \eqref{iden5}.\footnote{Equivalently, the set of functions (cf. equations \eqref{iden2}, \eqref{iden3}, and 
\eqref{iden4}, respectively) $h_{0}(\tau)$, $\hat{r}_{0}(\tau)$, and $\hat{u}_{0}(\tau)$.} Thus, the direct monodromy problem 
for equation \eqref{eq3.3} is solved asymptotically.
\begin{cccc} \label{linfnewlemm} 
Let $\widetilde{\Psi}_{k}(\widetilde{\mu},\tau)$, $k \! = \! \pm 1$, be the fundamental solution of equation \eqref{eq3.3} 
with asymptotics given in Lemma \ref{nprcl}, and let $\mathbb{Y}^{\infty}_{0}(\widetilde{\mu},\tau)$ be the canonical 
solution of equation \eqref{eq3.1}.\footnote{See Proposition \ref{prop1.4}.} Define {}\footnote{Since $\tau^{-\frac{1}{12} 
\sigma_{3}} \mathbb{Y}^{\infty}_{0}(\tau^{-1/6} \widetilde{\mu},\tau)$ (cf. equations \eqref{eq3.2}) is also a fundamental 
solution of equation \eqref{eq3.3}, it follows, therefore, that $\mathfrak{L}^{\infty}_{k}(\tau)$ is independent of 
$\widetilde{\mu}$.}
\begin{equation} \label{ellinfk1} 
\mathfrak{L}^{\infty}_{k}(\tau) \! := \! \big(\widetilde{\Psi}_{k}(\widetilde{\mu},\tau) \big)^{-1} \tau^{-\frac{1}{12} \sigma_{3}} 
\mathbb{Y}^{\infty}_{0}(\tau^{-1/6} \widetilde{\mu},\tau), \quad k \! = \! \pm 1.
\end{equation}
Assume that the parameters $\nu (k) \! + \! 1$ and $\delta_{k}$ satisfy the restrictions \eqref{pc4} and \eqref{restr1}, 
respectively, and, additionally, the following conditions are valid:\footnote{The conditions \eqref{iden5} and \eqref{restr1} 
are consistent with the conditions \eqref{ellinfk2a} and \eqref{ellinfk2b}.}
\begin{gather}
p_{k}(\tau) \mathfrak{B}_{k} \exp \! \left(-\mi \tau^{2/3}3 \sqrt{3}(\varepsilon b)^{1/3} \me^{\mi 2 \pi k/3} \right) 
\underset{\tau \to +\infty}{=} \mathcal{O} \! \left((\nu (k) \! + \! 1)^{\frac{1-k}{2}} \right), \label{ellinfk2a} \\
b(\tau) \tau^{\mi a/3} \exp \! \left(\mi \tau^{2/3}3(\varepsilon b)^{1/3} \me^{\mi 2 \pi k/3} \right) 
\underset{\tau \to +\infty}{=} \mathcal{O}(1), \label{ellinfk2b}
\end{gather}
where $p_{k}(\tau)$ and $\mathfrak{B}_{k}$ are defined in Lemma \ref{nprcl}.\footnote{{}From the results subsumed in the 
proof of Lemma \ref{ginversion}, it will be deduced \emph{a posteriori} that (cf. definition \eqref{prcy54}) $\mu_{k}(\tau)$ 
possesses the asymptotics $\mu_{k}(\tau) \genfrac{}{}{0pt}{3}{=}{\tau \to +\infty} \mi 4 \sqrtsign{3} \! + \! 
\sum_{\scriptscriptstyle \underset{m_{1}+m_{2}+m_{3} \geqslant 2}{m_{1},m_{2},m_{3} \in \mathbb{Z}_{\scriptscriptstyle +}}}
c_{m_{1},m_{2},m_{3}}(k)(\tilde{r}_{0}(\tau))^{m_{1}}(v_{0}(\tau))^{m_{2}}(\tau^{-1/3})^{m_{3}} \! + \! c_{\infty}(k) \tau^{-1/3} 
\me^{-(\beta (\tau)+\mi k \vartheta (\tau))} \big(1 \! + \! \mathcal{O}(\tau^{-1/3}) \big)$, $k \! = \! \pm 1$, where 
$c_{m_{1},m_{2},m_{3}}(k) \! \in \! \mathbb{C}$, and $\vartheta (\tau)$ and $\beta (\tau)$ are defined in equations 
\eqref{thmk12}; via this fact, and the definitions \eqref{prpr1}, \eqref{prcy22}, \eqref{prcy57}, and \eqref{prcy58}, a lengthy and 
circuitous calculation reveals that the asymptotic expansion of $\nu (k) \! + \! 1$, $k \! = \! \pm 1$, can be presented in the form
\begin{align*}
-(\nu (k) \! + \! 1) \underset{\tau \to +\infty}{=}& \, \frac{\mi}{8 \sqrt{3}} \! \left(\dfrac{-\alpha_{k}^{2}(8v_{0}^{2}(\tau) 
\! + \! 4v_{0}(\tau) \tilde{r}_{0}(\tau) \! - \! (\tilde{r}_{0}(\tau))^{2} \! - \! v_{0}(\tau)(\tilde{r}_{0}(\tau))^{2} \tau^{-1/3}) 
\! + \! 4(a \! - \! \mi/2)v_{0}(\tau) \tau^{-1/3}}{1 \! + \! v_{0}(\tau) \tau^{-1/3}} \right) \\
+& \, \frac{2 \mathfrak{p}_{k}(\tau)}{3 \sqrt{3} \alpha_{k}^{2}} \! + \! \sum_{m=2}^{\infty} \hat{\mu}_{m}^{\ast}(k)
(\tau^{-1/3})^{m} \! + \! \hat{c}_{\infty}(k) \tau^{-1/3} \me^{-(\beta (\tau)+\mi k \vartheta (\tau))} \big(1 \! + \! 
\mathcal{O}(\tau^{-1/3}) \big),
\end{align*}
where $\mathfrak{p}_{k}(\tau)$ is defined by equation \eqref{eqpeetee}. {}From the asymptotics \eqref{tr1} 
and \eqref{tr3}, and Propositions \ref{recursys} and \ref{proprr}, in conjunction with the formulae for the 
monodromy-data-dependent expansion coefficients $\mathrm{A}_{k}$, $k \! = \! \pm 1$, derived in the proof 
of Lemma \ref{ginversion} (see, in particular, equations \eqref{geek109} and \eqref{geek111}), the sum of the 
coefficients of each term $(\tau^{-1/3})^{j}$, $\mathbb{N} \! \ni \! j \! \geqslant \! 2$, and of the term $\tau^{-1/3} 
\me^{-(\beta (\tau)+\mi k \vartheta (\tau))}$ on the right-hand side of the latter asymptotic expansion for 
$\nu (k) \! + \! 1$ are equal to zero (e.g., $\hat{\mu}_{2}^{\ast}(k) \! = \! -\frac{\mi}{24 \sqrt{3} \alpha_{k}^{2}}((a \! - 
\! \mi/2)^{2} \! - \! 1/6))$, resulting, finally, in the asymptotics $\nu (k) \! + \! 1 \genfrac{}{}{0pt}{3}{=}{\tau \to +\infty} 
\mathcal{O}(\tau^{-2/3} \me^{-\beta (\tau)})$, $k \! = \! \pm 1$. The conditions \eqref{ellinfk2a} and \eqref{ellinfk2b} will 
be validated \emph{a posteriori\/} (see the proof of Lemma \ref{ginversion}) using the asymptotics $\nu (k) \! + \! 1 
\genfrac{}{}{0pt}{3}{=}{\tau \to +\infty} \mathcal{O}(\tau^{-2/3} \me^{-\beta (\tau)})$, $k \! = \! \pm 1$. Hereafter, whilst reading 
the text, the reader should be cognizant of the latter asymptotics for $\nu (k) \! + \! 1$, as all asymptotic expansions, 
estimates, orderings, etc., rely on this fact.} Then,
\begin{align} \label{ellinfk3} 
\mathfrak{L}^{\infty}_{k}(\tau) \underset{\tau \to +\infty}{=}& \, \mi (\mathcal{R}_{m_{\infty}}(k))^{-1} 
\me^{\tilde{\mathfrak{z}}_{k}^{0}(\tau) \sigma_{3}} \! \left(\dfrac{(\varepsilon b)^{1/4}(\sqrt{3} \! + \! 1)^{1/2}}{2^{1/4} 
\sqrt{\smash[b]{\mathfrak{B}_{k}}} \sqrt{\smash[b]{b(\tau)}}} \right)^{\sigma_{3}} \sigma_{2} 
\me^{-\Delta \tilde{\mathfrak{z}}_{k}(\tau) \sigma_{3}} 
\begin{pmatrix}
\hat{\mathbb{B}}_{0}^{\infty}(\tau) & 0 \\
0 & \hat{\mathbb{A}}_{0}^{\infty}(\tau)
\end{pmatrix} \nonumber \\
\times& \, \left(\mathrm{I} \! + \! \mathbb{E}_{{\scriptscriptstyle \mathcal{N}},k}^{\infty}(\tau) \right) \! \left(\mathrm{I} \! 
+ \! \mathcal{O} \! \left(\mathbb{E}^{\infty}_{k}(\tau) \right) \right),
\end{align}
where $\mathrm{M}_{2}(\mathbb{C}) \! \ni \! \mathcal{R}_{m_{\infty}}(k)$, $m_{\infty} \! \in \! \lbrace -1,0,1,2 \rbrace$, 
are defined in Remark \ref{leedasyue},
\begin{align} 
\tilde{\mathfrak{z}}_{k}^{0}(\tau) \! :=& \, -\dfrac{\mi a}{6} \ln \tau \! + \! \mi \tau^{2/3}3(\sqrt{3} \! - \! 1) 
\alpha_{k}^{2} \! + \! \mi (a \! - \! \mi/2) \ln ((\sqrt{3} \! + \! 1) \alpha_{k}/2), \label{ellinfk4} \\
\Delta \tilde{\mathfrak{z}}_{k}(\tau) \! :=& \, -\left(\dfrac{5 \! - \! \sqrt{3}}{6 \sqrt{3} \alpha_{k}^{2}} \right) \! 
\mathfrak{p}_{k}(\tau) \! + \! (\nu (k) \! + \! 1) \ln (2 \mu_{k}(\tau))^{1/2} \! + \! \dfrac{1}{3}(\nu (k) \! + \! 1) \ln \tau \nonumber \\
+& \, (\nu (k) \! + \! 1) \ln (6(\sqrt{3} \! + \! 1)^{-2} \alpha_{k}), \label{ellinfk5}
\end{align}
with $\mathfrak{p}_{k}(\tau)$ defined by equation \eqref{eqpeetee}, and $\mu_{k}(\tau)$ defined in Lemma \ref{nprcl},
\begin{gather}
\hat{\mathbb{A}}_{0}^{\infty}(\tau) \! := \! 1 \! + \! \dfrac{2^{1/4}(\Delta G_{k}^{\infty}(\tau))_{12}}{(\varepsilon b)^{1/4}
(\sqrt{3} \! + \! 1)^{1/2}}, \label{ellinfk6} \\
\hat{\mathbb{B}}_{0}^{\infty}(\tau) \! := \! 1 \! - \! \dfrac{(\varepsilon b)^{1/4}(\sqrt{3} \! + \! 1)^{1/2}}{2^{1/4}} \! \left(
(\Delta G_{k}^{\infty}(\tau))_{21} \! - \! \dfrac{\mathfrak{A}_{k}}{\mathfrak{B}_{k}} \! \left(\dfrac{\mi 4 \sqrt{3} 
\mathcal{Z}_{k}}{\chi_{k}(\tau)} \! - \! 1 \right) \! (\Delta G_{k}^{\infty}(\tau))_{11} \right), \label{ellinfk7}
\end{gather}
with $\mathcal{Z}_{k}$, $\mathfrak{A}_{k}$, and $\chi_{k}(\tau)$ defined in Lemma \ref{nprcl}, and
\begin{equation} \label{ellinfk8} 
\Delta G_{k}^{\infty}(\tau) \! := \! \dfrac{1}{(2 \sqrt{3}(\sqrt{3} \! + \! 1))^{1/2}} 
\begin{pmatrix}
(\Delta G_{k}^{\infty}(\tau))_{11} & (\Delta G_{k}^{\infty}(\tau))_{12} \\
(\Delta G_{k}^{\infty}(\tau))_{21} & (\Delta G_{k}^{\infty}(\tau))_{22}
\end{pmatrix},
\end{equation}
with
\begin{equation} \label{forinff}
\begin{aligned}
(\Delta G_{k}^{\infty}(\tau))_{11} \! =& \, (\sqrt{3} \! + \! 1)(\Delta \mathcal{G}_{0,k})_{22} \! + \! (2/\varepsilon b)^{1/2}
(\Delta \mathcal{G}_{0,k})_{12}, \\
(\Delta G_{k}^{\infty}(\tau))_{12} \! =& -(\sqrt{3} \! + \! 1)(\Delta \mathcal{G}_{0,k})_{12} \! + \! (2 \varepsilon b)^{1/2}
(\Delta \mathcal{G}_{0,k})_{22}, \\
(\Delta G_{k}^{\infty}(\tau))_{21} \! =& -(\sqrt{3} \! + \! 1)(\Delta \mathcal{G}_{0,k})_{21} \! - \! (2/\varepsilon b)^{1/2}
(\Delta \mathcal{G}_{0,k})_{11}, \\
(\Delta G_{k}^{\infty}(\tau))_{22} \! =& \, (\sqrt{3} \! + \! 1)(\Delta \mathcal{G}_{0,k})_{11} \! - \! (2 \varepsilon b)^{1/2}
(\Delta \mathcal{G}_{0,k})_{21},
\end{aligned}
\end{equation}
where $(\Delta \mathcal{G}_{0,k})_{i,j=1,2}$ are defined by equations \eqref{prcyg4}--\eqref{prcyg6},
\begin{align} \label{ellinfk9} 
\mathbb{E}_{{\scriptscriptstyle \mathcal{N}},k}^{\infty}(\tau) :=& \, 
\left(\dfrac{\me^{-\tilde{\beta}_{k}(\tau)}}{\sqrt{\smash[b]{b(\tau)}}} 
\right)^{\ad (\sigma_{3})} \! \left(-\dfrac{\mi 4 \sqrt{3} \mathcal{Z}_{k} \mathfrak{A}_{k} 
\ell_{0,k}^{+}}{\chi_{k}(\tau)} \! \left(\dfrac{(\varepsilon b)^{1/2}(\sqrt{3} \! + \! 1)
(\nu (k) \! + \! 1)}{\sqrt{2}p_{k}(\tau) \mathfrak{B}_{k}} \sigma_{+} \right. \right. \nonumber \\
+&\left. \left. \dfrac{p_{k}(\tau) \mathfrak{B}_{k}}{\sqrt{2}(\varepsilon b)^{1/2}(\sqrt{3} \! + \! 1) \mu_{k}(\tau)} 
\sigma_{-} \right) \! \sigma_{3} \! + \! \dfrac{1}{2 \sqrt{3}(\sqrt{3} \! + \! 1)} \right. \nonumber \\
\times&\left. 
\begin{pmatrix}
\frac{\mi 4 \sqrt{3} \mathcal{Z}_{k} \mathfrak{A}_{k} \ell_{0,k}^{+}}{\chi_{k}(\tau)} & 
-\frac{(\varepsilon b)^{1/2}(\sqrt{3}+1)}{\sqrt{2} \mathfrak{B}_{k}} \big((\frac{\mi 4 \sqrt{3} 
\mathcal{Z}_{k} \mathfrak{A}_{k}}{\chi_{k}(\tau)})^{2} \ell_{0,k}^{+} \! - \! \ell_{1,k}^{+} 
\! - \! \ell_{2,k}^{+} \big) \\
\frac{\sqrt{2} \mathfrak{B}_{k} \ell_{0,k}^{+}}{(\varepsilon b)^{1/2}(\sqrt{3}+1)} 
& -\frac{\mi 4 \sqrt{3} \mathcal{Z}_{k} \mathfrak{A}_{k} \ell_{0,k}^{+}}{\chi_{k}(\tau)} 
\end{pmatrix} \right. \nonumber \\
\times&\left. 
\begin{pmatrix}
\sqrt{3} \! + \! 1 & -(2 \varepsilon b)^{1/2} \\
(2/\varepsilon b)^{1/2} & \sqrt{3} \! + \! 1
\end{pmatrix} \! 
\begin{pmatrix}
\mathbb{T}_{11,k}(1;\tau) & \mathbb{T}_{12,k}(1;\tau) \\
\mathbb{T}_{21,k}(1;\tau) & \mathbb{T}_{22,k}(1;\tau)
\end{pmatrix} \right),
\end{align}
with $\ell_{0,k}^{+}$, $\ell_{1,k}^{+}$, and $\ell_{2,k}^{+}$ defined in Lemma \ref{nprcl}, $(\mathbb{T}_{ij,k}(1;\tau))_{i,j=1,2}$ 
defined in Proposition \ref{prop3.1.6}, and $\tilde{\beta}_{k}(\tau)$ defined by equation \eqref{ellinfk18}, and
\begin{equation} \label{ellinfk10} 
\mathcal{O} \big(\mathbb{E}^{\infty}_{k}(\tau) \big) \underset{\tau \to +\infty}{:=} 
\begin{pmatrix}
\mathcal{O} \big(\tau^{-\frac{1}{3}+3\delta_{k}} \big) & \mathcal{O} \big(\tau^{-\frac{1}{3}(\frac{1+k}{2})-\delta_{k}} \big) \\
\mathcal{O} \big(\tau^{-\frac{1}{3}(\frac{1-k}{2})-\delta_{k}} \big) & \mathcal{O} \big(\tau^{-\frac{1}{3}+3\delta_{k}} \big)
\end{pmatrix}.
\end{equation}
\end{cccc}

\emph{Proof}. Denote by $\widetilde{\Psi}_{\scriptscriptstyle \mathrm{WKB},k}(\widetilde{\mu},\tau)$, $k \! = \! \pm 1$, 
the solution of equation \eqref{eq3.3} that has leading-order asymptotics given by equations \eqref{eq3.16}--\eqref{eq3.18} 
in the canonical domain containing the Stokes curve approaching, for $k \! = \! +1$ (resp., $k \! = \! -1)$, the positive real 
$\widetilde{\mu}$-axis {}from above (resp., below) as $\widetilde{\mu} \! \to \! +\infty$. Let $\mathfrak{L}_{k}^{\infty}(\tau)$, 
$k \! = \! \pm 1$, be defined by equation \eqref{ellinfk1}; rewrite $\mathfrak{L}_{k}^{\infty}(\tau)$ in the following form:
\begin{equation} \label{ellinfk11} 
\mathfrak{L}_{k}^{\infty}(\tau) \! = \! \left(\big(\widetilde{\Psi}_{k}(\widetilde{\mu},\tau) \big)^{-1} 
\widetilde{\Psi}_{\scriptscriptstyle \mathrm{WKB},k}(\widetilde{\mu},\tau) \right) \! \left(
\big(\widetilde{\Psi}_{\scriptscriptstyle \mathrm{WKB},k}(\widetilde{\mu},\tau) \big)^{-1} 
\tau^{-\frac{1}{12} \sigma_{3}} \mathbb{Y}^{\infty}_{0}(\tau^{-1/6} \widetilde{\mu},\tau) \right).
\end{equation}
Taking note of the fact that $\widetilde{\Psi}_{k}(\widetilde{\mu},\tau)$, $\widetilde{\Psi}_{\scriptscriptstyle \mathrm{WKB},k}
(\widetilde{\mu},\tau)$, and $\tau^{-\frac{1}{12} \sigma_{3}} \mathbb{Y}^{\infty}_{0}(\tau^{-1/6} \widetilde{\mu},\tau)$ are 
all solutions of equation \eqref{eq3.3}, it follows that they differ on the right by non-degenerate, 
$\widetilde{\mu}$-independent, $\mathrm{M}_{2}(\mathbb{C})$-valued factors: via this observation, one evaluates, 
asymptotically, each of the factors appearing in equation \eqref{ellinfk11} by considering separate limits, namely, 
$\widetilde{\mu} \! \to \! \alpha_{k}$ and $\widetilde{\mu} \! \to \! +\infty$, respectively; more specifically, for $k \! = \! \pm 1$,
\begin{align} \label{ellinfk12} 
&\big(\widetilde{\Psi}_{k}(\widetilde{\mu},\tau) \big)^{-1} \widetilde{\Psi}_{\scriptscriptstyle \mathrm{WKB},k}
(\widetilde{\mu},\tau) \underset{\tau \to +\infty}{:=} \nonumber \\
&\underbrace{\left((b(\tau))^{-\frac{1}{2} \sigma_{3}} \mathcal{G}_{0,k} \mathfrak{B}_{k}^{\frac{1}{2} \sigma_{3}} 
\mathbb{F}_{k}(\tau) \Xi_{k}(\tau;\widetilde{\Lambda}) \hat{\chi}_{k}(\widetilde{\Lambda}) \Phi_{M,k}(\widetilde{\Lambda}) 
\right)^{-1}T(\widetilde{\mu}) \me^{\scriptscriptstyle \mathrm{W}_{k}(\widetilde{\mu},\tau)}}_{\widetilde{\mu}
=\widetilde{\mu}_{0,k}, \, \, \, \, \widetilde{\Lambda} \, \underset{\tau \to +\infty}{\thicksim} \, \mathcal{O}(\tau^{\delta_{k}}), 
\, \, 0< \delta < \delta_{k}<\frac{1}{24}, \, \, \, \, \arg (\widetilde{\Lambda})=\frac{\pi m_{\infty}}{2}+\frac{\pi}{4}-\frac{1}{2} 
\arg (\mu_{k}(\tau)), \, \, m_{\infty} \in \lbrace -1,0,1,2 \rbrace},
\end{align}
where (cf. Lemma~\ref{nprcl})
\begin{gather}
\mathbb{F}_{k}(\tau) \! = \! 
\begin{pmatrix}
1 & 0 \\
(\frac{\mi 4 \sqrt{3} \mathcal{Z}_{k}}{\chi_{k}(\tau)} \! - \! 1) \mathfrak{A}_{k} & 1
\end{pmatrix}, \label{ellinfk13} \\
\Xi_{k}(\tau;\widetilde{\Lambda}) \! = \! \mathrm{I} \! + \! \gimel_{\scriptscriptstyle A,k}(\tau) \widetilde{\Lambda} 
\! + \! \gimel_{\scriptscriptstyle B,k}(\tau) \widetilde{\Lambda}^{2}, \label{ellinfk14}
\end{gather}
and
\begin{equation} \label{ellinfk15} 
\hat{\chi}_{k}(\widetilde{\Lambda}) \underset{\tau \to +\infty}{=} \mathrm{I} \! + \! \mathcal{O} \! \left(\tilde{\mathfrak{C}}_{k}
(\tau) \lvert \nu (k) \! + \! 1 \rvert^{2} \lvert p_{k}(\tau) \rvert^{-2} \tau^{-\epsilon_{\mathrm{\scriptscriptstyle TP}}(k)} \right),
\end{equation}
with $\nu (k) \! + \! 1$, $p_{k}(\tau)$, $\widetilde{\mu}_{0,k}$, $\mathcal{G}_{0,k}$, $\mathfrak{A}_{k}$, $\mathfrak{B}_{k}$, 
$\mathcal{Z}_{k}$, $\gimel_{\scriptscriptstyle A,k}(\tau)$, $\gimel_{\scriptscriptstyle B,k}(\tau)$, $\mu_{k}(\tau)$, and 
$\chi_{k}(\tau)$ defined in Lemma \ref{nprcl}, $\mathrm{W}_{k}(\widetilde{\mu},\tau) \! := \! -\sigma_{3} \mi \tau^{2/3} 
\int_{\widetilde{\mu}_{0,k}}^{\widetilde{\mu}}l_{k}(\xi) \, \md \xi \! - \! \int_{\widetilde{\mu}_{0,k}}^{\widetilde{\mu}} \diag 
((T(\xi))^{-1} \partial_{\xi}T(\xi)) \, \md \xi$, $\epsilon_{\mathrm{\scriptscriptstyle TP}}(k) \! := \! \tfrac{1}{3} \! - \! 2(3 \! + 
\! \Re (\nu (k) \! + \! 1)) \delta_{k}$ $(> \! 0)$, and $\mathrm{M}_{2}(\mathbb{C}) \! \ni \! \tilde{\mathfrak{C}}_{k}(\tau) 
\genfrac{}{}{0pt}{3}{=}{\tau \to +\infty} \mathcal{O}(1)$, and
\begin{equation} \label{ellinfk16} 
\big(\widetilde{\Psi}_{\scriptscriptstyle \mathrm{WKB},k}(\widetilde{\mu},\tau) \big)^{-1} \tau^{-\frac{1}{12} \sigma_{3}} 
\mathbb{Y}^{\infty}_{0}(\tau^{-1/6} \widetilde{\mu},\tau) \underset{\tau \to +\infty}{:=} \lim_{\underset{\arg 
(\widetilde{\mu})=0}{\Omega_{0}^{\infty} \ni \widetilde{\mu} \to \infty}} \left(\big(T(\widetilde{\mu}) 
\me^{\scriptscriptstyle \mathrm{W}_{k}(\widetilde{\mu},\tau)} \big)^{-1} \tau^{-\frac{1}{12} \sigma_{3}} 
\mathbb{Y}^{\infty}_{0}(\tau^{-1/6} \widetilde{\mu},\tau) \right).
\end{equation}

One commences by considering the asymptotics subsumed in the definition~\eqref{ellinfk16}. {}From the asymptotics 
for $\mathbb{Y}^{\infty}_{0}(\tau^{-1/6} \widetilde{\mu},\tau)$ stated in Proposition \ref{prop1.4}, equations \eqref{iden3}, 
\eqref{iden4}, \eqref{expforeych}, \eqref{expforkapp}, \eqref{eq3.37}, \eqref{eq3.38}, \eqref{eqpeetee}, \eqref{eq3.41}, 
\eqref{asympforf7}, \eqref{eq3.51}, \eqref{prcy22}, and \eqref{prcyomg1}, one arrives at, via the conditions \eqref{iden5} 
and the asymptotics \eqref{terrbos6} and \eqref{asympforf3},
\begin{equation} \label{ellinfk17} 
\lim_{\underset{\arg (\widetilde{\mu})=0}{\Omega_{0}^{\infty} \ni \widetilde{\mu} \to \infty}} \left(\big(T(\widetilde{\mu}) 
\me^{\scriptscriptstyle \mathrm{W}_{k}(\widetilde{\mu},\tau)} \big)^{-1} \tau^{-\frac{1}{12} \sigma_{3}} \mathbb{Y}^{\infty}_{0}
(\tau^{-1/6} \widetilde{\mu},\tau) \right) \underset{\tau \to +\infty}{=} \exp \big(\tilde{\beta}_{k}(\tau) \sigma_{3} \big), 
\quad k \! = \! \pm 1,
\end{equation}
where
\begin{align} \label{ellinfk18} 
\tilde{\beta}_{k}(\tau) :=& \, \dfrac{\mi a}{6} \ln \tau \! - \! \mi \tau^{2/3}3(\sqrt{3} \! - \! 1) \alpha_{k}^{2} \! - \! \mi 2 \sqrt{3} 
\, \widetilde{\Lambda}^{2} \! - \! \mi (a \! - \! \mi/2) \ln ((\sqrt{3} \! + \! 1) \alpha_{k}/2) \! + \! \dfrac{(5 \! - \! \sqrt{3}) 
\mathfrak{p}_{k}(\tau)}{6 \sqrt{3} \alpha_{k}^{2}} \nonumber \\
+& \, \left(\dfrac{\mi}{2 \sqrt{3}} \! \left((a \! - \! \mi/2) \! + \! \alpha_{k}^{-2} \tau^{2/3} \hat{h}_{0}(\tau) \right) \! + \! 
\dfrac{2 \mathfrak{p}_{k}(\tau)}{3 \sqrt{3} \alpha_{k}^{2}} \right) \! \left(\dfrac{1}{3} \ln \tau \! - \! \ln \widetilde{\Lambda} 
\! + \! \ln \left(\dfrac{6 \alpha_{k}}{(\sqrt{3} \! + \! 1)^{2}} \right) \right) \nonumber \\
-& \, \dfrac{(\sqrt{3} \! - \! 1) \mathfrak{p}_{k}(\tau)}{\sqrt{3} \alpha_{k} \tau^{-1/3} \widetilde{\Lambda}} \! + \! \mathcal{O} 
\! \left(\left(\dfrac{\mathfrak{c}_{1,k} \tau^{-1/3} \! + \! \mathfrak{c}_{2,k} \tilde{r}_{0}(\tau)}{\widetilde{\Lambda}^{2}} \right) 
\! \left(\mathfrak{c}_{3,k} \tau^{-1/3} \! + \! \mathfrak{c}_{4,k}(\tilde{r}_{0}(\tau) \! + \! 4v_{0}(\tau)) \right) \right) \nonumber \\
+& \, \mathcal{O} \big(\tau^{-1/3} \widetilde{\Lambda}^{3} \big) \! + \! \mathcal{O} \big(\tau^{-1/3} \widetilde{\Lambda} \big) 
\! + \! \mathcal{O} \! \left(\dfrac{\tau^{-1/3}}{\widetilde{\Lambda}} \left(\mathfrak{c}_{5,k} \! + \! \mathfrak{c}_{6,k} \tau^{2/3} 
\hat{h}_{0}(\tau) \! + \! \mathfrak{c}_{7,k}(\tau^{2/3} \hat{h}_{0}(\tau))^{2} \right) \right) \nonumber \\
+& \, \mathcal{O} \! \left(\tau^{-2/3} \hat{d}_{0,k}(\tau) \! \left(\dfrac{1}{3} \ln \tau \! - \! \ln \widetilde{\Lambda} \right) \! \right),
\end{align}
$\mathfrak{c}_{m,k}$, $m \! = \! 1,2,\dotsc,7$, are $\mathcal{O}(1)$, and $\hat{d}_{0,k}(\tau)$ is defined in the proof of 
Proposition \ref{prop3.1.3}.

One now derives the asymptotics defined by equation \eqref{ellinfk12}. {}From the asymptotics \eqref{iden55} for $\varpi 
\! = \! +1$, equation \eqref{prcy2} for $\Phi_{M,k}(\widetilde{\Lambda})$ (in conjunction with its large-$\widetilde{\Lambda}$ 
asymptotics stated in Remark \ref{leedasyue}), the definitions \eqref{ellinfk13} and \eqref{ellinfk14} (concomitant with the 
fact that $\det (\Xi_{k}(\tau;\widetilde{\Lambda})) \! = \! 1)$, and the asymptotics \eqref{ellinfk15}, one shows, via the relation 
$(\mathrm{W}_{k}(\widetilde{\mu}_{0,k},\tau))_{i,j=1,2} \! = \! 0$ and the definition \eqref{ellinfk12}, that, for $k \! = \! \pm 1$,
{\fontsize{10pt}{11pt}\selectfont
\begin{align} \label{ellinfk19} 
\big(\widetilde{\Psi}_{k}(\widetilde{\mu},\tau) \big)^{-1} \widetilde{\Psi}_{\scriptscriptstyle 
\mathrm{WKB},k}(\widetilde{\mu},\tau) \underset{\tau \to +\infty}{:=}& \, 
\Phi_{M,k}^{-1}(\widetilde{\Lambda}) \hat{\chi}_{k}^{-1}(\widetilde{\Lambda}) 
\Xi_{k}^{-1}(\tau;\widetilde{\Lambda}) \mathbb{F}_{k}^{-1}(\tau) \mathfrak{B}_{k}^{
-\frac{1}{2} \sigma_{3}} \mathcal{G}_{0,k}^{-1}(b(\tau))^{\frac{1}{2} \sigma_{3}}
T(\widetilde{\mu}_{0,k}) \nonumber \\
\underset{\tau \to +\infty}{=}& \, (\mathcal{R}_{m_{\infty}}(k))^{-1} 
\me^{-\mathcal{P}_{0}^{\ast} \sigma_{3}} \mathfrak{Q}_{\infty,k}(\tau) \! \left(
\mathrm{I} \! + \! \dfrac{1}{\widetilde{\Lambda}} \mathfrak{Q}_{\infty,k}^{-1}
(\tau) \hat{\psi}_{1,k}^{-1}(\tau) \mathfrak{Q}_{\infty,k}(\tau) \right. \nonumber \\
+&\left. \, \dfrac{1}{\widetilde{\Lambda}^{2}} \mathfrak{Q}_{\infty,k}^{-1}
(\tau) \hat{\psi}_{2,k}^{-1}(\tau) \mathfrak{Q}_{\infty,k}(\tau) \! + \! \mathcal{O} 
\! \left(\dfrac{1}{\widetilde{\Lambda}^{3}} \mathfrak{Q}_{\infty,k}^{-1}(\tau) 
\hat{\psi}_{3,k}^{-1}(\tau) \mathfrak{Q}_{\infty,k}(\tau) \right) \right) \nonumber \\
\times& \, \left(\mathrm{I} \! + \! \mathcal{O} \! \left(\lvert \nu (k) \! + \! 1 \rvert^{2} 
\lvert p_{k}(\tau) \rvert^{-2} \tau^{-\epsilon_{\mathrm{\scriptscriptstyle TP}}(k)} 
\mathfrak{Q}_{\infty,k}^{-1}(\tau) \tilde{\mathfrak{C}}_{k}(\tau) \mathfrak{Q}_{\infty,k}
(\tau) \right) \right) \nonumber \\
\times& \, \left(\mathrm{I} \! + \! \widetilde{\Lambda} \mathfrak{Q}_{\infty,k}^{-1}
(\tau) \gimel_{\scriptscriptstyle A,k}^{-1}(\tau) \mathfrak{Q}_{\infty,k}(\tau) \! + \! 
\widetilde{\Lambda}^{2} \mathfrak{Q}_{\infty,k}^{-1}(\tau) 
\gimel_{\scriptscriptstyle B,k}^{-1}(\tau) \mathfrak{Q}_{\infty,k}(\tau) \right) 
\nonumber \\
\times& \, \left(\mathrm{I} \! + \! \widetilde{\Lambda} \tau^{-1/3} 
\mathbb{P}_{\infty,k}(\tau) \! + \! \dfrac{1}{\widetilde{\Lambda}} 
\widehat{\mathbb{E}}_{\infty,k}(\tau) \! + \! \mathcal{O} \! \left((\tau^{-1/3} 
\widetilde{\Lambda})^{2} \widetilde{\mathbb{E}}_{\infty,k}(\tau) \right) \right),
\end{align}}
where $\mathrm{M}_{2}(\mathbb{C}) \! \ni \! \mathcal{R}_{m_{\infty}}(k)$, $m_{\infty} \! \in \! \lbrace -1,0,1,2 \rbrace$, are 
defined in Remark \ref{leedasyue},
\begin{gather}
\mathcal{P}_{0}^{\ast} \! := \! \dfrac{1}{2} \mu_{k}(\tau) \widetilde{\Lambda}^{2} 
\! - \! (\nu (k) \! + \! 1) \ln \widetilde{\Lambda} \! - \! (\nu (k) \! + \! 1) 
\ln (2 \mu_{k}(\tau))^{1/2}, \label{ellinfk20} \\
\mathfrak{Q}_{\infty,k}(\tau) \! := \! \mathbb{F}_{k}^{-1}(\tau) \! \left(
\left(\dfrac{(\varepsilon b)^{1/4}(\sqrt{3} \! + \! 1)^{1/2}}{2^{1/4} 
\sqrt{\smash[b]{\mathfrak{B}_{k}}} \sqrt{\smash[b]{b(\tau)}}} \right)^{\sigma_{3}} 
\mi \sigma_{2} \! + \! \mathfrak{B}_{k}^{-\frac{1}{2} \sigma_{3}} \Delta G_{k}^{\infty}
(\tau)(b(\tau))^{\frac{1}{2} \sigma_{3}} \right), \label{ellinfk21}
\end{gather}
with $\Delta G_{k}^{\infty}(\tau)$ defined by equation \eqref{ellinfk8},
\begin{gather}
\hat{\psi}_{1,k}^{-1}(\tau) \! := \! \dfrac{1}{2 \mu_{k}(\tau)} \! 
\begin{pmatrix}
0 & p_{k}(\tau) \\
-q_{k}(\tau) & 0
\end{pmatrix}, \label{ellinfk22} \\
\hat{\psi}_{2,k}^{-1}(\tau) \! := \! \dfrac{(\nu (k) \! + \! 1)}{4 \mu_{k}(\tau)} \! 
\begin{pmatrix}
1 \! - \! (\nu (k) \! + \! 1) & 0 \\
0 & 1 \! + \! (\nu (k) \! + \! 1)
\end{pmatrix}, \label{ellinfk23} \\
\hat{\psi}_{3,k}^{-1}(\tau) \! := \! -\dfrac{1}{8(\mu_{k}(\tau))^{2}} \! 
\begin{pmatrix}
0 & (1 \! - \! (\nu (k) \! + \! 1))(2 \! - \! (\nu (k) \! + \! 1))p_{k}(\tau) \\
(1 \! + \! (\nu (k) \! + \! 1))(2 \! + \! (\nu (k) \! + \! 1))q_{k}(\tau) & 0 
\end{pmatrix}, \label{ellinfk24} \\
\mathbb{P}_{\infty,k}(\tau) \! := \! (b(\tau))^{-\frac{1}{2} \ad (\sigma_{3})} \! 
\begin{pmatrix}
0 & -\frac{(\varepsilon b)^{1/2}}{3 \sqrt{2} \alpha_{k}} \\
\frac{(\varepsilon b)^{-1/2}}{3 \sqrt{2} \alpha_{k}} & 0
\end{pmatrix}, \label{ellinfk25} \\
\widehat{\mathbb{E}}_{\infty,k}(\tau) \! := \! \dfrac{1}{2 \sqrt{3}(\sqrt{3} 
\! + \! 1)}(b(\tau))^{-\frac{1}{2} \ad (\sigma_{3})} \! 
\begin{pmatrix}
\sqrt{3} \! + \! 1 & -(2 \varepsilon b)^{1/2} \\
(2/\varepsilon b)^{1/2} & \sqrt{3} \! + \! 1
\end{pmatrix} \! 
\begin{pmatrix}
\mathbb{T}_{11,k}(1;\tau) & \mathbb{T}_{12,k}(1;\tau) \\
\mathbb{T}_{21,k}(1;\tau) & \mathbb{T}_{22,k}(1;\tau)
\end{pmatrix}, \label{ellinfk26} \\
\widetilde{\mathbb{E}}_{\infty,k}(\tau) \! := \! \dfrac{1}{2 \sqrt{3}(\sqrt{3} \! + \! 1)}(b(\tau))^{-\frac{1}{2} \ad (\sigma_{3})} \! 
\begin{pmatrix}
\sqrt{3} \! + \! 1 & -(2 \varepsilon b)^{1/2} \\
(2/\varepsilon b)^{1/2} & \sqrt{3} \! + \! 1
\end{pmatrix} \! \tilde{\mathfrak{C}}_{k}^{\lozenge}, \label{ellinfk27}
\end{gather}
$\mathrm{M}_{2}(\mathbb{C}) \! \ni \! \tilde{\mathfrak{C}}_{k}(\tau) \genfrac{}{}{0pt}{3}{=}{\tau \to +\infty} \mathcal{O}(1)$, 
$(\mathbb{T}_{ij,k}(1;\tau))_{i,j=1,2}$ defined in Proposition \ref{prop3.1.6}, and $\mathrm{M}_{2}(\mathbb{C}) \! \ni \! 
\tilde{\mathfrak{C}}_{k}^{\lozenge}$ is $\mathcal{O}(1)$.

Recalling the definitions \eqref{ellinfk12} and \eqref{ellinfk16}, and substituting the expansions \eqref{ellinfk17}, 
\eqref{ellinfk18}, and \eqref{ellinfk19} into equation \eqref{ellinfk11}, one shows, via the conditions \eqref{iden5}, 
the definition \eqref{prpr1}, the restrictions \eqref{pc4}, the asymptotics \eqref{prcyzeek1}, \eqref{prcychik1}, 
and \eqref{prcymuk1}, and (cf. step \pmb{(xi)} in the proof of Lemma \ref{nprcl}) $\arg (\mu_{k}(\tau)) 
\genfrac{}{}{0pt}{3}{=}{\tau \to +\infty} \tfrac{\pi}{2} \big(1 \! + \! \mathcal{O}(\tau^{-2/3}) \big)$, and the restriction \eqref{restr1}, that
\begin{align} \label{ellinfk28} 
\mathfrak{L}^{\infty}_{k}(\tau) \underset{\tau \to +\infty}{=}& \, \mi (\mathcal{R}_{m_{\infty}}(k))^{-1} 
\me^{\tilde{\mathfrak{z}}_{k}^{0}(\tau) \sigma_{3}} \! \left(\dfrac{(\varepsilon b)^{1/4}(\sqrt{3} \! + \! 1)^{1/2}}{2^{1/4} 
\sqrt{\smash[b]{\mathfrak{B}_{k}}} \sqrt{\smash[b]{b(\tau)}}} \right)^{\sigma_{3}} \sigma_{2} 
\me^{-\Delta \tilde{\mathfrak{z}}_{k}(\tau) \sigma_{3}} \nonumber \\
\times& \, \operatorname{diag} \! \Big(\hat{\mathbb{B}}_{0}^{\infty}(\tau),\hat{\mathbb{A}}_{0}^{\infty}(\tau) \Big) 
\overset{\Yup}{\mathbb{E}}_{\scriptscriptstyle \mathfrak{L}^{\infty}_{k}}^{\raise-6.75pt\hbox{$\scriptstyle \leftslice$}}
(\tau), \quad k \! = \! \pm 1,
\end{align}
where $\tilde{\mathfrak{z}}_{k}^{0}(\tau)$, $\Delta \tilde{\mathfrak{z}}_{k}(\tau)$, $\hat{\mathbb{A}}_{0}^{\infty}(\tau)$, 
and $\hat{\mathbb{B}}_{0}^{\infty}(\tau)$ are defined by equations \eqref{ellinfk4}--\eqref{ellinfk7}, respectively, and
\begin{align} \label{ellinfk29} 
\overset{\Yup}{\mathbb{E}}_{\scriptscriptstyle \mathfrak{L}^{\infty}_{k}}^{\raise-6.75pt\hbox{$\scriptstyle \leftslice$}}(\tau) 
\underset{\tau \to +\infty}{:=}& \, \left(\mathrm{I} \! + \! \mathcal{O}(\tau^{-1/3} \widetilde{\Lambda}^{3} \sigma_{3}) \right) \! 
\left(\mathrm{I} \! + \! \mathcal{O} \! \left(\frac{\hat{\mathbb{D}}_{0}^{\infty}(\tau)}{\hat{\mathbb{B}}_{0}^{\infty}(\tau)b(\tau) 
\me^{2 \widetilde{\beta}_{k}^{\ast}(\tau)}} \sigma_{+} \right) \! + \! \mathcal{O} \! \left(\frac{\hat{\mathbb{C}}_{0}^{\infty}
(\tau)b(\tau) \me^{2 \widetilde{\beta}_{k}^{\ast}(\tau)}}{\hat{\mathbb{A}}_{0}^{\infty}(\tau)} \sigma_{-} \right) \right) \nonumber \\ 
\times& \, \left(\mathrm{I} \! + \! \dfrac{1}{\widetilde{\Lambda}} \hat{\psi}_{1,k}^{-1,\sharp}(\tau) \! + \! 
\dfrac{1}{\widetilde{\Lambda}^{2}} \hat{\psi}_{2,k}^{-1,\sharp}(\tau) \! + \! \mathcal{O} \! \left(
\dfrac{1}{\widetilde{\Lambda}^{3}} \hat{\psi}_{3,k}^{-1,\sharp}(\tau) \right) \right) \nonumber \\
\times& \, \left(\mathrm{I} \! + \! \mathcal{O} \! \left(\dfrac{\lvert \nu (k) \! + \! 1 \rvert^{2} 
\tau^{-\epsilon_{\mathrm{\scriptscriptstyle TP}}(k)}}{\lvert p_{k}(\tau) \rvert^{2}} \me^{-\tilde{\beta}_{k}(\tau) \ad (\sigma_{3})} 
\mathfrak{Q}_{\infty,k}^{-1}(\tau) \tilde{\mathfrak{C}}_{k}(\tau) \mathfrak{Q}_{\infty,k}(\tau) \right) \right) \nonumber \\
\times& \, \left(\mathrm{I} \! + \! \widetilde{\Lambda} \gimel_{\scriptscriptstyle A,k}^{\sharp}(\tau) \! + \! 
\widetilde{\Lambda}^{2} \gimel_{\scriptscriptstyle B,k}^{\sharp}(\tau) \right) \! \left(\mathrm{I} \! + \! \widetilde{\Lambda} 
\tau^{-1/3} \mathbb{P}_{\infty,k}^{\sharp}(\tau) \! + \! \dfrac{1}{\widetilde{\Lambda}} \widehat{\mathbb{E}}_{\infty,k}^{\sharp}
(\tau) \right. \nonumber \\
+&\left. \, \mathcal{O} \! \left((\tau^{-1/3} \widetilde{\Lambda})^{2} \widetilde{\mathbb{E}}_{\infty,k}^{\sharp}(\tau) \right) \right),
\end{align}
where $\widetilde{\beta}_{k}^{\ast}(\tau) \! := \! \tfrac{\mi a}{6} \ln \tau \! + \! \mi 3 \alpha_{k}^{2} \tau^{2/3}$,
\begin{gather}
\hat{\mathbb{C}}_{0}^{\infty}(\tau) \! := \! (\Delta G_{k}^{\infty}(\tau))_{11}, \label{ellinfk30} \\
\hat{\mathbb{D}}_{0}^{\infty}(\tau) \! := \! (\Delta G_{k}^{\infty}(\tau))_{22} \! - \! 
\dfrac{\mathfrak{A}_{k}}{\mathfrak{B}_{k}} \! \left(\dfrac{\mi 4 \sqrt{3} \mathcal{Z}_{k}}{\chi_{k}(\tau)} 
\! - \! 1 \right) \! \left(\dfrac{(\varepsilon b)^{1/4}(\sqrt{3} \! + \! 1)^{1/2}}{2^{1/4}} \! + \! 
(\Delta G_{k}^{\infty}(\tau))_{12} \right), \label{ellinfk31} \\
\hat{\psi}_{m,k}^{-1,\sharp}(\tau) \! := \! \me^{-\tilde{\beta}_{k}(\tau) \ad (\sigma_{3})} 
\mathfrak{Q}_{\infty,k}^{-1}(\tau) \hat{\psi}_{m,k}^{-1}(\tau) \mathfrak{Q}_{\infty,k}(\tau), \quad 
m \! = \! 1,2,3, \label{ellinfk32} \\
\gimel_{\scriptscriptstyle A,k}^{\sharp}(\tau) \! := \! \me^{-\tilde{\beta}_{k}(\tau) \ad (\sigma_{3})} 
\mathfrak{Q}_{\infty,k}^{-1}(\tau) \gimel_{\scriptscriptstyle A,k}^{-1}(\tau) \mathfrak{Q}_{\infty,k}(\tau), 
\label{ellinfk33} \\
\gimel_{\scriptscriptstyle B,k}^{\sharp}(\tau) \! := \! \me^{-\tilde{\beta}_{k}(\tau) \ad (\sigma_{3})} 
\mathfrak{Q}_{\infty,k}^{-1}(\tau) \gimel_{\scriptscriptstyle B,k}^{-1}(\tau) \mathfrak{Q}_{\infty,k}(\tau), 
\label{ellinfk34} \\
\mathbb{P}_{\infty,k}^{\sharp}(\tau) \! := \! \me^{-\tilde{\beta}_{k}(\tau) \ad (\sigma_{3})} 
\mathbb{P}_{\infty,k}(\tau), \label{ellinfk35} \\
\widehat{\mathbb{E}}_{\infty,k}^{\sharp}(\tau) \! := \! \me^{-\tilde{\beta}_{k}(\tau) \ad (\sigma_{3})} 
\widehat{\mathbb{E}}_{\infty,k}(\tau), \label{ellinfk36} \\
\widetilde{\mathbb{E}}_{\infty,k}^{\sharp}(\tau) \! := \! \me^{-\tilde{\beta}_{k}(\tau) \ad (\sigma_{3})} 
\widetilde{\mathbb{E}}_{\infty,k}(\tau). \label{ellinfk37} 
\end{gather} 
Via the conditions \eqref{iden5}, the restrictions \eqref{pc4} and \eqref{restr1}, the definitions \eqref{eqpeetee}, 
\eqref{peekayity}, \eqref{prpr1}, \eqref{prpr3}, \eqref{prpr4}, \eqref{prcy57}, \eqref{prcy58}, 
\eqref{ellinfk6}--\eqref{ellinfk8}, \eqref{ellinfk13}, \eqref{ellinfk21}--\eqref{ellinfk27}, and 
\eqref{ellinfk30}--\eqref{ellinfk37}, and the asymptotics \eqref{tr1}, \eqref{tr3}, \eqref{asympforf3}, \eqref{prcyzeek1}, 
\eqref{prcyg4}--\eqref{prcybk1}, \eqref{prcyomg1}--\eqref{prcyell2k1}, and \eqref{ellinfk18}, upon imposing the 
conditions \eqref{ellinfk2a} and \eqref{ellinfk2b}, and defining
\begin{gather*}
J_{k}^{\infty} \! := \! \frac{1}{\sqrt{2 \sqrtsign{3}(\sqrtsign{3} \! + \! 1)}} 
\begin{pmatrix}
\sqrt{3} \! + \! 1 & -(2 \varepsilon b)^{1/2} \\
(2/\varepsilon b)^{1/2} & \sqrt{3} \! + \! 1
\end{pmatrix}, \, \, \, \qquad \, \, \, \pmb{\mathbb{T}}_{\infty,k}^{\sharp} \! 
:= \! (\mathbb{T}_{ij,k}(1;\tau))_{i,j=1,2}, \\
\mathbb{D}_{\infty,k}^{\sharp} \! := \! \mathfrak{B}_{k}^{-\frac{1}{2} \sigma_{3}} \! 
\begin{pmatrix}
0 & -\frac{(\varepsilon b)^{1/4}(\sqrt{3}+1)^{1/2}}{2^{1/4}} \\
\frac{2^{1/4}}{(\varepsilon b)^{1/4}(\sqrt{3}+1)^{1/2}} & 0
\end{pmatrix},
\end{gather*}
one shows that, for $k \! = \! \pm 1$,
\begin{align*} 
\overset{\Yup}{\mathbb{E}}_{\scriptscriptstyle 
\mathfrak{L}^{\infty}_{k}}^{\raise-6.75pt\hbox{$\scriptstyle \leftslice$}}(\tau) \underset{\tau \to +\infty}{=}& \, 
\left(\mathrm{I} \! + \! \mathcal{O}(\tau^{-1/3} \widetilde{\Lambda}^{3} \sigma_{3}) \right) \! \left(\mathrm{I} \! + \! \mathcal{O} \! 
\left(\frac{\hat{\mathbb{D}}_{0}^{\infty}(\tau)}{\hat{\mathbb{B}}_{0}^{\infty}(\tau)b(\tau) \me^{2 \widetilde{\beta}_{k}^{\ast}(\tau)}} 
\sigma_{+} \right) \! + \! \mathcal{O} \! \left(\frac{\hat{\mathbb{C}}_{0}^{\infty}(\tau)b(\tau) \me^{2 \widetilde{\beta}_{k}^{\ast}
(\tau)}}{\hat{\mathbb{A}}_{0}^{\infty}(\tau)} \sigma_{-} \right) \right) \nonumber \\ 
\times& \, \left(\mathrm{I} \! + \! \dfrac{1}{\widetilde{\Lambda}} 
\hat{\psi}_{1,k}^{-1,\sharp}(\tau) \! + \! \dfrac{1}{\widetilde{\Lambda}^{2}} 
\hat{\psi}_{2,k}^{-1,\sharp}(\tau) \! + \! \mathcal{O} \! \left(\dfrac{1}{\widetilde{
\Lambda}^{3}} \hat{\psi}_{3,k}^{-1,\sharp}(\tau) \right) \right) \nonumber \\
\times& \, \left(\mathrm{I} \! + \! \mathcal{O} \! \left(\dfrac{\lvert \nu (k) \! + \! 
1 \rvert^{2} \tau^{-\epsilon_{\mathrm{\scriptscriptstyle TP}}(k)}}{\lvert p_{k}
(\tau) \rvert^{2}} \me^{-\tilde{\beta}_{k}(\tau) \ad (\sigma_{3})} 
\mathfrak{Q}_{\infty,k}^{-1}(\tau) \tilde{\mathfrak{C}}_{k}(\tau) 
\mathfrak{Q}_{\infty,k}(\tau) \right) \right) \nonumber \\
\times& \, \left(\mathrm{I} \! + \! \gimel_{\scriptscriptstyle A,k}^{\sharp}(\tau) 
\widehat{\mathbb{E}}_{\infty,k}^{\sharp}(\tau) \! + \! \dfrac{1}{\widetilde{\Lambda}} 
\widehat{\mathbb{E}}_{\infty,k}^{\sharp}(\tau) \! + \! \widetilde{\Lambda} \! 
\left(\tau^{-1/3} \mathbb{P}_{\infty,k}^{\sharp}(\tau) \! + \! 
\gimel_{\scriptscriptstyle A,k}^{\sharp}(\tau) \right. \right. \nonumber \\
+&\left. \left. \, \gimel_{\scriptscriptstyle B,k}^{\sharp}(\tau) \widehat{\mathbb{E}}_{
\infty,k}^{\sharp}(\tau) \right) \! + \! \widetilde{\Lambda}^{2} \! \left(\tau^{-1/3} 
\gimel_{\scriptscriptstyle A,k}^{\sharp}(\tau) \mathbb{P}_{\infty,k}^{\sharp}(\tau) 
\! + \! \gimel_{\scriptscriptstyle B,k}^{\sharp}(\tau) \! + \! \mathcal{O} \! \left(
\tau^{-2/3} \widetilde{\mathbb{E}}_{\infty,k}^{\sharp}(\tau) \right) \right) 
\right. \nonumber \\
+&\left. \, \widetilde{\Lambda}^{3} \! \left(\tau^{-1/3} 
\gimel_{\scriptscriptstyle B,k}^{\sharp}(\tau) \mathbb{P}_{\infty,k}^{\sharp}(\tau) 
\! + \! \mathcal{O} \! \left(\tau^{-2/3} \gimel_{\scriptscriptstyle A,k}^{\sharp}(\tau) 
\widetilde{\mathbb{E}}_{\infty,k}^{\sharp}(\tau) \right) \right) \right) \nonumber \\
\underset{\tau \to +\infty}{=}& \, \left(\mathrm{I} \! + \! \mathcal{O}(\tau^{-1/3} \widetilde{\Lambda}^{3} \sigma_{3}) \right) \! 
\left(\mathrm{I} \! + \! \mathcal{O} \! \left(\frac{\hat{\mathbb{D}}_{0}^{\infty}(\tau)}{\hat{\mathbb{B}}_{0}^{\infty}(\tau)b(\tau) 
\me^{2 \widetilde{\beta}_{k}^{\ast}(\tau)}} \sigma_{+} \right) \! + \! \mathcal{O} \! \left(\frac{\hat{\mathbb{C}}_{0}^{\infty}
(\tau)b(\tau) \me^{2 \widetilde{\beta}_{k}^{\ast}(\tau)}}{\hat{\mathbb{A}}_{0}^{\infty}(\tau)} \sigma_{-} \right) \right) \nonumber \\ 
\times& \, \left(\mathrm{I} \! + \! \dfrac{1}{\widetilde{\Lambda}} 
\hat{\psi}_{1,k}^{-1,\sharp}(\tau) \! + \! \dfrac{1}{\widetilde{\Lambda}^{2}} 
\hat{\psi}_{2,k}^{-1,\sharp}(\tau) \! + \! \mathcal{O} \! \left(\dfrac{1}{\widetilde{
\Lambda}^{3}} \hat{\psi}_{3,k}^{-1,\sharp}(\tau) \right) \right) \nonumber \\
\times& \, \left(\mathrm{I} \! + \! \gimel_{\scriptscriptstyle A,k}^{\sharp}(\tau) 
\widehat{\mathbb{E}}_{\infty,k}^{\sharp}(\tau) \! + \! \dfrac{1}{\widetilde{\Lambda}} 
\dfrac{1}{2 \sqrt{3}(\sqrt{3} \! + \! 1)} \! \left(\dfrac{\me^{-\tilde{\beta}_{k}(\tau)}}{
\sqrt{\smash[b]{b(\tau)}}} \right)^{\ad (\sigma_{3})} J_{k}^{\infty} 
\pmb{\mathbb{T}}_{\infty,k}^{\sharp} \right. \nonumber \\
+&\left. \, \widetilde{\Lambda} \dfrac{\mi 4 \sqrt{3} \mathcal{Z}_{k} \mathfrak{A}_{k} 
\ell_{0,k}^{+}}{\chi_{k}(\tau)} \sigma_{3} \! + \! \mathcal{O} \! \left(\dfrac{\lvert \nu 
(k) \! + \! 1 \rvert^{2} \tau^{-\epsilon_{\mathrm{\scriptscriptstyle TP}}(k)}}{\lvert 
p_{k}(\tau) \rvert^{2}} \! \left(\dfrac{\me^{-\tilde{\beta}_{k}(\tau)}}{\sqrt{\smash[b]{
b(\tau)}}} \right)^{\ad (\sigma_{3})} \right. \right. \nonumber \\
\times&\left. \left. \, \mathbb{D}_{\infty,k}^{\sharp} \tilde{\mathfrak{C}}_{k}
(\tau)(\mathbb{D}_{\infty,k}^{\sharp})^{-1} \right) \right) \nonumber \\
\underset{\tau \to +\infty}{=}& \, \left(\mathrm{I} \! + \! \mathcal{O}(\tau^{-1/3} \widetilde{\Lambda}^{3} \sigma_{3}) \right) \! 
\left(\mathrm{I} \! + \! \mathcal{O} \! \left(\frac{\hat{\mathbb{D}}_{0}^{\infty}(\tau)}{\hat{\mathbb{B}}_{0}^{\infty}(\tau)b(\tau) 
\me^{2 \widetilde{\beta}_{k}^{\ast}(\tau)}} \sigma_{+} \right) \! + \! \mathcal{O} \! \left(\frac{\hat{\mathbb{C}}_{0}^{\infty}
(\tau)b(\tau) \me^{2 \widetilde{\beta}_{k}^{\ast}(\tau)}}{\hat{\mathbb{A}}_{0}^{\infty}(\tau)} \sigma_{-} \right) \right) \nonumber \\ 
\times& \, \left(\mathrm{I} \! + \! \gimel_{\scriptscriptstyle A,k}^{\sharp}(\tau) 
\widehat{\mathbb{E}}_{\infty,k}^{\sharp}(\tau) \! + \! \dfrac{\mi 4 \sqrt{3} 
\mathcal{Z}_{k} \mathfrak{A}_{k} \ell_{0,k}^{+}}{\chi_{k}(\tau)} \hat{\psi}_{1,k}^{
-1,\sharp}(\tau) \sigma_{3} \! + \! \widetilde{\Lambda} \dfrac{\mi 4 \sqrt{3} 
\mathcal{Z}_{k} \mathfrak{A}_{k} \ell_{0,k}^{+}}{\chi_{k}(\tau)} \sigma_{3} 
\right. \nonumber \\
+&\left. \, \dfrac{1}{\widetilde{\Lambda}} \! \left(\hat{\psi}_{1,k}^{-1,\sharp}
(\tau) \! + \! \dfrac{\mi 4 \sqrt{3} \mathcal{Z}_{k} \mathfrak{A}_{k} 
\ell_{0,k}^{+}}{\chi_{k}(\tau)} \hat{\psi}_{2,k}^{-1,\sharp}(\tau) \sigma_{3} 
\! + \! \dfrac{1}{2 \sqrt{3}(\sqrt{3} \! + \! 1)} \! \left(\dfrac{\me^{-\tilde{\beta}_{k}
(\tau)}}{\sqrt{\smash[b]{b(\tau)}}} \right)^{\ad (\sigma_{3})} \right. \right. 
\nonumber \\
\times&\left. \left. J_{k}^{\infty} \pmb{\mathbb{T}}_{\infty,k}^{\sharp} 
\right) \! + \! \dfrac{1}{\widetilde{\Lambda}^{2}} \! \left(\hat{\psi}_{2,k}^{-1,\sharp}(\tau) \! + \! 
\dfrac{1}{2 \sqrt{3}(\sqrt{3} \! + \! 1)} \hat{\psi}_{1,k}^{-1,\sharp}(\tau) \! \left(
\dfrac{\me^{-\tilde{\beta}_{k}(\tau)}}{\sqrt{\smash[b]{b(\tau)}}} \right)^{\ad 
(\sigma_{3})} J_{k}^{\infty} \pmb{\mathbb{T}}_{\infty,k}^{\sharp} \right) \right. \nonumber \\
+&\left. \, \mathcal{O} \! \left(\dfrac{\lvert \nu (k) \! + \! 1 \rvert^{2} \tau^{-
\epsilon_{\mathrm{\scriptscriptstyle TP}}(k)}}{\lvert p_{k}(\tau) \rvert^{2}} \! 
\left(\dfrac{\me^{-\tilde{\beta}_{k}(\tau)}}{\sqrt{\smash[b]{b(\tau)}}} \right)^{
\ad (\sigma_{3})} \mathbb{D}_{\infty,k}^{\sharp} \tilde{\mathfrak{C}}_{k}(\tau)
(\mathbb{D}_{\infty,k}^{\sharp})^{-1} \right) \right. \nonumber \\
+&\left. \, \mathcal{O} \! \left(\dfrac{1}{\widetilde{\Lambda}} \dfrac{\lvert \nu (k) 
\! + \! 1 \rvert^{2} \tau^{-\epsilon_{\mathrm{\scriptscriptstyle TP}}(k)}}{\lvert 
p_{k}(\tau) \rvert^{2}} \hat{\psi}_{1,k}^{-1,\sharp}(\tau) \! \left(\dfrac{\me^{-
\tilde{\beta}_{k}(\tau)}}{\sqrt{\smash[b]{b(\tau)}}} \right)^{\ad (\sigma_{3})} 
\mathbb{D}_{\infty,k}^{\sharp} \tilde{\mathfrak{C}}_{k}(\tau)(\mathbb{D}_{
\infty,k}^{\sharp})^{-1} \right) \right. \nonumber \\
+&\left. \, \mathcal{O} \! \left(\dfrac{1}{\widetilde{\Lambda}^{2}} \dfrac{\lvert \nu 
(k) \! + \! 1 \rvert^{2} \tau^{-\epsilon_{\mathrm{\scriptscriptstyle TP}}(k)}}{\lvert 
p_{k}(\tau) \rvert^{2}} \hat{\psi}_{2,k}^{-1,\sharp}(\tau) \! \left(\dfrac{\me^{-
\tilde{\beta}_{k}(\tau)}}{\sqrt{\smash[b]{b(\tau)}}} \right)^{\ad (\sigma_{3})} 
\mathbb{D}_{\infty,k}^{\sharp} \tilde{\mathfrak{C}}_{k}(\tau)(\mathbb{D}_{
\infty,k}^{\sharp})^{-1} \right) \right. \nonumber \\
+&\left. \, \mathcal{O} \! \left(\dfrac{1}{\widetilde{\Lambda}^{2}} 
\dfrac{\mi 4 \sqrt{3} \mathcal{Z}_{k} \mathfrak{A}_{k} \ell_{0,k}^{+}}{\chi_{k}(\tau)} 
\hat{\psi}_{3,k}^{-1,\sharp}(\tau) \sigma_{3} \right) \right) \nonumber \\
\underset{\tau \to +\infty}{=}& \, \mathrm{I} \! + \! \gimel_{\scriptscriptstyle A,k}^{\sharp}(\tau) 
\widehat{\mathbb{E}}_{\infty,k}^{\sharp}(\tau) \! + \! \dfrac{\mi 4 \sqrt{3} \mathcal{Z}_{k} 
\mathfrak{A}_{k} \ell_{0,k}^{+}}{\chi_{k}(\tau)} \hat{\psi}_{1,k}^{-1,\sharp}(\tau) \sigma_{3} 
\! + \! \mathcal{O} \big(\tau^{-\frac{1}{3}+3 \delta_{k}} \sigma_{3} \big) \nonumber \\
+& \, 
\begin{pmatrix}
0 & \mathcal{O}(\tau^{-2/3}) \\
\mathcal{O}(\tau^{-2/3}) & 0
\end{pmatrix} \! + \! 
\begin{pmatrix}
\mathcal{O}(\tau^{-\frac{1}{3}+\delta_{k}}) & 0 \\
0 & \mathcal{O}(\tau^{-\frac{1}{3}+\delta_{k}})
\end{pmatrix} \nonumber \\
+& \, 
\begin{pmatrix}
0 & \mathcal{O}(\tau^{-\delta_{k}}(\nu (k) \! + \! 1)^{\frac{1+k}{2}}) \\
\mathcal{O}(\tau^{-\delta_{k}}(\nu (k) \! + \! 1)^{\frac{1-k}{2}}) & 0
\end{pmatrix} \nonumber \\
+& \, 
\begin{pmatrix}
\mathcal{O}(\tau^{-\frac{1}{3}-\delta_{k}}(\nu (k) \! + \! 1)) & 0 \\
0 & \mathcal{O}(\tau^{-\frac{1}{3}-\delta_{k}}(\nu (k) \! + \! 1))
\end{pmatrix} \nonumber \\
+& \, 
\begin{pmatrix}
\mathcal{O}(\tau^{-\frac{1}{3}-\delta_{k}}) & \mathcal{O}(\tau^{-\frac{1}{3}
-\delta_{k}}) \\
\mathcal{O}(\tau^{-\frac{1}{3}-\delta_{k}}) & \mathcal{O}(\tau^{-\frac{1}{3}
-\delta_{k}})
\end{pmatrix} \! + \! 
\begin{pmatrix}
\mathcal{O}(\tau^{-2\delta_{k}}(\nu (k) \! + \! 1)) & 0 \\
0 & \mathcal{O}(\tau^{-2\delta_{k}}(\nu (k) \! + \! 1))
\end{pmatrix} \nonumber \\
+& \, 
\begin{pmatrix}
\mathcal{O}(\tau^{-\frac{1}{3}-2\delta_{k}}(\nu (k) \! + \! 1)^{\frac{1+k}{2}}) 
& \mathcal{O}(\tau^{-\frac{1}{3}-2\delta_{k}}(\nu (k) \! + \! 1)^{\frac{1+k}{2}}) \\
\mathcal{O}(\tau^{-\frac{1}{3}-2\delta_{k}}(\nu (k) \! + \! 1)^{\frac{1-k}{2}}) 
& \mathcal{O}(\tau^{-\frac{1}{3}-2\delta_{k}}(\nu (k) \! + \! 1)^{\frac{1-k}{2}}) 
\end{pmatrix} \nonumber \\
+& \, 
\begin{pmatrix}
\mathcal{O}(\tau^{-2-\epsilon_{\mathrm{\scriptscriptstyle TP}}(k)}) & \mathcal{O}
(\tau^{-1-\epsilon_{\mathrm{\scriptscriptstyle TP}}(k)}(\nu (k) \! + \! 1)^{\frac{1
+k}{2}}) \\
\mathcal{O}(\tau^{-3-\epsilon_{\mathrm{\scriptscriptstyle TP}}(k)}(\nu (k) \! + \! 
1)^{\frac{1-k}{2}}) & \mathcal{O}(\tau^{-2-\epsilon_{\mathrm{\scriptscriptstyle TP}}
(k)})
\end{pmatrix} \nonumber \\
+& \, 
\begin{pmatrix}
\mathcal{O}(\tau^{-3-\delta_{k}-\epsilon_{\mathrm{\scriptscriptstyle TP}}(k)}(\nu 
(k) \! + \! 1)) & \mathcal{O}(\tau^{-2-\delta_{k}-\epsilon_{\mathrm{\scriptscriptstyle TP}}
(k)}(\nu (k) \! + \! 1)^{\frac{1+k}{2}}) \\
\mathcal{O}(\tau^{-2-\delta_{k}-\epsilon_{\mathrm{\scriptscriptstyle TP}}(k)}
(\nu (k) \! + \! 1)^{\frac{1-k}{2}}) & \mathcal{O}(\tau^{-1-\delta_{k}-
\epsilon_{\mathrm{\scriptscriptstyle TP}}(k)}(\nu (k) \! + \! 1))
\end{pmatrix} \nonumber \\
+& \, 
\begin{pmatrix}
\mathcal{O}(\tau^{-2-2\delta_{k}-\epsilon_{\mathrm{\scriptscriptstyle TP}}(k)}(\nu (k) 
\! + \! 1)) & \mathcal{O}(\tau^{-1-2\delta_{k}-\epsilon_{\mathrm{\scriptscriptstyle TP}}
(k)}(\nu (k) \! + \! 1)^{\frac{3+k}{2}}) \\
\mathcal{O}(\tau^{-3-2\delta_{k}-\epsilon_{\mathrm{\scriptscriptstyle TP}}(k)}
(\nu (k) \! + \! 1)^{\frac{3-k}{2}}) & \mathcal{O}(\tau^{-2-2\delta_{k}-
\epsilon_{\mathrm{\scriptscriptstyle TP}}(k)}(\nu (k) \! + \! 1))
\end{pmatrix} \nonumber \\
+& \, 
\begin{pmatrix}
0 & \mathcal{O}(\tau^{-\frac{1}{3}-2\delta_{k}}(\nu (k) \! + \! 1)^{\frac{1+k}{2}}) \\
\mathcal{O}(\tau^{-\frac{1}{3}-2\delta_{k}}(\nu (k) \! + \! 1)^{\frac{1-k}{2}}) & 0
\end{pmatrix} \nonumber \\
\underset{\tau \to +\infty}{=}& \, \mathrm{I} \! + \! \gimel_{\scriptscriptstyle A,k}^{\sharp}(\tau) 
\widehat{\mathbb{E}}_{\infty,k}^{\sharp}(\tau) \! + \! \dfrac{\mi 4 \sqrt{3} \mathcal{Z}_{k} \mathfrak{A}_{k} 
\ell_{0,k}^{+}}{\chi_{k}(\tau)} \hat{\psi}_{1,k}^{-1,\sharp}(\tau) \sigma_{3} \! + \! 
\begin{pmatrix}
\mathcal{O} \big(\tau^{-\frac{1}{3}+3 \delta_{k}} \big) & \mathcal{O} \big(\tau^{-2/3} \big) \\
\mathcal{O} \big(\tau^{-2/3} \big) & \mathcal{O} \big(\tau^{-\frac{1}{3}+3 \delta_{k}} \big)
\end{pmatrix} \nonumber \\
+& \, 
\begin{pmatrix}
\mathcal{O} \big(\tau^{-\frac{1}{3}+\delta_{k}} \big) & \mathcal{O} \big(\tau^{-\frac{1}{3}(\frac{1+k}{2})-\delta_{k}} \big) \\
\mathcal{O} \big(\tau^{-\frac{1}{3}(\frac{1-k}{2})-\delta_{k}} \big) & \mathcal{O} \big(\tau^{-\frac{1}{3}+\delta_{k}} \big)
\end{pmatrix} \nonumber \\
\underset{\tau \to +\infty}{=}& \, \mathrm{I} \! + \! 
\underbrace{\gimel_{\scriptscriptstyle A,k}^{\sharp}(\tau) \widehat{\mathbb{E}}_{\infty,k}^{\sharp}(\tau) 
\! + \! \dfrac{\mi 4 \sqrt{3} \mathcal{Z}_{k} \mathfrak{A}_{k} \ell_{0,k}^{+}}{\chi_{k}(\tau)} \hat{\psi}_{1,k}^{-1,\sharp}(\tau) 
\sigma_{3}}_{=: \, \mathbb{E}_{{\scriptscriptstyle \mathcal{N}},k}^{\infty}(\tau)} \nonumber \\
+& \, \underbrace{\begin{pmatrix}
\mathcal{O}(\tau^{-\frac{1}{3}+3\delta_{k}}) & \mathcal{O}(\tau^{-\frac{1}{3}(\frac{1+k}{2})-\delta_{k}}) \\
\mathcal{O}(\tau^{-\frac{1}{3}(\frac{1-k}{2})-\delta_{k}}) & \mathcal{O}(\tau^{-\frac{1}{3}+3\delta_{k}})
\end{pmatrix}}_{=: \, \mathcal{O}(\mathbb{E}^{\infty}_{k}(\tau))} \nonumber \\
\underset{\tau \to +\infty}{=}& \, (\mathrm{I} \! + \! \mathbb{E}_{{\scriptscriptstyle \mathcal{N}},k}^{\infty}(\tau)) \! 
\left(\mathrm{I} \! + \! \underbrace{(\mathrm{I} \! + \! \mathbb{E}_{{\scriptscriptstyle \mathcal{N}},k}^{\infty}(\tau))^{-1}}_{= \, 
\mathcal{O}(1)} \, \mathcal{O}(\mathbb{E}^{\infty}_{k}(\tau)) \right) \quad \Rightarrow \nonumber
\end{align*}
\begin{equation} \label{ellinfk38} 
\overset{\Yup}{\mathbb{E}}_{\scriptscriptstyle \mathfrak{L}^{\infty}_{k}}^{\raise-6.75pt\hbox{$\scriptstyle \leftslice$}}
(\tau) \underset{\tau \to +\infty}{=} \big(\mathrm{I} \! + \! \mathbb{E}_{{\scriptscriptstyle \mathcal{N}},k}^{\infty}
(\tau) \big) \big(\mathrm{I} \! + \! \mathcal{O}(\mathbb{E}^{\infty}_{k}(\tau)) \big),
\end{equation} 
where $\mathbb{E}_{{\scriptscriptstyle \mathcal{N}},k}^{\infty}(\tau)$ and $\mathcal{O}(\mathbb{E}^{\infty}_{k}(\tau))$ 
are defined by equations \eqref{ellinfk9} and \eqref{ellinfk10}, respectively.\footnote{The asymptotics for the function 
$\mathbb{E}_{{\scriptscriptstyle \mathcal{N}},k}^{\infty}(\tau)$ is presented in the proof of Lemma \ref{ginversion} 
(see Section \ref{finalsec}).} Thus, via the asymptotics \eqref{ellinfk28} and \eqref{ellinfk38}, one arrives at the results 
stated in the lemma. \hfill $\qed$
\begin{cccc} \label{lzernewlemm} 
Let $\widetilde{\Psi}_{k}(\widetilde{\mu},\tau)$, $k \! = \! \pm 1$, be the fundamental solution of equation \eqref{eq3.3} 
with asymptotics given in Lemma \ref{nprcl}, and let $\mathbb{X}^{0}_{1-k}(\widetilde{\mu},\tau)$ be the canonical 
solution of equation \eqref{eq3.1}.\footnote{See Proposition \ref{prop1.4}.} Define {}\footnote{Since (cf. equations 
\eqref{eq3.2}) $\tau^{-\frac{1}{12} \sigma_{3}} \mathbb{X}^{0}_{1-k}(\tau^{-1/6} \widetilde{\mu},\tau)$, $k \! = \! \pm 1$, 
is also a fundamental solution of equation \eqref{eq3.3}, it follows, therefore, that $\mathfrak{L}^{0}_{k}(\tau)$ is 
independent of $\widetilde{\mu}$.}
\begin{equation} \label{ellohk1} 
\mathfrak{L}^{0}_{k}(\tau) \! := \! \big(\widetilde{\Psi}_{k}(\widetilde{\mu},\tau) \big)^{-1} \tau^{-\frac{1}{12} \sigma_{3}} 
\mathbb{X}^{0}_{1-k}(\tau^{-1/6} \widetilde{\mu},\tau), \quad k \! = \! \pm 1.
\end{equation}
Assume that the parameters $\nu (k) \! + \! 1$ and $\delta_{k}$ satisfy the restrictions \eqref{pc4} and \eqref{restr1}, respectively, 
and, additionally, the conditions \eqref{ellinfk2a} and \eqref{ellinfk2b} are valid. Then,
\begin{align} \label{ellohk3} 
\mathfrak{L}^{0}_{k}(\tau) \underset{\tau \to +\infty}{=}& \, (\mathcal{R}_{m_{0}}(k))^{-1} \me^{\hat{\mathfrak{z}}_{k}^{0}
(\tau) \sigma_{3}} \! \left(\dfrac{\mi 2^{1/4}}{(\sqrt{3} \! - \! 1)^{1/2} \sqrt{\smash[b]{\mathfrak{B}_{k}}}} 
\right)^{\sigma_{3}} \me^{\Delta \hat{\mathfrak{z}}_{k}(\tau) \sigma_{3}} 
\begin{pmatrix}
\hat{\mathbb{A}}_{0}^{0}(\tau) & 0 \\
0 & \hat{\mathbb{B}}_{0}^{0}(\tau)
\end{pmatrix} \nonumber \\
\times& \, \left(\mathrm{I} \! + \! \mathbb{E}_{{\scriptscriptstyle \mathcal{N}},k}^{0}(\tau) \right) \! \mathbb{S}_{k}^{\ast} 
\! \left(\mathrm{I} \! + \! \mathcal{O} \left(\mathbb{E}^{0}_{k}(\tau) \right) \right),
\end{align}
where $\mathrm{M}_{2}(\mathbb{C}) \! \ni \! \mathcal{R}_{m_{0}}(k)$, $m_{0} \! \in \! \lbrace -1,0,1,2 \rbrace$, are 
defined in Remark \ref{leedasyue},
\begin{align} 
\hat{\mathfrak{z}}_{k}^{0}(\tau) \! :=& \, \mi \tau^{2/3}3 \sqrt{3} \alpha_{k}^{2} \! + \! \mi (a \! - \! \mi/2) \ln (2^{-1/2}
(\sqrt{3} \! + \! 1)), \label{ellohk4} \\
\Delta \hat{\mathfrak{z}}_{k}(\tau) \! :=& \, -\left(\dfrac{5 \! + \! 9 \sqrt{3}}{6 \sqrt{3} \alpha_{k}^{2}} \right) \! 
\mathfrak{p}_{k}(\tau) \! + \! (\nu (k) \! + \! 1) \ln (2 \mu_{k}(\tau))^{1/2} \! + \! \dfrac{1}{3}(\nu (k) \! + \! 1) \ln \tau 
\nonumber \\
-& \, (\nu (k) \! + \! 1) \ln (\me^{\mi k \pi}/3 \alpha_{k}), \label{ellohk5}
\end{align}
with $\mathfrak{p}_{k}(\tau)$ defined by equation \eqref{eqpeetee}, and $\mathfrak{B}_{k}$ and $\mu_{k}(\tau)$ defined 
in Lemma \ref{nprcl},
\begin{gather}
\hat{\mathbb{A}}_{0}^{0}(\tau) \! := \! 1 \! + \! \dfrac{(\varepsilon b)^{1/4}
(\sqrt{3} \! - \! 1)^{1/2}(\Delta G_{k}^{0}(\tau))_{11}}{2^{1/4}}, \label{ellohk6} \\
\hat{\mathbb{B}}_{0}^{0}(\tau) \! := \! 1 \! + \! \dfrac{2^{1/4}}{(\varepsilon 
b)^{1/4}(\sqrt{3} \! - \! 1)^{1/2}} \! \left((\Delta G_{k}^{0}(\tau))_{22} \! - 
\! \dfrac{\mathfrak{A}_{k}}{\mathfrak{B}_{k}} \! \left(\dfrac{\mi 4 \sqrt{3} 
\mathcal{Z}_{k}}{\chi_{k}(\tau)} \! - \! 1 \right) \! (\Delta G_{k}^{0}(\tau))_{12} \right), 
\label{ellohk7}
\end{gather}
with $\mathcal{Z}_{k}$, $\mathfrak{A}_{k}$, and $\chi_{k}(\tau)$ defined in Lemma \ref{nprcl}, and
\begin{equation} \label{ellohk8} 
\Delta G_{k}^{0}(\tau) \! := \! \dfrac{1}{(2 \sqrt{3}(\sqrt{3} \! - \! 1))^{1/2}} 
\begin{pmatrix}
(\Delta G_{k}^{0}(\tau))_{11} & (\Delta G_{k}^{0}(\tau))_{12} \\
(\Delta G_{k}^{0}(\tau))_{21} & (\Delta G_{k}^{0}(\tau))_{22}
\end{pmatrix},
\end{equation}
with
\begin{equation} \label{forzilch}
\begin{aligned} 
(\Delta G_{k}^{0}(\tau))_{11} \! =& \, (\sqrt{3} \! - \! 1)(\Delta \mathcal{G}_{0,k})_{22} \! - \! (2/\varepsilon b)^{1/2}
(\Delta \mathcal{G}_{0,k})_{12}, \\
(\Delta G_{k}^{0}(\tau))_{12} \! =& -(\sqrt{3} \! - \! 1)(\Delta \mathcal{G}_{0,k})_{12} \! - \! (2 \varepsilon b)^{1/2}
(\Delta \mathcal{G}_{0,k})_{22}, \\
(\Delta G_{k}^{0}(\tau))_{21} \! =& -(\sqrt{3} \! - \! 1)(\Delta \mathcal{G}_{0,k})_{21} \! + \! (2/\varepsilon b)^{1/2}
(\Delta \mathcal{G}_{0,k})_{11}, \\
(\Delta G_{k}^{0}(\tau))_{22} \! =& \, (\sqrt{3} \! - \! 1)(\Delta \mathcal{G}_{0,k})_{11} \! + \! (2 \varepsilon b)^{1/2}
(\Delta \mathcal{G}_{0,k})_{21},
\end{aligned}
\end{equation}
where $(\Delta \mathcal{G}_{0,k})_{i,j=1,2}$ are defined by equations \eqref{prcyg4}--\eqref{prcyg6},
\begin{equation} \label{ellohk9} 
\mathbb{S}_{k}^{\ast} \! := \! 
\begin{pmatrix}
1 & -(1 \! + \! k)s_{0}^{0}/2 \\
(1 \! - \! k)s_{0}^{0}/2 & 1
\end{pmatrix},
\end{equation}
\begin{align} \label{ellohk10} 
\mathbb{E}_{{\scriptscriptstyle \mathcal{N}},k}^{0}(\tau) :=& \, 
\me^{-\hat{\beta}_{k}(\tau) \ad (\sigma_{3})} \! \left(\dfrac{\mi 4 \sqrt{3} 
\mathcal{Z}_{k} \mathfrak{A}_{k} \ell_{0,k}^{+}}{\chi_{k}(\tau)} \! \left(
\dfrac{(\sqrt{3} \! - \! 1)p_{k}(\tau) \mathfrak{B}_{k}}{2^{3/2} \mu_{k}(\tau)} 
\sigma_{+} \! + \! \dfrac{\sqrt{2}(\nu (k) \! + \! 1)}{(\sqrt{3} \! - \! 1)p_{k}
(\tau) \mathfrak{B}_{k}} \sigma_{-} \right) \! \sigma_{3} \right. \nonumber \\
+&\left. \, \dfrac{1}{2 \sqrt{3}(\sqrt{3} \! - \! 1)} \! 
\begin{pmatrix}
-\frac{\mi 4 \sqrt{3} \mathcal{Z}_{k} \mathfrak{A}_{k} \ell_{0,k}^{+}}{\chi_{k}(\tau)} & \frac{(\sqrt{3}-1) \mathfrak{B}_{k} 
\ell_{0,k}^{+}}{\sqrt{2}} \\
-\frac{\sqrt{2}}{(\sqrt{3}-1) \mathfrak{B}_{k}} \big((\frac{\mi 4 \sqrt{3} \mathcal{Z}_{k} \mathfrak{A}_{k}}{\chi_{k}(\tau)})^{2} 
\ell_{0,k}^{+} \! - \! \ell_{1,k}^{+} \! - \! \ell_{2,k}^{+} \big) & \frac{\mi 4 \sqrt{3} \mathcal{Z}_{k} \mathfrak{A}_{k} 
\ell_{0,k}^{+}}{\chi_{k}(\tau)} 
\end{pmatrix} \right. \nonumber \\
\times&\left. \, 
\begin{pmatrix}
\sqrt{3} \! - \! 1 & (2 \varepsilon b)^{1/2} \\
-(2/\varepsilon b)^{1/2} & \sqrt{3} \! - \! 1
\end{pmatrix} \! 
\begin{pmatrix}
\mathbb{T}_{11,k}(-1;\tau) & \mathbb{T}_{12,k}(-1;\tau) \\
\mathbb{T}_{21,k}(-1;\tau) & \mathbb{T}_{22,k}(-1;\tau)
\end{pmatrix} \right),
\end{align}
with $\ell_{0,k}^{+}$, $\ell_{1,k}^{+}$, and $\ell_{2,k}^{+}$ defined in Lemma \ref{nprcl}, $(\mathbb{T}_{ij,k}(-1;\tau))_{i,j=1,2}$ 
defined in Proposition \ref{prop3.1.6}, and $\hat{\beta}_{k}(\tau)$ defined by equation \eqref{ellohk16}, and
\begin{equation} \label{ellohk11} 
\mathcal{O} \big(\mathbb{E}^{0}_{k}(\tau) \big) \underset{\tau \to +\infty}{:=} 
\begin{pmatrix}
\mathcal{O} \big(\tau^{-\frac{1}{3}+3 \delta_{k}} \big) & \mathcal{O} \big(\tau^{-\frac{1}{3}(\frac{1-k}{2})-\delta_{k}} \big) \\
\mathcal{O} \big(\tau^{-\frac{1}{3}(\frac{1+k}{2})-\delta_{k}} \big) & \mathcal{O} \big(\tau^{-\frac{1}{3}+3\delta_{k}} \big)
\end{pmatrix}.
\end{equation}
\end{cccc}

\emph{Proof}. Denote by $\widetilde{\Psi}_{\scriptscriptstyle \mathrm{WKB},k}(\widetilde{\mu},\tau)$, $k \! = \! \pm 1$, the 
solution of equation \eqref{eq3.3} that has leading-order asymptotics given by equations \eqref{eq3.16}--\eqref{eq3.18} 
in the canonical domain containing the Stokes curve approaching, for $k \! = \! +1$ (resp., $k \! = \! -1)$, the real 
$\widetilde{\mu}$-axis {}from above (resp., below) as $\widetilde{\mu} \! \to \! 0$. Let $\mathfrak{L}_{k}^{0}(\tau)$, 
$k \! = \! \pm 1$, be defined by equation \eqref{ellohk1}; rewrite $\mathfrak{L}_{k}^{0}(\tau)$ in the following form:
\begin{equation} \label{ellohk12} 
\mathfrak{L}_{k}^{0}(\tau) \! = \! \left(\big(\widetilde{\Psi}_{k}(\widetilde{\mu},\tau) \big)^{-1} 
\widetilde{\Psi}_{\scriptscriptstyle \mathrm{WKB},k}(\widetilde{\mu},\tau) \right) \! \left(\big(
\widetilde{\Psi}_{\scriptscriptstyle \mathrm{WKB},k}(\widetilde{\mu},\tau) \big)^{-1} \tau^{-\frac{1}{12} \sigma_{3}} 
\mathbb{X}^{0}_{1}(\tau^{-1/6} \widetilde{\mu},\tau) \right) \! \mathbb{S}_{k}^{\ast},
\end{equation}
where $\mathbb{S}_{k}^{\ast}$ is defined by equation \eqref{ellohk9}. Noting that $\widetilde{\Psi}_{k}(\widetilde{\mu},
\tau)$, $\widetilde{\Psi}_{\scriptscriptstyle \mathrm{WKB},k}(\widetilde{\mu},\tau)$, and $\tau^{-\frac{1}{12} \sigma_{3}} 
\mathbb{X}^{0}_{1}(\tau^{-1/6} \widetilde{\mu},\tau)$ are all solutions of equation \eqref{eq3.3}, it follows that they differ 
on the right by non-degenerate, $\widetilde{\mu}$-independent, $\mathrm{M}_{2}(\mathbb{C})$-valued factors: via this 
observation, one evaluates, asymptotically, each of the factors appearing in equation \eqref{ellohk12} by considering 
separate limits, namely, $\widetilde{\mu} \! \to \! \alpha_{k}$ and $\widetilde{\mu} \! \to \! 0$, respectively; more precisely, 
for $k \! = \! \pm 1$,
\begin{align} \label{ellohk13} 
&\big(\widetilde{\Psi}_{k}(\widetilde{\mu},\tau) \big)^{-1} \widetilde{\Psi}_{\scriptscriptstyle 
\mathrm{WKB},k}(\widetilde{\mu},\tau) \underset{\tau \to +\infty}{:=} \nonumber \\
&\underbrace{\left((b(\tau))^{-\frac{1}{2} \sigma_{3}} \mathcal{G}_{0,k} 
\mathfrak{B}_{k}^{\frac{1}{2} \sigma_{3}} \mathbb{F}_{k}(\tau) \Xi_{k}(\tau;
\widetilde{\Lambda}) \hat{\chi}_{k}(\widetilde{\Lambda}) \Phi_{M,k}(\widetilde{\Lambda}) 
\right)^{-1}T(\widetilde{\mu}) \me^{\scriptscriptstyle \mathrm{W}_{k}(\widetilde{\mu},
\tau)}}_{\widetilde{\mu}=\widetilde{\mu}_{0,k}, \, \, \, \, \widetilde{\Lambda} \, 
\underset{\tau \to +\infty}{\thicksim} \, \mathcal{O}(\tau^{\delta_{k}}), \, \, 0< \delta < \delta_{k}<\frac{1}{24}, 
\, \, \, \, \arg (\widetilde{\Lambda})=\frac{\pi m_{0}}{2}+\frac{\pi}{4}-\frac{1}{2} \arg (\mu_{k}(\tau)), \, \, 
m_{0} \in \lbrace -1,0,1,2 \rbrace},
\end{align}
where (cf. Lemma \ref{linfnewlemm}) $\mathbb{F}_{k}(\tau)$ and $\Xi_{k}(\tau;\widetilde{\Lambda})$ are given 
in equations \eqref{ellinfk13} and \eqref{ellinfk14}, respectively, $\mathrm{W}_{k}(\widetilde{\mu},\tau) \! := \! 
-\sigma_{3} \mi \tau^{2/3} \int_{\widetilde{\mu}_{0,k}}^{\widetilde{\mu}}l_{k}(\xi) \, \md \xi \! - \! 
\int_{\widetilde{\mu}_{0,k}}^{\widetilde{\mu}} \diag ((T(\xi))^{-1} \partial_{\xi}T(\xi)) \, \md \xi$, and $\hat{\chi}_{k}
(\widetilde{\Lambda})$ has the asymptotics \eqref{ellinfk15}, and
\begin{equation} \label{ellohk14} 
\big(\widetilde{\Psi}_{\scriptscriptstyle \mathrm{WKB},k}(\widetilde{\mu},\tau) \big)^{-1} \tau^{-\frac{1}{12} \sigma_{3}} 
\mathbb{X}^{0}_{1}(\tau^{-1/6} \widetilde{\mu},\tau) \underset{\tau \to +\infty}{:=} \lim_{\underset{\arg (\widetilde{\mu})
=\pi}{\Omega_{1}^{0} \ni \widetilde{\mu} \to 0}} \left(\big(T(\widetilde{\mu}) \me^{\scriptscriptstyle \mathrm{W}_{k}
(\widetilde{\mu},\tau)} \big)^{-1} \tau^{-\frac{1}{12} \sigma_{3}} \mathbb{X}^{0}_{1}(\tau^{-1/6} \widetilde{\mu},\tau) \right).
\end{equation}

One commences by considering the asymptotics subsumed in the definition \eqref{ellohk14}. {}From the asymptotics 
for $\mathbb{X}^{0}_{1}(\tau^{-1/6} \widetilde{\mu},\tau)$ stated in Proposition \ref{prop1.4}, equations \eqref{iden3}, 
\eqref{iden4}, \eqref{expforeych}, \eqref{expforkapp}, \eqref{eq3.39}, \eqref{eq3.40}, \eqref{eqpeetee}, \eqref{eq3.43}, 
\eqref{asympforf8}, \eqref{eq3.52}, and \eqref{prcy22}, one arrives at, via the conditions \eqref{iden5} and the asymptotics 
\eqref{terrbos6}, \eqref{asympforf3}, and \eqref{prcyomg1}, 
\begin{equation} \label{ellohk15} 
\lim_{\underset{\arg (\widetilde{\mu})=\pi}{\Omega_{1}^{0} \ni \widetilde{\mu} \to 0}} \left(\big(T(\widetilde{\mu}) 
\me^{\scriptscriptstyle \mathrm{W}_{k}(\widetilde{\mu},\tau)} \big)^{-1} \tau^{-\frac{1}{12} \sigma_{3}} \mathbb{X}^{0}_{1}
(\tau^{-1/6} \widetilde{\mu},\tau) \right) \underset{\tau \to +\infty}{=} \left(\dfrac{\mi (\varepsilon b)^{1/4}}{\sqrt{\smash[b]{b(\tau)}}} 
\right)^{\sigma_{3}} \exp \big(\hat{\beta}_{k}(\tau) \sigma_{3} \big), \quad k \! = \! \pm 1,
\end{equation}
where
\begin{align} \label{ellohk16} 
\hat{\beta}_{k}(\tau) :=& \, \mi \tau^{2/3}3 \sqrt{3} \alpha_{k}^{2} \! + \! \mi 2 \sqrt{3} \, \widetilde{\Lambda}^{2} 
\! + \! \mi (a \! - \! \mi/2) \ln ((\sqrt{3} \! + \! 1)/\sqrt{2}) \! - \! \dfrac{(5 \! + \! 9 \sqrt{3}) \mathfrak{p}_{k}(\tau)}{6 
\sqrt{3} \alpha_{k}^{2}} \nonumber \\
+& \, \left(\dfrac{\mi}{2 \sqrt{3}} \! \left((a \! - \! \mi/2) \! + \! \alpha_{k}^{-2} \tau^{2/3} \hat{h}_{0}(\tau) \right) \! + \! 
\dfrac{2 \mathfrak{p}_{k}(\tau)}{3 \sqrt{3} \alpha_{k}^{2}} \right) \! \left(-\dfrac{1}{3} \ln \tau \! + \! \ln \widetilde{\Lambda} 
\! + \! \ln (\me^{\mi k \pi}/3 \alpha_{k}) \right) \nonumber \\
-& \, \dfrac{(\sqrt{3} \! + \! 1) \mathfrak{p}_{k}(\tau)}{\sqrt{3} \alpha_{k} \tau^{-1/3} \widetilde{\Lambda}} \! + \! 
\mathcal{O} \! \left(\left(\dfrac{\tilde{\mathfrak{c}}_{1,k} \tau^{-1/3} \! + \! \tilde{\mathfrak{c}}_{2,k} 
\tilde{r}_{0}(\tau)}{\widetilde{\Lambda}^{2}} \right) \! \left(\tilde{\mathfrak{c}}_{3,k} \tau^{-1/3} \! + \! 
\tilde{\mathfrak{c}}_{4,k}(\tilde{r}_{0}(\tau) \! + \! 4v_{0}(\tau)) \right) \right) \nonumber \\
+& \, \mathcal{O} \big(\tau^{-1/3} \widetilde{\Lambda}^{3} \big) \! + \! \mathcal{O} \big(\tau^{-1/3} \widetilde{\Lambda} \big) 
\! + \! \mathcal{O} \! \left(\dfrac{\tau^{-1/3}}{\widetilde{\Lambda}} \left(\tilde{\mathfrak{c}}_{5,k} \! + \! \tilde{\mathfrak{c}}_{6,k} 
\tau^{2/3} \hat{h}_{0}(\tau) \! + \! \tilde{\mathfrak{c}}_{7,k}(\tau^{2/3} \hat{h}_{0}(\tau))^{2} \right) \right) \nonumber \\
+& \, \mathcal{O} \! \left(\tau^{-2/3} \hat{d}_{0,k}(\tau) \! \left(-\dfrac{1}{3} \ln \tau \! + \! \ln \widetilde{\Lambda} \right) \! \right),
\end{align}
$\tilde{\mathfrak{c}}_{m,k}$, $m \! = \! 1,2,\dotsc,7$, are $\mathcal{O}(1)$, and $\hat{d}_{0,k}(\tau)$ is defined in the proof of 
Proposition \ref{prop3.1.3}.

One now derives the asymptotics defined by equation \eqref{ellohk13}. {}From the asymptotics \eqref{iden55} for $\varpi \! 
= \! -1$, equation \eqref{prcy2} for $\Phi_{M,k}(\widetilde{\Lambda})$ (in conjunction with its large-$\widetilde{\Lambda}$ 
symptotics stated in Remark \ref{leedasyue}), the definitions \eqref{ellinfk13} and \eqref{ellinfk14} (concomitant with the fact 
that $\det (\Xi_{k}(\tau;\widetilde{\Lambda})) \! = \! 1)$, and the asymptotics \eqref{ellinfk15}, one shows, via the relation 
$(\mathrm{W}_{k}(\widetilde{\mu}_{0,k},\tau))_{i,j=1,2} \! = \! 0$ and the definition \eqref{ellohk13}, that, for $k \! = \! \pm 1$,
{\fontsize{10pt}{11pt}\selectfont
\begin{align} \label{ellohk17} 
\big(\widetilde{\Psi}_{k}(\widetilde{\mu},\tau) \big)^{-1} \widetilde{\Psi}_{\scriptscriptstyle \mathrm{WKB},k}
(\widetilde{\mu},\tau) \underset{\tau \to +\infty}{:=}& \, \Phi_{M,k}^{-1}(\widetilde{\Lambda}) \hat{\chi}_{k}^{-1}
(\widetilde{\Lambda}) \Xi_{k}^{-1}(\tau;\widetilde{\Lambda}) \mathbb{F}_{k}^{-1}(\tau) \mathfrak{B}_{k}^{
-\frac{1}{2} \sigma_{3}} \mathcal{G}_{0,k}^{-1}(b(\tau))^{\frac{1}{2} \sigma_{3}}
T(\widetilde{\mu}_{0,k}) \nonumber \\
\underset{\tau \to +\infty}{=}& \, (\mathcal{R}_{m_{0}}(k))^{-1} 
\me^{-\mathcal{P}_{0}^{\ast} \sigma_{3}} \mathfrak{Q}_{0,k}(\tau) \! \left(
\mathrm{I} \! + \! \dfrac{1}{\widetilde{\Lambda}} \mathfrak{Q}_{0,k}^{-1}
(\tau) \hat{\psi}_{1,k}^{-1}(\tau) \mathfrak{Q}_{0,k}(\tau) \right. \nonumber \\
+&\left. \, \dfrac{1}{\widetilde{\Lambda}^{2}} \mathfrak{Q}_{0,k}^{-1}(\tau) 
\hat{\psi}_{2,k}^{-1}(\tau) \mathfrak{Q}_{0,k}(\tau) \! + \! \mathcal{O} \! 
\left(\dfrac{1}{\widetilde{\Lambda}^{3}} \mathfrak{Q}_{0,k}^{-1}(\tau) 
\hat{\psi}_{3,k}^{-1}(\tau) \mathfrak{Q}_{0,k}(\tau) \right) \right) \nonumber \\
\times& \, \left(\mathrm{I} \! + \! \mathcal{O} \! \left(\lvert \nu (k) \! + \! 1 \rvert^{2} 
\lvert p_{k}(\tau) \rvert^{-2} \tau^{-\epsilon_{\mathrm{\scriptscriptstyle TP}}(k)} 
\mathfrak{Q}_{0,k}^{-1}(\tau) \tilde{\mathfrak{C}}_{k}(\tau) 
\mathfrak{Q}_{0,k}(\tau) \right) \right) \nonumber \\
\times& \, \left(\mathrm{I} \! + \! \widetilde{\Lambda} \mathfrak{Q}_{0,k}^{-1}
(\tau) \gimel_{\scriptscriptstyle A,k}^{-1}(\tau) \mathfrak{Q}_{0,k}
(\tau) \! + \! \widetilde{\Lambda}^{2} \mathfrak{Q}_{0,k}^{-1}(\tau) 
\gimel_{\scriptscriptstyle B,k}^{-1}(\tau) \mathfrak{Q}_{0,k}(\tau) \right) 
\nonumber \\
\times& \, \left(\mathrm{I} \! + \! \widetilde{\Lambda} \tau^{-1/3} 
\mathbb{P}_{0,k}(\tau) \! + \! \dfrac{1}{\widetilde{\Lambda}} 
\widehat{\mathbb{E}}_{0,k}(\tau) \! + \! \mathcal{O} \! \left((\tau^{-1/3} 
\widetilde{\Lambda})^{2} \widetilde{\mathbb{E}}_{0,k}(\tau) \right) \right),
\end{align}}
where $\mathrm{M}_{2}(\mathbb{C}) \! \ni \! \mathcal{R}_{m_{0}}(k)$, $m_{0} \! \in \! \lbrace -1,0,1,2 \rbrace$, are 
defined in Remark \ref{leedasyue}, $\mathcal{P}_{0}^{\ast}$, $\hat{\psi}_{1,k}^{-1}(\tau)$, $\hat{\psi}_{2,k}^{-1}
(\tau)$, and $\hat{\psi}_{3,k}^{-1}(\tau)$ are defined by equations \eqref{ellinfk20}, \eqref{ellinfk22}, \eqref{ellinfk23}, 
and \eqref{ellinfk24}, respectively,
\begin{equation} \label{ellohk18} 
\mathfrak{Q}_{0,k}(\tau) \! := \! \mathbb{F}_{k}^{-1}(\tau) \! \left(\left(\dfrac{2^{1/4} 
\sqrt{\smash[b]{b(\tau)}}}{(\varepsilon b)^{1/4}(\sqrt{3} \! - \! 1)^{1/2} \sqrt{\smash[b]{\mathfrak{B}_{k}}}} \right)^{
\sigma_{3}} \! + \! \mathfrak{B}_{k}^{-\frac{1}{2} \sigma_{3}} \Delta G_{k}^{0}(\tau)(b(\tau))^{\frac{1}{2} \sigma_{3}} \right),
\end{equation}
with $\Delta G_{k}^{0}(\tau)$ defined by equation \eqref{ellohk8},
\begin{gather}
\mathbb{P}_{0,k}(\tau) \! := \! (b(\tau))^{-\frac{1}{2} \ad (\sigma_{3})} \! 
\begin{pmatrix}
0 & -\frac{(\varepsilon b)^{1/2}}{3 \sqrt{2} \alpha_{k}} \\
\frac{(\varepsilon b)^{-1/2}}{3 \sqrt{2} \alpha_{k}} & 0
\end{pmatrix}, \label{ellohk19} \\
\widehat{\mathbb{E}}_{0,k}(\tau) \! := \! \dfrac{1}{2 \sqrt{3}(\sqrt{3} \! - \! 1)}(b(\tau))^{-\frac{1}{2} \ad (\sigma_{3})} \! 
\begin{pmatrix}
\sqrt{3} \! - \! 1 & (2 \varepsilon b)^{1/2} \\
-(2/\varepsilon b)^{1/2} & \sqrt{3} \! - \! 1
\end{pmatrix} \! 
\begin{pmatrix}
\mathbb{T}_{11,k}(-1;\tau) & \mathbb{T}_{12,k}(-1;\tau) \\
\mathbb{T}_{21,k}(-1;\tau) & \mathbb{T}_{22,k}(-1;\tau)
\end{pmatrix}, \label{ellohk20} \\
\widetilde{\mathbb{E}}_{0,k}(\tau) \! := \! \dfrac{1}{2 \sqrt{3}(\sqrt{3} \! - \! 1)}(b(\tau))^{-\frac{1}{2} \ad (\sigma_{3})} \! 
\begin{pmatrix}
\sqrt{3} \! - \! 1 & (2 \varepsilon b)^{1/2} \\
-(2/\varepsilon b)^{1/2} & \sqrt{3} \! - \! 1
\end{pmatrix} \tilde{\mathfrak{C}}_{k}^{\lozenge}, \label{ellohk21}
\end{gather}
$\mathrm{M}_{2}(\mathbb{C}) \! \ni \! \tilde{\mathfrak{C}}_{k}(\tau) \genfrac{}{}{0pt}{3}{=}{\tau \to +\infty} \mathcal{O}(1)$, 
$(\mathbb{T}_{ij,k}(-1;\tau))_{i,j=1,2}$ defined in Proposition \ref{prop3.1.6}, and $\mathrm{M}_{2}(\mathbb{C}) \! \ni \! 
\tilde{\mathfrak{C}}_{k}^{\lozenge}$ is $\mathcal{O}(1)$.

Recalling the definitions \eqref{ellohk13} and \eqref{ellohk14}, and substituting the expansions \eqref{ellohk15}, \eqref{ellohk16}, 
and \eqref{ellohk17} into equation \eqref{ellohk12}, one shows, via the conditions \eqref{iden5}, the definition \eqref{prpr1}, the 
restrictions \eqref{pc4}, the asymptotics \eqref{prcyzeek1}, \eqref{prcychik1}, and \eqref{prcymuk1}, and (cf. step \pmb{(xi)} in the 
proof of Lemma \ref{nprcl}) $\arg (\mu_{k}(\tau)) \genfrac{}{}{0pt}{3}{=}{\tau \to +\infty} \tfrac{\pi}{2} \big(1 \! + \! \mathcal{O}
(\tau^{-2/3}) \big)$, and the restriction \eqref{restr1}, that
\begin{equation} \label{ellohk22} 
\mathfrak{L}^{0}_{k}(\tau) \underset{\tau \to +\infty}{=} (\mathcal{R}_{m_{0}}(k))^{-1} 
\me^{\hat{\mathfrak{z}}_{k}^{0}(\tau) \sigma_{3}} \! \left(\dfrac{\mi 2^{1/4}}{(\sqrt{3} \! - \! 1)^{1/2} 
\sqrt{\smash[b]{\mathfrak{B}_{k}}}} \right)^{\sigma_{3}} \me^{\Delta \hat{\mathfrak{z}}_{k}(\tau) \sigma_{3}} 
\operatorname{diag} \! \left(\hat{\mathbb{A}}_{0}^{0}(\tau),\hat{\mathbb{B}}_{0}^{0}(\tau) \right) \! 
\overset{\Ydown}{\mathbb{E}}_{\scriptscriptstyle \mathfrak{L}^{0}_{k}}^{\raise-6.75pt\hbox{$\scriptstyle \leftslice$}}(\tau) 
\mathbb{S}_{k}^{\ast}, \quad k \! = \! \pm 1,
\end{equation}
where $\hat{\mathfrak{z}}_{k}^{0}(\tau)$, $\Delta \hat{\mathfrak{z}}_{k}(\tau)$, $\hat{\mathbb{A}}_{0}^{0}(\tau)$, 
and $\hat{\mathbb{B}}_{0}^{0}(\tau)$ are defined by equations \eqref{ellohk4}--\eqref{ellohk7}, respectively, and
{\fontsize{10pt}{11pt}\selectfont
\begin{align} \label{ellohk23} 
\overset{\Ydown}{\mathbb{E}}_{\scriptscriptstyle \mathfrak{L}^{0}_{k}}^{\raise-6.75pt\hbox{$\scriptstyle \leftslice$}}
(\tau) \underset{\tau \to +\infty}{:=}& \, \left(\mathrm{I} \! + \! \mathcal{O}(\tau^{-1/3} \widetilde{\Lambda}^{3} \sigma_{3}) 
\right) \! \left(\mathrm{I} \! + \! \mathcal{O} \! \left(\frac{\hat{\mathbb{C}}_{0}^{0}(\tau) \sqrtsign{b(\tau)}}{\hat{\mathbb{A}}_{0}^{0}
(\tau) \me^{\widehat{\beta}_{k}^{\ast}(\tau)}} \sigma_{+} \right) \! + \! \mathcal{O} \! \left(\frac{\hat{\mathbb{D}}_{0}^{0}(\tau) 
\me^{\widehat{\beta}_{k}^{\ast}(\tau)}}{\hat{\mathbb{B}}_{0}^{0}(\tau) \sqrtsign{b(\tau)}} \sigma_{-} \right) \right) \nonumber \\
\times& \, \left(\mathrm{I} \! + \! \dfrac{1}{\widetilde{\Lambda}} \hat{\psi}_{1,k}^{-1,\natural}(\tau) \! + \! 
\dfrac{1}{\widetilde{\Lambda}^{2}} \hat{\psi}_{2,k}^{-1,\natural}(\tau) \! + \! \mathcal{O} \! \left(
\dfrac{1}{\widetilde{\Lambda}^{3}} \hat{\psi}_{3,k}^{-1,\natural}(\tau) \right) \right) \! \left(\mathrm{I} \! + \! \mathcal{O} 
\! \left(\dfrac{\lvert \nu (k) \! + \! 1 \rvert^{2} \tau^{-\epsilon_{\mathrm{\scriptscriptstyle TP}}(k)}}{\lvert p_{k}(\tau) \rvert^{2}} 
\right. \right. \nonumber \\
\times&\left. \left. \, \mathfrak{Q}_{\ast,k}^{-1}(\tau) \tilde{\mathfrak{C}}_{k}(\tau) \mathfrak{Q}_{\ast,k}(\tau) \right) 
\right) \! \left(\mathrm{I} \! + \! \widetilde{\Lambda} \gimel_{\scriptscriptstyle A,k}^{\natural}(\tau) \! + \! 
\widetilde{\Lambda}^{2} \gimel_{\scriptscriptstyle B,k}^{\natural}(\tau) \right) \! \left(\mathrm{I} \! + \! \widetilde{\Lambda} 
\tau^{-1/3} \mathbb{P}_{0,k}^{\natural}(\tau) 
\right. \nonumber \\
+&\left. \, \dfrac{1}{\widetilde{\Lambda}} \widehat{\mathbb{E}}_{0,k}^{\natural}(\tau) \! + \! \mathcal{O} \! \left((\tau^{-1/3} 
\widetilde{\Lambda})^{2} \widetilde{\mathbb{E}}_{0,k}^{\natural}(\tau) \right) \right),
\end{align}}
where $-\widehat{\beta}_{k}^{\ast}(\tau) \! := \! \tfrac{\mi a}{6} \ln \tau \! + \! \mi 3 \alpha_{k}^{2} \tau^{2/3}$,
\begin{gather}
\hat{\mathbb{C}}_{0}^{0}(\tau) \! := \! -\mi (\varepsilon b)^{-1/4}(\Delta G_{k}^{0}(\tau))_{12}, \label{ellohk24} \\
\hat{\mathbb{D}}_{0}^{0}(\tau) \! := \! \mi (\varepsilon b)^{1/4}(\Delta G_{k}^{0}(\tau))_{21} \! - \! 
\dfrac{\mathfrak{A}_{k}}{\mathfrak{B}_{k}} \! \left(\dfrac{\mi 4 \sqrt{3} \mathcal{Z}_{k}}{\chi_{k}(\tau)} \! - \! 1 \right) 
\! \left(\dfrac{\mi 2^{1/4}}{(\sqrt{3} \! - \! 1)^{1/2}} \! + \! \mi (\varepsilon b)^{1/4}(\Delta G_{k}^{0}(\tau))_{11} \right), 
\label{ellohk25} \\
\hat{\psi}_{m,k}^{-1,\natural}(\tau) \! := \! \mathfrak{Q}_{\ast,k}^{-1}(\tau) \hat{\psi}_{m,k}^{-1}(\tau) \mathfrak{Q}_{\ast,k}
(\tau), \quad m \! = \! 1,2,3, \label{ellohk26} \\
\mathfrak{Q}_{\ast,k}(\tau) \! := \! \mathfrak{Q}_{0,k}(\tau)(\mi (\varepsilon b)^{1/4})^{\sigma_{3}}(b(\tau))^{-\frac{1}{2} 
\sigma_{3}} \me^{\hat{\beta}_{k}(\tau) \sigma_{3}}, \label{ellohk27} \\
\gimel_{\scriptscriptstyle A,k}^{\natural}(\tau) \! := \! \mathfrak{Q}_{\ast,k}^{-1}(\tau) \gimel_{\scriptscriptstyle A,k}^{-1}
(\tau) \mathfrak{Q}_{\ast,k}(\tau), \label{ellohk28} \\
\gimel_{\scriptscriptstyle B,k}^{\natural}(\tau) \! := \! \mathfrak{Q}_{\ast,k}^{-1}(\tau) \gimel_{\scriptscriptstyle B,k}^{-1}
(\tau) \mathfrak{Q}_{\ast,k}(\tau), \label{ellohk29} \\
\mathbb{P}_{0,k}^{\natural}(\tau) \! := \! (\mi (\varepsilon b)^{1/4})^{-\ad (\sigma_{3})}(b(\tau))^{\frac{1}{2} \ad (\sigma_{3})} 
\me^{-\hat{\beta}_{k}(\tau) \ad (\sigma_{3})} \mathbb{P}_{0,k}(\tau), \label{ellohk30} \\
\widehat{\mathbb{E}}_{0,k}^{\natural}(\tau) \! := \! (\mi (\varepsilon b)^{1/4})^{-\ad (\sigma_{3})}(b(\tau))^{\frac{1}{2} 
\ad (\sigma_{3})} \me^{-\hat{\beta}_{k}(\tau) \ad (\sigma_{3})} \widehat{\mathbb{E}}_{0,k}(\tau), \label{ellohk31} \\
\widetilde{\mathbb{E}}_{0,k}^{\natural}(\tau) \! := \! (\mi (\varepsilon b)^{1/4})^{-\ad (\sigma_{3})}(b(\tau))^{\frac{1}{2} 
\ad (\sigma_{3})} \me^{-\hat{\beta}_{k}(\tau) \ad (\sigma_{3})} \widetilde{\mathbb{E}}_{0,k}(\tau). \label{ellohk32} 
\end{gather}
The calculations for the asymptotics (as $\tau \! \to \! +\infty$ with $\varepsilon b \! > \! 0$) of the error function 
$\overset{\Ydown}{\mathbb{E}}_{\scriptscriptstyle \mathfrak{L}^{0}_{k}}^{\raise-6.75pt\hbox{$\scriptstyle \leftslice$}}
(\tau)$ (cf. definition \eqref{ellohk23}) are similar to those for the error function 
$\overset{\Yup}{\mathbb{E}}_{\scriptscriptstyle \mathfrak{L}^{\infty}_{k}}^{\raise-6.75pt\hbox{$\scriptstyle \leftslice$}}(\tau)$ 
presented in the proof of Lemma \ref{linfnewlemm}; therefore, via the conditions \eqref{iden5}, the restrictions \eqref{pc4} 
and \eqref{restr1}, the definitions \eqref{eqpeetee}, \eqref{peekayity}, \eqref{prpr1}, \eqref{prpr3}, \eqref{prpr4}, \eqref{prcy57}, 
\eqref{prcy58}, \eqref{ellinfk13}, \eqref{ellinfk22}--\eqref{ellinfk24}, \eqref{ellohk6}--\eqref{ellohk8}, \eqref{ellohk18}--\eqref{ellohk21}, 
and \eqref{ellohk24}--\eqref{ellohk32}, and the asymptotics \eqref{tr1}, \eqref{tr3}, \eqref{asympforf3}, \eqref{prcyzeek1}, 
\eqref{prcyg4}--\eqref{prcybk1}, \eqref{prcyomg1}--\eqref{prcyell2k1}, and \eqref{ellohk16}, upon imposing the conditions 
\eqref{ellinfk2a} and \eqref{ellinfk2b} and proceeding as in the proof of Lemma \ref{linfnewlemm}, one shows that, for 
$k \! = \! \pm 1$,
\begin{equation} \label{ellohk33} 
\overset{\Ydown}{\mathbb{E}}_{\scriptscriptstyle \mathfrak{L}^{0}_{k}}^{\raise-6.75pt\hbox{$\scriptstyle \leftslice$}}
(\tau) \underset{\tau \to +\infty}{=} \big(\mathrm{I} \! + \! \mathbb{E}_{{\scriptscriptstyle \mathcal{N}},k}^{0}(\tau) \big) 
\big(\mathrm{I} \! + \! \mathcal{O}(\mathbb{E}^{0}_{k}(\tau)) \big),
\end{equation}
where $\mathbb{E}_{{\scriptscriptstyle \mathcal{N}},k}^{0}(\tau)$ and $\mathcal{O}(\mathbb{E}^{0}_{k}(\tau))$ are 
defined by equations \eqref{ellohk10} and \eqref{ellohk11}, respectively.\footnote{Note that $\mathcal{O} \! \left(
\frac{\hat{\mathbb{C}}_{0}^{0}(\tau) \sqrtsign{b(\tau)}}{\hat{\mathbb{A}}_{0}^{0}(\tau) \me^{\widehat{\beta}_{k}^{\ast}(\tau)}} 
\right) \underset{\tau \to +\infty}{=} \mathcal{O}(\tau^{-2/3})$ and $\mathcal{O} \! \left(\frac{\hat{\mathbb{D}}_{0}^{0}(\tau) 
\me^{\widehat{\beta}_{k}^{\ast}(\tau)}}{\hat{\mathbb{B}}_{0}^{0}(\tau) \sqrtsign{b(\tau)}} \right) 
\underset{\tau \to +\infty}{=} \mathcal{O}(\tau^{-2/3})$. The asymptotics for the function 
$\mathbb{E}_{{\scriptscriptstyle \mathcal{N}},k}^{0}(\tau)$ is presented in the proof of Lemma \ref{ginversion} (see 
Section \ref{finalsec}).} Thus, via the asymptotics \eqref{ellohk22} and \eqref{ellohk33}, one arrives at the results stated 
in the lemma. \hfill $\qed$
\begin{dddd} \label{theor3.3.1} 
Assume that the conditions \eqref{iden5}, \eqref{pc4}, \eqref{restr1}, \eqref{ellinfk2a}, and \eqref{ellinfk2b} are valid; 
then, the connection matrix has the following asymptotics:
\begin{equation} \label{jeek1} 
G_{k} \underset{\tau \to +\infty}{=} \widetilde{G}(k) \widehat{\mathscr{G}}(k) \big(\mathrm{I} \! + \! \mathcal{O}
(\mathbb{E}_{k}^{\scriptscriptstyle G_{k}}(\tau)) \big), \quad k \! = \! \pm 1,
\end{equation}
where
\begin{gather}
\widetilde{G}(k) \! := \! (\mathbb{S}_{k}^{\ast})^{-1} \mathrm{G}^{\ast}(k), \label{jeek2} \\
\widehat{\mathscr{G}}(k) \! := \! (\mathrm{G}^{\ast}(k))^{-1} \big(\mathrm{I} \! + \! \mathbb{E}_{{\scriptscriptstyle \mathcal{N}},k}^{0}
(\tau) \big)^{-1} \mathrm{G}^{\ast}(k) \big(\mathrm{I} \! + \! \mathbb{E}_{{\scriptscriptstyle \mathcal{N}},k}^{\infty}(\tau) \big), 
\label{jeek3}
\end{gather}
with $\mathbb{E}_{{\scriptscriptstyle \mathcal{N}},k}^{\infty}(\tau)$, $\mathbb{S}_{k}^{\ast}$, and 
$\mathbb{E}_{{\scriptscriptstyle \mathcal{N}},k}^{0}(\tau)$ defined by equations \eqref{ellinfk9}, \eqref{ellohk9}, 
and \eqref{ellohk10}, respectively, and
\begin{equation} \label{jeek4} 
\mathrm{G}^{\ast}(k) \! = \! 
\begin{pmatrix}
\frac{\hat{\mathbb{G}}_{11}(k) \hat{\mathbb{B}}_{0}^{\infty}(\tau)}{\hat{\mathbb{A}}_{0}^{0}(\tau)} 
\me^{-\Delta \tilde{\mathfrak{z}}_{k}(\tau)-\Delta \hat{\mathfrak{z}}_{k}(\tau)} & \frac{\hat{\mathbb{G}}_{12}(k) 
\hat{\mathbb{A}}_{0}^{\infty}(\tau)}{\hat{\mathbb{A}}_{0}^{0}(\tau)} \me^{\Delta \tilde{\mathfrak{z}}_{k}(\tau)
-\Delta \hat{\mathfrak{z}}_{k}(\tau)} \\
\frac{\hat{\mathbb{G}}_{21}(k) \hat{\mathbb{B}}_{0}^{\infty}(\tau)}{\hat{\mathbb{B}}_{0}^{0}(\tau)} 
\me^{-\Delta \tilde{\mathfrak{z}}_{k}(\tau)+\Delta \hat{\mathfrak{z}}_{k}(\tau)} & \frac{\hat{\mathbb{G}}_{22}(k) 
\hat{\mathbb{A}}_{0}^{\infty}(\tau)}{\hat{\mathbb{B}}_{0}^{0}(\tau)} \me^{\Delta \tilde{\mathfrak{z}}_{k}(\tau)
+\Delta \hat{\mathfrak{z}}_{k}(\tau)}
\end{pmatrix},
\end{equation}
where
\begin{gather} 
\hat{\mathbb{G}}_{11}(k) \! := \! -\dfrac{\mi \sqrt{2 \pi} \, p_{k}(\tau) \mathfrak{B}_{k} 
\sqrt{\smash[b]{b(\tau)}} \, \me^{\mi \pi (\nu (k)+1)}}{(\varepsilon b)^{1/4}
(2 \! + \! \sqrt{3})^{1/2}(2 \mu_{k}(\tau))^{1/2} \Gamma (-\nu (k))} 
\exp \! \left(-\tilde{\mathfrak{z}}_{k}^{0}(\tau) \! - \! \hat{\mathfrak{z}}_{k}^{0}
(\tau) \right), \label{jeek5} \\
\hat{\mathbb{G}}_{12}(k) \! := \! -\dfrac{\mi (\varepsilon b)^{1/4}}{\sqrt{
\smash[b]{b(\tau)}}} \exp \! \left(\tilde{\mathfrak{z}}_{k}^{0}(\tau) \! - \! 
\hat{\mathfrak{z}}_{k}^{0}(\tau) \right), \label{jeek6} \\
\hat{\mathbb{G}}_{21}(k) \! := \! -\dfrac{\mi \sqrt{\smash[b]{b(\tau)}} \, \me^{-2 
\pi \mi (\nu (k)+1)}}{(\varepsilon b)^{1/4}} \exp \! \left(-\tilde{\mathfrak{z}}_{k}^{0}
(\tau) \! + \! \hat{\mathfrak{z}}_{k}^{0}(\tau) \right), \label{jeek7} \\
\hat{\mathbb{G}}_{22}(k) \! := \! -\dfrac{\sqrt{2 \pi} \, (\varepsilon b)^{1/4}(2 \! 
+ \! \sqrt{3})^{1/2}(2 \mu_{k}(\tau))^{1/2} \me^{-2 \pi \mi (\nu (k)+1)}}{p_{k}
(\tau) \mathfrak{B}_{k} \sqrt{\smash[b]{b(\tau)}} \, \Gamma (\nu (k) \! + \! 1)} 
\exp \! \left(\tilde{\mathfrak{z}}_{k}^{0}(\tau) \! + \! \hat{\mathfrak{z}}_{k}^{0}
(\tau) \right), \label{jeek8}
\end{gather}
with $\tilde{\mathfrak{z}}_{k}^{0}(\tau)$, $\Delta \tilde{\mathfrak{z}}_{k}(\tau)$, $\hat{\mathbb{A}}_{0}^{\infty}(\tau)$, 
$\hat{\mathbb{B}}_{0}^{\infty}(\tau)$, $\hat{\mathfrak{z}}_{k}^{0}(\tau)$, $\Delta \hat{\mathfrak{z}}_{k}(\tau)$, 
$\hat{\mathbb{A}}_{0}^{0}(\tau)$, and $\hat{\mathbb{B}}_{0}^{0}(\tau)$ defined by equations \eqref{ellinfk4}, 
\eqref{ellinfk5}, \eqref{ellinfk6}, \eqref{ellinfk7}, \eqref{ellohk4}, \eqref{ellohk5}, \eqref{ellohk6}, and \eqref{ellohk7}, 
respectively, and
\begin{equation} \label{jeek9} 
\mathcal{O} \big(\mathbb{E}_{k}^{\scriptscriptstyle G_{k}}(\tau) \big) \underset{\tau \to +\infty}{:=} \mathcal{O} \big(
\mathbb{E}_{k}^{\infty}(\tau) \big) \! + \! \mathcal{O} \! \left(\big(\widetilde{G}(k) \widehat{\mathscr{G}}(k) \big)^{-1} 
\mathbb{E}_{k}^{0}(\tau) \widetilde{G}(k) \widehat{\mathscr{G}}(k) \right),
\end{equation}
with the asymptotics $\mathcal{O}(\mathbb{E}_{k}^{\infty}(\tau))$ and $\mathcal{O}(\mathbb{E}_{k}^{0}(\tau))$ defined 
by equations \eqref{ellinfk10} and \eqref{ellohk11}, respectively.
\end{dddd}

\emph{Proof}. Mimicking the calculations subsumed in the proof of Theorem 3.4.1 of \cite{av2}, one shows that
\begin{equation} \label{jeek10} 
G_{k} \! = \! \big(\mathfrak{L}_{k}^{0}(\tau) \big)^{-1} \mathfrak{L}_{k}^{\infty}(\tau), \quad k \! = \! \pm 1.
\end{equation}
{}From equations \eqref{ellinfk3}--\eqref{ellinfk10}, \eqref{ellohk3}--\eqref{ellohk11}, and \eqref{jeek10}, one arrives at
\begin{align} \label{jeek11} 
G_{k} \underset{\tau \to +\infty}{=}& \, \big(\mathrm{I} \! + \! \mathcal{O}(\mathbb{E}_{k}^{0}(\tau)) \big)
(\mathbb{S}_{k}^{\ast})^{-1} \big(\mathrm{I} \! + \! \mathbb{E}_{{\scriptscriptstyle \mathcal{N}},k}^{0}(\tau) \big)^{-1} 
\me^{-\Delta \hat{\mathfrak{z}}_{k}(\tau) \sigma_{3}} \operatorname{diag} \! \left((\hat{\mathbb{A}}_{0}^{0}(\tau))^{-1},
(\hat{\mathbb{B}}_{0}^{0}(\tau))^{-1} \right) \nonumber \\
\times& \left(\dfrac{\mi 2^{1/4}}{(\sqrt{3} \! - \! 1)^{1/2} \sqrt{\smash[b]{\mathfrak{B}_{k}}}} \right)^{-\sigma_{3}} 
\me^{-\hat{\mathfrak{z}}_{k}^{0}(\tau) \sigma_{3}} \mathcal{R}_{m_{0}}(k)(\mathcal{R}_{m_{\infty}}(k))^{-1} 
\me^{\tilde{\mathfrak{z}}_{k}^{0}(\tau) \sigma_{3}} \left(\dfrac{(\varepsilon b)^{1/4}(\sqrt{3} \! + \! 1)^{1/2}}{2^{1/4} 
\sqrt{\smash[b]{\mathfrak{B}_{k}}} \sqrt{\smash[b]{b(\tau)}}} \right)^{\sigma_{3}} \nonumber \\
\times& \, \mi \sigma_{2} \me^{-\Delta \tilde{\mathfrak{z}}_{k}(\tau) \sigma_{3}} \operatorname{diag} \! 
\left(\hat{\mathbb{B}}_{0}^{\infty}(\tau),\hat{\mathbb{A}}_{0}^{\infty}(\tau) \right) \! \big(\mathrm{I} \! + \! 
\mathbb{E}_{{\scriptscriptstyle \mathcal{N}},k}^{\infty}(\tau) \big) \big(\mathrm{I} \! + \! \mathcal{O}
(\mathbb{E}_{k}^{\infty}(\tau)) \big):
\end{align}
taking $(m_{\infty},m_{0}) \! = \! (0,2)$, that is, $\Delta \arg (\widetilde{\Lambda}) \! := \! \pi (m_{0} \! - \! m_{\infty})/2 
\! = \! \pi$, and using the definitions of $\mathcal{R}_{0}(k)$ and $\mathcal{R}_{2}(k)$ given in Remark \ref{leedasyue}, 
one arrives at, via equation \eqref{jeek11} and the reflection formula $\Gamma (z) \Gamma (1 \! - \! z) \! = \! \pi/\sin 
(\pi z)$, the result stated in the theorem. \hfill $\qed$
\section{The Inverse Monodromy Problem: Asymptotic Solution} \label{finalsec} 
In Subsection \ref{sec3.3}, the corresponding connection matrices, $G_{k}$, $k \! \in \! \lbrace \pm 1 \rbrace$, were 
calculated asymptotically (as $\tau \! \to \! +\infty$ with $\varepsilon b \! > \! 0)$ under the assumption of the validity 
of the conditions \eqref{iden5}, \eqref{pc4}, \eqref{restr1}, \eqref{ellinfk2a}, and \eqref{ellinfk2b}. Using these conditions, 
one can derive the $\tau$-dependent class(es) of functions $G_{k}$ belongs to: this, most general, approach will not 
be adopted here; rather, the isomonodromy condition will be evoked on $G_{k}$, that is, $g_{ij} \! := \! (G_{k})_{ij}$, 
$i,j \! \in \! \lbrace 1,2 \rbrace$, are $\mathcal{O}(1)$ constants, and then the formula for $G_{k}$ will be inverted in 
order to derive the coefficient functions of equation \eqref{eq3.3}, after which, it will be verified that they satisfy all of 
the imposed conditions for this isomonodromy case. The latter procedure gives rise to explicit asymptotic formulae 
for the coefficient functions of equation \eqref{eq3.3}, leading to asymptotics of the solution of the system of 
isomonodromy deformations \eqref{newlax8},\footnote{Via the definitions \eqref{newlax2}, also the asymptotics of the 
solution of the---original---system of isomonodromy deformations \eqref{eq1.4}.} and, in turn, defines asymptotics of the 
solution $u(\tau)$ of the DP3E \eqref{eq1.1} and the related, auxiliary functions $\mathcal{H}(\tau)$, $f_{\pm}(\tau)$, 
$\sigma (\tau)$,\footnote{See the definitions \eqref{eqh1}, \eqref{hatsoff7}, \eqref{pga3}, and \eqref{thmk23}, respectively.} 
and $\hat{\varphi}(\tau)$.
\begin{ccccc} \label{ginversion} 
Let $g_{ij} \! := \! (G_{k})_{ij}$, $i,j \! \in \! \lbrace 1,2 \rbrace$, $k \! = \! \pm 1$, denote the matrix elements of the 
corresponding connection matrices. Assume that all the conditions stated in Theorem \ref{theor3.3.1} are valid. For 
$k \! = \! +1$, let $g_{11}g_{12}g_{21} \! \neq \! 0$ and $g_{22} \! = \! 0$, and, for $k \! = \! -1$, let $g_{12}g_{21}g_{22} 
\! \neq \! 0$ and $g_{11} \! = \! 0$. Then, for $0 \! < \! \delta \! < \! \delta_{k} \! < \! 1/24$, $k \! = \! \pm 1$, the functions 
$v_{0}(\tau)$, $\tilde{r}_{0}(\tau)$,\footnote{See the asymptotics \eqref{tr1} and \eqref{tr3}, respectively.} and $b(\tau)$ 
have the following asymptotics:
\begin{align}
v_{0}(\tau) \! := \! v_{0,k}(\tau) \underset{\tau \to +\infty}{=} \sum_{m=0}^{\infty} \dfrac{\mathfrak{u}_{m}(k)}{(\tau^{1/3})^{m+1}} 
+ &\dfrac{\mi \me^{\mi \pi k/4} \me^{-\mi \pi k/3}(\mathscr{P}_{a})^{k}(s_{0}^{0} \! - \! \mi \me^{-\pi a})}{\sqrt{\smash[b]{2 \pi}} 
\, 3^{1/4}(\varepsilon b)^{1/6}} \me^{-(\beta (\tau)+\mi k \vartheta (\tau))} \nonumber \\
\times& \left(1 \! + \! \mathcal{O} \big(\tau^{-1/3} \big) \right), \label{geek1} \\
\tilde{r}_{0}(\tau) \! := \! \tilde{r}_{0,k}(\tau) \underset{\tau \to +\infty}{=} \sum_{m=0}^{\infty} \dfrac{\mathfrak{r}_{m}
(k)}{(\tau^{1/3})^{m+1}} + &\dfrac{\mi k(\sqrt{3} \! + \! 1)^{k} \me^{\mi \pi k/4} \me^{-\mi \pi k/3}
(\mathscr{P}_{a})^{k}(s_{0}^{0} \! - \! \mi \me^{-\pi a})}{\sqrt{\smash[b]{\pi}} \, 2^{(k-2)/2}3^{1/4}(\varepsilon b)^{1/6}} 
\nonumber \\
\times& \, \me^{-(\beta (\tau)+\mi k \vartheta (\tau))} \! \left(1 \! + \! \mathcal{O} \big(\tau^{-1/3} \big) \right), \label{geek2}
\end{align}
and
\begin{equation} 
\sqrtsign{\smash[b]{b(\tau)}} \underset{\tau \to +\infty}{=} \mathfrak{b}(k)(\varepsilon b)^{1/4} \exp \! \left(\mi 
(a \! - \! \mi/2) \ln \! \big(\alpha_{k}/\sqrtsign{2} \big) \! - \! \frac{\mi a}{6} \ln \tau \! + \! \frac{3k}{4}(\sqrtsign{3} 
\! + \! \mi k)(\varepsilon b)^{1/3} \tau^{2/3} \! + \! \mathcal{O} \big(\tau^{-\delta_{k}} \big) \right), \label{geek3}
\end{equation}
where $\vartheta (\tau)$ and $\beta (\tau)$ are defined in equations \eqref{thmk12},
\begin{gather} 
\mathscr{P}_{a} \! := \! (2 \! + \! \sqrt{3})^{\mi a}, \label{geek4} \\
\mathfrak{b}(k) \! = \! 
\begin{cases} 
g_{11} \me^{\pi a}, &\text{$k \! = \! +1$,} \\
-\big(g_{22} \me^{\pi a} \big)^{-1}, &\text{$k \! = \! -1$,}
\end{cases} \label{eequeb} 
\end{gather}
and the expansion coefficients $\mathfrak{u}_{m}(k)$ (resp., $\mathfrak{r}_{m}(k))$, $m \! \in \! \mathbb{Z}_{+}$, are 
given in equations \eqref{thmk2}--\eqref{thmk10} (resp., \eqref{thmk15} and \eqref{thmk16}$)$.\footnote{Trans-series 
asymptotics (as $\tau \! \to \! +\infty$ with $\varepsilon b \! > \! 0$) for $b(\tau)$ are given in the proof of Theorem 
\ref{pfeetotsa}.}
\end{ccccc}

\emph{Proof}. The scheme of the proof is, \emph{mutatis mutandis}, similar for both cases $(k \! = \! \pm 1)$; therefore, 
without loss of generality, the proof for the case $k \! = \! +1$ is presented: the case $k \! = \! -1$ is proved analogously.

It follows {}from the asymptotics \eqref{tr1}, \eqref{tr3}, and \eqref{prcybk1}, the conditions \eqref{ellinfk2a} 
and \eqref{ellinfk2b}, and the definitions \eqref{ellinfk4} and \eqref{ellohk4} that $p_{1}(\tau) 
\genfrac{}{}{0pt}{3}{=}{\tau \to +\infty} \mathcal{O}(\tau^{1/3} \me^{-\beta (\tau)})$ and $\sqrt{\smash[b]{b(\tau)}} 
\genfrac{}{}{0pt}{3}{=}{\tau \to +\infty} \mathcal{O}(\tau^{-\frac{\mi a}{6}} \me^{\frac{3 \sqrt{3}}{4}(\varepsilon b)^{1/3} 
\tau^{2/3}})$, where $\vartheta (\tau)$ and $\beta (\tau)$ are defined in equations \eqref{thmk12}. {}From the 
definitions \eqref{prpr1}, \eqref{prcy46}, \eqref{prcy54}, \eqref{prcy57}, and \eqref{prcy58}, and the asymptotics 
\eqref{tr1}, \eqref{tr3}, \eqref{prcyak1}, \eqref{prcyomg1}, and \eqref{prcychik1}--\eqref{prcyell2k1}, it follows, 
via a linearisation and inversion argument,\footnote{That is, retaining only those terms that are $\mathcal{O}
(\tau^{-1/3})$.} in conjunction with the latter asymptotics for $p_{1}(\tau)$, that, for $k \! = \! +1$,
\begin{align}
\mathfrak{r}_{0}(1) \tau^{-1/3} \! + \! \mathcal{O} \big(\tau^{-2/3} \big) \underset{\tau \to +\infty}{=}& \, \dfrac{1}{2 \sqrt{3}} 
\! \left(\dfrac{2(a \! - \! \mi/2) \tau^{-1/3}}{\sqrt{3} \alpha_{1}^{2}} \! - \! \dfrac{48 \sqrt{3}(p_{1}(\tau) \! - \! 1)
(\nu (1) \! + \! 1)}{p_{1}(\tau) \tau^{-1/3}} \right. \nonumber \\
-&\left. \, \dfrac{\mi p_{1}(\tau) \tau^{-1/3}}{3 \alpha_{1}^{2}(p_{1}(\tau) \! - \! 1)} \right), \label{geek5} \\
\mathfrak{u}_{0}(1) \tau^{-1/3} \! + \! \mathcal{O} \big(\tau^{-2/3} \big) \underset{\tau \to +\infty}{=}& \, \dfrac{1}{8 \sqrt{3}} 
\! \left(\dfrac{4(a \! - \! \mi/2) \tau^{-1/3}}{\sqrt{3} \alpha_{1}^{2}} \! + \! \dfrac{48 \sqrt{3}(\sqrt{3} \! + \! 1)(p_{1}(\tau) 
\! - \! 1)(\nu (1) \! + \! 1)}{p_{1}(\tau) \tau^{-1/3}} \right. \nonumber \\
+&\left. \, \dfrac{\mi \tau^{-1/3}}{3 \alpha_{1}^{2}} \! \left(\sqrt{3} \! + \! 1 \! - \! \dfrac{(\sqrt{3} \! - \! 1)}{p_{1}(\tau) 
\! - \! 1} \right) \right), \label{geek6}
\end{align}
where
\begin{equation} \label{geek7} 
-\dfrac{(\nu (1) \! + \! 1)}{p_{1}(\tau)} \! = \! \dfrac{q_{1}(\tau)}{2 \mu_{1}(\tau)},
\end{equation}
with
\begin{gather}
q_{1}(\tau) \underset{\tau \to +\infty}{=} c_{q}^{\ast}(1) \tau^{-2/3} \! + \! \mathcal{O} \big(\tau^{-1} \big), \label{geek8} \\
2 \mu_{1}(\tau) \underset{\tau \to +\infty}{=} \mi 8 \sqrt{3} \big(1 \! + \! \mathcal{O}(\tau^{-2/3}) \big), \label{geek9}
\end{gather}
where $c_{q}^{\ast}(1)$ is some to-be-determined coefficient. Recalling {}from Propositions \ref{recursys} and 
\ref{proprr}, respectively, that $\mathfrak{u}_{0}(1) \! = \! a/6 \alpha_{1}^{2}$ and $\mathfrak{r}_{0}(1) \! = \! 
(a \! - \! \mi/2)/3 \alpha_{1}^{2}$, it follows via the asymptotic relations \eqref{geek5} and \eqref{geek6}, equation 
\eqref{geek7}, the asymptotics \eqref{geek8} and \eqref{geek9}, and the asymptotics for $p_{1}(\tau)$ stated above that
\begin{gather}
\dfrac{(a \! - \! \mi/2) \tau^{-1/3}}{3 \alpha_{1}^{2}} \! + \! \mathcal{O} \big(\tau^{-2/3} \big) \underset{\tau \to +\infty}{=} 
\dfrac{\tau^{-1/3}}{2 \sqrt{3}} \! \left(\dfrac{2(a \! - \! \mi/2)}{\sqrt{3} \alpha_{1}^{2}} \! + \! \mi 6c_{q}^{\ast}(1) \right) 
\! + \! \mathcal{O} \big(\tau^{-2/3} \big), \label{geek10} \\
\dfrac{a \tau^{-1/3}}{6 \alpha_{1}^{2}} \! + \! \mathcal{O} \big(\tau^{-2/3} \big) \underset{\tau \to +\infty}{=} 
\dfrac{\tau^{-1/3}}{8 \sqrt{3}} \! \left(\dfrac{4a}{\sqrt{3} \alpha_{1}^{2}} \! - \! \mi 6 (\sqrt{3} \! + \! 1)c_{q}^{\ast}(1) \right) 
\! + \! \mathcal{O} \big(\tau^{-2/3} \big), \label{geek11}
\end{gather}
whence
\begin{equation} \label{geek12} 
c_{q}^{\ast}(1) \! = \! 0.
\end{equation}
Thus, {}from equation \eqref{geek7}, the asymptotics \eqref{geek8} and \eqref{geek9}, the relation \eqref{geek12}, 
and the asymptotics (see above) $p_{1}(\tau) \genfrac{}{}{0pt}{3}{=}{\tau \to +\infty} \mathcal{O}(\tau^{1/3} 
\me^{-\beta (\tau)})$, one deduces that, for $k \! = \! +1$,\footnote{Even though this realisation is not exploited 
in this work, it turns out that $\nu (k) \! + \! 1$ has the asymptotic trans-series expansion
\begin{equation*}
\nu (k) \! + \! 1 \underset{\tau \to +\infty}{=} \sum_{j \in \mathbb{Z}_{+}} \, \sum_{m \in \mathbb{N}} 
\hat{\mathfrak{s}}_{j,k}(m)(\tau^{-1/3})^{j} \big(\me^{-(\beta (\tau)+\mi k \vartheta (\tau))} \big)^{m}, \quad k \! = \! \pm 1,
\end{equation*}
for certain coefficients $\hat{\mathfrak{s}}_{j,k}(m) \colon \mathbb{Z}_{+} \times \{\pm 1\} \times \mathbb{N} \! \to \! 
\mathbb{C}$, where, in particular, $\hat{\mathfrak{s}}_{0,k}(1) \! = \! \hat{\mathfrak{s}}_{1,k}(1) \! = \! 0$.}
\begin{equation} \label{geek13} 
\nu (1) \! + \! 1 \underset{\tau \to +\infty}{=} \mathcal{O} \big(\tau^{-2/3} \me^{-\beta (\tau)} \big).
\end{equation}

{}From the corresponding $(k \! = \! +1)$ asymptotics \eqref{tr1} and \eqref{tr3}, the definitions \eqref{eqpeetee}, 
\eqref{ellinfk5}, and \eqref{ellohk5}, the expansion $\me^{z} \! = \! \sum_{m=0}^{\infty} \tfrac{z^{m}}{m!}$, and the 
leading-order asymptotics \eqref{geek9} and \eqref{geek13}, one shows that, for $k \! = \! +1$,
\begin{gather}
\me^{\pm \Delta \tilde{\mathfrak{z}}_{1}(\tau)} \underset{\tau \to +\infty}{=} 1 \! + \! \tau^{-2/3} \sum_{m=0}^{\infty} 
\tilde{\zeta}_{m}^{\pm}(1)(\tau^{-1/3})^{m} \! + \! \mathcal{O} \big(\tau^{-1/3} \me^{-\beta (\tau)} \big), \label{geek14} \\
\me^{\pm \Delta \hat{\mathfrak{z}}_{1}(\tau)} \underset{\tau \to +\infty}{=} 1 \! + \! \tau^{-2/3} \sum_{m=0}^{\infty} 
\hat{\zeta}_{m}^{\pm}(1)(\tau^{-1/3})^{m} \! + \! \mathcal{O} \big(\tau^{-1/3} \me^{-\beta (\tau)} \big), \label{geek15}
\end{gather}
for $\mathcal{O}(1)$ coefficients $\tilde{\zeta}_{m}^{\pm}(1)$ and $\hat{\zeta}_{m}^{\pm}(1)$. {}From the corresponding 
$(k \! = \! +1)$ asymptotics \eqref{tr1}, \eqref{tr3}, \eqref{prcyzeek1}, \eqref{prcyomg1}, \eqref{prcychik1}, and 
\eqref{prcymuk1}, the definition \eqref{prcy57}, and $p_{1}(\tau) \genfrac{}{}{0pt}{3}{=}{\tau \to +\infty} \mathcal{O}
(\tau^{1/3} \me^{-\beta (\tau)})$, it follows that, for $k \! = \! +1$,
\begin{equation} \label{geek16} 
\dfrac{1}{(2 \mu_{1}(\tau))^{1/2}} \underset{\tau \to +\infty}{=} \dfrac{\me^{-\mi \pi/4}}{2^{3/2}3^{1/4}} \! \left(
1 \! + \! \tau^{-2/3} \sum_{m=0}^{\infty} \alpha_{m}^{\sharp}(1)(\tau^{-1/3})^{m} \! + \! \mathcal{O} 
\big(\tau^{-1/3} \me^{-\beta (\tau)} \big) \right),
\end{equation}
for $\mathcal{O}(1)$ coefficients $\alpha_{m}^{\sharp}(1)$. {}From the corresponding $(k \! = \! +1)$ asymptotics \eqref{tr1}, 
\eqref{tr3}, \eqref{prcyzeek1}, \eqref{prcyg4}--\eqref{prcybk1}, \eqref{prcyomg1}, and \eqref{prcychik1}, and the definitions 
\eqref{prcy57}, \eqref{ellinfk6}--\eqref{ellinfk8}, and \eqref{ellohk6}--\eqref{ellohk8}, one shows that (cf. Lemmata 
\ref{linfnewlemm} and \ref{lzernewlemm}), for $k \! = \! \pm 1$, to leading order,
\begin{gather}
\mathbb{E}_{{\scriptscriptstyle \mathcal{N}},k}^{\infty}(\tau) \underset{\tau \to +\infty}{=} 
\begin{pmatrix}
\mathcal{O} \big(\tau^{-2/3} \big) & \mathcal{O} \! \left(\tau^{-1/3}(\me^{-\beta (\tau)})^{\frac{1+k}{2}} \right) \\
\mathcal{O} \! \left(\tau^{-1/3}(\me^{-\beta (\tau)})^{\frac{1-k}{2}} \right) & \mathcal{O} \big(\tau^{-2/3} \big)
\end{pmatrix}, \label{geek21} \\
\mathbb{E}_{{\scriptscriptstyle \mathcal{N}},k}^{0}(\tau) \underset{\tau \to +\infty}{=} 
\begin{pmatrix}
\mathcal{O} \big(\tau^{-2/3} \big) & \mathcal{O} \! \left(\tau^{-1/3}(\me^{-\beta (\tau)})^{\frac{1-k}{2}} \right) \\
\mathcal{O} \! \left(\tau^{-1/3}(\me^{-\beta (\tau)})^{\frac{1+k}{2}} \right) & \mathcal{O} \big(\tau^{-2/3} \big)
\end{pmatrix}, \label{geek22} 
\end{gather}
whence, via the asymptotics \eqref{geek13}, \eqref{geek21}, and \eqref{geek22}, and the above asymptotics for $p_{1}
(\tau)$, it follows via the relation $\det (\mathrm{I} \! + \! \mathbb{J}) \! = \! 1 \! + \! \tr (\mathbb{J}) \! + \! \det (\mathbb{J})$, 
$\mathbb{J} \! \in \! \mathrm{M}_{2}(\mathbb{C})$, that, for $k \! = \! +1$, to all orders,
\begin{gather}
\mathrm{I} \! + \! \mathbb{E}_{{\scriptscriptstyle \mathcal{N}},1}^{\infty}(\tau) \underset{\tau \to +\infty}{=} \mathrm{I} \! + \! 
\sum_{m=1}^{\infty} \zeta_{m}^{\flat}(1)(\tau^{-1/3})^{m} \! + \! \mathcal{O} \big(\tau^{-1/3} \me^{-\beta (\tau)} \sigma_{+} 
\big), \label{geek23} \\
(\mathrm{I} \! + \! \mathbb{E}_{{\scriptscriptstyle \mathcal{N}},1}^{0}(\tau))^{-1} \underset{\tau \to +\infty}{=} \mathrm{I} 
\! + \! \sum_{m=1}^{\infty} \zeta_{m}^{\natural}(1)(\tau^{-1/3})^{m} \! + \! \mathcal{O} \big(\tau^{-1/3} \me^{-\beta (\tau)} 
\sigma_{-} \big), \label{geek24}
\end{gather}
for $\mathrm{M}_{2}(\mathbb{C})$-valued, $\mathcal{O}(1)$ coefficients $\zeta_{m}^{\flat}(1)$ and $\zeta_{m}^{\natural}(1)$. 
It now follows {}from the corresponding $(k \! = \! +1)$ conditions \eqref{ellinfk2a} and \eqref{ellinfk2b}, that is, $p_{1}
(\tau) \mathfrak{B}_{1} \genfrac{}{}{0pt}{3}{=}{\tau \to +\infty} \mathcal{O}(\me^{2 \hat{\mathfrak{z}}_{1}^{0}(\tau)})$ and 
$\sqrt{\smash[b]{b(\tau)}} \genfrac{}{}{0pt}{3}{=}{\tau \to +\infty} \mathcal{O}(\me^{\tilde{\mathfrak{z}}_{1}^{0}(\tau)-
\hat{\mathfrak{z}}_{1}^{0}(\tau)})$, respectively, where $\tilde{\mathfrak{z}}_{1}^{0}(\tau)$ and $\hat{\mathfrak{z}}_{1}^{0}(\tau)$ 
are defined by equations \eqref{ellinfk4} and \eqref{ellohk4}, respectively, the expansion $\me^{z} \! = \! \sum_{m=0}^{\infty} 
\tfrac{z^{m}}{m!}$, the reflection formula $\Gamma (z) \Gamma (1 \! - \! z) \! = \! \pi/\sin \pi z$, the definitions 
\eqref{jeek5}--\eqref{jeek8}, and the asymptotics \eqref{geek13} and \eqref{geek16}, that, for $k \! = \! +1$,
\begin{equation} \label{geek25} 
\hat{\mathbb{G}}(1) \! := \! 
\begin{pmatrix}
\hat{\mathbb{G}}_{11}(1) & \hat{\mathbb{G}}_{12}(1) \\
\hat{\mathbb{G}}_{21}(1) & \hat{\mathbb{G}}_{22}(1)
\end{pmatrix} \underset{\tau \to +\infty}{=} 
\begin{pmatrix}
\mathcal{O}(1) & \mathcal{O}(1) \\
\mathcal{O}(1) & \mathcal{O}(\nu (1) \! + \! 1)
\end{pmatrix},
\end{equation}
and, {}from equation \eqref{jeek4}, the definitions \eqref{ellinfk6}, \eqref{ellinfk7}, \eqref{ellohk6}, and \eqref{ellohk7}, and 
the asymptotics \eqref{geek14}, \eqref{geek15}, and \eqref{geek25},
\begin{equation} \label{geek26} 
\mathrm{G}^{\ast}(1) \underset{\tau \to +\infty}{=} 
\begin{pmatrix}
\mathcal{O}(1) & \mathcal{O}(1) \\
\mathcal{O}(1) & \mathcal{O}(\nu (1) \! + \! 1)
\end{pmatrix}, 
\end{equation} 
whence, via the definitions \eqref{ellohk9}, \eqref{jeek2}, and \eqref{jeek3}, and the asymptotics \eqref{geek23} and 
\eqref{geek24},
\begin{gather}
\widetilde{G}(1) \underset{\tau \to +\infty}{=} 
\begin{pmatrix}
\mathcal{O}(1) & \mathcal{O}(1) \\
\mathcal{O}(1) & \mathcal{O}(\nu (1) \! + \! 1)
\end{pmatrix}, \label{geek27} \\
\widehat{\mathscr{G}}(1) \underset{\tau \to +\infty}{=} 
\begin{pmatrix}
\mathcal{O}(1) & \mathcal{O}(1) \\
\mathcal{O}(1) & \mathcal{O}(1)
\end{pmatrix}. \label{geek28}
\end{gather}
{}From the asymptotics \eqref{ellinfk10} and \eqref{ellohk11}, the definition \eqref{jeek9}, the asymptotics 
\eqref{geek27} and \eqref{geek28}, and the relations $\max \lbrace z_{1},z_{2} \rbrace \! = \! (z_{1} \! + \! z_{2} 
\! + \! \lvert z_{1} \! - \! z_{2} \rvert)/2$, $\min \lbrace z_{1},z_{2} \rbrace \! = \! (z_{1} \! + \! z_{2} \! - \! \lvert z_{1} 
\! - \! z_{2} \rvert)/2$, $z_{1},z_{2} \! \in \! \mathbb{R}$, and $\max_{k= \pm 1} \lbrace 3 \delta_{k} \! - \! 1/3,
-\delta_{k} \! - \! (1 \! + \! k)/6,-\delta_{k} \! - \! (1 \! - \! k)/6 \rbrace \! = \! -\delta_{k}$, it follows that, for $k \! = \! +1$,
\begin{equation} \label{geek29} 
\mathbb{E}_{1}^{\scriptscriptstyle G_{1}}(\tau) \underset{\tau \to +\infty}{=} \mathcal{O} 
\big(\tau^{-\delta_{1}} \big).
\end{equation}
Finally, {}from the asymptotics \eqref{jeek1} and \eqref{geek27}--\eqref{geek29}, one arrives at $(G_{1})_{i,j=1,2} 
\genfrac{}{}{0pt}{3}{=}{\tau \to +\infty} \mathcal{O}(1)$ (for $k \! = \! +1$), which is, in fact, the isomonodromy 
condition for the corresponding connection matrix.

{}From the definition \eqref{ellohk9}, the asymptotics \eqref{jeek1}, the definitions \eqref{jeek2} and \eqref{jeek3}, 
equation \eqref{jeek4}, the definitions \eqref{jeek5}--\eqref{jeek8}, the asymptotics \eqref{geek23}, \eqref{geek24}, 
and \eqref{geek29}, and the isomonodromy condition for the corresponding connection matrix $G_{1}$, it follows 
that, for $k \! = \! +1$, upon setting $g_{ij} \! := \! (G_{1})_{ij}$, $i,j \! \in \! \lbrace 1,2 \rbrace$,
\begin{equation} \label{geek30} 
\begin{pmatrix}
g_{11} & g_{12} \\
g_{21} & g_{22}
\end{pmatrix} \underset{\tau \to +\infty}{=} 
\begin{pmatrix}
1 & s_{0}^{0} \\
0 & 1
\end{pmatrix} \! 
\begin{pmatrix}
\mathrm{G}^{\ast}_{11}(1) & \mathrm{G}^{\ast}_{12}(1) \\
\mathrm{G}^{\ast}_{21}(1) & \mathrm{G}^{\ast}_{22}(1)
\end{pmatrix} \! 
\begin{pmatrix}
1 \! + \! \eta_{11}(\tau) & \eta_{12}(\tau) \\
\eta_{21}(\tau) & 1 \! + \! \eta_{22}(\tau)
\end{pmatrix} \! \left(\mathrm{I} \! + \! \mathcal{O}(\tau^{-\delta_{1}}) \right),
\end{equation}
where
\begin{equation} \label{geek31} 
\eta_{ij}(\tau) \underset{\tau \to +\infty}{:=} \sum_{m=1}^{\infty}(\mathbb{H}_{m}(1))_{ij}(\tau^{-1/3})^{m} 
\! + \! \mathcal{O} \big(\tau^{-1/3} \me^{-\beta (\tau)} \big), \quad i,j \! \in \! \lbrace 1,2 \rbrace,
\end{equation}
for $\mathcal{O}(1)$ coefficients $(\mathbb{H}_{m}(1))_{ij}$. It follows {}from the asymptotics \eqref{geek30} that
\begin{align} \label{geek32} 
g_{12}g_{21} \underset{\tau \to \infty}{=}& \, \left(\mathrm{G}_{21}^{\ast}(1)(1 \! + \! \eta_{11}(\tau)) 
\! + \! \mathrm{G}_{22}^{\ast}(1) \eta_{21}(\tau) \right) \! \left(\mathrm{G}_{12}^{\ast}(1) \! + \! 
s_{0}^{0} \mathrm{G}_{22}^{\ast}(1) \right. \nonumber \\
+&\left. \, (\mathrm{G}_{12}^{\ast}(1) \! + \! s_{0}^{0} \mathrm{G}_{22}^{\ast}(1)) \eta_{22}(\tau) \! + \! 
(\mathrm{G}_{11}^{\ast}(1) \! + \! s_{0}^{0} \mathrm{G}_{21}^{\ast}(1)) \eta_{12}(\tau) \right) \!
 \left(1 \! + \! \mathcal{O}(\tau^{-\delta_{1}}) \right).
\end{align}
{}From the corresponding $(k \! = \! +1)$ conditions~\eqref{ellinfk2a} and \eqref{ellinfk2b}, that is, $p_{1}(\tau) 
\mathfrak{B}_{1} \genfrac{}{}{0pt}{3}{=}{\tau \to +\infty} \mathcal{O}(\me^{2 \hat{\mathfrak{z}}_{1}^{0}(\tau)})$ and 
$\sqrt{\smash[b]{b(\tau)}} \genfrac{}{}{0pt}{3}{=}{\tau \to +\infty} \mathcal{O}(\me^{\tilde{\mathfrak{z}}_{1}^{0}
(\tau)-\hat{\mathfrak{z}}_{1}^{0}(\tau)})$, respectively, where $\tilde{\mathfrak{z}}_{1}^{0}(\tau)$ and 
$\hat{\mathfrak{z}}_{1}^{0}(\tau)$ are defined by equations \eqref{ellinfk4} and \eqref{ellohk4}, respectively, equation 
\eqref{jeek4}, the definitions \eqref{jeek5}--\eqref{jeek8}, the expansion $\me^{z} \! = \! \sum_{m=0}^{\infty} 
\tfrac{z^{m}}{m!}$, the asymptotics \eqref{geek13}--\eqref{geek16}, and the definitions \eqref{ellinfk6}, 
\eqref{ellinfk7}, \eqref{ellohk6}, and \eqref{ellohk7}, one shows that, for $k \! = \! +1$,
\begin{align}
\mathrm{G}_{21}^{\ast}(1) \eta_{11}(\tau) \! =&\, \eta_{11}(\tau) \dfrac{\hat{\mathbb{G}}_{21}(1) 
\hat{\mathbb{B}}_{0}^{\infty}(\tau)}{\hat{\mathbb{B}}_{0}^{0}(\tau)} \me^{-\Delta \tilde{\mathfrak{z}}_{1}
(\tau) + \Delta \hat{\mathfrak{z}}_{1}(\tau)} \underset{\tau \to +\infty}{=} \mathcal{O} \big(\tau^{-1/3} \big), 
\label{geek33} \\
\mathrm{G}_{22}^{\ast}(1) \eta_{21}(\tau) \! =& \, \eta_{21}(\tau) \dfrac{\hat{\mathbb{G}}_{22}(1) 
\hat{\mathbb{A}}_{0}^{\infty}(\tau)}{\hat{\mathbb{B}}_{0}^{0}(\tau)} \me^{\Delta \tilde{\mathfrak{z}}_{1}(\tau) + 
\Delta \hat{\mathfrak{z}}_{1}(\tau)} \underset{\tau \to +\infty}{=} \mathcal{O} \big(\tau^{-1} \me^{-\beta (\tau)} 
\big), \label{geek34}
\end{align}
\begin{align}
(\mathrm{G}_{12}^{\ast}(1) \! + \! s_{0}^{0} \mathrm{G}_{22}^{\ast}(1)) \eta_{22}(\tau) =& \, \eta_{22}(\tau) \! 
\left(\dfrac{\hat{\mathbb{G}}_{12}(1) \hat{\mathbb{A}}_{0}^{\infty}(\tau)}{\hat{\mathbb{A}}_{0}^{0}(\tau)} 
\me^{\Delta \tilde{\mathfrak{z}}_{1}(\tau)-\Delta \hat{\mathfrak{z}}_{1}(\tau)} \! + \! s_{0}^{0} \dfrac{
\hat{\mathbb{G}}_{22}(1) \hat{\mathbb{A}}_{0}^{\infty}(\tau)}{\hat{\mathbb{B}}_{0}^{0}(\tau)} \right. 
\nonumber \\
\times&\left. \me^{\Delta \tilde{\mathfrak{z}}_{1}(\tau)+\Delta \hat{\mathfrak{z}}_{1}(\tau)} \right) \! 
\underset{\tau \to +\infty}{=} \! \mathcal{O}(\tau^{-1/3}) \big(\mathcal{O}(1) \! + \! \mathcal{O}(\tau^{-2/3} 
\me^{-\beta (\tau)}) \big) \! \underset{\tau \to +\infty}{=} \! \mathcal{O} \big(\tau^{-1/3} \big), \label{geek35} \\
(\mathrm{G}_{11}^{\ast}(1) \! + \! s_{0}^{0} \mathrm{G}_{21}^{\ast}(1)) \eta_{12}(\tau) 
=& \, \eta_{12}(\tau) \! \left(\dfrac{\hat{\mathbb{G}}_{11}(1) \hat{\mathbb{B}}_{0}^{\infty}
(\tau)}{\hat{\mathbb{A}}_{0}^{0}(\tau)} \me^{-\Delta \tilde{\mathfrak{z}}_{1}(\tau)-
\Delta \hat{\mathfrak{z}}_{1}(\tau)} \! + \! s_{0}^{0} \dfrac{\hat{\mathbb{G}}_{21}(1) 
\hat{\mathbb{B}}_{0}^{\infty}(\tau)}{\hat{\mathbb{B}}_{0}^{0}(\tau)} \right. \nonumber \\
\times&\left. \me^{-\Delta \tilde{\mathfrak{z}}_{1}(\tau)+\Delta \hat{\mathfrak{z}}_{1}(\tau)} \right) \! 
\underset{\tau \to +\infty}{=} \! \mathcal{O}(\tau^{-1/3}) \big(\mathcal{O}(1) \! + \! \mathcal{O}(1)) \! 
\underset{\tau \to +\infty}{=} \! \mathcal{O} \big(\tau^{-1/3} \big), \label{geek36}
\end{align}
whence (cf. asymptotics~\eqref{geek32})
\begin{align}
g_{12}g_{21} \underset{\tau \to +\infty}{=}& \, \left(\mathrm{G}_{21}^{\ast}(1) \! + \! \mathcal{O}(\tau^{-1/3}) 
\! + \! \mathcal{O}(\tau^{-1} \me^{-\beta (\tau)}) \right) \! \left(1 \! + \! \mathcal{O}(\tau^{-\delta_{1}}) \right) 
\nonumber \\
\times& \, \left(\mathrm{G}_{12}^{\ast}(1) \! + \! \mathcal{O}(\tau^{-1/3}) \! + \! \mathcal{O}(\tau^{-2/3} 
\me^{-\beta (\tau)}) \right) \nonumber \\
\underset{\tau \to +\infty}{=}& \, \mathrm{G}_{12}^{\ast}(1) \mathrm{G}_{21}^{\ast}(1) \! \left(1 \! + \! \mathcal{O}
(\tau^{-\delta_{1}}) \right) \underset{\tau \to +\infty}{=} \hat{\mathbb{G}}_{12}(1) \hat{\mathbb{G}}_{21}(1) 
\dfrac{\hat{\mathbb{A}}_{0}^{\infty}(\tau) \hat{\mathbb{B}}_{0}^{\infty}(\tau)}{\hat{\mathbb{A}}_{0}^{0}(\tau) 
\hat{\mathbb{B}}_{0}^{0}(\tau)} \! \left(1 \! + \! \mathcal{O}(\tau^{-\delta_{1}}) \right) \nonumber \\
\underset{\tau \to +\infty}{=}& \, -\me^{-2 \pi \mi (\nu (1)+1)} \big(1 \! + \! \mathcal{O}(\tau^{-2/3}) \big) \big(
1 \! + \! \mathcal{O}(\tau^{-\delta_{1}}) \big) \underset{\tau \to +\infty}{=} -(1 \! + \! \mathcal{O}(\nu (1) \! + \! 1)) 
\nonumber \\
\times& \, \big(1 \! + \! \mathcal{O}(\tau^{-\delta_{1}}) \big) \underset{\tau \to +\infty}{=} -\big(1 \! + \! \mathcal{O}
(\tau^{-2/3} \me^{-\beta (\tau)}) \big) \big(1 \! + \! \mathcal{O}(\tau^{-\delta_{1}}) \big) \quad \Rightarrow \nonumber \\
-g_{12}g_{21} \underset{\tau \to +\infty}{=}& \, 1 \! + \! \mathcal{O} \big(\tau^{-\delta_{1}} \big); 
\label{geek37}
\end{align}
analogously,
\begin{align}
g_{21} \underset{\tau \to +\infty}{=}& \, \left(\mathrm{G}_{21}^{\ast}(1)(1 \! + \! \eta_{11}(\tau)) \! + \! 
\mathrm{G}_{22}^{\ast}(1) \eta_{21}(\tau) \right) \! \left(1 \! + \! \mathcal{O}(\tau^{-\delta_{1}}) \right) 
\nonumber \\
\underset{\tau \to +\infty}{=}& \, \left(\mathrm{G}_{21}^{\ast}(1) \! + \! \mathcal{O}(\tau^{-1/3}) \! + \! 
\mathcal{O}(\tau^{-1} \me^{-\beta (\tau)}) \right) \! \left(1 \! + \! \mathcal{O}(\tau^{-\delta_{1}}) \right) 
\nonumber \\
\underset{\tau \to +\infty}{=}& \, \mathrm{G}_{21}^{\ast}(1) \big(1 \! + \! \mathcal{O}(\tau^{-\delta_{1}}) \big) 
\underset{\tau \to +\infty}{=} \hat{\mathbb{G}}_{21}(1) \dfrac{\hat{\mathbb{B}}_{0}^{\infty}(\tau)}{
\hat{\mathbb{B}}_{0}^{0}(\tau)} \me^{-\Delta \tilde{\mathfrak{z}}_{1}(\tau)+\Delta \hat{\mathfrak{z}}_{1}(\tau)} 
\big(1 \! + \! \mathcal{O}(\tau^{-\delta_{1}}) \big) \nonumber \\
\underset{\tau \to +\infty}{=}& \, -\dfrac{\mi \sqrt{\smash[b]{b(\tau)}}}{(\varepsilon b)^{1/4}} 
\me^{-\tilde{\mathfrak{z}}_{1}^{0}(\tau)+\hat{\mathfrak{z}}_{1}^{0}(\tau)} \me^{-2 \pi \mi (\nu (1)+1)} \big(
1 \! + \! \mathcal{O}(\tau^{-2/3}) \big) \big(1 \! + \! \mathcal{O}(\tau^{-2/3}) \big) \big(1 \! + \! \mathcal{O}
(\tau^{-\delta_{1}}) \big) \nonumber \\
\underset{\tau \to +\infty}{=}& \, -\dfrac{\mi \sqrt{\smash[b]{b(\tau)}}}{(\varepsilon b)^{1/4}} 
\me^{-\tilde{\mathfrak{z}}_{1}^{0}(\tau)+\hat{\mathfrak{z}}_{1}^{0}(\tau)} \big(1 \! + \! \mathcal{O}(\nu (1) 
\! + \! 1) \big) \big(1 \! + \! \mathcal{O}(\tau^{-\delta_{1}}) \big) \nonumber \\
\underset{\tau \to +\infty}{=}& \, -\dfrac{\mi \sqrt{\smash[b]{b(\tau)}}}{(\varepsilon b)^{1/4}} 
\me^{-\tilde{\mathfrak{z}}_{1}^{0}(\tau)+\hat{\mathfrak{z}}_{1}^{0}(\tau)} \big(1 \! + \! \mathcal{O}(\tau^{-2/3} 
\me^{-\beta (\tau)}) \big) \big(1 \! + \! \mathcal{O}(\tau^{-\delta_{1}}) \big) \quad \Rightarrow \nonumber \\
g_{21} \underset{\tau \to +\infty}{=}& \, -\dfrac{\mi \sqrt{\smash[b]{b(\tau)}}}{(\varepsilon b)^{1/4}} 
\me^{-\tilde{\mathfrak{z}}_{1}^{0}(\tau)+\hat{\mathfrak{z}}_{1}^{0}(\tau)} \big(1 \! + \! \mathcal{O}
(\tau^{-\delta_{1}}) \big). \label{geek38}
\end{align}
It follows, upon inversion, {}from the asymptotics \eqref{geek37} and \eqref{geek38} that, for $k \! = \! +1$,
\begin{equation} \label{geek39} 
\sqrt{\smash[b]{b(\tau)}} \underset{\tau \to +\infty}{=} \mi g_{21}(\varepsilon b)^{1/4} \me^{\tilde{\mathfrak{z}}_{1}^{0}
(\tau)-\hat{\mathfrak{z}}_{1}^{0}(\tau)} \big(1 \! + \! \mathcal{O}(\tau^{-\delta_{1}}) \big) \underset{\tau \to +\infty}{=} 
-\mi g_{12}^{-1}(\varepsilon b)^{1/4} \me^{\tilde{\mathfrak{z}}_{1}^{0}(\tau)-\hat{\mathfrak{z}}_{1}^{0}(\tau)} \big(
1 \! + \! \mathcal{O}(\tau^{-\delta_{1}}) \big),
\end{equation}
whence, via equations \eqref{monok2} and the definitions \eqref{ellinfk4} and \eqref{ellohk4}, one arrives at the 
corresponding $(k \! = \! +1)$ asymptotics for $\sqrt{\smash[b]{b(\tau)}}$ stated in equation \eqref{geek3} of the 
lemma.\footnote{Note that the asymptotics \eqref{geek39} is consistent with the corresponding $(k \! = \! +1)$ 
condition \eqref{ellinfk2b}.}

Recall the following formula (cf. equations \eqref{monoeqns}), which is one of the defining relations for the manifold 
of the monodromy data, $\mathscr{M}$:
\begin{equation} \label{geek40} 
g_{21}g_{22} \! - \! g_{11}g_{12} \! + \! s_{0}^{0}g_{11}g_{22} \! = \! \mi \me^{-\pi a}.
\end{equation}
Let
\begin{equation} \label{geek52} 
x \! := \! \dfrac{\sqrt{2 \pi} \, p_{1}(\tau) \mathfrak{B}_{1} \me^{-2 \hat{\mathfrak{z}}_{1}^{0}(\tau)} 
\me^{\mi \pi (\nu (1)+1)}}{(2 \! + \! \sqrt{3})^{1/2}(2 \mu_{1}(\tau))^{1/2} \Gamma (-\nu (1))} 
\dfrac{\hat{\mathbb{A}}_{0}^{\infty}(\tau) \hat{\mathbb{B}}_{0}^{\infty}(\tau)}{\hat{\mathbb{A}}_{0}^{0}
(\tau) \hat{\mathbb{A}}_{0}^{0}(\tau)} \me^{-2 \Delta \hat{\mathfrak{z}}_{1}(\tau)}(1 \! + \! \eta_{11}
(\tau))(1 \! + \! \eta_{22}(\tau)):
\end{equation}
substituting equation \eqref{jeek4}, the definitions \eqref{jeek5}--\eqref{jeek8}, and the asymptotics \eqref{geek30} 
into equation \eqref{geek40}, an algebraic exercise reveals that, in terms of the newly-defined variable $x$, it can 
be recast in the form
\begin{equation} \label{geek53} 
y_{1}x^{-2} \! + \! (y_{2} \! + \! y_{3} \! + \! y_{4})x^{-1} \! + \! (1 \! + \! y_{5} \! + \! y_{6})x \! + \! y_{7}x^{2} \! + \! 
y_{8} \! + \! y_{9} \! + \! y_{10} \! + \! y_{11} \! - \! \mi \me^{-\pi a} \! + \! \mathcal{O} \big(\tau^{-\delta_{1}} \big) 
\underset{\tau \to +\infty}{=} 0,
\end{equation}
where
\begin{align}
y_{1} \! :=& \, \left(\mi 2g_{21}^{-1} \sin (\pi (\nu (1) \! + \! 1)) \me^{-\mi \pi (\nu (1)+1)} 
\dfrac{\hat{\mathbb{A}}_{0}^{\infty}(\tau) \hat{\mathbb{B}}_{0}^{\infty}(\tau)}{\hat{\mathbb{A}}_{0}^{0}(\tau) 
\hat{\mathbb{A}}_{0}^{0}(\tau)} \right)^{2} \! \left(\dfrac{\hat{\mathbb{A}}_{0}^{\infty}(\tau)}{\hat{\mathbb{B}}_{
0}^{0}(\tau)} \right)^{2} \me^{2(\Delta \tilde{\mathfrak{z}}_{1}(\tau)-\Delta \hat{\mathfrak{z}}_{1}(\tau))} \nonumber \\
\times& \, (1 \! + \! \eta_{11}(\tau))^{2}(1 \! + \! \eta_{22}(\tau))^{3} \eta_{21}(\tau) \big(1 \! + \! \mathcal{O}
(\tau^{-\delta_{1}}) \big), \label{geek59} \\
y_{2} \! :=& \, \mi 2 \sin (\pi (\nu (1) \! + \! 1)) \me^{-\mi 3 \pi (\nu (1)+1)} \! \left(\dfrac{\hat{\mathbb{A}}_{0}^{\infty}
(\tau) \hat{\mathbb{B}}_{0}^{\infty}(\tau)}{\hat{\mathbb{A}}_{0}^{0}(\tau) \hat{\mathbb{B}}_{0}^{0}(\tau)}
(1 \! + \! \eta_{11}(\tau))(1 \! + \! \eta_{22}(\tau)) \right)^{2}, \label{geek54} \\
y_{3} \! :=& \, \mi 2 s_{0}^{0}g_{21}^{-2} \sin (\pi (\nu (1) \! + \! 1)) \me^{-\mi \pi (\nu (1)+1)} \! \left(
\dfrac{\hat{\mathbb{A}}_{0}^{\infty}(\tau)}{\hat{\mathbb{A}}_{0}^{0}(\tau)} \right)^{3} \dfrac{\hat{\mathbb{B}}_{0}^{\infty}
(\tau)}{\hat{\mathbb{B}}_{0}^{0}(\tau)} \me^{2(\Delta \tilde{\mathfrak{z}}_{1}(\tau)-\Delta \hat{\mathfrak{z}}_{1}(\tau))} 
\nonumber \\
\times& \, (1 \! + \! \eta_{11}(\tau))(1 \! + \! \eta_{22}(\tau))^{2} \eta_{21}(\tau), \label{geek61} \\
y_{4} \! :=& \, \mi 2 \sin (\pi (\nu (1) \! + \! 1)) \me^{-\mi 3 \pi (\nu (1)+1)} \! \left(\dfrac{\hat{\mathbb{A}}_{0}^{\infty}
(\tau) \hat{\mathbb{B}}_{0}^{\infty}(\tau)}{\hat{\mathbb{A}}_{0}^{0}(\tau) \hat{\mathbb{B}}_{0}^{0}(\tau)} \right)^{2}
(1 \! + \! \eta_{11}(\tau)) (1 \! + \! \eta_{22}(\tau)) \eta_{12}(\tau) \eta_{21}(\tau), \label{geek62} \\
y_{5} \! :=& \, -s_{0}^{0}g_{21}^{2} \dfrac{\hat{\mathbb{A}}_{0}^{0}(\tau) \hat{\mathbb{B}}_{0}^{\infty}
(\tau)}{\hat{\mathbb{A}}_{0}^{\infty}(\tau) \hat{\mathbb{B}}_{0}^{0}(\tau)} \me^{-2(\Delta \tilde{\mathfrak{z}}_{1}
(\tau)-\Delta \hat{\mathfrak{z}}_{1}(\tau))} \dfrac{\eta_{12}(\tau)}{1 \! + \! \eta_{22}(\tau)} \big(1 \! + \! \mathcal{O}
(\tau^{-\delta_{1}}) \big), \label{geek58} \\
y_{6} \! :=& \, \dfrac{\eta_{12}(\tau) \eta_{21}(\tau)}{(1 \! + \! \eta_{11}(\tau))(1 \! + \! \eta_{22}(\tau))}, \label{geek63} \\
y_{7} \! :=& \, -g_{21}^{2} \! \left(\dfrac{\hat{\mathbb{A}}_{0}^{0}(\tau)}{\hat{\mathbb{A}}_{0}^{\infty}(\tau)} \right)^{2} 
\me^{-2(\Delta \tilde{\mathfrak{z}}_{1}(\tau)-\Delta \hat{\mathfrak{z}}_{1}(\tau))} \dfrac{\eta_{12}(\tau)}{(1 \! + \! 
\eta_{11}(\tau))(1 \! + \! \eta_{22}(\tau))^{2}} \big(1 \! + \! \mathcal{O}(\tau^{-\delta_{1}}) \big), \label{geek57} \\
y_{8} \! :=& \, s_{0}^{0} \me^{-\mi 2 \pi (\nu (1)+1)} \dfrac{\hat{\mathbb{A}}_{0}^{\infty}(\tau) 
\hat{\mathbb{B}}_{0}^{\infty}(\tau)}{\hat{\mathbb{A}}_{0}^{0}(\tau) \hat{\mathbb{B}}_{0}^{0}(\tau)}
(1 \! + \! \eta_{11}(\tau))(1 \! + \! \eta_{22}(\tau)), \label{geek55} \\
y_{9} \! :=& \, g_{21}^{2} \! \left(\dfrac{\hat{\mathbb{B}}_{0}^{\infty}(\tau)}{\hat{\mathbb{B}}_{0}^{0}(\tau)} \right)^{2} 
\me^{-2(\Delta \tilde{\mathfrak{z}}_{1}(\tau)-\Delta \hat{\mathfrak{z}}_{1}(\tau))}(1 \! + \! \eta_{11}(\tau)) \eta_{12}
(\tau) \big(1 \! + \! \mathcal{O}(\tau^{-\delta_{1}}) \big), \label{geek56} \\
y_{10} \! :=& \, -g_{21}^{-2} \! \left(\dfrac{\hat{\mathbb{A}}_{0}^{\infty}(\tau)}{\hat{\mathbb{A}}_{0}^{0}(\tau)} \right)^{2} 
\me^{2(\Delta \tilde{\mathfrak{z}}_{1}(\tau)-\Delta \hat{\mathfrak{z}}_{1}(\tau))}(1 \! + \! \eta_{22}(\tau)) \eta_{21}(\tau) 
\big(1 \! + \! \mathcal{O}(\tau^{-\delta_{1}}) \big), \label{geek60} \\
y_{11} \! :=& \, -\mi 2 s_{0}^{0} \sin (\pi (\nu (1) \! + \! 1)) \me^{-\mi \pi (\nu (1)+1)} \dfrac{\hat{\mathbb{A}}_{0}^{\infty}
(\tau) \hat{\mathbb{B}}_{0}^{\infty}(\tau)}{\hat{\mathbb{A}}_{0}^{0}(\tau) \hat{\mathbb{B}}_{0}^{0}(\tau)} \eta_{12}
(\tau) \eta_{21}(\tau). \label{geek64}
\end{align}
Via the asymptotics \eqref{geek13}--\eqref{geek16} and \eqref{geek31}, the definitions \eqref{ellinfk6}, \eqref{ellinfk7}, 
\eqref{ellohk6}, and \eqref{ellohk7}, and the expansion $\me^{z} \! = \! \sum_{m=0}^{\infty} \tfrac{z^{m}}{m!}$, it follows 
{}from the definitions \eqref{geek59}--\eqref{geek64} that
\begin{gather}
y_{1} \underset{\tau \to +\infty}{=} \mathcal{O} \big(\tau^{-5/3} \me^{-2 \beta (\tau)} \big), \qquad \quad y_{2} 
\underset{\tau \to +\infty}{=} \mathcal{O} \big(\tau^{-2/3} \me^{-\beta (\tau)} \big), \label{geek65} \\
y_{3} \underset{\tau \to +\infty}{=} \mathcal{O} \big(\tau^{-1} \me^{-\beta (\tau)} \big), \quad \quad y_{4} 
\underset{\tau \to +\infty}{=} \mathcal{O} \big(\tau^{-4/3} \me^{-\beta (\tau)} \big), \quad \quad y_{5} 
\underset{\tau \to +\infty}{=} \mathcal{O} \big(\tau^{-1/3} \big), \label{geek66} \\
y_{6} \underset{\tau \to +\infty}{=} \mathcal{O} \big(\tau^{-2/3} \big), \quad \quad y_{7} \underset{\tau \to +\infty}{=} 
\mathcal{O} \big(\tau^{-1/3} \big), \quad \quad y_{8} \underset{\tau \to +\infty}{=} s_{0}^{0} \big(1 \! + \! 
\mathcal{O}(\tau^{-1/3}) \big), \label{geek67} \\
y_{9} \underset{\tau \to +\infty}{=} \mathcal{O} \big(\tau^{-1/3} \big), \quad \quad y_{10} \underset{\tau \to +\infty}{=} 
\mathcal{O} \big(\tau^{-1/3} \big), \quad \quad y_{11} \underset{\tau \to +\infty}{=} \mathcal{O} \big(\tau^{-4/3} 
\me^{-\beta (\tau)} \big). \label{geek68}
\end{gather}
One notes that the asymptotic equation \eqref{geek53} is a quartic equation for the indeterminate $x$, which can 
be solved explicitly: via a study of the four solutions of the quartic equation (see, for example, \cite{antvsnon}), 
in conjunction with the asymptotics \eqref{geek65}--\eqref{geek68}, it can be shown that the sought-after solution, 
that is, the one for which $x \genfrac{}{}{0pt}{3}{=}{\tau \to +\infty} \mathcal{O}(1)$, can be extracted as one of the 
two solutions of the quadratic equation
\begin{equation} \label{geek76} 
(1 \! + \! \upsilon_{1}^{\ast})x^{2} \! + \! \left(y_{8} \! + \! \upsilon_{2}^{\ast} \! - \! \mi \me^{-\pi a} \! + \! \mathcal{O}
(\tau^{-\delta_{1}}) \right) \! x \! + \! \upsilon_{3}^{\ast} \underset{\tau \to +\infty}{=} 0,
\end{equation}
where
\begin{gather}
\upsilon_{1}^{\ast} \! := \! y_{5} \! + \! y_{6} \underset{\tau \to +\infty}{=} \mathcal{O} \big(\tau^{-1/3} \big), \, \, \, 
\quad \quad \, \, \, \upsilon_{2}^{\ast} \! := \! y_{9} \! + \! y_{10} \! + \! y_{11} \underset{\tau \to +\infty}{=} 
\mathcal{O} \big(\tau^{-1/3} \big), \label{geek77} \\
\upsilon_{3}^{\ast} \! := \! y_{2} \! + \! y_{3} \! + \! y_{4} \underset{\tau \to +\infty}{=} \mathcal{O} \big(\tau^{-2/3} 
\me^{-\beta (\tau)} \big). \label{geek79}
\end{gather}
The roots of the quadratic equation \eqref{geek76} are
\begin{equation} \label{geek80} 
x \underset{\tau \to +\infty}{=} \dfrac{-\big(y_{8} \! + \! \upsilon_{2}^{\ast} \! - \! \mi \me^{-\pi a} \! + \! \mathcal{O}
(\tau^{-\delta_{1}}) \big) \! \pm \! \sqrtsign{\big(y_{8} \! + \! \upsilon_{2}^{\ast} \! - \! \mi \me^{-\pi a} \! + \! \mathcal{O}
(\tau^{-\delta_{1}}) \big)^{2} \! - \! 4(1 \! + \! \upsilon_{1}^{\ast}) \upsilon_{3}^{\ast}}}{2(1 \! + \! \upsilon_{1}^{\ast})};
\end{equation}
of the two solutions given by equation \eqref{geek80}, the one that is consistent with the corresponding 
$(k \! = \! +1)$ condition \eqref{ellinfk2a} reads
\begin{equation} \label{geek81} 
x \underset{\tau \to +\infty}{=} \dfrac{-\big(y_{8} \! + \! \upsilon_{2}^{\ast} \! - \! \mi \me^{-\pi a} \! + \! \mathcal{O}
(\tau^{-\delta_{1}}) \big) \! - \! \sqrtsign{\big(y_{8} \! + \! \upsilon_{2}^{\ast} \! - \! \mi \me^{-\pi a} \! + \! \mathcal{O}
(\tau^{-\delta_{1}}) \big)^{2} \! - \! 4(1 \! + \! \upsilon_{1}^{\ast}) \upsilon_{3}^{\ast}}}{2(1 \! + \! \upsilon_{1}^{\ast})}:
\end{equation}
via the definition \eqref{geek52}, and the asymptotics \eqref{geek65}, \eqref{geek77}, and \eqref{geek79}, it 
follows {}from equation \eqref{geek81} and an application of the Binomial Theorem that, for $s_{0}^{0} \! 
\neq \! \mi \me^{-\pi a}$,
\begin{align}
&\dfrac{\sqrt{2 \pi} \, p_{1}(\tau) \mathfrak{B}_{1} \me^{-2 \hat{\mathfrak{z}}_{1}^{0}
(\tau)} \me^{\mi \pi (\nu (1)+1)}}{(2 \! + \! \sqrt{3})^{1/2}(2 \mu_{1}(\tau))^{1/2} 
\Gamma (-\nu (1))} \dfrac{\hat{\mathbb{A}}_{0}^{\infty}(\tau) \hat{\mathbb{B}}_{
0}^{\infty}(\tau)}{\hat{\mathbb{A}}_{0}^{0}(\tau) \hat{\mathbb{A}}_{0}^{0}(\tau)} 
\me^{-2 \Delta \hat{\mathfrak{z}}_{1}(\tau)}(1 \! + \! \eta_{11}(\tau))(1 \! + \! 
\eta_{22}(\tau)) \nonumber \\
&\underset{\tau \to +\infty}{=} -(s_{0}^{0} \! - \! \mi \me^{-\pi a}) \! + \! \mathcal{O} 
\big(\tau^{-\delta_{1}} \big). \label{geek82}
\end{align}
{}From the asymptotics \eqref{tr1}, \eqref{tr3}, \eqref{geek13}, \eqref{geek15}, and \eqref{geek31}, the definitions 
\eqref{ellinfk6}, \eqref{ellinfk7}, \eqref{ellohk6}, and \eqref{ellohk7}, the reflection formula $\Gamma (z) \Gamma 
(1 \! - \! z) \! = \! \pi/\sin \pi z$, the expansion $\me^{z} \! = \! \sum_{m=0}^{\infty} \tfrac{z^{m}}{m!}$, and the 
asymptotics $\tfrac{1}{\Gamma (-\nu (1))} \genfrac{}{}{0pt}{3}{=}{\tau \to +\infty} 1 \! + \! \mathcal{O}(\nu (1) \! + \! 1) 
\genfrac{}{}{0pt}{3}{=}{\tau \to +\infty} 1 \! + \! \mathcal{O}(\tau^{-2/3} \me^{-\beta (\tau)})$, one shows that, for 
$k \! = \! +1$,
\begin{gather}
\dfrac{\me^{\mi \pi (\nu (1)+1)}}{\Gamma (-\nu (1))} \dfrac{\hat{\mathbb{A}}_{0}^{\infty}(\tau) 
\hat{\mathbb{B}}_{0}^{\infty}(\tau)}{\hat{\mathbb{A}}_{0}^{0}(\tau) \hat{\mathbb{A}}_{0}^{0}(\tau)} 
\underset{\tau \to +\infty}{:=} 1 \! + \! \tau^{-2/3} \sum_{m=0}^{\infty} \alpha_{m}(1)(\tau^{-1/3})^{m} 
\! + \! \mathcal{O} \big(\tau^{-1/3} \me^{-\beta (\tau)} \big), \label{geek83} \\
\me^{-2 \Delta \hat{\mathfrak{z}}_{1}(\tau)} \underset{\tau \to +\infty}{:=} 1 \! + \! \tau^{-2/3} \sum_{m=0}^{\infty} 
\alpha_{m}^{\natural}(1)(\tau^{-1/3})^{m} \! + \! \mathcal{O} \big(\tau^{-1/3} \me^{-\beta (\tau)} \big), 
\label{geek84} \\
(1 \! + \! \eta_{11}(\tau))(1 \! + \! \eta_{22}(\tau)) \underset{\tau \to +\infty}{:=} 1 \! + \! \tau^{-1/3} \sum_{m=0}^{\infty} 
\alpha_{m}^{\flat}(1)(\tau^{-1/3})^{m} \! + \! \mathcal{O} \big(\tau^{-1/3} \me^{-\beta (\tau)} \big), \label{geek85}
\end{gather}
for $\mathcal{O}(1)$ coefficients $\alpha_{m}(1)$, $\alpha_{m}^{\natural}(1)$, and $\alpha_{m}^{\flat}(1)$. Via the 
asymptotics \eqref{geek16} and \eqref{geek83}--\eqref{geek85}, upon defining
\begin{align}
&\left(1 \! + \! \sum_{m_{1}=0}^{\infty} \dfrac{\alpha_{m_{1}}^{\flat}(1)}{(\tau^{1/3})^{m_{1}+1}} \! + \! \mathcal{O} 
\big(\tau^{-1/3} \me^{-\beta (\tau)} \big) \right) \! \left(1 \! + \! \sum_{m_{2}=0}^{\infty} \dfrac{\alpha_{m_{2}}(1)}{
(\tau^{1/3})^{m_{2}+2}} \! + \! \mathcal{O} \big(\tau^{-1/3} \me^{-\beta (\tau)} \big) \right) \nonumber \\
&\times \left(1 \! + \! \sum_{m_{3}=0}^{\infty} \dfrac{\alpha_{m_{3}}^{\natural}(1)}{(\tau^{1/3})^{m_{3}+2}} 
\! + \! \mathcal{O} \big(\tau^{-1/3} \me^{-\beta (\tau)} \big) \right) \! \left(1 \! + \! \sum_{m_{4}=0}^{\infty} 
\dfrac{\alpha_{m_{4}}^{\sharp}(1)}{(\tau^{1/3})^{m_{4}+2}} \! + \! \mathcal{O} \big(\tau^{-1/3} 
\me^{-\beta (\tau)} \big) \right) \nonumber \\
&\underset{\tau \to +\infty}{=:} 1 \! + \! \sum_{m=0}^{\infty} \dfrac{\hat{\epsilon}_{m}^{\sharp}(1)}{(\tau^{1/3})^{m+1}} 
\! + \! \mathcal{O} \big(\tau^{-1/3} \me^{-\beta (\tau)} \big), \label{geek86}
\end{align}
it follows {}from the corresponding $(k \! = \! +1)$ definition \eqref{ellohk4} and the asymptotics \eqref{geek82} 
and \eqref{geek86} that, for $s_{0}^{0} \! \neq \! \mi \me^{-\pi a}$,
\begin{align} \label{geek87} 
p_{1}(\tau) \mathfrak{B}_{1} \! \left(1 \! + \! \sum_{m=0}^{\infty} \frac{\hat{\epsilon}_{m}^{\sharp}(1)}{(\tau^{1/3})^{m+1}} 
\! + \! \mathcal{O} \big(\tau^{-1/3} \me^{-\beta (\tau)} \big) \right) \underset{\tau \to +\infty}{=}& \, -\frac{2^{3/2}3^{1/4} 
\me^{\mi \pi/4}(2 \! + \! \sqrt{3}) \mathscr{P}_{a}(s_{0}^{0} \! - \! \mi \me^{-\pi a})}{\sqrt{2 \pi}} \nonumber \\
\times& \, \me^{-(\beta (\tau)+\mi \vartheta (\tau))} \big(1 \! + \! \mathcal{O}(\tau^{-\delta_{1}}) \big),
\end{align}
where $\mathscr{P}_{a}$ is defined by equation \eqref{geek4}.\footnote{{}From the leading term of asymptotics for 
$\mathfrak{B}_{1}$ given in equation \eqref{prcybk1}, that is, $\mathfrak{B}_{1} \genfrac{}{}{0pt}{3}{=}{\tau \to +\infty} 
-\tfrac{(\sqrt{3}+1)}{6 \alpha_{1} \tau^{1/3}} \! + \! \mathcal{O}(\tau^{-1})$, and the asymptotics \eqref{geek87}, it 
follows that $p_{1}(\tau) \genfrac{}{}{0pt}{3}{=}{\tau \to +\infty} \mathfrak{D}_{1} \tau^{1/3} \me^{-(\beta (\tau)+\mi 
\vartheta (\tau))} \big(1 \! + \! \mathcal{O}(\tau^{-\delta_{1}}) \big)$, where $\mathfrak{D}_{1} \! := \! 6(\sqrt{3} \! + \! 
1)3^{1/4} \me^{\mi \pi/4} \alpha_{1} \mathscr{P}_{a}(s_{0}^{0} \! - \! \mi \me^{-\pi a})/\sqrt{\pi}$, whence $p_{1}(\tau) 
\mathfrak{B}_{1} \genfrac{}{}{0pt}{3}{=}{\tau \to +\infty} \mathcal{O}(\me^{-\beta (\tau)})$, which is consistent with the 
corresponding $(k \! = \! +1)$ condition~\eqref{ellinfk2a}.} Via the asymptotics \eqref{prcyellok1} and the definition 
\eqref{prcy57}, a multiplication argument shows that
\begin{align} \label{geek88} 
p_{1}(\tau) \mathfrak{B}_{1} \underset{\tau \to +\infty}{=}& \, -\dfrac{\mi \mathfrak{B}_{0,1}^{\sharp}}{8 \sqrt{3}} \! 
+ \! \mathfrak{B}_{1}(1 \! + \! \hat{\mathbb{L}}_{1}(\tau)) \! - \! \dfrac{\mi \tilde{r}_{0,1}(\tau) \tau^{-1/3}}{96 \sqrt{3}} 
\! \left(1 \! + \!  \mathcal{O} \big((\tilde{r}_{0,1}(\tau) \tau^{-1/3})^{2} \big) \right) \! \mathfrak{B}_{0,1}^{\sharp} 
\nonumber \\
+& \, \dfrac{\mi \omega_{0,1}^{2}}{(8 \sqrt{3})^{3}} \! \left(1 \! + \! \dfrac{\tilde{r}_{0,1}(\tau) \tau^{-1/3}}{12} \! + \! 
\mathcal{O} \big((\tilde{r}_{0,1}(\tau) \tau^{-1/3})^{3} \big) \right)^{3} \! \left(\dfrac{\mathfrak{B}_{0,1}^{\sharp}}{
\mathfrak{B}_{1}} \right)^{2} \! \mathfrak{B}_{1} \nonumber \\
+& \, \mathcal{O} \! \left(\omega_{0,1}^{4} \! \left(1 \! + \! \dfrac{\tilde{r}_{0,1}(\tau) \tau^{-1/3}}{12} \! + \! \mathcal{O} 
\big((\tilde{r}_{0,1}(\tau) \tau^{-1/3})^{3} \big) \right)^{5} \! \left(\dfrac{\mathfrak{B}_{0,1}^{\sharp}}{\mathfrak{B}_{1}} 
\right)^{3} \! \mathfrak{B}_{1} \right);
\end{align}
{}from the corresponding $(k \! = \! +1)$ asymptotics \eqref{tr1}, \eqref{tr3}, \eqref{prcybk1}, \eqref{prcyb0k1}, 
and \eqref{prcyomg1}, the various terms appearing in the asymptotics \eqref{geek88} can be presented as 
follows:\footnote{Note, in particular, that $\mathfrak{B}_{0,1}^{\sharp}/\mathfrak{B}_{1} 
\genfrac{}{}{0pt}{3}{=}{\tau \to +\infty} -\mi 8 \sqrt{3}(1 \! + \! o(1))$.}
{\fontsize{10pt}{11pt}\selectfont
\begin{gather}
-\dfrac{\mi \mathfrak{B}_{0,1}^{\sharp}}{8 \sqrt{3}} \underset{\tau \to +\infty}{=} \dfrac{(\sqrt{3} \! + \! 1) 
\tau^{-1/3}}{6 \alpha_{1}} \! + \! \sum_{m=0}^{\infty} \dfrac{\mathfrak{b}_{m}^{\flat}(1)}{(\tau^{1/3})^{m+3}} 
\! + \! \mathcal{O} \big(\tau^{-2/3} \me^{-\beta (\tau)} \big), \label{geek89} \\
-\dfrac{\mi \tilde{r}_{0,1}(\tau) \tau^{-1/3}}{96 \sqrt{3}} \! \left(1 \! + \!  \mathcal{O} \big((\tilde{r}_{0,1}(\tau) 
\tau^{-1/3})^{2} \big) \right) \! \mathfrak{B}_{0,1}^{\sharp} \underset{\tau \to +\infty}{=} \sum_{m=0}^{\infty} 
\dfrac{\mathfrak{b}_{m}^{\natural}(1)}{(\tau^{1/3})^{m+3}} \! + \! \mathcal{O} \big(\tau^{-2/3} \me^{-\beta (\tau)} \big), 
\label{geek90}
\end{gather}
\begin{align}
&\dfrac{\mi \omega_{0,1}^{2}}{(8 \sqrt{3})^{3}} \! \left(1 \! + \! \dfrac{\tilde{r}_{0,1}(\tau) \tau^{-1/3}}{12} \! + \! 
\mathcal{O} \big((\tilde{r}_{0,1}(\tau) \tau^{-1/3})^{3} \big) \right)^{3} \! \left(\dfrac{\mathfrak{B}_{0,1}^{\sharp}}{
\mathfrak{B}_{1}} \right)^{2} \! \mathfrak{B}_{1} \! + \! \mathcal{O} \! \left(\omega_{0,1}^{4} \! \left(1 \! + \! 
\dfrac{\tilde{r}_{0,1}(\tau) \tau^{-1/3}}{12} \right. \right. \nonumber \\
+&\left. \left. \mathcal{O} \big((\tilde{r}_{0,1}(\tau) \tau^{-1/3})^{3} \big) \right)^{5} \! \left(
\dfrac{\mathfrak{B}_{0,1}^{\sharp}}{\mathfrak{B}_{1}} \right)^{3} \! \mathfrak{B}_{1} \right) \underset{\tau \to +\infty}{=} 
\sum_{m=0}^{\infty} \dfrac{\mathfrak{b}_{m}^{\sharp}(1)}{(\tau^{1/3})^{m+3}} \! + \! \mathcal{O} \big(\tau^{-2/3} 
\me^{-\beta (\tau)} \big), \label{geek91} 
\end{align}}
for $\mathcal{O}(1)$ coefficients $\mathfrak{b}_{m}^{\flat}(1)$, $\mathfrak{b}_{m}^{\natural}(1)$, and 
$\mathfrak{b}_{m}^{\sharp}(1)$, whence (cf. asymptotics \eqref{geek88})
\begin{equation} \label{geek92} 
p_{1}(\tau) \mathfrak{B}_{1} \underset{\tau \to +\infty}{=} \mathfrak{B}_{1}(1 \! + \! \hat{\mathbb{L}}_{1}(\tau)) 
\! + \! \dfrac{(\sqrt{3} \! + \! 1) \tau^{-1/3}}{6 \alpha_{1}} \! + \! \sum_{m=0}^{\infty} 
\dfrac{\mathfrak{b}_{m}^{\dagger}(1)}{(\tau^{1/3})^{m+3}} \! + \! \mathcal{O} \big(\tau^{-2/3} 
\me^{-\beta (\tau)} \big),
\end{equation}
for $\mathcal{O}(1)$ coefficients $\mathfrak{b}_{m}^{\dagger}(1)$; for example,
\begin{equation} \label{nueeq1} 
\mathfrak{b}_{0}^{\dagger}(1) \! = \! \dfrac{\mi (\sqrt{3} \! + \! 1)}{48 \sqrt{3} \alpha_{1}} \! \left(6 \mi \mathfrak{r}_{0}(1) 
\! + \! 4(a \! - \! \mi/2) \mathfrak{u}_{0}(1) \! - \! \alpha_{1}^{2}(8 \mathfrak{u}_{0}^{2}(1) \! + \! 4 \mathfrak{u}_{0}(1) 
\mathfrak{r}_{0}(1) \! - \! \mathfrak{r}_{0}^{2}(1)) \right).
\end{equation}
One shows {}from the corresponding $(k \! = \! +1)$ asymptotics \eqref{tr1}, \eqref{tr3}, and \eqref{prcybk1} that
\begin{equation} \label{nueeq2} 
\mathfrak{B}_{1} \underset{\tau \to +\infty}{=} \llfloor \mathfrak{B}_{1} \rrfloor \! + \! \dfrac{\mi (\sqrt{3} \! + \! 1) 
\alpha_{1}}{2}(4 \mathrm{A}_{1} \! + \! (\sqrt{3} \! + \! 1) \mathrm{B}_{1}) \me^{-(\beta (\tau)+\mi \vartheta (\tau))} 
\big(1 \! + \! \mathcal{O}(\tau^{-1/3}) \big),
\end{equation}
where
\begin{equation} \label{nueeq3} 
\mathrm{B}_{1} \! := \! 2(1 \! + \! \sqrt{3}) \mathrm{A}_{1},
\end{equation}
and
\begin{equation} \label{nueeq4} 
\llfloor \mathfrak{B}_{1} \rrfloor \! := \! -\dfrac{(\sqrt{3} \! + \! 1) \tau^{-1/3}}{6 \alpha_{1}} \! + \! \sum_{m=0}^{\infty} 
\dfrac{b_{m}(1)}{(\tau^{1/3})^{m+3}},
\end{equation}
for $\mathcal{O}(1)$ coefficients $b_{m}(1)$; for example,
\begin{align}
b_{0}(1) \! =& \, \dfrac{\mi (\sqrt{3} \! + \! 1)^{2}}{2} \! \left(\alpha_{1} \mathfrak{r}_{2}(1) \! + \! \dfrac{1}{2 \sqrt{3}} 
\! \left(-\dfrac{\alpha_{1}}{2}(\mathfrak{r}_{0}^{2}(1) \! + \! 2(\sqrt{3} \! + \! 1) \mathfrak{r}_{0}(1) \mathfrak{u}_{0}
(1) \! + \! 8 \mathfrak{u}_{0}^{2}(1)) \right. \right. \nonumber \\
+&\left. \left. \frac{(a \! - \! \mi/2)}{6 \alpha_{1}}(12 \mathfrak{u}_{0}(1) \! + \! (2 \sqrt{3} \! - \! 1) \mathfrak{r}_{0}
(1)) \right) \right), \label{nueeq5} \\
b_{1}(1) \! =& \, 0. \label{nueeq6}
\end{align}
{}From the expansions \eqref{geek92} and \eqref{nueeq2}, and the definition \eqref{nueeq4}, it follows that
\begin{align} \label{nueeq7} 
p_{1}(\tau) \mathfrak{B}_{1} \underset{\tau \to +\infty}{=}& \, \tau^{-1} \sum_{m=0}^{\infty} 
\dfrac{d_{m}^{\ast}(1)}{(\tau^{1/3})^{m}} \! + \! \hat{\mathbb{L}}_{1}(\tau) \! \left(\llfloor \mathfrak{B}_{1} 
\rrfloor \! + \! \mathcal{O}(\me^{-\beta (\tau)}) \right) \nonumber \\
+& \, \dfrac{\mi (\sqrt{3} \! + \! 1) \alpha_{1}}{2}(4 \mathrm{A}_{1} \! + \! (\sqrt{3} \! + \! 1) \mathrm{B}_{1}) 
\me^{-(\beta (\tau)+\mi \vartheta (\tau))} \big(1 \! + \! \mathcal{O}(\tau^{-1/3}) \big),
\end{align}
for $\mathcal{O}(1)$ coefficients $d_{m}^{\ast}(1) \! := \! \mathfrak{b}_{m}^{\dagger}(1) \! + \! b_{m}(1)$, 
$m \! \in \! \mathbb{Z}_{+}$; for example,
\begin{align} \label{nueeq8} 
d_{0}^{\ast}(1) =& \, \dfrac{\mi (\sqrt{3} \! + \! 1)}{48 \sqrt{3} \alpha_{1}} \! \left(6 \mi \mathfrak{r}_{0}(1) \! + \! 
4(a \! - \! \mi/2) \mathfrak{u}_{0}(1) \! - \! \alpha_{1}^{2}(8 \mathfrak{u}_{0}^{2}(1) \! + \! 4 \mathfrak{u}_{0}(1) 
\mathfrak{r}_{0}(1) \! - \! \mathfrak{r}_{0}^{2}(1)) \right) \nonumber \\
+& \, \dfrac{\mi (\sqrt{3} \! + \! 1)^{2}}{2} \! \left(\alpha_{1} \mathfrak{r}_{2}(1) \! + \! \dfrac{1}{2 \sqrt{3}} \! 
\left(-\dfrac{\alpha_{1}}{2}(\mathfrak{r}_{0}^{2}(1) \! + \! 2(\sqrt{3} \! + \! 1) \mathfrak{r}_{0}(1) \mathfrak{u}_{0}
(1) \! + \! 8 \mathfrak{u}_{0}^{2}(1)) \right. \right. \nonumber \\
+&\left. \left. \dfrac{(a \! - \! \mi/2)}{6 \alpha_{1}}(12 \mathfrak{u}_{0}(1) \! + \! (2 \sqrt{3} \! - \! 1) 
\mathfrak{r}_{0}(1)) \right) \right).
\end{align}
Thus, via the asymptotics \eqref{geek87} and \eqref{nueeq7}, one arrives at
\begin{align} \label{nueeq9} 
&\left(\sum_{m=0}^{\infty} \dfrac{d_{m}^{\ast}(1)}{(\tau^{1/3})^{m+3}} \! + \! \hat{\mathbb{L}}_{1}(\tau) \! \left(
\llfloor \mathfrak{B}_{1} \rrfloor \! + \! \mathcal{O}(\me^{-\beta (\tau)}) \right) \! + \! \frac{\mi (\sqrt{3} \! + \! 1) 
\alpha_{1}}{2}(4 \mathrm{A}_{1} \! + \! (\sqrt{3} \! + \! 1) \mathrm{B}_{1}) \me^{-(\beta (\tau)+\mi \vartheta 
(\tau))} \right. \nonumber \\
&\times \left. \big(1 \! + \! \mathcal{O}(\tau^{-1/3}) \big) \right) \! \left(1 \! + \! \sum_{m=0}^{\infty} 
\dfrac{\hat{\epsilon}_{m}^{\sharp}(1)}{(\tau^{1/3})^{m+1}} \! + \! \mathcal{O} \big(\tau^{-1/3} \me^{-\beta (\tau)} 
\big) \right) \underset{\tau \to +\infty}{=} -\mathcal{Q}_{1} \me^{-(\beta (\tau)+\mi \vartheta (\tau))} \big(1 \! + 
\! \mathcal{O}(\tau^{-\delta_{1}})),
\end{align}
where
\begin{equation} \label{geek94} 
\mathcal{Q}_{1} \! := \! \dfrac{2^{3/2}3^{1/4} \me^{\mi \pi/4}(2 \! + \! \sqrt{3}) \mathscr{P}_{a}(s_{0}^{0} \! - \! 
\mi \me^{-\pi a})}{\sqrt{2 \pi}}.
\end{equation}
One now chooses $\hat{\mathbb{L}}_{1}(\tau)$ so that the---divergent---power series on the left-hand side of equation 
\eqref{nueeq9} is identically equal to zero:
\begin{equation} \label{nueeq10} 
\left(\tau^{-1} \sum_{m=0}^{\infty} \dfrac{d_{m}^{\ast}(1)}{(\tau^{1/3})^{m}} \! + \! \hat{\mathbb{L}}_{1}(\tau) 
\llfloor \mathfrak{B}_{1} \rrfloor \right) \! \left(1 \! + \! \tau^{-1/3} \sum_{m=0}^{\infty} \dfrac{\hat{\epsilon}_{m}^{\sharp}
(1)}{(\tau^{1/3})^{m}} \right) \! = \! 0;
\end{equation}
via the definition \eqref{nueeq4}, one solves equation \eqref{nueeq10} for $\hat{\mathbb{L}}_{1}(\tau)$ to arrive at
\begin{equation} \label{nueeq11} 
\hat{\mathbb{L}}_{1}(\tau) \! = \! \tau^{-2/3} \sum_{m=0}^{\infty} \dfrac{\hat{\mathfrak{l}}_{m+2}(1)}{(\tau^{1/3})^{m}},
\end{equation}
where the coefficients $\hat{\mathfrak{l}}_{m^{\prime}}(1)$, $m^{\prime} \! \in \! \mathbb{Z}_{+}$, are determined 
according to the recursion relation 
\begin{gather}
\hat{\mathfrak{l}}_{0}(1) \! = \! \hat{\mathfrak{l}}_{1}(1) \! = \! 0, \, \, \quad \quad \, \, \hat{\mathfrak{l}}_{2}(1) 
\! = \! \dfrac{6 \alpha_{1}d_{0}^{\ast}(1)}{\sqrt{3} \! + \! 1}, \label{nueeq12} \\
\hat{\mathfrak{l}}_{m+3}(1) \! = \! \dfrac{6 \alpha_{1}}{\sqrt{3} \! + \! 1} \! \left(d_{m+1}^{\ast}(1) \! + \! 
\sum_{p=0}^{m}d_{p}^{\ast}(1) \hat{\epsilon}_{m-p}^{\sharp}(1) \! + \! \sum_{j=0}^{m+2} \hat{\mathfrak{l}}_{j}(1) 
\hat{d}_{m+4-j}(1) \right), \quad m \! \in \! \mathbb{Z}_{+}, \label{nueeq13}
\end{gather}
with
\begin{gather}
\hat{d}_{0}(1) \! = \! 0, \! \quad \! \hat{d}_{1}(1) \! = \! -\dfrac{(\sqrt{3} \! + \! 1)}{6 \alpha_{1}}, \! \quad \! \hat{d}_{2}
(1) \! = \! -\dfrac{(\sqrt{3} \! + \! 1) \hat{\epsilon}_{0}^{\sharp}(1)}{6 \alpha_{1}}, \! \quad \! \hat{d}_{3}(1) \! = \! b_{0}
(1) \! - \! \dfrac{(\sqrt{3} \! + \! 1) \hat{\epsilon}_{1}^{\sharp}(1)}{6 \alpha_{1}}, \label{nueeq14} \\
\hat{d}_{m+4}(1) \! = \! b_{m+1}(1) \! - \! \dfrac{(\sqrt{3} \! + \! 1) \hat{\epsilon}_{m+2}^{\sharp}(1)}{6 \alpha_{1}} 
\! + \! \sum_{p=0}^{m}b_{p}(1) \hat{\epsilon}_{m-p}^{\sharp}(1), \quad m \! \in \! \mathbb{Z}_{+}. \label{nueeq15}
\end{gather}
{}From the condition \eqref{nueeq10}, equation \eqref{nueeq11}, and the asymptotics \eqref{nueeq9}, it follows that
\begin{equation} \label{nueeq16} 
\dfrac{\mi (\sqrt{3} \! + \! 1) \alpha_{1}}{2}(4 \mathrm{A}_{1} \! + \! (\sqrt{3} \! + \! 1) \mathrm{B}_{1}) 
\me^{-(\beta (\tau)+\mi \vartheta (\tau))} \underset{\tau \to +\infty}{=} -\mathcal{Q}_{1} 
\me^{-(\beta (\tau)+\mi \vartheta (\tau))} \big(1 \! + \! \mathcal{O}(\tau^{-\delta_{1}}) \big),
\end{equation}
whence, via the definitions \eqref{geek4}, \eqref{nueeq3}, and \eqref{geek94}, one arrives at
\begin{equation} \label{geek109} 
\mathrm{A}_{1} \! = \! \dfrac{\mi \me^{\mi \pi/4} \me^{-\mi \pi/3}(2 \! + \! \sqrt{3})^{\mi a}(s_{0}^{0} \! - \! 
\mi \me^{-\pi a})}{\sqrt{2 \pi} \, 3^{1/4}(\varepsilon b)^{1/6}}.
\end{equation}

Alternatively, one may proceed as follows. Substituting the asymptotics \eqref{geek92} and \eqref{nueeq2} into 
equation \eqref{geek87}, one shows, via the definition \eqref{nueeq4} and the definition $d_{m}^{\ast}(1) \! := \! 
\mathfrak{b}_{m}^{\dagger}(1) \! + \! b_{m}(1)$, $m \! \in \! \mathbb{Z}_{+}$, that
\begin{align}
&\mathfrak{B}_{1} \! + \! \dfrac{(\sqrt{3} \! + \! 1) \tau^{-1/3}}{6 \alpha_{1}} \! + \! \tau^{-1} \sum_{m=0}^{\infty} 
\dfrac{d_{m}(1)}{(\tau^{1/3})^{m}} \! + \! \hat{\mathbb{L}}_{1}(\tau) \mathfrak{B}_{1} \! \left(1 \! + \! \tau^{-1/3} 
\sum_{m=0}^{\infty} \dfrac{\hat{\epsilon}_{m}^{\sharp}(1)}{(\tau^{1/3})^{m}} \right. \nonumber \\
+&\left. \mathcal{O} \big(\tau^{-1/3} \me^{-\beta (\tau)} \big) \right) \! + \! \mathcal{O} \big(\tau^{-1/3} 
\me^{-\beta (\tau)} \big) \underset{\tau \to +\infty}{=} -\mathcal{Q}_{1} \me^{-(\beta (\tau)+\mi \vartheta (\tau))} 
\big(1 \! + \! \mathcal{O}(\tau^{-\delta_{1}}) \big), \label{nueeq17}
\end{align}
where $\mathcal{Q}_{1}$ is defined by equation \eqref{geek94},
\begin{equation} \label{nueeq18}
d_{0}(1) \! = \! \mathfrak{b}_{0}^{\dagger}(1), \quad \quad d_{m+1}(1) \! = \! \mathfrak{b}_{m+1}^{\dagger}(1) 
\! + \! \sum_{p=0}^{m}d_{p}^{\ast}(1) \hat{\epsilon}_{m-p}^{\sharp}(1), \quad m \! \in \! \mathbb{Z}_{+}.
\end{equation}
{}From the condition \eqref{nueeq10}, equation \eqref{nueeq11}, the asymptotics \eqref{nueeq17}, the 
definition $d_{m}^{\ast}(1) \! := \! \mathfrak{b}_{m}^{\dagger}(1) \! + \! b_{m}(1)$, $m \! \in \! \mathbb{Z}_{+}$, 
and equations \eqref{nueeq18}, it follows that
\begin{equation} \label{nueeq19} 
\mathfrak{B}_{1} \underset{\tau \to +\infty}{=} -\dfrac{(\sqrt{3} \! + \! 1) \tau^{-1/3}}{6 \alpha_{1}} 
\! + \! \tau^{-1} \sum_{m=0}^{\infty} \dfrac{b_{m}(1)}{(\tau^{1/3})^{m}} \! - \! \mathcal{Q}_{1} 
\me^{-(\beta (\tau)+\mi \vartheta (\tau))} \big(1 \! + \! \mathcal{O}(\tau^{-\delta_{1}}) \big).
\end{equation}
It follows {}from the corresponding $(k \! = \! +1)$ asymptotics \eqref{tr1}, \eqref{tr3}, and \eqref{prcybk1} that 
the function $\mathfrak{B}_{1}$ can also be presented in the form
\begin{align} \label{nueeq20} 
\mathfrak{B}_{1} \underset{\tau \to +\infty}{=}& \, \mi (\sqrt{3} \! + \! 1) \! \left(\dfrac{\alpha_{1}}{2}(4v_{0,1}
(\tau) \! + \! (\sqrt{3} \! + \! 1) \tilde{r}_{0,1}(\tau)) \! - \! \dfrac{(\sqrt{3} \! + \! 1)(a \! - \! \mi/2)}{2 \sqrt{3} \alpha_{1} 
\tau^{1/3}} \right) \! + \! \sum_{m=0}^{\infty} \dfrac{\hat{b}_{m}^{\ast}(1)}{(\tau^{1/3})^{m+3}} \nonumber \\
+& \, \mathcal{O} \big(\tau^{-2/3} \me^{-\beta (\tau)} \big),
\end{align}
for $\mathcal{O}(1)$ coefficients $\hat{b}_{m}^{\ast}(1)$ (see, for example, equations \eqref{nueeq22} and 
\eqref{nueeq23}); hence, {}from the asymptotics \eqref{nueeq19} and \eqref{nueeq20}, one deduces that
\begin{align} \label{geek97}  
4v_{0,1}(\tau) \! + \! (\sqrt{3} \! + \! 1) \tilde{r}_{0,1}(\tau) \underset{\tau \to +\infty}{=}& \, \dfrac{(\sqrt{3} \! + \! 1)
(\sqrt{3}a \! - \! \mi/2)}{3 \alpha_{1}^{2} \tau^{1/3}} \! + \! \sum_{m=0}^{\infty} \dfrac{\iota_{m}^{\ast}(1)}{(\tau^{1/3})^{m+3}} 
\nonumber \\
+& \, \dfrac{2 \mi \mathcal{Q}_{1} \me^{-(\beta (\tau)+\mi \vartheta (\tau))}}{(\sqrt{3} \! + \! 1) \alpha_{1}} 
\big(1 \! + \! \mathcal{O}(\tau^{-\delta_{1}}) \big),
\end{align}
where
\begin{equation} \label{nueeq21} 
\iota_{m}^{\ast}(1) \! := \! -\dfrac{\mi 2(b_{m}(1) \! - \! \hat{b}_{m}^{\ast}(1))}{(\sqrt{3} \! + \! 1) \alpha_{1}}, 
\quad m \! \in \! \mathbb{Z}_{+}.
\end{equation}
Combining the corresponding $(k \! = \! +1)$ equations \eqref{iden217} and \eqref{tr2}, it follows that, in terms 
of the corresponding $(k \! = \! +1)$ solution of the DP3E \eqref{eq1.1},
\begin{equation} \label{geek98}  
4v_{0,1}(\tau) \! + \! (\sqrt{3} \! + \! 1) \tilde{r}_{0,1}(\tau) \! = \! \dfrac{8 \me^{2 \pi \mi/3}u(\tau)}{\varepsilon 
(\varepsilon b)^{2/3}} \! - \! \dfrac{\mi (\sqrt{3} \! + \! 1) \me^{-\mi 2 \pi/3} \tau^{2/3}}{(\varepsilon b)^{1/3}} \! 
\left(\dfrac{u^{\prime}(\tau) \! - \! \mi b}{u(\tau)} \right) \! + \! 2(\sqrt{3} \! - \! 1) \tau^{1/3};
\end{equation}
finally, {}from the asymptotics \eqref{geek97} and equation \eqref{geek98}, one arrives at 
the---asymptotic---Riccati differential equation 
\begin{equation} \label{geek99} 
u^{\prime}(\tau) \underset{\tau \to +\infty}{=} \tilde{\mathfrak{a}}(\tau) \! + \! \tilde{\mathfrak{b}}(\tau)u(\tau) 
\! + \! \tilde{\mathfrak{c}}(\tau)(u(\tau))^{2},
\end{equation}
where
\begin{equation} \label{ricc1} 
\begin{gathered} 
\tilde{\mathfrak{a}}(\tau) \! := \! \mi b, \qquad \qquad \tilde{\mathfrak{c}}(\tau) \! := \! \frac{\mi 8 \sqrt{2} \varepsilon 
\alpha_{1} \tau^{-2/3}}{(\sqrt{3} \! + \! 1)(\varepsilon b)^{1/2}}, \\
\tilde{\mathfrak{b}}(\tau) \! := \! -\frac{8 \mi \alpha_{1}^{2} \tau^{-1/3}}{(\sqrt{3} \! + \! 1)^{2}} \! + \! 
\frac{2 \mi (\sqrt{3}a \! - \! \mi/2)}{3 \tau} \! + \! \frac{2 \mi \alpha_{1}^{2}}{(\sqrt{3} \! + \! 1)} \sum_{m=0}^{\infty} 
\frac{\iota_{m}^{\ast}(1)}{(\tau^{1/3})^{m+5}} \! - \! \frac{4 \alpha_{1} \mathcal{Q}_{1} 
\me^{-(\beta (\tau)+\mi \vartheta (\tau))}}{(\sqrt{3} \! + \! 1)^{2} \tau^{2/3}} \big(1 \! + \! \mathcal{O}
(\tau^{-\delta_{1}}) \big).
\end{gathered}
\end{equation}
Incidentally, changing the dependent variable according to $w(\tau) \! = \! \tfrac{1}{2} \tilde{\mathfrak{b}}(\tau) 
\! + \! \tfrac{1}{2} \tfrac{\tilde{\mathfrak{c}}^{\prime}(\tau)}{\tilde{\mathfrak{c}}(\tau)} \! + \! \tilde{\mathfrak{c}}
(\tau)u(\tau)$,\footnote{See Section 4.6 of \cite{ntsnmz}; see, also, Chapter 5 of \cite{eillh}.} it follows that the 
Riccati differential equation \eqref{geek99} transforms into
\begin{equation} \label{ricc2} 
w^{\prime}(\tau) \underset{\tau \to +\infty}{=} \Xi (\tau) \! + \! (w(\tau))^{2},
\end{equation}
where
\begin{equation} \label{ricc3} 
-\Xi (\tau) \! := \! -\tilde{\mathfrak{a}}(\tau) \tilde{\mathfrak{c}}(\tau) \! + \! \frac{1}{4}
(\tilde{\mathfrak{b}}(\tau))^{2} \! - \! \frac{1}{2} \tilde{\mathfrak{b}}^{\prime}(\tau) \! 
+ \! \frac{1}{2} \frac{\tilde{\mathfrak{b}}(\tau) \tilde{\mathfrak{c}}^{\prime}(\tau)}{
\tilde{\mathfrak{c}}(\tau)} \! - \! \frac{1}{2} \frac{\tilde{\mathfrak{c}}^{\prime \prime}
(\tau)}{\tilde{\mathfrak{c}}(\tau)} \! + \! \frac{3}{4} \! \left(\frac{\tilde{\mathfrak{c}}^{
\prime}(\tau)}{\tilde{\mathfrak{c}}(\tau)} \right)^{2}.
\end{equation}
Substituting the corresponding $(k \! = \! +1)$ differentiable asymptotics \eqref{recur15} into either the 
Riccati differential equation \eqref{geek99} or its dependent-variable-transformed variant~\eqref{ricc2}, and 
recalling that $c_{0,1} \! = \! \tfrac{1}{2} \varepsilon (\varepsilon b)^{2/3} \me^{-\mi 2 \pi/3}$, one shows that
\begin{align} \label{geek100} 
&\dfrac{8 \varepsilon \me^{\mi 2 \pi/3}}{(\varepsilon b)^{2/3}} \! \left(c_{0,1}^{2} \tau^{2/3} \! + \! 2c_{0,1}^{2} 
\sum_{m=0}^{\infty} \dfrac{\mathfrak{u}_{m}(1)}{(\tau^{1/3})^{m}} \! + \! c_{0,1}^{2} \tau^{-2/3} \sum_{m=0}^{\infty} 
\sum_{m_{1}=0}^{m} \mathfrak{u}_{m_{1}}(1) \mathfrak{u}_{m-m_{1}}(1)(\tau^{-1/3})^{m} \right. \nonumber \\
&+ \left. 2c_{0,1} \mathbb{P} \tau^{1/3} \me^{-(\beta (\tau)+\mi \vartheta (\tau))} \big(1 \! + \! \mathcal{O}
(\tau^{-1/3}) \big) \right) \! - \! \dfrac{\mi (\sqrt{3} \! + \! 1) \me^{-\mi 2 \pi/3} \tau^{2/3}}{(\varepsilon b)^{1/3}} 
\! \left(-\mi b \! + \! \dfrac{c_{0,1}}{3} \tau^{-2/3} \right. \nonumber \\
&- \left. \dfrac{c_{0,1}}{3} \sum_{m=0}^{\infty} \dfrac{(m \! + \! 1) \mathfrak{u}_{m}(1)}{(\tau^{1/3})^{m+4}} \! + 
\! \mi 2 \sqrt{3}(\varepsilon b)^{1/3} \me^{\mi 2 \pi/3} \mathbb{P} \tau^{-1/3} \me^{-(\beta (\tau)+\mi \vartheta 
(\tau))} \big(1 \! + \! \mathcal{O}(\tau^{-1/3}) \big) \right) \nonumber \\
&+ 2(\sqrt{3} \! - \! 1) \tau^{1/3} \! \left(c_{0,1} \tau^{1/3} \! + \! c_{0,1} \sum_{m=0}^{\infty} \dfrac{\mathfrak{u}_{m}
(1)}{(\tau^{1/3})^{m+1}} \! + \! \mathbb{P} \me^{-(\beta (\tau)+\mi \vartheta (\tau))} \big(1 \! + \! \mathcal{O}
(\tau^{-1/3}) \big) \right) \nonumber \\
& \underset{\tau \to +\infty}{=} \left(\dfrac{(\sqrt{3} \! + \! 1)(\sqrt{3}a \! - \! \mi/2)}{3 \alpha_{1}^{2} \tau^{1/3}} 
\! + \! \sum_{m=0}^{\infty} \frac{\iota_{m}^{\ast}(1)}{(\tau^{1/3})^{m+3}} \! + \! \dfrac{2 \mi \mathcal{Q}_{1} 
\me^{-(\beta (\tau)+\mi \vartheta (\tau))}}{(\sqrt{3} \! + \! 1) \alpha_{1}} \big(1 \! + \! \mathcal{O}(\tau^{-\delta_{1}}) 
\big) \right) \nonumber \\
& \, \times \left(c_{0,1} \tau^{1/3} \! + \! c_{0,1} \sum_{m=0}^{\infty} \dfrac{\mathfrak{u}_{m}(1)}{(\tau^{1/3})^{m+1}} 
\! + \! \mathbb{P} \me^{-(\beta (\tau)+\mi \vartheta (\tau))} \big(1 \! + \! \mathcal{O}(\tau^{-1/3}) \big) \right),
\end{align}
where
\begin{equation} \label{geek101} 
\mathbb{P} \! := \! c_{0,1} \mathrm{A}_{1}.
\end{equation}
Equating the coefficients of terms of order $\mathcal{O}(\tau^{1/3} \me^{-(\beta (\tau)+\mi \vartheta (\tau))})$, 
$\mathcal{O}(\tau^{2/3})$, $\mathcal{O}(1)$, $\mathcal{O}(\tau^{-1/3})$, $\mathcal{O}(\tau^{-2/3})$, and 
$\mathcal{O}(\tau^{-1})$, respectively, in equation \eqref{geek100}, one arrives at, in the indicated order:
\begin{align}
&\left(\dfrac{16 \me^{\mi 2 \pi/3}c_{0,1}}{\varepsilon (\varepsilon b)^{2/3}} \! + \! 2 \sqrt{3}(\sqrt{3} \! + \! 1) 
\! + \! 2(\sqrt{3} \! - \! 1) \right) \! \mathbb{P} \! = \! \dfrac{2 \mi \mathcal{Q}_{1}c_{0,1}}{(\sqrt{3} \! + \! 1) 
\alpha_{1}}, \label{geek102} \\
&\dfrac{8 \me^{\mi 2 \pi/3}c_{0,1}^{2}}{\varepsilon (\varepsilon b)^{2/3}} \! - \! \dfrac{(\sqrt{3} \! + \! 1)b 
\me^{-\mi 2 \pi/3}}{(\varepsilon b)^{1/3}} \! + \! 2(\sqrt{3} \! - \! 1)c_{0,1} \! = \! 0, \label{geek103} \\
&\dfrac{16 \me^{\mi 2 \pi/3}c_{0,1} \mathfrak{u}_{0}(1)}{\varepsilon (\varepsilon b)^{2/3}} \! - \! 
\dfrac{\mi (\sqrt{3} \! + \! 1) \me^{-\mi 2 \pi/3}}{3(\varepsilon b)^{1/3}} \! + \! 2(\sqrt{3} \! - \! 1) \mathfrak{u}_{0}(1) 
\! = \! \dfrac{(\sqrt{3} \! + \! 1)(\sqrt{3}a \! - \! \mi/2)}{3 \alpha_{1}^{2}}, \label{geek104} \\
&\left(\dfrac{16 \me^{\mi 2 \pi/3}c_{0,1}}{\varepsilon (\varepsilon b)^{2/3}} \! + \! 2(\sqrt{3} \! - \! 1) \right) \! 
\mathfrak{u}_{1}(1) \! = \! 0, \label{geek105} \\
&\dfrac{8 \me^{\mi 2 \pi/3}c_{0,1}}{\varepsilon (\varepsilon b)^{2/3}} \! \left(2 \mathfrak{u}_{2}(1) \! + \! 
\mathfrak{u}_{0}^{2}(1) \right) \! + \! \dfrac{\mi (\sqrt{3} \! + \! 1) \me^{-\mi 2 \pi/3} \mathfrak{u}_{0}(1)}{3
(\varepsilon b)^{1/3}} \! + \! 2(\sqrt{3} \! - \! 1) \mathfrak{u}_{2}(1) \nonumber \\
&= \! \dfrac{(\sqrt{3} \! + \! 1)(\sqrt{3}a \! - \! \mi/2)}{3 \alpha_{1}^{2}} \! + \! \iota_{0}^{\ast}(1), \label{geek106} \\
&\dfrac{16 \me^{\mi 2 \pi/3}c_{0,1}}{\varepsilon (\varepsilon b)^{2/3}} \! \left(\mathfrak{u}_{3}(1) \! + \! 
\mathfrak{u}_{0}(1) \mathfrak{u}_{1}(1) \right) \! + \! \dfrac{2 \mi (\sqrt{3} \! + \! 1) \me^{-\mi 2 \pi/3} 
\mathfrak{u}_{1}(1)}{3(\varepsilon b)^{1/3}} \! + \! 2(\sqrt{3} \! - \! 1) \mathfrak{u}_{3}(1) \nonumber \\
&= \! \dfrac{(\sqrt{3} \! + \! 1)(\sqrt{3}a \! - \! \mi/2) \mathfrak{u}_{1}(1)}{3 \alpha_{1}^{2}} \! + \! 
\iota_{1}^{\ast}(1). \label{geek107}
\end{align}
Using the corresponding $(k \! = \! +1)$ coefficients \eqref{thmk3}, in particular, $\mathfrak{u}_{0}(1) \! = \! a/6 
\alpha_{1}^{2}$ and $\mathfrak{u}_{1}(1) \! = \! \mathfrak{u}_{2}(1) \! = \! \mathfrak{u}_{3}(1) \! = \! 0$, one 
analyses equations \eqref{geek102}--\eqref{geek107}, in the indicated order, in order to arrive at the following 
conclusions: (i) solving equation \eqref{geek102} for $\mathbb{P}$, one gets that
\begin{equation} \label{geek108} 
\mathbb{P} \! = \! -\dfrac{\mi \varepsilon (\varepsilon b)^{1/2} \me^{\mi \pi/4} \mathscr{P}_{a}(s_{0}^{0} 
\! - \! \mi \me^{-\pi a})}{\sqrt{\pi} \, 2^{3/2}3^{1/4}},
\end{equation}
whence, {}from the definition \eqref{geek101}, one arrives, again, at equation \eqref{geek109}; (ii) equations 
\eqref{geek103}--\eqref{geek105} are identically true; and (iii) solving equations \eqref{geek106} and 
\eqref{geek107} for $\iota_{0}^{\ast}(1)$ and $\iota_{1}^{\ast}(1)$, respectively, one concludes that
\begin{equation} \label{geek110} 
\iota_{0}^{\ast}(1) \! = \! \dfrac{\mi a (1 \! + \! \mi a)(\sqrt{3} \! + \! 1)}{18 \alpha_{1}^{4}} \qquad 
\text{and} \qquad \iota_{1}^{\ast}(1) \! = \! 0;
\end{equation}
moreover, {}from equations \eqref{nueeq5} and \eqref{nueeq6}, the definition \eqref{nueeq21}, and equations 
\eqref{geek110}, it also follows that
\begin{align}
\hat{b}_{0}^{\ast}(1) \! =& \, \dfrac{\mi (\sqrt{3} \! + \! 1)^{2}}{4 \sqrt{3}} \! \left(-\dfrac{\alpha_{1}}{2}(\mathfrak{r}_{0}^{2}
(1) \! + \! 2(\sqrt{3} \! + \! 1) \mathfrak{r}_{0}(1) \mathfrak{u}_{0}(1) \! + \! 8 \mathfrak{u}_{0}^{2}(1)) \right. 
\nonumber \\
+&\left. \frac{(a \! - \! \mi/2)}{6 \alpha_{1}}(12 \mathfrak{u}_{0}(1) \! + \! (2 \sqrt{3} \! - \! 1) \mathfrak{r}_{0}(1)) 
\right), \label{nueeq22} \\
\hat{b}_{1}^{\ast}(1) \! =& \, 0. \label{nueeq23}
\end{align}
Finally, {}from the asymptotics \eqref{tr1} and \eqref{tr3} (for $k \! = \! +1)$ and equation \eqref{geek109}, one 
arrives at the corresponding asymptotics for $v_{0}(\tau) \! := \! v_{0,1}(\tau)$ and $\tilde{r}_{0}(\tau) \! := \! 
\tilde{r}_{0,1}(\tau)$ stated in equations \eqref{geek1} and \eqref{geek2}, respectively, of the lemma.

Similarly, proceeding as delineated above, one shows that, for $k \! = \! -1$,
\begin{equation} \label{geek111}
\mathrm{A}_{-1} \! = \! \dfrac{\mi \me^{-\mi \pi/4} \me^{\mi \pi/3}(2 \! + \! \sqrt{3})^{-\mi a}(s_{0}^{0} \! - \! 
\mi \me^{-\pi a})}{\sqrt{2 \pi} \, 3^{1/4}(\varepsilon b)^{1/6}};
\end{equation}
thus, {}from the asymptotics \eqref{tr1} and \eqref{tr3} (for $k \! = \! -1)$ and equation \eqref{geek111}, one arrives 
at the corresponding asymptotics for $v_{0}(\tau) \! := \! v_{0,-1}(\tau)$ and $\tilde{r}_{0}(\tau) \! := \! \tilde{r}_{0,-1}
(\tau)$ stated in equations \eqref{geek1} and \eqref{geek2}, respectively, of the lemma. \hfill $\qed$

{}From equation \eqref{iden217}, the asymptotics \eqref{geek1}, the definition \eqref{geek4}, and recalling that 
(cf. equation \eqref{thmk2}) $c_{0,k} \! = \! \tfrac{1}{2} \varepsilon (\varepsilon b)^{2/3} \me^{-\mi 2 \pi k/3}$, 
$k \! = \! \pm 1$, one arrives at the corresponding $(\varepsilon_{1},\varepsilon_{2},m(\varepsilon_{2}) \vert \ell) 
\! = \! (0,0,0 \vert 0)$ asymptotics (as $\tau \! \to \! +\infty$ with $\varepsilon b \! > \! 0)$ for the solution $u(\tau)$ 
of the DP3E \eqref{eq1.1} stated in Theorem \ref{theor2.1}.

Via the definitions \eqref{hatsoff7} and \eqref{pga3} and equations \eqref{pga5} and \eqref{tr2}, one deduces that, 
for $k \! = \! \pm 1$,
\begin{gather}
2f_{-}(\tau) \! = \! -\mi (a \! - \! \mi/2) \! + \! \frac{\mi (\varepsilon b)^{1/3} \me^{\mi 2 \pi k/3}}{2} \tau^{2/3} \! 
\left(-2 \! + \! \tilde{r}_{0}(\tau) \tau^{-1/3} \right), \label{textfeqn2} \\
\frac{4 \mi}{\varepsilon b}f_{+}(\tau) \! = \! \mi (a \! + \! \mi/2) \! + \! \frac{\mi (\varepsilon b)^{1/3} \me^{\mi 2 \pi k/3}}{2} 
\tau^{2/3} \! \left(-2 \! + \! \tilde{r}_{0}(\tau) \tau^{-1/3} \right) \! + \! \frac{\mi b \tau}{u(\tau)}; \label{efhpls5}
\end{gather}
thus, {}from the asymptotics \eqref{geek1} and \eqref{geek2}, the definition \eqref{geek4}, and equations 
\eqref{textfeqn2} and \eqref{efhpls5}, one arrives at the corresponding $(\varepsilon_{1},\varepsilon_{2},
m(\varepsilon_{2}) \vert \ell) \! = \! (0,0,0 \vert 0)$ asymptotics (as $\tau \! \to \! +\infty$ with $\varepsilon 
b \! > \! 0)$ for the principal auxiliary functions $f_{\pm}(\tau)$ (corresponding to $u(\tau)$) stated in 
Theorem \ref{theor2.1}.

It was shown in equation (4.25) of \cite{av2} that, in terms of the function $h_{0}(\tau)$, the Hamiltonian function 
$\mathcal{H}(\tau)$ (corresponding to $u(\tau))$ defined by equation \eqref{eqh1} is given by
\begin{equation} \label{eqnhtext1} 
\mathcal{H}(\tau) \! = \! 3(\varepsilon b)^{2/3} \tau^{1/3} \! + \! \dfrac{1}{2 \tau}(a \! - \! \mi/2)^{2} \! - \! 4 \tau^{-1/3}
h_{0}(\tau):
\end{equation}
via the definition \eqref{iden2}, and equation \eqref{eqnhtext1}, it follows that, in terms of the function $\hat{h}_{0}
(\tau) \! := \! \hat{h}_{0,k}(\tau)$,
\begin{equation} \label{eqnhtext2} 
\mathcal{H}(\tau) \! = \! 3(\varepsilon b)^{2/3} \me^{-\mi 2 \pi k/3} \tau^{1/3} \! + \! \dfrac{1}{2 \tau}(a \! - \! \mi/2)^{2} 
\! - \! 4 \tau^{1/3} \hat{h}_{0,k}(\tau), \quad k \! = \! \pm 1;
\end{equation}
consequently, {}from equation \eqref{expforeych}, the third relation of equations \eqref{expforkapp}, and equation 
\eqref{eqnhtext2}, upon recalling that $v_{0}(\tau) \! := \! v_{0,k}(\tau)$ and $\tilde{r}_{0}(\tau) \! := \! \tilde{r}_{0,k}(\tau)$, 
one shows that the Hamiltonian function is given by
\begin{align} \label{eqnhtext3} 
\mathcal{H}(\tau) \! =& \, 3(\varepsilon b)^{2/3} \me^{-\mi 2 \pi k/3} \tau^{1/3} \! + \! \dfrac{1}{2 \tau}(a \! - \! \mi/2)^{2} 
\! + \! \frac{\alpha_{k}^{2} \tau^{-1/3}}{1 \! + \! \tau^{-1/3}v_{0,k}(\tau)} \! \left(\vphantom{M^{M^{M^{M^{M}}}}} \! 
\alpha_{k}^{2} \left(8v_{0,k}^{2}(\tau) \! + \! (4v_{0,k}(\tau) \right. \right. \nonumber \\
-&\left. \left. \tilde{r}_{0,k}(\tau)) \tilde{r}_{0,k}(\tau) \! - \! \tau^{-1/3}v_{0,k}(\tau)(\tilde{r}_{0,k}(\tau))^{2} \right) 
\! + \! 4(a \! - \! \mi/2) \right), \quad k \! = \! \pm 1.
\end{align}
Finally, {}from the asymptotics \eqref{geek1} and \eqref{geek2}, the definition \eqref{geek4}, and equation \eqref{eqnhtext3}, 
one arrives at, after a lengthy, but otherwise straightforward, calculation, the corresponding $(\varepsilon_{1},\varepsilon_{2},
m(\varepsilon_{2}) \vert \ell) \! = \! (0,0,0 \vert 0)$ asymptotics (as $\tau \! \to \! +\infty$ with $\varepsilon b \! > \! 0)$ for the 
Hamiltonian function, $\mathcal{H}(\tau)$, stated in Theorem \ref{theor2.1}.

Via the definition \eqref{thmk23} and the asymptotics (as $\tau \! \to \! +\infty$ with $\varepsilon b \! > \! 0)$ for 
$f_{-}(\tau)$ and $\mathcal{H}(\tau)$ stated above, one arrives at the corresponding $(\varepsilon_{1},
\varepsilon_{2},m(\varepsilon_{2}) \vert \ell) \! = \! (0,0,0 \vert 0)$ asymptotics for the function $\sigma (\tau)$ 
stated in Theorem \ref{theor2.1}.
\begin{bbbbb} \label{isomonoabcd} 
Under the conditions of Lemma \ref{ginversion}, the functions $a(\tau)$, $b(\tau)$, $c(\tau)$, and $d(\tau)$, 
defining, via equations \eqref{eq3.2}, the solution of the corresponding system of isomonodromy deformations 
\eqref{newlax8}, have the following asymptotic representations: for $k \! = \! \pm 1$,
{\fontsize{10pt}{11pt}\selectfont
\begin{align}
\sqrt{\smash[b]{-a(\tau)b(\tau)}} \underset{\tau \to +\infty}{=}& \, \frac{(\varepsilon b)^{2/3} \me^{-\mi 2 \pi k/3}}{2} 
\! \left(1 \! + \! \sum_{m=0}^{\infty} \frac{\mathfrak{u}_{m}(k)}{(\tau^{1/3})^{m+2}} \right) \! - \! 
\frac{\mi (\varepsilon b)^{1/2} \me^{\mi \pi k/4}(\mathscr{P}_{a})^{k}(s_{0}^{0} \! - \! \mi \me^{-\pi a})}{\sqrt{\pi} 
\, 2^{3/2}3^{1/4} \tau^{1/3}} \nonumber \\
\times& \, \me^{-(\beta (\tau)+\mi k \vartheta (\tau))} \! \left(1 \! + \! \mathcal{O} \big(\tau^{-1/3} \big) \right), 
\label{isomk1} \\
a(\tau)d(\tau) \underset{\tau \to +\infty}{=}& \, -\frac{\mi (\varepsilon b)}{4} \! - \! \frac{\mi (\varepsilon b)^{2/3} 
\me^{-\mi 2 \pi k/3}}{4}(a \! - \! \mi/3) \tau^{-2/3} \! + \! \frac{\mi (\varepsilon b)}{8} \left(\mathfrak{r}_{1}(k) 
\! - \! 2 \mathfrak{u}_{1}(k) \right) \! \tau^{-1} \nonumber \\
+& \, (\tau^{-1/3})^{4} \sum_{m=0}^{\infty} \left(\frac{\mi (\varepsilon b)}{8} \left(\mathfrak{r}_{m+2}(k) \! - \! 2 
\mathfrak{u}_{m+2}(k) \right) \! - \! \frac{\mi (\varepsilon b)^{2/3} \me^{-\mi 2 \pi k/3}}{4}(a \! - \! \mi/2) 
\mathfrak{u}_{m}(k) \right. \nonumber \\
+&\left. \frac{\mi (\varepsilon b)}{8} \sum_{p=0}^{m} \mathfrak{u}_{p}(k) \mathfrak{r}_{m-p}(k) \right) \! 
(\tau^{-1/3})^{m} \! - \! \frac{k(\varepsilon b)^{5/6}3^{1/4} \me^{\mi \pi k/4}(\mathscr{P}_{a})^{k}(s_{0}^{0} 
\! - \! \mi \me^{-\pi a})}{4 \sqrt{2 \pi} \, \me^{\mi \pi k/3} \tau^{1/3}} \nonumber \\
\times& \, \me^{-(\beta (\tau)+\mi k \vartheta (\tau))} \! \left(1 \! + \! \mathcal{O} \big(\tau^{-1/3} \big) \right), 
\label{isomk2} \\
b(\tau)c(\tau) \underset{\tau \to +\infty}{=}& \, -\frac{\mi (\varepsilon b)}{4} \! - \! \frac{\mi (\varepsilon b)^{2/3} 
\me^{-\mi 2 \pi k/3}}{4}(a \! + \! \mi/3) \tau^{-2/3} \! - \! \frac{\mi (\varepsilon b)}{8} \left(\mathfrak{r}_{1}(k) 
\! - \! 2 \mathfrak{u}_{1}(k) \right) \! \tau^{-1} \nonumber \\
+& \, (\tau^{-1/3})^{4} \sum_{m=0}^{\infty} \left(-\frac{\mi (\varepsilon b)}{8} \left(\mathfrak{r}_{m+2}(k) \! - \! 2 
\mathfrak{u}_{m+2}(k) \right) \! - \! \frac{\mi (\varepsilon b)^{2/3} \me^{-\mi 2 \pi k/3}}{4}(a \! + \! \mi/2) 
\mathfrak{u}_{m}(k) \right. \nonumber \\
-&\left. \frac{\mi (\varepsilon b)}{8} \sum_{p=0}^{m} \mathfrak{u}_{p}(k) \mathfrak{r}_{m-p}(k) \right) \! 
(\tau^{-1/3})^{m} \! + \! \frac{k(\varepsilon b)^{5/6}3^{1/4} \me^{\mi \pi k/4}(\mathscr{P}_{a})^{k}(s_{0}^{0} 
\! - \! \mi \me^{-\pi a})}{4 \sqrt{2 \pi} \, \me^{\mi \pi k/3} \tau^{1/3}} \nonumber \\
\times& \, \me^{-(\beta (\tau)+\mi k \vartheta (\tau))} \! \left(1 \! + \! \mathcal{O} \big(\tau^{-1/3} \big) \right), 
\label{isomk3} \\
-c(\tau)d(\tau) \underset{\tau \to +\infty}{=}& \, \frac{(\varepsilon b)^{2/3} \me^{\mi \pi k/3}}{4} \! - \! 
\frac{a(\varepsilon b)^{1/3} \me^{\mi 2 \pi k/3}}{3} \tau^{-2/3} \! - \! \frac{(\varepsilon b)^{2/3} \me^{\mi \pi k/3}}{2} 
\mathfrak{u}_{1}(k) \tau^{-1} \nonumber \\
-& \, \left(\frac{1}{6}(a^{2} \! + \! 1/6) \! + \! \dfrac{(\varepsilon b)^{2/3} \me^{\mi \pi k/3}}{2} \mathfrak{u}_{2}(k) 
\right) \! (\tau^{-1/3})^{4} \! + \! (\tau^{-1/3})^{4} \sum_{m=1}^{\infty} \left(-\frac{(\varepsilon b)^{2/3} 
\me^{\mi \pi k/3}}{2} \right. \nonumber \\
\times&\left. \, \mathfrak{u}_{m+2}(k) \! + \! \frac{\mi (\varepsilon b)^{1/3} \me^{\mi 2 \pi k/3}}{8} \mathfrak{r}_{m}(k) 
\! - \! \dfrac{(\varepsilon b)^{1/3} \me^{\mi 2 \pi k/3}}{2}(a \! - \! \mi/2) \mathfrak{w}_{m}(k) \! - \! 
\frac{(\varepsilon b)^{2/3} \me^{\mi \pi k/3}}{2} \right. \nonumber \\
\times&\left. \, \sum_{p=0}^{m} \left(\left(\mathfrak{u}_{p}(k) \! + \! \frac{1}{2} \mathfrak{r}_{p}(k) \right) \! 
\mathfrak{w}_{m-p}(k) \! + \! \frac{1}{8} \mathfrak{r}_{p}(k) \mathfrak{r}_{m-p}(k) \right) \right) \! 
(\tau^{-1/3})^{m} \nonumber \\
-& \, \dfrac{\mi (\varepsilon b)^{1/2} \me^{\mi \pi k/4}(\mathscr{P}_{a})^{k}(s_{0}^{0} \! - \! \mi \me^{-\pi a})}{\sqrt{\pi} 
\, 2^{3/2}3^{1/4} \tau^{1/3}} \me^{-(\beta (\tau)+\mi k \vartheta (\tau))} \! \left(1 \! + \! \mathcal{O} \big(\tau^{-1/3} 
\big) \right), \label{isomk4}
\end{align}}
where the expansion coefficients $\mathfrak{u}_{m}(k)$ (resp., $\mathfrak{r}_{m}(k))$, $m \! \in \! \mathbb{Z}_{+}$, are 
given in equations \eqref{thmk2}--\eqref{thmk10} (resp., \eqref{thmk15} and \eqref{thmk16}$)$.
\end{bbbbb}

\emph{Proof}. If, for $k \! = \! \pm 1$, $g_{ij}$, $i,j \! \in \! \lbrace 1,2 \rbrace$, are $\tau$ dependent, then, functions 
whose asymptotics (as $\tau \! \to \! +\infty$ with $\varepsilon b \! > \! 0)$ are given by equations 
\eqref{geek1}--\eqref{geek3} satisfy the conditions \eqref{iden5}, \eqref{pc4}, \eqref{restr1}, \eqref{ellinfk2a}, and 
\eqref{ellinfk2b}; therefore, one can use the justification scheme suggested in \cite{a20} (see, also, \cite{a22}). {}From 
equations \eqref{iden4oldu}, \eqref{iden7}, \eqref{iden8}, and \eqref{iden10}, respectively, one shows, via the 
definitions \eqref{iden3} and \eqref{iden4}, that, for $k \! = \! \pm 1$,\footnote{Recall that (cf. Lemma \ref{ginversion}) 
$v_{0}(\tau) \! := \! v_{0,k}(\tau)$ and $\tilde{r}_{0}(\tau) \! := \! \tilde{r}_{0,k}(\tau)$, $k \! = \! \pm 1$.}
\begin{align}
\sqrt{\smash[b]{-a(\tau)b(\tau)}} =& \, \frac{(\varepsilon b)^{2/3} \me^{-\mi 2 \pi k/3}}{2} \big(1 \! + \! \tau^{-1/3}
v_{0,k}(\tau) \big), \label{isomk5} \\
a(\tau)d(\tau) =& \, \frac{\mi (\varepsilon b)}{8} \big(1 \! + \! \tau^{-1/3}v_{0,k}(\tau) \big)(-2 \! + \! \tau^{-1/3} 
\tilde{r}_{0,k}(\tau)) \nonumber \\
-& \, \frac{\mi (\varepsilon b)^{2/3} \me^{-\mi 2 \pi k/3}}{4}(a \! - \! \mi/2) \big(1 \! + \! \tau^{-1/3}v_{0,k}(\tau) \big) 
\tau^{-2/3}, \label{isomk6} \\
b(\tau)c(\tau) =& \, -\frac{\mi (\varepsilon b)}{2} \! - \! \frac{\mi (\varepsilon b)}{8} \big(1 \! + \! \tau^{-1/3}v_{0,k}
(\tau) \big)(-2 \! + \! \tau^{-1/3} \tilde{r}_{0,k}(\tau)) \nonumber \\
-& \, \frac{\mi (\varepsilon b)^{2/3} \me^{-\mi 2 \pi k/3}}{4}(a \! + \! \mi/2) \big(1 \! + \! \tau^{-1/3}v_{0,k}(\tau) \big) 
\tau^{-2/3}, \label{isomk7} \\
-c(\tau)d(\tau) =& \, -\frac{(\varepsilon b)^{2/3} \me^{\mi \pi k/3}}{4} \! \left(\frac{-2 \! + \! \tau^{-1/3} \tilde{r}_{0,k}
(\tau)}{1 \! + \! \tau^{-1/3}v_{0,k}(\tau)} \right) \! - \! \frac{(\varepsilon b)^{2/3} \me^{\mi \pi k/3}}{16}(-2 \! + \! 
\tau^{-1/3} \tilde{r}_{0,k}(\tau))^{2} \nonumber \\
-& \, \frac{1}{4}(a \! - \! \mi/2)(a \! + \! \mi/2) \tau^{-4/3} \! + \! \frac{(\varepsilon b)^{1/3} \me^{\mi 2 \pi k/3}}{2} \! 
\left(\mi (-2 \! + \! \tau^{-1/3} \tilde{r}_{0,k}(\tau))/4 \right. \nonumber \\
-&\left. \, \frac{(a \! - \! \mi/2)}{1 \! + \! \tau^{-1/3}v_{0,k}(\tau)} \right) \! \tau^{-2/3}. \label{isomk8}
\end{align}
Via the asymptotics \eqref{geek1} and \eqref{geek2}, and equations \eqref{isomk5}--\eqref{isomk8}, one arrives at the 
asymptotics (as $\tau \! \to \! +\infty$ with $\varepsilon b \! > \! 0)$ for the functions $\sqrt{\smash[b]{-a(\tau)b(\tau)}}$, 
$a(\tau)d(\tau)$, $b(\tau)c(\tau)$, and $-c(\tau)d(\tau)$ stated in equations \eqref{isomk1}--\eqref{isomk4}, respectively. 
\hfill $\qed$
\begin{eeeee} \label{intofmot} 
It is important to note that the asymptotics \eqref{isomk1}--\eqref{isomk4} are consistent with equation \eqref{iden6}$;$ 
moreover, via the definitions \eqref{newlax2}, equations \eqref{eq3.2}, and the asymptotics \eqref{geek3} and 
\eqref{isomk1}--\eqref{isomk4}, one arrives at the asymptotics (as $\tau \! \to \! +\infty$ with $\varepsilon b \! > \! 0)$ 
for the solution of the---original---system of isomonodromy deformations \eqref{eq1.4}. \hfill $\blacksquare$
\end{eeeee}
\appendix
\section{Appendix: Proof of Proposition \ref{recursys}} \label{apprecuru}
\textbf{\emph{Proof}.} 
As the exponentially small correction term does not contribute to the algebraic determination of the coefficients 
$\mathfrak{u}_{m}(k)$, $m \! \in \! \mathbb{Z}_{+}$, $k \! = \! \pm 1$, hereafter, only the following `truncated' 
(and differentiable) asymptotics of $u(\tau)$ will be considered (with abuse of notation, also denoted as $u(\tau))$:
\begin{equation} \label{recur16} 
u(\tau) \underset{\tau \to +\infty}{=} c_{0,k} \tau^{1/3} \! \left(1 \! + \! \tau^{-2/3} \sum_{m=0}^{\infty} 
\dfrac{\mathfrak{u}_{m}(k)}{(\tau^{1/3})^{m}} \right), \quad k \! = \! \pm 1.
\end{equation}
Via the asymptotics~\eqref{recur16}, one shows that
\begin{equation} \label{recur17} 
\dfrac{1}{u(\tau)} \underset{\tau \to +\infty}{=} \dfrac{\tau^{-1/3}}{c_{0,k}} 
\! \left(1 \! + \! \tau^{-2/3} \sum_{m=0}^{\infty} \dfrac{\mathfrak{w}_{m}
(k)}{(\tau^{1/3})^{m}} \right), \quad k \! = \! \pm 1,
\end{equation}
where $\mathfrak{w}_{m}(k)$, $m \! \in \! \mathbb{Z}_{+}$, are determined iteratively {}from equations 
\eqref{thmk8}; in particular (this will be required for the ensuing proof), for $k \! = \! \pm 1$,
\begin{align}
\mathfrak{w}_{0}(k) \! =& -\mathfrak{u}_{0}(k), \label{recur18} \\
\mathfrak{w}_{1}(k) \! =& -\mathfrak{u}_{1}(k), \label{recur19} \\
\mathfrak{w}_{2}(k) \! =& -\mathfrak{u}_{2}(k) \! + \! \mathfrak{u}_{0}^{2}
(k), \label{recur20} \\
\mathfrak{w}_{3}(k) \! =& -\mathfrak{u}_{3}(k) \! + \! 2 \mathfrak{u}_{0}
(k) \mathfrak{u}_{1}(k), \label{recur21} \\
\mathfrak{w}_{4}(k) \! =& -\mathfrak{u}_{4}(k) \! + \! 2 \mathfrak{u}_{0}(k) 
\mathfrak{u}_{2}(k) \! + \! \mathfrak{u}_{1}^{2}(k) \! - \! \mathfrak{u}_{0}^{3}
(k), \label{recur22} \\
\mathfrak{w}_{5}(k) \! =& -\mathfrak{u}_{5}(k) \! + \! 2 \mathfrak{u}_{0}
(k) \mathfrak{u}_{3}(k) \! + \! 2 \mathfrak{u}_{1}(k) \mathfrak{u}_{2}(k) 
\! - \! 3 \mathfrak{u}_{0}^{2}(k) \mathfrak{u}_{1}(k), \label{recur23} \\
\mathfrak{w}_{6}(k) \! =& -\mathfrak{u}_{6}(k) \! + \! 2 \mathfrak{u}_{0}(k) 
\mathfrak{u}_{4}(k) \! + \! 2 \mathfrak{u}_{1}(k) \mathfrak{u}_{3}(k) \! + \! 
\mathfrak{u}_{2}^{2}(k) \! - \! 3 \mathfrak{u}_{0}^{2}(k) \mathfrak{u}_{2}(k) 
\! - \! 3 \mathfrak{u}_{0}(k) \mathfrak{u}_{1}^{2}(k) \! + \! \mathfrak{u}_{0}^{4}(k), 
\label{recur24} \\
\mathfrak{w}_{7}(k) \! =& -\mathfrak{u}_{7}(k) \! + \! 2 \mathfrak{u}_{0}(k) 
\mathfrak{u}_{5}(k) \! + \! 2 \mathfrak{u}_{1}(k) \mathfrak{u}_{4}(k) \! + \! 
2 \mathfrak{u}_{2}(k) \mathfrak{u}_{3}(k) \! - \! 3 \mathfrak{u}_{3}(k) 
\mathfrak{u}_{0}^{2}(k) \! - \! 6 \mathfrak{u}_{0}(k) \mathfrak{u}_{1}(k) 
\mathfrak{u}_{2}(k) \nonumber \\
+& \, 4 \mathfrak{u}_{1}(k) \mathfrak{u}_{0}^{3}(k) \! - \! \mathfrak{u}_{1}^{3}(k). 
\label{recur25}
\end{align}
{}From equations \eqref{thmk8} and the asymptotics \eqref{recur16} and \eqref{recur17}, one shows that 
(cf. DP3E \eqref{eq1.1}), for $k \! = \! \pm 1$,
\begin{equation} \label{recur26} 
\dfrac{b^{2}}{u(\tau)} \underset{\tau \to +\infty}{=} \dfrac{b^{2} \tau^{-1/3}}{c_{0,k}} \! \left(1 \! - \! 
\mathfrak{u}_{0}(k) \tau^{-2/3} \! - \! \mathfrak{u}_{1}(k)(\tau^{-1/3})^{3} \! - \! (\tau^{-1/3})^{4} 
\sum_{m=0}^{\infty} \lambda_{m}(k)(\tau^{-1/3})^{m} \right),
\end{equation}
where $\lambda_{j}(k) \! := \! -\mathfrak{w}_{j+2}(k)$, $j \! \in \! \mathbb{Z}_{+}$,
{\fontsize{10pt}{11pt}\selectfont
\begin{align}
\dfrac{1}{\tau} \! \left(-8 \varepsilon u^{2}(\tau) \! + \! 2ab \right) 
\underset{\tau \to +\infty}{=}& \, -8 \varepsilon c_{0,k}^{2} \tau^{-1/3} \! + \! 
(2ab \! - \! 16 \varepsilon c_{0,k}^{2} \mathfrak{u}_{0}(k))(\tau^{-1/3})^{3} \! - 
\! 16 \varepsilon c_{0,k}^{2} \mathfrak{u}_{1}(k)(\tau^{-1/3})^{4} \nonumber \\
-& \, 8 \varepsilon c_{0,k}^{2}(\tau^{-1/3})^{5} \sum_{m=0}^{\infty} \! 
\left(2 \mathfrak{u}_{m+2}(k) \! + \! \sum_{p=0}^{m} \mathfrak{u}_{p}
(k) \mathfrak{u}_{m-p}(k) \right) \! (\tau^{-1/3})^{m}, \label{recur27}
\end{align}}
\begin{equation}
\dfrac{u^{\prime}(\tau)}{\tau} \underset{\tau \to +\infty}{=} 
\dfrac{1}{3}c_{0,k}(\tau^{-1/3})^{5} \! \left(1 \! - \! \tau^{-2/3} 
\sum_{m=0}^{\infty}(m \! + \! 1) \mathfrak{u}_{m}(k)(\tau^{-1/3})^{m} 
\right), \label{recur28}
\end{equation}
{\fontsize{10pt}{11pt}\selectfont 
\begin{align}
\dfrac{(u^{\prime}(\tau))^{2}}{u(\tau)} \underset{\tau \to +\infty}{=}& \, 
\dfrac{1}{9}c_{0,k}(\tau^{-1/3})^{5} \! \left(1 \! - \! 3 \mathfrak{u}_{0}
(k) \tau^{-2/3} \! - \! 5 \mathfrak{u}_{1}(k)(\tau^{-1/3})^{3} \! + \! 
(2 \mathfrak{u}_{0}^{2}(k) \! - \! \lambda_{0}(k) \! + \! \eta_{0}(k))
(\tau^{-1/3})^{4} \vphantom{M^{M^{M^{M^{M^{M^{M}}}}}}} \right. 
\nonumber \\
+&\left. \, (6 \mathfrak{u}_{0}(k) \mathfrak{u}_{1}(k) \! - \! \lambda_{1}
(k) \! + \! \eta_{1}(k))(\tau^{-1/3})^{5} \! + \! (4 \mathfrak{u}_{1}^{2}(k) 
\! - \! \lambda_{2}(k) \! + \! 2 \mathfrak{u}_{0}(k) \lambda_{0}(k) \! + \! 
\eta_{2}(k) \right. \nonumber \\
-&\left. \, \mathfrak{u}_{0}(k) \eta_{0}(k))(\tau^{-1/3})^{6} \! + \! 
(-\lambda_{3}(k) \! + \! 2 \mathfrak{u}_{0}(k) \lambda_{1}(k) \! + \! 
4 \mathfrak{u}_{1}(k) \lambda_{0}(k) \! + \! \eta_{3}(k) \! - \! 
\mathfrak{u}_{0}(k) \eta_{1}(k) \right. \nonumber \\
-&\left. \, \mathfrak{u}_{1}(k) \eta_{0}(k))(\tau^{-1/3})^{7} \! + \! 
(\tau^{-1/3})^{8} \sum_{m=0}^{\infty} \! \left(-\lambda_{m+4}(k) \! + 
\! 2 \mathfrak{u}_{0}(k) \lambda_{m+2}(k) \! + \! 4 \mathfrak{u}_{1}(k) 
\lambda_{m+1}(k) \vphantom{M^{M^{M^{M^{M^{M^{M}}}}}}} \right. 
\right. \nonumber \\
+&\left. \left. \eta_{m+4}(k) \! - \! \mathfrak{u}_{0}(k) \eta_{m+2}(k) 
\! - \! \mathfrak{u}_{1}(k) \eta_{m+1}(k) \! - \! \sum_{p=0}^{m} \eta_{p}
(k) \lambda_{m-p}(k) \right) \! (\tau^{-1/3})^{m} \right), \label{recur29}
\end{align}}
where $\eta_{m}(k)$ is defined by equation~\eqref{thmk10}, and
\begin{equation} \label{recur30} 
u^{\prime \prime}(\tau) \underset{\tau \to +\infty}{=} -\dfrac{2}{9}
c_{0,k}(\tau^{-1/3})^{5} \! \left(1 \! - \! \tau^{-2/3} \sum_{m=0}^{\infty} 
\dfrac{(m \! + \! 1)(m \! + \! 4)}{2} \mathfrak{u}_{m}(k)(\tau^{-1/3})^{m} 
\right).
\end{equation}
Substituting, now, the expansions \eqref{recur26}--\eqref{recur30} into the DP3E \eqref{eq1.1}, and 
equating coefficients of like powers of $(\tau^{-1/3})^{m}$, $m \! \in \! \mathbb{N}$, one arrives at, for 
$k \! = \! \pm 1$, the following system of recurrence relations for the expansion coefficients 
$\mathfrak{u}_{m^{\prime}}(k)$, $m^{\prime} \! \in \! \mathbb{Z}_{+}$:
\begin{align}
&\mathcal{O} \! \left(\tau^{-1/3} \right) \colon & \quad 0 \! =& 
-8 \varepsilon c_{0,k}^{2} \! + \! b^{2}c_{0,k}^{-1}, \label{recur31} \\
&\mathcal{O} \! \left((\tau^{-1/3})^{3} \right) \colon & \quad 0 \! =& 
-16 \varepsilon c_{0,k}^{2} \mathfrak{u}_{0}(k) \! + \! 2ab \! - \! 
b^{2}c_{0,k}^{-1} \mathfrak{u}_{0}(k), \label{recur32} \\
&\mathcal{O} \! \left((\tau^{-1/3})^{4} \right) \colon & \quad 0 \! =& 
-16 \varepsilon c_{0,k}^{2} \mathfrak{u}_{1}(k) \! - \! b^{2}c_{0,k}^{-1} 
\mathfrak{u}_{1}(k), \label{recur33} \\
&\mathcal{O} \! \left((\tau^{-1/3})^{5} \right) \colon & \quad 0 \! =& 
\mathfrak{t}_{k}(2,0), \label{recur34} \\
&\mathcal{O} \! \left((\tau^{-1/3})^{6} \right) \colon & \quad 0 \! =& 
\mathfrak{t}_{k}(3,1), \label{recur35} \\
&\mathcal{O} \! \left((\tau^{-1/3})^{7} \right) \colon & \quad 
\dfrac{4}{9}c_{0,k} \mathfrak{u}_{0}(k) \! =& \mathfrak{t}_{k}
(4,2), \label{recur36} \\
&\mathcal{O} \! \left((\tau^{-1/3})^{8} \right) \colon & \quad c_{0,k} 
\mathfrak{u}_{1}(k) \! =& \mathfrak{t}_{k}(5,3), \label{recur37} \\
&\mathcal{O} \! \left((\tau^{-1/3})^{9} \right) \colon & \quad c_{0,k} 
\mathfrak{u}_{2}(k) \! =& \dfrac{1}{9}c_{0,k} \! \left(2 \mathfrak{u}_{0}^{2}
(k) \! - \! \lambda_{0}(k) \! + \! \eta_{0}(k) \right) \nonumber \\
& & +& \, \mathfrak{t}_{k}(6,4), \label{recur38} \\
&\mathcal{O} \! \left((\tau^{-1/3})^{10} \right) \colon & \quad \left(
\dfrac{4}{3} \right)^{2}c_{0,k} \mathfrak{u}_{3}(k) \! =& \dfrac{1}{9}
c_{0,k} \! \left(6 \mathfrak{u}_{0}(k) \mathfrak{u}_{1}(k) \! - \! 
\lambda_{1}(k) \! + \! \eta_{1}(k) \right) \nonumber \\
& & +& \, \mathfrak{t}_{k}(7,5), \label{recur39} \\
&\mathcal{O} \! \left((\tau^{-1/3})^{11} \right) \colon & \quad \left(
\dfrac{5}{3} \right)^{2}c_{0,k} \mathfrak{u}_{4}(k) \! =& \dfrac{1}{9}
c_{0,k} \! \left(4 \mathfrak{u}_{1}^{2}(k) \! - \! \lambda_{2}(k) \! + \! 
2 \mathfrak{u}_{0}(k) \lambda_{0}(k) \right. \nonumber \\
& & +&\left. \eta_{2}(k) \! - \! \mathfrak{u}_{0}(k) \eta_{0}(k) \right) 
\! + \! \mathfrak{t}_{k}(8,6), \label{recur40} \\
&\mathcal{O} \! \left((\tau^{-1/3})^{12} \right) \colon & \quad \left(
\dfrac{6}{3} \right)^{2}c_{0,k} \mathfrak{u}_{5}(k) \! =& \dfrac{1}{9}
c_{0,k} \! \left(-\lambda_{3}(k) \! + \! 2 \mathfrak{u}_{0}(k) \lambda_{1}
(k) \! + \! 4 \mathfrak{u}_{1}(k) \lambda_{0}(k) \right. \nonumber \\
& & +&\left. \eta_{3}(k) \! - \! \mathfrak{u}_{0}(k) \eta_{1}(k) \! - \! 
\mathfrak{u}_{1}(k) \eta_{0}(k) \right) \! + \! \mathfrak{t}_{k}(9,7), 
\label{recur41} \\
&\mathcal{O} \! \left((\tau^{-1/3})^{m+13} \right) \colon & \quad \left(
\dfrac{m \! + \! 7}{3} \right)^{2}c_{0,k} \mathfrak{u}_{m+6}(k) \! =& 
\dfrac{1}{9}c_{0,k} \! \left(-\lambda_{m+4}(k) \! + \! 2 \mathfrak{u}_{0}
(k) \lambda_{m+2}(k) \right. \nonumber \\
& & +&\left. 4 \mathfrak{u}_{1}(k) \lambda_{m+1}(k) \! + \! \eta_{m+4}
(k) \! - \! \mathfrak{u}_{0}(k) \eta_{m+2}(k) \right. \nonumber \\
& & -&\left. \mathfrak{u}_{1}(k) \eta_{m+1}(k) \! - \! \sum_{p=0}^{m} 
\eta_{p}(k) \lambda_{m-p}(k) \right) \nonumber \\
& & +& \, \mathfrak{t}_{k}(m \! + \! 10,m \! + \! 8), \quad m \! \in \! \mathbb{Z}_{+}, 
\label{recur42}
\end{align}
where
\begin{equation} \label{teekayjayell} 
\mathfrak{t}_{k}(j,l) \! := \! -8 \varepsilon c_{0,k}^{2} \! \left(2 
\mathfrak{u}_{j}(k) \! + \! \sum_{p=0}^{l} \mathfrak{u}_{p}(k) 
\mathfrak{u}_{l-p}(k) \right) \! - \! b^{2}c_{0,k}^{-1} \lambda_{l}(k).
\end{equation}
Noting that (cf. definition \eqref{thmk2}) equation \eqref{recur31} is identically true, the algorithm, hereafter, is as follows: 
(i) one solves equation \eqref{recur32} for $\mathfrak{u}_{0}(k)$ in order to arrive at the first of equations \eqref{thmk3}; 
(ii) via the formula for $\mathfrak{u}_{0}(k)$, the definitions of $c_{0,k}$, $\lambda_{i}(k)$, and $\eta_{m}(k)$ given 
heretofore, and equations \eqref{recur18}--\eqref{recur25}, one solves equations \eqref{recur33}--\eqref{recur41}, in the 
indicated order, to arrive at the expressions for the coefficients $\mathfrak{u}_{j}(k)$, $j \! = \! 1,2,\dotsc,9$, given in 
equations \eqref{thmk3} and \eqref{thmk4}; and (iii) using the fact that $\mathfrak{u}_{1}(k) \! = \! 0$ (cf. equations \eqref{thmk3}), 
and the definition of $\lambda_{i}(k)$, one solves equation \eqref{recur42} for $\mathfrak{u}_{m+10}(k)$, $m \! \in \! \mathbb{Z}_{+}$, 
and, after a lengthy induction argument, arrives at equations \eqref{thmk6} and~\eqref{thmk7}. \hfill $\qed$
\appendix
\setcounter{section}{1}
\section{Appendix: Asymptotics as $\tau \! \to \! +\infty$ for $\mathcal{Z}_{k}$, $\mathcal{G}_{0,k}$, $\mathfrak{A}_{k}$, 
$\mathfrak{B}_{k}$, $\mathfrak{C}_{k}$, $\mathfrak{A}_{0,k}^{\sharp}$, $\mathfrak{B}_{0,k}^{\sharp}$, 
$\mathfrak{C}_{0,k}^{\sharp}$, $\omega_{0,k}^{2}$, $\ell_{0,k}^{+}$, $\chi_{k}(\tau)$, $\ell_{1,k}^{+}$, $\mu_{k}(\tau)$, 
and $\ell_{2,k}^{+}$, $k \! = \! \pm 1$} \label{formforleadasymps}
For the requisite estimates in step \pmb{(xi)} of the proof of Lemma \ref{nprcl}, the $\tau \! \to \! +\infty$ asymptotics for 
$\mathcal{Z}_{k}$, $\mathcal{G}_{0,k}$, $\mathfrak{A}_{k}$, $\mathfrak{B}_{k}$, $\mathfrak{C}_{k}$, 
$\mathfrak{A}_{0,k}^{\sharp}$, $\mathfrak{B}_{0,k}^{\sharp}$, $\mathfrak{C}_{0,k}^{\sharp}$, $\omega_{0,k}^{2}$, 
$\ell_{0,k}^{+}$, $\chi_{k}(\tau)$, $\ell_{1,k}^{+}$, $\mu_{k}(\tau)$, and $\ell_{2,k}^{+}$, $k \! = \! \pm 1$, are necessary. 
{}From the conditions \eqref{iden5}, the asymptotics \eqref{tr1} and \eqref{tr3}, the definitions \eqref{prcy5}, \eqref{prcy9}, 
\eqref{prcy10}, \eqref{prcy15}, \eqref{prcy22}, \eqref{prcy33}, \eqref{prcy34}, \eqref{prcy35}, \eqref{prcy40}, 
\eqref{prcy45}, \eqref{prcy46}, \eqref{prcy53}, and \eqref{prcy54}, and equations \eqref{prcy15}--\eqref{prcy20}, a lengthy, 
but otherwise straightforward, algebraic calculation shows that, in the indicated order:
\begin{equation} \label{prcyzeek1} 
\mathcal{Z}_{k} \underset{\tau \to +\infty}{=} 1 \! - \! \dfrac{\tilde{r}_{0}(\tau) \tau^{-1/3}}{12} \! + \! \left(\dfrac{\tilde{r}_{0}
(\tau) \tau^{-1/3}}{12} \right)^{2} \! + \! \mathcal{O} \! \left(\big(\tilde{r}_{0}(\tau) \tau^{-1/3} \big)^{3} \right),
\end{equation}
\begin{equation} \label{prcyg1} 
\mathcal{G}_{0,k} \underset{\tau \to +\infty}{=} \mathcal{G}_{0,k}^{\infty} \! + \! \Delta \mathcal{G}_{0,k}, \quad k \! = \! \pm 1,
\end{equation}
where
\begin{equation} \label{prcyg2} 
(6 \varepsilon b)^{1/4} \mathcal{G}_{0,k}^{\infty} \! = \! 
\begin{pmatrix}
\frac{(\varepsilon b)^{1/2}(\sqrt{3}-1)}{\sqrt{2}} & 
-\frac{(\varepsilon b)^{1/2}(\sqrt{3}+1)}{\sqrt{2}} \\
1 & 1
\end{pmatrix},
\end{equation}
and
\begin{equation} \label{prcyg3} 
\Delta \mathcal{G}_{0,k} \! := \! \mathcal{G}_{0,k} \! - \! \mathcal{G}_{0,k}^{\infty} 
\! = \! 
\begin{pmatrix}
(\Delta \mathcal{G}_{0,k})_{11} & (\Delta \mathcal{G}_{0,k})_{12} \\
(\Delta \mathcal{G}_{0,k})_{21} & (\Delta \mathcal{G}_{0,k})_{22}
\end{pmatrix},
\end{equation}
with
{\fontsize{10pt}{11pt}\selectfont
\begin{align}
(6 \varepsilon b)^{1/4}(\Delta \mathcal{G}_{0,k})_{11} \! =& \, \dfrac{(\varepsilon b)^{1/2}}{4 \sqrt{2}} \! 
\left(\dfrac{(\sqrt{3} \! - \! 1)(2 \sqrt{3} \! + \! 1)}{6} \tilde{r}_{0}(\tau) \tau^{-1/3} \! + \! \dfrac{1}{12 \sqrt{3}} 
\! \left(1 \! + \! \dfrac{(\sqrt{3} \! - \! 1)(4 \sqrt{3} \! - \! 1)}{8 \sqrt{3}} \right) \right. \nonumber \\
\times&\left. (\tilde{r}_{0}(\tau) \tau^{-1/3})^{2} \! + \! \mathcal{O}((\tilde{r}_{0}(\tau) \tau^{-1/3})^{3}) \right), 
\label{prcyg4} \\
(6 \varepsilon b)^{1/4}(\Delta \mathcal{G}_{0,k})_{12} \! =& \, \dfrac{(\varepsilon b)^{1/2}}{4 \sqrt{2}} \! 
\left(\dfrac{(\sqrt{3} \! + \! 1)(2 \sqrt{3} \! - \! 1)}{6} \tilde{r}_{0}(\tau) \tau^{-1/3} \! + \! \dfrac{1}{12 \sqrt{3}} 
\! \left(-1 \! + \! \dfrac{(\sqrt{3} \! + \! 1)(4 \sqrt{3} \! + \! 1)}{8 \sqrt{3}} \right) \right. \nonumber \\
\times&\left. (\tilde{r}_{0}(\tau) \tau^{-1/3})^{2} \! + \! \mathcal{O}((\tilde{r}_{0}(\tau) \tau^{-1/3})^{3}) \right), 
\label{prcyg5} \\
(6 \varepsilon b)^{1/4}(\Delta \mathcal{G}_{0,k})_{21} \! =& \, (6 \varepsilon b)^{1/4}(\Delta \mathcal{G}_{0,k})_{22} 
\! = \! \dfrac{1}{24} \tilde{r}_{0}(\tau) \tau^{-1/3} \! - \! \dfrac{1}{2(24)^{2}}(\tilde{r}_{0}(\tau) \tau^{-1/3})^{2} 
\nonumber \\
+& \, \mathcal{O}((\tilde{r}_{0}(\tau) \tau^{-1/3})^{3}), \label{prcyg6}
\end{align}}
{\fontsize{9pt}{10pt}\selectfont 
\begin{align}
\mathfrak{A}_{k} \underset{\tau \to +\infty}{=}& \, \dfrac{\mi (a \! - \! 
\mi/2) \tau^{-1/3}}{\sqrt{3} \alpha_{k}} \! + \! \dfrac{\mi \tau^{-1/3}}{4 
\sqrt{3}} \! \left(\alpha_{k}(4v_{0}(\tau)(\tilde{r}_{0}(\tau) \! + \! 2v_{0}
(\tau)) \! - \! (\tilde{r}_{0}(\tau))^{2}) \! - \! \dfrac{(a \! - \! \mi/2)(12
v_{0}(\tau) \! - \! \tilde{r}_{0}(\tau)) \tau^{-1/3}}{3 \alpha_{k}} \right) 
\nonumber \\
+& \, \mathcal{O} \! \left((6 \varepsilon b)^{-1/2} \! \left(-\mi (\varepsilon 
b)^{1/3}((\varepsilon b)^{1/3} \me^{\mi \pi k/3}(\tilde{r}_{0}(\tau) \! + \! 
2v_{0}(\tau)) \! + \! 2(a \! - \! \mi/2) \me^{\mi 2 \pi k/3} \tau^{-1/3})
(v_{0}(\tau) \tau^{-1/3})^{2} \right. \right. \nonumber \\
+&\left. \left. \dfrac{\mi (\varepsilon b)^{1/3}}{12} \! \left(-\dfrac{
(\varepsilon b)^{1/3} \me^{\mi \pi k/3}}{4}(\tilde{r}_{0}(\tau))^{2} 
\tau^{-1/3} \! + \! ((\varepsilon b)^{1/3} \me^{\mi \pi k/3}(\tilde{r}_{0}
(\tau) \! + \! 2v_{0}(\tau)) \! + \! 2(a \! - \! \mi/2) \me^{\mi 2 \pi k/3} 
\tau^{-1/3}) \right. \right. \right. \nonumber \\
\times&\left. \left. \left. v_{0}(\tau) \tau^{-1/3} \right) 
\tilde{r}_{0}(\tau) \tau^{-1/3} \right) \right), \label{prcyak1} \\
\mathfrak{B}_{k} \underset{\tau \to +\infty}{=}& \, \mi (\sqrt{3} \! + \! 1) \! 
\left(\dfrac{\alpha_{k}}{2}(4v_{0}(\tau) \! + \! (\sqrt{3} \! + \! 1) \tilde{r}_{0}
(\tau)) \! - \! \dfrac{(\sqrt{3} \! + \! 1)(a \! - \! \mi/2) \tau^{-1/3}}{2 
\sqrt{3} \alpha_{k}} \right) \! + \! \dfrac{\mi (\sqrt{3} \! + \! 1)^{2} 
\tau^{-1/3}}{4 \sqrt{3}} \nonumber \\
\times& \left(-\dfrac{\alpha_{k}}{2}((\tilde{r}_{0}(\tau))^{2} \! + \! 
2(\sqrt{3} \! + \! 1)v_{0}(\tau) \tilde{r}_{0}(\tau) \! + \! 8v_{0}^{2}(\tau)) 
\! + \! \dfrac{(a \! - \! \mi/2)(12v_{0}(\tau) \! + \! (2 \sqrt{3} \! - \! 1) 
\tilde{r}_{0}(\tau)) \tau^{-1/3}}{6 \alpha_{k}} \right) \nonumber \\
+& \, \mathcal{O} \! \left((6 \varepsilon b)^{-1/2} \! \left(-\dfrac{\mi (\sqrt{3} 
\! + \! 1)^{2}(a \! - \! \mi/2)(\varepsilon b)^{1/3} \me^{\mi 2 \pi k/3}}{12} 
\tilde{r}_{0}(\tau)(\tau^{-1/3})^{3} \left(v_{0}(\tau) \! + \! \tilde{r}_{0}
(\tau)/2 \sqrt{3} \right) \right. \right. \nonumber \\
-&\left. \left. \dfrac{\mi (\sqrt{3} \! + \! 1)(\varepsilon b)^{2/3} 
\me^{\mi \pi k/3}}{48 \sqrt{3}} \tilde{r}_{0}(\tau)(\tau^{-1/3})^{2}
((\tilde{r}_{0}(\tau))^{2} \! + \! (\sqrt{3} \! + \! 1)(\tilde{r}_{0}(\tau) \! + \! 
2 \sqrt{3} \, v_{0}(\tau))(\tilde{r}_{0}(\tau) \! + \! 2v_{0}(\tau))) \right. \right. 
\nonumber \\
+&\left. \left. \dfrac{\mi \alpha_{k}(\varepsilon b)^{1/2}}{24 \sqrt{6}}
(\tilde{r}_{0}(\tau))^{3}(\tau^{-1/3})^{2} \! + \! \left(\dfrac{\mi (\varepsilon 
b)^{1/3}(\sqrt{3} \! + \! 1)^{2}}{2}(v_{0}(\tau) \tau^{-1/3})^{2} \! + \! 
\dfrac{\mi (\varepsilon b)^{1/3}(3 \sqrt{3} \! + \! 4)}{48 \sqrt{3}}(\tilde{r}_{0}
(\tau) \tau^{-1/3})^{2} \right. \right. \right. \nonumber \\
+&\left. \left. \left. \dfrac{\mi (\varepsilon b)^{1/3}(2 \! + \! \sqrt{3})}{2 
\sqrt{3}}v_{0}(\tau) \tilde{r}_{0}(\tau)(\tau^{-1/3})^{2} \right) \! \left(
(\varepsilon b)^{1/3} \me^{\mi \pi k/3}(\tilde{r}_{0}(\tau) \! + \! 2v_{0}(\tau)) 
\! + \! 2(a \! - \! \mi/2) \me^{\mi 2 \pi k/3} \tau^{-1/3} \right) \right) 
\right), \label{prcybk1} \\
\mathfrak{C}_{k} \underset{\tau \to +\infty}{=}& \, -\mi (\sqrt{3} \! - \! 1) 
\! \left(\dfrac{\alpha_{k}}{2}(4v_{0}(\tau) \! - \! (\sqrt{3} \! - \! 1) \tilde{r}_{0}
(\tau)) \! - \! \dfrac{(\sqrt{3} \! - \! 1)(a \! - \! \mi/2) \tau^{-1/3}}{2 \sqrt{3} 
\alpha_{k}} \right) \! + \! \dfrac{\mi (\sqrt{3} \! - \! 1)^{2} \tau^{-1/3}}{4 
\sqrt{3}} \nonumber \\
\times& \left(\dfrac{\alpha_{k}}{2}((\tilde{r}_{0}(\tau))^{2} \! - \! 2
(\sqrt{3} \! - \! 1)v_{0}(\tau) \tilde{r}_{0}(\tau) \! + \! 8v_{0}^{2}(\tau)) 
\! - \! \dfrac{(a \! - \! \mi/2)(12v_{0}(\tau) \! - \! (2 \sqrt{3} \! + \! 1) 
\tilde{r}_{0}(\tau)) \tau^{-1/3}}{6 \alpha_{k}} \right) \nonumber \\
+& \, \mathcal{O} \! \left((6 \varepsilon b)^{-1/2} \! \left(\dfrac{\mi (\sqrt{3} 
\! - \! 1)^{2}(a \! - \! \mi/2)(\varepsilon b)^{1/3} \me^{\mi 2 \pi k/3}}{12} 
\tilde{r}_{0}(\tau)(\tau^{-1/3})^{3} \left(v_{0}(\tau) \! - \! \tilde{r}_{0}
(\tau)/2 \sqrt{3} \right) \right. \right. \nonumber \\
+&\left. \left. \dfrac{\mi (\sqrt{3} \! - \! 1)(\varepsilon b)^{2/3} 
\me^{\mi \pi k/3}}{48 \sqrt{3}} \tilde{r}_{0}(\tau)(\tau^{-1/3})^{2}
((\tilde{r}_{0}(\tau))^{2} \! + \! (\sqrt{3} \! - \! 1)(2 \sqrt{3} \, v_{0}(\tau) \! 
- \! \tilde{r}_{0}(\tau))(\tilde{r}_{0}(\tau) \! + \! 2v_{0}(\tau))) \right. \right. 
\nonumber \\
+&\left. \left. \dfrac{\mi \alpha_{k}(\varepsilon b)^{1/2}}{24 \sqrt{6}}
(\tilde{r}_{0}(\tau))^{3}(\tau^{-1/3})^{2} \! - \! \left(\dfrac{\mi (\varepsilon 
b)^{1/3}(\sqrt{3} \! - \! 1)^{2}}{2}(v_{0}(\tau) \tau^{-1/3})^{2} \! + \! 
\dfrac{\mi (\varepsilon b)^{1/3}(3 \sqrt{3} \! - \! 4)}{48 \sqrt{3}}
(\tilde{r}_{0}(\tau) \tau^{-1/3})^{2} \right. \right. \right. \nonumber \\
-&\left. \left. \left. \dfrac{\mi (\varepsilon b)^{1/3}(2 \! - \! \sqrt{3})}{2 
\sqrt{3}}v_{0}(\tau) \tilde{r}_{0}(\tau)(\tau^{-1/3})^{2} \right) \! \left(
(\varepsilon b)^{1/3} \me^{\mi \pi k/3}(\tilde{r}_{0}(\tau) \! + \! 2v_{0}
(\tau)) \! + \! 2(a \! - \! \mi/2) \me^{\mi 2 \pi k/3} \tau^{-1/3} \right) 
\right) \right), \label{prcyck1}
\end{align}}
\begin{gather}
\mathfrak{A}_{0,k}^{\sharp} \underset{\tau \to +\infty}{=} -\dfrac{14 \mi \tau^{-1/3}}{\sqrt{3} \alpha_{k}} \! - \! 
\dfrac{\mi \tilde{r}_{0}(\tau)(\tau^{-1/3})^{2}}{4 \sqrt{3} \alpha_{k}} \! \left(-\dfrac{4}{3} \! + \! \dfrac{1}{2} 
\tilde{r}_{0}(\tau) \tau^{-1/3} \! + \! \mathcal{O}((\tilde{r}_{0}(\tau) \tau^{-1/3})^{2}) \right), \label{prcya0k1} \\
\mathfrak{B}_{0,k}^{\sharp} \underset{\tau \to +\infty}{=} \dfrac{4 \mi (\sqrt{3} \! + \! 1) \tau^{-1/3}}{\sqrt{3} 
\alpha_{k}} \! + \! \dfrac{\mi \tilde{r}_{0}(\tau)(\tau^{-1/3})^{2}}{4 \sqrt{3} \alpha_{k}} \! \left(-\dfrac{2(3 \sqrt{3} 
\! + \! 7)}{3} \! + \! \mathcal{O}((\tilde{r}_{0}(\tau) \tau^{-1/3})^{2}) \right), \label{prcyb0k1} \\
\mathfrak{C}_{0,k}^{\sharp} \underset{\tau \to +\infty}{=} \dfrac{4 \mi (\sqrt{3} \! - \! 1) \tau^{-1/3}}{\sqrt{3} 
\alpha_{k}} \! + \! \dfrac{\mi \tilde{r}_{0}(\tau)(\tau^{-1/3})^{2}}{4 \sqrt{3} \alpha_{k}} \! \left(-\dfrac{2(3 \sqrt{3} 
\! - \! 7)}{3} \! + \! \mathcal{O}((\tilde{r}_{0}(\tau) \tau^{-1/3})^{2}) \right), \label{prcyc0k1}
\end{gather}
\begin{align}  
\omega_{0,k}^{2} \underset{\tau \to +\infty}{=}& \, -\alpha_{k}^{2}(8v_{0}^{2}(\tau) \! + \! 4v_{0}(\tau) 
\tilde{r}_{0}(\tau) \! - \! (\tilde{r}_{0}(\tau))^{2}) \! + \! 4(a \! - \! \mi/2)v_{0}(\tau) \tau^{-1/3} \nonumber \\
+& \, (4 \alpha_{k}^{2}v_{0}(\tau)(\tilde{r}_{0}(\tau) \! + \! 2v_{0}(\tau)) \! - \! 4(a \! - \! \mi/2)v_{0}(\tau) 
\tau^{-1/3})v_{0}(\tau) \tau^{-1/3} \nonumber \\
+& \, \mathcal{O} \! \left((-4 \alpha_{k}^{2}v_{0}(\tau)(\tilde{r}_{0}(\tau) \! + \! 2v_{0}(\tau)) \! + \! 
4(a \! - \! \mi/2)v_{0}(\tau) \tau^{-1/3})(v_{0}(\tau) \tau^{-1/3})^{2} \right), \label{prcyomg1} \\
\ell_{0,k}^{+} \underset{\tau \to +\infty}{=}& \, \dfrac{\mi}{8 \sqrt{3}} \! 
\left(1 \! + \! \dfrac{\tilde{r}_{0}(\tau) \tau^{-1/3}}{12} \! + \! \mathcal{O}
((\tilde{r}_{0}(\tau) \tau^{-1/3})^{3}) \right) \! \dfrac{\mathfrak{B}_{0,k}^{
\sharp}}{\mathfrak{B}_{k}} \nonumber \\
-& \, \dfrac{\mi \omega_{0,k}^{2}}{(8 \sqrt{3})^{3}} \! \left(1 \! + \! 
\dfrac{\tilde{r}_{0}(\tau) \tau^{-1/3}}{12} \! + \! \mathcal{O}((\tilde{r}_{0}
(\tau) \tau^{-1/3})^{3}) \right)^{3} \! \left(\dfrac{\mathfrak{B}_{0,k}^{
\sharp}}{\mathfrak{B}_{k}} \right)^{2} \nonumber \\
+& \, \mathcal{O} \! \left(\omega_{0,k}^{4} \! \left(1 \! + \! \dfrac{\tilde{r}_{0}
(\tau) \tau^{-1/3}}{12} \! + \! \mathcal{O}((\tilde{r}_{0}(\tau) \tau^{-1/3})^{3}) 
\right)^{5} \! \left(\dfrac{\mathfrak{B}_{0,k}^{\sharp}}{\mathfrak{B}_{k}} 
\right)^{3} \right), \label{prcyellok1} \\
\chi_{k}(\tau) \underset{\tau \to +\infty}{=}& \, \mi 4 \sqrt{3} \mathcal{Z}_{k} \! + \! 
\omega_{0,k}^{2} \ell_{0,k}^{+} \! + \! \dfrac{\mathfrak{R}_{0,k}^{\ast}(-\ell_{0,k}^{+} 
\! + \! 1)}{2(\mi 4 \sqrt{3} \mathcal{Z}_{k} \! + \! \omega_{0,k}^{2} \ell_{0,k}^{+})} \! 
- \! \dfrac{(\mathfrak{R}_{0,k}^{\ast}(-\ell_{0,k}^{+} \! + \! 1))^{2}}{8(\mi 4 \sqrt{3} 
\mathcal{Z}_{k} \! + \! \omega_{0,k}^{2} \ell_{0,k}^{+})^{3}} \nonumber \\
+& \, \mathcal{O} \! \left(\dfrac{(\mathfrak{R}_{0,k}^{\ast}(-\ell_{0,k}^{+} \! + \! 1))^{3}}{(\mi 4 \sqrt{3} 
\mathcal{Z}_{k} \! + \! \omega_{0,k}^{2} \ell_{0,k}^{+})^{5}} \right), \label{prcychik1} \\
\ell_{1,k}^{+} \underset{\tau \to +\infty}{=}& \, \dfrac{\mathfrak{R}_{0,k}^{\ast}}{2(\mi 4 \sqrt{3} \mathcal{Z}_{k} 
\! + \! \omega_{0,k}^{2} \ell_{0,k}^{+})} \! - \! \dfrac{(\mathfrak{R}_{0,k}^{\ast})^{2}(-\ell_{0,k}^{+} \! + \! 1)}{8
(\mi 4 \sqrt{3} \mathcal{Z}_{k} \! + \! \omega_{0,k}^{2} \ell_{0,k}^{+})^{3}} \! + \! \mathcal{O} \! \left(
\dfrac{(\mathfrak{R}_{0,k}^{\ast})^{3}(-\ell_{0,k}^{+} \! + \! 1)^{2}}{(\mi 4 \sqrt{3} \mathcal{Z}_{k} \! + \! 
\omega_{0,k}^{2} \ell_{0,k}^{+})^{5}} \right), \label{prcyell1k1} \\
\mu_{k}(\tau) \underset{\tau \to +\infty}{=}& \, \chi_{k}(\tau) \! - \! \dfrac{\mi 8 \sqrt{3} \mathcal{Z}_{k} 
\mathfrak{A}_{k}(-\ell_{0,k}^{+} \! + \! 1)(\mathfrak{P}_{0,k}^{\ast} \! - \! \mi 8 \sqrt{3} \mathcal{Z}_{k} 
\mathfrak{A}_{k} \ell_{0,k}^{+})}{2 \chi_{k}^{2}(\tau)} \nonumber \\
-& \, \dfrac{(\mi 8 \sqrt{3} \mathcal{Z}_{k} \mathfrak{A}_{k}(-\ell_{0,k}^{+} \! + \! 1)(\mathfrak{P}_{0,k}^{\ast} 
\! - \! \mi 8 \sqrt{3} \mathcal{Z}_{k} \mathfrak{A}_{k} \ell_{0,k}^{+}))^{2}}{8 \chi_{k}^{5}(\tau)} \nonumber \\
+& \, \mathcal{O} \! \left(\dfrac{(\mi 8 \sqrt{3} \mathcal{Z}_{k} \mathfrak{A}_{k}(-\ell_{0,k}^{+} \! + \! 1)
(\mathfrak{P}_{0,k}^{\ast} \! - \! \mi 8 \sqrt{3} \mathcal{Z}_{k} \mathfrak{A}_{k} \ell_{0,k}^{+}))^{3}}{
\chi_{k}^{8}(\tau)} \right), \label{prcymuk1} 
\end{align}
and
\begin{align} \label{prcyell2k1} 
\ell_{2,k}^{+} \underset{\tau \to +\infty}{=}& \, -\dfrac{\mi 4 \sqrt{3} \mathcal{Z}_{k} \mathfrak{A}_{k}
(\mathfrak{P}_{0,k}^{\ast} \! - \! \mi 8 \sqrt{3} \mathcal{Z}_{k} \mathfrak{A}_{k} \ell_{0,k}^{+})}{\chi_{k}^{2}(\tau)} 
\! - \! \dfrac{(-\ell_{0,k}^{+} \! + \! 1)(\mi 8 \sqrt{3} \mathcal{Z}_{k} \mathfrak{A}_{k}(\mathfrak{P}_{0,k}^{\ast} 
\! - \! \mi 8 \sqrt{3} \mathcal{Z}_{k} \mathfrak{A}_{k} \ell_{0,k}^{+}))^{2}}{8 \chi_{k}^{5}(\tau)} \nonumber \\
+& \, \mathcal{O} \! \left(\dfrac{(-\ell_{0,k}^{+} \! + \! 1)^{2}(\mi 8 \sqrt{3} \mathcal{Z}_{k} \mathfrak{A}_{k}
(\mathfrak{P}_{0,k}^{\ast} \! - \! \mi 8 \sqrt{3} \mathcal{Z}_{k} \mathfrak{A}_{k} \ell_{0,k}^{+}))^{3}}{
\chi_{k}^{8}(\tau)} \right).
\end{align}
\appendix
\setcounter{section}{2}
\section{Appendix: Asymptotic Estimates as $\tau \! \to \! +\infty$ for $\lvert (\Phi_{M,k}(\xi))_{ij} \rvert$, 
$k \! = \! \pm 1$, $i,j \! = \! 1,2$, on the Stokes Rays} \label{appenformoduli} 
Asymptotic estimates as $\tau \! \to \! +\infty$ for the moduli $\lvert (\Phi_{M,k}(\xi))_{ij} \rvert$, $k \! = \! 
\pm 1$, $i,j \! = \! 1,2$, on the Stokes rays $\pmb{\widehat{\mathcal{S}}}$ are: \textbf{(a)} for $\arg (\xi) 
\genfrac{}{}{0pt}{3}{=}{\tau \to +\infty} 0 \! + \! \mathcal{O} \big(\tau^{-2/3} \big)$,\footnote{The asymptotic estimate 
$\mathcal{O} \big(\tau^{-2/3} \big)$ appears on the Stokes rays because of the factor $(2 \mu_{k}(\tau))^{1/2}$ in the 
arguments of the various parabolic-cylinder functions in equation \eqref{prcy2} and the fact that (cf. expansions 
\eqref{prcyzeek1}, \eqref{prcychik1}, and \eqref{prcymuk1}) $\arg (\mu_{k}(\tau)) \genfrac{}{}{0pt}{3}{=}{\tau \to +\infty} 
\tfrac{\pi}{2} \big(1 \! + \! \mathcal{O}(\tau^{-2/3}) \big)$.}
{\fontsize{10pt}{11pt}\selectfont
\begin{align*}
\lvert (\Phi_{M,k}(\xi))_{11} \rvert \underset{\tau \to +\infty}{\leqslant}& \, 
\left(\dfrac{2^{3/2} \me^{\pi \Im (\nu (k)+1)/2}2^{\Re (\nu (k))/2} \cosh^{3}
(\frac{\pi}{2} \Im (\nu (k) \! + \! 1)) \Gamma (-\Re (\nu (k)))}{\Gamma 
(\frac{1}{2} \! - \! \frac{\Re (\nu (k))}{2}) \sin (-\frac{\pi}{2} \Re (\nu (k)))} \right. \\
+&\left. \, \dfrac{\sqrt{\pi} \, \me^{\pi \Im (\nu (k)+1)}2^{-\Re (\nu (k)+1)/2} 
\lvert \sin (\frac{\pi}{2}(\nu (k) \! + \! 1)) \rvert}{\Gamma (\frac{1}{2} \! + \! 
\frac{\Re (\nu (k)+1)}{2}) \sin (\frac{\pi}{2} \Re (\nu (k) \! + \! 1))} \right) \! 
\left(1 \! + \! \mathcal{O}(\tau^{-2/3}) \right), \\
\lvert (\Phi_{M,k}(\xi))_{12} \rvert \underset{\tau \to +\infty}{\leqslant}& \, 
\dfrac{\sqrt{\pi} \, 2^{\Re (\nu (k))/2} \cosh (\frac{\pi}{2} \Im (\nu (k) \! + \! 1))}{
\Gamma (\frac{1}{2} \! - \! \frac{\Re (\nu (k))}{2}) \sin (-\frac{\pi}{2} \Re (\nu (k)))} 
\! \left(1 \! + \! \mathcal{O}(\tau^{-2/3}) \right), \\
\lvert (\Phi_{M,k}(\xi))_{21} \rvert \underset{\tau \to +\infty}{\leqslant}& \, 
\dfrac{4 \sqrt{3} \, \lvert \xi \rvert \Re (\nu (k) \! + \! 1)}{\lvert p_{k}(\tau) \rvert} 
\! \left(\dfrac{\me^{\pi \Im (\nu (k)+1)}2^{\Re (\nu (k)+1)/2} \lvert \sin 
(\frac{\pi}{2}(\nu (k) \! + \! 1)) \rvert \Gamma (\frac{\Re (\nu (k)+1)}{2})}{\sin 
(\frac{\pi}{2} \Re (\nu (k) \! + \! 1)) \Gamma (\Re (\nu (k) \! + \! 1))} \right. \\
+&\left. \, \dfrac{2^{3/2} \me^{\pi \Im (\nu (k)+1)/2}2^{-\Re (\nu (k))/2} 
\cosh^{3}(\frac{\pi}{2} \Im (\nu (k) \! + \! 1)) \Gamma (-\frac{\Re (\nu (k))}{2})}{
\sqrt{\pi} \sin (-\frac{\pi}{2} \Re (\nu (k)))} \right) \\
\times& \left(1 \! + \! \mathcal{O}(\tau^{-2/3}) \right), \\
\lvert (\Phi_{M,k}(\xi))_{22} \rvert \underset{\tau \to +\infty}{\leqslant}& \, 
\dfrac{4 \sqrt{3} \, \lvert \xi \rvert \Re (\nu (k) \! + \! 1)2^{-\Re (\nu (k))/2} 
\cosh (\frac{\pi}{2} \Im (\nu (k) \! + \! 1)) \Gamma (-\frac{\Re (\nu (k))}{2})}{
\lvert p_{k}(\tau) \rvert \sin (-\frac{\pi}{2} \Re (\nu (k))) \Gamma (-\Re (\nu (k)))} \\
\times& \left(1 \! + \! \mathcal{O}(\tau^{-2/3}) \right);
\end{align*}}
\textbf{(b)} for $\arg (\xi) \genfrac{}{}{0pt}{3}{=}{\tau \to +\infty} -\pi/2 \! + \! \mathcal{O} \big(\tau^{-2/3} \big)$,
\begin{align}
\lvert (\Phi_{M,k}(\xi))_{11} \rvert \underset{\tau \to +\infty}{\leqslant}& \, 
\dfrac{\sqrt{\pi} \, 2^{-\Re (\nu (k)+1)/2} \lvert \sin (\frac{\pi}{2}(\nu (k) \! + \! 1)) 
\rvert}{\Gamma (\frac{1}{2} \! + \! \frac{\Re (\nu (k)+1)}{2}) \sin (\frac{\pi}{2} \Re 
(\nu (k) \! + \! 1))} \! \left(1 \! + \! \mathcal{O}(\tau^{-2/3}) \right) \nonumber \\
=:& \, \hat{\varrho}_{0}(k) \! \left(1 \! + \! \mathcal{O}(\tau^{-2/3}) \right), 
\label{prcy67} \\
\lvert (\Phi_{M,k}(\xi))_{12} \rvert \underset{\tau \to +\infty}{\leqslant}& \, 
\dfrac{\sqrt{\pi} \, 2^{\Re (\nu (k))/2} \cosh (\frac{\pi}{2} \Im (\nu (k) \! + \! 1))}{
\Gamma (\frac{1}{2} \! - \! \frac{\Re (\nu (k))}{2}) \sin (-\frac{\pi}{2} \Re (\nu (k)))} 
\! \left(1 \! + \! \mathcal{O}(\tau^{-2/3}) \right) \nonumber \\
=:& \, \hat{\varrho}_{1}(k) \! \left(1 \! + \! \mathcal{O}(\tau^{-2/3}) \right), 
\label{prcy68} \\
\lvert (\Phi_{M,k}(\xi))_{21} \rvert \underset{\tau \to +\infty}{\leqslant}& \, 
\dfrac{4 \sqrt{3} \, \lvert \xi \rvert \Re (\nu (k) \! + \! 1)2^{\Re (\nu (k)+1)/2} 
\Gamma (\frac{\Re (\nu (k)+1)}{2}) \lvert \sin (\frac{\pi}{2}(\nu (k) \! + \! 1)) 
\rvert}{\lvert p_{k}(\tau) \rvert \Gamma (\Re (\nu (k) \! + \! 1)) \sin 
(\frac{\pi}{2} \Re (\nu (k) \! + \! 1))} \nonumber \\
\times& \left(1 \! + \! \mathcal{O}(\tau^{-2/3}) \right) \! =: \! \hat{\varrho}_{2}
(k) \dfrac{\lvert \xi \rvert \Re (\nu (k) \! + \! 1)}{\lvert p_{k}(\tau) \rvert} \! \left(
1 \! + \! \mathcal{O}(\tau^{-2/3}) \right), \label{prcy69} \\
\lvert (\Phi_{M,k}(\xi))_{22} \rvert \underset{\tau \to +\infty}{\leqslant}& \, 
\dfrac{4 \sqrt{3} \, \lvert \xi \rvert \Re (\nu (k) \! + \! 1)2^{-\Re (\nu (k))/2} 
\cosh (\frac{\pi}{2} \Im (\nu (k) \! + \! 1)) \Gamma (-\frac{\Re (\nu (k))}{2})}{
\lvert p_{k}(\tau) \rvert \sin (-\frac{\pi}{2} \Re (\nu (k))) \Gamma (-\Re (\nu (k)))} 
\nonumber \\
\times& \left(1 \! + \! \mathcal{O}(\tau^{-2/3}) \right) \! =: \! \hat{\varrho}_{3}(k) 
\dfrac{\lvert \xi \rvert \Re (\nu (k) \! + \! 1)}{\lvert p_{k}(\tau) \rvert} \! \left(
1 \! + \! \mathcal{O}(\tau^{-2/3}) \right); \label{prcy70} 
\end{align}
\textbf{(c)} for $\arg (\xi) \genfrac{}{}{0pt}{3}{=}{\tau \to +\infty} -\pi \! + \! \mathcal{O} \big(\tau^{-2/3} \big)$,
{\fontsize{10pt}{11pt}\selectfont
\begin{align*}
\lvert (\Phi_{M,k}(\xi))_{11} \rvert \underset{\tau \to +\infty}{\leqslant}& \, 
\dfrac{\sqrt{\pi} \, 2^{-\Re (\nu (k)+1)/2} \lvert \sin (\frac{\pi}{2}(\nu (k) \! + \! 
1)) \rvert}{\Gamma (\frac{1}{2} \! + \! \frac{\Re (\nu (k)+1)}{2}) \sin (\frac{\pi}{2} 
\Re (\nu (k) \! + \! 1))} \! \left(1 \! + \! \mathcal{O}(\tau^{-2/3}) \right), \\
\lvert (\Phi_{M,k}(\xi))_{12} \rvert \underset{\tau \to +\infty}{\leqslant}& \, 
\left(\dfrac{2^{3/2} \me^{\pi \Im (\nu (k)+1)/2} \lvert \cos (\frac{\pi}{2}
(\nu (k) \! + \! 1)) \rvert \lvert \sin (\frac{\pi}{2}(\nu (k) \! + \! 1)) \rvert^{2} 
\Gamma (\Re (\nu (k) \! + \! 1))}{2^{\Re (\nu (k)+1)/2} \Gamma (\frac{1}{2} \! 
+ \! \frac{\Re (\nu (k)+1)}{2}) \sin (\frac{\pi}{2} \Re (\nu (k) \! + \! 1))} \right. \\
+&\left. \, \dfrac{\sqrt{\pi} \, \me^{\pi \Im (\nu (k)+1)}2^{\Re (\nu (k))/2} \cosh 
(\frac{\pi}{2} \Im (\nu (k)+1))}{\Gamma (\frac{1}{2} \! - \! \frac{\Re (\nu (k))}{2}) 
\sin (-\frac{\pi}{2} \Re (\nu (k)))} \right) \! \left(1 \! + \! \mathcal{O}
(\tau^{-2/3}) \right), \\
\lvert (\Phi_{M,k}(\xi))_{21} \rvert \underset{\tau \to +\infty}{\leqslant}& \, 
\dfrac{4 \sqrt{3} \, \lvert \xi \rvert \Re (\nu (k) \! + \! 1)2^{\Re (\nu (k)+1)/2} 
\lvert \sin (\frac{\pi}{2}(\nu (k) \! + \! 1)) \rvert \Gamma (\frac{\Re (\nu (k)
+1)}{2})}{\lvert p_{k}(\tau) \rvert \sin (\frac{\pi}{2} \Re (\nu (k) \! + \! 1)) 
\Gamma (\Re (\nu (k) \! + \! 1))} \\
\times& \left(1 \! + \! \mathcal{O}(\tau^{-2/3}) \right), \\
\lvert (\Phi_{M,k}(\xi))_{22} \rvert \underset{\tau \to +\infty}{\leqslant}& \, 
\dfrac{4 \sqrt{3} \, \lvert \xi \rvert \Re (\nu (k) \! + \! 1)}{\lvert p_{k}(\tau) 
\rvert} \! \left(\dfrac{\me^{\pi \Im (\nu (k)+1)} \cosh (\frac{\pi}{2} \Im 
(\nu (k) \! + \! 1)) \Gamma (-\frac{\Re (\nu (k))}{2})}{2^{\Re (\nu (k))/2} 
\sin (-\frac{\pi}{2} \Re (\nu (k))) \Gamma (-\Re (\nu (k)))} \right. \\
+&\left. \, \dfrac{2^{3/2} \me^{\pi \Im (\nu (k)+1)/2} \lvert \cos 
(\frac{\pi}{2}(\nu (k) \! + \! 1)) \rvert \lvert \sin (\frac{\pi}{2}(\nu (k) 
\! + \! 1)) \rvert^{2} \Gamma (\frac{\Re (\nu (k) +1)}{2})}{\sqrt{\pi} 
2^{-\Re (\nu (k)+1)/2} \sin (\frac{\pi}{2} \Re (\nu (k) \! + \! 1))} \right) \\
\times& \left(1 \! + \! \mathcal{O}(\tau^{-2/3}) \right);
\end{align*}}
and \textbf{(d)} for $\arg (\xi) \genfrac{}{}{0pt}{3}{=}{\tau \to +\infty} -3 \pi/2 \! + \! \mathcal{O} \big(\tau^{-2/3} \big)$,
{\fontsize{10pt}{11pt}\selectfont
\begin{align}
\lvert (\Phi_{M,k}(\xi))_{11} \rvert \underset{\tau \to +\infty}{\leqslant}& \, 
\left(\dfrac{2^{3/2} \me^{-\pi \Im (\nu (k)+1)/2}2^{\Re (\nu (k))/2} \cosh^{3}
(\frac{\pi}{2} \Im (\nu (k) \! + \! 1)) \Gamma (-\Re (\nu (k)))}{\Gamma (\frac{1}{2} 
\! - \! \frac{\Re (\nu (k))}{2}) \sin (-\frac{\pi}{2} \Re (\nu (k)))} \right. \nonumber \\
+&\left. \, \dfrac{\sqrt{\pi} \, \me^{-\pi \Im (\nu (k)+1)}2^{-\Re (\nu (k)+1)/2} 
\lvert \sin (\frac{\pi}{2}(\nu (k) \! + \! 1)) \rvert}{\Gamma (\frac{1}{2} \! + \! 
\frac{\Re (\nu (k)+1)}{2}) \sin (\frac{\pi}{2} \Re (\nu (k) \! + \! 1))} \right) 
\! \left(1 \! + \! \mathcal{O}(\tau^{-2/3}) \right) \nonumber \\
=:& \, \tilde{\varrho}_{0}(k) \! \left(1 \! + \! \mathcal{O}(\tau^{-2/3}) \right), 
\label{prcy71} \\
\lvert (\Phi_{M,k}(\xi))_{12} \rvert \underset{\tau \to +\infty}{\leqslant}& \, 
\left(\dfrac{2^{3/2} \me^{-\pi \Im (\nu (k)+1)/2} \lvert \cos (\frac{\pi}{2}
(\nu (k) \! + \! 1)) \rvert \lvert \sin (\frac{\pi}{2}(\nu (k) \! + \! 1)) \rvert^{2} 
\Gamma (\Re (\nu (k) \! + \! 1))}{2^{\Re (\nu (k)+1)/2} \Gamma (\frac{1}{2} 
\! + \! \frac{\Re (\nu (k)+1)}{2}) \sin (\frac{\pi}{2} \Re (\nu (k) \! + \! 1))} 
\right. \nonumber \\
+&\left. \, \dfrac{\sqrt{\pi} \, \me^{-\pi \Im (\nu (k)+1)}2^{\Re (\nu (k))/2} 
\cosh (\frac{\pi}{2} \Im (\nu (k) \! + \! 1))}{\Gamma (\frac{1}{2} \! - \! \frac{
\Re (\nu (k))}{2}) \sin (-\frac{\pi}{2} \Re (\nu (k)))} \right) \! \left(
1 \! + \! \mathcal{O}(\tau^{-2/3}) \right) \nonumber \\
=:& \, \tilde{\varrho}_{1}(k) \! \left(1 \! + \! \mathcal{O}(\tau^{-2/3}) \right), 
\label{prcy72} \\
\lvert (\Phi_{M,k}(\xi))_{21} \rvert \underset{\tau \to +\infty}{\leqslant}& \, 
\dfrac{4 \sqrt{3} \, \lvert \xi \rvert \Re (\nu (k) \! + \! 1)}{\lvert p_{k}(\tau) 
\rvert} \! \left(\dfrac{\me^{-\pi \Im (\nu (k)+1)/2} \lvert \sin (\frac{\pi}{2}
(\nu (k) \! + \! 1)) \rvert \Gamma (\frac{\Re (\nu (k)+1)}{2})}{2^{-\Re (\nu 
(k)+1)/2} \sin (\frac{\pi}{2} \Re (\nu (k) \! + \! 1)) \Gamma (\Re (\nu (k) 
\! + \! 1))} \right. \nonumber \\
+&\left. \, \dfrac{2^{3/2} \me^{-\pi \Im (\nu (k)+1)/2} \cosh^{3}(\frac{\pi}{2} 
\Im (\nu (k) \! + \! 1)) \Gamma (-\frac{\Re (\nu (k))}{2})}{2^{\Re (\nu (k))/2} 
\sin (-\frac{\pi}{2} \Re (\nu (k)))} \right) \! \left(1 \! + \! \mathcal{O}
(\tau^{-2/3}) \right) \nonumber \\
=:& \, \tilde{\varrho}_{2}(k) \dfrac{\lvert \xi \rvert \Re (\nu (k) \! + \! 1)}{
\lvert p_{k}(\tau) \rvert} \! \left(1 \! + \! \mathcal{O}(\tau^{-2/3}) \right), 
\label{prcy73} \\
\lvert (\Phi_{M,k}(\xi))_{22} \rvert \underset{\tau \to +\infty}{\leqslant}& \, 
\dfrac{4 \sqrt{3} \, \lvert \xi \rvert \Re (\nu (k) \! + \! 1)}{\lvert p_{k}(\tau) 
\rvert} \! \left(\dfrac{\me^{-\pi \Im (\nu (k)+1)} \cosh (\frac{\pi}{2} \Im (\nu 
(k) \! + \! 1)) \Gamma (-\frac{\Re (\nu (k))}{2})}{2^{\Re (\nu (k))/2} \sin 
(-\frac{\pi}{2} \Re (\nu (k))) \Gamma (-\Re (\nu (k)))} \right. \nonumber \\
+&\left. \, \dfrac{2^{3/2} \me^{-\pi \Im (\nu (k)+1)/2} \lvert \cos 
(\frac{\pi}{2}(\nu (k) \! + \! 1)) \rvert \lvert \sin (\frac{\pi}{2}(\nu (k) \! 
+ \! 1)) \rvert^{2} \Gamma (\frac{\Re (\nu (k)+1)}{2})}{\sqrt{\pi}2^{-\Re 
(\nu (k)+1)/2} \sin (\frac{\pi}{2} \Re (\nu (k) \! + \! 1))} \right) \nonumber \\
\times& \left(1 \! + \! \mathcal{O}(\tau^{-2/3}) \right) \! =: \! 
\tilde{\varrho}_{3}(k) \dfrac{\lvert \xi \rvert \Re (\nu (k) \! + \! 1)}{\lvert p_{k}
(\tau) \rvert} \! \left(1 \! + \! \mathcal{O}(\tau^{-2/3}) \right). \label{prcy74}
\end{align}}
\appendix
\setcounter{section}{3}
\section{Appendix: Symmetries and Transformations} \label{sectonsymm} 
It was shown in Proposition \ref{prop1.2} that (cf. system~\eqref{eq1.5}), given any solution $\hat{u}(\tau)$ of the DP3E 
\eqref{eq1.1}, the function $\hat{\varphi}(\tau)$ is defined as the general solution of the ODE $\hat{\varphi}^{\prime}(\tau) 
\! = \! 2a \tau^{-1} \! + \! b(\hat{u}(\tau))^{-1}$. {}From the latter ODE, it is clear that, given $\hat{u}(\tau)$, the function 
$\hat{\varphi}(\tau)$ is defined up to a $\tau$-independent ``additive parameter'', that is, $\hat{\varphi}(\tau) \! \to \! 
\hat{\varphi}(\tau) \! + \! \hat{\varphi}_{0}$, where $\hat{\varphi}_{0} \! \in \! \mathbb{C}$.\footnote{Of course, it also follows 
{}from the definitions \eqref{eq:ABCD} and \eqref{tempeq} that $\hat{\varphi}(\tau)$ is defined $\operatorname{mod}(2 \pi)$: 
similar statements apply, \emph{mutatis mutandis}, for the pair of functions $(u(\tau),\varphi (\tau))$ that solve the system 
\eqref{equu17}, where, in particular, $\varphi (\tau)$ is also defined $\operatorname{mod}(2 \pi)$ (cf. definitions \eqref{equu18} 
and \eqref{equu19}).} As the principal focus of the symmetry transformations derived in Section 6 of \cite{a1} was on the 
function $\hat{u}(\tau)$ and not the function $\hat{\varphi}(\tau)$, it must be noted that the additive parameter, 
$\hat{\varphi}_{0}$, appears non-uniformly (though correctly!) in those symmetries; for example, for the Transformation 6.2.1 
changing $\tau \! \to \! -\tau$, $\hat{\varphi}_{0} \! = \! -\pi \epsilon_{1}^{\ast}$, $\epsilon_{1}^{\ast} \! \in \! \lbrace \pm 1 
\rbrace$, whilst for the Transformation 6.2.3 changing $\tau \! \to \! \mi \tau$, $\hat{\varphi}_{0} \! = \! 0$. In order to, with 
abuse of nomenclature, ``uniformize'' the presentation of the final asymptotic results of the present work, this appendix 
considers the concomitant actions (see the brief discussion below) of the Lie-point symmetries for the DP3E \eqref{eq1.1} 
and the systems of isomonodromy deformations \eqref{eq1.4} and \eqref{newlax8} on the fundamental solutions of the 
systems \eqref{eqFGmain} and \eqref{newlax3} and the manifold of the monodromy data, $\mathscr{M}$,\footnote{The 
group of symmetries derived in this section preserve, in particular, the invariance of the system \eqref{monoeqns} defining 
$\mathscr{M}$.} under the strict caveat that, for every symmetry, the additive parameter is equal to zero; en route, novel 
sets of symmetry transformations not identified in \cite{a1} are obtained.

Before proceeding, however, some preamble regarding group actions on sets is necessary (see, for example, 
\cite{julzez}). The terms `function' and `transformation' will be used interchangeably throughout the following 
discussion. Let $\mathrm{G}$ be a group and $\mathrm{X}$ denote a set. An \emph{action\/} of $\mathrm{G}$ on 
$\mathrm{X}$ is a function {}from $\mathrm{G} \times \mathrm{X}$ to $\mathrm{X}$ if, for every pair $(\mathfrak{g},
\mathfrak{x}) \! \in \! \mathrm{G} \times \mathrm{X}$, there is an element $\mathfrak{g} \mathfrak{x} \! \in \! \mathrm{X}$ 
such that $(\mathfrak{g}_{1} \mathfrak{g}_{2}) \mathfrak{x} \! = \! \mathfrak{g}_{1}(\mathfrak{g}_{2} \mathfrak{x})$ and 
$\mathfrak{e} \mathfrak{x} \! = \! \mathfrak{x}$ ($\mathfrak{e}$ is the identity in $\mathrm{G}$). For fixed $\mathfrak{g} 
\! \in \! \mathrm{G}$, there is a function (transformation) $\aleph_{\mathfrak{g}} \colon \mathrm{X} \! \mapsto \! 
\mathfrak{g} \mathfrak{x}$ for $\mathfrak{x} \! \in \! \mathrm{X}$, that is, $\operatorname{Act}(\mathrm{G})_{
\scriptscriptstyle \mathrm{X}} \colon \mathrm{G} \times \mathrm{X} \! \to \! \mathrm{X}$, $(\mathfrak{g},\mathfrak{x}) 
\! \mapsto \! \aleph_{\mathfrak{g}}(\mathfrak{x}) \! := \! \mathfrak{g} \mathfrak{x}$. As $\aleph_{\mathfrak{g}_{1}} \circ 
\aleph_{\mathfrak{g}_{2}} \! = \! \aleph_{\mathfrak{g}_{1} \mathfrak{g}_{2}}$ and $\aleph_{\mathfrak{e}} \! = \! 
\id_{\scriptscriptstyle \mathrm{X}}$ (the identity mapping on $\mathrm{X}$), it follows that $\aleph_{\mathfrak{g}}$ is 
a bijection on $\mathrm{X}$, since $\aleph_{\mathfrak{g}} \circ \aleph_{\mathfrak{g}^{-1}} \! = \! \aleph_{\mathfrak{g} 
\mathfrak{g}^{-1}} \! = \! \aleph_{\mathfrak{e}} \! = \! \aleph_{\mathfrak{g}^{-1} \mathfrak{g}} \! = \! 
\aleph_{\mathfrak{g}^{-1}} \circ \aleph_{\mathfrak{g}}$, where $\aleph_{\mathfrak{g}^{-1}}$ denotes the inverse function 
of $\aleph_{\mathfrak{g}}$. All bijective functions $\aleph \colon \mathrm{X} \! \to \! \mathrm{X}$ form a group under 
composition of functions (the composition of functions is associative, the identity is the identity function $\id (\mathfrak{x}) 
\! = \! \mathfrak{x}$ for $\mathfrak{x} \! \in \! \mathrm{X}$, and the inverse of $\aleph$ is the inverse function $\aleph^{-1}$). 
Denoting by $\mathbb{B}(\mathfrak{x})$ the group of all bijections on $\mathrm{X}$, one defines a transformation group 
of $\mathrm{X}$ as any subgroup of $\mathbb{B}(\mathfrak{x})$.\footnote{In this work, the transformation group is a 
disjoint union of two subgroups of Lie-point symmetries for the DP3E \eqref{eq1.1} and the systems of isomonodromy 
deformations \eqref{eq1.4} and \eqref{newlax8}, and, in particular, the actions (symmetry transformations) of these 
subgroups on $\mathscr{M}$ is studied.} Any action of a group $\mathrm{G}$ on a set $\mathrm{X}$ defines a 
homomorphism {}from $\mathrm{G}$ to the transformation group $\mathbb{B}(\mathfrak{x})$ such that $\mathfrak{g} 
\! \in \! \mathrm{G}$ maps onto the transformation $\aleph_{\mathfrak{g}}$. Denoting such a homomorphism by 
$\daleth \colon \mathrm{G} \! \to \! \mathbb{B}(\mathfrak{x})$, it follows that $\daleth (\mathfrak{g}) \! = \! 
\aleph_{\mathfrak{g}}$; conversely, any homomorphism $\daleth \colon \mathrm{G} \! \to \! \mathbb{B}(\mathfrak{x})$ 
defines an action of $\mathrm{G}$ on $\mathrm{X}$ if one defines $\mathfrak{g} \mathfrak{x} \! := \! \daleth (\mathfrak{g})
(\mathfrak{x})$.\footnote{For $\mathfrak{g}_{1},\mathfrak{g}_{2} \! \in \! \mathrm{G}$ and $\mathfrak{x} \! \in \! \mathrm{X}$, 
the properties $\daleth (\mathfrak{g}_{1} \mathfrak{g}_{2}) \! = \! \daleth (\mathfrak{g}_{1}) \daleth (\mathfrak{g}_{2})$ and 
$\daleth (\mathfrak{e}) \! = \! \id$ imply that $(\mathfrak{g}_{1} \mathfrak{g}_{2}) \mathfrak{x} \! = \! \mathfrak{g}_{1}
(\mathfrak{g}_{2} \mathfrak{x})$ and $\mathfrak{e} \mathfrak{x} \! = \! \mathfrak{x}$.} For a group $\mathrm{G}$ acting 
on a set $\mathrm{X}$, the \emph{orbit} of $\mathfrak{x} \! \in \! \mathrm{X}$, denoted by $\mathrm{G} \mathfrak{x}$, 
is defined as $\mathrm{G} \mathfrak{x} \! := \! \lbrace \mathstrut \mathfrak{g} \mathfrak{x}, \, \forall \, \mathfrak{g} \! \in 
\! \mathrm{G} \rbrace$ (the set of all images of $\mathfrak{x}$ under the elements of $\mathrm{G}$).
\begin{eeeee} \label{remaction} 
In this work (see Appendix \ref{sectonsymmcomp} below for complete details), the group $\mathrm{G}$ of all (Lie-point) 
symmetries of interest is written as the disjoint union of two subgroups, $\mathrm{G} \! = \! \widetilde{\pmb{\mathcal{W}}} 
\cup \widehat{\pmb{\mathcal{W}}}$, where the elements of the subgroup $\widetilde{\pmb{\mathcal{W}}}$ are denoted by 
$\mathscr{F}^{\scriptscriptstyle \lbrace \ell \rbrace}_{\scriptscriptstyle \varepsilon_{1},\varepsilon_{2},m(\varepsilon_{2})}$, 
with $\varepsilon_{1} \! \in \! \lbrace 0,\pm 1 \rbrace$, $\varepsilon_{2} \! \in \! \lbrace 
0,\pm 1 \rbrace$, $m(\varepsilon_{2}) \! = \! 
\left\{
\begin{smallmatrix}
0, \, \, \varepsilon_{2}=0, \\
\pm \varepsilon_{2}, \, \, \varepsilon_{2} \in \lbrace \pm 1 \rbrace,
\end{smallmatrix}
\right.$ and $\ell \! \in \! \lbrace 0,1 \rbrace$, and the elements of the subgroup 
$\widehat{\pmb{\mathcal{W}}}$ are denoted by $\hat{\mathscr{F}}^{\scriptscriptstyle 
\lbrace \hat{\ell} \rbrace}_{\scriptscriptstyle \hat{\varepsilon}_{1},\hat{\varepsilon}_{2},
\hat{m}(\hat{\varepsilon}_{2})}$, with $\hat{\varepsilon}_{1} \! \in \! \lbrace \pm 1 
\rbrace$, $\hat{\varepsilon}_{2} \! \in \! \lbrace 0,\pm 1 \rbrace$, $\hat{m}
(\hat{\varepsilon}_{2}) \! = \! 
\left\{
\begin{smallmatrix}
0, \, \, \hat{\varepsilon}_{2} \in \lbrace \pm 1 \rbrace, \\
\pm \hat{\varepsilon}_{1}, \, \, \hat{\varepsilon}_{2}=0,
\end{smallmatrix}
\right.$ and $\hat{\ell} \! \in \! \lbrace 0,1 \rbrace$, and the action of the group elements 
$\mathscr{F}^{\scriptscriptstyle \lbrace \ell \rbrace}_{\scriptscriptstyle \varepsilon_{1},\varepsilon_{2},
m(\varepsilon_{2})}$ on $\mathscr{M}$,
\begin{align*}
\mathscr{F}^{\scriptscriptstyle \lbrace \ell \rbrace}_{\scriptscriptstyle \varepsilon_{1},
\varepsilon_{2},m(\varepsilon_{2})} \mathscr{M} :=& \, \left(\mathscr{F}^{\scriptscriptstyle 
\lbrace \ell \rbrace}_{\scriptscriptstyle \varepsilon_{1},\varepsilon_{2},m(\varepsilon_{2})}a,
\mathscr{F}^{\scriptscriptstyle \lbrace \ell \rbrace}_{\scriptscriptstyle \varepsilon_{1},
\varepsilon_{2},m(\varepsilon_{2})}s_{0}^{0},\mathscr{F}^{\scriptscriptstyle \lbrace \ell 
\rbrace}_{\scriptscriptstyle \varepsilon_{1},\varepsilon_{2},m(\varepsilon_{2})}s_{0}^{\infty},
\mathscr{F}^{\scriptscriptstyle \lbrace \ell \rbrace}_{\scriptscriptstyle \varepsilon_{1},
\varepsilon_{2},m(\varepsilon_{2})}s_{1}^{\infty},\mathscr{F}^{\scriptscriptstyle \lbrace 
\ell \rbrace}_{\scriptscriptstyle \varepsilon_{1},\varepsilon_{2},m(\varepsilon_{2})}g_{11}, 
\right. \nonumber \\
&\left. \, \mathscr{F}^{\scriptscriptstyle \lbrace \ell \rbrace}_{\scriptscriptstyle \varepsilon_{1},
\varepsilon_{2},m(\varepsilon_{2})}g_{12},\mathscr{F}^{\scriptscriptstyle \lbrace \ell 
\rbrace}_{\scriptscriptstyle \varepsilon_{1},\varepsilon_{2},m(\varepsilon_{2})}g_{21},
\mathscr{F}^{\scriptscriptstyle \lbrace \ell \rbrace}_{\scriptscriptstyle \varepsilon_{1},
\varepsilon_{2},m(\varepsilon_{2})}g_{22} \right),
\end{align*}
is given in equations \eqref{laxhat76}--\eqref{laxhat90} and \eqref{laxhat99}--\eqref{laxhat113} below, whilst the 
action of the group elements $\hat{\mathscr{F}}^{\scriptscriptstyle \lbrace \hat{\ell} \rbrace}_{\scriptscriptstyle 
\hat{\varepsilon}_{1},\hat{\varepsilon}_{2},\hat{m}(\hat{\varepsilon}_{2})}$ on $\mathscr{M}$,
\begin{align*}
\hat{\mathscr{F}}^{\scriptscriptstyle \lbrace \hat{\ell} \rbrace}_{\scriptscriptstyle 
\hat{\varepsilon}_{1},\hat{\varepsilon}_{2},\hat{m}(\hat{\varepsilon}_{2})} \mathscr{M} 
:=& \, \left(\hat{\mathscr{F}}^{\scriptscriptstyle \lbrace \hat{\ell} \rbrace}_{\scriptscriptstyle 
\hat{\varepsilon}_{1},\hat{\varepsilon}_{2},\hat{m}(\hat{\varepsilon}_{2})}a,
\hat{\mathscr{F}}^{\scriptscriptstyle \lbrace \hat{\ell} \rbrace}_{\scriptscriptstyle 
\hat{\varepsilon}_{1},\hat{\varepsilon}_{2},\hat{m}(\hat{\varepsilon}_{2})}s_{0}^{0},
\hat{\mathscr{F}}^{\scriptscriptstyle \lbrace \hat{\ell} \rbrace}_{\scriptscriptstyle 
\hat{\varepsilon}_{1},\hat{\varepsilon}_{2},\hat{m}(\hat{\varepsilon}_{2})}s_{0}^{\infty},
\hat{\mathscr{F}}^{\scriptscriptstyle \lbrace \hat{\ell} \rbrace}_{\scriptscriptstyle 
\hat{\varepsilon}_{1},\hat{\varepsilon}_{2},\hat{m}(\hat{\varepsilon}_{2})}s_{1}^{\infty},
\hat{\mathscr{F}}^{\scriptscriptstyle \lbrace \hat{\ell} \rbrace}_{\scriptscriptstyle 
\hat{\varepsilon}_{1},\hat{\varepsilon}_{2},\hat{m}(\hat{\varepsilon}_{2})}g_{11}, 
\right. \nonumber \\
&\left. \, \hat{\mathscr{F}}^{\scriptscriptstyle \lbrace \hat{\ell} \rbrace}_{\scriptscriptstyle 
\hat{\varepsilon}_{1},\hat{\varepsilon}_{2},\hat{m}(\hat{\varepsilon}_{2})}g_{12},
\hat{\mathscr{F}}^{\scriptscriptstyle \lbrace \hat{\ell} \rbrace}_{\scriptscriptstyle 
\hat{\varepsilon}_{1},\hat{\varepsilon}_{2},\hat{m}(\hat{\varepsilon}_{2})}g_{21},
\hat{\mathscr{F}}^{\scriptscriptstyle \lbrace \hat{\ell} \rbrace}_{\scriptscriptstyle 
\hat{\varepsilon}_{1},\hat{\varepsilon}_{2},\hat{m}(\hat{\varepsilon}_{2})}g_{22} \right),
\end{align*}
is given in equations \eqref{laxhat91}--\eqref{laxhat98} and \eqref{laxhat114}--\eqref{laxhat121} below. The orbit of 
$\mathrm{G}$ on $\mathscr{M}$ considered in this work reads:
\begin{equation*}
\mathrm{G} \mathscr{M} \! = \! \cup_{\mathfrak{g} \in \mathrm{G}} \cup_{\mathfrak{x} \in \mathscr{M}} 
\mathfrak{g} \mathfrak{x} \! = \! \cup_{\scriptscriptstyle \varepsilon_{1},\varepsilon_{2},m(\varepsilon_{2}),
\ell} \cup_{\mathfrak{x} \in \mathscr{M}} \lbrace \mathscr{F}^{\scriptscriptstyle \lbrace \ell 
\rbrace}_{\scriptscriptstyle \varepsilon_{1},\varepsilon_{2},m(\varepsilon_{2})} \mathfrak{x} 
\rbrace \bigcup \cup_{\scriptscriptstyle \hat{\varepsilon}_{1},\hat{\varepsilon}_{2},
\hat{m}(\hat{\varepsilon}_{2}),\hat{\ell}} \cup_{\mathfrak{x} \in \mathscr{M}} \lbrace 
\hat{\mathscr{F}}^{\scriptscriptstyle \lbrace \hat{\ell} \rbrace}_{\scriptscriptstyle 
\hat{\varepsilon}_{1},\hat{\varepsilon}_{2},\hat{m}(\hat{\varepsilon}_{2})} \mathfrak{x} \rbrace. 
\tag*{$\blacksquare$}
\end{equation*}
\end{eeeee}
\begin{eeeee} \label{simtation} 
Throughout this appendix, let $o$ denote ``old'' (or original) variables and let $n$ denote ``new'' (or transformed) variables, 
respectively. \hfill $\blacksquare$
\end{eeeee}
\subsection{The Transformation $\tau \! \to \! -\tau$} \label{sectonsymmt} 
Let $(\hat{u}_{o}(\tau_{o}),\hat{\varphi}_{o}(\tau_{o}))$ solve the system \eqref{eq1.5} for $\tau \! = \! \tau_{o}$, 
$\varepsilon \! = \! \varepsilon_{o} \! \in \! \lbrace \pm 1 \rbrace$, $a \! = \! a_{o}$, and $b \! = \! b_{o}$, and let 
the $4$-tuple of functions $(\hat{A}_{o}(\tau_{o}),\hat{B}_{o}(\tau_{o}),\hat{C}_{o}(\tau_{o}),\hat{D}_{o}(\tau_{o}))$, 
defined via equations \eqref{eq:ABCD} for $\hat{u}(\tau) \! = \! \hat{u}_{o}(\tau_{o})$, $\hat{\varphi}(\tau) \! = \! 
\hat{\varphi}_{o}(\tau_{o})$, $\tau \! = \! \tau_{o}$, and $\varepsilon \! = \! \varepsilon_{o}$, solve the system of 
isomonodromy deformations \eqref{eq1.4} for $\tau \! = \! \tau_{o}$ and $a \! = \! a_{o}$. Set $\hat{u}_{o}(\tau_{o}) 
\! = \! -\hat{u}_{n}(\tau_{n})$, $\hat{\varphi}_{o}(\tau_{o}) \! = \! \hat{\varphi}_{n}(\tau_{n})$, $\tau_{o} \! = \! 
\tau_{n} \me^{-\mi \pi \varepsilon_{1}}$, $\varepsilon_{1} \! \in \! \lbrace \pm 1 \rbrace$, $a_{o} \! = \! a_{n}$, 
$\varepsilon_{o} \! = \! \varepsilon_{n}$, $b_{o} \! = \! b_{n}$ (that is, $\varepsilon_{o}b_{o} \! = \! \varepsilon_{n}
b_{n}$), and $(\hat{A}_{o}(\tau_{o}),\hat{B}_{o}(\tau_{o}),\hat{C}_{o}(\tau_{o}),\hat{D}_{o}(\tau_{o})) \! = \! 
(\hat{A}_{n}(\tau_{n}),\hat{B}_{n}(\tau_{n}),-\hat{C}_{n}(\tau_{n}),-\hat{D}_{n}(\tau_{n}))$; then, $(\hat{u}_{n}
(\tau_{n}),\hat{\varphi}_{n}(\tau_{n}))$ solves the system \eqref{eq1.5} for $\tau \! = \! \tau_{n}$, $\varepsilon 
\! = \! \varepsilon_{n} \! \in \! \lbrace \pm 1 \rbrace$, $a \! = \! a_{n}$, and $b \! = \! b_{n}$, and the $4$-tuple of 
functions $(\hat{A}_{n}(\tau_{n}),\hat{B}_{n}(\tau_{n}),\hat{C}_{n}(\tau_{n}),\hat{D}_{n}(\tau_{n}))$, defined via 
equations \eqref{eq:ABCD} for $\hat{u}(\tau) \! = \! \hat{u}_{n}(\tau_{n})$, $\hat{\varphi}(\tau) \! = \! \hat{\varphi}_{n}
(\tau_{n})$, $\tau \! = \! \tau_{n}$, and $\varepsilon \! = \! \varepsilon_{n}$, solve the system~\eqref{eq1.4} for 
$\tau \! = \! \tau_{n}$, $a \! = \! a_{n}$, and $\sqrt{\smash[b]{-\hat{A}_{o}(\tau_{o}) \hat{B}_{o}(\tau_{o})}} \! = \! 
\sqrt{\smash[b]{-\hat{A}_{n}(\tau_{n}) \hat{B}_{n}(\tau_{n})}}$. Furthermore, let the functions $\hat{A}_{o}(\tau_{o})$, 
$\hat{B}_{o}(\tau_{o})$, $\hat{C}_{o}(\tau_{o})$, and $\hat{D}_{o}(\tau_{o})$ be the ones appearing in the definition 
\eqref{firstintegral} of $\hat{\alpha}(\tau)$ for $\tau \! = \! \tau_{o}$ and $a \! = \! a_{o}$, and in the first integral 
(cf. Remark \ref{alphwave}) for $\varepsilon \! = \! \varepsilon_{o} \! \in \! \lbrace \pm 1 \rbrace$ and $b \! = \! b_{o}$; 
then, under the above symmetry transformations, $\hat{\alpha}_{o}(\tau_{o}) \! = \! \hat{\alpha}_{n}(\tau_{n})$, 
where $\hat{\alpha}_{n}(\tau_{n}) \! := \! -2(\hat{B}_{n}(\tau_{n}))^{-1}(\mi a_{n} \sqrt{\smash[b]{-\hat{A}_{n}(\tau_{n}) 
\hat{B}_{n}(\tau_{n})}} \! + \! \tau_{n}(\hat{A}_{n}(\tau_{n}) \hat{D}_{n}(\tau_{n}) \! + \! \hat{B}_{n}(\tau_{n}) \hat{C}_{n}
(\tau_{n})))$, and $-\mi \hat{\alpha}_{n}(\tau_{n}) \hat{B}_{n}(\tau_{n}) \! = \! \varepsilon_{n}b_{n}$, $\varepsilon_{n} 
\! \in \! \lbrace \pm 1 \rbrace$. On the corresponding fundamental solution of the system \eqref{eqFGmain} (cf. 
equations \eqref{mpeea1} and \eqref{mpeea2}), the aforementioned transformations act as follows:
\begin{equation} \label{laxhat1} 
\mu_{o} \! = \! \mu_{n} \me^{\mi \pi l/2}, \quad l \! \in \! \lbrace \pm 1 \rbrace, \qquad \text{and} \qquad 
\widehat{\Psi}_{o}(\mu_{o},\tau_{o}) \! = \! \me^{-\frac{\mi \pi l}{4} \sigma_{3}} \widehat{\Psi}_{n}(\mu_{n},\tau_{n}).
\end{equation}

Let $(u_{o}(\tau_{o}),\varphi_{o}(\tau_{o}))$ solve the system~\eqref{equu17} for $\tau \! = \! \tau_{o}$, 
$\varepsilon \! = \! \varepsilon_{o} \! \in \! \lbrace \pm 1 \rbrace$, $a \! = \! a_{o}$, and $b \! = \! b_{o}$, and let the 
$4$-tuple of functions $(A_{o}(\tau_{o}),B_{o}(\tau_{o}),C_{o}(\tau_{o}),D_{o}(\tau_{o}))$, defined via equations 
\eqref{equu18} for $u(\tau) \! = \! u_{o}(\tau_{o})$, $\varphi (\tau) \! = \! \varphi_{o}(\tau_{o})$, $\tau \! = \! \tau_{o}$, 
and $\varepsilon \! = \! \varepsilon_{o}$, solve the corresponding system of isomonodromy deformations 
\eqref{newlax8} for $\tau \! = \! \tau_{o}$ and $a \! = \! a_{o}$. Set $u_{o}(\tau_{o}) \! = \! -u_{n}(\tau_{n})$, 
$\varphi_{o}(\tau_{o}) \! = \! \varphi_{n}(\tau_{n})$, $\tau_{o} \! = \! \tau_{n} \me^{-\mi \pi \varepsilon_{1}}$, 
$\varepsilon_{1} \! \in \! \lbrace \pm 1 \rbrace$, $a_{o} \! = \! a_{n}$, $\varepsilon_{o} \! = \! \varepsilon_{n}$, 
$b_{o} \! = \! b_{n}$ (that is, $\varepsilon_{o}b_{o} \! = \! \varepsilon_{n}b_{n}$), and $(A_{o}(\tau_{o}),B_{o}
(\tau_{o}),C_{o}(\tau_{o}),D_{o}(\tau_{o})) \! = \! (A_{n}(\tau_{n}),B_{n}(\tau_{n}),-C_{n}(\tau_{n}),-D_{n}(\tau_{n}))$; 
then, $(u_{n}(\tau_{n}),\varphi_{n}(\tau_{n}))$ solves the system \eqref{equu17} for $\tau \! = \! \tau_{n}$, 
$\varepsilon \! = \! \varepsilon_{n} \! \in \! \lbrace \pm 1 \rbrace$, $a \! = \! a_{n}$, and $b \! = \! b_{n}$, and the 
$4$-tuple of functions $(A_{n}(\tau_{n}),B_{n}(\tau_{n}),C_{n}(\tau_{n}),D_{n}(\tau_{n}))$, defined via equations 
\eqref{equu18} for $u(\tau) \! = \! u_{n}(\tau_{n})$, $\varphi (\tau) \! = \! \varphi_{n}(\tau_{n})$, $\tau \! = \! \tau_{n}$, 
and $\varepsilon \! = \! \varepsilon_{n}$, solve the system \eqref{newlax8} for $\tau \! = \! \tau_{n}$, $a \! = \! a_{n}$, 
and $\sqrt{\smash[b]{-A_{o}(\tau_{o})B_{o}(\tau_{o})}} \! = \! \sqrt{\smash[b]{-A_{n}(\tau_{n})B_{n}(\tau_{n})}}$. 
Furthermore, let the functions $A_{o}(\tau_{o})$, $B_{o}(\tau_{o})$, $C_{o}(\tau_{o})$, and $D_{o}(\tau_{o})$ be the 
ones appearing in the definition \eqref{aphnovij} of $\alpha (\tau)$ for $\tau \! = \! \tau_{o}$ and $a \! = \! a_{o}$, and 
in the first integral (cf. Remark \ref{newlax6}) for $\varepsilon \! = \! \varepsilon_{o} \! \in \! \lbrace \pm 1 \rbrace$ and 
$b \! = \! b_{o}$; then, under the above transformations, $\alpha_{o}(\tau_{o}) \! = \! \alpha_{n}(\tau_{n})$, where 
$\alpha_{n}(\tau_{n}) \! := \! -2(B_{n}(\tau_{n}))^{-1}(\mi a_{n} \sqrt{\smash[b]{-A_{n}(\tau_{n})B_{n}(\tau_{n})}} \! 
+ \! \tau_{n}(A_{n}(\tau_{n})D_{n}(\tau_{n}) \! + \! B_{n}(\tau_{n})C_{n}(\tau_{n})))$, and $-\mi \alpha_{n}(\tau_{n})
B_{n}(\tau_{n}) \! = \! \varepsilon_{n}b_{n}$, $\varepsilon_{n} \! \in \! \lbrace \pm 1 \rbrace$. On the corresponding 
fundamental solution of the system \eqref{newlax3} (cf. equations \eqref{nlxa} and \eqref{nlxb}), the aforementioned 
symmetry transformations act as follows:
\begin{equation} \label{laxhat2} 
\mu_{o} \! = \! \mu_{n} \me^{\mi \pi l/2}, \quad l \! \in \! \lbrace \pm 1 \rbrace, \qquad \text{and} \qquad 
\Psi_{o}(\mu_{o},\tau_{o}) \! = \! \me^{-\frac{\mi \pi l}{4} \sigma_{3}} \Psi_{n}(\mu_{n},\tau_{n}).
\end{equation}

In terms of the canonical solutions of the system \eqref{newlax3}, the actions \eqref{laxhat2} read: for 
$k \! \in \! \mathbb{Z}$ and $\varepsilon_{1},l \! \in \! \lbrace \pm 1 \rbrace$,
\begin{equation} \label{laxhat3} 
\mathbb{Y}_{o,k}^{\infty}(\mu_{o}) \! = \! \me^{-\frac{\mi \pi l}{4} \sigma_{3}} 
\mathbb{Y}_{n,k-l+\varepsilon_{1}}^{\infty}(\mu_{n}) \me^{\frac{\pi l a_{n}}{2} \sigma_{3}},
\end{equation}
and
\begin{equation} \label{laxhat4} 
\mathbb{X}_{o,k}^{0}(\mu_{o}) \! = \! 
\begin{cases} 
\me^{-\frac{\mi \pi l}{4} \sigma_{3}} \mathbb{X}_{n,k}^{0}(\mu_{n}), 
&\text{$\varepsilon_{1} \! = \! -l$,} \\
\mi l \me^{-\frac{\mi \pi l}{4} \sigma_{3}} \mathbb{X}_{n,k-l}^{0}(\mu_{n}) \sigma_{1}, 
&\text{$\varepsilon_{1} \! = \! l$.}
\end{cases}
\end{equation}

The transformations \eqref{laxhat3} and \eqref{laxhat4} for the canonical solutions of the system \eqref{newlax3} 
imply the following action on $\mathscr{M}$: for $k \! \in \! \mathbb{Z}$ and $\varepsilon_{1},l \! \in \! \lbrace 
\pm 1 \rbrace$,
\begin{gather}
S_{o,k}^{\infty} \! = \! \me^{-\frac{\pi l a_{n}}{2} \sigma_{3}}S_{n,k-l+\varepsilon_{1}}^{\infty} \me^{\frac{\pi l a_{n}}{2} 
\sigma_{3}}, \label{laxhat6} \\
S_{o,k}^{0} \! = \! 
\begin{cases} 
S_{n,k}^{0}, &\text{$\varepsilon_{1} \! = \! -l$,} \\
\sigma_{1}S_{n,k-l}^{0} \sigma_{1}, &\text{$\varepsilon_{1} \! = \! l$,}
\end{cases} \label{laxhat5} \\
G_{o} \! = \! 
\begin{cases}
-\mi S_{n,0}^{0} \sigma_{1}G_{n} \me^{\frac{\pi a_{n}}{2} \sigma_{3}}, 
&\text{$\varepsilon_{1} \! = \! 1$,} \\
\mi \sigma_{1}(S_{n,0}^{0})^{-1}G_{n} \me^{-\frac{\pi a_{n}}{2} \sigma_{3}}, 
&\text{$\varepsilon_{1} \! = \! -1$.} \label{laxhat7}
\end{cases}
\end{gather}

The actions \eqref{laxhat6}--\eqref{laxhat7} on $\mathscr{M}$ can be expressed in terms of an intermediate auxiliary 
mapping $\mathscr{F}_{\scriptscriptstyle \mathscr{M}}^{\, \rightslice}(\varepsilon_{1}) \colon \mathbb{C}^{8} \! \to \! 
\mathbb{C}^{8}$, $\varepsilon_{1} \! \in \! \lbrace \pm 1 \rbrace$, which is an isomorphism on $\mathscr{M}$:
\begin{align*} 
\mathscr{F}_{\scriptscriptstyle \mathscr{M}}^{\, \rightslice}(\varepsilon_{1}) \colon 
\mathscr{M} \! \to \! \mathscr{M}, \, &(a,s_{0}^{0},s_{0}^{\infty},s_{1}^{\infty},
g_{11},g_{12},g_{21},g_{22}) \! \mapsto \! \left(a,s_{0}^{0}(\varepsilon_{1}),
s_{0}^{\infty}(\varepsilon_{1}),s_{1}^{\infty}(\varepsilon_{1}), \right. \\
&\left. \, g_{11}(\varepsilon_{1}),g_{12}(\varepsilon_{1}),g_{21}(\varepsilon_{1}),
g_{22}(\varepsilon_{1}) \right),
\end{align*}
where, for $\varepsilon_{1} \! = \! -1$,
\begin{equation} \label{laxhat8} 
\begin{gathered}
s_{0}^{0}(-1) \! = \! s_{0}^{0}, \quad s_{0}^{\infty}(-1) \! = \! s_{0}^{\infty} 
\me^{\pi a}, \quad s_{1}^{\infty}(-1) \! = \! s_{1}^{\infty} \me^{-\pi a}, \quad 
g_{11}(-1) \! = \! -\mi (g_{21} \! + \! s_{0}^{0}g_{11}) \me^{\pi a/2}, \\
g_{12}(-1) \! = \! -\mi (g_{22} \! + \! s_{0}^{0}g_{12}) \me^{-\pi a/2}, \quad 
g_{21}(-1) \! = \! -\mi g_{11} \me^{\pi a/2}, \quad g_{22}(-1) \! = \! -\mi 
g_{12} \me^{-\pi a/2},
\end{gathered}
\end{equation}
and, for $\varepsilon_{1} \! = \! 1$,
\begin{equation} \label{laxhat9}
\begin{gathered}
s_{0}^{0}(1) \! = \! s_{0}^{0}, \quad s_{0}^{\infty}(1) \! = \! s_{0}^{\infty} 
\me^{-\pi a}, \quad s_{1}^{\infty}(1) \! = \! s_{1}^{\infty} \me^{\pi a}, 
\quad g_{11}(1) \! = \! \mi g_{21} \me^{-\pi a/2}, \\
g_{12}(1) \! = \! \mi g_{22} \me^{\pi a/2}, \quad g_{21}(1) \! = \! \mi 
(g_{11} \! - \! s_{0}^{0}g_{21}) \me^{-\pi a/2}, \quad g_{22}(1) \! = \! \mi 
(g_{12} \! - \! s_{0}^{0}g_{22}) \me^{\pi a/2}.
\end{gathered}
\end{equation}
One uses this transformation in order to arrive at asymptotics for $\tau \! < \! 0$ by using those for 
$\tau \! > \! 0$.\footnote{In Section 7, p. 45 of \cite{avlkv}, it is stated that the Lie-point symmetry $\tau \! \to \! -\tau$ 
in Subsection 6.2.1 of \cite{a1} requires correction. Keeping in mind the $\operatorname{mod}(2 \pi)$ arbitrariness 
inherent in the definition of the function $\hat{\varphi}(\tau)$ discussed in the Introduction to this appendix, the 
Lie-point symmetry $\tau \! \to \! -\tau$ alluded to in Section 7, p. 45 of \cite{avlkv} is the one for which the ``additive 
parameter'', denoted by $\hat{\varphi}_{0}$, is equal to zero: the transformation changing $\tau \! \to \! -\tau$ for which 
$\hat{\varphi}_{0} \! = \! 0$ is presented \textbf{here}, in Appendix \ref{sectonsymmt}, and \textbf{not} in Subsection 
6.2.1 of \cite{a1} wherein the Transformation 6.2.1 changing $\tau \! \to \! -\tau$ was derived under the condition 
$\hat{\varphi}_{o}(\tau_{o}) \! \to \! \hat{\varphi}_{o}(\tau_{o}) \! - \! \pi \epsilon_{1}^{\ast} \! =: \! \hat{\varphi}_{n}
(\tau_{n})$, $\epsilon_{1}^{\ast} \! \in \! \lbrace \pm 1 \rbrace$, that is, the additive parameter is equal to 
$-\pi \epsilon_{1}^{\ast}$ (unfortunately, the action of the symmetry $\tau \! \to \! -\tau$ on the function 
$\hat{\varphi}(\tau)$ was not emphasized in \cite{a1}).}
\subsection{The Transformation $\tau \! \to \! \tau$} \label{sectonsymmt2} 
Let $(\hat{u}_{o}(\tau_{o}),\hat{\varphi}_{o}(\tau_{o}))$ solve the system \eqref{eq1.5} for $\tau \! = \! \tau_{o}$, 
$\varepsilon \! = \! \varepsilon_{o} \! \in \! \lbrace \pm 1 \rbrace$, $a \! = \! a_{o}$, and $b \! = \! b_{o}$, and let the 
$4$-tuple of functions $(\hat{A}_{o}(\tau_{o}),\hat{B}_{o}(\tau_{o}),\hat{C}_{o}(\tau_{o}),\hat{D}_{o}(\tau_{o}))$, 
defined via equations \eqref{eq:ABCD} for $\hat{u}(\tau) \! = \! \hat{u}_{o}(\tau_{o})$, $\hat{\varphi}(\tau) \! = \! 
\hat{\varphi}_{o}(\tau_{o})$, $\tau \! = \! \tau_{o}$, and $\varepsilon \! = \! \varepsilon_{o}$, solve the system of 
isomonodromy deformations \eqref{eq1.4} for $\tau \! = \! \tau_{o}$ and $a \! = \! a_{o}$. Set $\hat{u}_{o}(\tau_{o}) 
\! = \! -\hat{u}_{n}(\tau_{n})$, $\hat{\varphi}_{o}(\tau_{o}) \! = \! \hat{\varphi}_{n}(\tau_{n})$, $\tau_{o} \! = \! 
\tau_{n}$, $a_{o} \! = \! a_{n}$, $\varepsilon_{o} \! = \! -\varepsilon_{n}$, $b_{o} \! = \! -b_{n}$ (that is, 
$\varepsilon_{o}b_{o} \! = \! \varepsilon_{n}b_{n}$), and $(\hat{A}_{o}(\tau_{o}),\hat{B}_{o}(\tau_{o}),\hat{C}_{o}
(\tau_{o}),\hat{D}_{o}(\tau_{o})) \! = \! (-\hat{A}_{n}(\tau_{n}),-\hat{B}_{n}(\tau_{n}),-\hat{C}_{n}(\tau_{n}),-\hat{D}_{n}
(\tau_{n}))$; then, $(\hat{u}_{n}(\tau_{n}),\hat{\varphi}_{n}(\tau_{n}))$ solves the system \eqref{eq1.5} for $\tau \! = \! 
\tau_{n}$, $\varepsilon \! = \! \varepsilon_{n} \! \in \! \lbrace \pm 1 \rbrace$, $a \! = \! a_{n}$, and $b \! = \! b_{n}$, and 
the $4$-tuple of functions $(\hat{A}_{n}(\tau_{n}),\hat{B}_{n}(\tau_{n}),\hat{C}_{n}(\tau_{n}),\hat{D}_{n}(\tau_{n}))$, 
defined via equations \eqref{eq:ABCD} for $\hat{u}(\tau) \! = \! \hat{u}_{n}(\tau_{n})$, $\hat{\varphi}(\tau) \! = \! 
\hat{\varphi}_{n}(\tau_{n})$, $\tau \! = \! \tau_{n}$, and $\varepsilon \! = \! \varepsilon_{n}$, solve the system 
\eqref{eq1.4} for $\tau \! = \! \tau_{n}$, $a \! = \! a_{n}$, and $\sqrt{\smash[b]{-\hat{A}_{o}(\tau_{o}) \hat{B}_{o}(\tau_{o})}} 
\! = \! \sqrt{\smash[b]{-\hat{A}_{n}(\tau_{n}) \hat{B}_{n}(\tau_{n})}}$. Moreover, let the functions $\hat{A}_{o}(\tau_{o})$, 
$\hat{B}_{o}(\tau_{o})$, $\hat{C}_{o}(\tau_{o})$, and $\hat{D}_{o}(\tau_{o})$ be the ones appearing in the definition 
\eqref{firstintegral} of $\hat{\alpha}(\tau)$ for $\tau \! = \! \tau_{o}$ and $a \! = \! a_{o}$, and in the first integral (cf. 
Remark \ref{alphwave}) for $\varepsilon \! = \! \varepsilon_{o} \! \in \! \lbrace \pm 1 \rbrace$ and $b \! = \! b_{o}$; then, 
under the above transformations, $\hat{\alpha}_{o}(\tau_{o}) \! = \! -\hat{\alpha}_{n}(\tau_{n})$, where $\hat{\alpha}_{n}
(\tau_{n}) \! := \! -2(\hat{B}_{n}(\tau_{n}))^{-1}(\mi a_{n} \sqrt{\smash[b]{-\hat{A}_{n}(\tau_{n}) \hat{B}_{n}(\tau_{n})}} \! 
+ \! \tau_{n}(\hat{A}_{n}(\tau_{n}) \hat{D}_{n}(\tau_{n}) \! + \! \hat{B}_{n}(\tau_{n}) \hat{C}_{n}(\tau_{n})))$, and 
$-\mi \hat{\alpha}_{n}(\tau_{n}) \hat{B}_{n}(\tau_{n}) \! = \! \varepsilon_{n}b_{n}$, $\varepsilon_{n} \! \in \! \lbrace \pm 1 
\rbrace$. On the corresponding fundamental solution of the system \eqref{eqFGmain} (cf. equations \eqref{mpeea1} 
and \eqref{mpeea2}), the aforementioned symmetry transformations act as follows:
\begin{equation} \label{laxhat10} 
\mu_{o} \! = \! \mu_{n} \me^{\mi \pi m}, \quad m \! \in \! \lbrace 0,1 \rbrace, \qquad \text{and} \qquad \widehat{\Psi}_{o}
(\mu_{o},\tau_{o}) \! = \! \me^{\frac{\mi \pi}{2}(m-1) \sigma_{3}} \widehat{\Psi}_{n}(\mu_{n},\tau_{n}).
\end{equation}

Let $(u_{o}(\tau_{o}),\varphi_{o}(\tau_{o}))$ solve the system~\eqref{equu17} for $\tau \! = \! \tau_{o}$, $\varepsilon \! 
= \! \varepsilon_{o} \! \in \! \lbrace \pm 1 \rbrace$, $a \! = \! a_{o}$, and $b \! = \! b_{o}$, and let the $4$-tuple of functions 
$(A_{o}(\tau_{o}),B_{o}(\tau_{o}),C_{o}(\tau_{o}),D_{o}(\tau_{o}))$, defined via equations \eqref{equu18} for $u(\tau) \! 
= \! u_{o}(\tau_{o})$, $\varphi (\tau) \! = \! \varphi_{o}(\tau_{o})$, $\tau \! = \! \tau_{o}$, and $\varepsilon \! = \! \varepsilon_{o}$, 
solve the corresponding system of isomonodromy deformations \eqref{newlax8} for $\tau \! = \! \tau_{o}$ and $a \! = \! 
a_{o}$. Set $u_{o}(\tau_{o}) \! = \! -u_{n}(\tau_{n})$, $\varphi_{o}(\tau_{o}) \! = \! \varphi_{n}(\tau_{n})$, $\tau_{o} \! = \! 
\tau_{n}$, $a_{o} \! = \! a_{n}$, $\varepsilon_{o} \! = \! -\varepsilon_{n}$, $b_{o} \! = \! -b_{n}$ (that is, $\varepsilon_{o}
b_{o} \! = \! \varepsilon_{n}b_{n}$), and $(A_{o}(\tau_{o}),B_{o}(\tau_{o}),C_{o}(\tau_{o}),D_{o}(\tau_{o})) \! = \! (-A_{n}
(\tau_{n}),-B_{n}(\tau_{n}),-C_{n}(\tau_{n}),-D_{n}(\tau_{n}))$; then, $(u_{n}(\tau_{n}),\varphi_{n}(\tau_{n}))$ solves the 
system \eqref{equu17} for $\tau \! = \! \tau_{n}$, $\varepsilon \! = \! \varepsilon_{n} \! \in \! \lbrace \pm 1 \rbrace$, $a \! 
= \! a_{n}$, and $b \! = \! b_{n}$, and the $4$-tuple of functions $(A_{n}(\tau_{n}),B_{n}(\tau_{n}),C_{n}(\tau_{n}),D_{n}
(\tau_{n}))$, defined via equations \eqref{equu18} for $u(\tau) \! = \! u_{n}(\tau_{n})$, $\varphi (\tau) \! = \! \varphi_{n}
(\tau_{n})$, $\tau \! = \! \tau_{n}$, and $\varepsilon \! = \! \varepsilon_{n}$, solve the system \eqref{newlax8} for $\tau \! 
= \! \tau_{n}$, $a \! = \! a_{n}$, and $\sqrt{\smash[b]{-A_{o}(\tau_{o})B_{o}(\tau_{o})}} \! = \! \sqrt{\smash[b]{-A_{n}(\tau_{n})
B_{n}(\tau_{n})}}$. Furthermore, let the functions $A_{o}(\tau_{o})$, $B_{o}(\tau_{o})$, $C_{o}(\tau_{o})$, and $D_{o}
(\tau_{o})$ be the ones appearing in the definition \eqref{aphnovij} of $\alpha (\tau)$ for $\tau \! = \! \tau_{o}$ and $a \! 
= \! a_{o}$, and in the first integral (cf. Remark \ref{newlax6}) for $\varepsilon \! = \! \varepsilon_{o} \! \in \! \lbrace \pm 1 
\rbrace$ and $b \! = \! b_{o}$; then, under the above transformations, $\alpha_{o}(\tau_{o}) \! = \! -\alpha_{n}(\tau_{n})$, 
where $\alpha_{n}(\tau_{n}) \! := \! -2(B_{n}(\tau_{n}))^{-1}(\mi a_{n} \sqrt{\smash[b]{-A_{n}(\tau_{n})B_{n}(\tau_{n})}} \! 
+ \! \tau_{n}(A_{n}(\tau_{n})D_{n}(\tau_{n}) \! + \! B_{n}(\tau_{n})C_{n}(\tau_{n})))$, and $-\mi \alpha_{n}(\tau_{n})B_{n}
(\tau_{n}) \! = \! \varepsilon_{n}b_{n}$, $\varepsilon_{n} \! \in \! \lbrace \pm 1 \rbrace$. On the corresponding fundamental 
solution of the system \eqref{newlax3} (cf. equations \eqref{nlxa} and \eqref{nlxb}), the aforementioned symmetry 
transformations act as follows:
\begin{equation} \label{laxhat11} 
\mu_{o} \! = \! \mu_{n} \me^{\mi \pi m}, \quad m \! \in \! \lbrace 0,1 \rbrace, \qquad \text{and} \qquad \Psi_{o}
(\mu_{o},\tau_{o}) \! = \! \me^{\frac{\mi \pi}{2}(m-1) \sigma_{3}} \Psi_{n}(\mu_{n},\tau_{n}).
\end{equation}

In terms of the canonical solutions of the system~\eqref{newlax3}, the actions \eqref{laxhat11} read: for $k \! \in \! 
\mathbb{Z}$, $m \! \in \! \lbrace 0,1 \rbrace$, and $\tilde{l} \! \in \! \lbrace \pm 1 \rbrace$,\footnote{As discussed 
in Remarks \ref{newrem12} and \ref{newrem13}, since the canonical solutions $\mathbb{X}^{0}_{k}(\mu)$, $k \! \in 
\! \mathbb{Z}$, are defined uniquely provided the branch of $(B(\tau))^{1/2}$ is fixed, it follows that, since the branch 
of $(B(\tau))^{1/2}$ is not fixed, the canonical solutions $\mathbb{X}^{0}_{k}(\mu)$, $k \! \in \! \mathbb{Z}$, are 
defined up to a sign (plus or minus), thus the appearance of the `sign parameter' $\tilde{l}$: this comment applies, 
\emph{mutatis mutandis}, throughout the remaining sub-appendices.}
\begin{equation} \label{laxhat12} 
\mathbb{Y}_{o,k}^{\infty}(\mu_{o}) \! = \! \me^{\frac{\mi \pi}{2}(m-1) \sigma_{3}} 
\mathbb{Y}_{n,k-2m}^{\infty}(\mu_{n}) \me^{-\frac{\mi \pi}{2}(m-1) \sigma_{3}} 
\me^{\pi m (a_{n}-\mi/2) \sigma_{3}},
\end{equation}
and
\begin{equation} \label{laxhat13} 
\mathbb{X}_{o,k}^{0}(\mu_{o}) \! = \! 
\begin{cases} 
-\tilde{l} \me^{-\frac{\mi \pi}{2} \sigma_{3}} \mathbb{X}_{n,k}^{0}(\mu_{n}), 
&\text{$m \! = \! 0$,} \\
\mi \tilde{l} \, \mathbb{X}_{n,k-1}^{0}(\mu_{n}) \sigma_{1}, &\text{$m \! = \! 1$.}
\end{cases}
\end{equation}

The transformations \eqref{laxhat12} and~\eqref{laxhat13} for the canonical solutions of the system \eqref{newlax3} 
imply the following action on $\mathscr{M}$: for $k \! \in \! \mathbb{Z}$, $m \! \in \! \lbrace 0,1 \rbrace$, and $\tilde{l} 
\! \in \! \lbrace \pm 1 \rbrace$,
\begin{gather}
S_{o,k}^{\infty} \! = \! \me^{\frac{\mi \pi}{2}(m-1) \sigma_{3}} \me^{-\pi m (a_{n}-\mi/2) \sigma_{3}}S_{n,k-2m}^{\infty} 
\me^{\pi m (a_{n}-\mi/2) \sigma_{3}} \me^{-\frac{\mi \pi}{2}(m-1) \sigma_{3}}, \label{laxhat15} \\
S_{o,k}^{0} \! = \! 
\begin{cases} 
S_{n,k}^{0}, &\text{$m \! = \! 0$,} \\
\sigma_{1}S_{n,k-1}^{0} \sigma_{1}, &\text{$m \! = \! 1$,}
\end{cases} \label{laxhat14} \\
G_{o} \! = \! -\tilde{l}G_{n} \me^{\frac{\mi \pi}{2} \sigma_{3}}. \label{laxhat16}
\end{gather}

The actions \eqref{laxhat15}--\eqref{laxhat16} on $\mathscr{M}$ can be expressed in terms of an intermediate auxiliary 
mapping $\mathscr{F}_{\scriptscriptstyle \mathscr{M}}^{\, \leftslice}(\tilde{l}) \colon \mathbb{C}^{8} \! \to \! \mathbb{C}^{8}$, 
$\tilde{l} \! \in \! \lbrace \pm 1 \rbrace$, which is an isomorphism on $\mathscr{M}$:
\begin{align*} 
\mathscr{F}_{\scriptscriptstyle \mathscr{M}}^{\, \leftslice}(\tilde{l}) \colon \mathscr{M} \! \to \! 
\mathscr{M}, \, &(a,s_{0}^{0},s_{0}^{\infty},s_{1}^{\infty},g_{11},g_{12},g_{21},
g_{22}) \! \mapsto \! \left(a,s_{0}^{0}(\tilde{l}),s_{0}^{\infty}(\tilde{l}),s_{1}^{\infty}
(\tilde{l}),\right. \\
&\left. \, g_{11}(\tilde{l}),g_{12}(\tilde{l}),g_{21}(\tilde{l}),g_{22}(\tilde{l}) \right),
\end{align*}
where
\begin{equation} \label{laxhat17} 
\begin{gathered}
s_{0}^{0}(\tilde{l}) \! = \! s_{0}^{0}, \quad s_{0}^{\infty}(\tilde{l}) \! = \! -s_{0}^{\infty}, 
\quad s_{1}^{\infty}(\tilde{l}) \! = \! -s_{1}^{\infty}, \quad g_{11}(\tilde{l}) 
\! = \! \mi \tilde{l}g_{11}, \quad g_{12}(\tilde{l}) \! = \! -\mi \tilde{l}g_{12}, \\
g_{21}(\tilde{l}) \! = \! \mi \tilde{l}g_{21}, \quad g_{22}(\tilde{l}) \! = \! -\mi \tilde{l}g_{22}.
\end{gathered}
\end{equation}
One uses this transformation in order to define an analogue of the identity map; see, in particular, Appendix 
\ref{sectonsymmcomp}, definitions~\eqref{laxhat52} and \eqref{laxhat53}.
\subsection{The Transformation $a \! \to \! -a$} \label{sectonsymmatoa}
Let $(\hat{u}_{o}(\tau_{o}),\hat{\varphi}_{o}(\tau_{o}))$ solve the system \eqref{eq1.5} 
for $\tau \! = \! \tau_{o}$, $\varepsilon \! = \! \varepsilon_{o} \! \in \! \lbrace 
\pm 1 \rbrace$, $a \! = \! a_{o}$, and $b \! = \! b_{o}$, and let the $4$-tuple 
of functions $(\hat{A}_{o}(\tau_{o}),\hat{B}_{o}(\tau_{o}),\hat{C}_{o}(\tau_{o}),
\hat{D}_{o}(\tau_{o}))$, defined via equations \eqref{eq:ABCD} for $\hat{u}(\tau) 
\! = \! \hat{u}_{o}(\tau_{o})$, $\hat{\varphi}(\tau) \! = \! \hat{\varphi}_{o}(\tau_{o})$, 
$\tau \! = \! \tau_{o}$, and $\varepsilon \! = \! \varepsilon_{o}$, solve the system 
of isomonodromy deformations~\eqref{eq1.4} for $\tau \! = \! \tau_{o}$ and 
$a \! = \! a_{o}$. Set $\hat{u}_{o}(\tau_{o}) \! = \! -\hat{u}_{n}(\tau_{n})$, 
$\hat{\varphi}_{o}(\tau_{o}) \! = \! -\hat{\varphi}_{n}(\tau_{n})$, $\tau_{o} \! 
= \! \tau_{n}$, $a_{o} \! = \! -a_{n}$, $\varepsilon_{o} \! = \! \varepsilon_{n} 
\me^{-\mi \pi \varepsilon_{2}}$, $\varepsilon_{2} \! \in \! \lbrace \pm 1 \rbrace$, 
$b_{o} \! = \! b_{n}$ (that is, $\varepsilon_{o}b_{o} \! = \! \varepsilon_{n}b_{n} 
\me^{-\mi \pi \varepsilon_{2}}$), and $(\hat{A}_{o}(\tau_{o}),\hat{B}_{o}(\tau_{o}),
\hat{C}_{o}(\tau_{o}),\hat{D}_{o}(\tau_{o})) \! = \! (\hat{B}_{n}(\tau_{n}),\hat{A}_{n}
(\tau_{n}),-\hat{D}_{n}(\tau_{n}),-\hat{C}_{n}(\tau_{n}))$; then, $(\hat{u}_{n}(\tau_{n}),
\hat{\varphi}_{n}(\tau_{n}))$ solves the system \eqref{eq1.5} for $\tau \! = \! \tau_{n}$, 
$\varepsilon \! = \! \varepsilon_{n} \! \in \! \lbrace \pm 1 \rbrace$, $a \! = \! a_{n}$, 
and $b \! = \! b_{n}$, and the $4$-tuple of functions $(\hat{A}_{n}(\tau_{n}),\hat{B}_{n}
(\tau_{n}),\hat{C}_{n}(\tau_{n}),\hat{D}_{n}(\tau_{n}))$, defined via equations \eqref{eq:ABCD} 
for $\hat{u}(\tau) \! = \! \hat{u}_{n}(\tau_{n})$, $\hat{\varphi}(\tau) \! = \! \hat{\varphi}_{n}
(\tau_{n})$, $\tau \! = \! \tau_{n}$, and $\varepsilon \! = \! \varepsilon_{n}$, solve 
the system \eqref{eq1.4} for $\tau \! = \! \tau_{n}$, $a \! = \! a_{n}$, and 
$\sqrt{\smash[b]{-\hat{A}_{o}(\tau_{o}) \hat{B}_{o}(\tau_{o})}} \! = \! 
\sqrt{\smash[b]{-\hat{A}_{n}(\tau_{n}) \hat{B}_{n}(\tau_{n})}}$. Furthermore, let the 
functions $\hat{A}_{o}(\tau_{o})$, $\hat{B}_{o}(\tau_{o})$, $\hat{C}_{o}(\tau_{o})$, 
and $\hat{D}_{o}(\tau_{o})$ be the ones appearing in the definition \eqref{firstintegral} 
of $\hat{\alpha}(\tau)$ for $\tau \! = \! \tau_{o}$ and $a \! = \! a_{o}$, and in the 
first integral (cf. Remark \ref{alphwave}) for $\varepsilon \! = \! \varepsilon_{o} \! 
\in \! \lbrace \pm 1 \rbrace$ and $b \! = \! b_{o}$; then, under the above symmetry 
transformations, $\hat{\alpha}_{o}(\tau_{o}) \! = \! -\hat{B}_{n}(\tau_{n})(\hat{A}_{n}
(\tau_{n}))^{-1} \hat{\alpha}_{n}(\tau_{n})$, where $\hat{\alpha}_{n}(\tau_{n}) \! := \! 
-2(\hat{B}_{n}(\tau_{n}))^{-1}(\mi a_{n} \sqrt{\smash[b]{-\hat{A}_{n}(\tau_{n}) \hat{B}_{n}
(\tau_{n})}} \! + \! \tau_{n}(\hat{A}_{n}(\tau_{n}) \hat{D}_{n}(\tau_{n}) \! + \! \hat{B}_{n}
(\tau_{n}) \hat{C}_{n}(\tau_{n})))$, and $-\mi \hat{\alpha}_{n}(\tau_{n}) \hat{B}_{n}
(\tau_{n}) \! = \! \varepsilon_{n}b_{n}$, $\varepsilon_{n} \! \in \! \lbrace \pm 1 \rbrace$. 
On the corresponding fundamental solution of the system \eqref{eqFGmain} (cf. 
equations \eqref{mpeea1} and \eqref{mpeea2}), the aforementioned transformations 
act as follows:
\begin{equation} \label{laxhat18} 
\mu_{o} \! = \! \mu_{n} \me^{\mi \pi m/2}, \quad m \! \in \! \lbrace \pm 1 \rbrace, \quad \text{and} 
\quad \widehat{\Psi}_{o}(\mu_{o},\tau_{o}) \! = \! \widehat{\mathcal{Q}}(\mu_{n},\tau_{n}) 
\widehat{\Psi}_{n}(\mu_{n},\tau_{n}),
\end{equation}
where
\begin{equation} \label{laxhat19} 
\widehat{\mathcal{Q}}(\mu_{n},\tau_{n}) \! := \! \left(\frac{\hat{B}_{n}(\tau_{n}) \me^{-\mi \pi m/4}}{
\sqrt{\smash[b]{-\hat{A}_{n}(\tau_{n}) \hat{B}_{n}(\tau_{n})}}} \right)^{\sigma_{3}} \! + \! \mu_{n} 
\me^{\mi \pi m/4} \sigma_{-}.
\end{equation}

Let $(u_{o}(\tau_{o}),\varphi_{o}(\tau_{o}))$ solve the system \eqref{equu17} for 
$\tau \! = \! \tau_{o}$, $\varepsilon \! = \! \varepsilon_{o} \! \in \! \lbrace \pm 1 
\rbrace$, $a \! = \! a_{o}$, and $b \! = \! b_{o}$, and let the $4$-tuple of functions 
$(A_{o}(\tau_{o}),B_{o}(\tau_{o}),C_{o}(\tau_{o}),D_{o}(\tau_{o}))$, defined via 
equations \eqref{equu18} for $u(\tau) \! = \! u_{o}(\tau_{o})$, $\varphi (\tau) \! = \! 
\varphi_{o}(\tau_{o})$, $\tau \! = \! \tau_{o}$, and $\varepsilon \! = \! \varepsilon_{o}$, 
solve the corresponding system of isomonodromy deformations \eqref{newlax8} for 
$\tau \! = \! \tau_{o}$ and $a \! = \! a_{o}$. Set $u_{o}(\tau_{o}) \! = \! -u_{n}
(\tau_{n})$, $\varphi_{o}(\tau_{o}) \! = \! -\varphi_{n}(\tau_{n})$, $\tau_{o} \! = \! 
\tau_{n}$, $a_{o} \! = \! -a_{n}$, $\varepsilon_{o} \! = \! \varepsilon_{n} \me^{-\mi \pi 
\varepsilon_{2}}$, $\varepsilon_{2} \! \in \! \lbrace \pm 1 \rbrace$, $b_{o} \! = \! b_{n}$ 
(that is, $\varepsilon_{o}b_{o} \! = \! \varepsilon_{n}b_{n} \me^{-\mi \pi \varepsilon_{2}}$), 
and $(A_{o}(\tau_{o}),B_{o}(\tau_{o}),C_{o}(\tau_{o}),D_{o}(\tau_{o})) \! = \! (B_{n}(\tau_{n}),
A_{n}(\tau_{n}),-D_{n}(\tau_{n}),-C_{n}(\tau_{n}))$; then, $(u_{n}(\tau_{n}),\varphi_{n}
(\tau_{n}))$ solves the system \eqref{equu17} for $\tau \! = \! \tau_{n}$, $\varepsilon \! = 
\! \varepsilon_{n} \! \in \! \lbrace \pm 1 \rbrace$, $a \! = \! a_{n}$, and $b \! = \! b_{n}$, 
and the $4$-tuple of functions $(A_{n}(\tau_{n}),B_{n}(\tau_{n}),C_{n}(\tau_{n}),D_{n}
(\tau_{n}))$, defined via equations \eqref{equu18} for $u(\tau) \! = \! u_{n}(\tau_{n})$, 
$\varphi (\tau) \! = \! \varphi_{n}(\tau_{n})$, $\tau \! = \! \tau_{n}$, and $\varepsilon 
\! = \! \varepsilon_{n}$, solve the system \eqref{newlax8} for $\tau \! = \! \tau_{n}$, 
$a \! = \! a_{n}$, and $\sqrt{\smash[b]{-A_{o}(\tau_{o})B_{o}(\tau_{o})}} \! = \! 
\sqrt{\smash[b]{-A_{n}(\tau_{n})B_{n}(\tau_{n})}}$. Furthermore, let the functions $A_{o}
(\tau_{o})$, $B_{o}(\tau_{o})$, $C_{o}(\tau_{o})$, and $D_{o}(\tau_{o})$ be the ones 
appearing in the definition \eqref{aphnovij} of $\alpha (\tau)$ for $\tau \! = \! \tau_{o}$ 
and $a \! = \! a_{o}$, and in the first integral (cf. Remark \ref{newlax6}) for $\varepsilon 
\! = \! \varepsilon_{o} \! \in \! \lbrace \pm 1 \rbrace$ and $b \! = \! b_{o}$; then, under the 
above transformations, $\alpha_{o}(\tau_{o}) \! = \! -B_{n}(\tau_{n})(A_{n}(\tau_{n}))^{-1} 
\alpha_{n}(\tau_{n})$, where $\alpha_{n}(\tau_{n}) \! := \! -2(B_{n}(\tau_{n}))^{-1}(\mi a_{n} 
\sqrt{\smash[b]{-A_{n}(\tau_{n})B_{n}(\tau_{n})}} \! + \! \tau_{n}(A_{n}(\tau_{n})D_{n}(\tau_{n}) 
\! + \! B_{n}(\tau_{n})C_{n}(\tau_{n})))$, and $-\mi \alpha_{n}(\tau_{n})B_{n}(\tau_{n}) 
\! = \! \varepsilon_{n}b_{n}$, $\varepsilon_{n} \! \in \! \lbrace \pm 1 \rbrace$. On the corresponding 
fundamental solution of the system \eqref{newlax3} (cf. equations \eqref{nlxa} and \eqref{nlxb}), 
the aforementioned symmetry transformations act as follows:
\begin{equation} \label{laxhat20} 
\mu_{o} \! = \! \mu_{n} \me^{\mi \pi m/2}, \quad m \! \in \! \lbrace \pm 1 \rbrace, \quad \text{and} 
\quad \Psi_{o}(\mu_{o},\tau_{o}) \! = \! \mathcal{Q}(\mu_{n},\tau_{n}) \Psi_{n}(\mu_{n},\tau_{n}),
\end{equation}
where
\begin{equation} \label{laxhat21} 
\mathcal{Q}(\mu_{n},\tau_{n}) \! := \! \left(\frac{B_{n}(\tau_{n}) \me^{-\mi \pi m/4}}{
\sqrt{\smash[b]{-A_{n}(\tau_{n})B_{n}(\tau_{n})}}} \right)^{\sigma_{3}} \! + \! \mu_{n} 
\me^{\mi \pi m/4} \sigma_{-}.
\end{equation}

In terms of the canonical solutions of the system \eqref{newlax3}, the actions \eqref{laxhat20} read: for 
$k \! \in \! \mathbb{Z}$ and $m,\varepsilon_{2},l \! \in \! \lbrace \pm 1 \rbrace$,
\begin{equation} \label{laxhat22} 
\mathbb{Y}_{o,k}^{\infty}(\mu_{o}) \! = \! \mathcal{Q}(\mu_{n},\tau_{n}) 
\mathbb{Y}_{n,k-m}^{\infty}(\mu_{n}) \me^{\frac{\pi m a_{n}}{2} \sigma_{3}} 
\sigma_{3} \sigma_{1},
\end{equation}
and
\begin{equation} \label{laxhat23} 
\mathbb{X}_{o,k}^{0}(\mu_{o}) \! = \! 
\begin{cases} 
l \mathcal{Q}(\mu_{n},\tau_{n}) \mathbb{X}_{n,k}^{0}(\mu_{n}), &\text{$m \! = \! -\varepsilon_{2}$,} \\
\mi l \mathcal{Q}(\mu_{n},\tau_{n}) \mathbb{X}_{n,k-m}^{0}(\mu_{n}) \sigma_{1}, 
&\text{$m \! = \! \varepsilon_{2}$.}
\end{cases}
\end{equation}

The transformations \eqref{laxhat22} and \eqref{laxhat23} for the canonical solutions of the system 
\eqref{newlax3} imply the following action on $\mathscr{M}$: for $k \! \in \! \mathbb{Z}$ and 
$m,\varepsilon_{2},l \! \in \! \lbrace \pm 1 \rbrace$,
\begin{gather}
S_{o,k}^{\infty} \! = \! \sigma_{1} \sigma_{3} \me^{-\frac{\pi m a_{n}}{2} \sigma_{3}}
S_{n,k-m}^{\infty} \me^{\frac{\pi m a_{n}}{2} \sigma_{3}} \sigma_{3} \sigma_{1}, 
\label{laxhat25} \\
S_{o,k}^{0} \! = \! 
\begin{cases} 
S_{n,k}^{0}, &\text{$m \! = \! -\varepsilon_{2}$,} \\
\sigma_{1}S_{n,k-m}^{0} \sigma_{1}, &\text{$m \! = \! \varepsilon_{2}$,}
\end{cases} \label{laxhat24} \\
G_{o} \! = \! 
\begin{cases}
-\mi l S_{o,0}^{0} \sigma_{1}G_{n} \me^{\pi (a_{n}-\mi/2) \sigma_{3}} \sigma_{3}
(S_{n,1}^{\infty})^{-1} \sigma_{3} \me^{-\pi (a_{n}-\mi/2) \sigma_{3}} 
\me^{\frac{\pi a_{n}}{2} \sigma_{3}} \sigma_{3} \sigma_{1}, &\text{$(m,
\varepsilon_{2}) \! = \! (1,1)$,} \\
lG_{n} \me^{\pi (a_{n}-\mi/2) \sigma_{3}} \sigma_{3}(S_{n,1}^{\infty})^{-1} \sigma_{3} 
\me^{-\pi (a_{n}-\mi/2) \sigma_{3}} \me^{\frac{\pi a_{n}}{2} \sigma_{3}} \sigma_{3} 
\sigma_{1}, &\text{$(m,\varepsilon_{2}) \! = \! (1,-1)$,} \\
lG_{n}S_{n,0}^{\infty} \me^{-\frac{\pi a_{n}}{2} \sigma_{3}} \sigma_{3} \sigma_{1}, 
&\text{$(m,\varepsilon_{2}) \! = \! (-1,1)$,} \\
-\mi l \sigma_{1}(S_{o,0}^{0})^{-1}G_{n}S_{n,0}^{\infty} \me^{-\frac{\pi a_{n}}{2} 
\sigma_{3}} \sigma_{3} \sigma_{1}, &\text{$(m,\varepsilon_{2}) \! = \! (-1,-1)$.} 
\label{laxhat26}
\end{cases}
\end{gather}

The actions \eqref{laxhat25}--\eqref{laxhat26} on $\mathscr{M}$ can be 
expressed in terms of an intermediate auxiliary mapping 
$\mathscr{F}_{\scriptscriptstyle \mathscr{M}}^{\, \leftrightarrowtriangle}(m,\varepsilon_{2}) \colon 
\mathbb{C}^{8} \! \to \! \mathbb{C}^{8}$, $m,\varepsilon_{2} \! \in \! \lbrace \pm 1 \rbrace$, which 
is an isomorphism on $\mathscr{M}$: for $l \! \in \! \lbrace \pm 1 \rbrace$,
\begin{align*} 
\mathscr{F}_{\scriptscriptstyle \mathscr{M}}^{\, \leftrightarrowtriangle}
(m,\varepsilon_{2}) \colon \mathscr{M} \! \to \! \mathscr{M}, \, &(a,s_{0}^{0},
s_{0}^{\infty},s_{1}^{\infty},g_{11},g_{12},g_{21},g_{22}) \! \mapsto \! \left(
-a,s_{0}^{0}(m,\varepsilon_{2}),s_{0}^{\infty}(m,\varepsilon_{2}),s_{1}^{\infty}
(m,\varepsilon_{2}),\right. \\
&\left. \, g_{11}(m,\varepsilon_{2}),g_{12}(m,\varepsilon_{2}),g_{21}
(m,\varepsilon_{2}),g_{22}(m,\varepsilon_{2}) \right),
\end{align*}
where, for $(m,\varepsilon_{2}) \! = \! (1,1)$,
\begin{equation} \label{laxhat27} 
\begin{gathered}
s_{0}^{0}(1,1) \! = \! s_{0}^{0}, \quad s_{0}^{\infty}(1,1) \! = \! -s_{1}^{\infty} 
\me^{\pi a}, \quad s_{1}^{\infty}(1,1) \! = \! -s_{0}^{\infty} \me^{\pi a}, 
\quad g_{11}(1,1) \! = \! \mi lg_{22} \me^{\pi a/2}, \\
g_{12}(1,1) \! = \! -\mi 
l(g_{21} \! + \! s_{0}^{\infty}g_{22}) \me^{-\pi a/2}, \quad g_{21}(1,1) 
\! = \! \mi l(g_{12} \! - \! s_{0}^{0}g_{22}) \me^{\pi a/2}, \\
g_{22}(1,1) \! = \! \mi l(-g_{11} \! - \! s_{0}^{\infty}g_{12} \! + \! 
s_{0}^{0}(g_{21} \! + \! s_{0}^{\infty}g_{22})) \me^{-\pi a/2},
\end{gathered}
\end{equation}
for $(m,\varepsilon_{2}) \! = \! (1,-1)$,
\begin{equation} \label{laxhat28} 
\begin{gathered}
s_{0}^{0}(1,-1) \! = \! s_{0}^{0}, \quad s_{0}^{\infty}(1,-1) \! = \! -s_{1}^{\infty} 
\me^{\pi a}, \quad s_{1}^{\infty}(1,-1) \! = \! -s_{0}^{\infty} \me^{\pi a}, \\
g_{11}(1,-1) \! = \! lg_{12} \me^{\pi a/2}, \quad g_{12}(1,-1) \! = \! -l(g_{11} 
\! + \! s_{0}^{\infty}g_{12}) \me^{-\pi a/2}, \\
g_{21}(1,-1) \! = \! lg_{22} \me^{\pi a/2}, \quad g_{22}(1,-1) \! = \! -l(g_{21} 
\! + \! s_{0}^{\infty}g_{22}) \me^{-\pi a/2},
\end{gathered}
\end{equation}
for $(m,\varepsilon_{2}) \! = \! (-1,1)$,
\begin{equation} \label{laxhat29} 
\begin{gathered}
s_{0}^{0}(-1,1) \! = \! s_{0}^{0}, \quad s_{0}^{\infty}(-1,1) \! = \! -s_{1}^{\infty} 
\me^{\pi a}, \quad s_{1}^{\infty}(-1,1) \! = \! -s_{0}^{\infty} \me^{\pi a}, \\
g_{11}(-1,1) \! = \! l(g_{12} \! - \! s_{1}^{\infty}g_{11} \me^{2 \pi a}) 
\me^{-\pi a/2}, \quad g_{12}(-1,1) \! = \! -lg_{11} \me^{\pi a/2}, \\
g_{21}(-1,1) \! = \! l(g_{22} \! - \! s_{1}^{\infty}g_{21} \me^{2 \pi a}) 
\me^{-\pi a/2}, \quad g_{22}(-1,1) \! = \! -lg_{21} \me^{\pi a/2},
\end{gathered}
\end{equation}
and, for $(m,\varepsilon_{2}) \! = \! (-1,-1)$,
\begin{equation} \label{laxhat30} 
\begin{gathered}
s_{0}^{0}(-1,-1) \! = \! s_{0}^{0}, \quad s_{0}^{\infty}(-1,-1) \! = \! -s_{1}^{\infty} 
\me^{\pi a}, \quad s_{1}^{\infty}(-1,-1) \! = \! -s_{0}^{\infty} \me^{\pi a}, \\
g_{11}(-1,-1) \! = \! \mi l(g_{22} \! - \! s_{1}^{\infty}g_{21} \me^{2 \pi a} \! + \! 
s_{0}^{0}(g_{12} \! - \! s_{1}^{\infty}g_{11} \me^{2 \pi a})) \me^{-\pi a/2}, \\
g_{12}(-1,-1) \! = \! -\mi l(g_{21} \! + \! s_{0}^{0}g_{11}) \me^{\pi a/2}, \quad 
g_{21}(-1,-1) \! = \! \mi l(g_{12} \! - \! s_{1}^{\infty}g_{11} \me^{2 \pi a}) 
\me^{-\pi a/2}, \\
g_{22}(-1,-1) \! = \! -\mi lg_{11} \me^{\pi a/2}.
\end{gathered}
\end{equation}
One uses this transformation in order to arrive at asymptotics for $\varepsilon b \! < \! 0$ by using 
those for $\varepsilon b \! > \! 0$.
\subsection{The Transformation $\tau \! \to \! \pm \mi \tau$} \label{sectonsymmtit} 
Let $(\hat{u}_{o}(\tau_{o}),\hat{\varphi}_{o}(\tau_{o}))$ solve the system \eqref{eq1.5} 
for $\tau \! = \! \tau_{o}$, $\varepsilon \! = \! \varepsilon_{o} \! \in \! \lbrace 
\pm 1 \rbrace$, $a \! = \! a_{o}$, and $b \! = \! b_{o}$, and let the $4$-tuple of 
functions $(\hat{A}_{o}(\tau_{o}),\hat{B}_{o}(\tau_{o}),\hat{C}_{o}(\tau_{o}),\hat{D}_{o}
(\tau_{o}))$, defined via equations \eqref{eq:ABCD} for $\hat{u}(\tau) \! = \! 
\hat{u}_{o}(\tau_{o})$, $\hat{\varphi}(\tau) \! = \! \hat{\varphi}_{o}(\tau_{o})$, 
$\tau \! = \! \tau_{o}$, and $\varepsilon \! = \! \varepsilon_{o}$, solve the system 
of isomonodromy deformations \eqref{eq1.4} for $\tau \! = \! \tau_{o}$ and 
$a \! = \! a_{o}$. Set $\hat{u}_{o}(\tau_{o}) \! = \! \hat{u}_{n}(\tau_{n}) 
\me^{\mi \pi \tilde{\varepsilon}_{1}/2}$, $\tilde{\varepsilon}_{1} \! \in \! \lbrace 
\pm 1 \rbrace$, $\hat{\varphi}_{o}(\tau_{o}) \! = \! \hat{\varphi}_{n}(\tau_{n})$, 
$\tau_{o} \! = \! \tau_{n} \me^{-\mi \pi \tilde{\varepsilon}_{1}/2}$, $a_{o} \! 
= \! a_{n}$, $\varepsilon_{o} \! = \! \varepsilon_{n}$, and $b_{o} \! = \! b_{n} 
\me^{-\mi \pi \tilde{\varepsilon}_{2}}$, $\tilde{\varepsilon}_{2} \! \in \! \lbrace 
\pm 1 \rbrace$ (that is, $\varepsilon_{o}b_{o} \! = \! \varepsilon_{n}b_{n} 
\me^{-\mi \pi \tilde{\varepsilon}_{2}}$), and $(\hat{A}_{o}(\tau_{o}),\hat{B}_{o}
(\tau_{o}),\hat{C}_{o}(\tau_{o}),\hat{D}_{o}(\tau_{o})) \! = \! (\hat{A}_{n}(\tau_{n}) 
\me^{\mi \pi \tilde{\varepsilon}_{1}},\hat{B}_{n}(\tau_{n}) \me^{\mi \pi 
\tilde{\varepsilon}_{1}},\hat{C}_{n}(\tau_{n}) \me^{\mi \pi \tilde{\varepsilon}_{1}/2},
\hat{D}_{n}(\tau_{n}) \me^{\mi \pi \tilde{\varepsilon}_{1}/2})$; then, $(\hat{u}_{n}
(\tau_{n}),\hat{\varphi}_{n}(\tau_{n}))$ solves the system \eqref{eq1.5} for $\tau \! 
= \! \tau_{n}$, $\varepsilon \! = \! \varepsilon_{n} \! \in \! \lbrace \pm 1 \rbrace$, 
$a \! = \! a_{n}$, and $b \! = \! b_{n}$, and the $4$-tuple of functions $(\hat{A}_{n}
(\tau_{n}),\hat{B}_{n}(\tau_{n}),\hat{C}_{n}(\tau_{n}),\hat{D}_{n}(\tau_{n}))$, defined 
via equations~\eqref{eq:ABCD} for $\hat{u}(\tau) \! = \! \hat{u}_{n}(\tau_{n})$, 
$\hat{\varphi}(\tau) \! = \! \hat{\varphi}_{n}(\tau_{n})$, $\tau \! = \! \tau_{n}$, 
and $\varepsilon \! = \! \varepsilon_{n}$, solve the system \eqref{eq1.4} for 
$\tau \! = \! \tau_{n}$, $a \! = \! a_{n}$, and $\sqrt{\smash[b]{-\hat{A}_{o}
(\tau_{o}) \hat{B}_{o}(\tau_{o})}} \! = \! \me^{\mi \pi \tilde{\varepsilon}_{1}} 
\sqrt{\smash[b]{-\hat{A}_{n}(\tau_{n}) \hat{B}_{n}(\tau_{n})}}$. Moreover, let the 
functions $\hat{A}_{o}(\tau_{o})$, $\hat{B}_{o}(\tau_{o})$, $\hat{C}_{o}(\tau_{o})$, and 
$\hat{D}_{o}(\tau_{o})$ be the ones appearing in the definition \eqref{firstintegral} 
of $\hat{\alpha}(\tau)$ for $\tau \! = \! \tau_{o}$ and $a \! = \! a_{o}$, and in the 
first integral (cf. Remark \ref{alphwave}) for $\varepsilon \! = \! \varepsilon_{o} \! 
\in \! \lbrace \pm 1 \rbrace$ and $b \! = \! b_{o}$; then, under the above symmetry 
transformations, $\hat{\alpha}_{o}(\tau_{o}) \! = \! \hat{\alpha}_{n}(\tau_{n})$, 
where $\hat{\alpha}_{n}(\tau_{n}) \! := \! -2(\hat{B}_{n}(\tau_{n}))^{-1}(\mi a_{n} 
\sqrt{\smash[b]{-\hat{A}_{n}(\tau_{n}) \hat{B}_{n}(\tau_{n})}} \! + \! \tau_{n}
(\hat{A}_{n}(\tau_{n}) \hat{D}_{n}(\tau_{n}) \! + \! \hat{B}_{n}(\tau_{n}) \hat{C}_{n}
(\tau_{n})))$, and $-\mi \hat{\alpha}_{n}(\tau_{n}) \hat{B}_{n}(\tau_{n}) \! = \! 
\varepsilon_{n}b_{n}$, $\varepsilon_{n} \! \in \! \lbrace \pm 1 \rbrace$. On 
the corresponding fundamental solution of the system~\eqref{eqFGmain} (cf. equations 
\eqref{mpeea1} and \eqref{mpeea2}), the aforementioned transformations act as follows:
\begin{equation} \label{maga1} 
\mu_{o} \! = \! \mu_{n} \me^{\mi \pi \tilde{\varepsilon}_{1}/4}, \quad \tilde{\varepsilon}_{1} \! \in \! 
\lbrace \pm 1 \rbrace, \qquad \text{and} \qquad \widehat{\Psi}_{o}(\mu_{o},\tau_{o}) \! = \! 
\me^{-\frac{\mi \pi \tilde{\varepsilon}_{1}}{8} \sigma_{3}} \widehat{\Psi}_{n}(\mu_{n},\tau_{n}).
\end{equation}

Let $(u_{o}(\tau_{o}),\varphi_{o}(\tau_{o}))$ solve the system \eqref{equu17} for 
$\tau \! = \! \tau_{o}$, $\varepsilon \! = \! \varepsilon_{o} \! \in \! \lbrace \pm 1 
\rbrace$, $a \! = \! a_{o}$, and $b \! = \! b_{o}$, and let the $4$-tuple of functions 
$(A_{o}(\tau_{o}),B_{o}(\tau_{o}),C_{o}(\tau_{o}),D_{o}(\tau_{o}))$, defined via 
equations \eqref{equu18} for $u(\tau) \! = \! u_{o}(\tau_{o})$, $\varphi (\tau) \! = \! 
\varphi_{o}(\tau_{o})$, $\tau \! = \! \tau_{o}$, and $\varepsilon \! = \! \varepsilon_{o}$, 
solve the corresponding system of isomonodromy deformations \eqref{newlax8} 
for $\tau \! = \! \tau_{o}$ and $a \! = \! a_{o}$. Set $u_{o}(\tau_{o}) \! = \! u_{n}
(\tau_{n}) \me^{\mi \pi \tilde{\varepsilon}_{1}/2}$, $\tilde{\varepsilon}_{1} \! \in 
\! \lbrace \pm 1 \rbrace$, $\varphi_{o}(\tau_{o}) \! = \! \varphi_{n}(\tau_{n})$, 
$\tau_{o} \! = \! \tau_{n} \me^{-\mi \pi \tilde{\varepsilon}_{1}/2}$, $a_{o} \! 
= \! a_{n}$, $\varepsilon_{o} \! = \! \varepsilon_{n}$, and $b_{o} \! = \! b_{n} 
\me^{-\mi \pi \tilde{\varepsilon}_{2}}$, $\tilde{\varepsilon}_{2} \! \in \! \lbrace 
\pm 1 \rbrace$ (that is, $\varepsilon_{o}b_{o} \! = \! \varepsilon_{n}b_{n} 
\me^{-\mi \pi \tilde{\varepsilon}_{2}}$), and $(A_{o}(\tau_{o}),B_{o}(\tau_{o}),
C_{o}(\tau_{o}),D_{o}(\tau_{o})) \! = \! (A_{n}(\tau_{n}) \me^{\mi \pi 
\tilde{\varepsilon}_{1}},B_{n}(\tau_{n}) \me^{\mi \pi \tilde{\varepsilon}_{1}},
C_{n}(\tau_{n}) \me^{\mi \pi \tilde{\varepsilon}_{1}/2},D_{n}(\tau_{n}) 
\me^{\mi \pi \tilde{\varepsilon}_{1}/2})$; then, $(u_{n}(\tau_{n}),\varphi_{n}(\tau_{n}))$ 
solves the system \eqref{equu17} for $\tau \! = \! \tau_{n}$, $\varepsilon \! = \! 
\varepsilon_{n} \! \in \! \lbrace \pm 1 \rbrace$, $a \! = \! a_{n}$, and $b \! = \! b_{n}$, 
and the $4$-tuple of functions $(A_{n}(\tau_{n}),B_{n}(\tau_{n}),C_{n}(\tau_{n}),D_{n}
(\tau_{n}))$, defined via equations \eqref{equu18} for $u(\tau) \! = \! u_{n}(\tau_{n})$, 
$\varphi (\tau) \! = \! \varphi_{n}(\tau_{n})$, $\tau \! = \! \tau_{n}$, and $\varepsilon 
\! = \! \varepsilon_{n}$, solve the system \eqref{newlax8} for $\tau \! = \! \tau_{n}$, 
$a \! = \! a_{n}$, and $\sqrt{\smash[b]{-A_{o}(\tau_{o})B_{o}(\tau_{o})}} \! = \! 
\me^{\mi \pi \tilde{\varepsilon}_{1}} \sqrt{\smash[b]{-A_{n}(\tau_{n})B_{n}(\tau_{n})}}$. 
Furthermore, let the functions $A_{o}(\tau_{o})$, $B_{o}(\tau_{o})$, $C_{o}(\tau_{o})$, 
and $D_{o}(\tau_{o})$ be the ones appearing in the definition \eqref{aphnovij} of 
$\alpha (\tau)$ for $\tau \! = \! \tau_{o}$ and $a \! = \! a_{o}$, and in the first 
integral (cf. Remark \ref{newlax6}) for $\varepsilon \! = \! \varepsilon_{o} \! \in \! 
\lbrace \pm 1 \rbrace$ and $b \! = \! b_{o}$; then, under the above transformations, 
$\alpha_{o}(\tau_{o}) \! = \! \alpha_{n}(\tau_{n})$, where $\alpha_{n}(\tau_{n}) \! := \! 
-2(B_{n}(\tau_{n}))^{-1}(\mi a_{n} \sqrt{\smash[b]{-A_{n}(\tau_{n})B_{n}(\tau_{n})}} \! + 
\! \tau_{n}(A_{n}(\tau_{n})D_{n}(\tau_{n}) \! + \! B_{n}(\tau_{n})C_{n}(\tau_{n})))$, and 
$-\mi \alpha_{n}(\tau_{n})B_{n}(\tau_{n}) \! = \! \varepsilon_{n}b_{n}$, $\varepsilon_{n} 
\! \in \! \lbrace \pm 1 \rbrace$. On the corresponding fundamental solution of the system 
\eqref{newlax3} (cf. equations \eqref{nlxa} and \eqref{nlxb}), the aforementioned symmetry 
transformations act as follows:
\begin{equation} \label{maga2} 
\mu_{o} \! = \! \mu_{n} \me^{\mi \pi \tilde{\varepsilon}_{1}/4}, \quad \tilde{\varepsilon}_{1} \! 
\in \! \lbrace \pm 1 \rbrace, \qquad \text{and} \qquad \Psi_{o}(\mu_{o},\tau_{o}) \! = \! 
\me^{-\frac{\mi \pi \tilde{\varepsilon}_{1}}{8} \sigma_{3}} \Psi_{n}(\mu_{n},\tau_{n}).
\end{equation}

In terms of the canonical solutions of the system \eqref{newlax3}, the actions \eqref{maga2} read: for 
$k \! \in \! \mathbb{Z}$ and $\tilde{\varepsilon}_{1},\tilde{\varepsilon}_{2} \! \in \! \lbrace \pm 1 \rbrace$,
\begin{equation} \label{maga3} 
\mathbb{Y}_{o,k}^{\infty}(\mu_{o}) \! = \! \me^{-\frac{\mi \pi \tilde{\varepsilon}_{1}}{8} \sigma_{3}} 
\mathbb{Y}_{n,k}^{\infty}(\mu_{n}) \me^{\frac{a_{n} \pi \tilde{\varepsilon}_{1}}{4} \sigma_{3}},
\end{equation}
and
\begin{equation} \label{maga4} 
\mathbb{X}_{o,k}^{0}(\mu_{o}) \! = \! 
\begin{cases} 
\me^{-\frac{\mi \pi \tilde{\varepsilon}_{1}}{8} \sigma_{3}} \mathbb{X}_{n,k}^{0}
(\mu_{n}), &\text{$\tilde{\varepsilon}_{1} \! = \! -\tilde{\varepsilon}_{2}$,} \\
-\mi \tilde{\varepsilon}_{1} \me^{-\frac{\mi \pi \tilde{\varepsilon}_{1}}{8} 
\sigma_{3}} \mathbb{X}_{n,k-\tilde{\varepsilon}_{1}}^{0}(\mu_{n}) \sigma_{1}, 
&\text{$\tilde{\varepsilon}_{1} \! = \! \tilde{\varepsilon}_{2}$.}
\end{cases}
\end{equation}

The transformations \eqref{maga3} and \eqref{maga4} for the canonical solutions of the system 
\eqref{newlax3} imply the following action on $\mathscr{M}$: for $k \! \in \! \mathbb{Z}$ and 
$\tilde{\varepsilon}_{1},\tilde{\varepsilon}_{2} \! \in \! \lbrace \pm 1 \rbrace$,
\begin{gather}
S_{o,k}^{\infty} \! = \! \me^{-\frac{a_{n} \pi \tilde{\varepsilon}_{1}}{4} \sigma_{3}}
S_{n,k}^{\infty} \me^{\frac{a_{n} \pi \tilde{\varepsilon}_{1}}{4} \sigma_{3}}, 
\label{maga6} \\
S_{o,k}^{0} \! = \! 
\begin{cases} 
S_{n,k}^{0}, &\text{$\tilde{\varepsilon}_{1} \! = \! -\tilde{\varepsilon}_{2}$,} \\
\sigma_{1}S_{n,k-\tilde{\varepsilon}_{1}}^{0} \sigma_{1}, &\text{$\tilde{\varepsilon}_{1} 
\! = \! \tilde{\varepsilon}_{2}$,}
\end{cases} \label{maga5} \\
G_{o} \! = \! 
\begin{cases} 
\mi S_{o,0}^{0} \sigma_{1} G_{n} \me^{\frac{a_{n} \pi}{4} \sigma_{3}}, 
&\text{$(\tilde{\varepsilon}_{1},\tilde{\varepsilon}_{2}) \! = \! (1,1)$,} \\
G_{n} \me^{\frac{a_{n} \pi}{4} \sigma_{3}}, &\text{$(\tilde{\varepsilon}_{1},
\tilde{\varepsilon}_{2}) \! = \! (1,-1)$,} \\
G_{n} \me^{-\frac{a_{n} \pi}{4} \sigma_{3}}, &\text{$(\tilde{\varepsilon}_{1},
\tilde{\varepsilon}_{2}) \! = \! (-1,1)$,} \\
-\mi \sigma_{1}(S_{o,0}^{0})^{-1}G_{n} \me^{-\frac{a_{n} \pi}{4} \sigma_{3}}, 
&\text{$(\tilde{\varepsilon}_{1},\tilde{\varepsilon}_{2}) \! = \! (-1,-1)$.} \label{maga7}
\end{cases}
\end{gather}

The actions \eqref{maga6}--\eqref{maga7} on $\mathscr{M}$ can be expressed in terms of an intermediate 
auxiliary mapping $\mathscr{F}_{\scriptscriptstyle \mathscr{M}}^{\, \leftrightsquigarrow}(\tilde{\varepsilon}_{1},
\tilde{\varepsilon}_{2}) \colon \mathbb{C}^{8} \! \to \! \mathbb{C}^{8}$, $\tilde{\varepsilon}_{1},
\tilde{\varepsilon}_{2}\! \in \! \lbrace \pm 1 \rbrace$, which is an isomorphism on $\mathscr{M}$:
\begin{align*} 
\mathscr{F}_{\scriptscriptstyle \mathscr{M}}^{\, \leftrightsquigarrow}(\tilde{\varepsilon}_{1},
\tilde{\varepsilon}_{2}) \colon \mathscr{M} \! \to \! \mathscr{M}, \, &(a,s_{0}^{0},
s_{0}^{\infty},s_{1}^{\infty},g_{11},g_{12},g_{21},g_{22}) \! \mapsto \! \left(a,s_{0}^{0}
(\tilde{\varepsilon}_{1},\tilde{\varepsilon}_{2}),s_{0}^{\infty}(\tilde{\varepsilon}_{1},
\tilde{\varepsilon}_{2}),s_{1}^{\infty}(\tilde{\varepsilon}_{1},\tilde{\varepsilon}_{2}), 
\right. \\
&\left. \, g_{11}(\tilde{\varepsilon}_{1},\tilde{\varepsilon}_{2}),g_{12}
(\tilde{\varepsilon}_{1},\tilde{\varepsilon}_{2}),g_{21}(\tilde{\varepsilon}_{1},
\tilde{\varepsilon}_{2}),g_{22}(\tilde{\varepsilon}_{1},\tilde{\varepsilon}_{2}) \right),
\end{align*}
where, for $(\tilde{\varepsilon}_{1},\tilde{\varepsilon}_{2}) \! = \! (1,1)$,
\begin{equation} \label{laxhat31} 
\begin{gathered}
s_{0}^{0}(1,1) \! = \! s_{0}^{0}, \quad s_{0}^{\infty}(1,1) \! = \! s_{0}^{\infty} 
\me^{-\pi a/2}, \quad s_{1}^{\infty}(1,1) \! = \! s_{1}^{\infty} \me^{\pi a/2}, \\
g_{11}(1,1) \! = \! -\mi g_{21} \me^{-\pi a/4}, \quad g_{12}(1,1) \! = \! 
-\mi g_{22} \me^{\pi a/4}, \\
g_{21}(1,1) \! = \! -\mi (g_{11} \! - \! s_{0}^{0}g_{21}) \me^{-\pi a/4}, \quad 
g_{22}(1,1) \! = \! -\mi (g_{12} \! - \! s_{0}^{0}g_{22}) \me^{\pi a/4},
\end{gathered}
\end{equation}
for $(\tilde{\varepsilon}_{1},\tilde{\varepsilon}_{2}) \! = \! (1,-1)$,
\begin{equation} \label{laxhat32} 
\begin{gathered}
s_{0}^{0}(1,-1) \! = \! s_{0}^{0}, \quad s_{0}^{\infty}(1,-1) \! = \! s_{0}^{\infty} 
\me^{-\pi a/2}, \quad s_{1}^{\infty}(1,-1) \! = \! s_{1}^{\infty} \me^{\pi a/2}, \\
g_{11}(1,-1) \! = \! g_{11} \me^{-\pi a/4}, \quad g_{12}(1,-1) \! = \! g_{12} 
\me^{\pi a/4}, \quad g_{21}(1,-1) \! = \! g_{21} \me^{-\pi a/4}, \\
g_{22}(1,-1) \! = \! g_{22} \me^{\pi a/4},
\end{gathered}
\end{equation}
for $(\tilde{\varepsilon}_{1},\tilde{\varepsilon}_{2}) \! = \! (-1,1)$,
\begin{equation} \label{laxhat33} 
\begin{gathered}
s_{0}^{0}(-1,1) \! = \! s_{0}^{0}, \quad s_{0}^{\infty}(-1,1) \! = \! s_{0}^{\infty} 
\me^{\pi a/2}, \quad s_{1}^{\infty}(-1,1) \! = \! s_{1}^{\infty} \me^{-\pi a/2}, \\
g_{11}(-1,1) \! = \! g_{11} \me^{\pi a/4}, \quad g_{12}(-1,1) \! = \! g_{12} 
\me^{-\pi a/4}, \quad g_{21}(-1,1) \! = \! g_{21} \me^{\pi a/4}, \\
g_{22}(-1,1) \! = \! g_{22} \me^{-\pi a/4},
\end{gathered}
\end{equation}
and, for $(\tilde{\varepsilon}_{1},\tilde{\varepsilon}_{2}) \! = \! (-1,-1)$,
\begin{equation} \label{laxhat34} 
\begin{gathered}
s_{0}^{0}(-1,-1) \! = \! s_{0}^{0}, \quad s_{0}^{\infty}(-1,-1) \! = \! s_{0}^{\infty} 
\me^{\pi a/2}, \quad s_{1}^{\infty}(-1,-1) \! = \! s_{1}^{\infty} \me^{-\pi a/2}, \\
g_{11}(-1,-1) \! = \! \mi (g_{21} \! + \! s_{0}^{0}g_{11}) \me^{\pi a/4}, \quad 
g_{12}(-1,-1) \! = \! \mi (g_{22} \! + \! s_{0}^{0}g_{12}) \me^{-\pi a/4}, \\
g_{21}(-1,-1) \! = \! \mi g_{11} \me^{\pi a/4}, \quad g_{22}(-1,-1) \! = \! 
\mi g_{12} \me^{-\pi a/4}.
\end{gathered}
\end{equation}
One uses this transformation in order to arrive at asymptotics for pure-imaginary $\tau$ by using those 
for real $\tau$.
\subsection{Composed Symmetries and Asymptotics} \label{sectonsymmcomp}
In order to derive the complete set of requisite transformations, one considers the actions~\eqref{laxhat8}, 
\eqref{laxhat9}, \eqref{laxhat17}, \eqref{laxhat27}--\eqref{laxhat30}, and \eqref{laxhat31}--\eqref{laxhat34} 
as a group of basis symmetries, the compositions of whose elements yield the remaining isomorphisms 
on $\mathscr{M}$.

In order to do so, however, additional notation is necessary. For symmetries related to real $\tau$, introduce 
the auxiliary parameters $\varepsilon_{1} \! \in \! \lbrace 0,\pm 1 \rbrace$, $\varepsilon_{2} \! \in \! \lbrace 0,\pm 1 
\rbrace$, $m(\varepsilon_{2}) \! = \! 
\left\{
\begin{smallmatrix}
0, \, \, \varepsilon_{2}=0, \\
\pm \varepsilon_{2}, \, \, \varepsilon_{2} \in \lbrace \pm 1 \rbrace,
\end{smallmatrix}
\right.$ and $\ell \! \in \! \lbrace 0,1 \rbrace$, and consider the $4$-tuple 
$(\varepsilon_{1},\varepsilon_{2},m(\varepsilon_{2}) \vert \ell)$ concomitant with 
its associated isomorphism(s) on $\mathscr{M}$ denoted by $\mathscr{F}^{
\scriptscriptstyle \lbrace \ell \rbrace}_{\scriptscriptstyle \varepsilon_{1},
\varepsilon_{2},m(\varepsilon_{2})} \colon \mathbb{C}^{8} \! \to \! \mathbb{C}^{8}$, where
\begin{align} \label{laxhat35} 
\mathscr{F}^{\lbrace \ell \rbrace}_{\varepsilon_{1},\varepsilon_{2},m(\varepsilon_{2})} 
\colon \mathscr{M} \! \to \! \mathscr{M}, \, &(a,s_{0}^{0},s_{0}^{\infty},s_{1}^{\infty},
g_{11},g_{12},g_{21},g_{22}) \! \mapsto \! \left((-1)^{\varepsilon_{2}}a,s_{0}^{0}
(\varepsilon_{1},\varepsilon_{2},m(\varepsilon_{2}) \vert \ell),\right. \nonumber \\
&\left. \, s_{0}^{\infty}(\varepsilon_{1},\varepsilon_{2},m(\varepsilon_{2}) \vert \ell),
s_{1}^{\infty}(\varepsilon_{1},\varepsilon_{2},m(\varepsilon_{2}) \vert \ell),g_{11}
(\varepsilon_{1},\varepsilon_{2},m(\varepsilon_{2}) \vert \ell), \right. \nonumber \\
&\left. \, g_{12}(\varepsilon_{1},\varepsilon_{2},m(\varepsilon_{2}) \vert \ell),g_{21}
(\varepsilon_{1},\varepsilon_{2},m(\varepsilon_{2}) \vert \ell),g_{22}(\varepsilon_{1},
\varepsilon_{2},m(\varepsilon_{2}) \vert \ell) \right);
\end{align}
and, for symmetries related to pure-imaginary $\tau$, introduce the auxiliary parameters $\hat{\varepsilon}_{1} 
\! \in \! \lbrace \pm 1 \rbrace$, $\hat{\varepsilon}_{2} \! \in \! \lbrace 0,\pm 1 \rbrace$, $\hat{m}(\hat{\varepsilon}_{2}) 
\! = \! 
\left\{
\begin{smallmatrix}
0, \, \, \hat{\varepsilon}_{2} \in \lbrace \pm 1 \rbrace, \\
\pm \hat{\varepsilon}_{1}, \, \, \hat{\varepsilon}_{2}=0,
\end{smallmatrix}
\right.$ and $\hat{\ell} \! \in \! \lbrace 0,1 \rbrace$, and consider the $4$-tuple 
$(\hat{\varepsilon}_{1},\hat{\varepsilon}_{2},\hat{m}(\hat{\varepsilon}_{2}) \vert 
\hat{\ell})$ concomitant with its associated isomorphism(s) on $\mathscr{M}$ denoted 
by $\hat{\mathscr{F}}^{\scriptscriptstyle \lbrace \hat{\ell} \rbrace}_{\scriptscriptstyle 
\hat{\varepsilon}_{1},\hat{\varepsilon}_{2},\hat{m}(\hat{\varepsilon}_{2})} \colon 
\mathbb{C}^{8} \! \to \! \mathbb{C}^{8}$, where
\begin{align} \label{laxhat36} 
\hat{\mathscr{F}}^{\lbrace \hat{\ell} \rbrace}_{\hat{\varepsilon}_{1},\hat{\varepsilon}_{2},
\hat{m}(\hat{\varepsilon}_{2})} \colon \mathscr{M} \! \to \! \mathscr{M}, \, &(a,
s_{0}^{0},s_{0}^{\infty},s_{1}^{\infty},g_{11},g_{12},g_{21},g_{22}) \! \mapsto \! \left(
(-1)^{1+\hat{\varepsilon}_{2}}a,\hat{s}_{0}^{0}(\hat{\varepsilon}_{1},\hat{\varepsilon}_{2},
\hat{m}(\hat{\varepsilon}_{2}) \vert \hat{\ell}),\right. \nonumber \\
&\left. \, \hat{s}_{0}^{\infty}(\hat{\varepsilon}_{1},\hat{\varepsilon}_{2},\hat{m}
(\hat{\varepsilon}_{2}) \vert \hat{\ell}),\hat{s}_{1}^{\infty}(\hat{\varepsilon}_{1},
\hat{\varepsilon}_{2},\hat{m}(\hat{\varepsilon}_{2}) \vert \hat{\ell}),\hat{g}_{11}
(\hat{\varepsilon}_{1},\hat{\varepsilon}_{2},\hat{m}(\hat{\varepsilon}_{2}) \vert 
\hat{\ell}),\right. \nonumber \\
&\left. \, \hat{g}_{12}(\hat{\varepsilon}_{1},\hat{\varepsilon}_{2},\hat{m}
(\hat{\varepsilon}_{2}) \vert \hat{\ell}),\hat{g}_{21}(\hat{\varepsilon}_{1},
\hat{\varepsilon}_{2},\hat{m}(\hat{\varepsilon}_{2}) \vert \hat{\ell}),\hat{g}_{22}
(\hat{\varepsilon}_{1},\hat{\varepsilon}_{2},\hat{m}(\hat{\varepsilon}_{2}) \vert 
\hat{\ell}) \right).
\end{align}
Let
\begin{equation} \label{laxhat37} 
\mathscr{F}^{\lbrace 0 \rbrace}_{0,0,0} \colon \mathscr{M} \! \to \! \mathscr{M}, \, 
(a,s_{0}^{0},s_{0}^{\infty},s_{1}^{\infty},g_{11},g_{12},g_{21},g_{22}) \! \mapsto \! 
(a,s_{0}^{0},s_{0}^{\infty},s_{1}^{\infty},g_{11},g_{12},g_{21},g_{22})
\end{equation}
denote the \emph{identity map},\footnote{That is, $s_{0}^{0}(0,0,0 \vert 0) \! = \! 
s_{0}^{0}$, $s_{0}^{\infty}(0,0,0 \vert 0) \! = \! s_{0}^{\infty}$, $s_{1}^{\infty}(0,0,0 
\vert 0) \! = \! s_{1}^{\infty}$, and $g_{ij}(0,0,0 \vert 0) \! = \! g_{ij}$, $i,j \! \in \! 
\lbrace 1,2 \rbrace$.} and, for $\ell \! = \! 0$, set
\begin{equation} \label{laxhat38} 
\mathscr{F}^{\lbrace 0 \rbrace}_{\varepsilon_{1},\varepsilon_{2},m(\varepsilon_{2})} 
\! := \! 
\begin{cases}
\mathscr{F}_{\scriptscriptstyle \mathscr{M}}^{\, \rightslice}(1), &\text{$(\varepsilon_{1},
\varepsilon_{2},m(\varepsilon_{2}) \vert \ell) \! = \! (1,0,0 \vert 0)$,} \\
\mathscr{F}_{\scriptscriptstyle \mathscr{M}}^{\, \rightslice}(-1), &\text{$(\varepsilon_{1},
\varepsilon_{2},m(\varepsilon_{2}) \vert \ell) \! = \! (-1,0,0 \vert 0)$,} \\
\mathscr{F}_{\scriptscriptstyle \mathscr{M}}^{\, \leftrightarrowtriangle}(1,1), &\text{$
(\varepsilon_{1},\varepsilon_{2},m(\varepsilon_{2}) \vert \ell) \! = \! (0,1,1 \vert 0)$,} \\
\mathscr{F}_{\scriptscriptstyle \mathscr{M}}^{\, \leftrightarrowtriangle}(1,-1), &\text{$
(\varepsilon_{1},\varepsilon_{2},m(\varepsilon_{2}) \vert \ell) \! = \! (0,-1,1 \vert 0)$,} \\
\mathscr{F}_{\scriptscriptstyle \mathscr{M}}^{\, \leftrightarrowtriangle}(-1,1), &\text{$
(\varepsilon_{1},\varepsilon_{2},m(\varepsilon_{2}) \vert \ell) \! = \! (0,1,-1 \vert 0)$,} \\
\mathscr{F}_{\scriptscriptstyle \mathscr{M}}^{\, \leftrightarrowtriangle}(-1,-1), &\text{$
(\varepsilon_{1},\varepsilon_{2},m(\varepsilon_{2}) \vert \ell) \! = \! (0,-1,-1 \vert 0)$,}
\end{cases}
\end{equation}
and, for $\hat{\ell} \! = \! 0$, set
\begin{equation} \label{laxhat39} 
\hat{\mathscr{F}}^{\lbrace 0 \rbrace}_{\hat{\varepsilon}_{1},\hat{\varepsilon}_{2},\hat{m}
(\hat{\varepsilon}_{2})} \! := \! 
\begin{cases}
\mathscr{F}_{\scriptscriptstyle \mathscr{M}}^{\, \leftrightsquigarrow}(1,1), &\text{$
(\hat{\varepsilon}_{1},\hat{\varepsilon}_{2},\hat{m}(\hat{\varepsilon}_{2}) \vert 
\hat{\ell}) \! = \! (1,1,0 \vert 0)$,} \\
\mathscr{F}_{\scriptscriptstyle \mathscr{M}}^{\, \leftrightsquigarrow}(1,-1), &\text{$
(\hat{\varepsilon}_{1},\hat{\varepsilon}_{2},\hat{m}(\hat{\varepsilon}_{2}) \vert 
\hat{\ell}) \! = \! (1,-1,0 \vert 0)$,} \\
\mathscr{F}_{\scriptscriptstyle \mathscr{M}}^{\, \leftrightsquigarrow}(-1,1), &\text{$
(\hat{\varepsilon}_{1},\hat{\varepsilon}_{2},\hat{m}(\hat{\varepsilon}_{2}) \vert 
\hat{\ell}) \! = \! (-1,1,0 \vert 0)$,} \\
\mathscr{F}_{\scriptscriptstyle \mathscr{M}}^{\, \leftrightsquigarrow}(-1,-1), &\text{$
(\hat{\varepsilon}_{1},\hat{\varepsilon}_{2},\hat{m}(\hat{\varepsilon}_{2}) \vert 
\hat{\ell}) \! = \! (-1,-1,0 \vert 0)$.}
\end{cases}
\end{equation}
Via the definitions \eqref{laxhat37}--\eqref{laxhat39}, define the following compositions (isomorphisms on 
$\mathscr{M}$): for $\ell \! = \! 0$,\footnote{Recall {}from Remarks \ref{newrem12} and \ref{newrem13} that 
$G_{1} \! \equiv \! G_{2} \! \Leftrightarrow \! (G_{1})_{ij} \! = \! -(G_{2})_{ij}$, $i,j \! \in \! \lbrace 1,2 \rbrace$.} set
\begin{alignat}{2}
\mathscr{F}^{\lbrace 0 \rbrace}_{-1,-1,-1} \! :=& \, \mathscr{F}^{\lbrace 0 \rbrace}_{0,-1,-1} 
\circ \mathscr{F}^{\lbrace 0 \rbrace}_{-1,0,0}, & \qquad &(\varepsilon_{1},\varepsilon_{2},
m(\varepsilon_{2}) \vert \ell) \! = \! (-1,-1,-1 \vert 0), \label{laxhat40} \\
\mathscr{F}^{\lbrace 0 \rbrace}_{1,-1,-1} \! :=& \, \mathscr{F}^{\lbrace 0 \rbrace}_{0,-1,-1} 
\circ \mathscr{F}^{\lbrace 0 \rbrace}_{1,0,0}, & \qquad &(\varepsilon_{1},\varepsilon_{2},
m(\varepsilon_{2}) \vert \ell) \! = \! (1,-1,-1 \vert 0), \label{laxhat41} \\
\mathscr{F}^{\lbrace 0 \rbrace}_{-1,-1,1} \! :=& \, \mathscr{F}^{\lbrace 0 \rbrace}_{0,-1,1} 
\circ \mathscr{F}^{\lbrace 0 \rbrace}_{-1,0,0}, & \qquad &(\varepsilon_{1},\varepsilon_{2},
m(\varepsilon_{2}) \vert \ell) \! = \! (-1,-1,1 \vert 0), \label{laxhat42} \\
\mathscr{F}^{\lbrace 0 \rbrace}_{1,-1,1} \! :=& \ \mathscr{F}^{\lbrace 0 \rbrace}_{0,-1,1} 
\circ \mathscr{F}^{\lbrace 0 \rbrace}_{1,0,0}, & \qquad &(\varepsilon_{1},\varepsilon_{2},
m(\varepsilon_{2}) \vert \ell) \! = \! (1,-1,1 \vert 0), \label{laxhat43} \\
\mathscr{F}^{\lbrace 0 \rbrace}_{-1,1,-1} \! :=& \, \mathscr{F}^{\lbrace 0 \rbrace}_{0,1,-1} 
\circ \mathscr{F}^{\lbrace 0 \rbrace}_{-1,0,0}, & \qquad &(\varepsilon_{1},\varepsilon_{2},
m(\varepsilon_{2}) \vert \ell) \! = \! (-1,1,-1 \vert 0), \label{laxhat44} \\
\mathscr{F}^{\lbrace 0 \rbrace}_{1,1,-1} \! :=& \, \mathscr{F}^{\lbrace 0 \rbrace}_{0,1,-1} 
\circ \mathscr{F}^{\lbrace 0 \rbrace}_{1,0,0}, & \qquad &(\varepsilon_{1},\varepsilon_{2},
m(\varepsilon_{2}) \vert \ell) \! = \! (1,1,-1 \vert 0), \label{laxhat45} \\
\mathscr{F}^{\lbrace 0 \rbrace}_{-1,1,1} \! :=& \, \mathscr{F}^{\lbrace 0 \rbrace}_{0,1,1} 
\circ \mathscr{F}^{\lbrace 0 \rbrace}_{-1,0,0}, & \qquad &(\varepsilon_{1},\varepsilon_{2},
m(\varepsilon_{2}) \vert \ell) \! = \! (-1,1,1 \vert 0), \label{laxhat46} \\
\mathscr{F}^{\lbrace 0 \rbrace}_{1,1,1} \! :=& \, \mathscr{F}^{\lbrace 0 \rbrace}_{0,1,1} 
\circ \mathscr{F}^{\lbrace 0 \rbrace}_{1,0,0}, & \qquad &(\varepsilon_{1},\varepsilon_{2},
m(\varepsilon_{2}) \vert \ell) \! = \! (1,1,1 \vert 0), \label{laxhat47}
\end{alignat}
and, for $\hat{\ell} \! = \! 0$, set
\begin{alignat}{2}
\hat{\mathscr{F}}^{\lbrace 0 \rbrace}_{1,0,-1} \! :=& \, \mathscr{F}^{\lbrace 0 \rbrace}_{
0,-1,-1} \circ \hat{\mathscr{F}}^{\lbrace 0 \rbrace}_{1,1,0}, & \qquad &(\hat{\varepsilon}_{1},
\hat{\varepsilon}_{2},\hat{m}(\hat{\varepsilon}_{2}) \vert \hat{\ell}) \! = \! (1,0,-1 \vert 0), 
\label{laxhat48} \\
\hat{\mathscr{F}}^{\lbrace 0 \rbrace}_{-1,0,-1} \! :=& \, \mathscr{F}^{\lbrace 0 \rbrace}_{
0,-1,-1} \circ \hat{\mathscr{F}}^{\lbrace 0 \rbrace}_{-1,1,0}, & \qquad &(\hat{\varepsilon}_{1},
\hat{\varepsilon}_{2},\hat{m}(\hat{\varepsilon}_{2}) \vert \hat{\ell}) \! = \! (-1,0,-1 \vert 0), 
\label{laxhat49} \\
\hat{\mathscr{F}}^{\lbrace 0 \rbrace}_{1,0,1} \! :=& \, \mathscr{F}^{\lbrace 0 \rbrace}_{
0,1,1} \circ \hat{\mathscr{F}}^{\lbrace 0 \rbrace}_{1,-1,0}, & \qquad &(\hat{\varepsilon}_{1},
\hat{\varepsilon}_{2},\hat{m}(\hat{\varepsilon}_{2}) \vert \hat{\ell}) \! = \! (1,0,1 \vert 0), 
\label{laxhat50} \\
\hat{\mathscr{F}}^{\lbrace 0 \rbrace}_{-1,0,1} \! :=& \, \mathscr{F}^{\lbrace 0 \rbrace}_{
0,1,1} \circ \hat{\mathscr{F}}^{\lbrace 0 \rbrace}_{-1,-1,0}, & \qquad &(\hat{\varepsilon}_{1},
\hat{\varepsilon}_{2},\hat{m}(\hat{\varepsilon}_{2}) \vert \hat{\ell}) \! = \! (-1,0,1 \vert 0). 
\label{laxhat51}
\end{alignat}

The cases $\ell,\hat{\ell} \! = \! 1$ are a bit more subtle, because there is no analogue, \emph{per se}, of 
the---standard---identity map \eqref{laxhat37}; rather, the r\^{o}le of the identity map for $\ell,\hat{\ell} \! = \! 1$ 
is mimicked by the endomorphism $\mathscr{F}_{\scriptscriptstyle \mathscr{M}}^{\, \leftslice}(\tilde{l})$, 
$\tilde{l} \! \in \! \lbrace \pm 1 \rbrace$, given in Appendix \ref{sectonsymmt2} (cf. equations \eqref{laxhat17}); 
with conspicuous changes in notation (which are in line with the notations introduced in this subsection), it 
reads (for $\ell \! = \! 1)$:
\begin{align} \label{laxhat52} 
\mathscr{F}^{\lbrace 1 \rbrace}_{0,0,0} \colon \mathscr{M} \! \to \! \mathscr{M}, \, 
&(a,s_{0}^{0},s_{0}^{\infty},s_{1}^{\infty},g_{11},g_{12},g_{21},g_{22}) \! \mapsto \! 
\left(a,s_{0}^{0}(0,0,0 \vert 1),s_{0}^{\infty}(0,0,0 \vert 1),s_{1}^{\infty}(0,0,0 \vert 1), 
\right. \nonumber \\
&\left. \, g_{11}(0,0,0 \vert 1),g_{12}(0,0,0 \vert 1),g_{21}(0,0,0 \vert 1),g_{22}
(0,0,0 \vert 1) \right),
\end{align}
where, for $\tilde{l} \! \in \! \lbrace \pm 1 \rbrace$,
\begin{equation} \label{laxhat53} 
\begin{gathered}
s_{0}^{0}(0,0,0 \vert 1) \! := \! s_{0}^{0}(\tilde{l}), \quad s_{0}^{\infty}(0,0,0 \vert 1) 
\! := \! s_{0}^{\infty}(\tilde{l}), \quad s_{1}^{\infty}(0,0,0 \vert 1) \! := \! 
s_{1}^{\infty}(\tilde{l}), \\
g_{ij}(0,0,0 \vert 1) \! := \! g_{ij}(\tilde{l}), \quad i,j \! \in \! \lbrace 1,2 \rbrace.
\end{gathered}
\end{equation}
To complete the list of the remaining $\ell,\hat{\ell} \! = \! 1$ mappings, define, in analogy with the definitions 
\eqref{laxhat38}--\eqref{laxhat51}, the following compositions (isomorphisms) on $\mathscr{M}$: for $\ell \! = \! 1$,
\begin{alignat}{2}
\mathscr{F}^{\lbrace 1 \rbrace}_{-1,0,0} \! :=& \, \mathscr{F}^{\lbrace 0 \rbrace}_{-1,0,0} 
\circ \mathscr{F}^{\lbrace 1 \rbrace}_{0,0,0}, & \qquad &(\varepsilon_{1},\varepsilon_{2},
m(\varepsilon_{2}) \vert \ell) \! = \! (-1,0,0 \vert 1), \label{laxhat54} \\
\mathscr{F}^{\lbrace 1 \rbrace}_{1,0,0} \! :=& \, \mathscr{F}^{\lbrace 0 \rbrace}_{1,0,0} 
\circ \mathscr{F}^{\lbrace 1 \rbrace}_{0,0,0}, & \qquad &(\varepsilon_{1},\varepsilon_{2},
m(\varepsilon_{2}) \vert \ell) \! = \! (1,0,0 \vert 1), \label{laxhat55} \\
\mathscr{F}^{\lbrace 1 \rbrace}_{0,-1,-1} \! :=& \, \mathscr{F}^{\lbrace 0 \rbrace}_{0,-1,-1} 
\circ \mathscr{F}^{\lbrace 1 \rbrace}_{0,0,0}, & \qquad &(\varepsilon_{1},\varepsilon_{2},
m(\varepsilon_{2}) \vert \ell) \! = \! (0,-1,-1 \vert 1), \label{laxhat56} \\
\mathscr{F}^{\lbrace 1 \rbrace}_{0,-1,1} \! :=& \ \mathscr{F}^{\lbrace 0 \rbrace}_{0,-1,1} 
\circ \mathscr{F}^{\lbrace 1 \rbrace}_{0,0,0}, & \qquad &(\varepsilon_{1},\varepsilon_{2},
m(\varepsilon_{2}) \vert \ell) \! = \! (0,-1,1 \vert 1), \label{laxhat57} \\
\mathscr{F}^{\lbrace 1 \rbrace}_{0,1,-1} \! :=& \, \mathscr{F}^{\lbrace 0 \rbrace}_{0,1,-1} 
\circ \mathscr{F}^{\lbrace 1 \rbrace}_{0,0,0}, & \qquad &(\varepsilon_{1},\varepsilon_{2},
m(\varepsilon_{2}) \vert \ell) \! = \! (0,1,-1 \vert 1), \label{laxhat58} \\
\mathscr{F}^{\lbrace 1 \rbrace}_{0,1,1} \! :=& \, \mathscr{F}^{\lbrace 0 \rbrace}_{0,1,1} 
\circ \mathscr{F}^{\lbrace 1 \rbrace}_{0,0,0}, & \qquad &(\varepsilon_{1},\varepsilon_{2},
m(\varepsilon_{2}) \vert \ell) \! = \! (0,1,1 \vert 1), \label{laxhat59} \\
\mathscr{F}^{\lbrace 1 \rbrace}_{-1,-1,-1} \! :=& \, \mathscr{F}^{\lbrace 1 \rbrace}_{0,-1,-1} 
\circ \mathscr{F}^{\lbrace 0 \rbrace}_{-1,0,0}, & \qquad &(\varepsilon_{1},\varepsilon_{2},
m(\varepsilon_{2}) \vert \ell) \! = \! (-1,-1,-1 \vert 1), \label{laxhat60} \\
\mathscr{F}^{\lbrace 1 \rbrace}_{1,-1,-1} \! :=& \, \mathscr{F}^{\lbrace 1 \rbrace}_{0,-1,-1} 
\circ \mathscr{F}^{\lbrace 0 \rbrace}_{1,0,0}, & \qquad &(\varepsilon_{1},\varepsilon_{2},
m(\varepsilon_{2}) \vert \ell) \! = \! (1,-1,-1 \vert 1), \label{laxhat61} \\
\mathscr{F}^{\lbrace 1 \rbrace}_{-1,-1,1} \! :=& \, \mathscr{F}^{\lbrace 1 \rbrace}_{0,-1,1} 
\circ \mathscr{F}^{\lbrace 0 \rbrace}_{-1,0,0}, & \qquad &(\varepsilon_{1},\varepsilon_{2},
m(\varepsilon_{2}) \vert \ell) \! = \! (-1,-1,1 \vert 1), \label{laxhat62} \\
\mathscr{F}^{\lbrace 1 \rbrace}_{1,-1,1} \! :=& \, \mathscr{F}^{\lbrace 1 \rbrace}_{0,-1,1} 
\circ \mathscr{F}^{\lbrace 0 \rbrace}_{1,0,0}, & \qquad &(\varepsilon_{1},\varepsilon_{2},
m(\varepsilon_{2}) \vert \ell) \! = \! (1,-1,1 \vert 1), \label{laxhat63} \\
\mathscr{F}^{\lbrace 1 \rbrace}_{-1,1,-1} \! :=& \, \mathscr{F}^{\lbrace 1 \rbrace}_{0,1,-1} 
\circ \mathscr{F}^{\lbrace 0 \rbrace}_{-1,0,0}, & \qquad &(\varepsilon_{1},\varepsilon_{2},
m(\varepsilon_{2}) \vert \ell) \! = \! (-1,1,-1 \vert 1), \label{laxhat64} \\
\mathscr{F}^{\lbrace 1 \rbrace}_{1,1,-1} \! :=& \ \mathscr{F}^{\lbrace 1 \rbrace}_{0,1,-1} 
\circ \mathscr{F}^{\lbrace 0 \rbrace}_{1,0,0}, & \qquad &(\varepsilon_{1},\varepsilon_{2},
m(\varepsilon_{2}) \vert \ell) \! = \! (1,1,-1 \vert 1), \label{laxhat65} \\
\mathscr{F}^{\lbrace 1 \rbrace}_{-1,1,1} \! :=& \, \mathscr{F}^{\lbrace 1 \rbrace}_{0,1,1} 
\circ \mathscr{F}^{\lbrace 0 \rbrace}_{-1,0,0}, & \qquad &(\varepsilon_{1},\varepsilon_{2},
m(\varepsilon_{2}) \vert \ell) \! = \! (-1,1,1 \vert 1), \label{laxhat66} \\
\mathscr{F}^{\lbrace 1 \rbrace}_{1,1,1} \! :=& \, \mathscr{F}^{\lbrace 1 \rbrace}_{0,1,1} 
\circ \mathscr{F}^{\lbrace 0 \rbrace}_{1,0,0}, & \qquad &(\varepsilon_{1},\varepsilon_{2},
m(\varepsilon_{2}) \vert \ell) \! = \! (1,1,1 \vert 1); \label{laxhat67}
\end{alignat}
and, for $\hat{\ell} \! = \! 1$,
\begin{alignat}{2}
\hat{\mathscr{F}}^{\lbrace 1 \rbrace}_{1,1,0} \! :=& \, \hat{\mathscr{F}}^{\lbrace 0 
\rbrace}_{1,1,0} \circ \mathscr{F}^{\lbrace 1 \rbrace}_{0,0,0}, & \qquad 
&(\hat{\varepsilon}_{1},\hat{\varepsilon}_{2},\hat{m}(\hat{\varepsilon}_{2}) \vert 
\hat{\ell}) \! = \! (1,1,0 \vert 1), \label{laxhat68} \\
\hat{\mathscr{F}}^{\lbrace 1 \rbrace}_{1,-1,0} \! :=& \, \hat{\mathscr{F}}^{\lbrace 0 
\rbrace}_{1,-1,0} \circ \mathscr{F}^{\lbrace 1 \rbrace}_{0,0,0}, & \qquad 
&(\hat{\varepsilon}_{1},\hat{\varepsilon}_{2},\hat{m}(\hat{\varepsilon}_{2}) \vert 
\hat{\ell}) \! = \! (1,-1,0 \vert 1), \label{laxhat69} \\
\hat{\mathscr{F}}^{\lbrace 1 \rbrace}_{-1,1,0} \! :=& \, \hat{\mathscr{F}}^{\lbrace 0 
\rbrace}_{-1,1,0} \circ \mathscr{F}^{\lbrace 1 \rbrace}_{0,0,0}, & \qquad 
&(\hat{\varepsilon}_{1},\hat{\varepsilon}_{2},\hat{m}(\hat{\varepsilon}_{2}) \vert 
\hat{\ell}) \! = \! (-1,1,0 \vert 1), \label{laxhat70} \\
\hat{\mathscr{F}}^{\lbrace 1 \rbrace}_{-1,-1,0} \! :=& \, \hat{\mathscr{F}}^{\lbrace 0 
\rbrace}_{-1,-1,0} \circ \mathscr{F}^{\lbrace 1 \rbrace}_{0,0,0}, & \qquad 
&(\hat{\varepsilon}_{1},\hat{\varepsilon}_{2},\hat{m}(\hat{\varepsilon}_{2}) \vert 
\hat{\ell}) \! = \! (-1,-1,0 \vert 1), \label{laxhat71} \\
\hat{\mathscr{F}}^{\lbrace 1 \rbrace}_{1,0,-1} \! :=& \, \mathscr{F}^{\lbrace 0 
\rbrace}_{0,1,-1} \circ \hat{\mathscr{F}}^{\lbrace 1 \rbrace}_{1,-1,0}, & \qquad 
&(\hat{\varepsilon}_{1},\hat{\varepsilon}_{2},\hat{m}(\hat{\varepsilon}_{2}) \vert 
\hat{\ell}) \! = \! (1,0,-1 \vert 1), \label{laxhat72} \\
\hat{\mathscr{F}}^{\lbrace 1 \rbrace}_{1,0,1} \! :=& \, \mathscr{F}^{\lbrace 0 
\rbrace}_{0,1,1} \circ \hat{\mathscr{F}}^{\lbrace 1 \rbrace}_{1,-1,0}, & \qquad 
&(\hat{\varepsilon}_{1},\hat{\varepsilon}_{2},\hat{m}(\hat{\varepsilon}_{2}) \vert 
\hat{\ell}) \! = \! (1,0,1 \vert 1), \label{laxhat73} \\
\hat{\mathscr{F}}^{\lbrace 1 \rbrace}_{-1,0,-1} \! :=& \, \mathscr{F}^{\lbrace 0 
\rbrace}_{0,1,-1} \circ \hat{\mathscr{F}}^{\lbrace 1 \rbrace}_{-1,-1,0}, & \qquad 
&(\hat{\varepsilon}_{1},\hat{\varepsilon}_{2},\hat{m}(\hat{\varepsilon}_{2}) \vert 
\hat{\ell}) \! = \! (-1,0,-1 \vert 1), \label{laxhat74} \\
\hat{\mathscr{F}}^{\lbrace 1 \rbrace}_{-1,0,1} \! :=& \, \mathscr{F}^{\lbrace 0 
\rbrace}_{0,1,1} \circ \hat{\mathscr{F}}^{\lbrace 1 \rbrace}_{-1,-1,0}, & \qquad 
&(\hat{\varepsilon}_{1},\hat{\varepsilon}_{2},\hat{m}(\hat{\varepsilon}_{2}) \vert 
\hat{\ell}) \! = \! (-1,0,1 \vert 1). \label{laxhat75}
\end{alignat}

Via the elementary symmetries \eqref{laxhat8}, \eqref{laxhat9}, \eqref{laxhat17}, \eqref{laxhat27}--\eqref{laxhat30}, 
and \eqref{laxhat31}--\eqref{laxhat34}, and the definitions \eqref{laxhat37}--\eqref{laxhat75}, one arrives at the 
following explicit list of actions on $\mathscr{M}$ of the isomorphisms (cf. definition \eqref{laxhat35}) 
$\mathscr{F}^{\scriptscriptstyle \lbrace \ell \rbrace}_{\scriptscriptstyle \varepsilon_{1},
\varepsilon_{2},m(\varepsilon_{2})}$, relevant for real $\tau$, and (cf. definition \eqref{laxhat36}) 
$\hat{\mathscr{F}}^{\scriptscriptstyle \lbrace \hat{\ell} \rbrace}_{\scriptscriptstyle \hat{\varepsilon}_{1},
\hat{\varepsilon}_{2},\hat{m}(\hat{\varepsilon}_{2})}$, relevant for pure-imaginary $\tau$: for $\tilde{l},
l^{\prime} \! \in \! \lbrace \pm 1 \rbrace$,
\begin{enumerate}
\item[(1)] $\mathscr{F}^{\lbrace 0 \rbrace}_{0,0,0} \Rightarrow$
\begin{equation} \label{laxhat76} 
\begin{gathered}
s_{0}^{0}(0,0,0 \vert 0) \! = \! s_{0}^{0}, \quad s_{0}^{\infty}(0,0,0 \vert 0) \! = \! 
s_{0}^{\infty}, \quad s_{1}^{\infty}(0,0,0 \vert 0) \! = \! s_{1}^{\infty}, \quad g_{ij}
(0,0,0 \vert 0) \! = \! g_{ij}, \quad i,j \! \in \! \lbrace 1,2 \rbrace;
\end{gathered}
\end{equation}
\item[(2)] $\mathscr{F}^{\lbrace 0 \rbrace}_{-1,0,0} \Rightarrow$
\begin{equation} \label{laxhat77} 
\begin{gathered}
s_{0}^{0}(-1,0,0 \vert 0) \! = \! s_{0}^{0}, \quad s_{0}^{\infty}(-1,0,0 \vert 0) \! = \! 
s_{0}^{\infty} \me^{\pi a}, \quad s_{1}^{\infty}(-1,0,0 \vert 0) \! = \! s_{1}^{\infty} 
\me^{-\pi a}, \\
g_{11}(-1,0,0 \vert 0) \! = \! -\mi (g_{21} \! + \! s_{0}^{0}g_{11}) \me^{\pi a/2}, 
\quad g_{12}(-1,0,0 \vert 0) \! = \! -\mi (g_{22} \! + \! s_{0}^{0}g_{12}) 
\me^{-\pi a/2}, \\
g_{21}(-1,0,0 \vert 0) \! = \! -\mi g_{11} \me^{\pi a/2}, \quad g_{22}
(-1,0,0 \vert 0) \! = \! -\mi g_{12} \me^{-\pi a/2};
\end{gathered}
\end{equation}
\item[(3)] $\mathscr{F}^{\lbrace 0 \rbrace}_{1,0,0} \Rightarrow$
\begin{equation} \label{laxhat78} 
\begin{gathered}
s_{0}^{0}(1,0,0 \vert 0) \! = \! s_{0}^{0}, \quad s_{0}^{\infty}(1,0,0 \vert 0) \! = \! 
s_{0}^{\infty} \me^{-\pi a}, \quad s_{1}^{\infty}(1,0,0 \vert 0) \! = \! s_{1}^{\infty} 
\me^{\pi a}, \\
g_{11}(1,0,0 \vert 0) \! = \! \mi g_{21} \me^{-\pi a/2}, \quad g_{12}(1,0,0 \vert 0) 
\! = \! \mi g_{22} \me^{\pi a/2}, \\
g_{21}(1,0,0 \vert 0) \! = \! \mi (g_{11} \! - \! s_{0}^{0}g_{21}) \me^{-\pi a/2}, \quad 
g_{22}(1,0,0 \vert 0) \! = \! \mi (g_{12} \! - \! s_{0}^{0}g_{22}) \me^{\pi a/2};
\end{gathered}
\end{equation}
\item[(4)] $\mathscr{F}^{\lbrace 0 \rbrace}_{0,-1,-1} \Rightarrow$
\begin{equation} \label{laxhat79} 
\begin{gathered}
s_{0}^{0}(0,-1,-1 \vert 0) \! = \! s_{0}^{0}, \quad s_{0}^{\infty}(0,-1,-1 \vert 0) \! = \! 
-s_{1}^{\infty} \me^{\pi a}, \quad s_{1}^{\infty}(0,-1,-1 \vert 0) \! = \! -s_{0}^{\infty} 
\me^{\pi a}, \\
g_{11}(0,-1,-1 \vert 0) \! = \! \mi l^{\prime}(g_{22} \! - \! g_{21}s_{1}^{\infty} 
\me^{2 \pi a} \! + \! s_{0}^{0}(g_{12} \! - \! g_{11}s_{1}^{\infty} \me^{2 \pi a})) 
\me^{-\pi a/2}, \\
g_{12}(0,-1,-1 \vert 0) \! = \! -\mi l^{\prime}(g_{21} \! + \! s_{0}^{0}g_{11}) 
\me^{\pi a/2}, \quad g_{21}(0,-1,-1 \vert 0) \! = \! \mi l^{\prime}(g_{12} 
\!  - \! g_{11}s_{1}^{\infty} \me^{2 \pi a}) \me^{-\pi a/2}, \\
g_{22}(0,-1,-1 \vert 0) \! = \! -\mi l^{\prime}g_{11} \me^{\pi a/2};
\end{gathered}
\end{equation}
\item[(5)] $\mathscr{F}^{\lbrace 0 \rbrace}_{0,-1,1} \Rightarrow$
\begin{equation} \label{laxhat80} 
\begin{gathered}
s_{0}^{0}(0,-1,1 \vert 0) \! = \! s_{0}^{0}, \quad s_{0}^{\infty}(0,-1,1 \vert 0) \! = \! 
-s_{1}^{\infty} \me^{\pi a}, \quad s_{1}^{\infty}(0,-1,1 \vert 0) \! = \! -s_{0}^{\infty} 
\me^{\pi a}, \\
g_{11}(0,-1,1 \vert 0) \! = \! l^{\prime}g_{12} \me^{\pi a/2}, \quad g_{12}(0,-1,1 
\vert 0) \! = \! -l^{\prime}(g_{11} \! + \! s_{0}^{\infty}g_{12}) \me^{-\pi a/2}, \\
g_{21}(0,-1,1 \vert 0) \! = \! l^{\prime}g_{22} \me^{\pi a/2}, \quad g_{22}
(0,-1,1 \vert 0) \! = \! -l^{\prime}(g_{21} \! + \! s_{0}^{\infty}g_{22}) 
\me^{-\pi a/2};
\end{gathered}
\end{equation}
\item[(6)] $\mathscr{F}^{\lbrace 0 \rbrace}_{0,1,-1} \Rightarrow$
\begin{equation} \label{laxhat81} 
\begin{gathered}
s_{0}^{0}(0,1,-1 \vert 0) \! = \! s_{0}^{0}, \quad s_{0}^{\infty}(0,1,-1 \vert 0) \! = \! 
-s_{1}^{\infty} \me^{\pi a}, \quad s_{1}^{\infty}(0,1,-1 \vert 0) \! = \! -s_{0}^{\infty} 
\me^{\pi a}, \\
g_{11}(0,1,-1 \vert 0) \! = \! l^{\prime}(g_{12} \! - \! g_{11}s_{1}^{\infty} 
\me^{2 \pi a}) \me^{-\pi a/2}, \quad g_{12}(0,1,-1 \vert 0) \! = \! -l^{\prime}
g_{11} \me^{\pi a/2}, \\
g_{21}(0,1,-1 \vert 0) \! = \! l^{\prime}(g_{22} \! - \! g_{21}s_{1}^{\infty} 
\me^{2 \pi a}) \me^{-\pi a/2}, \quad g_{22}(0,1,-1 \vert 0) \! = \! -l^{\prime}
g_{21} \me^{\pi a/2};
\end{gathered}
\end{equation}
\item[(7)] $\mathscr{F}^{\lbrace 0 \rbrace}_{0,1,1} \Rightarrow$
\begin{equation} \label{laxhat82} 
\begin{gathered}
s_{0}^{0}(0,1,1 \vert 0) \! = \! s_{0}^{0}, \quad s_{0}^{\infty}(0,1,1 \vert 0) \! = \! 
-s_{1}^{\infty} \me^{\pi a}, \quad s_{1}^{\infty}(0,1,1 \vert 0) \! = \! -s_{0}^{\infty} 
\me^{\pi a}, \\
g_{11}(0,1,1 \vert 0) \! = \! \mi l^{\prime}g_{22} \me^{\pi a/2}, \quad g_{12}
(0,1,1 \vert 0) \! = \! -\mi l^{\prime}(g_{21} \! + \! s_{0}^{\infty}g_{22}) 
\me^{-\pi a/2}, \\
g_{21}(0,1,1 \vert 0) \! = \! \mi l^{\prime}(g_{12} \! - \! s_{0}^{0}g_{22}) 
\me^{\pi a/2}, \quad g_{22}(0,1,1 \vert 0) \! = \! \mi l^{\prime}(-g_{11} \! 
- \! g_{12}s_{0}^{\infty} \! + \! s_{0}^{0}(g_{21} \! + \! s_{0}^{\infty}g_{22})) 
\me^{-\pi a/2};
\end{gathered}
\end{equation}
\item[(8)] $\mathscr{F}^{\lbrace 0 \rbrace}_{-1,-1,-1} \Rightarrow$
\begin{equation} \label{laxhat83} 
\begin{gathered}
s_{0}^{0}(-1,-1,-1 \vert 0) \! = \! s_{0}^{0}, \quad s_{0}^{\infty}(-1,-1,-1 \vert 0) 
\! = \! -s_{1}^{\infty}, \quad s_{1}^{\infty}(-1,-1,-1 \vert 0) \! = \! -s_{0}^{\infty} 
\me^{2 \pi a}, \\
g_{11}(-1,-1,-1 \vert 0) \! = \! l^{\prime}((g_{12} \! - \! g_{11}s_{1}^{\infty} 
\me^{2 \pi a})(1 \! + \! (s_{0}^{0})^{2}) \! + \! s_{0}^{0}(g_{22} \! - \! g_{21}
s_{1}^{\infty} \me^{2 \pi a})) \me^{-\pi a}, \\
g_{12}(-1,-1,-1 \vert 0) \! = \! -l^{\prime}(g_{11}(1 \! + \! (s_{0}^{0})^{2}) 
\! + \! s_{0}^{0}g_{21}) \me^{\pi a}, \\
g_{21}(-1,-1,-1 \vert 0) \! = \! l^{\prime}(g_{22} \! - \! g_{21}s_{1}^{\infty} 
\me^{2 \pi a} \! + \! s_{0}^{0}(g_{12} \! - \! g_{11}s_{1}^{\infty} \me^{2 \pi a})) 
\me^{-\pi a}, \\
g_{22}(-1,-1,-1 \vert 0) \! = \! -l^{\prime}(g_{21} \! + \! s_{0}^{0}g_{11}) \me^{\pi a};
\end{gathered}
\end{equation}
\item[(9)] $\mathscr{F}^{\lbrace 0 \rbrace}_{1,-1,-1} \Rightarrow$
\begin{equation} \label{laxhat84} 
\begin{gathered}
s_{0}^{0}(1,-1,-1 \vert 0) \! = \! s_{0}^{0}, \quad s_{0}^{\infty}(1,-1,-1 \vert 0) \! = \! 
-s_{1}^{\infty} \me^{2 \pi a}, \quad s_{1}^{\infty}(1,-1,-1 \vert 0) \! = \! -s_{0}^{\infty}, \\
g_{11}(1,-1,-1 \vert 0) \! = \! -l^{\prime}(g_{12} \! - \! g_{11}s_{1}^{\infty} \me^{2 \pi a}), 
\quad g_{12}(1,-1,-1 \vert 0) \! = \! l^{\prime}g_{11}, \\
g_{21}(1,-1,-1 \vert 0) \! = \! -l^{\prime}(g_{22} \! - \! g_{21}s_{1}^{\infty} \me^{2 \pi a}), 
\quad g_{22}(1,-1,-1 \vert 0) \! = \! l^{\prime}g_{21};
\end{gathered}
\end{equation}
\item[(10)] $\mathscr{F}^{\lbrace 0 \rbrace}_{-1,-1,1} \Rightarrow$
\begin{equation} \label{laxhat85} 
\begin{gathered}
s_{0}^{0}(-1,-1,1 \vert 0) \! = \! s_{0}^{0}, \quad s_{0}^{\infty}(-1,-1,1 \vert 0) \! = \! 
-s_{1}^{\infty}, \quad s_{1}^{\infty}(-1,-1,1 \vert 0) \! = \! -s_{0}^{\infty} \me^{2 \pi a}, \\
g_{11}(-1,-1,1 \vert 0) \! = \! -\mi l^{\prime}(g_{22} \! + \! s_{0}^{0}g_{12}), \quad 
g_{12}(-1,-1,1 \vert 0) \! = \! \mi l^{\prime}(g_{21} \! + \! s_{0}^{\infty}g_{22} \! + \! 
s_{0}^{0}(g_{11} \! + \! s_{0}^{\infty}g_{12})), \\
g_{21}(-1,-1,1 \vert 0) \! = \! -\mi l^{\prime}g_{12}, \quad g_{22}(-1,-1,1 \vert 0) \! 
= \! \mi l^{\prime}(g_{11} \! + \! s_{0}^{\infty}g_{12});
\end{gathered}
\end{equation}
\item[(11)] $\mathscr{F}^{\lbrace 0 \rbrace}_{1,-1,1} \Rightarrow$
\begin{equation} \label{laxhat86} 
\begin{gathered}
s_{0}^{0}(1,-1,1 \vert 0) \! = \! s_{0}^{0}, \quad s_{0}^{\infty}(1,-1,1 \vert 0) \! = \! 
-s_{1}^{\infty} \me^{2 \pi a}, \quad s_{1}^{\infty}(1,-1,1 \vert 0) \! = \! -s_{0}^{\infty}, \\
g_{11}(1,-1,1 \vert 0) \! = \! \mi l^{\prime}g_{22} \me^{\pi a}, \quad g_{12}
(1,-1,1 \vert 0) \! = \! -\mi l^{\prime}(g_{21} \! + \! s_{0}^{\infty}g_{22}) 
\me^{-\pi a}, \\
g_{21}(1,-1,1 \vert 0) \! = \! \mi l^{\prime}(g_{12} \! - \! s_{0}^{0}g_{22}) \me^{\pi a}, 
\quad g_{22}(1,-1,1 \vert 0) \! = \! -\mi l^{\prime}(g_{11} \! + \! s_{0}^{\infty}g_{12} 
\! - \! s_{0}^{0}(g_{21} \! + \! s_{0}^{\infty}g_{22})) \me^{-\pi a};
\end{gathered}
\end{equation}
\item[(12)] $\mathscr{F}^{\lbrace 0 \rbrace}_{-1,1,-1} \Rightarrow$
\begin{equation} \label{laxhat87} 
\begin{gathered}
s_{0}^{0}(-1,1,-1 \vert 0) \! = \! s_{0}^{0}, \quad s_{0}^{\infty}(-1,1,-1 \vert 0) \! = \! 
-s_{1}^{\infty}, \quad s_{1}^{\infty}(-1,1,-1 \vert 0) \! = \! -s_{0}^{\infty} \me^{2 \pi a}, \\
g_{11}(-1,1,-1 \vert 0) \! = \! -\mi l^{\prime}(g_{22} \! - \! g_{21}s_{1}^{\infty} 
\me^{2 \pi a} \! + \! s_{0}^{0}(g_{12} \! - \! g_{11}s_{1}^{\infty} \me^{2 \pi a})) 
\me^{-\pi a}, \\
g_{12}(-1,1,-1 \vert 0) \! = \! \mi l^{\prime}(g_{21} \! + \! s_{0}^{0}g_{11}) \me^{\pi a}, 
\quad g_{21}(-1,1,-1 \vert 0) \! = \! -\mi l^{\prime}(g_{12} \!  - \! g_{11}s_{1}^{\infty} 
\me^{2 \pi a}) \me^{-\pi a}, \\
g_{22}(-1,1,-1 \vert 0) \! = \! \mi l^{\prime}g_{11} \me^{\pi a};
\end{gathered}
\end{equation}
\item[(13)] $\mathscr{F}^{\lbrace 0 \rbrace}_{1,1,-1} \Rightarrow$
\begin{equation} \label{laxhat88} 
\begin{gathered}
s_{0}^{0}(1,1,-1 \vert 0) \! = \! s_{0}^{0}, \quad s_{0}^{\infty}(1,1,-1 \vert 0) \! = \! 
-s_{1}^{\infty} \me^{2 \pi a}, \quad s_{1}^{\infty}(1,1,-1 \vert 0) \! = \! -s_{0}^{\infty}, \\
g_{11}(1,1,-1 \vert 0) \! = \! \mi l^{\prime}(g_{22} \! - \! g_{21}s_{1}^{\infty} 
\me^{2 \pi a}), \quad g_{12}(1,1,-1 \vert 0) \! = \! -\mi l^{\prime}g_{21}, \\
g_{21}(1,1,-1 \vert 0) \! = \! \mi l^{\prime}(g_{12} \!  - \! g_{11}s_{1}^{\infty} 
\me^{2 \pi a} \! - \! s_{0}^{0}(g_{22} \! - \! g_{21}s_{1}^{\infty} \me^{2 \pi a})), 
\quad g_{22}(1,1,-1 \vert 0) \! = \! -\mi l^{\prime}(g_{11} \! - \! s_{0}^{0}g_{21});
\end{gathered}
\end{equation}
\item[(14)] $\mathscr{F}^{\lbrace 0 \rbrace}_{-1,1,1} \Rightarrow$
\begin{equation} \label{laxhat89} 
\begin{gathered}
s_{0}^{0}(-1,1,1 \vert 0) \! = \! s_{0}^{0}, \quad s_{0}^{\infty}(-1,1,1 \vert 0) \! = \! 
-s_{1}^{\infty}, \quad s_{1}^{\infty}(-1,1,1 \vert 0) \! = \! -s_{0}^{\infty} 
\me^{2 \pi a}, \\
g_{11}(-1,1,1 \vert 0) \! = \! l^{\prime}g_{12}, \quad g_{12}(-1,1,1 \vert 0) \! = \! 
-l^{\prime}(g_{11} \! + \! s_{0}^{\infty}g_{12}), \\
g_{21}(-1,1,1 \vert 0) \! = \! l^{\prime}g_{22}, \quad g_{22}(-1,1,1 \vert 0) \! = \! 
-l^{\prime}(g_{21} \! + \! s_{0}^{\infty}g_{22});
\end{gathered}
\end{equation}
\item[(15)] $\mathscr{F}^{\lbrace 0 \rbrace}_{1,1,1} \Rightarrow$
\begin{equation} \label{laxhat90} 
\begin{gathered}
s_{0}^{0}(1,1,1 \vert 0) \! = \! s_{0}^{0}, \quad s_{0}^{\infty}(1,1,1 \vert 0) \! = \! 
-s_{1}^{\infty} \me^{2 \pi a}, \quad s_{1}^{\infty}(1,1,1 \vert 0) \! = \! 
-s_{0}^{\infty}, \\
g_{11}(1,1,1 \vert 0) \! = \! -l^{\prime}(g_{12} \! - \! s_{0}^{0}g_{22}) \me^{\pi a}, 
\quad g_{12}(1,1,1 \vert 0) \! = \! -l^{\prime}(-g_{11} \! - \! g_{12}s_{0}^{\infty} 
\! + \! s_{0}^{0}(g_{21} \! + \! g_{22}s_{0}^{\infty})) \me^{-\pi a}, \\
g_{21}(1,1,1 \vert 0) \! = \! -l^{\prime}(g_{22} \! - \! s_{0}^{0}(g_{12} \! - 
\! s_{0}^{0}g_{22})) \me^{\pi a}, \\
g_{22}(1,1,1 \vert 0) \! = \! l^{\prime}((g_{21} \! + \! g_{22}s_{0}^{\infty})
(1 \! + \! (s_{0}^{0})^{2}) \! - \! s_{0}^{0}(g_{11} \! + \! s_{0}^{\infty}g_{12})) 
\me^{-\pi a};
\end{gathered}
\end{equation}
\item[(16)] $\hat{\mathscr{F}}^{\lbrace 0 \rbrace}_{1,1,0} \Rightarrow$
\begin{equation} \label{laxhat91} 
\begin{gathered}
\hat{s}_{0}^{0}(1,1,0 \vert 0) \! = \! s_{0}^{0}, \quad \hat{s}_{0}^{\infty}(1,1,0 \vert 0) 
\! = \! s_{0}^{\infty} \me^{-\pi a/2}, \quad \hat{s}_{1}^{\infty}(1,1,0 \vert 0) \! = \! 
s_{1}^{\infty} \me^{\pi a/2}, \\
\hat{g}_{11}(1,1,0 \vert 0) \! = \! -\mi g_{21} \me^{-\pi a/4}, \quad \hat{g}_{12}
(1,1,0 \vert 0) \! = \! -\mi g_{22} \me^{\pi a/4}, \\
\hat{g}_{21}(1,1,0 \vert 0) \! = \! -\mi (g_{11} \! - \! s_{0}^{0}g_{21}) \me^{-\pi a/4}, 
\quad \hat{g}_{22}(1,1,0 \vert 0) \! = \! -\mi (g_{12} \! - \! s_{0}^{0}g_{22}) 
\me^{\pi a/4};
\end{gathered}
\end{equation}
\item[(17)] $\hat{\mathscr{F}}^{\lbrace 0 \rbrace}_{1,-1,0} \Rightarrow$
\begin{equation} \label{laxhat92} 
\begin{gathered}
\hat{s}_{0}^{0}(1,-1,0 \vert 0) \! = \! s_{0}^{0}, \quad \hat{s}_{0}^{\infty}(1,-1,0 \vert 0) 
\! = \! s_{0}^{\infty} \me^{-\pi a/2}, \quad \hat{s}_{1}^{\infty}(1,-1,0 \vert 0) \! = \! 
s_{1}^{\infty} \me^{\pi a/2}, \\
\hat{g}_{11}(1,-1,0 \vert 0) \! = \! g_{11} \me^{-\pi a/4}, \quad \hat{g}_{12}
(1,-1,0 \vert 0) \! = \! g_{12} \me^{\pi a/4}, \\
\hat{g}_{21}(1,-1,0 \vert 0) \! = \! g_{21} \me^{-\pi a/4}, \quad \hat{g}_{22}
(1,-1,0 \vert 0) \! = \! g_{22} \me^{\pi a/4};
\end{gathered}
\end{equation}
\item[(18)] $\hat{\mathscr{F}}^{\lbrace 0 \rbrace}_{-1,1,0} \Rightarrow$
\begin{equation} \label{laxhat93} 
\begin{gathered}
\hat{s}_{0}^{0}(-1,1,0 \vert 0) \! = \! s_{0}^{0}, \quad \hat{s}_{0}^{\infty}(-1,1,0 \vert 0) 
\! = \! s_{0}^{\infty} \me^{\pi a/2}, \quad \hat{s}_{1}^{\infty}(-1,1,0 \vert 0) \! = \! 
s_{1}^{\infty} \me^{-\pi a/2}, \\
\hat{g}_{11}(-1,1,0 \vert 0) \! = \! g_{11} \me^{\pi a/4}, \quad \hat{g}_{12}
(-1,1,0 \vert 0) \! = \! g_{12} \me^{-\pi a/4}, \\
\hat{g}_{21}(-1,1,0 \vert 0) \! = \! g_{21} \me^{\pi a/4}, \quad \hat{g}_{22}
(-1,1,0 \vert 0) \! = \! g_{22} \me^{-\pi a/4};
\end{gathered}
\end{equation}
\item[(19)] $\hat{\mathscr{F}}^{\lbrace 0 \rbrace}_{-1,-1,0} \Rightarrow$
\begin{equation} \label{laxhat94} 
\begin{gathered}
\hat{s}_{0}^{0}(-1,-1,0 \vert 0) \! = \! s_{0}^{0}, \quad \hat{s}_{0}^{\infty}(-1,-1,0 \vert 0) 
\! = \! s_{0}^{\infty} \me^{\pi a/2}, \quad \hat{s}_{1}^{\infty}(-1,-1,0 \vert 0) \! = \! 
s_{1}^{\infty} \me^{-\pi a/2}, \\
\hat{g}_{11}(-1,-1,0 \vert 0) \! = \! \mi (g_{21} \! + \! s_{0}^{0}g_{11}) \me^{\pi a/4}, 
\quad \hat{g}_{12}(-1,-1,0 \vert 0) \! = \! \mi (g_{22} \! + \! s_{0}^{0}g_{12}) 
\me^{-\pi a/4}, \\
\hat{g}_{21}(-1,-1,0 \vert 0) \! = \! \mi g_{11} \me^{\pi a/4}, \quad \hat{g}_{22}
(-1,-1,0 \vert 0) \! = \! \mi g_{12} \me^{-\pi a/4};
\end{gathered}
\end{equation}
\item[(20)] $\hat{\mathscr{F}}^{\lbrace 0 \rbrace}_{1,0,-1} \Rightarrow$
\begin{equation} \label{laxhat95} 
\begin{gathered}
\hat{s}_{0}^{0}(1,0,-1 \vert 0) \! = \! s_{0}^{0}, \quad \hat{s}_{0}^{\infty}(1,0,-1 \vert 
0) \! = \! -s_{1}^{\infty} \me^{3 \pi a/2}, \quad \hat{s}_{1}^{\infty}(1,0,-1 \vert 0) 
\! = \! -s_{0}^{\infty} \me^{\pi a/2}, \\
\hat{g}_{11}(1,0,-1 \vert 0) \! = \! l^{\prime}(g_{12} \! - \! g_{11}s_{1}^{\infty} 
\me^{2 \pi a}) \me^{-\pi a/4}, \quad \hat{g}_{12}(1,0,-1 \vert 0) \! = \! -l^{\prime}
g_{11} \me^{\pi a/4}, \\
\hat{g}_{21}(1,0,-1 \vert 0) \! = \! l^{\prime}(g_{22} \! - \! g_{21}s_{1}^{\infty} 
\me^{2 \pi a}) \me^{-\pi a/4}, \quad \hat{g}_{22}(1,0,-1 \vert 0) \! = \! -l^{\prime}
g_{21} \me^{\pi a/4};
\end{gathered}
\end{equation}
\item[(21)] $\hat{\mathscr{F}}^{\lbrace 0 \rbrace}_{-1,0,-1} \Rightarrow$
\begin{equation} \label{laxhat96} 
\begin{gathered}
\hat{s}_{0}^{0}(-1,0,-1 \vert 0) \! = \! s_{0}^{0}, \quad \hat{s}_{0}^{\infty}(-1,0,-1 \vert 0) 
\! = \! -s_{1}^{\infty} \me^{\pi a/2}, \quad \hat{s}_{1}^{\infty}(-1,0,-1 \vert 0) \! = \! 
-s_{0}^{\infty} \me^{3 \pi a/2}, \\
\hat{g}_{11}(-1,0,-1 \vert 0) \! = \! \mi l^{\prime}(g_{22} \! - \! g_{21}s_{1}^{\infty} 
\me^{2 \pi a} \! + \! s_{0}^{0}(g_{12} \! - \! g_{11}s_{1}^{\infty} \me^{2 \pi a})) 
\me^{-3 \pi a/4}, \\
\hat{g}_{12}(-1,0,-1 \vert 0) \! = \! -\mi l^{\prime}(g_{21} \! + \! s_{0}^{0}g_{11}) 
\me^{3 \pi a/4}, \quad \hat{g}_{21}(-1,0,-1 \vert 0) \! = \! \mi l^{\prime}(g_{12} 
\! - \! g_{11}s_{1}^{\infty} \me^{2 \pi a}) \me^{-3 \pi a/4}, \\
\hat{g}_{22}(-1,0,-1 \vert 0) \! = \! -\mi l^{\prime}g_{11} \me^{3 \pi a/4};
\end{gathered}
\end{equation}
\item[(22)] $\hat{\mathscr{F}}^{\lbrace 0 \rbrace}_{1,0,1} \Rightarrow$
\begin{equation} \label{laxhat97} 
\begin{gathered}
\hat{s}_{0}^{0}(1,0,1 \vert 0) \! = \! s_{0}^{0}, \quad \hat{s}_{0}^{\infty}(1,0,1 \vert 0) 
\! = \! -s_{1}^{\infty} \me^{3 \pi a/2}, \quad \hat{s}_{1}^{\infty}(1,0,1 \vert 0) \! = \! 
-s_{0}^{\infty} \me^{\pi a/2}, \\
\hat{g}_{11}(1,0,1 \vert 0) \! = \! \mi l^{\prime}g_{22} \me^{3 \pi a/4}, \quad 
\hat{g}_{12}(1,0,1 \vert 0) \! = \! -\mi l^{\prime}(g_{21} \! + \! s_{0}^{\infty}g_{22}) 
\me^{-3 \pi a/4}, \\
\hat{g}_{21}(1,0,1 \vert 0) \! = \! \mi l^{\prime}(g_{12} \! - \! s_{0}^{0}g_{22}) 
\me^{3 \pi a/4}, \quad \hat{g}_{22}(1,0,1 \vert 0) \! = \! \mi l^{\prime}(-g_{11} 
\! - \! s_{0}^{\infty}g_{12} \! + \! s_{0}^{0}(g_{21} \! + \! s_{0}^{\infty}g_{22})) 
\me^{-3 \pi a/4};
\end{gathered}
\end{equation}
\item[(23)] $\hat{\mathscr{F}}^{\lbrace 0 \rbrace}_{-1,0,1} \Rightarrow$
\begin{equation} \label{laxhat98} 
\begin{gathered}
\hat{s}_{0}^{0}(-1,0,1 \vert 0) \! = \! s_{0}^{0}, \quad \hat{s}_{0}^{\infty}(-1,0,1 \vert 0) 
\! = \! -s_{1}^{\infty} \me^{\pi a/2}, \quad \hat{s}_{1}^{\infty}(-1,0,1 \vert 0) \! = \! 
-s_{0}^{\infty} \me^{3 \pi a/2}, \\
\hat{g}_{11}(-1,0,1 \vert 0) \! = \! -l^{\prime}g_{12} \me^{\pi a/4}, \quad \hat{g}_{12}
(-1,0,1 \vert 0) \! = \! l^{\prime}(g_{11} \! + \! s_{0}^{\infty}g_{12}) \me^{-\pi a/4}, \\
\hat{g}_{21}(-1,0,1 \vert 0) \! = \! -l^{\prime}g_{22} \me^{\pi a/4}, \quad \hat{g}_{22}
(-1,0,1 \vert 0) \! = \! l^{\prime}(g_{21} \! + \! s_{0}^{\infty}g_{22}) \me^{-\pi a/4};
\end{gathered}
\end{equation}
\item[(24)] $\mathscr{F}^{\lbrace 1 \rbrace}_{0,0,0} \Rightarrow$
\begin{equation} \label{laxhat99} 
\begin{gathered}
s_{0}^{0}(0,0,0 \vert 1) \! = \! s_{0}^{0}, \quad s_{0}^{\infty}(0,0,0 \vert 1) \! = \! 
-s_{0}^{\infty}, \quad s_{1}^{\infty}(0,0,0 \vert 1) \! = \! -s_{1}^{\infty}, \\
g_{11}(0,0,0 \vert 1) \! = \! \mi \tilde{l}g_{11}, \quad g_{12}(0,0,0 \vert 1) \! = \! 
-\mi \tilde{l}g_{12}, \quad g_{21}(0,0,0 \vert 1) \! = \! \mi \tilde{l}g_{21}, \\
g_{22}(0,0,0 \vert 1) \! = \! -\mi \tilde{l}g_{22};
\end{gathered}
\end{equation}
\item[(25)] $\mathscr{F}^{\lbrace 1 \rbrace}_{-1,0,0} \Rightarrow$
\begin{equation} \label{laxhat100} 
\begin{gathered}
s_{0}^{0}(-1,0,0 \vert 1) \! = \! s_{0}^{0}, \quad s_{0}^{\infty}(-1,0,0 \vert 1) \! = \! 
-s_{0}^{\infty} \me^{\pi a}, \quad s_{1}^{\infty}(-1,0,0 \vert 1) \! = \! -s_{1}^{\infty} 
\me^{-\pi a}, \\
g_{11}(-1,0,0 \vert 1) \! = \! \tilde{l}(g_{21} \! + \! s_{0}^{0}g_{11}) \me^{\pi a/2}, 
\quad g_{12}(-1,0,0 \vert 1) \! = \! -\tilde{l}(g_{22} \! + \! s_{0}^{0}g_{12}) 
\me^{-\pi a/2}, \\
g_{21}(-1,0,0 \vert 1) \! = \! \tilde{l}g_{11} \me^{\pi a/2}, \quad g_{22}(-1,0,0 \vert 1) 
\! = \! -\tilde{l}g_{12} \me^{-\pi a/2};
\end{gathered}
\end{equation}
\item[(26)] $\mathscr{F}^{\lbrace 1 \rbrace}_{1,0,0} \Rightarrow$
\begin{equation} \label{laxhat101} 
\begin{gathered}
s_{0}^{0}(1,0,0 \vert 1) \! = \! s_{0}^{0}, \quad s_{0}^{\infty}(1,0,0 \vert 1) \! = \! 
-s_{0}^{\infty} \me^{-\pi a}, \quad s_{1}^{\infty}(1,0,0 \vert 1) \! = \! 
-s_{1}^{\infty} \me^{\pi a}, \\
g_{11}(1,0,0 \vert 1) \! = \! -\tilde{l}g_{21} \me^{-\pi a/2}, \quad g_{12}
(1,0,0 \vert 1) \! = \! \tilde{l}g_{22} \me^{\pi a/2}, \\
g_{21}(1,0,0 \vert 1) \! = \! -\tilde{l}(g_{11} \! - \! s_{0}^{0}g_{21}) \me^{-\pi a/2}, 
\quad g_{22}(1,0,0 \vert 1) \! = \! \tilde{l}(g_{12} \! - \! s_{0}^{0}g_{22}) \me^{\pi a/2};
\end{gathered}
\end{equation}
\item[(27)] $\mathscr{F}^{\lbrace 1 \rbrace}_{0,-1,-1} \Rightarrow$
\begin{equation} \label{laxhat102} 
\begin{gathered}
s_{0}^{0}(0,-1,-1 \vert 1) \! = \! s_{0}^{0}, \quad s_{0}^{\infty}(0,-1,-1 \vert 1) \! = \! 
s_{1}^{\infty} \me^{\pi a}, \quad s_{1}^{\infty}(0,-1,-1 \vert 1) \! = \! s_{0}^{\infty} 
\me^{\pi a}, \\
g_{11}(0,-1,-1 \vert 1) \! = \! -\tilde{l}l^{\prime}(g_{22} \! - \! g_{21}s_{1}^{\infty} 
\me^{2 \pi a} \! + \! s_{0}^{0}(g_{12} \! - \! g_{11}s_{1}^{\infty} \me^{2 \pi a})) 
\me^{-\pi a/2}, \\
g_{12}(0,-1,-1 \vert 1) \! = \! -\tilde{l}l^{\prime}(g_{21} \! + \! s_{0}^{0}g_{11}) 
\me^{\pi a/2}, \quad g_{21}(0,-1,-1 \vert 1) \! = \! -\tilde{l}l^{\prime}(g_{12} \!  
- \! g_{11}s_{1}^{\infty} \me^{2 \pi a}) \me^{-\pi a/2}, \\
g_{22}(0,-1,-1 \vert 1) \! = \! -\tilde{l}l^{\prime}g_{11} \me^{\pi a/2};
\end{gathered}
\end{equation}
\item[(28)] $\mathscr{F}^{\lbrace 1 \rbrace}_{0,-1,1} \Rightarrow$
\begin{equation} \label{laxhat103} 
\begin{gathered}
s_{0}^{0}(0,-1,1 \vert 1) \! = \! s_{0}^{0}, \quad s_{0}^{\infty}(0,-1,1 \vert 1) 
\! = \! s_{1}^{\infty} \me^{\pi a}, \quad s_{1}^{\infty}(0,-1,1 \vert 1) \! = \! 
s_{0}^{\infty} \me^{\pi a}, \\
g_{11}(0,-1,1 \vert 1) \! = \! \mi \tilde{l}l^{\prime}g_{12} \me^{\pi a/2}, \quad 
g_{12}(0,-1,1 \vert 1) \! = \! \mi \tilde{l}l^{\prime}(g_{11} \! + \! s_{0}^{\infty}
g_{12}) \me^{-\pi a/2}, \\ 
g_{21}(0,-1,1 \vert 1) \! = \! \mi \tilde{l}l^{\prime}g_{22} \me^{\pi a/2}, \quad 
g_{22}(0,-1,1 \vert 1) \! = \! \mi \tilde{l}l^{\prime}(g_{21} \! + \! s_{0}^{\infty}
g_{22}) \me^{-\pi a/2};
\end{gathered}
\end{equation}
\item[(29)] $\mathscr{F}^{\lbrace 1 \rbrace}_{0,1,-1} \Rightarrow$
\begin{equation} \label{laxhat104} 
\begin{gathered}
s_{0}^{0}(0,1,-1 \vert 1) \! = \! s_{0}^{0}, \quad s_{0}^{\infty}(0,1,-1 \vert 1) 
\! = \! s_{1}^{\infty} \me^{\pi a}, \quad s_{1}^{\infty}(0,1,-1 \vert 1) \! = \! 
s_{0}^{\infty} \me^{\pi a}, \\
g_{11}(0,1,-1 \vert 1) \! = \! \mi \tilde{l}l^{\prime}(g_{12} \! - \! g_{11}s_{1}^{\infty} 
\me^{2 \pi a}) \me^{-\pi a/2}, \quad g_{12}(0,1,-1 \vert 1) \! = \! \mi \tilde{l}
l^{\prime}g_{11} \me^{\pi a/2}, \\ 
g_{21}(0,1,-1 \vert 1) \! = \! \mi \tilde{l}l^{\prime}(g_{22} \! - \! g_{21}s_{1}^{\infty} 
\me^{2 \pi a}) \me^{-\pi a/2}, \quad g_{22}(0,1,-1 \vert 1) \! = \! \mi \tilde{l}
l^{\prime}g_{21} \me^{\pi a/2};
\end{gathered}
\end{equation}
\item[(30)] $\mathscr{F}^{\lbrace 1 \rbrace}_{0,1,1} \Rightarrow$
\begin{equation} \label{laxhat105} 
\begin{gathered}
s_{0}^{0}(0,1,1 \vert 1) \! = \! s_{0}^{0}, \quad s_{0}^{\infty}(0,1,1 \vert 1) \! = \! 
s_{1}^{\infty} \me^{\pi a}, \quad s_{1}^{\infty}(0,1,1 \vert 1) \! = \! s_{0}^{\infty} 
\me^{\pi a}, \\
g_{11}(0,1,1 \vert 1) \! = \! -\tilde{l}l^{\prime}g_{22} \me^{\pi a/2}, \quad g_{12}
(0,1,1 \vert 1) \! = \! -\tilde{l}l^{\prime}(g_{21} \! + \! s_{0}^{\infty}g_{22}) 
\me^{-\pi a/2}, \\
g_{21}(0,1,1 \vert 1) \! = \! -\tilde{l}l^{\prime}(g_{12} \! - \! s_{0}^{0}g_{22}) 
\me^{\pi a/2}, \quad g_{22}(0,1,1 \vert 1) \! = \! \tilde{l}l^{\prime}(-g_{11} \! 
- \! s_{0}^{\infty}g_{12} \! + \! s_{0}^{0}(g_{21} \! + \! s_{0}^{\infty}g_{22})) 
\me^{-\pi a/2};
\end{gathered}
\end{equation}
\item[(31)] $\mathscr{F}^{\lbrace 1 \rbrace}_{-1,-1,-1} \Rightarrow$
\begin{equation} \label{laxhat106} 
\begin{gathered}
s_{0}^{0}(-1,-1,-1 \vert 1) \! = \! s_{0}^{0}, \quad s_{0}^{\infty}(-1,-1,-1 \vert 1) 
\! = \! s_{1}^{\infty}, \quad s_{1}^{\infty}(-1,-1,-1 \vert 1) \! = \! s_{0}^{\infty} 
\me^{2 \pi a}, \\
g_{11}(-1,-1,-1 \vert 1) \! = \! \mi \tilde{l}l^{\prime}((g_{12} \! - \! g_{11}
s_{1}^{\infty} \me^{2 \pi a})(1 \! + \! (s_{0}^{0})^{2}) \! + \! s_{0}^{0}(g_{22} 
\! - \! g_{21}s_{1}^{\infty} \me^{2 \pi a})) \me^{-\pi a}, \\
g_{12}(-1,-1,-1 \vert 1) \! = \! \mi \tilde{l}l^{\prime}(g_{11}(1 \! + \! 
(s_{0}^{0})^{2}) \! + \! s_{0}^{0}g_{21}) \me^{\pi a}, \\
g_{21}(-1,-1,-1 \vert 1) \! = \! \mi \tilde{l}l^{\prime}(g_{22} \! - \! g_{21}
s_{1}^{\infty} \me^{2 \pi a} \! + \! s_{0}^{0}(g_{12} \! - \! g_{11}s_{1}^{\infty} 
\me^{2 \pi a})) \me^{-\pi a}, \\
g_{22}(-1,-1,-1 \vert 1) \! = \! \mi \tilde{l}l^{\prime}(g_{21} \! + \! s_{0}^{0}
g_{11}) \me^{\pi a};
\end{gathered}
\end{equation}
\item[(32)] $\mathscr{F}^{\lbrace 1 \rbrace}_{1,-1,-1} \Rightarrow$
\begin{equation} \label{laxhat107} 
\begin{gathered}
s_{0}^{0}(1,-1,-1 \vert 1) \! = \! s_{0}^{0}, \quad s_{0}^{\infty}(1,-1,-1 \vert 1) 
\! = \! s_{1}^{\infty} \me^{2 \pi a}, \quad s_{1}^{\infty}(1,-1,-1 \vert 1) \! = \! 
s_{0}^{\infty}, \\
g_{11}(1,-1,-1 \vert 1) \! = \! -\mi \tilde{l}l^{\prime}(g_{12} \! - \! g_{11}
s_{1}^{\infty} \me^{2 \pi a}), \quad g_{12}(1,-1,-1 \vert 1) \! = \! -\mi \tilde{l}
l^{\prime}g_{11}, \\
g_{21}(1,-1,-1 \vert 1) \! = \! -\mi \tilde{l}l^{\prime}(g_{22} \! - \! g_{21}
s_{1}^{\infty} \me^{2 \pi a}), \quad g_{22}(1,-1,-1 \vert 1) \! = \! -\mi \tilde{l}
l^{\prime}g_{21};
\end{gathered}
\end{equation}
\item[(33)] $\mathscr{F}^{\lbrace 1 \rbrace}_{-1,-1,1} \Rightarrow$
\begin{equation} \label{laxhat108} 
\begin{gathered}
s_{0}^{0}(-1,-1,1 \vert 1) \! = \! s_{0}^{0}, \quad s_{0}^{\infty}(-1,-1,1 \vert 1) 
\! = \! s_{1}^{\infty}, \quad s_{1}^{\infty}(-1,-1,1 \vert 1) \! = \! s_{0}^{\infty} 
\me^{2 \pi a}, \\
g_{11}(-1,-1,1 \vert 1) \! = \! \tilde{l}l^{\prime}(g_{22} \! + \! s_{0}^{0}g_{12}), 
\quad g_{12}(-1,-1,1 \vert 1) \! = \! \tilde{l}l^{\prime}(g_{21} \! + \! s_{0}^{\infty}
g_{22} \! + \! s_{0}^{0}(g_{11} \! + \! s_{0}^{\infty}g_{12})), \\
g_{21}(-1,-1,1 \vert 1) \! = \! \tilde{l}l^{\prime}g_{12}, \quad g_{22}
(-1,-1,1 \vert 1) \! = \! \tilde{l}l^{\prime}(g_{11} \! + \! s_{0}^{\infty}g_{12});
\end{gathered}
\end{equation}
\item[(34)] $\mathscr{F}^{\lbrace 1 \rbrace}_{1,-1,1} \Rightarrow$
\begin{equation} \label{laxhat109} 
\begin{gathered}
s_{0}^{0}(1,-1,1 \vert 1) \! = \! s_{0}^{0}, \quad s_{0}^{\infty}(1,-1,1 \vert 1) 
\! = \! s_{1}^{\infty} \me^{2 \pi a}, \quad s_{1}^{\infty}(1,-1,1 \vert 1) \! = \! 
s_{0}^{\infty}, \\
g_{11}(1,-1,1 \vert 1) \! = \! -\tilde{l}l^{\prime}g_{22} \me^{\pi a}, \quad 
g_{12}(1,-1,1 \vert 1) \! = \! -\tilde{l}l^{\prime}(g_{21} \! + \! s_{0}^{\infty}
g_{22}) \me^{-\pi a}, \\
g_{21}(1,-1,1 \vert 1) \! = \! -\tilde{l}l^{\prime}(g_{12} \!  - \! s_{0}^{0}g_{22}) 
\me^{\pi a}, \\
g_{22}(1,-1,1 \vert 1) \! = \! -\tilde{l}l^{\prime}(g_{11} \! + \! s_{0}^{\infty}
g_{12} \! - \! s_{0}^{0}(g_{21} \! + \! s_{0}^{\infty}g_{22})) \me^{-\pi a};
\end{gathered}
\end{equation}
\item[(35)] $\mathscr{F}^{\lbrace 1 \rbrace}_{-1,1,-1} \Rightarrow$
\begin{equation} \label{laxhat110} 
\begin{gathered}
s_{0}^{0}(-1,1,-1 \vert 1) \! = \! s_{0}^{0}, \quad s_{0}^{\infty}(-1,1,-1 \vert 1) 
\! = \! s_{1}^{\infty}, \quad s_{1}^{\infty}(-1,1,-1 \vert 1) \! = \! s_{0}^{\infty} 
\me^{2 \pi a}, \\
g_{11}(-1,1,-1 \vert 1) \! = \! \tilde{l}l^{\prime}(g_{22} \! - \! g_{21}s_{1}^{\infty} 
\me^{2 \pi a} \! + \! s_{0}^{0}(g_{12} \! - \! g_{11}s_{1}^{\infty} \me^{2 \pi a})) 
\me^{-\pi a}, \\
g_{12}(-1,1,-1 \vert 1) \! = \! \tilde{l}l^{\prime}(g_{21} \! + \! s_{0}^{0}g_{11}) 
\me^{\pi a}, \quad g_{21}(-1,1,-1 \vert 1) \! = \! \tilde{l}l^{\prime}
(g_{12} \!  - \! g_{11}s_{1}^{\infty} \me^{2 \pi a}) \me^{-\pi a}, \\
g_{22}(-1,1,-1 \vert 1) \! = \! \tilde{l}l^{\prime}g_{11} \me^{\pi a};
\end{gathered}
\end{equation}
\item[(36)] $\mathscr{F}^{\lbrace 1 \rbrace}_{1,1,-1} \Rightarrow$
\begin{equation} \label{laxhat111} 
\begin{gathered}
s_{0}^{0}(1,1,-1 \vert 1) \! = \! s_{0}^{0}, \quad s_{0}^{\infty}(1,1,-1 \vert 1) 
\! = \! s_{1}^{\infty} \me^{2 \pi a}, \quad s_{1}^{\infty}(1,1,-1 \vert 1) \! = \! 
s_{0}^{\infty}, \\
g_{11}(1,1,-1 \vert 1) \! = \! -\tilde{l}l^{\prime}(g_{22} \! - \! g_{21}s_{1}^{\infty} 
\me^{2 \pi a}), \quad g_{12}(1,1,-1 \vert 1) \! = \! -\tilde{l}l^{\prime}g_{21}, \\ 
g_{21}(1,1,-1 \vert 1) \! = \! -\tilde{l}l^{\prime}(g_{12} \! - \! g_{11}s_{1}^{\infty} 
\me^{2 \pi a} \! - \! s_{0}^{0}(g_{22} \! - \! g_{21}s_{1}^{\infty} \me^{2 \pi a})), \\
g_{22}(1,1,-1 \vert 1) \! = \! -\tilde{l}l^{\prime}(g_{11} \!  - \! s_{0}^{0}g_{21});
\end{gathered}
\end{equation}
\item[(37)] $\mathscr{F}^{\lbrace 1 \rbrace}_{-1,1,1} \Rightarrow$
\begin{equation} \label{laxhat112} 
\begin{gathered}
s_{0}^{0}(-1,1,1 \vert 1) \! = \! s_{0}^{0}, \quad s_{0}^{\infty}(-1,1,1 \vert 1) 
\! = \! s_{1}^{\infty}, \quad s_{1}^{\infty}(-1,1,1 \vert 1) \! = \! s_{0}^{\infty} 
\me^{2 \pi a}, \\
g_{11}(-1,1,1 \vert 1) \! = \! \mi \tilde{l}l^{\prime}g_{12}, \quad g_{12}
(-1,1,1 \vert 1) \! = \! \mi \tilde{l}l^{\prime}(g_{11} \!  + \! s_{0}^{\infty}
g_{12}), \\
g_{21}(-1,1,1 \vert 1) \! = \! \mi \tilde{l}l^{\prime}g_{22}, \quad g_{22}
(-1,1,1 \vert 1) \! = \! \mi \tilde{l}l^{\prime}(g_{21} \!  + \! s_{0}^{\infty}
g_{22});
\end{gathered}
\end{equation}
\item[(38)] $\mathscr{F}^{\lbrace 1 \rbrace}_{1,1,1} \Rightarrow$
\begin{equation} \label{laxhat113} 
\begin{gathered}
s_{0}^{0}(1,1,1 \vert 1) \! = \! s_{0}^{0}, \quad s_{0}^{\infty}(1,1,1 \vert 1) \! = \! 
s_{1}^{\infty} \me^{2 \pi a}, \quad s_{1}^{\infty}(1,1,1 \vert 1) \! = \! s_{0}^{\infty}, \\
g_{11}(1,1,1 \vert 1) \! = \! -\mi \tilde{l}l^{\prime}(g_{12} \! - \! s_{0}^{0}g_{22}) 
\me^{\pi a}, \quad g_{12}(1,1,1 \vert 1) \! = \! \mi \tilde{l}l^{\prime}(-g_{11} 
\! - \! s_{0}^{\infty}g_{12} \! + \! s_{0}^{0}(g_{21} \! + \! s_{0}^{\infty}g_{22})) 
\me^{-\pi a}, \\
g_{21}(1,1,1 \vert 1) \! = \! -\mi \tilde{l}l^{\prime}(g_{22} \! - \! s_{0}^{0}
(g_{12} \! - \! s_{0}^{0}g_{22})) \me^{\pi a}, \\
g_{22}(1,1,1 \vert 1) \! = \! -\mi \tilde{l}l^{\prime}((g_{21} \! + \! s_{0}^{\infty}
g_{22})(1 \! + \! (s_{0}^{0})^{2}) \! - \! s_{0}^{0}(g_{11} \! + \! s_{0}^{\infty}g_{12})) 
\me^{-\pi a};
\end{gathered}
\end{equation}
\item[(39)] $\hat{\mathscr{F}}^{\lbrace 1 \rbrace}_{1,1,0} \Rightarrow$
\begin{equation} \label{laxhat114} 
\begin{gathered}
\hat{s}_{0}^{0}(1,1,0 \vert 1) \! = \! s_{0}^{0}, \quad \hat{s}_{0}^{\infty}(1,1,0 \vert 1) 
\! = \! -s_{0}^{\infty} \me^{-\pi a/2}, \quad \hat{s}_{1}^{\infty}(1,1,0 \vert 1) \! = \! 
-s_{1}^{\infty} \me^{\pi a/2}, \\
\hat{g}_{11}(1,1,0 \vert 1) \! = \! \tilde{l}g_{21} \me^{-\pi a/4}, \quad \hat{g}_{12}
(1,1,0 \vert 1) \! = \! -\tilde{l}g_{22} \me^{\pi a/4}, \\
\hat{g}_{21}(1,1,0 \vert 1) \! = \! \tilde{l}(g_{11} \! - \! s_{0}^{0}g_{21}) \me^{-\pi a/4}, 
\quad \hat{g}_{22}(1,1,0 \vert 1) \! = \! -\tilde{l}(g_{12} \! - \! s_{0}^{0}g_{22}) 
\me^{\pi a/4};
\end{gathered}
\end{equation}
\item[(40)] $\hat{\mathscr{F}}^{\lbrace 1 \rbrace}_{1,-1,0} \Rightarrow$
\begin{equation} \label{laxhat115} 
\begin{gathered}
\hat{s}_{0}^{0}(1,-1,0 \vert 1) \! = \! s_{0}^{0}, \quad \hat{s}_{0}^{\infty}(1,-1,0 \vert 1) 
\! = \! -s_{0}^{\infty} \me^{-\pi a/2}, \quad \hat{s}_{1}^{\infty}(1,-1,0 \vert 1) \! = \! 
-s_{1}^{\infty} \me^{\pi a/2}, \\
\hat{g}_{11}(1,-1,0 \vert 1) \! = \! \mi \tilde{l}g_{11} \me^{-\pi a/4}, \quad \hat{g}_{12}
(1,-1,0 \vert 1) \! = \! -\mi \tilde{l}g_{12} \me^{\pi a/4}, \\
\hat{g}_{21}(1,-1,0 \vert 1) \! = \! \mi \tilde{l}g_{21} \me^{-\pi a/4}, \quad \hat{g}_{22}
(1,-1,0 \vert 1) \! = \! -\mi \tilde{l}g_{22} \me^{\pi a/4};
\end{gathered}
\end{equation}
\item[(41)] $\hat{\mathscr{F}}^{\lbrace 1 \rbrace}_{-1,1,0} \Rightarrow$
\begin{equation} \label{laxhat116} 
\begin{gathered}
\hat{s}_{0}^{0}(-1,1,0 \vert 1) \! = \! s_{0}^{0}, \quad \hat{s}_{0}^{\infty}(-1,1,0 \vert 1) 
\! = \! -s_{0}^{\infty} \me^{\pi a/2}, \quad \hat{s}_{1}^{\infty}(-1,1,0 \vert 1) \! = \! 
-s_{1}^{\infty} \me^{-\pi a/2}, \\
\hat{g}_{11}(-1,1,0 \vert 1) \! = \! \mi \tilde{l}g_{11} \me^{\pi a/4}, \quad \hat{g}_{12}
(-1,1,0 \vert 1) \! = \! -\mi \tilde{l}g_{12} \me^{-\pi a/4}, \\
\hat{g}_{21}(-1,1,0 \vert 1) \! = \! \mi \tilde{l}g_{21} \me^{\pi a/4}, \quad \hat{g}_{22}
(-1,1,0 \vert 1) \! = \! -\mi \tilde{l}g_{22} \me^{-\pi a/4};
\end{gathered}
\end{equation}
\item[(42)] $\hat{\mathscr{F}}^{\lbrace 1 \rbrace}_{-1,-1,0} \Rightarrow$
\begin{equation} \label{laxhat117} 
\begin{gathered}
\hat{s}_{0}^{0}(-1,-1,0 \vert 1) \! = \! s_{0}^{0}, \quad \hat{s}_{0}^{\infty}(-1,-1,0 \vert 1) 
\! = \! -s_{0}^{\infty} \me^{\pi a/2}, \quad \hat{s}_{1}^{\infty}(-1,-1,0 \vert 1) \! = \! 
-s_{1}^{\infty} \me^{-\pi a/2}, \\
\hat{g}_{11}(-1,-1,0 \vert 1) \! = \! -\tilde{l}(g_{21} \! + \! s_{0}^{0}g_{11}) \me^{\pi a/4}, 
\quad \hat{g}_{12}(-1,-1,0 \vert 1) \! = \! \tilde{l}(g_{22} \! + \! s_{0}^{0}g_{12}) 
\me^{-\pi a/4}, \\
\hat{g}_{21}(-1,-1,0 \vert 1) \! = \! -\tilde{l}g_{11} \me^{\pi a/4}, \quad \hat{g}_{22}
(-1,-1,0 \vert 1) \! = \! \tilde{l}g_{12} \me^{-\pi a/4};
\end{gathered}
\end{equation}
\item[(43)] $\hat{\mathscr{F}}^{\lbrace 1 \rbrace}_{1,0,-1} \Rightarrow$
\begin{equation} \label{laxhat118} 
\begin{gathered}
\hat{s}_{0}^{0}(1,0,-1 \vert 1) \! = \! s_{0}^{0}, \quad \hat{s}_{0}^{\infty}(1,0,-1 \vert 1) 
\! = \! s_{1}^{\infty} \me^{3 \pi a/2}, \quad \hat{s}_{1}^{\infty}(1,0,-1 \vert 1) \! = \! 
s_{0}^{\infty} \me^{\pi a/2}, \\
\hat{g}_{11}(1,0,-1 \vert 1) \! = \! \mi \tilde{l}l^{\prime}(g_{12} \!  - \! g_{11}s_{1}^{\infty} 
\me^{2 \pi a}) \me^{-\pi a/4}, \quad \hat{g}_{12}(1,0,-1 \vert 1) \! = \! \mi \tilde{l}
l^{\prime}g_{11} \me^{\pi a/4}, \\
\hat{g}_{21}(1,0,-1 \vert 1) \! = \! \mi \tilde{l}l^{\prime}(g_{22} \! - \! g_{21}s_{1}^{\infty} 
\me^{2 \pi a}) \me^{-\pi a/4}, \quad \hat{g}_{22}(1,0,-1 \vert 1) \! = \! \mi \tilde{l}
l^{\prime}g_{21} \me^{\pi a/4};
\end{gathered}
\end{equation}
\item[(44)] $\hat{\mathscr{F}}^{\lbrace 1 \rbrace}_{1,0,1} \Rightarrow$
\begin{equation} \label{laxhat119} 
\begin{gathered}
\hat{s}_{0}^{0}(1,0,1 \vert 1) \! = \! s_{0}^{0}, \quad \hat{s}_{0}^{\infty}(1,0,1 \vert 1) 
\! = \! s_{1}^{\infty} \me^{3 \pi a/2}, \quad \hat{s}_{1}^{\infty}(1,0,1 \vert 1) \! = \! 
s_{0}^{\infty} \me^{\pi a/2}, \\
\hat{g}_{11}(1,0,1 \vert 1) \! = \! -\tilde{l}l^{\prime}g_{22} \me^{3 \pi a/4}, \quad 
\hat{g}_{12}(1,0,1 \vert 1) \! = \! -\tilde{l}l^{\prime}(g_{21} \! + \! g_{22}s_{0}^{\infty}) 
\me^{-3 \pi a/4}, \\
\hat{g}_{21}(1,0,1 \vert 1) \! = \! -\tilde{l}l^{\prime}(g_{12} \! - \! s_{0}^{0}g_{22}) 
\me^{3 \pi a/4}, \\
\hat{g}_{22}(1,0,1 \vert 1) \! = \! \tilde{l}l^{\prime}(-g_{11} \! - \! s_{0}^{\infty}
g_{12} \! + \! s_{0}^{0}(g_{21} \! + \! s_{0}^{\infty}g_{22})) \me^{-3 \pi a/4};
\end{gathered}
\end{equation}
\item[(45)] $\hat{\mathscr{F}}^{\lbrace 1 \rbrace}_{-1,0,-1} \Rightarrow$
\begin{equation} \label{laxhat120} 
\begin{gathered}
\hat{s}_{0}^{0}(-1,0,-1 \vert 1) \! = \! s_{0}^{0}, \quad \hat{s}_{0}^{\infty}(-1,0,-1 \vert 1) 
\! = \! s_{1}^{\infty} \me^{\pi a/2}, \quad \hat{s}_{1}^{\infty}(-1,0,-1 \vert 1) \! = \! 
s_{0}^{\infty} \me^{3 \pi a/2}, \\
\hat{g}_{11}(-1,0,-1 \vert 1) \! = \! -\tilde{l}l^{\prime}(g_{22} \! - \! g_{21}s_{1}^{\infty} 
\me^{2 \pi a} \! + \! s_{0}^{0}(g_{12} \! - \! g_{11}s_{1}^{\infty} \me^{2 \pi a})) 
\me^{-3 \pi a/4}, \\
\hat{g}_{12}(-1,0,-1 \vert 1) \! = \! -\tilde{l}l^{\prime}(g_{21} \! + \! s_{0}^{0}g_{11}) 
\me^{3 \pi a/4}, \quad \hat{g}_{21}(-1,0,-1 \vert 1) \! = \! -\tilde{l}l^{\prime}
(g_{12} \! - \! g_{11}s_{1}^{\infty} \me^{2 \pi a}) \me^{-3 \pi a/4}, \\
\hat{g}_{22}(-1,0,-1 \vert 1) \! = \! -\tilde{l}l^{\prime}g_{11} \me^{3 \pi a/4};
\end{gathered}
\end{equation}
\item[(46)] $\hat{\mathscr{F}}^{\lbrace 1 \rbrace}_{-1,0,1} \Rightarrow$
\begin{equation} \label{laxhat121} 
\begin{gathered}
\hat{s}_{0}^{0}(-1,0,1 \vert 1) \! = \! s_{0}^{0}, \quad \hat{s}_{0}^{\infty}(-1,0,1 \vert 1) 
\! = \! s_{1}^{\infty} \me^{\pi a/2}, \quad \hat{s}_{1}^{\infty}(-1,0,1 \vert 1) \! = \! 
s_{0}^{\infty} \me^{3 \pi a/2}, \\
\hat{g}_{11}(-1,0,1 \vert 1) \! = \! -\mi \tilde{l}l^{\prime}g_{12} \me^{\pi a/4}, \quad 
\hat{g}_{12}(-1,0,1 \vert 1) \! = \! -\mi \tilde{l}l^{\prime}(g_{11} \! + \! s_{0}^{\infty}
g_{12}) \me^{-\pi a/4}, \\
\hat{g}_{21}(-1,0,1 \vert 1) \! = \! -\mi \tilde{l}l^{\prime}g_{22} \me^{\pi a/4}, \quad 
\hat{g}_{22}(-1,0,1 \vert 1) \! = \! -\mi \tilde{l}l^{\prime}(g_{21} \! + \! s_{0}^{\infty}
g_{22}) \me^{-\pi a/4}.
\end{gathered}
\end{equation}
\end{enumerate}

Finally, applying the isomorphism $\mathscr{F}^{\scriptscriptstyle \lbrace \ell \rbrace}_{\scriptscriptstyle 
\varepsilon_{1},\varepsilon_{2},m(\varepsilon_{2})}$ (resp., $\hat{\mathscr{F}}^{\scriptscriptstyle \lbrace 
\hat{\ell} \rbrace}_{\scriptscriptstyle \hat{\varepsilon}_{1},\hat{\varepsilon}_{2},\hat{m}(\hat{\varepsilon}_{2})}$), 
whose action on $\mathscr{M}$ is given by equations \eqref{laxhat76}--\eqref{laxhat90} and 
\eqref{laxhat99}--\eqref{laxhat113} (resp., equations \eqref{laxhat91}--\eqref{laxhat98} and 
\eqref{laxhat114}--\eqref{laxhat121}), to the corresponding $(\varepsilon_{1},\varepsilon_{2},m(\varepsilon_{2}) 
\vert \ell) \! = \! (0,0,0 \vert 0)$ (resp., $(\hat{\varepsilon}_{1},\hat{\varepsilon}_{2},\hat{m}(\hat{\varepsilon}_{2}) 
\vert \hat{\ell}) \! = \! (0,0,0 \vert 0)$) asymptotics (as $\tau \! \to \! +\infty$ with $\varepsilon b \! > \! 0$) for 
$u(\tau)$, $f_{\pm}(\tau)$, $\mathcal{H}(\tau)$, and $\sigma (\tau)$ derived in Section \ref{finalsec}, one arrives 
at the asymptotics as $\tau \! \to \! \pm \infty$ (resp., $\tau \! \to \! \pm \mi \infty$) for $u(\tau)$, $f_{\pm}(\tau)$, 
$\mathcal{H}(\tau)$, and $\sigma (\tau)$ stated in Theorem \ref{theor2.1} (resp., Theorem \ref{appen}).\footnote{In 
Section 3 (resp., Section 2), p. 1174 (resp., p. 7) of \cite{a1} (resp., \cite{av2}), for item (9) in the definition of the 
mapping $\mathscr{F}_{1,1}$, the formula for $g_{21}(1,1)$ is missing: it reads $g_{21}(1,1) = \mi g_{12} \me^{\pi a}$.}
\appendix
\setcounter{section}{4} 
\section{Appendix: Asymptotics of $\hat{\varphi}(\tau)$ as $\tau \! \to \! \pm \infty$ and $\tau \! \to \! \pm \mi \infty$} 
\label{feetics} 
In this appendix, asymptotics as $\tau \! \to \! \pm \infty$ (resp., $\tau \! \to \! \pm \mi \infty$) for $\pm \varepsilon b 
\! > \! 0$ of the function $\hat{\varphi}(\tau)$ (cf. Proposition \ref{prop1.2}) are presented in Theorem \ref{pfeetotsa} 
(resp., Theorem \ref{pfeetotsb}). The results of this appendix are seminal for an upcoming series of works on 
asymptotics of integrals of solutions to the {\rm DP3E} \eqref{eq1.1} and related functions.
\begin{eeeee} \label{aboutfemoed} 
Since the function $\hat{\varphi}(\tau)$ is defined $\operatorname{mod}(2 \pi)$, the reader should be cognizant of the 
fact that the asymptotics for $\hat{\varphi}(\tau)$ stated in Theorems \ref{pfeetotsa} and \ref{pfeetotsb} are defined 
$\operatorname{mod}(2 \pi)$; this $\operatorname{mod}(2 \pi)$ arbitrariness, however, is not important, because the 
requisite functions are $u(\tau)$ and $\exp (\mi \hat{\varphi}(\tau))$.\hfill $\blacksquare$
\end{eeeee} 
\begin{eeeee} \label{noteaboutrealvarp} 
If one is only interested in the asymptotics as $\tau \! \to \! +\infty$ for $\varepsilon b \! > \! 0$ of the function 
$\hat{\varphi}(\tau)$, then, in Theorem \ref{pfeetotsa}, one sets $(\varepsilon_{1},\varepsilon_{2},
m(\varepsilon_{2}) \vert \ell) \! = \! (0,0,0 \vert 0)$ and uses the fact that (cf. Appendix \ref{sectonsymmcomp}, 
the identity map \eqref{laxhat76}) $s_{0}^{0}(0,0,0 \vert 0) \! = \! s_{0}^{0}$, $s_{0}^{\infty}(0,0,0 \vert 0) \! = \! 
s_{0}^{\infty}$, $s_{1}^{\infty}(0,0,0 \vert 0) \! = \! s_{1}^{\infty}$, and $g_{ij}(0,0,0 \vert 0) \! = \! g_{ij}$, $i,j \! \in 
\! \lbrace 1,2 \rbrace$. \hfill $\blacksquare$
\end{eeeee}
\begin{ddddd} \label{pfeetotsa} 
Let $u(\tau)$ be a solution of the {\rm DP3E} \eqref{eq1.1} and $\hat{\varphi}(\tau)$ be the general solution of the 
{\rm ODE} $\hat{\varphi}^{\prime}(\tau) \! = \! 2a \tau^{-1} \! + \! b(u(\tau))^{-1}$ for $\varepsilon b \! > \! 0$ 
corresponding to the monodromy data $(a,s^{0}_{0},s^{\infty}_{0},s^{\infty}_{1},g_{11},g_{12},g_{21},g_{22})$. Let 
$\varepsilon_{1},\varepsilon_{2} \! \in  \! \lbrace 0,\pm 1 \rbrace$, 
$m(\varepsilon_{2}) \! = \! 
\left\{
\begin{smallmatrix}
0, \, \, \varepsilon_{2}=0, \\
\pm \varepsilon_{2}, \, \, \varepsilon_{2} \in \lbrace \pm 1 \rbrace,
\end{smallmatrix} 
\right.$ $\ell \! \in \! \lbrace 0,1 \rbrace$, and $\varepsilon b \! = \! \vert \varepsilon b \vert \me^{\mi \pi \varepsilon_{2}}$. 
For $k \! = \! +1$, let
\begin{equation*}
g_{11}(\varepsilon_{1},\varepsilon_{2},m(\varepsilon_{2}) \vert \ell)g_{12}
(\varepsilon_{1},\varepsilon_{2},m(\varepsilon_{2}) \vert \ell)g_{21}(\varepsilon_{1},
\varepsilon_{2},m(\varepsilon_{2}) \vert \ell) \! \neq \! 0 \quad \text{and} \quad 
g_{22}(\varepsilon_{1},\varepsilon_{2},m(\varepsilon_{2}) \vert \ell) \! = \! 0,
\end{equation*}
and, for $k \! = \! -1$, let
\begin{equation*}
g_{11}(\varepsilon_{1},\varepsilon_{2},m(\varepsilon_{2}) \vert \ell) \! = \! 0 \quad 
\text{and} \quad g_{12}(\varepsilon_{1},\varepsilon_{2},m(\varepsilon_{2}) \vert 
\ell)g_{21}(\varepsilon_{1},\varepsilon_{2},m(\varepsilon_{2}) \vert \ell)g_{22}
(\varepsilon_{1},\varepsilon_{2},m(\varepsilon_{2}) \vert \ell) \! \neq \! 0,
\end{equation*}
where explicit expressions for $g_{ij}(\varepsilon_{1},\varepsilon_{2},m(\varepsilon_{2}) \vert \ell)$, $i,j \! \in \! 
\lbrace 1,2 \rbrace$, are given in Appendix \ref{sectonsymm}, equations \eqref{laxhat76}--\eqref{laxhat90} and 
\eqref{laxhat99}--\eqref{laxhat113}. Then, for $s_{0}^{0}(\varepsilon_{1},\varepsilon_{2},m(\varepsilon_{2}) \vert 
\ell) \! \neq \! \mi \me^{(-1)^{1+\varepsilon_{2}} \pi a}$,\footnote{Recall that (cf. Remark \ref{stokmxx}) $s_{0}^{0}
(\varepsilon_{1},\varepsilon_{2},m(\varepsilon_{2}) \vert \ell) \! = \! s_{0}^{0}$. For $s_{0}^{0}(\varepsilon_{1},
\varepsilon_{2},m(\varepsilon_{2}) \vert \ell) \! = \! \mi \me^{(-1)^{1+\varepsilon_{2}} \pi a}$, the exponentially 
small correction term in the asymptotics \eqref{totsa1} is absent.}
\begin{align} \label{totsa1} 
(-1)^{\varepsilon_{2}} \hat{\varphi}(\tau) \underset{\tau \to +\infty \me^{\mi \pi \varepsilon_{1}}}{=}& 
\, 3 \me^{\mi 2 \pi k/3}(-1)^{\varepsilon_{2}}(\varepsilon b)^{1/3} \tau^{2/3} \! + \! 
2(-1)^{\varepsilon_{2}}a \ln \! \left(\frac{2 \me^{-\mi \pi k/3} \tau^{2/3}}{(\varepsilon b 
\me^{-\mi \pi \varepsilon_{2}})^{1/6}} \right) \! + \! \mi \mathsterling_{k}(\varepsilon_{1},
\varepsilon_{2},m(\varepsilon_{2}) \vert \ell) \nonumber \\
-& \, k \pi \! - \! \mi \, \mathlarger{\sum_{m=2}^{\infty}} \! \left(2 \widetilde{\nu}_{m}(k) \! + \! 
\sum_{\substack{n,\mathfrak{l} \in \mathbb{N}\\\mathfrak{l} \geqslant n\\n+ \mathfrak{l}=m}} \, 
\sum_{\substack{\mathfrak{i}_{1}+2 \mathfrak{i}_{2} + \dotsb + \mathfrak{l} \mathfrak{i}_{
\mathfrak{l}}=\mathfrak{l}\\\mathfrak{i}_{1} + \mathfrak{i}_{2} + \dotsb + 
\mathfrak{i}_{\mathfrak{l}}=n}} \, (-1)^{n-1}(n \! - \! 1)! \prod_{\mathfrak{j}=1}^{\mathfrak{l}} 
\frac{(\mathfrak{u}_{\mathfrak{j}-1}(k))^{\mathfrak{i}_{\mathfrak{j}}}}{\mathfrak{i}_{\mathfrak{j}}!} 
\right) \nonumber \\
\times& \, \big((-1)^{\varepsilon_{1}} \tau^{-1/3} \big)^{m} \! - \! \frac{k(-1)^{\varepsilon_{1}} 
\me^{-\mi \pi k/3} \me^{\mi \pi k/4}(2 \! + \! \sqrt{3})^{\mi k(-1)^{\varepsilon_{2}}a}}{\sqrt{2 \pi} 
\, 3^{3/4}(\varepsilon b \me^{-\mi \pi \varepsilon_{2}})^{1/6} \tau^{1/3}} \! \left(
\vphantom{M^{M^{M}}}s_{0}^{0}(\varepsilon_{1},\varepsilon_{2},m(\varepsilon_{2}) \vert \ell) 
\right. \nonumber \\
-&\left. \, \mi \me^{(-1)^{1+\varepsilon_{2}} \pi a} \right) \! \me^{-\frac{3 \sqrt{3}}{2}(\sqrt{3}+\mi k)(-1)^{\varepsilon_{2}}
(\varepsilon b)^{1/3} \tau^{2/3}} \! \left(1 \! + \! \mathcal{O} \big(\tau^{-1/3} \big) \right), \quad k \! \in \! \lbrace \pm 1 \rbrace,
\end{align}
where
\begin{equation} \label{totsa2} 
\mathsterling_{k}(\varepsilon_{1},\varepsilon_{2},m(\varepsilon_{2}) \vert \ell) \! = \! 
\begin{cases}
\ln \! \left(g_{11}(\varepsilon_{1},\varepsilon_{2},m(\varepsilon_{2}) \vert \ell) 
\me^{(-1)^{\varepsilon_{2}} \pi a} \right)^{2}, &\text{$k \! = \! +1$,} \\
\ln \! \left(g_{22}(\varepsilon_{1},\varepsilon_{2},m(\varepsilon_{2}) \vert \ell) 
\me^{(-1)^{\varepsilon_{2}} \pi a} \right)^{-2}, &\text{$k \! = \! -1$,}
\end{cases}
\end{equation}
\begin{equation} \label{totsa3}
\begin{gathered}
\widetilde{\nu}_{1}(k) \! = \! 0, \qquad \qquad \widetilde{\nu}_{2}(k) \! = \! \frac{a(1 \! + \! 
\mi (-1)^{\varepsilon_{2}}a) \me^{\mi \pi k/3}}{6(\varepsilon b)^{1/3}}, \qquad \qquad 
\widetilde{\nu}_{3}(k) \! = \! 0, \\
\widetilde{\nu}_{4}(k) \! = \! -\frac{\mi (-1)^{\varepsilon_{2}}a \me^{\mi 2 \pi k/3}}{36(\varepsilon b)^{2/3}} 
\! \left(\frac{1 \! - \! 2a^{2}}{3} \! + \! \mi (-1)^{\varepsilon_{2}}a \right),
\end{gathered}
\end{equation}
and
\begin{align} \label{totsa4} 
(m \! + \! 5) \widetilde{\nu}_{m+5}(k) \! =& \, \frac{3 \mi}{2} \me^{-\mi \pi k/3}(-1)^{\varepsilon_{2}}
(\varepsilon b)^{1/3} \mathfrak{u}_{m+5}(k) \! + \! \frac{\mi (-1)^{\varepsilon_{2}} \me^{\mi \pi k/3}
(1 \! + \! 2 \mi (-1)^{\varepsilon_{2}}a)}{12(\varepsilon b)^{1/3}} \mu_{m+1}^{\ast}(k) \! + \! 
\frac{1}{4} \mu_{m+3}^{\ast}(k) \nonumber \\
-& \, \frac{\mi (-1)^{\varepsilon_{2}} \me^{\mi \pi k/3}}{12(\varepsilon b)^{1/3}} \! \left(
\vphantom{M^{M^{M^{M}}}}(m \! + \! 3)(m \! + \! 5 \! + \! 2 \mi (-1)^{\varepsilon_{2}}a) 
\widetilde{\nu}_{m+3}(k) \! - \! \frac{\mi (-1)^{\varepsilon_{2}}2a^{2} \me^{\mi \pi k/3}}{3
(\varepsilon b)^{1/3}}(m \! + \! 1) \widetilde{\nu}_{m+1}(k) \right. \nonumber \\
+&\left. \, \sum_{j=0}^{m-1}(j \! + \! 1) \widetilde{\nu}_{j+1}(k)(\mu_{m-j}^{\ast}(k) \! - \! 
2(m \! + \! 2 \! - \! j) \widetilde{\nu}_{m+2-j}(k)) \right), \quad m \! \in \! \mathbb{Z}_{+},
\end{align} 
with
\begin{equation} \label{totsa5} 
\begin{gathered} 
\mu_{0}^{\ast}(k) \! = \! \frac{2a \me^{\mi \pi k/3}}{3(\varepsilon b)^{1/3}}, \qquad 
\quad \mu_{1}^{\ast}(k) \! = \! 0, \\
\mu_{m_{1}+2}^{\ast}(k) \! = \! -2 \! \left(\mathrm{P}_{m_{1}+2}^{\ast}(k) \! + \! 
\mathfrak{w}_{m_{1}+2}(k) \! + \! \sum_{j=0}^{m_{1}} \mathrm{P}_{j}^{\ast}(k) 
\mathfrak{w}_{m_{1}-j}(k) \right), \quad m_{1} \! \in \! \mathbb{Z}_{+},
\end{gathered}
\end{equation}
and
\begin{equation} \label{totsa6} 
\begin{gathered} 
\mathrm{P}_{0}^{\ast}(k) \! = \! -\frac{2a \me^{\mi \pi k/3}}{3(\varepsilon b)^{1/3}}, 
\quad \qquad \mathrm{P}_{1}^{\ast}(k) \! = \! 0, \\
\mathrm{P}_{j}^{\ast}(k) \! = \! \frac{3}{2} \! \left(\mathfrak{u}_{j}(k) \! - \! \mi 
(-1)^{\varepsilon_{2}} \me^{\mi 2 \pi k/3}(\varepsilon b)^{1/3} \! \left(\mathfrak{r}_{j+2}
(k) \! - \! 2 \mathfrak{u}_{j+2}(k) \! + \! \sum_{m_{2}=0}^{j} \mathfrak{u}_{m_{2}}(k) 
\mathfrak{r}_{j-m_{2}}(k) \right) \! \right), \quad \mathbb{N} \! \ni \! j \! \geqslant \! 2,
\end{gathered}
\end{equation}
where the expansion coefficients $\mathfrak{u}_{m}(k)$ and $\mathfrak{w}_{m}(k)$ (resp., $\mathfrak{r}_{m}
(k))$, $m \! \in \! \mathbb{Z}_{+}$, $k \! \in \! \lbrace \pm 1 \rbrace$, are given in equations 
\eqref{thmk2}--\eqref{thmk10} (resp., \eqref{thmk15} and \eqref{thmk16}$)$.\footnote{Note: 
$\sum_{j=0}^{-1} \mathlarger{\ast} \! := \! 0$.}
\end{ddddd}

\emph{Proof}. The proof is presented for the case $\tau \! \to \! +\infty$ with $\varepsilon b \! > \! 0$, that 
is, $(\varepsilon_{1},\varepsilon_{2},m(\varepsilon_{2}) \vert \ell) \! = \! (0,0,0 \vert 0)$ (cf. Appendix 
\ref{sectonsymm}). Recall {}from Proposition \ref{prop1.2} that, given any solution $u(\tau)$ of the {\rm DP3E} 
\eqref{eq1.1}, the function $\hat{\varphi}(\tau)$ is defined as the general solution of the {\rm ODE} 
$\hat{\varphi}^{\prime}(\tau) \! = \! 2a \tau^{-1} \! + \! b(u(\tau))^{-1}$. {}From Propositions 1.2 and~4.1.1 of 
\cite{a1} (see, also, Section 1 of \cite{avkavint}), it can be shown that, for $\varepsilon \! \in \! \lbrace \pm 1 \rbrace$,
\begin{equation} \label{totsa7}
\hat{\varphi}(\tau) \! = \! -\mi \ln \! \left(\frac{\varepsilon \tau^{\mi a}u(\tau)}{\tau^{1/3}b(\tau)} \right):
\end{equation}
the trans-series asymptotics (as $\tau \! \to \! +\infty$ with $\varepsilon b \! > \! 0$) for $u(\tau)$ is given in Theorem 
\ref{theor2.1}, whilst only the leading-order asymptotics for the function $b(\tau)$ is derived in Lemma \ref{ginversion} 
(cf. equations \eqref{geek3}--\eqref{eequeb}); therefore, in order to proceed with the proof, trans-series asymptotics 
for $b(\tau)$ must be derived.

Commencing with the asymptotics \eqref{geek1} and \eqref{geek2}, and repeating, \emph{verbatim}, the 
asymptotic analysis of Section \ref{finalsec}, one shows that the asymptotic representation (as $\tau \! \to \! +\infty$ 
with $\varepsilon b \! > \! 0$) for the function $b(\tau)$ reads
\begin{equation} \label{totsa8} 
b(\tau) \! \underset{\tau \to +\infty}{=} \! \mathfrak{b}_{0}^{\ast}(k) \exp (-2 \mathscr{B}_{k}(\tau)), \quad 
k \! \in \! \lbrace \pm 1 \rbrace,
\end{equation}
where
\begin{equation} \label{totsa9} 
\mathfrak{b}_{0}^{\ast}(k) \! := \! (\mathfrak{b}(k))^{2}(\varepsilon b)^{1/2} \exp \! \left(2 \mi (a \! - \! \mi/2) \ln 
\big((\varepsilon b)^{1/6} \me^{\mi \pi k/3}/2 \big) \right),
\end{equation}
with $\mathfrak{b}(k)$ given in equation \eqref{eequeb}, and
\begin{align} \label{totsa10} 
\mathscr{B}_{k}(\tau) \! :=& \, \frac{\mi a}{6} \ln \tau \! - \! \frac{3k}{4}(\sqrt{3} \! + \! \mi k)(\varepsilon b)^{1/3} 
\tau^{2/3} \! + \! \sum_{m=1}^{\infty} \widetilde{\nu}_{m}(k) \big(\tau^{-1/3} \big)^{m} \nonumber \\
+& \, \left(\sum_{m=0}^{\infty} \frac{\upsilon_{m}(k)}{(\tau^{1/3})^{m}} \! + \! \mathcal{O} \! \left(
\me^{-\frac{3 \sqrt{3}}{2}(\sqrt{3}+\mi k)(\varepsilon b)^{1/3} \tau^{2/3}} \right) \! \right) 
\me^{-\frac{3 \sqrt{3}}{2}(\sqrt{3}+\mi k)(\varepsilon b)^{1/3} \tau^{2/3}};
\end{align}
it remains to determine the expansion coefficients $\lbrace \widetilde{\nu}_{m}(k) \rbrace_{m=1}^{\infty}$ and the 
first non-zero coefficient $\upsilon_{m}(k)$. Via the definitions \eqref{newlax2}, the isomonodromy deformations 
\eqref{newlax8}, the definitions \eqref{equu18}, \eqref{equu19}, and \eqref{eq3.2}, and equation \eqref{totsa8}, one 
shows that the function $\mathscr{B}_{k}(\tau)$ solves the following inhomogeneous second-order non-linear {\rm ODE}:
\begin{equation} \label{totsa11} 
\mathscr{B}_{k}^{\prime \prime}(\tau) \! - \! 2(\mathscr{B}_{k}^{\prime}(\tau))^{2} \! - \! \left(\frac{\md}{\md \tau} 
\ln \! \left(\frac{u(\tau)}{\tau^{2/3}} \right) \! \right) \! \mathscr{B}_{k}^{\prime}(\tau) \! = \! \frac{1}{2 \tau} \! \left(
\frac{2}{3} \frac{\md}{\md \tau} \ln \! \left(\frac{u(\tau)}{\tau^{1/3}} \right) \! + \! \mi a \frac{\md}{\md \tau} \ln \! 
\left(\frac{u(\tau)}{\tau^{1+ \mi a}} \right) \! + \! 8 \varepsilon u(\tau) \right),
\end{equation}
where (cf. equation \eqref{iden217}) $u(\tau) \! = \! \tfrac{1}{2} \varepsilon (\varepsilon b)^{2/3} \me^{-\mi 2 \pi k/3} \big(
\tau^{1/3} \! + \! v_{0,k}(\tau) \big)$, with $v_{0,k}(\tau)$ given in the asymptotics \eqref{geek1}. {}From the expression 
for $u^{\prime}(\tau)$ given in the proof of Proposition 5.7 in \cite{a1} and the definitions \eqref{newlax2} and 
\eqref{eq3.2}, it follows that
\begin{equation} \label{totsa12} 
\frac{\md}{\md \tau} \ln (u(\tau)) \! = \! \frac{u^{\prime}(\tau)}{u(\tau)} \! = \! \frac{1}{\tau} \! + \! 2 \varepsilon 
\! \left(\frac{a(\tau)d(\tau) \! - \! b(\tau)c(\tau)}{u(\tau)} \right);
\end{equation}
via equation \eqref{iden217}, the asymptotics \eqref{geek1}, \eqref{geek2}, \eqref{isomk2}, and \eqref{isomk3}, 
and equation \eqref{totsa12}, one shows that, for $k \! \in \! \lbrace \pm 1 \rbrace$,
\begin{equation} \label{totsa13} 
\frac{\md}{\md \tau} \ln (u(\tau)) \underset{\tau \to +\infty}{=} \frac{1}{3 \tau} \! \left(1 \! + \! \sum_{m=0}^{\infty} 
\frac{\mu_{m}^{\ast}(k)}{(\tau^{1/3})^{m+2}} \right) \! - \! \mathbb{V}_{0}(k) \tau^{-2/3} \me^{-\frac{3 \sqrt{3}}{2}
(\sqrt{3}+\mi k)(\varepsilon b)^{1/3} \tau^{2/3}} \! \left(1 \! + \! \mathcal{O} \big(\tau^{-1/3} \big) \right),
\end{equation}
where the expansion coefficients $\lbrace \mu_{m}^{\ast}(k) \rbrace_{m=0}^{\infty}$ are given in equations 
\eqref{totsa5} and \eqref{totsa6}, and
\begin{equation} \label{totsa14} 
\mathbb{V}_{0}(k) \! := \! \frac{k2^{1/2}3^{1/4} \me^{\mi \pi k/3} \me^{\mi \pi k/4}(\varepsilon b)^{1/6}
(s_{0}^{0} \! - \! \mi \me^{-\pi a})}{\sqrt{\pi}(2 \! + \! \sqrt{3})^{-\mi ka}}.
\end{equation}
Substituting the asymptotic expansions \eqref{thmk11}, \eqref{totsa10}, and \eqref{totsa13} into the second-order 
non-linear {\rm ODE} \eqref{totsa11}, and equating coefficients of terms of orders $\mathcal{O} \big((\tau^{-1/3})^{m_{1}} 
\exp (-\tfrac{3 \sqrt{3}}{2}(\sqrt{3} \! + \! \mi k)(\varepsilon b)^{1/3} \tau^{2/3}) \big)$, $m_{1} \! = \! 2,3$, and $\mathcal{O} 
\big((\tau^{-1/3})^{m_{2}} \big)$, $\mathbb{N} \! \ni \! m_{2} \! \geqslant \! 2$, one arrives at, after simplification, for $k \! 
\in \! \lbrace \pm 1 \rbrace$, in the indicated order: (i) $\mathcal{O} \big(\tau^{-2/3} \exp (-\tfrac{3 \sqrt{3}}{2}(\sqrtsign{3} 
\! + \! \mi k)(\varepsilon b)^{1/3} \tau^{2/3}) \big)$ $\Rightarrow$
\begin{equation} \label{totsa15} 
\sqrt{3}(\sqrt{3} \! + \! \mi k)^{2}(\sqrt{3} \! - \! 2k)(\varepsilon b)^{2/3} \upsilon_{0}(k) \! = \! 0;
\end{equation}
(ii) $\mathcal{O} \big(\tau^{-1} \exp (-\tfrac{3 \sqrt{3}}{2}(\sqrt{3} \! + \! \mi k)(\varepsilon b)^{1/3} \tau^{2/3}) \big)$ $\Rightarrow$
\begin{equation} \label{totsa16} 
\sqrt{3}(\sqrt{3} \! + \! \mi k)^{2}(\sqrt{3} \! - \! 2k)(\varepsilon b)^{2/3} \upsilon_{1}(k) \! = \! \frac{\left(-\mi 2 \! + \! 
\sqrt{3}(\sqrt{3} \! + \! \mi k) \me^{\mi \pi k/3} \right) \! \me^{\mi \pi k/4}(\varepsilon b)^{1/2}(s_{0}^{0} \! - \! \mi 
\me^{-\pi a})}{\sqrt{2 \pi} \, 3^{1/4}(2 \! + \! \sqrt{3})^{-\mi ka}};
\end{equation}
(iii) $\mathcal{O} \big(\tau^{-2/3} \big)$ $\Rightarrow$
\begin{equation} \label{totsa17} 
-4 \me^{-\mi 2 \pi k/3} \! = \! (k \sqrt{3} \! + \! \mi)^{2};
\end{equation}
(iv) $\mathcal{O} \big(\tau^{-4/3} \big)$ $\Rightarrow$
\begin{equation} \label{totsa18} 
2 \mi \me^{-\mi \pi k/3} \! = \! k \sqrt{3} \! + \! \mi;
\end{equation}
(v) $\mathcal{O} \big(\tau^{-5/3} \big)$ $\Rightarrow$
\begin{equation} \label{totsa19} 
\widetilde{\nu}_{1}(k) \! = \! 0;
\end{equation}
(vi) $\mathcal{O} \big(\tau^{-2} \big)$ $\Rightarrow$
\begin{equation} \label{totsa20} 
4 \widetilde{\nu}_{2}(k) \! - \! \frac{a \me^{\mi \pi k/3}}{3(\varepsilon b)^{1/3}} \! = \! 
\frac{2 \mi a(a \! - \! \mi/2) \me^{\mi \pi k/3}}{3(\varepsilon b)^{1/3}};
\end{equation}
(vii) $\mathcal{O} \big(\tau^{-7/3} \big)$ $\Rightarrow$
\begin{equation} \label{totsa21} 
\widetilde{\nu}_{3}(k) \! = \! 0;
\end{equation}
(viii) $\mathcal{O} \big(\tau^{-8/3} \big)$ $\Rightarrow$
\begin{equation} \label{totsa22} 
4 \mi \me^{-\mi \pi k/3}(\varepsilon b)^{1/3} \widetilde{\nu}_{4}(k) \! = \! \frac{a \me^{\mi \pi k/3}}{9(\varepsilon b)^{1/3}} 
\! \left(\frac{1 \! - \! 2a^{2}}{3} \! + \! \mi a \right);
\end{equation}
and (ix) $\mathcal{O} \big(\tau^{-(m+9)/3} \big)$, $m \! \in \! \mathbb{Z}_{+}$, $\Rightarrow$
\begin{align} \label{totsa23} 
4 \mi \me^{-\mi \pi k/3}(\varepsilon b)^{1/3}(m \! + \! 5) \widetilde{\nu}_{m+5}(k) 
\! =& \, -6 \me^{-\mi 2 \pi k/3}(\varepsilon b)^{2/3} \mathfrak{u}_{m+5}(k) \! - \! 
\frac{(1 \! + \! \mi 2a)}{3} \mu_{m+1}^{\ast}(k) \nonumber \\
+& \, \mi \me^{-\mi \pi k/3}(\varepsilon b)^{1/3} \mu_{m+3}^{\ast}(k) \! + \! 
\frac{1}{3} \! \left(\vphantom{M^{M^{M^{M}}}}(m \! + \! 3)(m \! + \! 5 \! + \! \mi 2a) 
\widetilde{\nu}_{m+3}(k) \right. \nonumber \\
+&\left. \, \sum_{j=0}^{m-1}(j \! + \! 1) \widetilde{\nu}_{j+1}(k)(\mu_{m-j}^{\ast}(k) 
\! - \! 2(m \! + \! 2 \! - \! j) \widetilde{\nu}_{m+2-j}(k)) \right. \nonumber \\
-&\left. \frac{2 \mi a^{2} \me^{\mi \pi k/3}}{3(\varepsilon b)^{1/3}}(m \! + \! 1) 
\widetilde{\nu}_{m+1}(k) \right),
\end{align}
with the convention $\sum_{j=0}^{-1} \mathlarger{\ast} \! := \! 0$. Solving equations \eqref{totsa15} and \eqref{totsa16} 
for $\upsilon_{0}(k)$ and $\upsilon_{1}(k)$, $k \! \in \! \lbrace \pm 1 \rbrace$, respectively, one shows that
\begin{equation}
\upsilon_{0}(k) \! = \! 0 \quad \qquad \text{and} \qquad \quad \upsilon_{1}(k) \! = \! -\frac{\mi \me^{-\mi \pi k/3} 
\me^{\mi \pi k/4}(2 \! + \! \sqrt{3})^{\mi ka}(s_{0}^{0} \! - \! \mi \me^{-\pi a})}{\sqrt{2 \pi} \, 3^{3/4}(\sqrt{3} \! - \! k)
(\varepsilon b)^{1/6}}.
\end{equation}
Equations \eqref{totsa17} and \eqref{totsa18} are identities. Solving equations \eqref{totsa19}--\eqref{totsa23} for the 
coefficients $\widetilde{\nu}_{1}(k)$, $\widetilde{\nu}_{2}(k)$, $\widetilde{\nu}_{3}(k)$, $\widetilde{\nu}_{4}(k)$, and 
$\widetilde{\nu}_{m+5}(k)$, $k \! \in \! \lbrace \pm 1 \rbrace$, $m \! \in \! \mathbb{Z}_{+}$, respectively, one arrives 
at equations \eqref{totsa3}--\eqref{totsa6}; therefore, the trans-series asymptotics for the function $b(\tau)$ is now 
established via equations \eqref{totsa8}--\eqref{totsa10}; in particular, for $k \! \in \! \lbrace \pm 1 \rbrace$,
\begin{align} \label{totsa25} 
\mathscr{B}_{k}(\tau) \! \underset{\tau \to +\infty}{=}& \, \frac{\mi a}{6} \ln \tau \! - \! \frac{3k}{4}(\sqrt{3} \! + \! \mi k)
(\varepsilon b)^{1/3} \tau^{2/3} \! + \! \sum_{m=1}^{\infty} \widetilde{\nu}_{m}(k) \big(\tau^{-1/3} \big)^{m} \nonumber \\
-& \, \frac{\mi \me^{-\mi \pi k/3} \me^{\mi \pi k/4}(2 \! + \! \sqrt{3})^{\mi ka}(s_{0}^{0} \! - \! \mi \me^{-\pi a})}{\sqrt{2 \pi} 
\, 3^{3/4}(\sqrt{3} \! - \! k)(\varepsilon b)^{1/6} \tau^{1/3}} \me^{-\frac{3 \sqrt{3}}{2}(\sqrt{3}+\mi k)
(\varepsilon b)^{1/3} \tau^{2/3}} \! \left(1 \! + \! \mathcal{O} \big(\tau^{-1/3} \big) \right).
\end{align}

Via equation \eqref{iden217}, the asymptotics \eqref{geek1} and \eqref{geek2}, equation \eqref{totsa7}, the definition 
\eqref{totsa9} (cf. equation \eqref{eequeb}), the asymptotics \eqref{totsa25}, and the expansion
\begin{equation} \label{totsa26} 
\ln \! \left(1 \! + \! \sum_{m=0}^{\infty} \frac{\mathfrak{u}_{m}(k)}{(\tau^{1/3})^{m+2}} \right) 
\underset{\tau \to +\infty}{=} \mathlarger{\sum_{m=2}^{\infty}} \sum_{\substack{n,\mathfrak{l} 
\in \mathbb{N}\\\mathfrak{l} \geqslant n\\n+ \mathfrak{l}=m}} \sum_{\substack{\mathfrak{i}_{1}+2 
\mathfrak{i}_{2} + \dotsb + \mathfrak{l} \mathfrak{i}_{\mathfrak{l}}=\mathfrak{l}\\
\mathfrak{i}_{1} + \mathfrak{i}_{2} + \dotsb + \mathfrak{i}_{\mathfrak{l}}=n}} \, 
\frac{\mathrm{S}^{(1)}_{n}(\mathfrak{u}_{0}(k))^{\mathfrak{i}_{1}}(\mathfrak{u}_{1}
(k))^{\mathfrak{i}_{2}} \dotsb (\mathfrak{u}_{\mathfrak{l}-1}(k))^{\mathfrak{i}_{
\mathfrak{l}}}}{\mathfrak{i}_{1}! \mathfrak{i}_{2}! \dotsb \mathfrak{i}_{\mathfrak{l}}!}
(\tau^{-1/3})^{m},
\end{equation}
where $\mathrm{S}^{(1)}_{n} \! = \! (-1)^{n-1}(n \! - \! 1)!$ is a special value of the Stirling number of the first 
kind \cite{a24}, one arrives at, for $k \! \in \! \lbrace \pm 1 \rbrace$, the $(\varepsilon_{1},\varepsilon_{2},
m(\varepsilon_{2}) \vert \ell) \! = \! (0,0,0 \vert 0)$ trans-series asymptotics (as $\tau \! \to \! +\infty$ with 
$\varepsilon b \! > \! 0$) for the function $\hat{\varphi}(\tau)$:
\begin{align} \label{totsa27} 
\hat{\varphi}(\tau) \underset{\tau \to +\infty}{=}& \, \mi \mathsterling_{k}(0,0,0 \vert 0) \! - \! k \pi \! + \! 
\frac{3k \mi}{2}(\sqrt{3} \! + \! \mi k)(\varepsilon b)^{1/3} \tau^{2/3} \! + \! 2a \ln \! \left(\frac{2 \me^{-\mi \pi k/3} 
\tau^{2/3}}{(\varepsilon b)^{1/6}} \right) \nonumber \\
-& \, \mi \, \mathlarger{\sum_{m=2}^{\infty}} \! \left(2 \widetilde{\nu}_{m}(k) \! + 
\! \sum_{\substack{n,\mathfrak{l} \in \mathbb{N}\\\mathfrak{l} \geqslant n\\n+ 
\mathfrak{l}=m}} \sum_{\substack{\mathfrak{i}_{1}+2 \mathfrak{i}_{2} + \dotsb + 
\mathfrak{l} \mathfrak{i}_{\mathfrak{l}}=\mathfrak{l}\\\mathfrak{i}_{1} + 
\mathfrak{i}_{2} + \dotsb + \mathfrak{i}_{\mathfrak{l}}=n}} \, 
(-1)^{n-1}(n \! - \! 1)! \prod_{\mathfrak{j}=1}^{\mathfrak{l}} \frac{(\mathfrak{u}_{
\mathfrak{j}-1}(k))^{\mathfrak{i}_{\mathfrak{j}}}}{\mathfrak{i}_{\mathfrak{j}}!} \right) \! \big(\tau^{-1/3} \big)^{m} 
\nonumber \\
-& \, \frac{k \me^{-\mi \pi k/3}\me^{\mi \pi k/4}(2 \! + \! \sqrt{3})^{\mi ka}(s_{0}^{0} \! - \! \mi \me^{-\pi a})}{\sqrt{2 \pi} 
\, 3^{3/4}(\varepsilon b)^{1/6} \tau^{1/3}} \me^{-\frac{3 \sqrt{3}}{2}(\sqrt{3}+\mi k)(\varepsilon b)^{1/3} \tau^{2/3}} \! 
\left(1 \! + \! \mathcal{O} \big(\tau^{-1/3} \big) \right),
\end{align}
where
\begin{equation} \label{totsa28} 
\mathsterling_{k}(0,0,0 \vert 0) \! = \! 
\begin{cases}
\ln \! \left(g_{11} \me^{\pi a} \right)^{2}, &\text{$k \! = \! +1$,} \\
\ln \! \left(g_{22} \me^{\pi a} \right)^{-2}, &\text{$k \! = \! -1$.}
\end{cases}
\end{equation}

Finally, applying the (map) isomorphism (cf. Appendix \ref{sectonsymm}) 
$\mathscr{F}^{\scriptscriptstyle \lbrace \ell \rbrace}_{\scriptscriptstyle \varepsilon_{1},
\varepsilon_{2},m(\varepsilon_{2})}$, whose action on $\mathscr{M}$ is given by equations 
\eqref{laxhat76}--\eqref{laxhat90} and \eqref{laxhat99}--\eqref{laxhat113}, to the corresponding 
$(\varepsilon_{1},\varepsilon_{2},m(\varepsilon_{2}) \vert \ell) \! = \! (0,0,0 \vert 0)$ asymptotics 
\eqref{totsa27} for $\hat{\varphi}(\tau)$, one arrives at the trans-series asymptotics \eqref{totsa1} 
(and equations \eqref{totsa2}--\eqref{totsa6}) stated in the theorem. \hfill $\qed$
\begin{eeeee} \label{isoasybeeh} 
Via equation \eqref{totsa8}, the definition \eqref{totsa9} (cf. equation \eqref{eequeb}), and the asymptotics 
\eqref{totsa25}, one arrives at, {}from the asymptotics \eqref{isomk2}, \eqref{isomk3}, and \eqref{isomk4}, 
respectively, the trans-series asymptotics (as $\tau \! \to \! +\infty$ with $\varepsilon b \! > \! 0$) for the 
functions $c(\tau)$, $d(\tau)$, and $a(\tau)$. \hfill $\blacksquare$
\end{eeeee} 
\begin{eeeee} \label{mer2mus} 
It is instructive to illustrate the first few contributions of the multi-indexed double summation in equation \eqref{totsa26} 
to the asymptotics of $\hat{\varphi}(\tau)$ for various values of the index $m$: (i) for $m \! = \! 2$ (that is, $\mathcal{O}
(\tau^{-2/3})$), $(n,\mathfrak{l}) \! = \! (1,1)$ $\Rightarrow$ $\mathfrak{i}_{1} \! = \! 1$, thus, for $k \! \in \! \lbrace \pm 1 
\rbrace$,\footnote{Recall that the expansion coefficients $\lbrace \mathfrak{u}_{j}(k) \rbrace_{j=0}^{\infty}$, $k \! \in \! 
\lbrace \pm 1 \rbrace$, are given in equations \eqref{thmk2}--\eqref{thmk10}.}
\begin{equation*} 
\sum_{\substack{n,\mathfrak{l} \in \mathbb{N}\\\mathfrak{l} \geqslant n\\n+ \mathfrak{l}=2}} 
\sum_{\substack{\mathfrak{i}_{1}+2 \mathfrak{i}_{2} + \dotsb + \mathfrak{l} \mathfrak{i}_{
\mathfrak{l}}=\mathfrak{l}\\\mathfrak{i}_{1} + \mathfrak{i}_{2} + \dotsb + \mathfrak{i}_{
\mathfrak{l}}=n}} \, (-1)^{n-1}(n \! - \! 1)! \prod_{\mathfrak{j}=1}^{\mathfrak{l}} 
\frac{(\mathfrak{u}_{\mathfrak{j}-1}(k))^{\mathfrak{i}_{\mathfrak{j}}}}{\mathfrak{i}_{\mathfrak{j}}!} 
\! = \! \mathfrak{u}_{0}(k) \! = \! \frac{a \me^{-\mi 2 \pi k/3}}{3(\varepsilon b)^{1/3}};
\end{equation*}
(ii) for $m \! = \! 3$ (that is, $\mathcal{O}(\tau^{-1})$), $(n,\mathfrak{l}) \! = \! (1,2)$ 
$\Rightarrow$ $(\mathfrak{i}_{1},\mathfrak{i}_{2}) \! = \! (0,1)$, thus, for $k \! \in \! \lbrace 
\pm 1 \rbrace$,
\begin{equation*} 
\sum_{\substack{n,\mathfrak{l} \in \mathbb{N}\\\mathfrak{l} \geqslant n\\n+ \mathfrak{l}=3}} 
\sum_{\substack{\mathfrak{i}_{1}+2 \mathfrak{i}_{2} + \dotsb + \mathfrak{l} \mathfrak{i}_{
\mathfrak{l}}=\mathfrak{l}\\\mathfrak{i}_{1} + \mathfrak{i}_{2} + \dotsb + \mathfrak{i}_{
\mathfrak{l}}=n}} \, (-1)^{n-1}(n \! - \! 1)! \prod_{\mathfrak{j}=1}^{\mathfrak{l}} 
\frac{(\mathfrak{u}_{\mathfrak{j}-1}(k))^{\mathfrak{i}_{\mathfrak{j}}}}{\mathfrak{i}_{\mathfrak{j}}!} 
\! = \! \mathfrak{u}_{1}(k) \! = \! 0;
\end{equation*}
(iii) for $m \! = \! 4$ (that is, $\mathcal{O}(\tau^{-4/3})$), $(n,\mathfrak{l}) \! = \! (2,2)$ 
$\Rightarrow$ $(\mathfrak{i}_{1},\mathfrak{i}_{2}) \! = \! (2,0)$, and $(n,\mathfrak{l}) \! = \! 
(1,3)$ $\Rightarrow$ $(\mathfrak{i}_{1},\mathfrak{i}_{2},\mathfrak{i}_{3}) \! = \! (0,0,1)$, thus, 
for $k \! \in \! \lbrace \pm 1 \rbrace$,
\begin{equation*} 
\sum_{\substack{n,\mathfrak{l} \in \mathbb{N}\\\mathfrak{l} \geqslant n\\n+ \mathfrak{l}=4}} 
\sum_{\substack{\mathfrak{i}_{1}+2 \mathfrak{i}_{2} + \dotsb + \mathfrak{l} \mathfrak{i}_{
\mathfrak{l}}=\mathfrak{l}\\\mathfrak{i}_{1} + \mathfrak{i}_{2} + \dotsb + \mathfrak{i}_{
\mathfrak{l}}=n}} \, (-1)^{n-1}(n \! - \! 1)! \prod_{\mathfrak{j}=1}^{\mathfrak{l}} 
\frac{(\mathfrak{u}_{\mathfrak{j}-1}(k))^{\mathfrak{i}_{\mathfrak{j}}}}{\mathfrak{i}_{
\mathfrak{j}}!} \! = \! \mathfrak{u}_{2}(k) \! - \! \frac{(\mathfrak{u}_{0}(k))^{2}}{2} 
\! = \! \frac{a^{2} \me^{-\mi \pi k/3}}{18(\varepsilon b)^{2/3}};
\end{equation*}
(iv) for $m \! = \! 5$ (that is, $\mathcal{O}(\tau^{-5/3})$), $(n,\mathfrak{l}) \! = \! (2,3)$ 
$\Rightarrow$ $(\mathfrak{i}_{1},\mathfrak{i}_{2},\mathfrak{i}_{3}) \! = \! (1,1,0)$, and 
$(n,\mathfrak{l}) \! = \! (1,4)$ $\Rightarrow$ $(\mathfrak{i}_{1},\mathfrak{i}_{2},\mathfrak{i}_{3},
\mathfrak{i}_{4}) \! = \! (0,0,0,1)$, thus, for $k \! \in \! \lbrace \pm 1 \rbrace$,
\begin{equation*} 
\sum_{\substack{n,\mathfrak{l} \in \mathbb{N}\\\mathfrak{l} \geqslant n\\n+ \mathfrak{l}=5}} 
\sum_{\substack{\mathfrak{i}_{1}+2 \mathfrak{i}_{2} + \dotsb + \mathfrak{l} \mathfrak{i}_{
\mathfrak{l}}=\mathfrak{l}\\\mathfrak{i}_{1} + \mathfrak{i}_{2} + \dotsb + \mathfrak{i}_{
\mathfrak{l}}=n}} \, (-1)^{n-1}(n \! - \! 1)! \prod_{\mathfrak{j}=1}^{\mathfrak{l}} 
\frac{(\mathfrak{u}_{\mathfrak{j}-1}(k))^{\mathfrak{i}_{\mathfrak{j}}}}{\mathfrak{i}_{\mathfrak{j}}!} 
\! = \! \mathfrak{u}_{3}(k) \! - \! \mathfrak{u}_{0}(k) \mathfrak{u}_{1}(k) \! = \! 0;
\end{equation*}
and (v) for $m \! = \! 6$ (that is, $\mathcal{O}(\tau^{-2})$), $(n,\mathfrak{l}) \! = \! (3,3)$ 
$\Rightarrow$ $(\mathfrak{i}_{1},\mathfrak{i}_{2},\mathfrak{i}_{3}) \! = \! (3,0,0)$, 
$(n,\mathfrak{l}) \! = \! (2,4)$ $\Rightarrow$ $(\mathfrak{i}_{1},\mathfrak{i}_{2},\mathfrak{i}_{3},
\mathfrak{i}_{4}) \! \in \! \lbrace (1,0,1,0),(0,2,0,0) \rbrace$, and $(n,\mathfrak{l}) \! = \! (1,5)$ 
$\Rightarrow$ $(\mathfrak{i}_{1},\mathfrak{i}_{2},\mathfrak{i}_{3},\mathfrak{i}_{4},\mathfrak{i}_{5}) 
\! = \! (0,0,0,0,1)$, thus, for $k \! \in \! \lbrace \pm 1 \rbrace$,
\begin{align*} 
\sum_{\substack{n,\mathfrak{l} \in \mathbb{N}\\\mathfrak{l} \geqslant n\\n+ \mathfrak{l}=6}} 
\sum_{\substack{\mathfrak{i}_{1}+2 \mathfrak{i}_{2} + \dotsb + \mathfrak{l} \mathfrak{i}_{
\mathfrak{l}}=\mathfrak{l}\\\mathfrak{i}_{1} + \mathfrak{i}_{2} + \dotsb + \mathfrak{i}_{
\mathfrak{l}}=n}} \, (-1)^{n-1}(n \! - \! 1)! \prod_{\mathfrak{j}=1}^{\mathfrak{l}} 
\frac{(\mathfrak{u}_{\mathfrak{j}-1}(k))^{\mathfrak{i}_{\mathfrak{j}}}}{\mathfrak{i}_{
\mathfrak{j}}!} =& \, \mathfrak{u}_{4}(k) \! - \! \mathfrak{u}_{0}(k) \mathfrak{u}_{2}(k) \! + \! 
\frac{(\mathfrak{u}_{0}(k))^{3}}{3} \! - \! \frac{(\mathfrak{u}_{1}(k))^{2}}{2} \nonumber \\
=& -\frac{a}{3^{4}(\varepsilon b)}. \tag*{$\blacksquare$}
\end{align*}
\end{eeeee} 
\begin{ddddd} \label{pfeetotsb} 
Let $u(\tau)$ be a solution of the {\rm DP3E} \eqref{eq1.1} and $\hat{\varphi}(\tau)$ be the general solution of 
the {\rm ODE} $\hat{\varphi}^{\prime}(\tau) \! = \! 2a \tau^{-1} \! + \! b(u(\tau))^{-1}$ for $\varepsilon b \! > \! 0$ 
corresponding to the monodromy data $(a,s^{0}_{0},s^{\infty}_{0},s^{\infty}_{1},g_{11},g_{12},g_{21},g_{22})$. 
Let $\hat{\varepsilon}_{1} \! \in \! \lbrace \pm 1 \rbrace$, $\hat{\varepsilon}_{2} \! \in  \! \lbrace 0,\pm 1 \rbrace$, 
$\hat{m}(\hat{\varepsilon}_{2}) \! = \! 
\left\{
\begin{smallmatrix}
0, \, \, \hat{\varepsilon}_{2} \in \lbrace \pm 1 \rbrace, \\
\pm \hat{\varepsilon}_{1}, \, \, \hat{\varepsilon}_{2}=0,
\end{smallmatrix}
\right.$ $\hat{\ell} \! \in \! \lbrace 0,1 \rbrace$, and $\varepsilon b \! = \! \vert 
\varepsilon b \vert \me^{\mi \pi \hat{\varepsilon}_{2}}$. For $k \! = \! +1$, let
\begin{equation*}
\hat{g}_{11}(\hat{\varepsilon}_{1},\hat{\varepsilon}_{2},\hat{m}(\hat{\varepsilon}_{2}) 
\vert \hat{\ell}) \hat{g}_{12}(\hat{\varepsilon}_{1},\hat{\varepsilon}_{2},\hat{m}
(\hat{\varepsilon}_{2}) \vert \hat{\ell}) \hat{g}_{21}(\hat{\varepsilon}_{1},
\hat{\varepsilon}_{2},\hat{m}(\hat{\varepsilon}_{2}) \vert \hat{\ell}) \! \neq \! 0 \quad 
\text{and} \quad \hat{g}_{22}(\hat{\varepsilon}_{1},\hat{\varepsilon}_{2},\hat{m}
(\hat{\varepsilon}_{2}) \vert \hat{\ell}) \! = \! 0,
\end{equation*}
and, for $k \! = \! -1$, let
\begin{equation*}
\hat{g}_{11}(\hat{\varepsilon}_{1},\hat{\varepsilon}_{2},\hat{m}(\hat{\varepsilon}_{2}) 
\vert \hat{\ell}) \! = \! 0 \quad \text{and} \quad \hat{g}_{12}(\hat{\varepsilon}_{1},
\hat{\varepsilon}_{2},\hat{m}(\hat{\varepsilon}_{2}) \vert \hat{\ell}) \hat{g}_{21}
(\hat{\varepsilon}_{1},\hat{\varepsilon}_{2},\hat{m}(\hat{\varepsilon}_{2}) \vert 
\hat{\ell}) \hat{g}_{22}(\hat{\varepsilon}_{1},\hat{\varepsilon}_{2},\hat{m}
(\hat{\varepsilon}_{2}) \vert \hat{\ell}) \! \neq \! 0,
\end{equation*}
where explicit expressions for $\hat{g}_{ij}(\hat{\varepsilon}_{1},\hat{\varepsilon}_{2},
\hat{m}(\hat{\varepsilon}_{2}) \vert \hat{\ell})$, $i,j \! \in \! \lbrace 1,2 \rbrace$, are 
given in Appendix \ref{sectonsymm}, equations \eqref{laxhat91}--\eqref{laxhat98} and 
\eqref{laxhat114}--\eqref{laxhat121}. Then, for $\hat{s}_{0}^{0}(\hat{\varepsilon}_{1},
\hat{\varepsilon}_{2},\hat{m}(\hat{\varepsilon}_{2}) \vert \hat{\ell}) \! \neq \! \mi 
\me^{(-1)^{\hat{\varepsilon}_{2}} \pi a}$,\footnote{Recall that (cf. Remark \ref{stokmxx}) 
$\hat{s}_{0}^{0}(\hat{\varepsilon}_{1},\hat{\varepsilon}_{2},\hat{m}(\hat{\varepsilon}_{2}) \vert 
\hat{\ell}) \! = \! s_{0}^{0}$. For $\hat{s}_{0}^{0}(\hat{\varepsilon}_{1},\hat{\varepsilon}_{2},
\hat{m}(\hat{\varepsilon}_{2}) \vert \hat{\ell}) \! = \! \mi \me^{(-1)^{\hat{\varepsilon}_{2}} 
\pi a}$, the exponentially small correction term in the asymptotics \eqref{totsb1} is absent.}
\begin{align} \label{totsb1} 
(-1)^{1+\hat{\varepsilon}_{2}} \hat{\varphi}(\tau) \underset{\tau \to +\infty \me^{\mi \pi 
\hat{\varepsilon}_{1}/2}}{=}& \, 3 \me^{\frac{\mi 2 \pi k}{3}}(-1)^{\hat{\varepsilon}_{2}}(\varepsilon b)^{1/3} 
\tau_{\ast}^{2/3} \! + \! 2(-1)^{1+\hat{\varepsilon}_{2}}a \ln \! \left(\frac{2 \me^{-\frac{\mi \pi k}{3}} 
\tau_{\ast}^{2/3}}{(\varepsilon b \me^{-\mi \pi \hat{\varepsilon}_{2}})^{1/6}} \right) \! + \! \mi 
\hat{\mathsterling}_{k}(\hat{\varepsilon}_{1},\hat{\varepsilon}_{2},\hat{m}(\hat{\varepsilon}_{2}) 
\vert \hat{\ell}) \nonumber \\
-& \, k \pi  \! - \! \mi \, \mathlarger{\sum_{m=2}^{\infty}} \! \left(2 \widehat{\nu}_{m}(k) \! + \! 
\sum_{\substack{n,\mathfrak{l} \in \mathbb{N}\\\mathfrak{l} \geqslant n\\n+ \mathfrak{l}=m}} 
\, \sum_{\substack{\mathfrak{i}_{1}+2 \mathfrak{i}_{2} + \dotsb + \mathfrak{l} 
\mathfrak{i}_{\mathfrak{l}}=\mathfrak{l}\\\mathfrak{i}_{1} + \mathfrak{i}_{2} + \dotsb + 
\mathfrak{i}_{\mathfrak{l}}=n}} \, (-1)^{n-1}(n \! - \! 1)! \prod_{\mathfrak{j}=1}^{\mathfrak{l}} 
\frac{(\hat{\mathfrak{u}}_{\mathfrak{j}-1}(k))^{\mathfrak{i}_{\mathfrak{j}}}}{\mathfrak{i}_{\mathfrak{j}}!} 
\right) \nonumber \\
\times& \, \big(\tau_{\ast}^{-1/3} \big)^{m} - \! \frac{k \me^{-\frac{\mi \pi k}{3}} \me^{\frac{\mi \pi k}{4}}
(2 \! + \! \sqrt{3})^{\mi k(-1)^{1+\hat{\varepsilon}_{2}}a}}{\sqrt{2 \pi} \, 3^{3/4}(\varepsilon b \me^{-\mi \pi 
\hat{\varepsilon}_{2}})^{1/6} \tau_{\ast}^{1/3}} \! \left(\hat{s}_{0}^{0}(\hat{\varepsilon}_{1},\hat{\varepsilon}_{2},
\hat{m}(\hat{\varepsilon}_{2}) \vert \hat{\ell}) \! - \! \mi \me^{(-1)^{\hat{\varepsilon}_{2}} \pi a} \right) 
\nonumber \\
\times& \, \me^{-\frac{3 \sqrt{3}}{2}(\sqrt{3}+\mi k)(-1)^{\hat{\varepsilon}_{2}}(\varepsilon b)^{1/3} 
\tau_{\ast}^{2/3}} \! \left(1 \! + \! \mathcal{O} \big(\tau^{-1/3} \big) \right), \quad k \! \in \! \lbrace \pm 1 \rbrace,
\end{align}
where $\tau_{\ast}$ is defined by equation \eqref{teas},
\begin{equation} \label{totsb2} 
\hat{\mathsterling}_{k}(\hat{\varepsilon}_{1},\hat{\varepsilon}_{2},\hat{m}(\hat{\varepsilon}_{2}) 
\vert \hat{\ell}) \! = \! 
\begin{cases}
\ln \! \left(\hat{g}_{11}(\hat{\varepsilon}_{1},\hat{\varepsilon}_{2},\hat{m}(\hat{\varepsilon}_{2}) 
\vert \hat{\ell}) \me^{(-1)^{1+\hat{\varepsilon}_{2}} \pi a} \right)^{2}, &\text{$k \! = \! +1$,} \\
\ln \! \left(\hat{g}_{22}(\hat{\varepsilon}_{1},\hat{\varepsilon}_{2},\hat{m}(\hat{\varepsilon}_{2}) 
\vert \hat{\ell}) \me^{(-1)^{1+\hat{\varepsilon}_{2}} \pi a} \right)^{-2}, &\text{$k \! = \! -1$,}
\end{cases}
\end{equation}
\begin{equation} \label{totsb3}
\begin{gathered}
\widehat{\nu}_{1}(k) \! = \! 0, \quad \qquad \widehat{\nu}_{2}(k) \! = \! -\frac{a(1 \! + \! \mi 
(-1)^{1+\hat{\varepsilon}_{2}}a) \me^{\mi \pi k/3}}{6(\varepsilon b)^{1/3}}, \qquad \quad 
\widehat{\nu}_{3}(k) \! = \! 0, \\
\widehat{\nu}_{4}(k) \! = \! \frac{\mi (-1)^{\hat{\varepsilon}_{2}}a \me^{\mi 2 \pi k/3}}{36(\varepsilon b)^{2/3}} 
\! \left(\frac{1 \! - \! 2a^{2}}{3} \! + \! \mi (-1)^{1+\hat{\varepsilon}_{2}}a \right),
\end{gathered}
\end{equation}
and 
\begin{align} \label{totsb4} 
(m \! + \! 5) \widehat{\nu}_{m+5}(k) \! =& \, \frac{3 \mi}{2} \me^{-\mi \pi k/3}(-1)^{\hat{\varepsilon}_{2}}
(\varepsilon b)^{1/3} \hat{\mathfrak{u}}_{m+5}(k) \! + \! \frac{\mi (-1)^{\hat{\varepsilon}_{2}} \me^{\mi \pi k/3}
(1 \! + \! 2 \mi (-1)^{1+ \hat{\varepsilon}_{2}}a)}{12(\varepsilon b)^{1/3}} \hat{\mu}_{m+1}^{\ast}(k) \! + \! 
\frac{1}{4} \hat{\mu}_{m+3}^{\ast}(k) \nonumber \\
-& \, \frac{\mi (-1)^{\hat{\varepsilon}_{2}} \me^{\mi \pi k/3}}{12(\varepsilon b)^{1/3}} \! \left(
\vphantom{M^{M^{M^{M}}}}(m \! + \! 3)(m \! + \! 5 \! + \! 2 \mi (-1)^{1+\hat{\varepsilon}_{2}}a) 
\widehat{\nu}_{m+3}(k) \! - \! \frac{\mi (-1)^{\hat{\varepsilon}_{2}}2a^{2} \me^{\mi \pi k/3}}{3
(\varepsilon b)^{1/3}}(m \! + \! 1) \widehat{\nu}_{m+1}(k) \right. \nonumber \\
+&\left. \sum_{j=0}^{m-1}(j \! + \! 1) \widehat{\nu}_{j+1}(k)(\hat{\mu}_{m-j}^{\ast}(k) \! - \! 
2(m \! + \! 2 \! - \! j) \widehat{\nu}_{m+2-j}(k)) \right), \quad m \! \in \! \mathbb{Z}_{+},
\end{align} 
with
\begin{equation} \label{totsb5} 
\begin{gathered} 
\hat{\mu}_{0}^{\ast}(k) \! = \! -\frac{2a \me^{\mi \pi k/3}}{3(\varepsilon b)^{1/3}}, \quad 
\qquad \hat{\mu}_{1}^{\ast}(k) \! = \! 0, \\
\hat{\mu}_{m_{1}+2}^{\ast}(k) \! = \! -2 \! \left(\widehat{\mathrm{P}}_{m_{1}+2}^{\ast}(k) \! 
+ \! \hat{\mathfrak{w}}_{m_{1}+2}(k) \! + \! \sum_{j=0}^{m_{1}} \widehat{\mathrm{P}}_{j}^{
\ast}(k) \hat{\mathfrak{w}}_{m_{1}-j}(k) \right), \quad m_{1} \! \in \! \mathbb{Z}_{+},
\end{gathered}
\end{equation}
and
\begin{equation} \label{totsb6} 
\begin{gathered} 
\widehat{\mathrm{P}}_{0}^{\ast}(k) \! = \! \frac{2a \me^{\mi \pi k/3}}{3(\varepsilon b)^{1/3}}, 
\quad \qquad \widehat{\mathrm{P}}_{1}^{\ast}(k) \! = \! 0, \\
\widehat{\mathrm{P}}_{j}^{\ast}(k) \! = \! \frac{3}{2} \! \left(\hat{\mathfrak{u}}_{j}(k) 
\! - \! \mi (-1)^{\hat{\varepsilon}_{2}} \me^{\mi 2 \pi k/3}(\varepsilon b)^{1/3} \! \left(
\hat{\mathfrak{r}}_{j+2}(k) \! - \! 2 \hat{\mathfrak{u}}_{j+2}(k) \! + \! \sum_{m_{2}=0}^{j} 
\hat{\mathfrak{u}}_{m_{2}}(k) \hat{\mathfrak{r}}_{j-m_{2}}(k) \right) \! \right), \quad 
\mathbb{N} \! \ni \! j \! \geqslant \! 2,
\end{gathered}
\end{equation}
where the expansion coefficients $\hat{\mathfrak{u}}_{m}(k)$ and $\hat{\mathfrak{w}}_{m}(k)$ (resp., 
$\hat{\mathfrak{r}}_{m}(k))$, $m \! \in \! \mathbb{Z}_{+}$, $k \! \in \! \lbrace \pm 1 \rbrace$, are given in 
equations \eqref{appen2}--\eqref{appen8} (resp., \eqref{appen12} and \eqref{appen13}$)$.
\end{ddddd}

\emph{Proof}. Applying the (map) isomorphism (cf. Appendix \ref{sectonsymm}) 
$\hat{\mathscr{F}}^{\lbrace \hat{\ell} \rbrace}_{\hat{\varepsilon}_{1},\hat{\varepsilon}_{2},\hat{m}(\hat{\varepsilon}_{2})}$, 
whose action on $\mathscr{M}$ is given by equations \eqref{laxhat91}--\eqref{laxhat98} and 
\eqref{laxhat114}--\eqref{laxhat121}, to the $(\varepsilon_{1},\varepsilon_{2},m(\varepsilon_{2}) \vert \ell) 
\! = \! (0,0,0 \vert 0)$ asymptotics \eqref{totsa27} (as $\tau \! \to \! +\infty$ with $\varepsilon b \! > \! 0$) 
for $\hat{\varphi}(\tau)$, one arrives at the trans-series asymptotics \eqref{totsb1} (and equations 
\eqref{totsb2}--\eqref{totsb6}) stated in the theorem. \hfill $\qed$
\appendix
\setcounter{section}{5} 
\section{Appendix: Literature Survey of the DP3E} \label{litsurvdp3e} 
The interested reader will find representative samples of the ubiquitous manifestations of the DP3E \eqref{eq1.1} in this appendix.
\begin{enumerate}
\item[\textbf{(i)}] It was shown in \cite{sulei} that a variant of the DP3E \eqref{eq1.1} appears in the characterisation 
of the effect of the small dispersion on the self-focusing of solutions of the fundamental equations of non-linear optics 
in the one-dimensional case, where the main order of the influence of this effect is described via a universal special 
monodromic solution of the non-linear Schr\"{o}dinger equation (NLSE); in particular, the author studies the asymptotics 
of a function that can be identified as a solution (the so-called `Suleimanov solution') of a slightly modified, yet equivalent, 
version of the DP3E \eqref{eq1.1} for the parameter values $a \! = \! \mi/2$ and $b \! = \! 64k^{-3}$, where $k \! > \! 0$ 
is a physical variable.
\item[\textbf{(ii)}] In \cite{avlkv}, an extensive number-theoretic and asymptotic analysis of the universal special 
monodromic solution considered in \cite{sulei} is presented: the author studies a particular meromorphic solution of 
the DP3E \eqref{eq1.1} that vanishes at the origin; more specifically, it is proved that, for $-2 \mi a \! \in \! \mathbb{Z}$, 
the aforementioned solution exists and is unique, and, for the case $a \! - \! \mi/2 \! \in \! \mathbb{Z}$, this solution 
exists and is unique provided that $u(\tau) \! = \! -u(-\tau)$. The bulk of the analysis presented in \cite{avlkv} focuses 
on the study of the Taylor-series expansion coefficients of the solution to the DP3E \eqref{eq1.1} that is holomorphic at 
$\tau \! = \! 0$; in particular, upon invoking the `normalisation condition' $b \! = \! a$ and taking $\varepsilon \! = \! +1$, 
it is shown that, for general values of the parameter $a$, these coefficients are rational functions of $a^{2}$ that possess 
remarkable number-theoretic properties: en route, novel notions such as super-generating functions and quasi-periodic 
fences are introduced. The author also studies the connection problem for the ``Suleimanov solution'' \cite{kitavzap2023} 
of the DP3E \eqref{eq1.1}.
\item[\textbf{(iii)}] Unlike the physical optics context adopted in \cite{sulei}, the authors of \cite{peetsdb} provide a colossal 
Riemann-Hilbert problem (RHP) asymptotic analysis of the solution of the focusing NLSE, $\mi \partial_{\mathrm{T}} 
\Psi \! + \! \tfrac{1}{2} \partial_{\mathrm{X}}^{2} \Psi \! + \! \lvert \Psi \rvert^{2} \Psi \! = \! 0$, by considering the rogue 
wave solution $\Psi (\mathrm{X},\mathrm{T})$ of infinite order, that is, a scaling limit of a sequence of particular solutions 
of the focusing NLSE modelling so-called rogue waves of ever-increasing amplitude, and show that, in the regime of 
large variables $\mathbb{R}^{2} \! \ni \! (\mathrm{X},\mathrm{T})$ when $\lvert \mathrm{X} \rvert \! \to \! +\infty$ in such 
a way that $\mathrm{T} \lvert \mathrm{X} \rvert^{-3/2} \! - \! 54^{-1/2} \! = \! \mathcal{O}(\lvert \mathrm{X} \rvert^{-1/3})$, 
the rogue wave of infinite order $\Psi (\mathrm{X},\mathrm{T})$ can be expressed explicitly in terms of a function 
$\mathcal{V}(y)$ extracted {}from the solution of the Jimbo-Miwa Painlev\'{e} II ($\mathrm{P} \mathrm{II}$) RHP for 
parameters $p \! = \! \ln (2)/2 \pi$ and $\tau \! = \! 1$;\footnote{Not to be confused with the independent variable $\tau$ 
that appears in the DP3E \eqref{eq1.1} and throughout this work.} in particular, Corollary 6 of \cite{peetsdb} presents 
the leading term of the $\mathrm{T} \! \to \! +\infty$ asymptotics of the rogue wave of infinite order $\Psi (0,\mathrm{T})$ 
(see, also, Theorem 2 and Section 4 of \cite{peetsdq}),\footnote{For the rogue wave of infinite order \cite{peetsdb}, one 
needs to consider asymptotics of tronqu\'{e}e/tritronqu\'{e}e solutions of the inhomogeneous $\mathrm{P} \mathrm{II}$ 
equation, $\tfrac{\md^{2}u(x;\alpha)}{\md x^{2}} \! = \! 2(u(x;\alpha))^{3} \! + \! xu(x;\alpha) \! - \! \alpha$, for the special 
complex value of $\alpha \! = \! \tfrac{1}{2} \! + \! \mi \tfrac{\ln (2)}{2 \pi}$ (asymptotics for tronqu\'{e}e/tritronqu\'{e}e solutions 
of the $\mathrm{P} \mathrm{II}$ equation with $\alpha \! = \! 0$ are given in the monograph \cite{a5}), and to know that the 
increasing tritronqu\'{e}e solution, denoted $u_{\mathrm{TT}}^{-}(x;\alpha)$ in \cite{peetsdc}, is void of poles on $\mathbb{R}$; 
furthermore, for the function $\mathcal{V}(y)$ to have sense as a meaningful asymptotic representation of the rogue wave 
of infinite order $\Psi (\mathrm{X},\mathrm{T})$, it is, additionally, necessary that $u_{\mathrm{TT}}^{-}(x;\alpha)$ be a 
global solution (analytic $\forall$ $x \! \in \! \mathbb{R}$) of the $\mathrm{P} \mathrm{II}$ equation for $\alpha \! = \! 
\tfrac{1}{2} \! + \! \mi \tfrac{\ln (2)}{2 \pi}$. In \cite{peetsdc}, the author provides a complete RHP asymptotic analysis of 
the global nature of tritronqu\'{e}e solutions of the $\mathrm{P} \mathrm{II}$ equation for various complex values of 
$\alpha$, including the particular value $\alpha \! = \! \tfrac{1}{2} \! + \! \mi \tfrac{\ln (2)}{2 \pi}$, and relates the function 
$\mathcal{V}(y)$ to the $\mathrm{P} \mathrm{II}$ equation, subsequently identifying the particular solution that is requisite 
in order to construct $\mathcal{V}(y)$ as the increasing tritronqu\'{e}e solution $u_{\mathrm{TT}}^{-}(x;\alpha)$ for the special 
parameter value $\alpha \! = \! \tfrac{1}{2} \! + \! \mi \tfrac{\ln (2)}{2 \pi}$; moreover, the value of the total, regularised integral 
over $\mathbb{R}$ for the increasing tritronqu\'{e}e solution is evaluated.} which, in the context of the DP3E \eqref{eq1.1}, 
coincides, up to a scalar, $\tau$-independent factor, with $\exp (\mi \hat{\varphi}(\tau))$, $\mathrm{T} \! = \! \tau^{2}$, where, 
given the solution, denoted by $\hat{u}(\tau)$, say, of the DP3E \eqref{eq1.1} studied in \cite{avlkv} for the monodromy data 
corresponding to $a \! = \! \mi/2$ (and a suitable choice for the parameter $b$), $\hat{\varphi}(\tau)$ is the general solution 
of the ODE $\hat{\varphi}^{\prime}(\tau) \! = \! 2a \tau^{-1} \! + \! b(\hat{u}(\tau))^{-1}$. 
\item[\textbf{(iv)}] The authors of \cite{BuckinghamMiller2022} present an expansive study of algebraic solutions (rational 
functions of $\tau^{1/3}$) of the DP3E \eqref{eq1.1} for the parameter values $\varepsilon \! = \! -1$, $b \! = \! \mi$, and 
$a \! = \! -\mi n$, $n \! \in \! \mathbb{Z}$. By considering the Lax-pair equations associated with the DP3E \eqref{eq1.1}, 
the authors construct their simultaneous solutions, called the `seed' lax-pair solutions, corresponding to the simplest 
algebraic solution of the DP3E \eqref{eq1.1}, $u(\tau) \! := \! u_{0}(\tau) \! = \! \tfrac{1}{2} \tau^{1/3}$, for $\varepsilon \! 
= \! -1$, $b \! = \! \mi$, and $a \! = \! 0$ in terms of Airy functions, and then formulate, as Riemann-Hilbert Problem 1 
(RHP1), the inverse monodromy problem for the rational solution $u(\tau) \! := \! u_{n}(\tau)$ for $a \! = \! -\mi n$, 
$n \! \in \! \mathbb{Z} \setminus \lbrace 0 \rbrace$ (the case $a \! = \! -\mi n$ for $n \! = \! 0$ is solved via the `seed' Lax-pair 
solutions); in particular, the authors show that, if RHP1 is solvable for $\tau \! > \! 0$ and $n \! \in \! \mathbb{Z}$, then the 
function $u_{n}(\tau)$ defined by equation~(101) in \cite{BuckinghamMiller2022} is the unique solution of the DP3E 
\eqref{eq1.1} with $\varepsilon \! = \! -1$, $b \! = \! \mi$, and $a \! = \! -\mi n$, $n \! \in \! \mathbb{Z}$, that is a rational 
function of $\tau^{1/3}$ (see Theorem~1 of \cite{BuckinghamMiller2022}). The authors then use the RHP1 representation 
for the algebraic solution $u_{n}(\tau)$ of the DP3E \eqref{eq1.1} to consider the large-positive-$n$ asymptotic behaviour 
of the solution (as a consequence of an inherent symmetry of the DP3E \eqref{eq1.1} that is discussed at the beginning 
of Subsection 4.1 of \cite{BuckinghamMiller2022}, it is sufficient to consider large $n \! \in \! \mathbb{N}$); in particular, 
after a rescaling argument for both the independent variable and the spectral parameter, the authors present a rigorous 
asymptotic analysis of RHP1 and derive $\mathbb{N} \! \ni \! n \! \to \! \infty$ (for sufficiently large rescaled $\tau \! > \! 0$) 
asymptotics of the function $u_{n}(\tau)$ (see, in particular, Theorems 2 and 3 of \cite{BuckinghamMiller2022}). (In this 
context, see, also, \cite{buckmill24}.)
\item[\textbf{(v)}] Introducing the substitution $\varepsilon \tau u \! = \! (x/3)^{2}y$, $\varepsilon b \tau^{2} \! = \! 2(x/3)^{3}$, 
the author of \cite{shimdeepee3} transforms the DP3E \eqref{eq1.1} into the second-order non-linear ODE $y^{\prime 
\prime}(x) \! = \! \tfrac{(y^{\prime}(x))^{2}}{y(x)} \! - \! \tfrac{y^{\prime}(x)}{x} \! - \! 2(y(x))^{2} \! + \! \tfrac{3a}{x} \! + \! 
\tfrac{1}{y(x)}$, where the prime denotes differentiation with respect to $x$, and then, via additional auxiliary changes of 
variables, shows that, with $x \! = \! t \me^{\mi \phi}$, the latter ODE for $y$ governs the isomonodromy deformation of 
a $2 \! \times \! 2$ linear system $\partial_{\lambda} \Psi (\lambda,t) \! = \! \tfrac{t}{3} \mathcal{B}(\lambda,t) \Psi (\lambda,t)$, 
where $\mathrm{M}_{2}(\mathbb{C}) \! \ni \! \mathcal{B}(\lambda,t)$ is given in equation (1.4), or, equivalently, equation 
(3.2), of \cite{shimdeepee3}. By applying the isomonodromy deformation method \cite{a2}, the author demonstrates the 
Boutroux ansatz (near the point at infinity) by deriving an elliptic asymptotic representation of the general solution $y(x)$ 
in terms of the Weierstrass $\wp$-function as $x \! = \! t \me^{\mi \phi} \! \to \! \infty$ in cheese-like strip domains along 
generic directions; see, in particular, the leading-order asymptotics of $y(x)$ stated in Theorems 2.1--2.3 of \cite{shimdeepee3}.
\item[\textbf{(vi)}] In \cite{zamol3}, the authors study the eigenvalue correlation kernel, denoted by $K_{n}(x,y,t)$, for the singularly 
perturbed Laguerre unitary ensemble (pLUE)\footnote{The pLUE and its relation to the Painlev\'{e} III ($\mathrm{P} \mathrm{III}$) 
equation was introduced and studied in \cite{schtin}.} on the space $\mathscr{H}_{n}^{+}$ of $n \times n$ positive-definite 
Hermitian matrices $M \! = \! (M)_{i,j=1}^{n}$ defined by the probability measure $Z_{n}^{-1}(\det M)^{\alpha} \exp 
(-\tr V_{t}(M)) \, \md M$, $n \! \in \! \mathbb{N}$, $\alpha \! > \! 0$, $t \! > \! 0$, where $Z_{n} \! := \! \int_{\mathscr{H}_{n}^{+}}
(\det M)^{\alpha} \me^{-\tr V_{t}(M)} \, \md M$ is the normalisation constant, $\md M \! := \! \prod_{i=1}^{n} \md M_{ii} 
\prod_{j=1}^{n-1} \prod_{k=j+1}^{n} \md \Re (M_{jk}) \, \md \Im (M_{jk})$, and $V_{t}(x) \! := \! x \! + \! t/x$, $x \! \in \! (0,+\infty)$. 
By considering, for example, a variety of double-scaling limits such as $n \! \to \! \infty$ and $(0,d] \! \ni \! t \! \to \! 0^{+}$, 
$d \! > \! 0$, such that $s \! := \! 2nt$ belongs to compact subsets of $(0,+\infty)$, or $n \! \to \! \infty$ and $t \! \to \! 0^{+}$ 
such that $s \! \to \! 0^{+}$, or $n \! \to \! \infty$ and $(0,d] \! \ni \! t$ such that $s \! \to \! +\infty$, the authors derive the 
corresponding limiting behaviours of the eigenvalue correlation kernel by studying the large-$n$ asymptotics of the 
orthogonal polynomials associated with the singularly perturbed Laguerre weight $w(x;t,\alpha) \! = \! x^{\alpha} 
\me^{-V_{t}(x)}$, and, en route, demonstrate that some of the limiting kernels involve certain functions related to a special 
solution of $(P_{\mathrm{III}^{\prime}})_{D_{7}}$ \eqref{dpee3d7}; moreover, in the follow-up work \cite{zamol4} on the 
pLUE, the authors derive the large-$n$ asymptotic formula (uniformly valid for $(0,d] \! \ni \! t$, $d \! > \! 0$ and fixed) for 
the Hankel determinant $D_{n}[w;t] \! := \! \det (\smallint_{0}^{+\infty} x^{j+k}w(x;t,\alpha) \, \md x)_{j,k=0}^{n-1}$ associated 
with the singularly perturbed Laguerre weight $w(x;t,\alpha)$, and show that the asymptotic representation for $D_{n}[w;t]$ 
involves a function related to a particular solution of $(P_{\mathrm{III}^{\prime}})_{D_{7}}$ \eqref{dpee3d7}. In the study of 
the Hankel determinant $D_{n}(t,\alpha,\beta) \! := \! \det (\smallint_{0}^{1} \xi^{j+k}w(\xi;t,\alpha,\beta) \, \md \xi)_{j,k=0}^{n-1}$ 
generated by the Pollaczek-Jacobi-type weight $w(x;t,\alpha,\beta) \! = \! x^{\alpha}(1 \! - \! x)^{\beta} \me^{-t/x}$, $x \! \in \! [0,1]$, 
$t \! \geqslant \! 0$, $\alpha,\beta \! > \! 0$, which is a fundamental object in unitary random matrix theory, under a double-scaling 
limit where $n$, the dimension of the Hankel matrix, tends to $\infty$ and $t \! \to \! 0^{+}$ in such a way that $s \! := \! 2n^{2}t$ 
remains bounded, the authors of \cite{mchychefa} show that the double-scaled Hankel determinant has an integral 
representation in terms of particular asymptotic solutions of a scaled version of the DP3E \eqref{eq1.1} (or, equivalently, 
$(P_{\mathrm{III}^{\prime}})_{D_{7}}$ \eqref{dpee3d7}). In \cite{zamol6}, the authors study singularly perturbed unitary 
invariant random matrix ensembles on $\mathscr{H}_{n}^{+}$ defined by the probability measure $C_{n}^{-1}(\det M)^{\alpha} 
\exp (-n \tr V_{k}(M)) \, \md M$, $n,k \! \in \! \mathbb{N}$, $\alpha \! > \! -1$, where $C_{n} \! := \! \int_{\mathscr{H}_{n}^{+}}
(\det M)^{\alpha} \me^{-n\tr V_{k}(M)} \, \md M$, and the---perturbed---potential $V_{k}(x)$ has a pole of order $k$ at the 
origin, $V_{k}(x) \! := \! V(x) \! + \! (t/x)^{k}$, $t \! > \! 0$, with the regular part, $V$, of the potential being real analytic on 
$[0,+\infty)$ and satisfying certain constraints; in particular, for the pLUE, the authors obtain, in various double-scaling 
limits when the size of the matrix $n \! \to \! \infty$ (at an appropriately adjusted rate) and the ``strength'' of the perturbation 
$t \! \to \! 0$, asymptotics of the associated eigenvalue correlation kernel and partition function, which are characterised 
in terms of special, pole-free solutions of a hierarchy (indexed by $k$) of higher-order analogues of the $\mathrm{P} 
\mathrm{III}$ equation: the first ($k \! = \! 1$) member of this $\mathrm{P} \mathrm{III}$ hierarchy, denoted by $\ell_{1}(s)$, 
$s \! > \! 0$, solves a rescaled version of the DP3E \eqref{eq1.1}. (Analogous results for the singularly perturbed Gaussian 
unitary ensemble (pGUE) on the set $\mathscr{H}_{n}$ of $n \times n$ Hermitian matrices are also obtained in \cite{zamol6}.) 
For the pLUE with perturbed potential $V_{k}(x) \! := \! V(x) \! + \! (t/x)^{k}$, $k \! \in \! \mathbb{N}$, $x \! \in \! (0,+\infty)$, 
$t \! > \! 0$, studied in \cite{zamol6}, the authors of \cite{zamol5} consider a related Fredholm determinant of an integral 
operator, denoted by $\mathcal{K}_{\mathrm{P} \mathrm{III}}$, acting on the space $L^{2}((0,+\infty))$, whose kernel is 
constructed {}from a certain $\mathrm{M}_{2}(\mathbb{C})$-valued function associated with a hierarchy (indexed by $k$) 
of higher-order analogues of the $\mathrm{P} \mathrm{III}$ equation; more precisely, for the Fredholm determinant 
$F(s;\lambda) \! := \! \ln \det (\mathrm{I} \! - \! \mathcal{K}_{\mathrm{P} \mathrm{III}})$, $s,\lambda \! > \! 0$, the authors 
of \cite{zamol5} obtain $s \! \to \! +\infty$ asymptotics of $F(s;\lambda)$ characterised in terms of an explicit integral 
representation of a special, pole-free solution for the first $(k \! = \! 1$) member of the corresponding $\mathrm{P} \mathrm{III}$ 
hierarchy: this solution is denoted by $\ell_{1}(\lambda)$, and it solves a rescaled version of the DP3E \eqref{eq1.1}.
\item[\textbf{(vii)}] In \cite{tradom}, the authors compute small-$t$ asymptotics of a class of solutions to the two-dimensional 
cylindrical Toda equations (2DCTE),\footnote{See, also, its generalisations \cite{gutsin1,gutsin2,gutsin3,gutsin4}.} 
$q_{k}^{\prime \prime}(t) \! + \! t^{-1}q_{k}^{\prime}(t) \! = \! 4(\me^{q_{k}(t)-q_{k-1}(t)} \! - \! \me^{q_{k+1}(t)-q_{k}(t)})$, 
$k \! \in \! \mathbb{Z}$, satisfying the periodicity conditions $q_{k+n}(t) \! = \! q_{k}(t)$, where the integer $n$ is arbitrary 
but fixed. Solutions that are valid for all $t \! > \! 0$ have the representation $q_{k}(t) \! = \! \log \det (\mathrm{I} \! - \! 
\lambda \mathcal{K}_{k}) \! - \! \log \det (\mathrm{I} \! - \! \lambda \mathcal{K}_{k-1})$, where $\mathcal{K}_{k}$ is the 
integral operator on $\mathbb{R}_{+}$ with kernel $\sum_{\lbrace \omega^{n}=1 \rbrace \setminus \lbrace 1 \rbrace} 
\omega^{k}c_{\omega} \tfrac{\me^{-t((1-\omega)u+(1-\omega^{-1})u^{-1})}}{-\omega u+v}$, for some coefficients 
$c_{\omega}$, and $\lambda$ is a free parameter. For $n \! = \! 3$ and the imposition of an additional constraint, which 
implies $q_{1}(t) \! = \! 0$ and $q_{2}(t) \! = \! -q_{3}(t)$, the 2DCTE gives rise to the radial Bullough-Dodd equation (for 
$q_{3}(t)$), $q_{3}^{\prime \prime}(t) \! + \! t^{-1}q_{3}^{\prime}(t) \! = \! 4(\me^{2q_{3}(t)} \! - \! \me^{-q_{3}(t)})$, which, 
via the dependent-variable transformation $w(t) \! = \! \me^{-q_{3}(t)}$, reduces to the non-linear ODE $w^{\prime \prime}(t) 
\! = \! \tfrac{(w^{\prime}(t))^{2}}{w(t)} \! - \! \tfrac{w^{\prime}(t)}{t} \! + \! 4(w(t))^{2} \! - \! \tfrac{4}{w(t)}$; by making one more 
change of variables, namely, $t \! = \! \lambda^{2/3}$ and $w(t) \! = \! \lambda^{-1/3} \mathcal{W}(\lambda)$, this ODE can, 
in turn, be transformed to the $\mathrm{P} \mathrm{III}$ equation with parameter values $(16/9,0,0,-16/9)$,
\begin{equation*} 
\mathcal{W}^{\prime \prime}(\lambda) \! = \! \frac{(\mathcal{W}^{\prime}(\lambda))^{2}}{\mathcal{W}
(\lambda)} \! - \! \frac{\mathcal{W}^{\prime}(\lambda)}{\lambda} \! + \! \frac{16}{9} \frac{(\mathcal{W}
(\lambda))^{2}}{\lambda} \! - \! \frac{16}{9} \frac{1}{\mathcal{W}(\lambda)},
\end{equation*}
where the prime denotes differentiation with respect to $\lambda$, which can be identified as a special reduction of the 
DP3E \eqref{eq1.1} for $a \! = \! 0$. The small-$t$ asymptotics of $q_{k}(t)$ are derived by computing the asymptotics 
$\det (\mathrm{I} \! - \! \lambda \mathcal{K}_{k}) \genfrac{}{}{0pt}{3}{\thicksim}{t \to 0^{+}} \! b_{k}(t/n)^{a_{k}}$, 
$n \! = \! 2,3$, where explicit expressions for the coefficients $a_{k}$ and $b_{k}$ are presented in \cite{tradom}.
\item[\textbf{(viii)}] The DP3E \eqref{eq1.1} also plays a prominent r\^{o}le in the description of surfaces with constant 
negative Gaussian curvature ($K$-surfaces) and two straight asymptotic lines (\emph{Amsler surfaces\/}) \cite{peetsde}. 
A non-degenerate surface in $\mathbb{R}^{3}$ is called an \emph{affine sphere\/} if all affine normal directions intersect 
at a point: this class of surfaces is described by an integrable equation first derived by Tzitz\'{e}ica. As discussed in 
\cite{peetsde}, for affine spheres characterised by the property that they possess two intersecting straight affine lines, 
the corresponding Tzitz\'{e}ica equation reduces to the $\mathrm{P} \mathrm{III}$ equation with parameter values 
$(1,0,0,-1)$,
\begin{equation*} 
y^{\prime \prime}(t) \! = \! \frac{(y^{\prime}(t))^{2}}{y(t)} \! - \! \dfrac{y^{\prime}(t)}{t} \! + \! \frac{(y(t))^{2}}{t} \! - \! 
\dfrac{1}{y(t)},
\end{equation*}
where the prime denotes differentiation with respect to $t$, with $y(t) \! = \! t^{1/3}H(r)$ and $t \! = \! \frac{8}{3^{3/2}}
r^{3/4}$, and where $H(r)$, with $r \! := \! xy$, is a Lorentz invariant solution of the Tzitz\'{e}ica equation that satisfies 
the second-order non-linear ODE
\begin{equation*} 
H^{\prime \prime}(r) \! = \! \frac{(H^{\prime}(r))^{2}}{H(r)} \! - \! \frac{H^{\prime}(r)}{r} \! + \! \frac{1}{r} 
\! \left((H(r))^{2} \! - \! \frac{1}{H(r)} \right),
\end{equation*}
where the prime denotes differentiation with respect to $r$; in fact, the ODE for the function $y(t)$ can be identified 
as a special reduction of the DP3E \eqref{eq1.1} for $a \! = \! 0$: letting $\tau \! = \! 2^{-3/2} \me^{\mi (2m+1) \pi/4}t$ 
and $u(\tau) \! = \! -2^{-3/2} \me^{-\mi (2m+1) \pi/4}y(t)$, $m \! = \! 0,1,2,3$, and choosing the---external---parameter 
values $\varepsilon \! = \! b \! = \! +1$ and $a \! = \! 0$, it follows that the DP3E \eqref{eq1.1} reduces to the ODE 
for $y(t)$. The algebroid theory for solutions of the ODE for $H(r)$ is presented in \cite{avkavalgeb}.
\item[\textbf{(ix)}] Let $\mathcal{X}$ be a six-dimensional Calabi-Yau (CY) manifold (a complex K\"{a}hler 
three-fold with covariantly constant holomorphic three-form $\Omega$). The Strominger-Yau-Zaslow (SYZ) 
conjecture (see \cite{DP} for details) states that, near the large complex structure limit, both $\mathcal{X}$ 
and its mirror should be the fibrations over the moduli space of special Lagrangian tori (submanifolds admitting 
a unitary flat connection). As an examination of the SYZ conjecture, Loftin-Yau-Zaslow (LYZ) set out to prove 
the existence of the metric of Hessian form $g_{B} \! = \! \tfrac{\partial^{2} \phi}{\partial x^{j} \partial x^{k}} \, 
\md x^{j} \otimes \md x^{k}$, where $x^{j}$, $j \! = \! 1,2,3$, are local coordinates on a real three-dimensional 
manifold, and $\phi$ (a K\"{a}hler potential) is homogeneous of degree two in $x^{j}$ and satisfies the real 
Monge-Amp\'{e}re equation $\det \! \left(\tfrac{\partial^{2} \phi}{\partial x^{j} \partial x^{k}} \right) \! = \! 1$: 
LYZ showed that the construction of the metric is tantamount to searching for solutions of the definite affine 
sphere equation (DASE) $\psi_{z \overline{z}} \! + \! \tfrac{1}{2} \me^{\psi} \! + \! \lvert U \rvert^{2} \me^{-2 \psi} 
\! = \! 0$, $U_{\overline{z}} \! = \! 0$, where $\psi$ and $U$ are real- and complex-valued functions, respectively, 
on an open subset of $\mathbb{C}$. For $U \! = \! z^{-2}$, LYZ proved the existence of the radially symmetric 
solution $\psi$ of the DASE with a prescribed behaviour near the singularity $z \! = \! 0$, and established the 
existence of the global solution to the coordinate-independent version of the DASE on $\mathbb{S}^{2}$ with 
three points excised. In \cite{DP}, the authors show that the DASE, and a closely related equation called the 
Tzitz\'{e}ica equation, arise as reductions of anti-self-dual Yang-Mills (ASDYM) system by two translations; 
moroever, they show that the ODE characterising its radial solutions give rise to an isomonodromy problem 
described by the $\mathrm{P} \mathrm{III}$ equation for special values of its parameters. In particular (see 
Proposition~1.3 of \cite{DP}), the authors show that, for $U \! = \! z^{-2}$, solutions of the DASE that are invariant 
under the group of rotations (rotational symmetry) $z \! \to \! \me^{\mi \mathfrak{c}}z$, $\mathfrak{c} \! \in \! 
\mathbb{R}$, are of the form $\psi (z,\overline{z}) \! = \! \ln (\mathscr{H}(s)) \! - \! 3 \ln (s)$, with $s \! := \! \lvert 
z \rvert^{1/2}$, where $\mathscr{H}(s)$ solves the $\mathrm{P} \mathrm{III}$ equation with parameter values 
$(-8,0,0,-16)$,
\begin{equation*} 
\mathscr{H}^{\prime \prime}(s) \! = \! \frac{(\mathscr{H}^{\prime}(s))^{2}}{\mathscr{H}(s)} \! - \! 
\frac{\mathscr{H}^{\prime}(s)}{s} \! - \! \frac{8(\mathscr{H}(s))^{2}}{s} \! - \! \frac{16}{\mathscr{H}(s)},
\end{equation*}
where the prime denotes differentiation with respect to $s$, which can be identified as a special reduction of the 
DP3E \eqref{eq1.1} for $a \! = \! 0$. The authors of \cite{DP} demonstrate that the existence theorem for Hessian 
metrics with prescribed monodromy reduces to the study of the $\mathrm{P} \mathrm{III}$ equation with parameters 
$(-8,0,0,-16)$, that is, a class of semi-flat CY metrics is obtained in terms of real solutions of the DP3E \eqref{eq1.1} 
for $a \! = \! 0$ (see, also, \cite{contdor15,contdor17,Dun2012,Dun2023}).
\item[\textbf{(x)}] In \cite{peetsdd}, the author introduces affine spheres as immersions of a manifold $\mathcal{M}$ 
as a hypersurface in $\mathbb{R}^{n}$ with certain properties and defines the affine metric $h$ and the cubic form 
$C$ on $\mathcal{M}$. By identifying, for $3$-dimensional cones and, correspondingly, affine $2$-spheres, the 
manifold $\mathcal{M}$ with a non-compact, simply-connected domain in $\mathbb{C}$, one can introduce complex 
isothermal co-ordinates $z$ on $\mathcal{M}$, in terms of which the affine metric $h$ may equivalently be described 
by a real conformal factor $u(z)$ and the cubic form $C$ by a holomorphic function $U(z)$ on $\mathcal{M}$, the 
relations being $h \! = \! \me^{u} \lvert \md z \rvert^{2}$ and $C \! = \! 2 \Re (U(z)) \md z^{3}$: the compatibility 
condition of the pair $(u,U)$ is referred to as \emph{Wang's equation}, $\me^{u} \! = \! \tfrac{1}{2} \Delta u \! + \! 
2 \lvert U \rvert^{2} \me^{-2u}$, where $\Delta u \! = \! u_{xx} \! + \! u_{yy} \! = \! 4u_{z \overline{z}}$ is the Laplacian 
of $u$, $\partial_{z} \! := \! \tfrac{1}{2}(\partial_{x} \! - \! \mi \partial_{y})$, and $\partial_{\overline{z}} \! := \! 
\tfrac{1}{2}(\partial_{x} \! + \! \mi \partial_{y})$. By classifying pairs $(\psi,U)$, where $\psi$ is a vector field on 
$\mathcal{M}$ generating a one-parameter group of conformal automorphisms on $\mathcal{M}$ which multiply 
$U$ by unimodular complex constants, the author finds, for every pair $(\psi,U)$, a unique solution $u$ of Wang's 
equation such that the corresponding affine metric $h$ is complete on $\mathcal{M}$ and $\psi$ is a Killing vector 
field for $h$: this latter property permits Wang's equation to be reduced to a second-order non-linear ODE that is 
equivalent to the DP3E \eqref{eq1.1}, a detailed qualitative study for which is presented in \cite{peetsdd}. The 
author presents a complete classification of self-associated cones (one calls a cone self-associated if it is linearly 
isomorphic to all its associated cones, with two cones said to be associated with each other if the Blaschke metrics 
on the corresponding affine spheres are related by an orientation-preserving isometry) and computes isothermal 
parametrisations of the corresponding affine spheres, the solution(s) of which can be expressed in terms of 
degenerate $\mathrm{P} \mathrm{III}$ transcendents (solutions of the DP3E \eqref{eq1.1}).
\end{enumerate}
Whilst not directly relevant to the DP3E \eqref{eq1.1}, the following facts are worth mentioning: (1) elliptic asymptotic 
representations in terms of the Jacobi $\mathrm{sn}$-function in cheese-like strip domains along generic directions are obtained 
for the general solution of the `complete' $\mathrm{P} \mathrm{III}$ equation in \cite{shimcomplp3}; (2) a detailed study of the 
$\mathrm{P} \mathrm{III}$ monodromy maps under the $D_{6} \! \to \! D_{8}$ confluence has recently been presented in 
\cite{barovyyerov}; (3) parametric Stokes phenomena for the $D_{6}$ and $D_{7}$ cases of the $\mathrm{P} \mathrm{III}$ 
equation are studied in \cite{iwak1}; (4) application of the $\mathrm{P} \mathrm{III}$ equation to the study of transformation 
phenomena for parametric Painlev\'{e} equations for the $D_{6}$ and $D_{7}$ cases is considered in \cite{iwak2}, whilst the 
$D_{8}$ case is studied in \cite{ayota,hwakytak}; (5) the monograph \cite{magch} studies the relation of the $\mathrm{P} 
\mathrm{III}$ equation of type $(P_{\mathrm{III}})_{D_{6}}$ to isomonodromic families of vector bundles on $\mathbb{P}^{1}$ 
with meromorphic connections; (6) in \cite{gavyy}, the $\pmb{\pmb{\tau}}$-function associated with the degenerate $\mathrm{P} 
\mathrm{III}$ equation of type $D_{8}$ is shown to admit a Fredholm determinant representation in terms of a generalised 
Bessel kernel; and (7) by using the universal example of the Gross-Witten-Wadia (GWW) third-order phase transition in the 
unitary matrix model, concomitant with the explicit Tracy-Widom mapping of the GWW partition function to a solution of a 
$\mathrm{P} \mathrm{III}$ equation, the transmutation (change in the resurgent asymptotic properties) of a trans-series in 
two parameters (a coupling $g^{2}$ and a gauge index $N$) at all coupling and all finite $N$ is studied in \cite{ahmedunne} 
(see, also, \cite{eund}).

\vspace*{0.35cm}
\noindent
\textbf{\Large Acknowledgements} 

\vspace*{0.25cm}
\noindent 
The author is grateful to A. V. Kitaev for perspicacious comments and criticisms related to preliminary results of this work, and to the St. 
Petersburg Branch of the Steklov Mathematical Institute (POMI) for hospitality. The author is grateful to the Referees for a cornucopia 
of constructive recommendations which were instrumental in improving the content of the paper. \emph{The author also wishes to honour 
the memories of L. D. Faddeev and P. P. Kulish\/}.

\end{document}